\documentclass[nospthms]{svmono}
\usepackage[latin1]{inputenc}
\usepackage[OT2,T1]{fontenc}
\usepackage[main=francais,english]{babel}
\usepackage{lmodern}
\usepackage{amsmath}
\usepackage{amssymb}
\usepackage{mathrsfs}
\usepackage[ps,dvips,arrow,matrix,tips]{xy}
\usepackage{makeidx}
\usepackage{multicol}
\makeatletter
\usepackage[paperwidth=14.8cm,paperheight=22cm,centering,textwidth=28pc,textheight=540\p@,marginparwidth=90\p@,headsep=12\p@]{geometry}
\makeatother

\entrymodifiers={+!!<0pt,\the\fontdimen22\textfont2>}
\SelectTips{cm}{10}
\SilentMatrices
\newcommand{\myxyhook}{
}
\newdir^{ (}{{}*!/-5pt/\dir^{(}}
\newdir_{ (}{{}*!/-5pt/\dir_{(}}

\newcount\tmpstorecat
\def\csarg#1#2{\expandafter#1\csname#2\endcsname}
\def\storecat#1{\tmpstorecat\escapechar \escapechar=-1\csarg\edef{restorecat\string#1}{\catcode`\string#1=\the\catcode\expandafter`\string#1}\catcode\expandafter`\string#1=12\relax\escapechar\tmpstorecat}
\def\restorecat#1{\tmpstorecat\escapechar \escapechar=-1\csname restorecat\string#1\endcsname\escapechar\tmpstorecat}
\def\myxyin{\storecat!\storecat;}
\def\myxyout{\restorecat!\restorecat;}

\newcounter{enonce}[chapter]
\renewcommand{\chaptermark}[1]{\markboth{\chaptername\ \thechapter.\ #1}{}}

\makeatletter
\def\@chapter[#1]#2{\ifnum \c@secnumdepth >\m@ne
                       \refstepcounter{chapter}%
                       \if@mainmatter
                         \typeout{\@chapapp\space\thechapter.}%
                       \fi
                    \addcontentsline{toc}{chapter}{\chaptername\ \thechapter.\ #1}%
                    \fi
                    \chaptermark{#1}%
                    \addtocontents{lof}{\protect\addvspace{10\p@}}%
                    \addtocontents{lot}{\protect\addvspace{10\p@}}%
                    \if@twocolumn
                      \@topnewpage[\@makechapterhead{#2}]%
                    \else
                      \@makechapterhead{#2}%
                      \@afterheading
                    \fi}
\def\@makechapterhead#1{{\parindent\z@\raggedright\normalfont
  \hyphenpenalty \@M
  \interlinepenalty\@M
  \chapsize\chapstyle
  \chap@hangfrom{\chaptername\ \thechapter\thechapterend\hskip\betweenumberspace}%!!!
  \ignorespaces#1\par\nobreak
  \processchapstarthook
  \ifdim\pagetotal>157\p@
     \vskip 11\p@
  \else
     \@tempdima=168\p@\advance\@tempdima by-\pagetotal
     \vskip\@tempdima
  \fi}}
\newcommand\chapternotoc{\startnewpage
                    \@ifundefined{thispagecropped}{}{\thispagecropped}
                    \thispagestyle{empty}%
                    \global\@topnum\z@
                    \@afterindentfalse
                    \secdef\@chapternotoc\@schapter}
\def\@chapternotoc[#1]#2{\ifnum \c@secnumdepth >\m@ne
                       \refstepcounter{chapter}%
                       \if@mainmatter
                         \typeout{\@chapapp\space\thechapter.}%
                       \fi
                    \fi
                    \chaptermark{#1}%
                    \addtocontents{lof}{\protect\addvspace{10\p@}}%
                    \addtocontents{lot}{\protect\addvspace{10\p@}}%
                    \if@twocolumn
                      \@topnewpage[\@makechapterhead{#2}]%
                    \else
                      \@makechapterhead{#2}%
                      \@afterheading
                    \fi}
\newcommand\sectionnotoc{\@startsectionnotoc{section}{1}{\z@}%
                       {-24\p@ \@plus -4\p@ \@minus -4\p@}%
                       {12\p@ \@plus 4\p@ \@minus 4\p@}%
                       {\normalfont\secsize\secstyle
                        \rightskip=\z@ \@plus 8em\pretolerance=10000 }}
\def\@startsectionnotoc#1#2#3#4#5#6{%
  \if@noskipsec \leavevmode \fi
  \par
  \@tempskipa #4\relax
  \@afterindenttrue
  \ifdim \@tempskipa <\z@
    \@tempskipa -\@tempskipa \@afterindentfalse
  \fi
  \if@nobreak
    \everypar{}%
  \else
    \addpenalty\@secpenalty\addvspace\@tempskipa
  \fi
  \@ifstar
    {\@ssect{#3}{#4}{#5}{#6}}%
    {\@dblarg{\@sectnotoc{#1}{#2}{#3}{#4}{#5}{#6}}}}
\def\@sectnotoc#1#2#3#4#5#6[#7]#8{%
  \ifnum #2>\c@secnumdepth
    \let\@svsec\@empty
  \else
    \refstepcounter{#1}%
    \protected@edef\@svsec{\@seccntformat{#1}\relax}%
  \fi
  \@tempskipa #5\relax
  \ifdim \@tempskipa>\z@
    \begingroup
      #6{%
        \@hangfrom{\hskip #3\relax\@svsec}%
          \interlinepenalty \@M #8\@@par}%
    \endgroup
    \csname #1mark\endcsname{#7}%
  \else
    \def\@svsechd{%
      #6{\hskip #3\relax
      \@svsec #8}%
      \csname #1mark\endcsname{#7}}%
  \fi
  \@xsect{#5}}
\let\if@onecolind\iffalse
\def\onecolindex{\let\if@onecolind\iftrue}
\renewenvironment{theindex}
               {\if@twocolumn
                  \@restonecolfalse
                \else
                  \@restonecoltrue
                \fi
                \columnseprule \z@
                \columnsep 1cc
                \@nobreaktrue
                \if@onecolind
                   \chapter*{\indexname}%
                \else
                   \begin{multicols}{2}[\chapter*{\indexname}%
                \fi
                {\csname indexst@rthook\endcsname}%
                \if@onecolind\else]\fi
                \global\let\indexst@rthook=\undefined
                \markboth{\indexname}{\indexname}%
                \addcontentsline{toc}{chapter}{\indexname}%
                \flushbottom
                \parindent\z@
                \rightskip\z@ \@plus 40\p@
                \parskip\z@ \@plus .3\p@\relax
                \flushbottom
                \let\item\@idxitem
                \def\,{\relax\ifmmode\mskip\thinmuskip
                             \else\hskip0.2em\ignorespaces\fi}%
                \normalfont\small}
               {\if@onecolind\else\end{multicols}\fi
                \global\let\if@threecolind\iffalse
                \if@restonecol\onecolumn\else\clearpage\fi}
\def\theenonce{\thechapter.\@arabic\c@enonce}
\xdef\longdashrightarrow{\mathrel{\dabar@\dabar@\dabar@
                              \mathchar"0\hexnumber@\symAMSa 4B}}%
\def\myrightarrow{{\setbox\z@\hbox{$\rightarrow$}\dimen0\ht\z@\multiply\dimen0 6\divide\dimen0 10\ht\z@\dimen0\box\z@}}
\def\myrightarrowfill@{\arrowfill@\relbar\relbar\myrightarrow}
\newcommand{\myxrightarrow}[2][]{\ext@arrow 0359\myrightarrowfill@{#1}{#2}}
\makeatother
{\setbox0\hbox{$ $}}\fontdimen16\textfont2=\fontdimen17\textfont2
\lineskip=\lineskiplimit

\renewcommand{\Gamma}{\varGamma}
\renewcommand{\Delta}{\varDelta}
\renewcommand{\Theta}{\varTheta}
\renewcommand{\Lambda}{\varLambda}
\renewcommand{\Phi}{\varPhi}
\renewcommand{\Pi}{\varPi}
\renewcommand{\Omega}{\varOmega}
\mathcode`A="7041 \mathcode`B="7042 \mathcode`C="7043 \mathcode`D="7044
\mathcode`E="7045 \mathcode`F="7046 \mathcode`G="7047 \mathcode`H="7048
\mathcode`I="7049 \mathcode`J="704A \mathcode`K="704B \mathcode`L="704C
\mathcode`M="704D \mathcode`N="704E \mathcode`O="704F \mathcode`P="7050
\mathcode`Q="7051 \mathcode`R="7052 \mathcode`S="7053 \mathcode`T="7054
\mathcode`U="7055 \mathcode`V="7056 \mathcode`W="7057 \mathcode`X="7058
\mathcode`Y="7059 \mathcode`Z="705A
\newcommand{\ft}[2]{#1_{\mkern-1mu#2}}
\newcommand{\fy}[2]{#1_{\mkern-1mu#2}}
\newcommand{\ff}[2]{#1_{\mkern-1mu#2}}
\newcommand{\fp}[2]{#1_{\mkern-1.5mu#2}}
\newcommand{\fv}[2]{#1_{\mkern-1mu#2}}

\makeatletter
\def\glossaryentry#1#2{\item #1\leaders\hbox{$\m@th \mkern \@dotsep mu\hbox{ }\mkern \@dotsep mu$}\hfill\nobreak\hbox to 1.5em{\hfil\normalfont \normalcolor #2}\hspace*{.5cm}}
\def\printnotation{{%
\def\indexname{Index des notations}
\begin{theindex}
\@input{\jobname.ntn}
\end{theindex}
}}
\makeatother

\makeglossary
\makeindex

\newenvironment{enonce}[1]{\noindent{\textbf{#1}}---\!\, \begin{itshape}}{\end{itshape}}
\newenvironment{proposition}[1][]{\refstepcounter{enonce}\begin{enonce}{Proposition \theenonce{}#1 }}{\end{enonce}}
\newenvironment{theoreme*}[1][]{\refstepcounter{enonce}\begin{enonce}{Théorème #1}}{\end{enonce}}
\newenvironment{definition}[1][]{\refstepcounter{enonce}\begin{enonce}{Définition \theenonce{}#1 }}{\end{enonce}}
\newenvironment{theoreme}[1][]{\refstepcounter{enonce}\begin{enonce}{Théorème \theenonce{}#1 }}{\end{enonce}}
\newenvironment{corollaire}[1][]{\refstepcounter{enonce}\begin{enonce}{Corollaire \theenonce{}#1 }}{\end{enonce}}
\newenvironment{lemme}[1][]{\refstepcounter{enonce}\begin{enonce}{Lemme \theenonce{}#1 }}{\end{enonce}}
\newenvironment{souslemme}[1][]{\refstepcounter{enonce}\begin{enonce}{Sous-lemme \theenonce{}#1 }}{\end{enonce}}
\newenvironment{conjectures}[1][]{\refstepcounter{enonce}\begin{enonce}{Conjectures \theenonce{}#1 }}{\end{enonce}}
\newenvironment{exemple}[1][]{\refstepcounter{enonce}\noindent{\textbf{Exemple \theenonce{}#1 }}---\!\, }{}
\newenvironment{remarque}[1][]{\noindent{\textbf{Remarque{}#1 }}---\!\, }{}
\newenvironment{remarques}[1][]{\noindent{\textbf{Remarques{}#1 }}---\!\, }{}
\newenvironment{demo}[1][]{\noindent{\it{}Démonstration{}#1 }--- }{\hfill $\square$}
\newenvironment{premdemo}[1][]{\noindent{\it{}Première démonstration{}#1 }--- }{\hfill $\square$}
\newenvironment{secdemo}[1][]{\noindent{\it{}Seconde démonstration{}#1 }--- }{\hfill $\square$}

\renewcommand{\emptyset}{\varnothing}
\newcommand{\souligne}[1]{\underline{\smash[b]{#1}}}
\renewcommand{\phi}{\varphi}
\renewcommand{\epsilon}{\varepsilon}
\newcommand{\limind}{\varinjlim}
\newcommand{\longhookrightarrow}{\lhook\joinrel\longrightarrow}
\newcommand{\Aut}{{\mathrm{Aut}}}
\newcommand{\FEt}{\text{FÉt}}
\newcommand{\Ens}{\text{Ens}}
\newcommand{\rg}{{\mathrm{rg}}}
\newcommand{\PGL}{{\mathrm{PGL}}}
\newcommand{\Hom}{{\mathrm{Hom}}}
\newcommand{\Br}{{\mathrm{Br}}}
\newcommand{\Res}{{\mathrm{Res}}}
\newcommand{\Irr}{{\mathrm{Irr}}}
\newcommand{\dgras}{\mathbf{d}}
\renewcommand{\vert}{{\mathrm{vert}}}
\newcommand{\reg}{{\mathrm{reg}}}
\newcommand{\hor}{{\mathrm{hor}}}
\newcommand{\gnr}{{\mathrm{gnr}}}
\newcommand{\sh}{{\mathrm{sh}}}
\newcommand{\h}{{\mathrm{h}}}
\newcommand{\nr}{{\mathrm{nr}}}
\newcommand{\Arond}{\mathscr{A}}
\newcommand{\Crond}{\mathscr{C}}
\newcommand{\Drond}{\mathscr{D}}
\newcommand{\Erond}{\mathscr{E}}
\newcommand{\Frond}{\mathscr{F}}

\newcommand{\Irond}{\mathscr{I}}

\newcommand{\Lrond}{\mathscr{L}}
\newcommand{\Mrond}{\mathscr{M}}
\newcommand{\Orond}{\mathscr{O}}
\newcommand{\Rrond}{\mathscr{R}}
\newcommand{\Srond}{\mathscr{S}}
\newcommand{\Trond}{\mathscr{T}}
\newcommand{\Urond}{\mathscr{U}}

\newcommand{\Wrond}{\mathscr{W}}
\newcommand{\Xrond}{\mathscr{X}}
\newcommand{\Yrond}{\mathscr{Y}}

\newcommand{\Tinfty}{\ft{T}\infty}
\newcommand{\C}{\mathbf{C}}
\newcommand{\Q}{\mathbf{Q}}
\newcommand{\Z}{\mathbf{Z}}
\newcommand{\F}{\mathbf{F}}
\newcommand{\N}{\mathbf{N}}
\renewcommand{\P}{\mathbf{P}}
\newcommand{\A}{\mathbf{A}}
\newcommand{\PPic}{\mathbf{Pic}}

\newcommand{\tors}[2]{{\vphantom{#2}}_{#1}{#2}}

\newcommand{\ensemble}[2]{\{ #1 \, ; \, #2 \}}
\newcommand{\bigensemble}[2]{\left\{ #1 \, ; \, #2 \right\}}
\newcommand{\uplet}[2]{#1, \mskip2.5mu \ldots \mskip-1mu, \mskip2.5mu #2}
\newcommand{\cD}{$(\mathrm{D})$}
\newcommand{\cDg}{$(\mathrm{D_g})$}
\newcommand{\cDgp}{$(\mathrm{D'_g})$}
\newcommand{\cDgpp}{$(\mathrm{D''_g})$}
\newcommand{\cE}{$(\mathrm{E})$}
\newcommand{\cDz}{$(\mathrm{D}_0)$}
\newcommand{\condD}{{condition~\cD}}
\newcommand{\condE}{{condition~\cE}}
\newcommand{\condDz}{{condition~\cDz}}
\newcommand{\condDC}[1]{{condition~$(\mathrm{D}/#1)$}}
\newcommand{\CondD}{Condition~$(\mathrm{D})$}
\newcommand{\CondE}{Condition~$(\mathrm{E})$}
\newcommand{\textcyr}[1]{{\fontencoding{OT2}\fontfamily{wncyr}\fontseries{m}\fontshape{n}\selectfont #1}}
\newcommand{\Sha}{{\mbox{\textcyr{Sh}}}}
\newcommand{\Ba}{\mbox{\textcyr{B}}}
\newcommand{\SG}{{\mathfrak{S}}}
\newcommand{\Gal}{\mathrm{Gal}}
\newcommand{\Pic}{\mathrm{Pic}}
\newcommand{\SDC}[1]{{\mathfrak{S}_{\mathrm{D}/#1}}}
\newcommand{\SD}{{\mathfrak{S}_{\mathrm{D}}}}
\newcommand{\SDz}{{\mathfrak{S}_{\mathrm{D}_0}}}
\newcommand{\SDCp}[1]{{\mathfrak{S}'_{\mathrm{D}/#1}}}
\newcommand{\SDp}{{\mathfrak{S}'_{\mathrm{D}}}}
\newcommand{\SDzp}{{\mathfrak{S}'_{\mathrm{D}_0}}}
\newcommand{\SDSzp}[1]{{\mathfrak{S}'_{\mathrm{D}_0,#1}}}
\newcommand{\SDCpp}[1]{{\mathfrak{S}''_{\mathrm{D}/#1}}}
\newcommand{\SDpp}{{\mathfrak{S}''_{\mathrm{D}}}}
\newcommand{\SDzpp}{{\mathfrak{S}''_{\mathrm{D}_0}}}
\newcommand{\SDSzpp}[1]{{\mathfrak{S}''_{\mathrm{D}_0,#1}}}
\newcommand{\TDC}[1]{{\mathfrak{T}_{\mathrm{D}/#1}}}
\newcommand{\TD}{{\mathfrak{T}_{\mathrm{D}}}}
\newcommand{\TDz}{{\mathfrak{T}_{\mathrm{D}_0}}}
\newcommand{\TDCp}[1]{{\mathfrak{T}'_{\mathrm{D}/#1}}}
\newcommand{\TDp}{{\mathfrak{T}'_{\mathrm{D}}}}
\newcommand{\TDzp}{{\mathfrak{T}'_{\mathrm{D}_0}}}
\newcommand{\TDCpp}[1]{{\mathfrak{T}''_{\mathrm{D}/#1}}}
\newcommand{\TDpp}{{\mathfrak{T}''_{\mathrm{D}}}}
\newcommand{\TDzpp}{{\mathfrak{T}''_{\mathrm{D}_0}}}
\newcommand{\RA}{{\Rrond_{\mathrm{A}}}}
\newcommand{\RO}{{\Rrond_{\mathrm{O}}}}
\newcommand{\RDC}[1]{{\Rrond_{\mathrm{D}/#1}}}
\newcommand{\RDCS}[2]{{\Rrond_{\mathrm{D}/#1,#2}}}
\newcommand{\RD}{{\Rrond_{\mathrm{D}}}}
\newcommand{\RDz}{{\Rrond_{\mathrm{D}_0}}}
\newcommand{\RDS}[1]{{\Rrond_{\mathrm{D},#1}}}
\newcommand{\RDSz}[1]{{\Rrond_{\mathrm{D}_0,#1}}}
\newcommand{\Sel}{\mathrm{Sel}}

\newcommand{\Gm}{\mathbf{G}_\mathrm{m}}
\newcommand{\Gmsur}[1]{\mathbf{G}_{\mathrm{m},#1}}
\newcommand{\Ga}{\mathbf{G}_\mathrm{a}}
\newcommand{\R}{\mathbf{R}}
\newcommand{\Qp}{\mathbf{Q}_p}

\newcommand{\longisoto}{\myxrightarrow{\,\,\sim\,\,}}
\newcommand{\verylongisoto}{\myxrightarrow{\,\,\,\sim\,\,\,}}
\newcommand{\isoto}{\myxrightarrow{\,\sim\,}}
\newcommand{\kappabarre}{{\overline{\kappa}}}
\newcommand{\ksep}{{\mkern1mu\overline{\mkern-1mu{}k\mkern-1mu}\mkern1mu}}
\newcommand{\Ksep}{{\mkern1mu\overline{\mkern-1mu{}K\mkern-1mu}\mkern1mu}}
\newcommand{\Erondsouligne}{{\souligne{\Erond\mkern-1.8mu}\mkern1.8mu}}
\newcommand{\cbarre}{{\overline{c}}}
\newcommand{\sbarre}{{\overline{s}}}
\newcommand{\inftytilde}{{\widetilde{\infty}}}
\newcommand{\atilde}{{\widetilde{a}}}
\newcommand{\xtilde}{{\widetilde{x}}}
\newcommand{\Ctilde}{{\widetilde{C}}}
\newcommand{\Dtilde}{{\widetilde{D}}}
\newcommand{\Mtilde}{{\widetilde{M}}}
\newcommand{\Htilde}{{\widetilde{H}}}
\newcommand{\Ltilde}{{\widetilde{L}}}
\newcommand{\mtilde}{{\widetilde{m}}}
\newcommand{\Mtildez}{{\widetilde{{\setbox0\hbox{$M_0$}\dimen0\ht0\advance\dimen0-.03em\ht0\dimen0\box0}}}}
\newcommand{\stilde}{{\widetilde{{\setbox0\hbox{$s$}\dimen0\ht0\advance\dimen0-.03em\ht0\dimen0\box0}}}}
\newcommand{\Xtilde}{{\widetilde{X}}}
\newcommand{\Ytilde}{{\widetilde{Y}}}
\newcommand{\Ztilde}{{\widetilde{Z}}}
\newcommand{\Hgoth}{{\mathfrak{H}}}
\newcommand{\mgoth}{{\mathfrak{m}}}
\newcommand{\tgoth}{{\mathfrak{t}}}

\DeclareMathOperator{\Spec}{Spec}
\DeclareMathOperator{\inv}{inv}
\DeclareMathOperator{\Card}{Card}
\DeclareMathOperator{\Ker}{Ker}
\DeclareMathOperator{\Coker}{Coker}
\newcommand{\Cores}{\mathrm{Cores}}

\renewcommand{\Im}{\operatorname{Im}}

\begin{document}
\author{Olivier Wittenberg}
\date{}
\title{Intersections de deux quadriques et pinceaux de courbes de genre~$1$}
\subtitle{Intersections of two quadrics and pencils of curves of genus~$1$}
\maketitle
\frontmatter
\selectlanguage{english}
\makeatletter
\startnewpage
\@ifundefined{thispagecropped}{}{\thispagecropped}
\thispagestyle{empty}%
\global\@topnum\z@
\@afterindentfalse
{\parindent \z@ \raggedright\normalfont
\hyphenpenalty \@M
\interlinepenalty\@M
\chapsize\chapstyle
{\noindent\phantom{\vrule height 34pt width 0pt depth 0pt}
\rlap{\smash{\lower 5pt\hbox to\textwidth{{\hfill}}}}\hbox{}
\vskip10pt}
\ignorespaces{Note to the English reader}\par\nobreak
\processchapstarthook
\vskip 15\p@}\@afterheading
\makeatother

In an effort to make the contents of this French-language volume accessible to a wider audience, we have supplemented the original text
with an extensive introduction written in English.  The idea was the Springer Lecture Notes in Mathematics series editors',
and I am grateful to them for suggesting it.
As a side effect, this monograph contains two introductions.  They are not translations of each other; rather, they
are complementary, in spite of their more than substantial overlap.  The emphasis in the English-language introduction is on
the context, while the French-language introduction more systematically describes precise results and refers to the existing literature.

\selectlanguage{francais}
\makeatletter
\global\@topnum\z@
\@afterindentfalse
{\parindent \z@ \raggedright\normalfont
\hyphenpenalty \@M
\interlinepenalty\@M
\chapsize\chapstyle
{\noindent\phantom{\vrule height 34pt width 0pt depth 0pt}
\rlap{\smash{\lower 5pt\hbox to\textwidth{{\hfill}}}}\hbox{}
\vskip10pt}
\ignorespaces{Remerciements}\par\nobreak
\processchapstarthook
\vskip 15\p@}\@afterheading
\makeatother

Le contenu de cet ouvrage reprend celui de ma thèse de doctorat de l'Université Paris-Sud, et j'ai plaisir à remercier
encore une fois tous ceux qui ont contribué d'une manière ou d'une autre à son élaboration.
Ma gratitude va en premier lieu à Jean-Louis Colliot-Thélène, qui a dirigé mon travail et a toujours trouvé le temps de répondre à mes questions.
Je suis reconnaissant envers Peter Swinnerton-Dyer pour ses encouragements et pour les nombreux échanges mathématiques que nous avons eus.
L'influence sur le contenu de cet ouvrage
des idées qu'il a introduites et publiées depuis une douzaine d'années saurait difficilement être surestimée.
Je tiens enfin à remercier David Harari de m'avoir transmis~\cite{harprepub}, sur lequel les résultats du troisième chapitre reposent,
et Per Salberger, Jean-Pierre Serre et Alexei Skorobogatov pour l'intérêt qu'ils portent
à ce travail.

\vspace{.6cm}
\begin{flushright}\noindent
Houston, octobre 2006\hfill {\it Olivier Wittenberg}
\end{flushright}

\tableofcontents
\mainmatter
\selectlanguage{english}
\mathcode`A="7141 \mathcode`B="7142 \mathcode`C="7143 \mathcode`D="7144
\mathcode`E="7145 \mathcode`F="7146 \mathcode`G="7147 \mathcode`H="7148
\mathcode`I="7149 \mathcode`J="714A \mathcode`K="714B \mathcode`L="714C
\mathcode`M="714D \mathcode`N="714E \mathcode`O="714F \mathcode`P="7150
\mathcode`Q="7151 \mathcode`R="7152 \mathcode`S="7153 \mathcode`T="7154
\mathcode`U="7155 \mathcode`V="7156 \mathcode`W="7157 \mathcode`X="7158
\mathcode`Y="7159 \mathcode`Z="715A

\chapter*{General introduction\markboth{General introduction}{General introduction}}
\addcontentsline{toc}{chapter}{General introduction}
\label{pageintroeng}
{
\makeatletter
\def\theequation{\@arabic\c@equation}
\makeatother

This monograph is devoted to the study, from the point of view of arithmetic
geometry, of algebraic surfaces endowed with a pencil of curves of genus~$1$,
and of intersections of two quadrics in projective space.

Some of the foundational questions of arithmetic geometry can be traced back
at least to the work of Diophantus.  Quite generally, given a
system
\begin{equation}
\label{eifeq}
\left\{
\begin{aligned}
f_1(\uplet{x_0}{x_n})&=0\\
\vdots\hspace*{.9cm}\\
f_m(\uplet{x_0}{x_n})&=0
\end{aligned}
\right.
\end{equation}
of polynomial equations with integer (or rational) coefficients, one ultimately wants to understand
the set of tuples $(\uplet{x_0}{x_n})$ in~$\Z^{n+1}$ (or in~$\Q^{n+1}$) which are solutions of the system.
The simplest specific questions that may be asked about this set are whether it is empty or not,
whether it is finite or infinite; if it is finite, whether one can list its elements, and if it is
infinite, whether one can still describe it in a satisfactory manner.
Answers to these questions for particular systems of equations have been known for centuries.
Suffice it to mention the theorems of Fermat, Gauss, and Lagrange on the representation of integers by sums of two, three, or four squares,
and the work of Brahmagupta and of Lagrange on Pell's equation.

The proper language for discussing more recent work on these problems turns out to be that of algebraic geometry.
Indeed, as was first recognised by Weil, the behaviour displayed by the set of integer (or rational) solutions of a system of polynomial equations with
integer (or rational) coefficients depends very closely on the topology of the set of complex solutions, and hence on the
geometry of the underlying algebraic variety.  Naïve measures of the complexity of a system of equations such as the number of unknowns, the
number of equations, or the degrees of the equations, are only of secondary relevance.

We now restrict attention to systems of homogeneous equations, in which case studying integer solutions and studying rational solutions are tantamount to
each other.
Let us fix a smooth and proper algebraic variety~$X$ over a number field~$k$ (for instance $k=\Q$).  If~$X$ is the subvariety of
$n$\nobreakdash-dimensional projective space over~$k$ defined by the system~(\ref{eifeq}), where the~$f_i$ are homogeneous polynomials, then the set
of nonzero solutions $(\uplet{x_0}{x_n}) \in k^{n+1}$ up to scaling is nothing but the set of rational points~$X(k)$ of~$X$ over~$k$.
It is most natural to embed $X(k)$ diagonally into the topological space $X(\A_k)$ of adelic points of~$X$.
The latter may be defined as the product $\prod_{v \in \Omega}X(k_v)$, since~$X$ is proper.
Here $\Omega$ denotes the set of places of~$k$ and~$k_v$ is the completion of~$k$ at a place $v \in \Omega$ (so that if $k=\Q$, the fields $k_v$ for $v \in \Omega$
are the field~$\Qp$ of $p$\nobreakdash-adic numbers for each prime~$p$, and the field~$\R$ of real numbers).
We note that unlike $X(k)$, the set $X(\A_k)$ is fairly easy to understand from a qualitative point of view.  Indeed, the tools of $p$\nobreakdash-adic
(resp.~real) analysis can be applied to the study of $X(k_v)$ for each finite (resp.~real) place~$v$.

It should be clear to the reader from the outset that even determining whether $X(k)$ is empty or not is a notoriously difficult question.
It might be undecidable in general (the analogous question for integral points if $k=\Q$ and~$X$ is not assumed to be proper is known to be
undecidable, by Matiyasevich's negative answer to Hilbert's tenth problem).  Understanding the set $X(k)$ for various classes
of varieties whose geometry is not too intricate nevertheless remains a topic of very high current
interest, as the huge literature on the arithmetic of elliptic curves clearly exemplifies.

Quadric hypersurfaces in projective space (hereafter called quadrics) were the first class of varieties whose arithmetic came to be elucidated.
Work of Legendre, Minkowski, Hilbert, and Hasse resulted in the proof by Hasse in~1924 of the celebrated
Hasse-Minkowski theorem, according to which $X(k)$ is nonempty as soon
as $X(\A_k)$ is nonempty, if~$X$ is a quadric.  Hence in this case
there is a simple criterion for deciding whether $X(k)$ is empty; moreover, if $X(k)$ is nonempty,
then projection from a rational point of~$X$ determines a $k$\nobreakdash-birational equivalence between~$X$ and projective space of dimension $\dim(X)$,
so that $X(k)$ is infinite, is dense in $X(\A_k)$, and can be described straightforwardly.
Note that this also settles the case of curves of genus~$0$, since such curves are canonically isomorphic to plane conics.

Once the Hasse-Minkowski theorem was known, Hasse, and in his wake many more mathematicians than would be fitting to list here,
studied the validity of the implication $X(\A_k)\neq\emptyset\Rightarrow X(k)\neq\emptyset$ in a systematic fashion.  When this implication holds
for a given variety~$X$, one says that~$X$ satisfies the Hasse principle.  A number of classes of varieties were thus shown to
satisfy the Hasse principle, but many counterexamples were also found, in particular among curves of genus~$1$ (by Reichardt and Lind,
then by Selmer), toric varieties (by Hasse and Witt),
compactifications of principal homogeneous spaces under semisimple algebraic groups (by Serre), smooth cubic surfaces (by Swinnerton-Dyer,
then by Cassels and Guy), surfaces endowed with a pencil of curves of genus~$0$ (by Iskovskikh), and smooth intersections of two quadrics in $4$\nobreakdash-dimensional
projective space (by Birch and Swinnerton-Dyer).

The next major step towards the understanding of the Hasse principle was taken by Manin, who presented in his influential 1970 ICM address
a general obstruction to the existence of rational points on arbitrary varieties.  All of the hitherto known counterexamples to the Hasse principle
thus found themselves given a common explanation (except for the cubic surface of Cassels and Guy, which was only properly dealt with fifteen years later,
by Colliot-Thélène, Kanevsky, and Sansuc).  The obstruction considered by Manin relies on the interplay between global class field theory on the one hand and
Brauer groups of schemes, which had been introduced in the sixties by Grothendieck, on the other hand.  We shall be content with simply mentioning that
Manin defines a closed subset $X(\A_k)^\Br \subseteq X(\A_k)$ in the space of adelic points of~$X$ and proves that it must contain $X(k)$.  When $X(\A_k)^\Br$ is empty,
one says that there is a Brauer-Manin obstruction to the existence of a rational point on~$X$.  Of course, if $X(\A_k)$ is nonempty and there is a Brauer-Manin obstruction to
the existence of a rational point on~$X$, the Hasse principle can only fail.

The discovery, by Reichardt, Lind, and Selmer, of counterexamples to the Hasse principle among curves of genus~$1$ over~$\Q$ triggered rapid progress in the arithmetic
theory of elliptic curves, and most notably led to the definition of the Tate-Shafarevich\index{groupe de Tate-Shafarevich} group of an elliptic curve
over a number field.  This group measures the failure of the Hasse principle for principal homogeneous spaces under the given elliptic curve. It is conjectured to always
be finite.  Assuming this conjecture, a classical algorithm, known as descent and developed in particular by Mordell, Weil, Selmer, Cassels, and Tate,
gives (at least in theory) a way to describe the set $X(k)$ if~$X$ is a curve of genus~$1$ over~$k$.
Namely, by applying descent, one first finds out whether~$X(k)$ is empty; if it is not, the choice of a rational point of~$X$ determines,
according to the Mordell-Weil theorem, the structure of a finitely generated abelian
group on~$X(k)$, and descent then provides an explicit presentation of this group.
Already in 1970, Manin saw that the theory of descent on elliptic curves could be recast in the framework of the Brauer-Manin obstruction.
He could thus prove, assuming the finiteness\index{groupe de Tate-Shafarevich!finitude}
of Tate-Shafarevich groups of elliptic curves, that if~$X$ is a curve of genus~$1$, then $X(\A_k)^\Br$ is nonempty if
and only if $X(k)$ is nonempty.

If, for a given variety~$X$, one knows that the set $X(k)$ is dense in $X(\A_k)^\Br$, then many qualitative questions about $X(k)$ can be reduced to questions
about $X(\A_k)^\Br$.  The latter are often much more tractable (and sometimes entirely trivial: it can happen that $X(\A_k)^\Br=X(\A_k)$ for purely algebraic reasons,
for instance this equality holds if~$X$ is a smooth complete intersection of dimension~$\geq 3$ in projective space).
General theoretical considerations, as well as positive results for specific families of varieties,
led Colliot-Thélène and Sansuc to conjecture in~1979 that $X(k)$ is dense in $X(\A_k)^\Br$ as soon as~$X$ is a rational surface (\emph{i.e.}, a surface which becomes birational to the
projective plane after an extension of scalars).
This conjecture still stands today as one of the most important unsolved problems in the arithmetic of surfaces.

Rational surfaces are really the simplest of all surfaces from the geometric point of view.  Despite unexpected recent progress concerning
surfaces endowed with a pencil of curves of genus~$1$, the arithmetic of nonrational surfaces is so poorly understood on the whole
that it does not even fit into a conjectural picture. (There are however a few important but specific conjectures about nonrational surfaces, most notably Lang's conjectures.)
The situation for curves of arbitrary genus is only slightly better: if~$X$ is a curve of genus~$\geq 2$, then $X(k)$ is known to be finite by Faltings' theorem, but it is
an open question whether $X(k)$ is dense in $X(\A_k)^\Br$ (for nonrational surfaces, $X(k)$ need not be dense
in $X(\A_k)^\Br$; Skorobogatov gave an example of this in 1998).

Several methods have been developed to attack the question of the density of~$X(k)$ in $X(\A_k)^\Br$ for various classes of varieties.
One of them is the theory of descent on rational varieties of Colliot-Thélène and Sansuc.  It mimics (and actually generalises) the classical
theory of descent on elliptic curves.  When applied to a rational variety~$X$, it reduces the study of the Brauer-Manin obstruction on~$X$ to the study
of the Hasse principle on certain auxiliary varieties of higher dimension whose arithmetic is hopefully easier to understand.
One of the complete successes of this method was the proof in~1987, by
Colliot-Thélène, Sansuc, and Swinnerton-Dyer, of the density of~$X(k)$ in $X(\A_k)^\Br$ when~$X$ is a Châtelet surface, and of the Hasse
principle for smooth intersections of two quadrics in projective space of dimension~$\geq 8$ over~$k$.
The Hasse principle for smooth intersections of two quadrics in projective space of dimension~$n$ over a number field~$k$ had previously been
established for $k=\Q$ and $n \geq 12$ by Mordell in~1959, and for $k=\Q$ and $n \geq 10$ by Swinnerton-Dyer in~1964.
(Conjecturally it holds for any number field~$k$ and any $n \geq 5$.)
Châtelet surfaces are the simplest nontrivial conic bundle surfaces (\emph{i.e.}, the simplest nontrivial rational surfaces
which carry a pencil of curves of genus~$0$).
The number of singular geometric fibres in the conic bundle structure of a Châtelet surface
is exactly four.
By~1991, thanks to a number of other results of Colliot-Thélène, Skorobogatov, and Salberger,
it was known that $X(k)$ is dense in $X(\A_k)^\Br$ if~$X$ is a conic bundle surface with at most five singular geometric fibres.
More complicated conic bundle surfaces have so far resisted all attempts at studying their rational points, except
for a specific family of conic bundles with six singular geometric fibres, by Swinnerton-Dyer in~1999.
As soon as the number of singular geometric fibres exceeds~$7$, if~$X$ is a conic bundle surface such that $X(k)\neq\emptyset$, one does not even
know that $X(k)$ is dense in~$X$ for the Zariski topology.

It is possible to obtain much stronger results about the arithmetic of conic bundle surfaces if one accepts\index{hypothèse de Schinzel} Schinzel's hypothesis~(H).
This hypothesis is a conjecture on the existence of prime values of polynomials with integer coefficients, which subsumes in particular the twin prime conjecture;
its quantitative version, due to Bateman and Horn, is supported by compelling numerical evidence.
The relevance of Schinzel's hypothesis to the study of rational points on conic bundle surfaces, at least over~$\Q$,
was first noticed by Colliot-Thélène and Sansuc, around~1978.
In~1992, Serre proved the density of $X(k)$ in $X(\A_k)^\Br$ for arbitrary conic bundle surfaces over
arbitrary number fields, under~(H).
Thus the conjecture of Colliot-Thélène and Sansuc mentioned above is known to hold, under~(H), for a substantial class of rational surfaces.
What are the remaining cases~?  According to Enriques, Manin, and Iskovskikh, any rational surface is, up to $k$\nobreakdash-birational equivalence,
either a conic bundle surface or a del Pezzo surface.  With each del Pezzo surface is associated an integer comprised between~$1$ and~$9$, called the degree;
the lower the degree, the more complex the surface.  From a qualitative point of view, the arithmetic of del Pezzo surfaces of degree~$\geq 5$ is featureless.
Indeed, for such surfaces, the set $X(k)$ is dense in $X(\A_k)$ (in particular the Hasse principle holds), by theorems
of Manin and Swinnerton-Dyer (proven in~1966 and~1970).
Del Pezzo surfaces of degree~$3$ or~$4$ are very familiar objects: the del Pezzo surfaces of degree~$3$ (resp.~$4$) are exactly the smooth cubic surfaces
(resp.~the smooth intersections of two quadrics in $4$\nobreakdash-dimensional projective space).
They do not satisfy the Hasse principle in general.  Some del Pezzo surfaces of degree~$4$ are at the same time conic bundle surfaces.  These conic bundle
surfaces necessarily have four singular geometric fibres, so that the results discussed in the previous paragraph are applicable to them.
However, most del Pezzo surfaces of degree~$4$ do not enjoy such a structure, and their arithmetic seems to be entirely unapproachable by
the method of descent.

One of the main goals of this monograph is to establish, for the first time, positive results on the Hasse principle for general del Pezzo surfaces of degree~$4$.
We shall take advantage of the abundance of curves of genus~$1$ lying on del Pezzo surfaces of degree~$4$.
Because of this, we have to assume the finiteness of Tate-Shafarevich groups of elliptic curves over number fields.  We also have to assume~(H).
Hopefully these assumptions are somewhat compensated for by the generality of our conclusions:
only finitely many explicit exceptional families of del Pezzo surfaces
of degree~$4$ evade our analysis.  As a corollary of these results, we prove, again under~(H) and the finiteness of
Tate-Shafarevich groups, but this time without exceptions, that \emph{smooth intersections of two quadrics in $n$\nobreakdash-dimensional projective space satisfy the Hasse
principle, over
any number field and for any $n \geq 5$}.

Swinnerton-Dyer introduced in~1993 a new technique for dealing with surfaces which carry a pencil of curves of genus~$1$.
It was greatly refined by Colliot-Thélène, Skorobogatov, and Swinnerton-Dyer in~1998, and has already been the key to a number of spectacular results.
Swinnerton-Dyer's initial purpose was the study, under~(H) and assuming the finiteness of Tate-Shafarevich groups, of
del Pezzo surfaces of degree~$4$ defined by simultaneously diagonal quadratic forms
(these constitute an exceptional family of del Pezzo surfaces of degree~$4$).  In~2001, Swinnerton-Dyer was able to prove the Hasse principle
for most diagonal cubic surfaces over~$\Q$ (these form an exceptional family of del Pezzo surfaces of degree~$3$),
assuming the finiteness of Tate-Shafarevich groups but not~(H).  Even the arithmetic of some nonrational surfaces could be tackled.  The first results in this connection
were obtained by Colliot-Thélène, Skorobogatov, and Swinnerton-Dyer in~1998.  Shortly afterwards,
Swinnerton-Dyer successfully studied a family of diagonal quartic surfaces (these are~$K3$ surfaces), under~(H) and the finiteness of Tate-Shafarevich groups,
and in~2005, Skorobogatov and Swinnerton-Dyer established the Hasse principle for a class of Kummer surfaces, assuming the finiteness of Tate-Shafarevich groups but not~(H).

We shall need to develop this technique further before we can look into del Pezzo surfaces of degree~$4$.
The first two chapters of this volume are thus devoted to (not necessarily rational) surfaces endowed with a pencil of curves of genus~$1$;
we refer the reader to their respective introductions for more details about their contents and their originality.
The third chapter deals with del Pezzo surfaces of degree~$4$ and intersections of two quadrics.  The results of the second chapter
are a major ingredient in the proof of the main theorems of the third chapter, but quite a sizeable amount of additional work is required.
Indeed the results of the first two chapters cannot be applied to a general del Pezzo surface of degree~$4$.
We shall apply them to certain nonrational surfaces of Kodaira dimension~$1$ built out of del Pezzo surfaces of degree~$4$,
following an ingenious construction of Swinnerton-Dyer that he introduced in~2000.
Another ingredient is a theorem of general interest (also established in the third chapter) about the so-called ``vertical'' Brauer-Manin obstruction for fibrations
over projective space, which extends a previous result of Colliot-Thélène, Skorobogatov, and Swinnerton-Dyer about fibrations over the projective line.
The passage from dimension~$1$ to dimension~$>1$ is not a formality; our proof rests on a recent theorem of Harari.
}

\selectlanguage{francais}
\renewcommand{\leq}{\leqslant}
\renewcommand{\geq}{\geqslant}
\mathcode`A="7041 \mathcode`B="7042 \mathcode`C="7043 \mathcode`D="7044
\mathcode`E="7045 \mathcode`F="7046 \mathcode`G="7047 \mathcode`H="7048
\mathcode`I="7049 \mathcode`J="704A \mathcode`K="704B \mathcode`L="704C
\mathcode`M="704D \mathcode`N="704E \mathcode`O="704F \mathcode`P="7050
\mathcode`Q="7051 \mathcode`R="7052 \mathcode`S="7053 \mathcode`T="7054
\mathcode`U="7055 \mathcode`V="7056 \mathcode`W="7057 \mathcode`X="7058
\mathcode`Y="7059 \mathcode`Z="705A

\chapter*{Introduction générale\markboth{Introduction générale}{Introduction générale}}
\addcontentsline{toc}{chapter}{Introduction générale}

Le thème de cet ouvrage est l'étude qualitative des points
rationnels des variétés algébriques définies sur un corps de nombres.
Dans les deux premiers chapitres nous considérons les surfaces munies d'un pinceau de courbes de genre~$1$.
Le troisième chapitre, qui repose sur le second, est quant à lui consacré aux intersections de deux quadriques et en particulier
aux surfaces de del Pezzo de degré~$4$.
Il est possible de lire le troisième chapitre après
avoir seulement pris connaissance des énoncés du second.
La lecture détaillée du second chapitre présuppose néanmoins celle du premier chapitre.

Les quelques paragraphes qui suivent ont pour but d'énoncer les deux résultats principaux de ce travail,
après avoir rappelé le contexte dans lequel ils se placent.
Il~est impossible de donner une vue d'ensemble raisonnablement complète
des questions, conjectures, théorèmes et méthodes qui interviennent dans l'étude
qualitative des points rationnels des variétés algébriques sur les corps de nombres
sans que l'exposé ne prenne de proportions incongrues.
Nous nous sommes donc borné à décrire les résultats directement liés aux problèmes
abordés dans cet ouvrage.
Le~lecteur trouvera de très nombreux compléments dans les articles de synthèse
\cite{cttoulouse}, \cite{cttrichy}, \cite{ctbudapest}, \cite{harariweakapp},
\cite{peyrebour}, \cite{swdprogress} ainsi que dans le livre~\cite{skotorsors}.

\bigskip
Soient~$k$ un corps de nombres et~$X$ une variété algébrique propre, lisse et
géométriquement connexe sur~$k$.
Les premières questions qualitatives qui viennent à l'esprit lorsqu'on
considère l'ensemble~$X(k)$ des points rationnels de~$X$ sont les suivantes:
cet ensemble est-il non vide~? est-il infini~? est-il dense
dans~$X$ pour la topologie de Zariski~?
existe-t-il un algorithme (au sens de Turing), aussi inefficace soit-il, qui puisse répondre à chacune
de ces questions en temps fini pour tout~$X$~?
Rappelons que
les travaux de Matiyasevich et de ses prédécesseurs sur le dixième
problème de Hilbert concernaient les solutions \emph{entières}
à des systèmes d'équations polynomiales à coefficients entiers; ils ont établi
qu'aucun algorithme ne peut déterminer systématiquement si de
telles solutions existent.
Le problème analogue avec des coefficients
rationnels est une question ouverte (l'avis selon lequel sa réponse est négative n'est d'ailleurs
pas universellement partagé).

Pour qu'il existe un point rationnel sur~$X$, il est évidemment nécessaire
que $X(k_v)\neq\emptyset$ pour toute place $v\in \Omega$, où~$\Omega$
désigne l'ensemble des places de~$k$ et~$k_v$ le complété de~$k$ en~$v$.
Lorsque cette condition est suffisante, c'est\nobreakdash-à-\-dire lorsque
$\prod_{v \in \Omega} X(k_v) \neq \emptyset \Rightarrow X(k) \neq \emptyset$,
on dit que \emph{la variété~$X$ satisfait au principe de Hasse}.
C'est en effet Hasse qui, le premier, étudia systématiquement la validité de
cette implication.  Il démontra notamment, dans la lignée des travaux de Legendre
et de Minkowski, le célèbre théorème de Hasse-Minkowski: toute quadrique sur~$k$ satisfait
au principe de Hasse.

Pendant la cinquantaine d'années qui suivirent cette découverte, le principe de Hasse
fut établi pour de nombreuses classes de variétés, par des auteurs encore plus nombreux;
dans le même temps, on vit apparaître
toute une série d'exemples de variétés pour lesquelles le principe de Hasse
est en défaut.  Le premier fut donné par Hasse lui-même.
Il y en eut notamment parmi les courbes de genre~$1$ (Reichardt et Lind,
puis Selmer),
parmi les compactifications lisses d'espaces principaux homogènes sous des groupes
algé\-briques semi-simples (Serre),
parmi les surfaces cubiques lisses (Swinnerton-Dyer, puis Cassels et Guy),
parmi les surfaces fibrées en coniques au-dessus de~$\P^1_k$ (Iskovskikh)
et parmi les intersections lisses de deux quadriques dans~$\P^4_k$ (Birch et
Swinnerton-Dyer).

Manin mit un terme à cette situation inconfortable en 1970: il
dégagea une condition nécessaire toute générale et cependant non triviale à la validité
du principe de Hasse, depuis lors connue sous le nom d'\emph{obstruction de Brauer-Manin}
(cf.~\cite{maninicm}).  Cette obstruction permit d'expliquer tous les contre-exemples
au principe de Hasse connus à cette date (excepté celui de Cassels et Guy, qui dut attendre
une quinzaine d'années).

L'obstruction de Brauer-Manin se définit comme suit. Soit $\Br(X)$ le groupe
de Brauer cohomologique de~$X$, défini par Grothendieck (cf.~\cite{grothbr2});
si~$\kappa(X)$ désigne le corps des fonctions de~$X$, le groupe~$\Br(X)$
n'est autre que le sous-groupe du groupe de Brauer usuel $\Br(\kappa(X))$
constitué des classes d'algèbres centrales simples
sur~$\kappa(X)$ qui admettent (en un certain sens)
une bonne spécialisation en tout point de~$X$.
Si~$\ell$ est une extension de~$k$
et $P \in X(\ell)$, on note $A \mapsto A(P)$ la flèche
$\Br(X) \rightarrow \Br(\ell)$
de spécialisation (ou encore: d'évaluation) en~$P$, où $\Br(\ell)$ est le groupe de Brauer
usuel du corps~$\ell$.
La théorie du corps de classes local fournit une injection $\inv_v \colon
\Br(k_v) \hookrightarrow \Q/\Z$
pour toute place $v\in \Omega$.  Grâce à la propreté de~$X$,
on peut montrer que pour tout $A \in \Br(X)$,
il existe un ensemble fini $S \subset \Omega$ tel que
pour tout $v \in \Omega \setminus S$ et tout $P \in X(k_v)$, on ait
$A(P)=0$.  Ceci donne un sens à la somme $\sum_{v \in \Omega} \inv_v A(\fp{P}v) \in \Q/\Z$
pour tout $A \in \Br(X)$ et tout $(\fp{P}v)_{v \in \Omega} \in \prod_{v \in \Omega}X(k_v)$.
Lorsque, pour une famille $(\fp{P}v)_{v \in \Omega}$ donnée, cette somme est nulle
pour tout $A \in B$, où~$B$ est un sous-groupe de~$\Br(X)$, on dit que
\emph{la famille $(\fp{P}v)_{v \in \Omega}$ est orthogonale au groupe~$B$}.
L'observation de Manin fut la suivante: l'image de tout point rationnel de~$X$
par l'inclusion diagonale
$X(k) \subset \prod_{v \in \Omega}X(k_v)$ est orthogonale
à~$\Br(X)$.  Cela résulte immédiatement de ce que la suite exacte\index{loi de réciprocité globale}
$$
\xymatrix@C=6.5ex{
0 \ar[r] & \Br(k) \ar[r] & \displaystyle\bigoplus_{v \in \Omega} \Br(k_v) \ar[r]^(.57){\sum \inv_v}
& \Q/\Z \ar[r] & 0 \rlap{\text{,}}
}
$$
issue de la théorie du corps de classes global, est un complexe.
Autrement dit, notant~$\A_k$ l'algèbre des adèles de~$k$ (de sorte que
$\prod_{v \in \Omega}X(k_v)=X(\A_k)$ puisque~$X$ est propre) et
$X(\A_k)^B$ le sous-ensemble de $X(\A_k)$ constitué des familles orthogonales à~$B$,
on a toujours $X(k) \subset X(\A_k)^B$.  En particulier, si $X(\A_k)^B=\emptyset$,
alors $X(k)=\emptyset$; on dit dans ce cas
qu'\emph{il y a une obstruction de Brauer-Manin à l'existence
d'un point rationnel sur~$X$ (associée à~$B$)}.

Une fois cette obstruction définie, il est naturel de se demander pour quelles
variétés l'absence d'obstruction de Brauer-Manin suffit à assurer l'existence
d'un point rationnel.  Autrement dit, il s'agit d'étudier la validité de l'implication
$X(\A_k)^{\Br(X)}\neq \emptyset \Rightarrow X(k)\neq\emptyset$.
Notons que celle-ci
nous renseigne non seulement sur l'existence d'un point rationnel sur~$X$,
mais aussi sur la question algorithmique mentionnée précédemment; en effet,
comme le remarque Poonen (cf.~\cite[Remark~5.3]{poonencurves}),
il existe un algorithme permettant de déterminer si~$X(k)$ est vide
pour toutes les variétés~$X$ pour lesquelles
$X(\A_k)^{\Br(X)}\neq \emptyset \Rightarrow X(k)\neq\emptyset$.

Il n'est pas difficile de voir que $X(\A_k)^{\Br(X)}$ est un fermé de~$X(\A_k)$.
L'adhérence $\overline{X(k)}$ de~$X(k)$ dans~$X(\A_k)$ est donc incluse dans~$X(\A_k)^{\Br(X)}$.
La question la plus pertinente semble être la suivante: pour quelles variétés
a-t-on $\overline{X(k)}=X(\A_k)^{\Br(X)}$~?
Il a toujours été clair que cette égalité ne devait pas être satisfaite en général
(cela contredirait la conjecture selon laquelle les points rationnels ne
sont jamais Zariski-denses sur une variété de type général) mais ce n'est
que relativement récemment que le premier contre-exemple fut exhibé (cf.~\cite{skoinvent};
il s'agit d'un contre-exemple tel que $X(k)=\emptyset$ mais $X(\A_k)^{\Br(X)}\neq\emptyset$).
Il existe néanmoins plusieurs classes intéressantes de variétés pour lesquelles
il est conjecturé que $\overline{X(k)}=X(\A_k)^{\Br(X)}$; nous en citerons quelques-unes ci-dessous.

Voici ce qui est connu au sujet des points rationnels
sur les courbes, du point de vue qualitatif.  Si~$X$ est une courbe
de genre~$0$, c'est une conique plane, de sorte que $\overline{X(k)}=X(\A_k)$.
En particulier, le principe de Hasse est satisfait, l'ensemble $X(k)$ est infini
s'il est non vide, et l'on peut décider algorithmiquement l'existence d'un point
rationnel sur~$X$.
Si~$X$ est une courbe de genre~$1$, c'est un torseur sous une courbe elliptique.
Manin a reformulé la théorie de la descente sur les courbes elliptiques en termes de
l'obstruction de Brauer-Manin (\emph{via} un résultat de Cassels), ce qui lui a permis d'établir
l'implication $X(\A_k)^{\Br(X)}\neq\emptyset \Rightarrow X(k)\neq\emptyset$
en supposant le groupe de Tate-Shafarevich de la jacobienne de~$X$ fini.
(Plus généralement, Wang a étudié la validité de l'égalité
$\overline{X(k)}=X(\A_k)^{\Br(X)}$.)
Si~$X$ est une courbe de genre~\mbox{$\geq 2$}, le théorème de Faltings
assure que l'ensemble $X(k)$ est fini.  Divers auteurs conjecturent
que $X(\A_k)^{\Br(X)}\neq\emptyset \Rightarrow X(k)\neq\emptyset$ (et en particulier
que l'on peut décider algorithmiquement l'existence d'un point
rationnel sur~$X$).  D'autre part, si l'on admet une version effective de la conjecture $abc$, on dispose d'un
algorithme permettant d'établir la liste de \emph{tous} les éléments de~$X(k)$.

Après les courbes viennent les surfaces.  Compte tenu du principe général, remontant à Weil, selon lequel
l'arithmétique d'une variété est gouvernée par sa géométrie, il est raisonnable de se limiter dans un premier
temps aux surfaces dont la géométrie est la plus simple, à savoir celles qui sont
rationnelles (sous-entendu, géométriquement).  Il s'avère que les trois
conditions $X(k)\neq\emptyset$,
$X(\A_k)^{\Br(X)}\neq\emptyset$ et $\overline{X(k)}=X(\A_k)^{\Br(X)}$
sont des invariants $k$\nobreakdash-birationnels des variétés
propres et lisses sur~$k$.
Par ailleurs, Enriques, Iskovskikh et Manin ont
établi une classification des surfaces rationnelles à équivalence
$k$\nobreakdash-birationnelle près (cf.~\cite{iskomin}): toute surface rationnelle sur~$k$
est $k$\nobreakdash-birationnelle à une surface de del Pezzo (c'est-à-dire à une
surface propre et lisse de faisceau anti-canonique ample)
ou à un fibré en coniques au-dessus
d'une conique.  Il suffit donc de s'intéresser à ces deux familles de surfaces rationnelles.

Supposons que~$X$ soit une surface rationnelle de l'un de ces deux types.
On appelle \emph{degré} de~$X$, et l'on note~$d$, le nombre d'auto-intersection
du faisceau canonique de~$X$.
Si~$X$ est une surface de del Pezzo, alors $1 \leq d \leq 9$; si~$X$ est un fibré en coniques
au-dessus d'une conique avec~$r$ fibres géométriques singulières, alors $d=8-r$ (de sorte
que $d \leq 8$ et qu'il n'y a pas de borne inférieure pour~$d$).
Manin et Swinnerton-Dyer
ont résolu toutes les questions
qualitatives concernant l'ensemble~$X(k)$ lorsque $d \geq 5$
(cf.~\cite{maninratsurf} et~\cite{swddelpezzo5}).
Il~résulte en effet de leurs
travaux que l'on a toujours $\overline{X(k)}=X(\A_k)$ si $d \geq 5$; en particulier, le principe
de Hasse est vérifié.

Restent les surfaces fibrées en coniques au-dessus d'une conique avec \mbox{$r \geq 4$} fibres
géométriques singulières et les surfaces de del Pezzo de degré \mbox{$d\leq 4$}.
Le meilleur résultat connu à ce jour concernant les points rationnels
des surfaces fibrées en coniques au-dessus d'une conique
est dû à Colliot-Thélène, Salberger, Sansuc, Skorobogatov et Swinnerton-Dyer,
qui établirent l'égalité $\overline{X(k)}=X(\A_k)^{\Br(X)}$ dès que $r \leq 5$
(c'est-à-dire $d \geq 3$).
La démonstration en est extrêmement subtile (cf.~\cite{ctsansdi}, \cite{ctsansdii},
\cite{ctsemtn}, \cite{salpezzoconique}, \cite{salskoweak});
notamment, le cas particulier
des surfaces de Châtelet (cas où~$X$ est de degré~$4$ mais n'est pas
une surface de del Pezzo)
fut l'un des grands succès de la théorie de la descente
de Colliot-Thélène et Sansuc.
Swinnerton-Dyer a par ailleurs obtenu des résultats pour certaines surfaces
fibrées en coniques au-dessus d'une conique avec six fibres géométriques singulières
(cf.~\cite[§7.4]{skotorsors}).

Soit~$X$ une surface de del Pezzo de degré $d \leq 4$. Si $d \geq 3$, le faisceau anti-canonique
de~$X$ est très ample et permet de voir~$X$ comme une surface de degré~$d$ dans~$\P^d_k$;
les surfaces de del Pezzo de degré~$3$ sont exactement les surfaces cubiques
lisses et les surfaces de del Pezzo de degré~$4$ sont les intersections lisses et de
codimension~$2$
de deux quadriques dans~$\P^4_k$ (cf.~\cite[p.~96]{maninratsurf}).
Ainsi les exemples de Swinnerton-Dyer, Cassels et Guy puis Birch et Swinnerton-Dyer
montrent-ils que le principe de Hasse peut être en défaut pour les surfaces de del
Pezzo de degré~$3$ et~$4$ (et \emph{a fortiori} pour les surfaces de del Pezzo de degré~$2$:
il suffit d'éclater un point de degré~$2$ sur une surface de del Pezzo de degré~$4$
qui est un contre-exemple au principe de Hasse).
Salberger et Skorobogatov~\cite{salskoweak}
ont établi que si $d=4$ et si $X(k)\neq\emptyset$, alors
$\overline{X(k)}=X(\A_k)^{\Br(X)}$.  C'est le seul résultat véritablement
général dont on dispose à l'heure
actuelle pour les surfaces de del Pezzo de degré~$\leq 4$; il ne dit rien sur l'implication
$X(\A_k)^{\Br(X)}\neq\emptyset \Rightarrow X(k)\neq\emptyset$, qui n'est connue
que dans très peu de cas (dont nous dresserons la liste plus bas lorsque $d=4$).

Suite à leurs travaux sur la descente pour les variétés rationnelles, Colliot-Thélène
et Sansuc ont conjecturé en~1979 que l'égalité $\overline{X(k)}=X(\A_k)^{\Br(X)}$ vaut
pour toute surface rationnelle~$X$ (et même pour toute variété rationnelle, cf.~\cite[p.~319]{cttoulouse}).
Comme il ressort des deux paragraphes qui précèdent,
on est encore très loin de savoir établir cette conjecture.

La situation est déjà plus favorable si l'on admet l'hypothèse de Schinzel\index{hypothèse de Schinzel}
(cf.~\cite{schinzelsierp}).  Il s'agit d'une conjecture hardie englobant notamment
la conjecture des nombres premiers jumeaux.  Elle s'énonce ainsi:
si $\uplet{p_1}{p_s}\in \Z[t]$ sont des polynômes irréductibles de coefficients
dominants positifs et s'il n'existe pas d'entier~$n>1$ divisant
$\prod_{i=1}^s p_i(m)$ pour tout $m \in \Z$, alors il existe une infinité
de $m \in \Z$ tels que les entiers $\uplet{p_1(m)}{p_s(m)}$ soient tous des nombres premiers.
Le cas où $s=1$ et $\deg(p_1)=1$ est un théorème bien connu de Dirichlet.
Hasse s'en servit pour démontrer le théorème de Hasse-Minkowski;
en~1978, Colliot-Thélène et Sansuc remarquèrent qu'un
argument similaire à celui employé par Hasse permet d'établir l'existence de points rationnels
sur certaines surfaces fibrées en coniques au-dessus de~$\P^1_\Q$, à condition de remplacer
le théorème de Dirichlet par
l'hypothèse de Schinzel (cf.~\cite{ctsanschinzel}).
En~1992, Serre généralisa ce résultat et prouva: si~$X$ est une surface fibrée en coniques
au-dessus d'une conique et si l'on admet l'hypothèse de Schinzel,
alors $\overline{X(k)}=X(\A_k)^{\Br(X)}$ (cf.~\cite[Chapitre~2, Annexe, Théorème~7.6]{serrecg}
et \cite[§4]{cscrelle94}; le corps de nombres~$k$ est ici quelconque).
Ainsi, sous l'hypothèse de Schinzel, les seules surfaces
rationnelles pour lesquelles il reste à établir
que $\overline{X(k)}=X(\A_k)^{\Br(X)}$ sont les surfaces de del Pezzo de degré $1$, $2$, $3$ et~$4$.

La preuve de l'implication $X(\A_k)^{\Br(X)}\neq\emptyset \Rightarrow X(k)\neq\emptyset$
lorsque~$X$ est un fibré en coniques au-dessus d'une conique, sous l'hypothèse
de Schinzel, consiste à trouver un point rationnel de la conique de base
au-dessus duquel la fibre de~$X$ possède un $k_v$\nobreakdash-point pour toute place $v\in \Omega$.
Il apparut bientôt que la démonstration de l'existence de telles fibres
à partir de l'hypothèse de Schinzel fonctionnait pour de nombreuses fibrations
au-dessus de~$\P^1_k$, sans aucune restriction sur la géométrie des fibres lisses;
ceci fut formalisé par Colliot-Thélène, Skorobogatov
et Swinnerton-Dyer dans~\cite{csscrelle98}.
Lorsque les fibres lisses satisfont au principe de Hasse et que les
fibres singulières possèdent une composante irréductible de multiplicité~$1$ déployée
par une extension abélienne du corps de base, on obtient ainsi
l'implication $X(\A_k)^{\Br(X)}\neq\emptyset \Rightarrow X(k)\neq\emptyset$ pour
l'espace total de la fibration considérée.

En 1993, Swinnerton-Dyer réussit le tour de force de combiner l'argument que l'on vient d'évoquer
avec un processus de $2$\nobreakdash-descente sur une courbe elliptique variable afin d'établir
le principe de Hasse pour certaines surfaces de del Pezzo de degré~$4$ sur~$\Q$
munies d'une fibration
en courbes de genre~$1$ au-dessus de~$\P^1_\Q$, en admettant l'hypothèse
de Schinzel et la finitude des groupes de Tate-Shafarevich des courbes elliptiques sur~$\Q$
(cf.~\cite{swdegloffstein}).
Si \mbox{$\pi \colon X \rightarrow \P^1_\Q$} est une fibration en courbes de genre~$1$ et
de période~$2$
et si $b \in \P^1(\Q)$ est tel que la fibre $X_b$ soit lisse et admette un $k_v$\nobreakdash-point
pour tout $v \in \Omega$, il se peut très bien que~$X_b$ n'ait aucun point rationnel;
mais Swinnerton-Dyer parvint à trouver de tels~$b$ pour lesquels
le sous-groupe de $2$\nobreakdash-torsion du groupe de Tate-Shafarevich de la jacobienne de~$X_b$
est nul, ce qui force~$X_b$ à contenir un point rationnel.
Colliot-Thélène, Skorobogatov et Swinnerton-Dyer~\cite{css} généralisèrent l'argument
de Swinnerton-Dyer~\cite{swdegloffstein}
à tout corps de nombres et le formulèrent de manière que la donnée
de départ soit une fibration en courbes de genre~$1$ abstraite
(alors que Swinnerton-Dyer considérait une famille particulière de fibrations).
Dans le cas où~$X$ est une surface de del Pezzo de degré~$4$, ils purent ainsi
établir le principe de Hasse pour~$X$, sous l'hypothèse de Schinzel et la finitude des groupes
de Tate-Shafarevich des courbes elliptiques sur~$k$, lorsque~$X$ est
une intersection suffisamment générale (en un sens explicite)
de deux quadriques \emph{simultanément diagonales} dans~$\P^4_k$.

Des idées similaires (mais plus délicates)
permirent à Swinnerton-Dyer de démontrer le principe de Hasse
pour une large classe de surfaces cubiques diagonales sur~$\Q$, en supposant la finitude
des groupes de Tate-Shafarevich des courbes elliptiques mais non l'hypothèse
de Schinzel (cf.~\cite{sdcubdiag}; l'hypothèse de Schinzel est remplacée par
le théorème de Dirichlet, qui suffit dans le cas considéré).

De manière assez inattendue, le théorème principal de~\cite{css} et ses variantes
s'appliquent également à des surfaces qui ne sont pas rationnelles, par exemple
des surfaces~$K3$.  Swinnerton-Dyer prouva ainsi le principe de Hasse pour une famille
de surfaces quartiques diagonales, sous l'hypothèse de Schinzel et la finitude des groupes
de Tate-Shafarevich (cf.~\cite{sddiagquartic}).
Plus récemment, et malgré des difficultés techniques formidables,
Skorobogatov et Swinnerton-Dyer~\cite{ssdkummer}
ont pu démontrer le principe de Hasse pour certaines surfaces de Kummer, en admettant
la finitude des groupes de Tate-Shafarevich des courbes elliptiques sur~$k$ mais non l'hypothèse
de Schinzel (remplacée par le théorème de Dirichlet).

Pour ce qui est des variétés de dimension~$>2$, il est difficile à l'heure
actuelle d'entamer une étude systématique des propriétés qualitatives de l'ensemble
de leurs points
rationnels, pour la simple raison que la classification géométrique
des variétés à étudier n'est pas encore assez aboutie (et même dans les cas où elle est
comprise, elle nécessiterait d'être précisée du point de vue $k$\nobreakdash-birationnel
pour~$k$ non algébriquement clos).
On peut néanmoins considérer des variétés qui sont «~simples~» d'un certain point de vue;
par exemple, de nombreux résultats ont été obtenus par divers auteurs
pour les compactifications lisses d'espaces homogènes sous des groupes algébriques linéaires
ou, dans une autre direction, pour les espaces totaux des fibrations au-dessus de~$\P^1_k$
dont les fibres vérifient l'égalité $\overline{X(k)}=X(\A_k)^{\Br(X)}$.
Le lecteur trouvera dans~\cite{cttoulouse}
une conjecture prédisant l'égalité $\overline{X(k)}=X(\A_k)^{\Br(X)}$ pour une large
classe de variétés non nécessairement
rationnelles qui englobe notamment les surfaces fibrées en courbes
de genre~$1$ étudiées dans~\cite{css}.
Une autre possibilité est de considérer les variétés projectives définies par
les équations les plus simples possibles: après les quadriques, qui satisfont au principe de Hasse,
se trouvent les intersections de deux quadriques dans~$\P^n_k$, puis
les hypersurfaces cubiques, etc.  Pour chaque tel type de variété projective,
la théorie analytique des nombres permet d'établir le principe de Hasse
si~$n$ est assez grand, avec une borne explicite,
au moins lorsque $k=\Q$.  Toute la difficulté est de s'approcher de la valeur minimale de~$n$.

Soient $n \geq 4$ et~$X$ une intersection lisse de deux quadriques dans $\P^n_k$.
Comme on l'a vu, il est conjecturé que $\overline{X(k)}=X(\A_k)^{\Br(X)}$ lorsque~$n=4$;
cette conjecture
implique que $\overline{X(k)}=X(\A_k)$ pour $n\geq 5$ (cf.~\cite{harduke}).
Les deux assertions sont connues lorsque $X(k)\neq\emptyset$, de sorte que seul le principe de Hasse
reste à étudier. Voici ce que l'on en sait.
En~1959, Mordell~\cite{mordell12}
démontra le principe de Hasse pour~$X$ lorsque $n \geq 12$ et $k=\Q$.
En~1964, Swinnerton-Dyer~\cite{swd10} l'établit pour $n \geq 10$ et $k=\Q$
(cf.~également~\cite[Remark~10.5.2]{ctsansdii}).
En~1971, Cook~\cite{cook} y parvint pour $n \geq 8$ et $k=\Q$ lorsque~$X$ est une intersection lisse
de deux quadriques \emph{simultanément diagonales} dans~$\P^n_k$.
En~1987, pour~$k$ arbitraire, Colliot-Thélène, Sansuc et Swinnerton-Dyer obtinrent le principe
de Hasse pour~$X$ lorsque $n \geq 8$
ainsi que dans quelques cas particuliers concernant
les intersections de deux quadriques
qui contiennent soit deux droites gauches conjuguées, soit une quadrique de dimension~$2$
(cf.~\cite{ctsansdi}, \cite{ctsansdii}).
Par la suite, Debbache (non publié) démontra le principe de Hasse pour~$X$
lorsque~$n\geq 7$ et que~$X$ est inclus dans un cône quadratique
et Salberger (non publié) prouva
que $X(\A_k)^{\Br(X)}\neq\emptyset \Rightarrow X(k)\neq\emptyset$
pour $n \geq 4$
si~$X$ contient une conique (ce qui revient, pour $n=4$, à demander que~$X$ admette une structure de fibré en coniques sur~$\P^1_k$).
Si $n=4$ et que le groupe de Picard de~$X$ est de rang~$1$ (c'est-à-dire lorsque~$X$ n'est ni l'éclaté d'une surface de del Pezzo de degré~$\geq 5$
(auquel cas il n'y aurait rien à démontrer) et que~$X$ n'est pas un fibré en coniques sur~$\P^1_k$),
le seul autre résultat connu est celui, déjà cité, de Colliot-Thélène, Skorobogatov et
Swinnerton-Dyer, qui établissent le principe de Hasse pour les intersections
suffisamment générales de deux quadriques \emph{simultanément diagonales} dans~$\P^4_k$,
sous l'hypothèse de Schinzel et la finitude des groupes de Tate-Shafarevich.

Les deux résultats principaux du présent travail sont les suivants.

\bigskip
\begin{theoreme*}[(cf.~théorème~\ref{ch3introp5}) ]%
Admettons l'hypothèse de Schinzel et la finitude des groupes de Tate-Shafarevich des courbes
elliptiques sur les corps de nombres.  Soit $n \geq 5$.
Toute intersection lisse de deux quadriques dans~$\P^n_k$ satisfait au principe de Hasse.
\end{theoreme*}

\bigskip
Si~$X$ est une surface de del Pezzo de degré~$4$ sur~$k$, choisissons
des formes quadratiques homogènes~$q_1$ et~$q_2$ en cinq variables telles
que~$X$ soit isomorphe à la variété projective d'équations $q_1=q_2=0$.
Numérotons les cinq racines du polynôme homogène $f(\lambda,\mu)=\det(\lambda q_1 + \mu q_2)$
dans un corps de décomposition $k'/k$
et notons $G \subset \mathfrak{S}_5$ le groupe de Galois de~$k'$ sur~$k$.
La classe de conjugaison du sous-groupe $G\subset \mathfrak{S}_5$ ne dépend
que de~$X$.  De plus, il est possible de diagonaliser d'autant plus de variables
simultanément dans~$q_1$ et~$q_2$ que le polynôme~$f$ admet de racines $k$\nobreakdash-rationnelles.
À l'extrême, le polynôme~$f$ est scindé si et seulement s'il existe une base
dans laquelle~$q_1$ et~$q_2$ sont simultanément diagonales.

\bigskip
\begin{theoreme*}[(cf.~théorème~\ref{ch3introthdp4}) ]%
Admettons l'hypothèse de Schinzel et la finitude des groupes de Tate-Shafarevich des courbes
elliptiques sur les corps de nombres.  Soit~$X$ une surface de del Pezzo de degré~$4$ sur~$k$.
Dans chacun des cas suivants, $X$ satisfait au principe de Hasse:
\begin{itemize}
\item le sous-groupe $G \subset \mathfrak{S}_5$ est $3$\nobreakdash-transitif (\emph{i.e.} $G=\mathfrak{A}_5$
ou $G=\mathfrak{S}_5$);
\item le polynôme~$f$ admet exactement deux racines $k$\nobreakdash-rationnelles et d'autre part $\Br(X)/\Br(k)=0$;
\item le polynôme~$f$ est scindé et $\Br(X)/\Br(k)=0$.
\end{itemize}
\end{theoreme*}

\bigskip
Ce théorème est le premier résultat positif concernant l'arithmétique des
surfaces de del Pezzo de degré~$4$ générales
(l'égalité $G=\mathfrak{S}_5$ est en effet satisfaite «~en général~»).
Néanmoins son intérêt ne se limite pas à de telles surfaces:
même dans le cas simultanément diagonal, ce théorème
généralise strictement le résultat de Colliot-Thélène, Skorobogatov et Swinnerton-Dyer mentionné précédemment.

Les conditions suffisantes qui figurent dans l'énoncé du théorème ne sont
qu'un échantillon; en réalité, pour toute classe de conjugaison~$\Crond$
de sous-groupes de $\mathfrak{S}_5$, nous obtenons le principe de Hasse
dès que~$X$ est suffisamment générale (en un sens explicite) parmi les surfaces
de del Pezzo de degré~$4$ pour lesquelles $G \in \Crond$.

Nous déduisons le premier théorème du second.
La preuve de celui-ci repose entre autres sur une généralisation convenable de la méthode introduite par Swinnerton-Dyer~\cite{swdegloffstein} et développée
par Colliot-Thélène, Skorobogatov et Swinnerton-Dyer~\cite{css} permettant d'étudier l'arithmétique des surfaces munies d'un pinceau de courbes de genre~$1$,
sur une construction suggérée par Swinnerton-Dyer~\cite[§6]{bsd} et sur un théorème récent d'Harari~\cite{harprepub}.

\chapter*{Conventions\markboth{Conventions}{Conventions}}
\addcontentsline{toc}{chapter}{Conventions}

Si~$M$ est un groupe abélien (ou un objet en groupes abéliens dans une certaine catégorie)
et~$n$ un entier naturel, on note respectivement $\tors{n}M$ et $M/n$ le noyau et le conoyau
de l'endomorphisme de multiplication par~$n$.  Plus généralement, si $f \colon M \rightarrow N$
est un homomorphisme de groupes, on note parfois $\tors{f}M$ le noyau de~$f$.
Si~$p$ est un nombre premier, $M\{p\}$ désigne le sous-groupe
de torsion $p$\nobreakdash-primaire de~$M$.  Le sous-groupe de~$M$ engendré par $\uplet{x_1}{x_n}$
est noté $\langle\uplet{x_1}{x_n}\rangle$; les~$x_i$ peuvent être des éléments ou des parties de~$M$.

Si~$k$ est un corps et $k'$, $k''$ sont deux extensions quadratiques ou triviales de~$k$
(ou plus généralement multiquadratiques), on note~$k'k''$ l'extension composée de $k'/k$
et~$k''/k$ (dans un corps arbitraire contenant~$k'$ et~$k''$);
elle est bien définie à isomorphisme près.  Étant donnés un corps~$k$ et
une propriété~$(\mathrm{P})$ des extensions de~$k$, on dira qu'une extension $\ell/k$ est
\emph{la plus petite extension satisfaisant à~$(\mathrm{P})$} si elle satisfait à~$(\mathrm{P})$ et
si d'autre part elle se plonge dans toute extension de~$k$ satisfaisant à~$(\mathrm{P})$.

Soit~$X$ un schéma.  Le corps résiduel d'un point $x \in X$ est
noté~$\kappa(x)$.
Si~$X$ est intègre, le corps des fonctions rationnelles sur~$X$ est
noté~$\kappa(X)$.
On emploiera parfois la notation $\{x\}$ pour désigner le schéma $\Spec(\kappa(x))$.
Un point $x \in X$ est régulier si l'anneau local $\Orond_{X,x}$ est régulier;
il est \emph{singulier} sinon. De même, le schéma~$X$ est dit \emph{singulier} s'il n'est
pas régulier.  Si $R \rightarrow S$ est un morphisme d'anneaux (que l'on suppose toujours
commutatifs et unitaires) et si~$X$ est un $R$\nobreakdash-schéma, on notera $X \otimes_R S$
ou parfois $X_S$ (si aucune confusion n'est possible) le $S$\nobreakdash-schéma $X\times_{\Spec(R)}\Spec(S)$.
Si $f \colon X \rightarrow S$ est un morphisme de schémas et si $s \in S$, on notera
$X_s = f^{-1}(s)$ la fibre de~$f$ en~$s$.
L'ensemble des points de codimension~$1$ d'un schéma~$X$ est
noté~$X^{(1)}$.
Une \emph{variété} sur un corps~$k$ est, par définition, un $k$\nobreakdash-schéma de type fini.
Une \emph{courbe} (resp.~\emph{surface}) est une variété de dimension~$1$ (resp.~$2$).
Une variété~$X$ est dite \emph{rationnelle} si elle devient birationnelle à l'espace
projectif après une extension des scalaires; pour signifier que~$X$ est birationnelle
à l'espace projectif sur le corps de base, on dira qu'elle est $k$\nobreakdash-rationnelle.
Enfin, lorsqu'on parlera de pinceaux sur une variété, il s'agira de pinceaux
linéaires (c'est-à-dire de systèmes linéaires de dimension~$1$).

Pour tout schéma intègre~$S$, on note $0,1,\infty \in \P^1_S$ les
points de $\P^1(\kappa(S))$ de coordonnées homogènes respectives
$[1:0]$, $[1:1]$, $[0:1]$, ce qui détermine une immersion ouverte
$\A^1_S \hookrightarrow \P^1_S$ fonctorielle en~$S$ ainsi qu'un
isomorphisme $\kappa(\P^1_k)=k(t)$ pour tout corps~$k$.

L'expression «~faisceau étale sur~$X$~» désigne un faisceau sur le petit
site étale du schéma~$X$.  La cohomologie employée sera toujours la cohomologie étale
(ou la cohomologie galoisienne), de sorte qu'on ne la précisera pas en indice.
Si $f \colon X \rightarrow Y$ est un morphisme de schémas et~$\Frond$ un
faisceau étale sur~$Y$, on notera parfois $H^n(X,\Frond)$ le groupe
$H^n(X,f^\star \Frond)$ (en particulier lorsque~$f$ est une immersion ouverte).
Le foncteur de Picard relatif d'un morphisme
de schémas $f \colon X \rightarrow Y$ est noté~$\PPic_{X/Y}$; on ne le considérera
que comme un faisceau étale sur~$Y$, et en tant que tel, il coïncide avec
$R^1 f_\star \Gm$
calculé pour la topologie étale, lorsque~$f$ est propre (cf.~\cite[p.~203]{blr}).
Si~$S$ est un schéma et $G$ un $S$\nobreakdash-schéma en groupes, on réserve
le terme «~torseur~» pour désigner les espaces principaux homogènes sous~$G$
qui sont des $S$\nobreakdash-schémas.  Enfin, si~$G$ est un schéma en groupes lisse et
de type fini sur~$S$, on appelle \emph{composante neutre de~$G$} l'unique
sous-schéma en groupes ouvert de~$G$ à fibres connexes (cf.~\cite[15.6.5]{ega43}).

Le \emph{groupe de Brauer}\index{groupe de Brauer} d'un schéma~$X$ (ou d'un anneau~$X$), noté $\Br(X)$\glossary{$\Br(X)$, $\Br_1(X)$}, est son groupe de
Brauer cohomologique $H^2(X,\Gm)$ (cf.~\cite{grothbr2}, \cite{grothbr3}).
Nous renvoyons à~\cite[§1]{cscrelle94} pour les propriétés de base du groupe de Brauer
et pour les notions associées de ramification et de résidus.
Si~$X$ est une variété sur~$k$, nous commettrons l'abus consistant à noter $\Br(X)/\Br(k)$
le conoyau de la flèche natu\-relle $\Br(k) \rightarrow \Br(X)$, qui n'est pas nécessairement
injective.  Étant donné un morphisme de schémas $f \colon X \rightarrow Y$ avec~$Y$ intègre,
le \emph{groupe de Brauer vertical\index{groupe de Brauer!vertical} de~$X$}, noté
$\Br_\vert(X)$\glossary{$\Br_\vert(X)$, $\Br_\hor(X)$}, est le sous-groupe
de~$\Br(X)$ constitué des classes dont la restriction à la fibre générique de~$f$
appartient à l'image de la flèche naturelle $\Br(\eta)\rightarrow\Br(X_\eta)$,
où~$\eta$ est le point générique de~$Y$; et le \emph{groupe de Brauer horizontal\index{groupe de Brauer!horizontal} de~$X$}
est $\Br_\hor(X)=\Br(X)/\Br_\vert(X)$.
Il existe donc une injection
naturelle $\Br_\hor(X) \hookrightarrow \Br_\hor(X_\eta)=\Br(X_\eta)/\Br(\eta)$.
Soient enfin~$X$ une variété sur un corps~$k$ de caractéristique~$0$ et~$\ksep$ une
clôture algébrique de~$k$.  Le \emph{groupe de Brauer\index{groupe de Brauer!algébrique} algébrique} de~$X$,
noté $\Br_1(X)$,
est le noyau de l'application naturelle $\Br(X) \rightarrow \Br(X_\ksep)$.
Par ailleurs, si~$X$ est lisse et connexe, une classe $A \in \Br(\kappa(X))$ sera dite
\emph{géométriquement non ramifiée} si son image dans $\Br(\kappa(X) \otimes_k \ksep)$
appartient au sous-groupe $\Br(X_\ksep)$.  Ces deux notions ne dépendent pas du choix de~$\ksep$.

Soient~$k$ un corps de nombres et~$X$ une variété sur~$k$. Notons~$\Omega$ l'ensemble
des places de~$k$ et~$\A_k$ l'anneau de ses adèles,
de sorte que $X(\A_k)=\prod_{v \in \Omega}X(k_v)$ si~$X$ est propre sur~$k$.
On dit que~\emph{$X$ satisfait au principe de Hasse} (resp.~\emph{à l'approximation faible})
si $X(\A_k)\neq \emptyset \Rightarrow X(k)\neq\emptyset$ (resp.~si $X(k)$ est dense dans~$X(\A_k)$
pour la topologie adélique). Pour $B \subset \Br(X)$, on note $X(\A_k)^B$\glossary{$X(\A_k)^B$,
$X(\A_k)^\Br$} l'ensemble
des $(\fp{P}v)_{v \in \Omega} \in X(\A_k)$ tels que $\sum_{v \in \Omega} \inv_v A(\fp{P}v)=0$
pour tout \mbox{$A \in B$}, où $\inv_v \colon \Br(k_v) \hookrightarrow \Q/\Z$ désigne l'invariant
de la théorie du corps de classes local. On pose
de plus $X(\A_k)^{\Br}=X(\A_k)^{\Br(X)}$\glossary{$X(\A_k)^{\Br_\vert}$, $X(\A_k)^{\Br_1}$}, $X(\A_k)^{\Br_1}=X(\A_k)^{\Br_1(X)}$,
et si l'on dispose d'un morphisme $f \colon X \rightarrow Y$
avec~$Y$ intègre, $X(\A_k)^{\Br_\vert}=X(\A_k)^{\Br_\vert(X)}$.
On dit qu'il y a obstruction de Brauer-Manin à l'existence d'un point rationnel
(resp.~à l'approximation faible) sur~$X$ si $X(\A_k)^\Br=\emptyset$
(resp.~si $X(\A_k)^\Br \neq X(\A_k)$). On définit de même l'obstruction de Brauer-Manin
algébrique (resp.~verticale, le cas échéant) en remplaçant~$\Br$ par~$\Br_1$ (resp.~par $\Br_\vert$).
Il résulte de la théorie du corps de classes global que ce sont bien des obstructions
(cf.~\cite[§3]{cscrelle94}).

Nous ferons souvent appel au théorème des fonctions\index{théorème des fonctions implicites} implicites; nous entendons par là \cite[Theorem~2, Ch.~III, §9]{serreharvard},
qui s'applique sur tout corps muni d'une valeur absolue pour laquelle il est complet.

Un autre outil dont nous aurons souvent besoin est un théorème dû à Harari et communément appelé\index{lemme formel} «~lemme formel~».
L'énoncé auquel cette expression fera référence est le suivant. Il s'agit d'une variante de~\cite[Corollaire~2.6.1]{harduke} que nous n'avons pas trouvée dans la littérature;
par souci de complétude et pour la commodité du lecteur, nous en donnons une preuve, mais celle-ci est directement adaptée de celle donnée par Harari dans~\cite{harduke}.

\bigskip
\begin{theoreme*}[(«~lemme formel~») ]%
Soient~$k$ un corps de nombres et~$X$ une variété propre, lisse et géométriquement
connexe sur~$k$.  Notons~$\Omega$ l'ensemble des places de~$k$ et~$\A_k$ l'anneau
de ses adèles. Pour tout ouvert $U \subset X$, tout sous-groupe fini $B \subset \Br(U)$,
tout sous-ensemble fini $S \subset \Omega$ et toute famille $(\fp{P}v)_{v \in \Omega} \in X(\A_k)^{\Br(X) \cap B}$
telle que $(\fp{P}v)_{v\in S} \in \prod_{v \in S}U(k_v)$,
il existe un sous-ensemble fini $S_1 \subset \Omega$ contenant~$S$ et
une famille $(Q_v)_{v \in S_1} \in \prod_{v \in S_1}U(k_v)$ tels
que $Q_v=\fp{P}v$ pour tout $v \in S$ et
$\sum_{v \in S_1} \inv_v A(Q_v)=0$ pour tout $A \in B$.
\end{theoreme*}

\bigskip
\begin{demo}
Notons $B^\star=\Hom(B,\Q/\Z)$ le dual de Pontrjagin de~$B$ et $\phi_v \colon U(k_v) \rightarrow B^\star$, pour $v \in \Omega$, l'application qui à $Q \in U(k_v)$
associe le caractère $A \mapsto \inv_v A(Q)$.
Soit $\Gamma \subset B^\star$ le sous-groupe engendré par l'ensemble des éléments de~$B^\star$ qui appartiennent à l'image de~$\phi_v$ pour une
infinité de~$v$.

La finitude du groupe $B \cap \Br(X)$, le théorème des fonctions implicites et la continuité de l'évaluation des classes de~$\Br(X)$
sur~$X(k_v)$
permettent de supposer que~$\fp{P}v$ appartient à~$U(k_v)$ pour tout $v \in \Omega$.
Quitte à agrandir~$S$, on peut aussi supposer
que pour $v \in \Omega \setminus S$, l'application $\phi_v$ est à valeurs dans~$\Gamma$.

Notons $w=\sum_{v \in S} \phi_v(\fp{P}v) \in B^\star$ et prouvons que $w \in \Gamma$.
La dualité de \mbox{Pontrjagin} des groupes abéliens finis étant une dualité parfaite,
il suffit de vérifier que l'orthogonal de~$\Gamma$ pour l'accouplement canonique
$B \times B^\star \rightarrow \Q/\Z$
est inclus dans l'orthogonal de~$w$.
Si $A \in B$ est orthogonal à~$\Gamma$,
on a $\inv_v A(Q)=0$ pour tout $v \in \Omega \setminus S$ et tout $Q \in U(k_v)$.
Par un théorème d'Harari, cela entraîne que $A \in \Br(X)$
(cf.~\cite[Théorème~2.1.1]{harduke} et \cite[Théorème~1.3]{ctbudapest}).
Il s'ensuit que~$A$ est bien orthogonal à~$w$, compte tenu
que $(\fp{P}v)_{v \in \Omega} \in X(\A_k)^{\Br(X)\cap B}$ et que $\phi_v(P_v)\in\Gamma$ pour $v \in \Omega\setminus S$.

L'appartenance de~$w$ (ou plutôt de~$-w$) à~$\Gamma$ se traduit par l'existence
d'un ensemble fini $T \subset \Omega \setminus S$ et d'une famille $(Q_v)_{v \in T} \in \prod_{v \in T}U(k_v)$
tels que $-w=\sum_{v \in T} \phi_v(Q_v)$.
Posant $S_1=S \cup T$ et $Q_v=\fp{P}v$ pour $v \in S$, on obtient ainsi une famille $(Q_v)_{v \in S_1}\in \prod_{v \in S_1}U(k_v)$
qui vérifie la conclusion du théorème.
\end{demo}

\bigskip
Enfin, si~$E$ et~$E'$ sont des courbes elliptiques sur un corps de nombres~$k$
et $f \colon E \rightarrow E'$ est une isogénie, on convient d'appeler \emph{groupe
de $f$\nobreakdash-Selmer de~$E$}, et de noter $\Sel_f(k,E)$,
le\index{groupe de Selmer} sous-groupe de $H^1(k,\tors{f}E)$ constitué des classes
dont l'image dans $H^1(k_v,E)$ est nulle pour toute place~$v$ de~$k$.
Cette définition semble être la plus répandue, mais il y a lieu de la préciser
car elle n'est pas standard; certains auteurs nomment ce groupe
le groupe de $f$\nobreakdash-Selmer de~$E'$ (ce qui se justifie par
l'interprétation des éléments de $H^1(k,\tors{f}E)$
comme des $f$\nobreakdash-revêtements de~$E'$).

\chapternotoc[Arithmétique des pinceaux semi-stables I]{Arithmétique des pinceaux semi-stables de courbes de genre~$1$ dont
les jacobiennes ont leur $2$\nobreakdash-torsion rationnelle}
\addcontentsline{toc}{chapter}{\chaptername\ \thechapter.\ Arithmétique des pinceaux semi-stables de courbes de genre~$1$ (première partie)}%

\section{Introduction}

Nous nous intéressons dans ce chapitre aux questions d'existence et de Zariski-densité
des points rationnels pour certaines surfaces propres et lisses, fibrées en courbes de
genre~$1$ au-dessus de la droite projective, définies sur un corps de nombres.
Le tout premier théorème significatif concernant une question de ce type
fut démontré par Swinnerton-Dyer~\cite{swdegloffstein}.
C'est une famille spécifique de surfaces, de surcroît définies sur~$\Q$, que Swinnerton-Dyer
étudiait.
Ses techniques furent développées par
Colliot-Thélène, Skorobogatov et Swinnerton-Dyer, qui
obtinrent ainsi dans~\cite{css} les premiers résultats généraux d'existence
de points rationnels pour des surfaces fibrées en courbes de genre~$1$ au-dessus de la droite
projective.
Un certain nombre
d'hypothèses apparaissent néanmoins dans leur théorème principal, dont notamment les
trois suivantes:
les fibres géométriques singulières de la fibration considérée sont
des réunions de deux courbes rationnelles intègres se rencontrant transversalement
en deux points distincts,
tous les points d'ordre~$2$ de la jacobienne de la fibre générique sont rationnels,
et soit la fibration ne possède pas de section et
le rang de Mordell-Weil de la jacobienne de la fibre générique est nul,
soit la fibration possède une section et le rang de la fibre générique est exactement~$1$.

Les idées introduites dans~\cite{css} furent réutilisées à plusieurs
reprises dans des situations où le théorème principal de~\cite{css} ne s'appliquait pas.
Swinnerton-Dyer~\cite{sddiagquartic} a en particulier obtenu des conditions suffisantes explicites
génériquement vraies
pour qu'une surface quartique diagonale
$$
a_0 x_0^4 + a_1 x_1^4 + a_2 x_2^4 + a_3 x_3^4 = 0
$$
avec~$\uplet{a_0}{a_3}\in\Q$ tels que~$a_0a_1a_2a_3$ soit un carré
dans~$\Q^\star$ satisfasse au principe de Hasse, en admettant l'hypothèse de Schinzel et
la finitude des groupes de Tate-Shafarevich des courbes elliptiques sur~$\Q$.  Pour ce faire,
il s'est servi de l'existence sur de telles surfaces de pinceaux de courbes de genre~$1$
dont la jacobienne générique vérifie les conditions ci-dessus (rang de Mordell-Weil nul
et points d'ordre~$2$ rationnels) mais dont les fibres singulières géométriques possèdent
quatre composantes irréductibles, réduites et
organisées en quadrilatère avec intersections transverses.

Seules certaines configurations de composantes irréductibles peuvent apparaître dans les
fibres géométriques des pinceaux de courbes de genre~$1$ dont l'espace total est régulier.
Dans le cas d'une fibration relativement minimale dont chaque fibre possède une composante
irréductible de multiplicité~$1$, les configurations possibles ont été classifiées
par Kodaira et Néron (cf.~\cite{kodaira}, \cite{neronihes}, \cite{silvaec2}).
Il apparaît dans cette classification qu'un rôle crucial est joué par la propriété
qu'a la jacobienne de la fibre générique d'être à réduction semi-stable (c'est-à-dire
multiplicative ou bonne) ou non (réduction additive).
Lorsqu'elle est à réduction semi-stable, les fibres géométriques singulières
sont réduites et leurs composantes irréductibles sont des courbes rationnelles
organisées en polygone à~$n$ côtés, avec intersections transverses.  Le théorème
principal de~\cite{css} traite le cas~$n=2$, et Swinnerton-Dyer dans~\cite{sddiagquartic}
applique la méthode à un cas particulier pour lequel~$n=4$.

L'un des buts du présent chapitre est d'obtenir des résultats dans le
cas général de réduction semi-stable, sans hypothèse sur~$n$, afin
notamment de couvrir simultanément les théorèmes principaux de~\cite{css}
et de~\cite{sddiagquartic}.  Nous en profitons pour supprimer l'hypothèse
sur le rang générique, qui s'avère inessentielle.
Ainsi prouvons-nous:

\bigskip
\begin{theoreme}[ (cf.~théorème~\ref{ch1appobm2})]\label{ch1introth1}Soient~$k$
un corps de nombres et~$X$ une surface propre, lisse et connexe sur~$k$
munie d'un morphisme~$\pi \colon X \rightarrow \P^1_k$ dont la fibre générique~$X_\eta$
est une courbe lisse et géométriquement connexe de genre~$1$ et de période~$2$,
et dont les fibres sont toutes réduites.
Supposons que la jacobienne de~$X_\eta$ soit à réduction
semi-stable en tout point fermé de~$\P^1_k$, que sa $2$\nobreakdash-torsion soit
rationnelle et que~$\pi$ satisfasse à la \condD{}.
Admettons l'hypothèse de Schinzel et la finitude des groupes de Tate-Shafarevich des
courbes elliptiques sur~$k$. Alors $X(k)\neq\emptyset$
dès que $X(\A_k)^\Br \neq \emptyset$ (et~$X(k)$ est même Zariski-dense dans~$X$).
\end{theoreme}

\bigskip
La \condD{} est une hypothèse technique étroitement liée au groupe de Brauer
de~$X$. Elle sera définie précisément au paragraphe~\ref{ch1hypnot}.  Qu'il suffise
pour le moment de préciser qu'une hypothèse analogue apparaît dans~\cite{css}
et dans~\cite{sddiagquartic}, et que celle-ci en est une généralisation commune.

L'hypothèse de Schinzel et la finitude des groupes de Tate-Shafarevich sont deux
conjectures qui semblent actuellement incontournables pour obtenir des résultats
généraux sur les points rationnels d'une surface fibrée en courbes de genre~$1$
sur un corps de nombres
au moyen des techniques introduites dans~\cite{css}.
Au mieux peut-on espérer remplacer dans certains cas l'hypothèse de Schinzel
par le théorème de Dirichlet sur les nombres premiers dans une progression arithmétique,
comme l'ont fait Swinnerton-Dyer~\cite{sdcubdiag} pour les surfaces cubiques diagonales
et Skorobogatov et Swinnerton-Dyer~\cite{ssdkummer} pour certaines surfaces de Kummer.

\bigskip
Le théorème~\ref{ch1introth1} est obtenu comme corollaire du
théorème~\ref{ch1thprin}, d'énoncé plus technique.
Ce chapitre contient deux autres applications du théorème~\ref{ch1thprin}.

\bigskip
La première concerne les pinceaux de courbes de genre~$1$ qui ne satisfont pas à
la\index{condition (D)@\condD{}} \condD.
Soit~$\pi \colon X \rightarrow \P^1_k$ un pinceau vérifiant toutes les hypothèses du
théorème~\ref{ch1introth1} autres que la \condD.
Admettant l'hypothèse de Schinzel, nous prouvons au paragraphe~\ref{ch1parsecdesc}
(théorème~\ref{ch1thsecdesc})
qu'en l'absence d'obstruction de Brauer-Manin à l'approximation faible sur~$X$
et moyennant une légère hypothèse technique (qui est en tout cas
toujours satisfaite lorsque les fibres géométriques singulières de~$\pi$ possèdent
deux composantes irréductibles, c'est-à-dire dans la situation envisagée dans~\cite{css};
il s'agit de l'hypothèse «~$L_M=\kappa(M)$ pour tout $M \in \Mrond$~»
du théorème~\ref{ch1thsecdesc}),
il existe nécessairement des points rationnels~$x \in \P^1(k)$ au-dessus desquels
la fibre~$X_x$ est lisse, possède un~$k_v$\nobreakdash-point pour toute place~$v$ de~$k$
et vérifie la condition suivante: notant~$E_x$ la jacobienne
de~$X_x$, la classe de~$X_x$ dans le groupe de Tate-Shafarevich~$\Sha(E_x)$ est orthogonale
au sous-groupe~$\tors{2}\Sha(E_x)$ pour l'accouplement\index{accouplement de Cassels-Tate} de Cassels-Tate; et l'ensemble
de ces~$x \in \P^1(k)$ est même dense dans~$\P^1(\A_k)$ pour la topologie adélique.

L'intérêt de ce résultat est double.  Tout d'abord, il a une conséquence
concrète concernant l'approximation faible.  Swinnerton-Dyer donne
dans~\cite[§8]{sddiagquartic} un exemple de famille de surfaces~$K3$ pour
lesquelles il exhibe un défaut d'approximation faible à l'aide d'une seconde
descente sur des pinceaux de courbes de genre~$1$, cette seconde\index{seconde descente} descente
étant effectuée explicitement grâce à l'algorithme de Cassels~\cite{casselssecond}.
Notre résultat montre que lorsque l'hypothèse technique du théorème~\ref{ch1thsecdesc}
est satisfaite,
un défaut d'approximation faible mis en évidence par cette méthode est nécessairement
expliqué par l'obstruction de Brauer-Manin, si l'on admet l'hypothèse de Schinzel.
Il se trouve que l'hypothèse technique du théorème~\ref{ch1thsecdesc} n'est pas vérifiée
dans l'exemple de~\cite[§8]{sddiagquartic}, mais la méthode de Swinnerton-Dyer est
toute générale et il ne fait pas de doute qu'elle conduit aussi à des défauts d'approximation
faible notamment pour des fibrations du type envisagé dans~\cite{css}.

Nous pensons cependant que l'intérêt du théorème~\ref{ch1thsecdesc} réside plus dans
sa preuve que dans son énoncé.
La \condD{} implique que le sous-groupe de torsion $2$\nobreakdash-primaire du groupe de Brauer de~$X$ est
entièrement vertical,
et la véritable raison pour laquelle jusqu'à présent
cette condition ou une condition analogue est apparue dans tous les articles
élaborant les idées de~\cite{css} est que l'on ne sait pas utiliser l'hypothèse
qu'une classe de~$\Br(X)$ ne fournit pas d'obstruction de Brauer-Manin si cette
classe n'est pas verticale.  Le théorème~\ref{ch1thsecdesc}, qui s'applique que
la \condD{} soit satisfaite ou non, représente donc un premier pas dans cette
direction: sa preuve prend effectivement en compte les classes non verticales de~$\Br(X)$.

\bigskip
Nous appliquons enfin le théorème~\ref{ch1thprin} à la notion de courbe
elliptique «~de rang élevé~» sur~$k(t)$, où~$k$ est un\index{courbe elliptique de rang élevé} corps de nombres.
Ce sont les courbes elliptiques~$E/k(t)$ dont le rang de Mordell-Weil est strictement
inférieur à celui de la spécialisation~$E_x/k$ en tout point rationnel~$x \in \P^1(k)$
hors d'un ensemble fini. On n'en connaît aucun exemple inconditionnel;
admettant la conjecture de parité pour les courbes elliptiques,
Cassels et Schinzel~\cite{cassschin} puis Rohrlich~\cite{rohrlich}
ont donné des exemples de courbes elliptiques de rang élevé sur~$\Q(t)$.
Leurs exemples sont tous isotriviaux.
Conrad, Conrad et Helfgott~\cite{cchrootnumbers} ont récemment
montré qu'une courbe elliptique de rang élevé sur~$k(t)$ est nécessairement isotriviale
lorsque~$k=\Q$, en admettant trois conjectures arithmétiques plus ou moins classiques,
dont une conjecture de densité concernant les rangs des courbes elliptiques~$E_x$ lorsque~$x$
varie.
Nous déduisons du théorème~\ref{ch1thprin} un résultat similaire, également conditionnel,
nécessitant de plus quelques hypothèses sur les courbes elliptiques considérées,
mais valable sur tout corps de nombres et ne dépendant d'aucune
conjecture qui concerne les courbes elliptiques (cf.~théorème~\ref{ch1rangeleve}).

\bigskip
Le plan du chapitre est le suivant.
Nous fixons les notations au paragraphe~\ref{ch1hypnot}, nous y
énonçons le théorème~\ref{ch1thprin}, nous y décrivons les grandes étapes
de sa démonstration et nous y soulignons les difficultés qui n'existaient pas
dans la situation de~\cite{css}.
Pour les applications, il est important
de disposer d'une forme explicite de la \condD{}; nous en donnons une
au paragraphe~\ref{ch1parexplicited}.
La preuve du théorème~\ref{ch1thprin} occupe les paragraphes~\ref{ch1parpreuvedebut}
et~\ref{ch1parpreuvefin}.
Au paragraphe~\ref{ch1parcondd},
nous discutons les liens entre la \condD{}, le
groupe de Brauer et le groupe~$\Drond$ introduit dans~\cite[§4]{css}.
Ce groupe, défini en toute généralité, permit aux auteurs de~\cite{css}
de formuler une condition abstraite qui équivaut à la \condD{} sous les hypothèses
de leur théorème principal (voir~\cite[§4.7]{css}).
Nous verrons que cette équivalence n'est plus valable dans la situation générale
de réduction semi-stable; c'est précisément de cette difficulté que naît
l'hypothèse technique du théorème~\ref{ch1thsecdesc} sur la seconde descente.
Nous consacrons ensuite le paragraphe~\ref{ch1parapprat} aux applications
du théorème~\ref{ch1thprin} à l'existence de points rationnels:
d'une part nous établissons le théorème~\ref{ch1introth1}, qui généralise
le théorème principal de~\cite{css}, et d'autre part nous déduisons du théorème~\ref{ch1introth1}
les résultats de Swinnerton-Dyer sur les surfaces quartiques diagonales énoncés
dans~\cite[§3]{sddiagquartic}.
Enfin, les paragraphes~\ref{ch1parsecdesc} et~\ref{ch1parrangel} contiennent
respectivement l'application aux défauts d'approximation
faible mis en évidence par une seconde descente, et l'application
aux courbes elliptiques de rang élevé.
Ils font tous deux usage des résultats obtenus au paragraphe~\ref{ch1parcondd}.

\section{Hypothèses et notations}
\label{ch1hypnot}

Soient~$k$ un corps de caractéristique~$0$ et~$C$ un schéma de Dedekind connexe sur~$k$,
de point générique~$\eta$\glossary{$k$, $C$, $\eta$, $X$, $\pi$}.
Supposons donnée une surface~$X$ lisse et géométriquement connexe sur~$k$, munie
d'un morphisme propre et plat $\pi \colon X \rightarrow C$
dont la fibre générique~$X_\eta$
est une courbe lisse de genre~$1$ sur~$K=\kappa(C)$ et dont toutes les fibres sont réduites.
Supposons de plus que la période de la courbe~$X_\eta$ divise~$2$, c'est-à-dire
que la classe de $H^1(K, E_\eta)$ définie par le torseur~$X_\eta$ est tuée par~$2$,
en notant~$E_\eta$ la jacobienne de~$X_\eta$.  Supposons enfin que la courbe
elliptique~$E_\eta$ soit à réduction semi-stable en tout point fermé de~$C$ et
que ses points d'ordre~$2$ soient tous $K$\nobreakdash-rationnels.

Notons~$\Erond$\glossary{$E_\eta$, $\Erond$, $\Erond^0$} le modèle de Néron\index{modèle de Néron} de~$E_\eta$ sur~$C$,
$\Erond^0 \subset \Erond$ sa composante neutre,
$\Mrond \subset C$ l'ensemble des points fermés de mauvaise réduction
pour~$E_\eta$ et~$U$\glossary{$\Mrond$, $U$}
un ouvert dense de~$C$ au-dessus duquel~$\pi$ est lisse.
Un tel ouvert est nécessairement disjoint de~$\Mrond$.
Pour~$M \in \Mrond$, notons~$\ff{F}M$\glossary{$\ff{F}M$}
le $\kappa(M)$\nobreakdash-schéma en groupes fini étale~$\Erond_M/\Erond^0_M$.

\bigskip
\begin{lemme}
\label{ch1lemmefm}
Soit~$M \in \Mrond$.  Si le~$\kappa(M)$\nobreakdash-groupe $\ff{F}M$ n'est pas constant, il existe
une extension quadratique~$L_M/\kappa(M)$ telle que $\ff{F}M \otimes_{\kappa(M)} L_M$
soit un~$L_M$\nobreakdash-groupe constant, et $\Gal(L_M/\kappa(M))$ agit alors sur~$\ff{F}M(L_M)$
par multiplication par~$-1$.
\end{lemme}

\bigskip
\begin{demo}
Ce lemme est prouvé en annexe, cf.~proposition~\ref{annpropcyclpair}.
\end{demo}

\bigskip
Lorsque~$\ff{F}M$ est constant, on notera $L_M=\kappa(M)$. Ainsi a-t-on défini dans tous
les cas une extension au plus quadratique $L_M/\kappa(M)$\glossary{$L_M$}.
Il résulte des hypothèses
de semi-stabilité de~$E_\eta$ et de $K$\nobreakdash-rationalité de ses points d'ordre~$2$
que les groupes $\ff{F}M(L_M)$ pour $M \in \Mrond$ sont cycliques d'ordre pair (cf.~proposition~\ref{annpropcyclpair}).

\bigskip
Pour~$M \in C$, notons~$\Orond_M^\sh$\glossary{$\Orond_M^\sh$, $K_M^\sh$}
l'anneau strictement local de~$C$ en~$M$
et~$K_M^\sh$ son corps des fractions.
La suite spectrale de Leray associée au faisceau étale~$E_\eta$
sur~$K$ et au morphisme canonique $\eta \rightarrow C$ montre que la
suite
\begin{equation}
\label{ch1setsg}
\myxyhook\xymatrix{
0 \ar[r] & H^1(C, \Erond) \ar[r] & H^1(K, E_\eta) \ar[r] & \displaystyle \prod_{M \in C}
H^1(K_M^\sh, E_\eta)
}
\end{equation}
est exacte.  Le groupe $H^1(C, \Erond)$ joue donc le rôle\index{groupe de Tate-Shafarevich!géométrique}
d'un \emph{groupe de Tate-Shafarevich géométrique} (cf.~\cite[§4]{grothbr3} et~\cite[§4.1]{css}).

La suite exacte de Kummer
\begin{equation}
\label{ch1sekummer}
\myxyhook\xymatrix{
0 \ar[r] & \tors{2}E_\eta \ar[r] & E_\eta \ar[r]^{\times 2} & E_\eta \ar[r] & 0
}
\end{equation}
induit une suite exacte
\begin{equation*}
\myxyhook\xymatrix{
0 \ar[r] & E_\eta(K)/2 \ar[r] & H^1(K, \tors{2}E_\eta) \ar[r]^\alpha &
\tors{2}H^1(K, E_\eta) \ar[r] & 0\rlap{\text{.}}
}
\end{equation*}
Notons $\SG_2(C,\Erond) \subset H^1(K, \tors{2}E_\eta)$\glossary{$\SG_2(C,\Erond)$}
l'image réciproque du
sous-groupe $\tors{2}H^1(C, \Erond)$ de $\tors{2}H^1(K,E_\eta)$ par~$\alpha$,
de sorte que l'on obtient une suite exacte
\begin{equation}
\label{ch1sesgtsg}
\myxyhook\xymatrix{
0 \ar[r] & E_\eta(K)/2 \ar[r] & \SG_2(C,\Erond) \ar[r] & \tors{2}H^1(C, \Erond) \ar[r] & 0
\rlap{\text{.}}
}
\end{equation}
Il y a lieu de nommer\index{groupe de Selmer!géométrique} $\SG_2(C,\Erond)$
\emph{groupe de $2$\nobreakdash-Selmer géométrique} (cf.~\cite[§4.2]{css}).

\bigskip
\begin{lemme}
\label{ch1lemmepeutevaluer}
On a l'inclusion $\SG_2(C,\Erond) \subset H^1(U,\tors{2}\Erond)$ de
sous-groupes de $H^1(K, \tors{2}E_\eta)$.
\end{lemme}

\bigskip
\begin{demo}
Combinant la suite exacte~(\ref{ch1setsg}) pour~$C=U$ avec la suite exacte obtenue
de manière analogue à partir du faisceau~$\tors{2}E_\eta$, on obtient le diagramme
commutatif
$$
\myxyhook\xymatrix{
0 \ar[r] & H^1(U, \tors{2}\Erond) \ar[d] \ar[r] & H^1(K, \tors{2}E_\eta) \ar[d]
\ar[r] & \displaystyle \prod_{M \in U} H^1(K_M^\sh, \tors{2}E_\eta) \ar[d] \\
0 \ar[r] & H^1(U, \Erond) \ar[r] & H^1(K, E_\eta) \ar[r] & \displaystyle
\prod_{M \in U} H^1(K_M^\sh, E_\eta)\rlap{\text{,}}
}
$$
dont les lignes sont exactes. Ainsi suffit-il de prouver que pour tout $M \in U$,
la flèche naturelle $H^1(K_M^\sh, \tors{2}E_\eta) \rightarrow H^1(K_M^\sh, E_\eta)$
est injective, autrement dit que $E_\eta(K_M^\sh)/2=0$;
ceci découle de la proposition~\ref{annpropekshsurdeux}, compte tenu
que~$E_\eta$ a bonne réduction en~$M$.
\end{demo}

\bigskip
Comme les fibres de~$\pi$ sont réduites, on a~$X_M(K_M^\sh)\neq \emptyset$
pour tout $M \in C$. La suite exacte~(\ref{ch1setsg}) permet d'en déduire
que la classe de $H^1(K, E_\eta)$ définie par le torseur~$X_\eta$ appartient
au sous-groupe $H^1(C, \Erond)$.  Cette classe correspond à un faisceau
représentable d'après~\cite[Theorem 4.3]{milneec}, d'où l'existence
d'un torseur $\Xrond \rightarrow C$\glossary{$\Xrond$}
sous~$\Erond$ dont la fibre générique est égale à~$X_\eta$.  Sa classe
dans $H^1(C, \Erond)$ sera notée~$[\Xrond]$.

Pour $M \in \Mrond$, notons
$\delta_M \colon H^1(C, \Erond) \rightarrow H^1(L_M, \ff{F}M)$\glossary{$\delta_M$}
la composée de la flèche induite par le morphisme de faisceaux $\Erond \rightarrow
i_{M \star} \ff{F}M$, où~$i_M$ désigne l'inclusion canonique $i_M\colon\Spec(\kappa(M)) \rightarrow C$,
et de la flèche de restriction $H^1(\kappa(M), \ff{F}M) \rightarrow H^1(L_M, \ff{F}M)$.
Soient $\TDC{C}$\glossary{$\SDC{C}$, $\TDC{C}$, $(\mathrm{D}/C)$} le noyau de l'application composée
$$
\myxyhook\xymatrix@C=8ex{
\tors{2}H^1(C, \Erond) \ar[r]^(0.41){\prod \delta_M} & \displaystyle
\prod_{M \in \Mrond} H^1(L_M, \ff{F}M) \ar[r] &
\displaystyle \prod_{M \in \Mrond} \frac{H^1(L_M, \ff{F}M)}{\langle \delta_M([\Xrond]) \rangle}
}
$$
et $\SDC{C} \subset \SG_2(C,\Erond)$ l'image réciproque de $\TDC{C}$ par la\index{condition (D)@\condD{}}
flèche de droite de la suite exacte~(\ref{ch1sesgtsg}).  On dira que \emph{la
\condDC{C} est satisfaite} si $\TDC{C}$ est engendré par~$[\Xrond]$.

Supposons maintenant que~$k$ soit un corps de nombres.  On note~$\Omega$ l'ensemble
de ses places, $\Omega_f \subset \Omega$ l'ensemble de ses places finies
et~$\A_k$ l'anneau des adèles de~$k$.  Soit\glossary{$\RA$}
$$
\RA = \bigensemble{x \in U(k)}{X_x(\A_k)\neq\emptyset}
$$
et soit~$\RDC{C}$\glossary{$\RDC{C}$, $\RD$, $\RDz$} l'ensemble des~$x \in \Rrond_A$ tels que tout élément
du groupe de $2$\nobreakdash-Selmer de~$\Erond_x$ appartienne à l'image de la composée
$$
\SDC{C} \subset \SG_2(C,\Erond) \subset H^1(U, \tors{2}\Erond) \rightarrow
H^1(k, \tors{2}\Erond_x)
$$
(dans laquelle l'avant-dernière flèche est donnée par le lemme~\ref{ch1lemmepeutevaluer}
et la dernière flèche est l'évaluation en~$x$) et tels que la restriction
de cette composée à l'image réciproque de $\{0,[\Xrond]\}$ par la flèche
de droite de~(\ref{ch1sesgtsg}) soit injective.

Lorsque~$C=\P^1_k$, on prend pour~$U$ le plus grand ouvert de~$\A^1_k$ au-dessus
duquel~$\pi$ est lisse et l'on\glossary{$\TD$, $\TDz$, $\SD$, $\SDz$}
note\glossary{\cD{}, \cDz{}} d'une part~$\RD$, $\TD$, $\SD$ et~\cD{} les ensembles~$\RDC{\P^1_k}$, $\TDC{\P^1_k}$,
$\SDC{\P^1_k}$ et la \condDC{\P^1_k}, et
d'autre part~$\RDz$, $\TDz$, $\SDz$ et~\cDz{} les ensembles~$\RDC{\A^1_k}$, $\TDC{\A^1_k}$,
$\SDC{\A^1_k}$ et la \condDC{\A^1_k},
étant entendu que $\RDC{\A^1_k}$, $\TDC{\A^1_k}$,
$\SDC{\A^1_k}$ et la \condDC{\A^1_k} désignent les ensembles et la condition
obtenus en appliquant les définitions ci-dessus avec $C=\A^1_k$
après avoir restreint~$\pi$ au-dessus de $\A^1_k \subset \P^1_k$.

\bigskip
Le théorème dont la démonstration occupera les paragraphes~\ref{ch1parpreuvedebut}
et~\ref{ch1parpreuvefin} est le suivant.

\bigskip
\begin{theoreme}\label{ch1thprin}\label{CH1THPRIN}Admettons
l'hypothèse de Schinzel. Supposons que $C=\P^1_k$ et que la fibre de~$\pi$
au-dessus du point $\infty \in \P^1(k)$ soit lisse.
Il existe alors un ensemble fini $S_0 \subset \Omega$ et un sous-groupe fini
$B_0 \subset \Br(U)$ tels que l'assertion suivante soit vérifiée.
Soient un ensemble $S_1 \subset \Omega$
fini contenant~$S_0$ et une famille $(x_v)_{v\in S_1}\in \prod_{v \in S_1} U(k_v)$.
Supposons que $X_{x_v}(k_v)\neq \emptyset$ pour tout $v \in S_1 \cap \Omega_f$ et que
$$
\sum_{v \in S_1} \inv_v A(x_v)=0
$$
pour tout $A \in B_0$. Supposons aussi que
pour toute place $v \in S_1$ réelle, on ait $X_\infty(k_v)\neq\emptyset$ et $x_v$
appartienne à la composante connexe non majorée de~$U(k_v)$.  Alors
\begin{enumerate}
\item[a)] si $\Mrond \neq \emptyset$, il existe un élément de~$\RD$ arbitrairement
proche de~$x_v$ en chaque place $v \in S_1 \cap \Omega_f$ et arbitrairement
grand en chaque place archimédienne de~$k$;
\item[b)] il existe un élément de~$\RDz$ arbitrairement
proche de~$x_v$ en chaque place $v \in S_1 \cap \Omega_f$,
arbitrairement grand en chaque place archimédienne de~$k$ et entier hors de~$S_1$.
\end{enumerate}
\end{theoreme}

\bigskip
L'appartenance d'un $x \in U(k)$ à~$\RD$ (resp.~à~$\RDz$)
est une condition arithmétique forte sur la courbe
elliptique~$\Erond_x$.  Elle entraîne par exemple que l'ordre du groupe de $2$\nobreakdash-Selmer
de~$\Erond_x$ est majoré par celui de~$\SD$ (resp.~$\SDz$),
groupe que l'on peut calculer explicitement et qui
ne dépend pas de~$x$.  Nous renvoyons le lecteur aux paragraphes~\ref{ch1parapprat}
à~\ref{ch1parrangel} pour des conséquences plus concrètes de ce théorème.

\bigskip
Le principe général de la preuve du théorème~\ref{ch1thprin} ne présente pas d'originalité par
rapport à~\cite{css}.  Contentons-nous donc de le rappeler brièvement.
Étant donné un ensemble fini de places~$T$
contenant~$S_1$, l'hypothèse de Schinzel fournit un $x \in U(k)$ satisfaisant des conditions
locales prescrites aux places de~$T$ et tel que la courbe elliptique~$\Erond_x$ ait bonne
réduction hors de~$T$ sauf en~$\Card(\Mrond)$ places, qui ne sont pas contrôlées mais en
lesquelles la mauvaise réduction de~$\Erond_x$ l'est.  Si~$B_0$ contient un système de
représentants modulo $\Br(k)$ des classes de $\Br(U)$ qui deviennent non ramifiées sur~$X$ et
si les conditions locales sur~$x$ aux places de~$T$ ont été choisies correctement, on sait
déduire de la loi de réciprocité globale l'existence de points adéliques dans la \mbox{fibre~$X_x$ --- c'est} maintenant un procédé standard, voir~\cite[Theorem 1.1]{csscrelle98}.
La difficulté de la
preuve du théorème~\ref{ch1thprin} réside dans le contrôle du groupe de $2$\nobreakdash-Selmer\index{groupe de Selmer}
de~$\Erond_x$.  Celui-ci s'effectue en trois étapes; l'hypothèse de semi-stabilité
est cruciale pour les deux dernières.  La première consiste à exprimer
le groupe de $2$\nobreakdash-Selmer de~$\Erond_x$ comme noyau d'un certain accouplement symétrique.
Dans la seconde étape,
on établit un théorème de comparaison entre ces accouplements pour
différents~$T$ et~$x$ comme ci-dessus.  Celui-ci permet, si l'on fixe une fois
pour toutes un point~$x_0$ comme ci-dessus associé à l'ensemble~$T=S_1$,
de prévoir pour d'autres choix de~$T$ et de conditions locales aux places
de~$T \setminus S_1$ ce que sera le groupe de $2$\nobreakdash-Selmer de~$\Erond_x$
quel que soit le point~$x$ donné par l'hypothèse de Schinzel.
Le but de la troisième et dernière étape est d'en déduire que l'on peut
construire un ensemble~$T$ et des conditions locales
permettant d'obtenir la conclusion du théorème.

Les difficultés liées à la généralisation du cas de réduction~$I_2$ au cas
général de réduction semi-stable sont essentiellement dues aux trois
phénomènes suivants:
\begin{itemize}
\item[(i)] les $\kappa(M)$\nobreakdash-groupes~$\ff{F}M$ ne sont plus nécessairement constants;
\item[(ii)] il n'est plus envisageable de travailler sur des modèles propres et réguliers
explicites des courbes elliptiques considérées et de leurs espaces principaux
homogènes, le nombre d'éclatements nécessaires pour obtenir de tels modèles pouvant
être arbitrairement grand;
\item[(iii)] les groupes $E_\eta(K_M^\sh)/2$ pour $M \in \Mrond$
ne sont plus
nécessairement engendrés par les classes des points d'ordre~$2$.
\end{itemize}

La structure des groupes~$\ff{F}M$ intervient dans toutes les questions
d'existence ou de non existence de points locaux sur les $2$\nobreakdash-revêtements
des courbes elliptiques~$\Erond_x$ pour $x \in U(k)$.  Qu'ils ne soient pas
constants n'est pas un problème si, pour les $x \in U(k)$ considérés,
les groupes de composantes connexes des fibres du modèle de Néron
de~$\Erond_x$ au-dessus de $\Spec(\Orond_{S_1})$ sont, eux, constants.
Il est possible de forcer cette condition à être satisfaite aux places
de $T \setminus S_1$ en choisissant convenablement ces dernières;
on peut ensuite faire en sorte que cette condition soit nécessairement satisfaite en toute place de
mauvaise réduction de~$\Erond_x$ à l'aide d'un argument de réciprocité similaire à celui permettant
de trouver des points adéliques dans les fibres de~$\pi$.

Les conséquences des phénomènes~(ii) et~(iii) sont multiples.
À titre d'exemple, l'explicitation des groupes de $2$\nobreakdash-Selmer
géométriques $\SG_2(\A^1_k,\Erond)$ et $\SG_2(\P^1_k,\Erond)$ en termes d'une
équation de Weierstrass pour~$E_\eta$ n'a plus rien d'immédiat (sans parler
des applications~$\delta_M$); voir la preuve de la proposition~\ref{ch1explsg}. Ou encore,
pour $x \in U(k)$, si l'on cherche à montrer qu'un $2$\nobreakdash-revêtement de~$\Erond_x$ n'admet
pas de $k_v$\nobreakdash-point pour un $v \in \Omega_f$ donné, on ne peut plus nécessairement
exhiber un point d'ordre~$2$ de~$\Erond_x$ dont la classe dans $\Erond_x(k_v)/2
\subset H^1(k_v,\tors{2}\Erond_x)$ n'est pas orthogonale à celle du $2$\nobreakdash-revêtement
en question (pour l'accouplement induit par l'accouplement de Weil).
On comparera ainsi la preuve du lemme~\ref{ch1lemmepeuttuer}
et celle de \cite[Lemma~2.5.2]{css}.
La différence la plus notable réside cependant dans la démonstration de l'indépendance
par rapport à~$x$ et~$T$ convenables de l'accouplement symétrique qui permet de calculer
le groupe de $2$\nobreakdash-Selmer de~$\Erond_x$: un ingrédient nouveau s'avère indispensable
(à savoir, le lemme~\ref{ch1exisgm}).

De manière générale, afin de contourner les difficultés que l'on vient d'évoquer,
nous avons remplacé dans la plupart des arguments le recours à une équation de Weierstrass
par l'utilisation des modèles de\index{modèle de Néron} Néron.

\section{Explicitation de la \condD{}}
\label{ch1parexplicited}

Supposons\index{condition (D)@\condD{}|(}
que $C=\P^1_k$ et que la courbe elliptique~$E_\eta$ ait bonne réduction
à l'infini.  Par souci de simplicité, supposons de plus que la fibration $\pi \colon X
\rightarrow \P^1_k$ soit relativement minimale.
Conformément aux conventions du paragraphe~\ref{ch1hypnot},
l'ouvert $U \subset \P^1_k$ désigne alors le complémentaire de $\Mrond \cup \{\infty\}$.
Nous reformulons ici les conditions~\cD{} et~\cDz{} sous une forme qui se prête plus facilement
aux calculs dans des cas concrets (cf.~notamment le paragraphe~\ref{ch1parquartiques}).

La courbe elliptique $E_\eta/K$ admet une équation de Weierstrass minimale de la forme
\begin{equation}
\label{ch1expleqw}
Y^2=(X-e_1)(X-e_2)(X-e_3)
\end{equation}
avec $e_1,e_2,e_3\in k[t]$. Posons $p_i=e_j-e_k$ pour toute permutation
cyclique $(i,j,k)$ de $(1,2,3)$ et $r=p_1p_2p_3$.  L'équation de Weierstrass~(\ref{ch1expleqw})
a pour discriminant $\Delta=16r^2$, de sorte que l'ensemble~$\Mrond$ s'identifie au lieu
d'annulation du polynôme~$r$ sur $\P^1_k$.
Comme dans le théorème~\ref{ch1thprin},
nous supposons que la courbe elliptique~$E_\eta$ a bonne réduction à l'infini, ce
qui revient à demander que les polynômes $p_1$, $p_2$ et~$p_3$ soient du même degré
et que ce degré soit pair.
L'hypothèse de réduction semi-stable en
tout point fermé de~$\P^1_k$ équivaut à ce que les polynômes $p_1$, $p_2$ et~$p_3$
soient deux à deux premiers entre eux.  Pour $M \in \Mrond$,
le nombre de composantes connexes de la fibre
géométrique de~$\Erond$ en~$M$ se lit alors comme l'ordre d'annulation de~$\Delta$
en~$M$ (cf.~\cite[Step~2, p.~366]{silvaec2}).

Il sera commode de noter $\Hgoth(S)=\Gm(S)/2 \times \Gm(S)/2$\glossary{$\Hgoth(S)$}
pour tout schéma ou anneau~$S$.
Fixons un $K$\nobreakdash-isomorphisme $\tors{2}E_\eta \isoto (\Z/2)^2$ en choisissant
d'envoyer le point de coordonnées $(X,Y)=(e_1,0)$ sur $(0,1)$ et le
point de coordonnées $(X,Y)=(e_2,0)$ sur $(1,0)$.
Des isomorphismes $H^1(K, \tors{2}E_\eta)=\Hgoth(K)$
et $H^1(U, \tors{2}\Erond)=\Hgoth(U)=\Hgoth(k[t][1/r])$
s'en déduisent
par additivité de la cohomologie.
La définition suivante permet de préserver la symétrie entre les~$e_i$ dans les énoncés à venir:
pour $\mgoth \in \Hgoth(K)$,
nous dirons que \emph{le triplet $(m_1,m_2,m_3)$ représente~$\mgoth$} si c'est
un triplet de polynômes séparables de~$k[t]$ dont le produit est un carré (non nul)
et si la classe dans $\Hgoth(K)$
du couple $(m_1,m_2)$ est égale à~$\mgoth$.

Les trois propositions suivantes explicitent respectivement
les groupes\index{groupe de Selmer!géométrique|(} de $2$\nobreakdash-Selmer géométriques, les extensions $L_M/\kappa(M)$
et les applications $\delta_M$; autrement dit, tous les objets qui apparaissent
dans la définition de la \condD{}.

\bigskip
\begin{proposition}
\label{ch1explsg}
Le sous-groupe $\SG_2(\A^1_k,\Erond) \subset H^1(U,\tors{2}\Erond)=\Hgoth(k[t][1/r])$
(cf.~lemme~\ref{ch1lemmepeutevaluer}) est égal à l'ensemble
des classes $\mgoth \in \Hgoth(k[t][1/r])$ représentées par
des triplets $(m_1,m_2,m_3)$ tels que pour tout $i \in \{1,2,3\}$,
les polynômes~$m_i$ et~$p_i$ soient premiers entre eux.
Le sous-groupe $\SG_2(\P^1_k,\Erond) \subset \SG_2(\A^1_k,\Erond)$
est égal à l'ensemble obtenu en ajoutant la condition que les polynômes~$m_i$ soient
tous de degré pair.
\end{proposition}

\bigskip
Pour $M \in \Mrond$, notons $v_M \colon K^\star \rightarrow \Z$ la valuation normalisée en~$M$.

\bigskip
\begin{demo}
Si~$M$ est un point fermé de~$\P^1_k$,
l'image de~$\mgoth$ dans $H^1(K_M^\sh,\tors{2}E_\eta)=\Z/2 \times \Z/2$
est égale à la classe du couple $(v_M(m_1), v_M(m_2))$.
Notamment, pour $M=\infty$, il en résulte que
les~$m_i$ sont tous de degré pair si et seulement si l'image de~$\mgoth$
dans $H^1(K_\infty^\sh,\tors{2}E_\eta)$ est nulle.  Comme la courbe elliptique~$E_\eta$
a bonne réduction à l'infini, on a $E_\eta(K_\infty^\sh)/2=0$
(cf.~proposition~\ref{annpropekshsurdeux});
par conséquent,
l'image de~$\mgoth$ dans $H^1(K_\infty^\sh,\tors{2}E_\eta)$ est nulle si et seulement
si l'image de~$\mgoth$ dans $H^1(K_\infty^\sh,E_\eta)$ est nulle.  La seconde assertion
de la proposition découle donc de la première.

Pour prouver celle-ci, il suffit de vérifier que pour tout $i \in \{1,2,3\}$
et tout facteur irréductible~$p$ de~$p_i$, notant $M \in \Mrond$ le point où~$p$
s'annule,
l'image de~$\mgoth$ dans $H^1(K_M^\sh,\tors{2}E_\eta)$ appartient au sous-groupe
$E_\eta(K_M^\sh)/2$ si et seulement si les polynômes~$m_i$ et~$p$ sont premiers entre eux.

\bigskip
\begin{lemme}
\label{ch1expldlemme}
Soit $M \in \Mrond$. L'unique point d'ordre~$2$ de $E_\eta$ qui se spécialise
sur $\Erond^0_M$
est celui de coordonnées $(X,Y)=(e_i,0)$, où $i \in \{1,2,3\}$
est tel que $p_i(M)=0$ (cf.~lemme~\ref{annlemmeuniquespec}).
\end{lemme}

\bigskip
\begin{demo}
Comme l'équation de Weierstrass~(\ref{ch1expleqw})
est minimale,
l'ouvert de lissité sur~$\A^1_k$ du sous-schéma fermé de $\P^2_k \times_k \A^1_k$
défini par~(\ref{ch1expleqw})
s'identifie à $\Erond^0_{\A^1_k}$
(cf.~\cite[Ch.~IV, §9, Cor.~9.1]{silvaec2}).
Il suffit donc de lire l'équation~(\ref{ch1expleqw})
modulo~$M$ pour déterminer quel point d'ordre~$2$
se spécialise sur $\Erond^0_M$: c'est celui qui ne se spécialise pas sur
le point singulier de la cubique réduite.
\end{demo}

\bigskip
Soient $i \in \{1,2,3\}$ et~$p$ un facteur irréductible de~$p_i$, s'annulant
en $M \in \Mrond$.
D'après le lemme~\ref{ch1expldlemme} et la proposition~\ref{annpropcaracteksh},
le sous-groupe $E_\eta(K_M^\sh)/2 \subset H^1(K_M^\sh,\tors{2}E_\eta)=\tors{2}E_\eta(K)$
est égal au sous-groupe de $\tors{2}E_\eta(K)$ engendré par le
point de coordonnées $(X,Y)=(e_i,0)$.
L'image de~$\mgoth$ dans $H^1(K_M^\sh,\tors{2}E_\eta)$ appartient donc
à $E_\eta(K_M^\sh)/2$ si et seulement si $v_M(m_i)=0$, d'où la\index{groupe de Selmer!géométrique|)} proposition.
\end{demo}

\bigskip
\begin{proposition}
\label{ch1expldlm}
Soit $M \in \Mrond$.  Si le polynôme~$r$ s'annule en~$M$ avec multiplicité~$1$,
alors $L_M=\kappa(M)$.  Sinon, la classe dans $\kappa(M)^\star/\kappa(M)^{\star 2}$
de l'extension quadratique ou triviale $L_M/\kappa(M)$ est égale
à la classe de $p_j(M)$, où $(i,j,k)$ est l'unique permutation cyclique
de $(1,2,3)$ telle que $p_i(M)=0$.
\end{proposition}

\bigskip
\begin{demo}
Notons~$n$ l'ordre du groupe $\ff{F}M(L_M)$.  Comme $n=2v_M(r)$,
soit~$r$ s'annule en~$M$ avec multiplicité~$1$, auquel cas $n=2$
et donc $L_M=\kappa(M)$,
soit $n>2$ et~$L_M$ est alors la plus petite extension de~$\kappa(M)$
sur laquelle sont définies les pentes des tangentes au point singulier de la cubique
obtenue en réduisant l'équation~(\ref{ch1expleqw}) modulo~$M$
(cf.~proposition~\ref{annpropcorfauxrat} et lemme~\ref{annlemmeextreme}).
La classe de cette dernière extension est bien celle de $p_j(M)$.
\end{demo}

\bigskip
Pour $M \in \Mrond$, notons $\gamma_M \in \kappa(M)^\star/\kappa(M)^{\star 2}$\glossary{$\gamma_M$}
la classe de $L_M/\kappa(M)$.
Comme le $L_M$\nobreakdash-groupe $\ff{F}M \otimes_{\kappa(M)} L_M$ est constant cyclique d'ordre
pair (cf.~proposition~\ref{annpropcyclpair}), les suites exactes
\begin{equation}
\label{ch1explse1}
\xymatrix{
0 \ar[r] & \ff{\tors{2}F}M \ar[r] & \ff{F}M \ar[r] & 2\ff{F}M \ar[r] & 0
}
\end{equation}
et
\begin{equation}
\label{ch1explse2}
\xymatrix{
0 \ar[r] & 2\ff{F}M \ar[r] & \ff{F}M \ar[r] & \ff{F}M/2 \ar[r] & 0
}
\end{equation}
permettent de voir que
$\tors{2}H^1(L_M,\ff{F}M)=H^1(L_M,\ff{\tors{2}F}M)=H^1(L_M,\Z/2)=L_M^\star/L_M^{\star 2}$,
d'où une injection canonique
$\kappa(M)^\star/\langle \kappa(M)^{\star 2}, \gamma_M\rangle
\hookrightarrow H^1(L_M,\ff{F}M)$.

\bigskip
\begin{proposition}
\label{ch1expldelta}
Soit $M \in \Mrond$.
L'image de $\delta_M \colon \SG_2(\A^1_k,\Erond)
\rightarrow H^1(L_M,\ff{F}M)$ est incluse dans le sous-groupe
$\kappa(M)^\star/\langle \kappa(M)^{\star 2}, \gamma_M\rangle \subset H^1(L_M,\ff{F}M)$.
Pour tout $\mgoth \in \SG_2(\A^1_k,\Erond)$, on a
$$
\delta_M(\mgoth)=m_i(M) \left(p_j(M)\right)^{v_M(m_j)}
\in \kappa(M)^\star/\langle \kappa(M)^{\star 2}, \gamma_M\rangle
\rlap{\text{,}}
$$
où le triplet $(m_1,m_2,m_3)$ représente~$\mgoth$ et $(i,j,k)$ est l'unique
permutation cyclique de $(1,2,3)$ telle que $p_i(M)=0$.
\end{proposition}

\bigskip
(Par un léger abus d'écriture, nous notons encore $\delta_M$ la composée
de la flèche naturelle $\SG_2(\A^1_k,\Erond) \rightarrow \tors{2}H^1(\A^1_k,\Erond)$
(cf.~suite exacte~(\ref{ch1sesgtsg}))
et de l'application $\delta_M$ du paragraphe~\ref{ch1hypnot}.)

\bigskip
\begin{demo}
Posons $G=\Gal(L_M/\kappa(M))$.
Les suites spectrales de Hochschild-Serre
$H^p(G,H^q(L_M,A))\Longrightarrow H^{p+q}(\kappa(M),A)$ pour $A=\ff{F}M$ et $A=\ff{\tors{2}F}M$
fournissent le diagramme commutatif
$$
\xymatrix{
H^1(\kappa(M), \ff{\tors{2}F}M) \ar[d] \ar[r] & H^1(L_M,\ff{\tors{2}F}M)^G \ar[d] \ar[r] &
H^2(G,\ff{\tors{2}F}M(L_M)) \ar[d] \\
\tors{2}H^1(\kappa(M),\ff{F}M) \ar[r] & \tors{2}H^1(L_M,\ff{F}M)^G \ar[r] & \tors{2}H^2(G,\ff{F}M(L_M))
\rlap{\text{,}}
}
$$
dont la ligne supérieure est exacte et la ligne inférieure est un complexe.
On~a déjà remarqué que la flèche verticale du milieu
est un isomorphisme.  La~flèche verticale de droite est par ailleurs
injective; en effet, si le groupe~$G$ n'est pas trivial, il est isomorphe à~$\Z/2$
et agit par multiplication par~$-1$ sur $\ff{F}M(L_M)$ (cf.~lemme~\ref{ch1lemmefm}),
auquel cas la périodicité de la cohomologie modifiée des groupes cycliques permet
d'identifier la flèche qui nous intéresse à l'inclusion
naturelle $H^0(G,\ff{\tors{2}F}M(L_M)) \rightarrow \tors{2}H^0(G,\ff{F}M(L_M))$.
Une chasse au diagramme montre maintenant que l'image de la flèche
de restriction $\tors{2}H^1(\kappa(M),\ff{F}M) \rightarrow \tors{2}H^1(L_M,\ff{F}M)$
est incluse dans l'image de la composée
$H^1(\kappa(M),\ff{\tors{2}F}M) \rightarrow H^1(L_M,\ff{\tors{2}F}M)=
\tors{2}H^1(L_M,\ff{F}M)$, autrement dit dans le
sous-groupe $\kappa(M)^\star/\langle \kappa(M)^{\star 2},\gamma_M \rangle$.
La première assertion de la proposition est donc établie.

Pour la seconde, quitte à étendre les scalaires de~$k$ à~$\kappa(M)$ puis
à~$L_M$, on peut supposer que le point~$M$ est $k$\nobreakdash-rationnel
et que $L_M=\kappa(M)$, de sorte que $\gamma_M=1$.
Notons alors $\Orond_M^\h$ le hensélisé de l'anneau local de~$\P^1_k$ en~$M$,
$K_M^\h$ son corps des fractions et~$v$ la valuation normalisée associée.
Posons $C=\Spec(\Orond_M^h)$.
Si $\mgoth \in \SG_2(C,\Erond)$, nous dirons que \emph{le triplet
$(m_1,m_2,m_3)$ représente~$\mgoth$}
si c'est un triplet d'éléments de~$K_M^{\h}$ dont le produit est un carré,
si $v(m_1),v(m_2),v(m_3)\in\{0,1\}$
et si la classe dans $\Hgoth(K_M^\h)$ du couple $(m_1,m_2)$ est égale à~$\mgoth$.
Soit $(i,j,k)$ l'unique permutation cyclique de $(1,2,3)$
telle que $p_i(M)=0$.
Comme dans la proposition~\ref{ch1explsg}, le sous-groupe $\SG_2(C,\Erond) \subset \Hgoth(K_M^\h)$
est égal à l'ensemble des classes représentées par un triplet $(m_1,m_2,m_3)$
tel que $v(m_i)=0$.

Soit $\delta_C$ la composée
$\SG_2(C,\Erond) \rightarrow \tors{2}H^1(C,\Erond) \rightarrow \tors{2}H^1(\kappa(M),\ff{F}M)$,
où la première flèche est issue de la suite exacte~(\ref{ch1sesgtsg}) et
la seconde est induite par le morphisme $\Erond \rightarrow i_{M\star}\ff{F}M$.
Vu le triangle commutatif
$$
\xymatrix@R=3ex{
\SG_2(\A^1_k,\Erond) \ar[dd] \ar[dr]^{\delta_M} \\
& \tors{2}H^1(\kappa(M),\ff{F}M)\rlap{\text{,}} \\
\SG_2(C,\Erond) \ar[ur]_{\delta_C}
}
$$
il suffit, pour conclure, d'établir que
l'on a $\delta_C(\mgoth)=m_i(M)(p_j(M))^{v(m_j)}$
dans $\kappa(M)^\star/\kappa(M)^{\star 2}$
pour tout $\mgoth \in \SG_2(C,\Erond)$, où le triplet $(m_1,m_2,m_3)$ représente~$\mgoth$.

Vérifions d'abord cette égalité lorsque $\mgoth \in E_\eta(K_M^\h)/2 \subset
\SG_2(C,\Erond)$
et que $v(m_j)=1$.  Comme $\delta_C(\mgoth)=1$ (cf.~suite exacte~(\ref{ch1sesgtsg})),
il s'agit de montrer que $m_i(M)=p_j(M)$ à un carré de~$\kappa(M)$ près.
Choisissons un relèvement $P \in E_\eta(K_M^\h)$ de~$\mgoth$.
Soit $\phi_\eta \colon E_\eta \rightarrow E''_\eta$ le quotient de~$E_\eta$
par l'unique point d'ordre~$2$ de~$E_\eta(K)$ qui se spécialise sur~$\Erond^0_M$.
Notons $\phi''_\eta$ l'isogénie duale de~$\phi_\eta$ et considérons le diagramme
commutatif
\begin{equation}
\label{ch1expldiagcom}
\begin{aligned}
\myxyhook\xymatrix{
0 \ar[r] & \tors{2}E_\eta \ar[d] \ar[r] & E_\eta \ar[r]^2 \ar[d]^{\phi_\eta}
& E_\eta \ar@{=}[d] \ar[r] & 0 \\
0 \ar[r] & \Z/2 \ar[r] & E''_\eta \ar[r]^{\phi''_\eta} & E_\eta \ar[r] & 0 \rlap{\text{,}}
}
\end{aligned}
\end{equation}
dont les lignes sont exactes.
Le lemme~\ref{ch1expldlemme} montre que la classe de~$m_i$
dans $K_M^{\h\star}/K_M^{\h\star 2} = H^1(K_M^\h,\Z/2)$ est égale à l'image de~$P$ par le
bord de la suite exacte inférieure de ce diagramme, puisque celui-ci commute.
Soient $\Erond''$ le modèle de Néron de~$E''_\eta$ sur~$C$, $\Erond''^0 \subset \Erond''$
sa composante neutre, $\ff{F''}M$ la fibre spéciale de $\Erond''/\Erond''^0$
et $i_M \colon \Spec(\kappa(M)) \rightarrow C$ l'immersion fermée canonique.
L'isogénie~$\phi''_\eta$ induit un morphisme de suites exactes de faisceaux
étales sur~$C$
$$
\xymatrix{
0 \ar[r] & \Erond''^0 \ar[d] \ar[r] & \Erond'' \ar[d] \ar[r] &
i_{M\star} \ff{F''}M \ar[d] \ar[r] & 0 \\
0 \ar[r] & \Erond^0 \ar[r] & \Erond \ar[r] & i_{M\star} \ff{F}M \ar[r] & 0 \rlap{\text{,}}
}
$$
où l'on convient de noter encore~$\Erond$ la restriction de~$\Erond$ à~$C$.
La flèche verticale de gauche est surjective d'après le lemme~\ref{annlemmesurje0},
la flèche verticale de droite l'est d'après la proposition~\ref{annpropppp}; la flèche verticale
du milieu est donc elle aussi surjective.
Compte tenu que le noyau de la flèche verticale
de droite est isomorphe à~$i_{M\star} \Z/2$ (cf.~proposition~\ref{annpropppp}), il en résulte
que les lignes du diagramme commutatif
\begin{equation}
\label{ch1expldiagcom2}
\begin{aligned}
\xymatrix{
0 \ar[r] & \Z/2 \ar[r] \ar[d] & \Erond'' \ar[d] \ar[r] & \Erond \ar[d] \ar[r] & 0 \\
0 \ar[r] & i_{M\star}\Z/2 \ar[r] & i_{M\star}\ff{F''}M \ar[r] & i_{M\star}\ff{F}M \ar[r] & 0
}
\end{aligned}
\end{equation}
sont exactes.  Notons $\alpha \in \ff{F}M(\kappa(M))$ l'image de~$P$ par la flèche
verticale de droite.  Comme~(\ref{ch1expldiagcom2}) commute, l'image~$\beta$
de~$\alpha$
dans $H^1(\kappa(M),\Z/2)=\kappa(M)^\star/\kappa(M)^{\star 2}$ par le bord
de la suite exacte inférieure de ce diagramme est égale à la classe de $m_i(M)$;
il reste donc seulement à vérifier que~$\beta$ coïncide avec la classe de~$p_j(M)$.
Soit $L''_M/\kappa(M)$ l'extension quadratique ou triviale minimale telle
que le $L''_M$\nobreakdash-groupe $\ff{F''}M \otimes_{\kappa(M)} L''_M$ soit
constant (cf.~proposition~\ref{annpropcyclpair}).
L'application $\ff{F''}M(L''_M) \rightarrow \ff{F}M(L''_M)$ induite par~$\phi''_\eta$
est surjective
puisque le $L''_M$\nobreakdash-groupe $\ff{F''}M \otimes_{\kappa(M)} L''_M$ est constant;
la classe~$\beta$ appartient donc au noyau de la flèche
de restriction $H^1(\kappa(M),\Z/2)\rightarrow H^1(L''_M,\Z/2)$.
D'autre part, comme $\ff{F}M(L_M)$ est d'ordre
pair, la suite exacte inférieure du diagramme~(\ref{ch1expldiagcom2}) montre
que $\ff{F''}M(L''_M)$ est d'ordre~$>2$.  Compte tenu de la proposition~\ref{annpropcorfauxrat}
et du lemme~\ref{annlemmeextreme}, il en résulte que la classe de l'extension
$L''_M/\kappa(M)$
dans $\kappa(M)^\star/\kappa(M)^{\star 2}$ est égale à la classe de~$p_j(M)$.
Ainsi a-t-on nécessairement $\beta=1$ ou $\beta=p_j(M)$ dans $\kappa(M)^\star/\kappa(M)^{\star 2}$.
Supposons, par l'absurde, que $\beta \neq p_j(M)$, auquel cas $\beta=1$. Cela
signifie que~$\alpha$ appartient à l'image de la flèche
$\ff{F''}M(\kappa(M)) \rightarrow \ff{F}M(\kappa(M))$ induite par~$\phi''_\eta$.
Comme $\beta \neq p_j(M)$ et que $\beta=1$, on a $p_j(M) \neq 1$ (toujours
dans $\kappa(M)^\star /\kappa(M)^{\star 2}$), de sorte que $L''_M \neq \kappa(M)$.
Autrement dit, le $\kappa(M)$\nobreakdash-groupe $\ff{F''}M$ n'est pas constant.
Vu la proposition~\ref{annpropcyclpair} et la ligne inférieure du diagramme~(\ref{ch1expldiagcom2}),
il s'ensuit que la flèche $\ff{F''}M(\kappa(M)) \rightarrow \ff{F}M(\kappa(M))$
induite par~$\phi''_\eta$ est nulle, et donc que $\alpha=0$.
La proposition~\ref{annpropekshsurdeux} permet d'en déduire que l'image
de~$\mgoth$ dans $E_\eta(K_M^\sh)/2$ est nulle, ce qui contredit l'hypothèse
selon laquelle $v(m_j)=1$.

Passons maintenant au cas général. Fixons $\mgoth \in \SG_2(C,\Erond)$ et
considérons le diagramme commutatif
$$
\xymatrix{
E_\eta(K_M^\h)/2 \ar[r] \ar[d] & \ff{F}M(\kappa(M))/2 \ar[d] \\
E_\eta(K_M^\sh)/2 \ar[r] & \ff{F}M\!\left(\overline{\kappa(M)}\right)\!/2 \rlap{\text{,}}
}
$$
où $\overline{\kappa(M)}$ est le corps résiduel de $\Orond_M^\sh$.
La proposition~\ref{annpropekshsurdeux} montre que la flèche horizontale supérieure
est surjective et que la flèche horizontale inférieure est bijective.  D'autre part,
la flèche verticale de droite est bijective puisque le $\kappa(M)$\nobreakdash-groupe $\ff{F}M$ est constant.
Il en résulte que la flèche verticale de gauche est surjective.

L'image de~$\mgoth$ dans $H^1(K_M^\sh,\tors{2}E_\eta)=\Z/2 \times \Z/2$ est
égale à $(v(m_1),v(m_2))$ et appartient par ailleurs au sous-groupe
$E_\eta(K_M^\sh)/2 \subset H^1(K_M^\sh,\tors{2}E_\eta)$, par définition du groupe
de Selmer géométrique.
D'après la propriété de surjectivité que l'on vient d'établir,
il existe un élément $\mgoth' \in E_\eta(K_M^\h)/2$ ayant même image que~$\mgoth$
dans $H^1(K_M^\sh,\tors{2}E_\eta)$.
Comme l'égalité que l'on veut démontrer est multiplicative et qu'elle a déjà été
établie pour~$\mgoth'$ si l'image de~$\mgoth$ dans $H^1(K_M^\sh,\tors{2}E_\eta)$
n'est pas nulle, on peut supposer,
quitte à remplacer~$\mgoth$ par $\mgoth+\mgoth'$,
que $v(m_1)=v(m_2)=v(m_3)=0$, auquel
cas $\mgoth \in H^1(C,\tors{2}\Erond) \subset \SG_2(C,\Erond)$.
Le carré commutatif
$$
\xymatrix{
H^1(C,\tors{2}\Erond) \ar[r] \ar[d] & H^1(C,\Erond) \ar[d] \\
H^1(\kappa(M),\ff{\tors{2}F}M) \ar[r] & H^1(\kappa(M), \ff{F}M)
}
$$
montre maintenant
que $\delta_C(\mgoth) \in \kappa(M)^\star/\kappa(M)^{\star 2}=H^1(\kappa(M),\ff{\tors{2}F}M)$
est égal à l'image de~$\mgoth$ par la flèche $H^1(C,\tors{2}\Erond)
\rightarrow H^1(\kappa(M),\ff{\tors{2}F}M)$ induite par le
morphisme de faisceaux $\tors{2}\Erond \rightarrow i_{M\star}
(\ff{\tors{2}F}M)$. Celui-ci s'identifie à la composée
$\tors{2}\Erond \rightarrow i_{M\star}(\tors{2}\Erond_M)
\rightarrow i_{M\star}(\tors{2}\Erond_M/\tors{2}\Erond_M^0)$,
compte tenu que $\tors{2}\Erond_M=\Z/2 \times \Z/2$, $\tors{2}\Erond_M^0=\Z/2$
et $\ff{\tors{2}F}M=\Z/2$.
Comme $H^1(C,\tors{2}\Erond)=\Hgoth(\Orond_M^\h)=\Hgoth(\kappa(M))$
et que ces isomorphismes canoniques appliquent $\mgoth \in H^1(C,\tors{2}\Erond)$
sur $(m_1(M),m_2(M)) \in \Hgoth(\kappa(M))$, on en conclut
que $\delta_C(\mgoth)$ est égal à l'image
de $(m_1(M),m_2(M))$ par l'application
$
\Hgoth(\kappa(M)) \longrightarrow \kappa(M)^\star/\kappa(M)^{\star 2}
$
obtenue en tensorisant par $\kappa(M)^\star/\kappa(M)^{\star 2}$
la flèche $\tors{2}E_\eta(K)=\Z/2\times\Z/2 \longrightarrow \Z/2$ de quotient par
l'unique point d'ordre~$2$ de $E_\eta(K)$ qui se spécialise sur~$\Erond^0_M$.
D'où le résultat, grâce au lemme~\ref{ch1expldlemme}.
\end{demo}

\bigskip
L'énoncé de la proposition~\ref{ch1expldelta} se simplifie dans le cas de mauvaise réduction en~$M$ de type~$I_n$
avec $n>2$.

\bigskip
\begin{corollaire}
\label{ch1expldeltai4}
Soit $M \in \Mrond$ tel que le polynôme~$r$ s'annule avec multiplicité~$>1$ en~$M$.
Pour tout $\mgoth \in \SG_2(\A^1_k,\Erond)$, on a
$\delta_M(\mgoth)=m_i(M) \in \kappa(M)^\star/\langle \kappa(M)^{\star 2}, \gamma_M\rangle$,
où le triplet $(m_1,m_2,m_3)$ représente~$\mgoth$ et $i \in \{1,2,3\}$
est tel que $p_i(M)=0$.
\end{corollaire}

\bigskip
\begin{demo}
Vu la proposition~\ref{ch1expldelta}, il suffit de vérifier
que la classe de $p_j(M)$ dans $\kappa(M)^\star/\kappa(M)^{\star 2}$
est égale à~$\gamma_M$; ceci résulte de la proposition~\ref{annpropcorfauxrat}
et du lemme~\ref{annlemmeextreme}.
\end{demo}

\bigskip
L'explicitation de la \condD{} permet de donner une preuve courte du résultat
de finitude suivant, qui sera utile à la fois pour la preuve du théorème~\ref{ch1thprin}
et pour son application aux secondes descentes (théorème~\ref{ch1thsecdesc}).

\bigskip
\begin{corollaire}
\label{ch1finitude}
Si $\Mrond \neq \emptyset$, les groupes $\SD$ et $\SDz$ sont finis.
\end{corollaire}

\bigskip
\begin{demo}
Le sous-groupe $\Hgoth(k) \subset \SG_2(\A^1_k,\Erond)$ étant
d'indice fini (cf.~proposition~\ref{ch1explsg}), il suffit de prouver
que le noyau de la flèche
$$
\Hgoth(k) \xrightarrow{\prod\delta_M} \prod_{M \in \Mrond} \kappa(M)^\star/
\langle \kappa(M)^{\star 2}, \gamma_M \rangle
$$
est fini.
L'hypothèse $\Mrond \neq \emptyset$ signifie que le polynôme~$r$
n'est pas constant.  Il~en résulte qu'aucun des~$p_i$ n'est constant,
puisqu'ils sont tous du même degré.  Par conséquent, il existe $M_1,M_2 \in \Mrond$
tels que $p_1(M_1)=0$ et $p_2(M_2)=0$.
D'après la proposition~\ref{ch1expldelta}, l'application
$$
\Hgoth(k) \xrightarrow{\delta_{M_1} \times \delta_{M_2}}
\kappa(M_1)^\star/\langle \kappa(M_1)^{\star 2}, \gamma_{M_1} \rangle
\times
\kappa(M_2)^\star/\langle \kappa(M_2)^{\star 2}, \gamma_{M_2} \rangle
$$
envoie $(m_1,m_2) \in \Hgoth(k)$ sur la classe de $(m_1,m_2)$.
On conclut alors en choisissant une extension finie galoisienne $\ell/k$
dans laquelle se plongent $L_{M_1}$ et $L_{M_2}$ et en remarquant
que le noyau de la flèche de restriction $\Hgoth(k) \rightarrow \Hgoth(\ell)$
est fini puisqu'il s'identifie au groupe de cohomologie\index{condition (D)@\condD{}|)} $H^1(\Gal(\ell/k),\Z/2 \times \Z/2)$.
\end{demo}

\section{Symétrisation du calcul de Selmer}
\label{ch1parpreuvedebut}
\label{ch1parsymetrisation}

Nous\index{groupe de Selmer|(} détaillons ici une construction due à Colliot-Thélène, Skorobogatov et
Swinnerton-Dyer (cf.~\cite[§1]{css}) permettant d'exprimer le groupe de $2$\nobreakdash-Selmer
d'une courbe elliptique sur un corps de nombres comme le noyau d'une certaine forme bilinéaire symétrique,
l'intérêt de cette opération étant que le comportement en famille de ladite forme bilinéaire
sera plus facile à étudier que celui du groupe de $2$\nobreakdash-Selmer.
Les résultats de ce paragraphe ne sont pas nouveaux, à l'exception de
la proposition~\ref{ch1wrvneron}; tout au plus leur présentation est-elle quelque
peu simplifiée (notamment, les espaces~$\Wrond_v$ n'interviennent que dans
la proposition~\ref{ch1wrvneron}).

La symétrisation du calcul de Selmer nécessite trois types d'ingrédients: les théorèmes
de dualité locale (pour les modules finis et pour les courbes elliptiques,
cf.~\cite[Ch.~I, §2 et §3]{milneadt}), un peu de théorie du corps de classes global,
et de l'algèbre linéaire sur~$\ff{\F}2$. Les preuves sont légèrement plus simples
lorsque les points d'ordre~$2$ de la courbe elliptique considérée sont tous rationnels
(et nous ferons donc cette hypothèse); le type de réduction ne joue en revanche
aucun rôle.

\bigskip
Les notations de ce paragraphe sont indépendantes de celles du reste du chapitre.
Soit~$k$ un corps de nombres.  Notons~$\Omega$ l'ensemble de ses places, et fixons
un sous-ensemble fini~$S_0 \subset \Omega$ contenant les places archimédiennes, les
places dyadiques et un système de générateurs du groupe de classes de~$k$.
Soit~$E$ une courbe elliptique sur~$k$ dont tous les points d'ordre~$2$ sont rationnels.
Pour $v \in \Omega$, on définit les $\ff{\F}2$\nobreakdash-espaces vectoriels suivants:
$$
\fv{V}v=H^1(k_v,\tors{2}E) \quad;\quad W_v(E)=E(k_v)/2 \quad;\quad
\ft{T}v=\Ker\left(\fv{V}v \rightarrow H^1(k_v^\nr,\tors{2}E)\right)\!\rlap{\text{,}}
$$
\glossary{$\fv{V}v$, $W_v(E)$, $\ft{T}v$, $\fv{V}S$, $I^S$}où~$k_v^\nr$ désigne une extension non ramifiée maximale de~$k_v$.
La suite exacte de Kummer permet de voir~$W_v(E)$ naturellement comme un sous-espace de~$\fv{V}v$.
Pour $S \subset \Omega$ fini contenant~$S_0$, on pose:
$$
\fv{V}S=\bigoplus_{v \in S} \fv{V}v \quad;\quad
I^S=\Ker\left(H^1(k,\tors{2}E) \rightarrow \prod_{v \in \Omega \setminus S} H^1(k_v^\nr,
\tors{2}E)\right)\!\rlap{\text{.}}
$$

\bigskip
\begin{lemme}
La flèche naturelle $I^S \rightarrow \fv{V}S$ est injective.
\end{lemme}

\bigskip
\begin{demo}
Le choix d'un isomorphisme de $k$\nobreakdash-groupes $\tors{2}E\isoto (\Z/2)^2$ nous ramène
à prouver que la flèche naturelle
$\Orond_S^\star/\Orond_S^{\star 2} \rightarrow \prod_{v \in S} k_v^\star/k_v^{\star 2}$
est injective, en notant~$\Orond_S$ l'anneau des $S$\nobreakdash-entiers de~$k$.
Il s'agit là d'une conséquence de la théorie du corps de classes (cf.~\cite[Chapter~VII, Lemma~9.2]{casselsfrohlich}; c'est ici que l'hypothèse que
$S_0$ contient un système de générateurs du groupe de classes de~$k$ sert).
\end{demo}

\bigskip
L'accouplement de Weil $\tors{2}E \times \tors{2}E \rightarrow \Z/2$ et le
cup-produit induisent une forme bilinéaire alternée non dégénérée $\fv{V}v \times
\fv{V}v \rightarrow \Z/2$ pour tout $v \in \Omega$ (cf.~\cite[Ch.~I, Cor.~2.3 et
Th.~2.13]{milneadt}), notée par la suite $\langle \cdot, \cdot \rangle_v$.  Le
sous-espace $W_v(E) \subset \fv{V}v$ est totalement isotrope maximal (cf.~\cite[Ch~I,
Cor.~3.4 et Rem.~3.7]{milneadt}).  De plus, il résulte du théorème de F.~K.~Schmidt
que $W_v(E)=\ft{T}v$ pour toute place
$v \in \Omega \setminus S_0$ de bonne réduction pour~$E$.

Munissons~$\fv{V}S$ de la somme orthogonale des accouplements $\langle \cdot, \cdot \rangle_v$.
La loi de réciprocité globale\index{loi de réciprocité globale} montre que le sous-espace $I^S \subset \fv{V}S$ est totalement
isotrope; un calcul de dimensions utilisant le théorème des unités de Dirichlet et la
formule du produit permet de voir qu'il est même totalement isotrope maximal
(cf.~\cite[Prop.~1.1.1]{css}).

\bigskip
\begin{lemme}
\label{ch1existekv}
Il existe des sous-espaces totalement isotropes maximaux
$K_v \subset \fv{V}v$\glossary{$K_v$}
pour $v \in S_0$ tels
que $\fv{V}{S_0} = I^{S_0} \oplus \bigoplus_{v \in S_0} K_v$.
\end{lemme}

\bigskip
\begin{demo}
Ce lemme est prouvé dans~\cite[Lemma~1.1.3]{css}. En voici une démonstration plus
courte, due à Swinnerton-Dyer.
Les espaces~$\fv{V}v$ s'écrivent comme des sommes orthogonales de plans hyperboliques.
Il suffit donc de voir qu'étant donnés une famille $(\fv{V}i)_{1 \leq i \leq n}$ de plans
hyperboliques sur~$\ff{\F}2$ et un sous-espace totalement isotrope maximal~$I$ de leur somme
orthogonale $V=\bigoplus \fv{V}i$, il existe des droites $K_i \subset \fv{V}i$
telles que $V=I\oplus \bigoplus K_i$. Prouvons-le par récurrence sur~$n$.
Si $n=0$, il n'y a rien à démontrer.
Supposons donc $n>0$ et prenons pour~$K_n$ une droite quelconque de~$\fv{V}n$ non incluse dans~$I$.
Posons $V^-=\bigoplus_{i=1}^{n-1}\fv{V}i$ et $I^-=V^- \cap (I+K_n)$. Il est facile de voir
que $I^- \subset V^-$ est totalement isotrope. De plus, $\dim(V^-)=2n-2$ et
$\dim(I^-) \geq \dim(V^-) + \dim(I+K_n)-\dim(V) = n-1$, ce qui
montre que~$I^-$ est totalement isotrope maximal. L'hypothèse de récurrence
fournit des droites~$K_i \subset \fv{V}i$ pour~$i<n$ telles que
$V^-=I^-\oplus \bigoplus_{i=1}^{n-1}K_i$. On a alors $(I+K_n)\cap \bigoplus_{i=1}^{n-1} K_i
\subset I^- \cap \bigoplus_{i=1}^n K_i=0$, donc $I \cap \bigoplus_{i=1}^n K_i \subset I
\cap K_n = 0$, d'où finalement $V=I\oplus \bigoplus_{i=1}^n K_i$.
\end{demo}

\bigskip
Supposons choisis des $K_v \subset \fv{V}v$ pour $v \in S_0$ comme dans le lemme ci-dessus
et notons $K_v=\ft{T}v$ pour $v \in \Omega\setminus S_0$. On dispose maintenant
d'un sous-espace totalement isotrope maximal $K_v \subset \fv{V}v$ pour chaque $v \in \Omega$.
Posons $K_S=\bigoplus_{v\in S}K_v$ pour $S \subset \Omega$.

\bigskip
\begin{lemme}
\label{ch1vsisks}
On a $\fv{V}S=I^S \oplus K_S$
pour tout $S \subset \Omega$ fini contenant~$S_0$.
\end{lemme}

\bigskip
\begin{demo}
En effet on a l'inclusion $I^S \cap K_S \subset I^{S_0} \cap K_{S_0}$
de sous-groupes de $H^1(k,\tors{2}E)$, or $I^{S_0} \cap K_{S_0}=0$ par définition des~$K_v$.
\end{demo}

\bigskip
Le groupe de $2$\nobreakdash-Selmer de~$E$ est par définition le sous-groupe $\Sel_2(k,E)$
de $H^1(k,\tors{2}E)$
constitué des classes dont l'image dans~$\fv{V}v$ appartient à~$W_v(E)$ pour tout $v \in \Omega$.
Il s'identifie à $I^S \cap W_S(E)$ si~$S$ contient l'ensemble des places de mauvaise
réduction pour~$E$ puisqu'alors $W_v(E)=\ft{T}v$ pour $v \not\in S$.
Nous allons maintenant nous intéresser à un sous-groupe intermédiaire entre
$\Sel_2(k,E)$ et~$I^S$, à savoir
$$\Irond^S(E) = I^S \cap (W_S(E) + K_S)\rlap{\text{.}}$$\glossary{$\Irond^S(E)$}%
Cette définition est à mettre en parallèle avec celle du groupe de $2$\nobreakdash-Selmer \emph{géométrique}
d'une courbe elliptique sur le corps des fonctions de la droite projective,
dans laquelle intervenaient
les hensélisés stricts des anneaux locaux de la base (par opposition aux hensélisés).

Soient $a,b\in \Irond^S(E)$.  Par hypothèse, il existe $\alpha_v$, $\alpha'_v$
pour $v \in S$ avec $\alpha_v \in W_v(E)$, $\alpha'_v \in K_v$ et
$a=\sum_{v\in S}(\alpha_v+\alpha'_v)$. Posons
\begin{equation}
\label{ch1defacc}
\langle a, b \rangle = \sum_{v \in S} \langle \alpha_v, b \rangle_v \rlap{\text{.}}
\end{equation}
Étant donné que $W_v(E)$ et~$K_v$ sont totalement isotropes,
cet élément de $\Z/2$ ne dépend que du couple $(a,b)$.
La formule ci-dessus définit donc une
forme bilinéaire $\Irond^S(E) \times \Irond^S(E) \rightarrow \Z/2$.

Soient $\beta_v \in W_v(E)$, $\beta'_v\in K_v$ pour $v \in S$ tels que
$b=\sum_{v \in S}(\beta_v+\beta'_v)$. La loi de réciprocité globale\index{loi de réciprocité globale}
entraîne que $\sum_{v\in\Omega}\langle a, b\rangle_v=0$. Par ailleurs,
$\langle a, b\rangle_v=0$ pour $v \not\in S$ car les images de~$a$ et~$b$
dans~$\fv{V}v$ pour $v \not\in S$ appartiennent au sous-espace totalement isotrope $W_v(E)=\ft{T}v$;
d'où $\sum_{v \in S} \langle a,b\rangle_v=0$. Utilisant à nouveau la propriété
qu'ont $W_v(E)$ et~$K_v$ d'être totalement isotropes, on en déduit immédiatement:

\bigskip
\begin{proposition}
\label{ch1estsymetrique}
La forme bilinéaire $\Irond^S(E) \times \Irond^S(E) \rightarrow \Z/2$ définie
par~(\ref{ch1defacc}) est symétrique.
\end{proposition}

\bigskip
Il est évident sur~(\ref{ch1defacc}) que $\langle a,b\rangle=0$ si $b \in \Sel_2(k,E)$;
le groupe de $2$\nobreakdash-Selmer est donc inclus dans le noyau de cette forme bilinéaire symétrique.
Montrons maintenant:

\bigskip
\begin{proposition}
\label{ch1noyauselmer}
Si $S$ contient l'ensemble des places de mauvaise réduction pour~$E$,
le noyau de la forme bilinéaire symétrique
définie par~(\ref{ch1defacc}) est égal
au groupe de $2$\nobreakdash-Selmer de~$E$.
\end{proposition}

\bigskip
\begin{demo}
Supposons que $b\in \Irond^S(E)$ appartienne au noyau. On veut montrer que $b\in I^S \cap W_S(E)$;
pour cela, il suffit que $b$ soit orthogonal à~$W_S(E)$ dans~$\fv{V}S$, puisque~$W_S(E)$ est
totalement isotrope maximal. Soit $w \in W_S(E)$. D'après le lemme~\ref{ch1vsisks},
il existe $a \in I^S$ tel que $a-w\in K_S$.  On a alors $a\in\Irond^S(E)$, et par
hypothèse $\langle a, b\rangle=0$, ce qui signifie précisément que~$w$ et~$b$ sont
orthogonaux\index{groupe de Selmer|)} dans~$\fv{V}S$.
\end{demo}

\bigskip
Notons $\Wrond_v(E)$\glossary{$\Wrond_v(E)$} l'image de
$W_v(E)$ par la flèche de restriction $H^1(k_v, \tors{2}E)
\rightarrow H^1(k_v^\nr, \tors{2}E)$,
pour $v \in \Omega \setminus S_0$.  Le sous-espace $W_v(E) +
K_v$ de~$\fv{V}v$ est exactement l'image réciproque de $\Wrond_v(E)$ par
cette flèche, puisque $K_v=\ft{T}v$ pour de tels~$v$.  Autrement dit, le
sous-espace $\Wrond_v(E)$ de $H^1(k_v^\nr, \tors{2}E)$
suffit à déterminer la condition en une place $v \in S \setminus S_0$ pour qu'un élément
de $I^S$ appartienne à $\Irond^S(E)$.
La proposition suivante nous sera utile
par la suite: elle permet d'exprimer cette condition en termes de modèles de Néron.

\bigskip
\begin{proposition}
\label{ch1wrvneron}
Soit\index{modèle de Néron} $v \in\Omega \setminus S_0$. Notons $\kappa$ le corps résiduel de $v$,
$\kappabarre$ une clôture algébrique de $\kappa$ et $F$ le $\kappa$\nobreakdash-schéma en
groupes fini étale des composantes connexes de la fibre spéciale du modèle de
Néron de $E \otimes_k k_v$ sur $\Orond_v$.  On a alors canoniquement
$\Wrond_v(E) = \Im\left(F(\kappa) \rightarrow F(\kappabarre)/2\right)$.
\end{proposition}

\bigskip
\begin{demo}
Vu le carré commutatif
$$
\myxyhook\xymatrix{
E(k_v)/2 \ar@{^{ (}->}[r] \ar[d] & H^1(k_v, \tors{2}E) \ar[d] \\
E(k_v^\nr)/2 \ar@{^{ (}->}[r] & H^1(k_v^\nr,\tors{2}E) \rlap{\text{,}}
}
$$
le groupe~$\Wrond_v(E)$ s'identifie à l'image de la flèche naturelle
$E(k_v)/2 \rightarrow E(k_v^\nr)/2$.
Considérons maintenant le carré commutatif
$$
\xymatrix{
E(k_v)/2 \ar[d] \ar[r]^r & F(\kappa)/2 \ar[d] \\
E(k_v^\nr)/2 \ar[r]^{r'} & F(\kappabarre)/2 \rlap{\text{,}}
}
$$
où~$r$ et~$r'$ sont induites par les flèches de spécialisation.
Il résulte de la proposition~\ref{annpropekshsurdeux} que~$r'$ est un isomorphisme
et que~$r$ est surjective,
d'où le résultat.
\end{demo}

\section{Preuve du théorème~\ref{ch1thprin}}
\label{ch1parpreuvefin}

Pour prouver le théorème~\ref{ch1thprin}, on peut évidemment supposer
la fibration $\pi \colon X \rightarrow \P^1_k$ relativement minimale.
L'intérêt de cette hypothèse est essentiellement qu'elle permet d'écrire
que $U=\A^1_k \setminus \Mrond$.
On désigne par~$\Orond$ (resp.~$\Orond_S$, pour $S \subset \Omega$) l'anneau des entiers
(resp.~des $S$\nobreakdash-entiers) de~$k$.  D'autre part, on
conserve la notation $\Hgoth(X)=\Gm(X)/2\times\Gm(X)/2$ du paragraphe~\ref{ch1parexplicited}.

Les paragraphes~\ref{ch1premierparpreuve} à~\ref{ch1dernierparavantpreuve} contiennent
des définitions et résultats préliminaires à la preuve proprement dite du théorème~\ref{ch1thprin},
qui occupe le paragraphe~\ref{ch1parpreuve}. 

\subsection{Triplets préadmissibles, admissibles}
\label{ch1premierparpreuve}

Si $x \in \P^1_k$ est un point fermé, on note~$\xtilde$ son adhérence
schématique dans~$\P^1_{\Orond}$.  C'est un $\Orond$\nobreakdash-schéma fini
génériquement étale.  On définit de même l'adhérence~$\xtilde$ de~$x$
dans~$\P^1_{\Orond_v}$ lorsque~$x$ est un point fermé de~$\P^1_{k_v}$
et que $v \in \Omega_f$.  Si $S \subset \Omega$ est assez grand pour
que $\xtilde \cap \P^1_{\Orond_S}$ soit étale sur $\Spec(\Orond_S)$,
l'ensemble des points fermés de $\xtilde \cap \P^1_{\Orond_S}$
s'identifie à l'ensemble des places de~$\kappa(x)$ dont la trace
sur~$k$ n'appartient pas à~$S$.  On utilisera librement cette
identification par la suite.

Soit~$S$ un ensemble fini de places finies de $k$ contenant les places dyadiques.
Soit $\Trond$ une famille $(\ft{T'}M)_{M \in \Mrond}$, où $\ft{T'}M$ est un ensemble fini de
places finies de~$\kappa(M)$.
Notons $\ft{T}M \subset \Omega_f$ l'ensemble des traces sur $k$ des
places de $\ft{T'}M$.
On dit que le couple $(S, \Trond)$
est \emph{préadmissible}\index{préadmissible!couple} si la condition suivante est satisfaite:
\begin{itemize}
\medskip
\item[(\refstepcounter{equation}\theequation{})]
les sous-ensembles $\ft{T}M \subset \Omega$ pour $M \in \Mrond$
sont deux à deux disjoints et disjoints de~$S$; le $\Orond_S$\nobreakdash-schéma
$\P^1_{\Orond_S} \cap \bigcup_{M \in \Mrond \cup \{\infty\}} \Mtilde$
est étale;
pour tout $M \in \Mrond$, l'application trace $\ft{T'}M \rightarrow \ft{T}M$ est bijective.
\end{itemize}
\medskip

Soit $x \in U(k)$. On dit que le triplet $(S, \Trond, x)$ est
\emph{préadmissible}\index{préadmissible!triplet} si le couple $(S,\Trond)$ l'est et si de plus la
condition suivante est satisfaite:
\begin{itemize}
\medskip
\item[(\refstepcounter{equation}\theequation{})] $x$ est entier hors de~$S$ (\emph{i.e.}
$\xtilde \cap \widetilde{\infty} \cap \P^1_{\Orond_S}=\emptyset$);
pour tout $M \in \Mrond$, le schéma $\xtilde \cap \Mtilde \cap \P^1_{\Orond_S}$ est réduit, et son
ensemble sous-jacent est la réunion de $\ft{T'}M$ et
d'une place de $\kappa(M)$ hors de $\ft{T'}M$, que l'on note $w_M$.
\end{itemize}

\bigskip
Étant donné un tel triplet $(S,\Trond, x)$, on notera $v_M$ la trace de $w_M$
sur $k$, et l'on posera\glossary{$v_M$, $w_M$, $\ft{T}M$, $T(x)$}
$$
T = S \cup \bigcup_{M \in \Mrond} \ft{T}M
$$
et $$T(x) = T \cup \ensemble{v_M}{M \in \Mrond}\rlap{\text{.}}$$

\bigskip
Les conditions de préadmissibilité sont de nature géométrique.  Nous voudrons aussi imposer
une condition arithmétique sur les triplets considérés afin d'assurer l'existence de
points locaux sur~$X_x$ (cf.~condition~(\ref{ch1conditionadmiss}) ci-dessous, et proposition~\ref{ch1propexistadel}).
Le lemme suivant, valable sur tout corps de caractéristique~$0$, nous sera utile pour définir cette condition.

\bigskip
\begin{lemme}
\label{ch1existekmd}
Soit $M \in \Mrond$.
Pour tout $d \in \tors{2}H^1(\kappa(M), \ff{F}M)$, il existe une extension quadratique ou triviale
minimale $K_{M,d}$ de~$\kappa(M)$ telle que~$d$ appartienne au noyau de
la restriction $H^1(\kappa(M), \ff{F}M) \rightarrow H^1(L_M K_{M,d}, \ff{F}M)$
et que $L_M K_{M,d}$ se plonge $L_M$\nobreakdash-linéairement dans toute extension
$\ell/L_M$ telle que l'image de~$d$ dans~$H^1(\ell,\ff{F}M)$
soit nulle.
\end{lemme}

\bigskip
\begin{demo}
Choisissons un isomorphisme $\ff{F}M(L_M) \isoto \Z/2n$ pour un $n \in \N^\star$.
Cet isomorphisme induit une action de $G=\Gal(L_M/\kappa(M))$ sur $\Z/2n$ par transport de structure.
Munissons $\Z/n$ et~$\Z/4n$ de l'action triviale de~$G$ si $G=1$,
de l'action de~$G$ par multiplication par~$-1$ si $G=\Z/2$,
et considérons le diagramme commutatif de groupes abéliens suivant, dont les
lignes et les colonnes sont exactes:
\begin{equation}
\label{ch1existekmddiag}
\begin{aligned}
\xymatrix{
&& 0 \ar[d] & 0  \ar[d] \\
0 \ar[r] & \Z/2 \ar@{=}[d] \ar[r] & \Z/2n \ar[d] \ar[rd] \ar[r] & \Z/n \ar[d] \ar[r] & 0 \\
0 \ar[r] & \Z/2 \ar[r] & \Z/4n \ar[r] \ar[d] & \Z/2n \ar[d] \ar[r] & 0 \rlap{\text{.}} \\
&& \Z/2 \ar@{=}[r] \ar[d] & \Z/2 \ar[d] \\
&& 0 & 0
}
\end{aligned}
\end{equation}
Compte tenu du lemme~\ref{ch1lemmefm}, toutes les flèches de ce diagramme
sont $G$\nobreakdash-équi\-variantes, de sorte que des suites exactes de cohomologie s'en déduisent par
passage aux invariants sous~$G$.
Notons $u \colon H^1(\kappa(M),\Z/2) \rightarrow H^1(\kappa(M),\ff{F}M)$ l'application induite
par l'inclusion $\Z/2 \hookrightarrow \ff{F}M$
et $v \colon \Z/2 \rightarrow H^1(\kappa(M),\ff{F}M)$ le bord
de la suite exacte verticale de gauche du diagramme ci-dessus.
Remarquant que la flèche oblique est l'endomorphisme de multiplication par~$2$ dans $\Z/2n$,
une chasse au diagramme permet de voir qu'il existe
$\alpha \in H^1(\kappa(M),\Z/2)$ et $\beta \in \Z/2$
tels que $d=u(\alpha)+v(\beta)$, grâce à l'hypothèse selon laquelle $2d=0$.
Soit $K_{M,d}$ l'extension quadratique ou triviale de~$\kappa(M)$ définie
par~$\alpha$.
L'image de~$v(\beta)$ dans $H^1(L_M,\ff{F}M)$ est nulle
puisque l'application $H^0(L_M,\Z/4n) \rightarrow H^0(L_M,\Z/2)$ issue de
la suite exacte verticale du diagramme~(\ref{ch1existekmddiag})
est surjective.  Il en résulte que l'extension $K_{M,d}/\kappa(M)$
satisfait bien aux conditions voulues.
\end{demo}

\bigskip
On notera par la\glossary{$K_{M,d}$, $K_M$} suite~$K_M$ le corps $K_{M,d}$ donné par le lemme~\ref{ch1existekmd}
en prenant pour~$d$ l'image de~$[\Xrond]$ par la flèche naturelle $H^1(C,\Erond)\rightarrow H^1(\kappa(M),\ff{F}M)$.

\bigskip
Le couple $(S,\Trond)$ sera dit \emph{admissible}\index{admissible!couple} s'il est préadmissible
et que la condition suivante est vérifiée:
\begin{itemize}
\medskip
\item[(\refstepcounter{equation}\theequation{})\label{ch1conditionadmiss}]
pour tout $M \in \Mrond$, les places de~$\ft{T'}M$ sont totalement décomposées dans $L_M K_M$.
\end{itemize}
\medskip

Si $x \in U(k)$, on dit enfin que le triplet $(S,\Trond, x)$ est
\emph{admissible}\index{admissible!triplet}
s'il est préadmissible, si le couple $(S, \Trond)$ est admissible et si pour tout $M \in \Mrond$,
la place $w_M$ est totalement décomposée dans $L_M K_M$.

\subsection{Calcul symétrique des groupes de Selmer en famille}
\label{ch1parsymfam}

Pour chaque $x \in U(k)$, les résultats du paragraphe~\ref{ch1parsymetrisation}
fournissent un $\ff{\F}2$\nobreakdash-espace vectoriel de dimension finie muni d'une forme bilinéaire
symétrique dont le noyau est canoniquement isomorphe au groupe de $2$\nobreakdash-Selmer
de~$\Erond_x$.  Si l'on cherche à étudier la variation des groupes de $2$\nobreakdash-Selmer
des fibres de $\Erond_U \rightarrow U$ au-dessus de points rationnels, il convient
de s'intéresser d'abord à la variation de cette forme bilinéaire et de l'espace sur lequel
elle est définie.

Fixons un isomorphisme de $K$\nobreakdash-groupes $\tors{2}E_\eta \isoto (\Z/2)^2$.
Il s'en déduit un isomorphisme $\tors{2}\Erond \isoto (\Z/2)^2$ de $\P^1_k$\nobreakdash-schémas en
groupes (par exemple en appliquant aux faisceaux étales associés
le foncteur d'image directe par l'inclusion du point générique), et donc un
isomorphisme de $k$\nobreakdash-schémas en groupes
$\tors{2}\Erond_x \isoto (\Z/2)^2$ pour tout $x \in \P^1_k$.
Utilisant l'additivité de la cohomologie, on obtient des isomorphismes
$H^1(K, \tors{2}E_\eta) \isoto \Hgoth(K)$, $H^1(k, \tors{2}\Erond_x)\isoto \Hgoth(k)$,~etc.
Ils~seront par la suite sous-entendus.

Posons $\Urond = \P^1_{\Orond} \setminus \left( \bigcup_{M \in \Mrond \cup
\{\infty\}} \Mtilde \right)$\glossary{$\Urond$}.  Conformément aux conventions de l'introduction,
nous noterons $\Urond_{\Orond_S}$ le $\Orond_S$\nobreakdash-schéma $\Urond \otimes_\Orond \Orond_S$,
pour $S \subset \Omega$.
Il~est égal à $\Urond \cap \P^1_{\Orond_S}$.
Soit $S_0 \subset \Omega$ un ensemble fini contenant l'ensemble~$S_0$ du paragraphe~\ref{ch1parsymetrisation},
assez grand pour que l'image de $\tors{2}E_\eta(K)$ dans $E_\eta(K)/2 \subset H^1(U,\Z/2)=\Hgoth(U)$
(cf.~lemme~\ref{ch1lemmepeutevaluer})
soit incluse dans $\Hgoth(\Urond_{\Orond_{S_0}})$ et pour que
la conclusion du lemme suivant soit satisfaite.

\bigskip
\begin{lemme}
\label{ch1symfams0}
Il existe un ensemble fini $S_0 \subset \Omega$ tel que pour tout
triplet admissible $(S,\Trond,x)$ avec $S_0 \subset S$,
tout $M \in \Mrond$ et tout $v \in \ft{T}M \cup \{ v_M \}$,
le $\kappa(v)$\nobreakdash-groupe des composantes connexes de la fibre spéciale
du modèle de Néron\index{modèle de Néron} de $\Erond_x\otimes_k k_v$ au-dessus de~$\Orond_v$ soit constant
isomorphe à~$\ff{F}M(L_M)$
et l'isomorphisme canonique $\tors{2}\Erond_x(k)=\tors{2}E_\eta(K)$ fasse se correspondre
l'unique point de $\tors{2}\Erond_x(k) \setminus \{0\}$ qui se spécialise
sur la composante neutre du modèle de Néron de $\Erond_x \otimes_k k_v$
au-dessus de~$\Orond_v$ et l'unique point de $\tors{2}E_\eta(K) \setminus \{0\}$ qui
se spécialise sur~$\Erond^0_M$ (cf.~lemme~\ref{annlemmeuniquespec}).
\end{lemme}

\bigskip
\begin{demo}
La courbe elliptique $E_\eta/K$ admet une équation de Weierstrass minimale de la
forme
\begin{equation}
\label{ch1symfameqw}
Y^2=(X-e_1)(X-e_2)(X-e_3)
\end{equation}
avec $e_1,e_2,e_3 \in k[t]$.  Soit $S_0 \subset \Omega$ l'ensemble
des places finies en lesquelles l'un des coefficients des polynômes~$e_i$ pour $i\in\{1,2,3\}$
n'est pas entier
ou le coefficient dominant de l'un des polynômes $e_i-e_j$ pour $i,j\in\{1,2,3\}$ distincts
n'est pas une unité.  Soit $(S,\Trond,x)$ un triplet admissible tel que $S_0 \subset S$.
La courbe elliptique $\Erond_x/k$ a pour équation de Weierstrass
\begin{equation}
\label{ch1symfameqwx}
Y^2=(X-e_1(x))(X-e_2(x))(X-e_3(x)) \rlap{\text{.}}
\end{equation}
Soient $M \in \Mrond$ et $v \in \ft{T}M \cup \{ v_M \}$.
Comme~$E_\eta$ est à réduction multiplicative en~$M$, un et un seul
des polynômes $e_1-e_2$, $e_2-e_3$ et $e_3-e_1$ s'annule en~$M$.
Quitte à renuméroter les~$e_i$, on peut supposer que $(e_1-e_2)(M)=0$.
Notons~$v$ (resp.~$v_M$) la valuation normalisée de~$k$ (resp.~de~$K$)
associée à la place~$v$ (resp.~au point~$M$).
La définition de~$S_0$ et la préadmissibilité du triplet $(S,\Trond,x)$ entraînent que
les coefficients de l'équation~(\ref{ch1symfameqwx}) sont des entiers $v$\nobreakdash-adiques
et que $v(e_i(x)-e_j(x))=v_M(e_i-e_j)$ pour tous $i,j$ distincts.
Comme $(e_1-e_2)(M)=0$, $(e_2-e_3)(M)\neq 0$
et $(e_3-e_1)(M)\neq 0$,
il en résulte que l'équation de Weierstrass~(\ref{ch1symfameqwx}) est minimale en~$v$,
que la courbe elliptique~$\Erond_x$ est à réduction multiplicative en~$v$,
qu'elle est à réduction multiplicative déployée si et seulement si
l'image de $e_1(x)-e_3(x)$ dans~$\kappa(v)$ est un carré,
que le $\kappa(v)$\nobreakdash-groupe des composantes connexes de la fibre spéciale du
modèle de Néron de $\Erond \otimes_k k_v$ au-dessus de~$\Orond_v$ est constant si $e_1-e_2$ s'annule en~$M$
avec multiplicité~$1$
et enfin qu'il est isomorphe à $\ff{F}M(L_M)$ s'il est constant.
Vérifions maintenant que ce groupe est constant si $v_M(e_1-e_2)>1$.
Sous cette hypothèse,
la classe dans $\kappa(M)^\star/\kappa(M)^{\star 2}$ de l'extension quadratique
ou triviale $L_M/\kappa(M)$ est égale à celle de $(e_1-e_3)(M)$ (cf.~proposition~\ref{ch1expldlm}).
Comme le triplet $(S,\Trond,x)$ est admissible,
l'unique place~$w$ de $\kappa(M)$ divisant~$v$ et appartenant à l'intersection $\xtilde \cap \Mtilde$
est totalement décomposée dans~$L_M$; autrement dit, l'image de $(e_1-e_3)(M)$ dans $\kappa(M)_w^\star$ est un carré.
La courbe elliptique~$\Erond_x$ est donc à réduction
multiplicative déployée en~$v$, compte tenu que $(e_1-e_3)(M) \in \kappa(M)$ est
une unité $w$\nobreakdash-adique et que son image modulo~$w$ est égale
à l'image de $e_1(x)-e_3(x)$ dans~$\kappa(v)$.  L'assertion voulue s'ensuit
(cf.~corollaire~\ref{anncorcstdep}).

L'équation de Weierstrass~(\ref{ch1symfameqw}) (resp.~(\ref{ch1symfameqwx})) étant
minimale, il suffit de la réduire modulo~$M$ (resp.~modulo~$v$) pour
déterminer quel point d'ordre~$2$ se spécialise sur la composante neutre du modèle
de Néron (cf.~preuve du lemme~\ref{ch1expldlemme}).
On voit ainsi que c'est le point de coordonnée $X=e_3(M)$ (resp.~$X=e_3(x)$),
d'où le résultat.
\end{demo}

\bigskip
Choisissons une fois pour toutes une famille
de sous-espaces $(K_v)_{v\in S_0}$ vérifiant la conclusion du lemme~\ref{ch1existekv}.

Soit $(S,\Trond, x)$ un triplet admissible tel que $S_0 \subset S$.
Notons $N = \Irond^{T(x)}(\Erond_x)$, vu comme sous-groupe de $H^1(k,\tors{2}\Erond_x) = \Hgoth(k)$.
Par un léger abus de langage, pour $v \in \Omega_f$, on parlera de la valuation en~$v$
d'un élément de~$\Hgoth(k)$ pour désigner son image dans $(\Z/2)^2$ par
la flèche induite par la valuation normalisée associée à~$v$; de même pour la valuation
en $M \in \Mrond$ d'un élément de~$\Hgoth(U)$.

La valuation d'un élément de~$\Hgoth(k)$ est égale à son
image par la flèche de restriction $\Hgoth(k)=H^1(k, \tors{2}\Erond_x) \rightarrow
H^1(k_v^\nr, \tors{2}\Erond_x)=\Hgoth(k_v^\nr)=(\Z/2)^2$.
La condition en une place $v \in T(x) \setminus S$ sur un $a \in \Hgoth(k)$
pour qu'il appartienne à~$N$ est donc que sa valuation en~$v$
appartienne au sous-groupe $\Wrond_v(\Erond_x) \subset (\Z/2)^2$
introduit au paragraphe~\ref{ch1parsymetrisation}.

Soit $N_1 \subset N$ le sous-groupe constitué des vecteurs $e\in N$
ayant mêmes valuations en~$v$ et en~$v_M$ pour tout $M \in \Mrond$ et
tout $v \in \ft{T}M$.
Nous allons
maintenant montrer que l'on peut injecter~$N$ dans le groupe de $2$\nobreakdash-Selmer
géométrique de~$\Erond$ au-dessus de la droite affine, ce qui permettra ensuite de
comparer les sous-espaces~$N_1$ pour diverses valeurs de~$x$ et~$\Trond$.

\bigskip
\begin{lemme}
\label{ch1evxiso}
L'application $\Hgoth(\Urond_{\Orond_T}) \rightarrow
\Hgoth(\Orond_{T(x)})$ d'évaluation en~$x$ est un isomorphisme.
\end{lemme}

\bigskip
\begin{demo}
Remarquons d'abord que l'application en question est bien définie, car
$\xtilde \cap \P^1_{\Orond_{T(x)}} \subset \Urond_{\Orond_T}$.
Elle s'inscrit dans le diagramme
\begin{equation}
\label{ch1diagpsi}
\begin{aligned}
\myxyhook\xymatrix{
0 \ar[r] & \Hgoth(\Orond_T) \ar@{=}[d] \ar[r] &
\Hgoth(\Urond_{\Orond_T}) \ar[d] \ar[r]^(.45)\alpha
& \displaystyle\prod_{M \in \Mrond} (\Z/2)^2
\ar@{=}[d] \ar[r] & 0 \\
0 \ar[r] & \Hgoth(\Orond_T) \ar[r] & \Hgoth(\Orond_{T(x)})
\ar[r]^(.45)\beta & \displaystyle\prod_{M \in \Mrond} (\Z/2)^2 \ar[r] & 0\rlap{\text{,}}
}
\end{aligned}
\end{equation}
où~$\alpha$ associe à un couple de classes de fonctions rationnelles sur~$\P^1_k$ la famille des
classes modulo~$2$ de leurs ordres d'annulation aux points de~$\Mrond$,
et~$\beta$ est induite par les valuations~$v_M$ sur~$k$.  Ce diagramme commute
en vertu de l'hypothèse de transversalité dans la définition d'un couple préadmissible,
et ses lignes sont exactes parce que
les groupes $\Pic(\A^1_{\Orond_T})$ et~$\Pic(\Orond_T)$ sont nuls (par construction de~$S_0$).
Le lemme des cinq permet de conclure.
\end{demo}

\bigskip
On notera $\psi \colon \Hgoth(\Orond_{T(x)}) \rightarrow
\Hgoth(\Urond_{\Orond_T})$\glossary{$\psi$} l'isomorphisme inverse de l'évaluation en~$x$.
Le groupe $I^{T(x)} \subset H^1(k,\tors{2}E)=\Hgoth(k)$ s'identifie
à $\Hgoth(\Orond_{T(x)})$ puisque $T(x)$ contient un système de générateurs
du groupe de classes de~$k$.  On peut donc considérer l'image de~$N$ par~$\psi$.

\bigskip
\begin{proposition}
\label{ch1nsg}
Pour tout $M\in \Mrond$
et tout $a \in \Hgoth(\Urond_{\Orond_T})$, l'image de~$a$ dans $H^1(K_M^\sh,
\tors{2}E_\eta)=(\Z/2)^2$ appartient au sous-groupe $E_\eta(K_M^\sh)/2$ si et
seulement si l'image de~$a(x)$ par la
valuation~$v_M$ appartient au sous-groupe $\Wrond_{v_M}(\Erond_x) \subset (\Z/2)^2$.
\end{proposition}

\bigskip
\begin{demo}
Vu la commutativité du diagramme~(\ref{ch1diagpsi}),
il suffit de montrer que les sous-groupes $\Wrond_{v_M}(\Erond_x)$
et $E_\eta(K_M^\sh)/2$ de $(\Z/2)^2$ sont égaux, ce qui résulte du lemme suivant.
\end{demo}

\bigskip
\begin{lemme}
\label{ch1lemmekvnrkmsh}
Pour tout $M \in \Mrond$ et tout $v \in \ft{T}M \cup \{ v_M \}$, les sous-groupes $\Wrond_v(\Erond_x)$ et $E_\eta(K_M^\sh)/2$ de
$(\Z/2)^2$ sont égaux et d'ordre~$2$.
\end{lemme}

\bigskip
\begin{demo}
La proposition~\ref{annpropcaracteksh} montre que
$E_\eta(K_M^\sh)/2=\Erond_x(k_v^\nr)/2$ (comme sous-groupes de $(\Z/2)^2$),
compte tenu de la conclusion du lemme~\ref{ch1symfams0}.
Par ailleurs, le groupe $E_\eta(K_M^\sh)/2$ est d'ordre~$2$ d'après
les propositions~\ref{annpropekshsurdeux} et~\ref{annpropcyclpair}.
Comme $\Wrond_v(\Erond_x)$ est inclus dans $\Erond_x(k_v^\nr)/2$ (en tant que sous-groupe de $(\Z/2)^2$),
il reste seulement à vérifier que $\Wrond_v(\Erond_x)\neq 0$;
or ceci découle de la proposition~\ref{ch1wrvneron} puisque
le $\kappa(v)$\nobreakdash-groupe des composantes connexes de la fibre spéciale
du modèle de Néron de $\Erond_x \otimes_k k_v$ sur $\Orond_v$
est un $\kappa(v)$\nobreakdash-groupe constant cyclique d'ordre pair (cf.~proposition~\ref{annpropcyclpair}
et conclusion du lemme~\ref{ch1symfams0}).
\end{demo}

\bigskip
\begin{corollaire}
\label{ch1corpsinsg}
Le sous-espace $\psi(N)$ de $\Hgoth(K)=H^1(K, \tors{2}E_\eta)$ est inclus
dans $\SG_2(\A^1_k, \Erond)$.
\end{corollaire}

\bigskip
\begin{demo}
Étant donné que $\psi(N) \subset \Hgoth(U) \subset H^1(U,\tors{2}\Erond)$,
il suffit d'appliquer la proposition~\ref{ch1nsg}.
\end{demo}

\bigskip
\begin{proposition}
\label{ch1propindep}
Soit $(x_v)_{v \in S} \in \prod_{v \in S} (\P^1_k\setminus \Mrond)(k_v)$. Il existe une famille
$(\Arond_v)_{v \in S}$ de voisinages $v$\nobreakdash-adiques des~$x_v$ telle que
ni le sous-espace $\psi(N_1)$ de $\SG_2(\A^1_k, \Erond)$
ni la forme bilinéaire symétrique
\begin{equation}
\psi(N_1) \times \psi(N_1) \longrightarrow \Z/2
\end{equation}
induite par la restriction à $N_1$ de la forme bilinéaire
sur~$N$ définie par~(\ref{ch1defacc})
ne varient lorsque $(\Trond, x)$ parcourt l'ensemble des couples tels que
le triplet $(S,\Trond,x)$ soit admissible et que $x \in \Arond_v$ pour tout $v \in S$.
\end{proposition}

\bigskip
\begin{remarque}
L'énoncé ci-dessus ne sous-entend pas que l'ensemble parcouru par le couple $(\Trond, x)$
n'est pas vide.
\end{remarque}

\bigskip
\begin{demo}
La préadmissibilité du triplet $(S,\Trond,x)$ entraîne que
$$
\Hgoth(\Urond_{\Orond_S}) = \bigensemble{a \in \Hgoth(\Urond_{\Orond_T})}
{\forall M  \in \Mrond, \; \forall v \in \ft{T}M, \; v(a(x))=v_M(a(x))\in (\Z/2)^2}\!\rlap{\text{.}}
$$
De cette égalité, de la définition de~$N_1$ et du lemme~\ref{ch1lemmekvnrkmsh}, on déduit:
\begin{equation*}
\psi(N_1)= \SG_2(\A^1_k, \Erond) \cap
\bigensemble{a \in \Hgoth(\Urond_{\Orond_S})}
{\forall v\in S, \; a(x)|_{\fv{V}v} \in W_v(\Erond_x) + K_v}\!\rlap{\text{,}}
\end{equation*}
où $a(x)|_{\fv{V}v}$ désigne l'image de~$a(x)$ par la flèche
naturelle $\Hgoth(k) \rightarrow \Hgoth(k_v)=\fv{V}v$.
Le membre de droite est évidemment indépendant
de~$\Trond$, et il ne dépend pas non plus de~$x$ si l'on choisit les
voisinages~$\Arond_v$ assez petits, comme le montre le lemme suivant.

\bigskip
\begin{lemme}
\label{ch1wvconst}
Pour tout $v \in \Omega$, le sous-groupe $W_v(\Erond_{x_v})$ de~$\Hgoth(k_v)$ est une fonction
localement constante de $x_v \in U'(k_v)$, où $U' = \P^1_k\setminus \Mrond$.
\end{lemme}

\bigskip
\begin{demo}
Comme~$\Hgoth(k_v)$ est fini, il suffit de voir que pour tout \mbox{$\alpha \in \Hgoth(k_v)$},
l'ensemble des $x_v \in U'(k_v)$ tels que~$\alpha$ appartienne à $W_v(\Erond_{x_v})$ est
ouvert et fermé dans $U'(k_v)$. Autrement dit, si l'on note
$Y \rightarrow U'_{k_v}$ un torseur sous~$\Erond_{U'_{k_v}}$
dont la classe
dans $H^1(U'_{k_v}, \Erond)$ est l'image de~$\alpha$ par la composée
$$
\Hgoth(k_v) \subset \Hgoth(U'_{k_v}) \subset H^1(U'_{k_v}, \tors{2}\Erond) \longrightarrow
H^1(U'_{k_v}, \Erond) \rlap{\text{,}}
$$
on doit prouver que l'image de l'application $Y(k_v) \rightarrow U'(k_v)$ est ouverte
et fermée.  Ces deux propriétés découlent respectivement de la lissité et de la
propreté du morphisme $Y \rightarrow U'_{k_v}$ (l'image est ouverte
d'après le théorème des fonctions implicites; pour montrer qu'elle est fermée,
on peut remarquer que le morphisme $Y \rightarrow U'_{k_v}$ est projectif, auquel
cas l'assertion est évidente, ou se servir du lemme de Chow pour se ramener à cette
situation).
\end{demo}

\bigskip
Passons maintenant à l'étude de la forme bilinéaire symétrique $\langle \cdot, \cdot \rangle$
sur~$\psi(N_1)$ induite par la forme bilinéaire $N \times N \rightarrow \Z/2$
définie par~(\ref{ch1defacc}).

\bigskip
\begin{lemme}
\label{ch1exisgm}
Pour tout $M \in \Mrond$, il existe $g_M \in \SG_2(\A^1_k, \Erond) \cap \Hgoth(\Urond_{\Orond_S})$
tel que pour tous~$\Trond$
et~$x$ tels que~$(S, \Trond, x)$ soit admissible et tout $v \in \ft{T}M \cup \{v_M\}$,
l'image de $g_M(x)$ dans~$\fv{V}v$ appartienne à~$W_v(\Erond_x)$ mais pas à~$K_v$.
\end{lemme}

\bigskip
\begin{demo}
Soit $M \in \Mrond$.  Supposons d'abord que le groupe $\ff{F}M(L_M)$ ne soit pas d'ordre~$2$.
Comme $\Pic(\Orond_S)=0$, il existe une fonction $f \in \Gm(\Urond_{\Orond_S})$
de valuation~$1$ en~$M$ et de valuation nulle en tout autre point de~$\Mrond$.
Notons $P \in \tors{2}E_\eta(K)$ l'unique point d'ordre~$2$ qui se spécialise
sur~$\Erond^0_M$ (cf.~lemme~\ref{annlemmeuniquespec}) et $g_M \in \Hgoth(\Urond_{\Orond_S})$ l'image de~$f$ par l'application
obtenue
en tensorisant
par $\Gm(\Urond_{\Orond_S})/2$
l'inclusion $\iota \colon
\Z/2 \hookrightarrow \Z/2 \times \Z/2$ de~$P$ dans $\tors{2}E_\eta(K)=\Z/2 \times \Z/2$.

Montrons que $g_M \in \SG_2(\A^1_k, \Erond)$.  La condition aux points de $\Mrond \setminus \{M\}$
est satisfaite puisque l'image de~$g_M$ dans $H^1(K_{M'}^\sh,\tors{2}E_\eta)$ est même nulle
pour $M' \neq M$.  Pour vérifier que l'image de~$g_M$ dans $H^1(K_M^\sh,\tors{2}E_\eta)$ appartient
au sous-groupe $E_\eta(K_M^\sh)/2$, il suffit de constater qu'elle appartient
à l'image de la flèche $H^1(K_M^\sh, \Z/2) \rightarrow H^1(K_M^\sh,\tors{2}E_\eta)$
induite par~$\iota$ (c'est l'image de la classe de~$f$) et d'appliquer la proposition~\ref{annpropcaracteksh}.

Soient~$\Trond$ et~$x$ tels que le triplet $(S, \Trond, x)$ soit admissible.  Soit $v \in \ft{T}M\cup\{v_M\}$.
La commutativité du diagramme~(\ref{ch1diagpsi}) entraîne que les images respectives de~$g_M$ et de~$g_M(x)$
dans $(\Z/2)^2$ par les valuations en~$M$ et en~$v_M$ coïncident.  Compte tenu du lemme~\ref{ch1lemmekvnrkmsh}
et de la définition de~$g_M$, il en résulte que l'image de~$g_M(x)$ dans~$\fv{V}v$
appartient à $W_v(\Erond_x)+K_v$ mais pas à~$K_v$.
Il reste seulement à prouver qu'elle appartient à $W_v(\Erond_x)$,
ou encore, ce qui revient au même, que son image~$\gamma$ dans $H^1(k_v, \Erond_x)$ est nulle.
Notons $\fp{P}x \in \Erond_x(k)$ l'image de~$P$ par l'isomorphisme
canonique $\tors{2}E_\eta(K)=\tors{2}\Erond_x(k)$ et $\phi \colon \Erond_x \rightarrow \Erond''_x$
le quotient de~$\Erond_x$ par~$\fp{P}x$.
Soient $\Erondsouligne_x$ (resp.~$\Erondsouligne''_x$) le modèle de Néron de $\Erond_x \otimes_k k_v$
(resp.~$\Erond''_x \otimes_k k_v$)
au-dessus de~$\Orond_v$
et~$F$ (resp.~$F''$) le $\kappa(v)$\nobreakdash-groupe fini étale des composantes connexes de la fibre spéciale de $\Erondsouligne_x$
(resp.~$\Erondsouligne''_x$).  D'après la propriété universelle des modèles de Néron, l'isogénie~$\phi$
s'étend en un morphisme $\souligne{\phi} \colon \Erondsouligne_x \rightarrow \Erondsouligne''_x$.
L'appartenance de l'image de~$g_M(x)$ dans~$\fv{V}v$ au sous-groupe $W_v(\Erond_x)+K_v$ se traduit
par celle de~$\gamma$ à $H^1(\Orond_v,\Erondsouligne_x)$, d'après la suite spectrale de Leray
pour le faisceau étale $\Erond_x \otimes_k k_v$ sur $\Spec(k_v)$ et l'inclusion du point
générique de $\Spec(\Orond_v)$.
Comme $\ff{F}M(L_M)$ est d'ordre pair (cf.~proposition~\ref{annpropcyclpair}) différent de~$2$,
il résulte de la conclusion
du lemme~\ref{ch1symfams0} que la courbe elliptique~$\Erond_x$ est à réduction multiplicative
déployée en~$v$ (cf.~corollaire~\ref{anncorcstdep}).
La flèche de spécialisation $H^1(\Orond_v,\Erondsouligne_x)
\rightarrow H^1(\kappa(v),F)$ est donc injective (d'après le théorème de Hilbert~90
et le lemme de Hensel, cf.~preuve de la proposition~\ref{annpropekshsurdeux} pour
les détails),
de sorte qu'il suffit, pour conclure,
de prouver que l'image de~$\gamma$ dans $H^1(\kappa(v),F)$ est nulle.

La classe~$\gamma$ appartient à l'image de la flèche $H^1(k_v,\Z/2)\rightarrow H^1(k_v,\Erond_x)$
induite par l'inclusion de~$\fp{P}x$ dans~$\Erond_x$.  Son image par la flèche
$H^1(\Orond_v,\Erondsouligne_x) \rightarrow H^1(\Orond_v,\Erondsouligne''_x)$
induite par~$\souligne{\phi}$ est donc nulle.  Notant~$i$ l'inclusion du point
générique de $\Spec(\Orond_v)$ et $\souligne{\phi}_0 \colon F \rightarrow F''$ le morphisme
déduit de~$\souligne{\phi}$, le carré commutatif
$$
\xymatrix@C=6ex{
\Erondsouligne_x \ar[d] \ar[r]^{\souligne{\phi}} & \Erondsouligne''_x \ar[d] \\
i_\star F \ar[r]^{i_\star\souligne{\phi}_0} & i_\star F''
}
$$
permet d'en déduire que l'image de~$\gamma$ dans $H^1(\kappa(v),F)$ appartient au noyau
de la flèche $H^1(\kappa(v),F) \rightarrow H^1(\kappa(v),F'')$ induite par $\souligne{\phi}_0$.
D'après la proposition~\ref{annpropppp} et la conclusion du lemme~\ref{ch1symfams0},
on a une suite exacte
$$
\xymatrix{
0 \ar[r] & F \ar[r] & F'' \ar[r] & \Z/2 \ar[r] & 0
}
$$
de $\kappa(v)$\nobreakdash-groupes.
Comme $\Erond_x$ est à réduction multiplicative déployée en~$v$ et
de type~$I_n$ avec $n>2$,
les $\kappa(v)$\nobreakdash-groupes de cette suite exacte sont tous constants
(cf.~corollaire~\ref{anncorcorfaux}).
La flèche $H^1(\kappa(v),F) \rightarrow H^1(\kappa(v),F'')$ est donc injective,
ce qui prouve finalement que $\gamma=0$.

Il reste à traiter le cas où $\ff{F}M=\Z/2$.  On choisit alors
$Q \in \tors{2}E_\eta(K) \setminus \{0,P\}$ et l'on définit~$g_M$
comme l'image de~$Q$ dans $E_\eta(K)/2 \subset \SG_2(\A^1_k,\Erond)$.
On a $g_M \in \Hgoth(\Urond_{\Orond_S})$ par définition de~$S_0$.
Il est évident que l'image de~$g_M(x)$ dans $\fv{V}v$ appartient à $W_v(\Erond_x)$
pour tout $x \in U(k)$ et toute place~$v$. Supposons qu'il existe~$\Trond$
et~$x$ tels que le triplet $(S,\Trond,x)$ soit admissible, et une place $v \in \ft{T}M \cup \{v_M\}$
telle que l'image de $g_M(x)$ dans~$\fv{V}v$ appartienne à~$K_v$, autrement dit telle
que l'image de~$g_M(x)$ dans $H^1(k_v^\nr,\tors{2}\Erond_x)$ soit nulle.
Compte tenu de l'inclusion $\Erond_x(k_v^\nr)/2 \subset H^1(k_v^\nr,\tors{2}\Erond_x)$,
l'image de~$Q_x$ dans $\Erond_x(k_v^\nr)/2$ est alors nulle, où
$Q_x \in \tors{2}\Erond_x(k)$ désigne l'image de~$Q$ par l'isomorphisme
canonique $\tors{2}E_\eta(K)=\tors{2}\Erond_x(k)$.
D'après le corollaire~\ref{annpropekshsurdeux},
cela signifie que l'image de~$Q_x$ dans $F(\kappa(v))=F(\kappa(v))/2=\Z/2$ est nulle;
en d'autres termes, le point~$Q_x$ partage avec~$\fp{P}x$ la propriété de
se spécialiser dans $\souligne{\Erond}^0$ (cf.~conclusion du lemme~\ref{ch1symfams0}).
Vu le lemme~\ref{annlemmeuniquespec}, il en résulte que $Q_x \in \{0,\fp{P}x\}$, d'où une contradiction.
\end{demo}

\bigskip
Soient $a, b \in \psi(N_1)$
et $\alpha_v \in W_v(\Erond_x)$, $\alpha'_v \in K_v$ pour $v \in T(x)$ tels
que l'image de $a(x)$ dans $\fv{V}{T(x)}$ s'écrive $\sum_{v \in T(x)}(\alpha_v+\alpha'_v)$.
Posons:
$$
\Mrond(a) = \bigensemble{M \in \Mrond}{a \not\in \Ker\left(\Hgoth(U) \rightarrow
\Hgoth(K_M^\sh)\right)}\!\rlap{\text{.}}
$$
Pour tout $M \in \Mrond$ et tout $v \in \ft{T}M \cup \{v_M\}$, on
a $M \in \Mrond(a)$ si et seulement si $\alpha_{v_M}\not\in K_{v_M}$
(commutativité du diagramme~(\ref{ch1diagpsi})), si et seulement
si $\alpha_v \not\in K_v$ (définition de~$N_1$),
si et seulement si $\alpha_v - g_M(x) \in K_v$
(propriétés de $g_M(x)$ et lemme~\ref{ch1lemmekvnrkmsh}, qui affirme que
$\Wrond_v(\Erond_x)$ est d'ordre~$2$).
Quitte à modifier les~$\alpha'_v$, on peut donc supposer que $\alpha_v=g_M(x)$
pour tout $v \in \bigcup_{M \in \Mrond(a)}(\ft{T}M \cup \{v_M\})$
et que $\alpha_v=0$ pour tout autre $v \in T(x) \setminus S$.
On a alors
\begin{eqnarray*}
\langle a, b \rangle &=& \sum_{v \in T(x)} \langle \alpha_v, b(x) \rangle_v\\
&=& \sum_{v \in S} \langle \alpha_v, b(x) \rangle_v + \sum_{M \in \Mrond(a)} \left(
\sum_{v \in \ft{T}M \cup \{v_M\}} \langle g_M(x), b(x) \rangle_v \right) \!\rlap{\text{.}}
\end{eqnarray*}
Le terme $\sum_{v \in S} \langle \alpha_v, b(x) \rangle_v$ est indépendant de
$(\Trond,x)$ si les~$\Arond_v$ sont assez petits, d'après le
lemme~\ref{ch1wvconst}. Le lemme suivant montre qu'il en va de même pour les
autres termes de la somme
ci-dessus, ce qui termine de prouver la proposition~\ref{ch1propindep}.
\end{demo}

\bigskip
\begin{lemme}
Soient $f, g\in  \SG_2(\A^1_k, \Erond)\cap\Hgoth(\Urond_{\Orond_S})$.  Pour tout $M\in\Mrond$,
$$
\sum_{v \in \ft{T}M \cup \{v_M\}} \langle f(x), g(x) \rangle_v \in \Z/2
$$
est indépendant du couple $(\Trond,x)$ tel que le triplet $(S,\Trond,x)$ soit
admissible et que $x \in \Arond_v$ pour tout $v \in S$.
\end{lemme}

\bigskip
\begin{demo}
La surjectivité de l'application
$$
\Gm(\Urond_{\Orond_S})/2 \longrightarrow \prod_{M \in \Mrond} \Z/2 \rlap{\text{,}}
$$
qui est une conséquence de la nullité de $\Pic(\A^1_{\Orond_S})$,
permet de supposer
que la valuation de~$f$ est nulle en tout point de~$\Mrond$ sauf au plus un,
et de même pour~$g$.
Quitte à remplacer~$f$ ou~$g$ par~$fg$, ce qui ne modifie pas la somme
ci-dessus puisque l'accouplement $\langle \cdot, \cdot \rangle_v$ est alterné,
on peut supposer qu'au moins l'un de~$f$ et de~$g$ est de valuation nulle en~$M$;
en effet, les valuations de~$f$ et de~$g$ en~$M$ sont
égales si elles ne sont pas nulles
puisque $f,g\in\SG_2(\A^1_k,\Erond)$ et que $E_\eta(K_M^\sh)/2$ est d'ordre~$2$.
Enfin, quitte à échanger~$f$ et~$g$, on peut supposer la valuation de~$f$ en~$M$ nulle.

Tout $h \in \Gm(\A^1_{\Orond_S} \setminus \Mtilde)/2$ est représenté par une fonction
régulière sur~$\A^1_{\Orond_S}$, autrement dit par un polynôme $p(t) \in \Orond_S[t]$.
Fixons~$h$ et supposons~$p$ non constant. Alors $\kappa(M)=k[t]/(p(t))$,
et si l'on note $\theta \in \kappa(M)$ l'image de la classe de~$t$
et $u \in \Orond_S$ le coefficient dominant de~$p$,
on a \mbox{$p(t)=u N_{\kappa(M)/k}(t-\theta)$}.  Comme~$\Mtilde$ ne rencontre
pas la section à l'infini au-dessus de $\Orond_S$ (par définition de la
préadmissibilité), le polynôme $N_{\kappa(M)/k}(t-\theta)$ est de valuation
nulle en toute place $v \not\in S$. Par conséquent $u \in \Orond_S^\star$.
Ainsi a-t-on prouvé que \mbox{$\Gm(\A^1_{\Orond_S} \setminus \Mtilde)/2$} est engendré
par $\Gm(\Orond_S)/2$ et par le sous-groupe des normes de~$\kappa(M)$ à~$k$.

Appliquons ceci aux deux composantes de~$g$: il existe $g_1 \in
\Hgoth(\A^1_{\kappa(M)} \setminus \{M\})$ et $u \in \Hgoth(\Orond_S)$ tels que
$g=u N_{\kappa(M)/k}(g_1)$.  On a alors $g(x)=u N_{\kappa(M)/k}(g_1(x))$.
Notons $f(x) \cup g(x)$ la classe de $\tors{2}\Br(k)$ obtenue par
l'accouplement de Weil et le cup-produit des classes de~$f(x)$ et de~$g(x)$
dans $H^1(k,\tors{2}\Erond_x)$.
Par définition de l'invariant de la théorie du corps de classes local,
on a
$\inv_v(f(x) \cup g(x))= \langle f(x),g(x)\rangle_v$. D'autre part,
pour $v \in \ft{T}M \cup \{v_M\}$, on a $\inv_v(f(x) \cup u)=0$ puisque
$\ft{T}v \subset \fv{V}v$ est totalement isotrope et que~$u$ et~$f(x)$ sont de
valuation nulle
en~$v$ (nullité de la valuation de~$f$ en~$M$ et préadmissibilité).  On en
déduit, à l'aide de la formule de projection (cf.~\cite[Ch. XIV,
Ex. 4]{serrecl}) et de la compatibilité de l'invariant à la
corestriction:
$$
\langle f(x), g(x) \rangle_v
= \inv_v \Cores_{\kappa(M)/k}(f(x) \cup g_1(x))
= \sum_{w|v} \inv_w(f(x)\cup g_1(x)) \rlap{\text{,}}
$$
où la seconde somme porte sur les places~$w$ de~$\kappa(M)$ divisant~$v$.
Comme la valuation de~$f$ en~$M$ est nulle et que~$\xtilde$ rencontre
transversalement~$\Mtilde$ au-dessus de~$v$, la classe de $f(x)f(M)$
dans $\Hgoth(\kappa(M)_w)$ est triviale pour toute place~$w$ de~$\kappa(M)$
divisant~$v$. D'où
$$
 \langle f(x), g(x) \rangle_v = \sum_{w|v} \inv_w(f(M) \cup g_1(x))
= \inv_v \Cores_{\kappa(M)/k}(f(M) \cup g_1(x))
$$
pour $v \in \ft{T}M \cup \{v_M\}$.
Les hypothèses sur~$f$, $g$ et~$x$ impliquent que~$f(M)$ et~$g_1(x)$ sont
de valuation nulle en toute place de~$\kappa(M)$ dont la trace sur~$k$
n'appartient pas à $S \cup \ft{T}M \cup \{v_M\}$.  La loi de réciprocité globale\index{loi de réciprocité globale}
permet donc finalement de déduire de l'égalité précédente que
$$
\sum_{v \in \ft{T}M \cup \{v_M\}} \langle f(x), g(x) \rangle_v =
\sum_{v \in S} \inv_v \Cores_{\kappa(M)/k}(f(M) \cup g_1(x)) \rlap{\text{,}}
$$
et cette quantité est clairement indépendante du couple $(\Trond,x)$ choisi
si les voisinages~$\Arond_v$ pour $v \in S$ sont suffisamment petits.
\end{demo}

\bigskip
Nous énonçons pour terminer une conséquence immédiate de la définition de~$N_1$
et du lemme~\ref{ch1lemmekvnrkmsh}:

\bigskip
\begin{lemme}
\label{ch1codimn1n}
La codimension du sous-$\ff{\F}2$\nobreakdash-espace vectoriel~$N_1$ de~$N$ est inférieure ou égale au cardinal
de $T \setminus S$.
\end{lemme}

\bigskip
(On peut prouver que cette inégalité est en fait une égalité, mais cela ne nous servira
pas.)

\subsection{Admissibilité et existence de points adéliques}
\label{ch1parexistadel}
\label{ch1dernierparavantpreuve}

Soit $(S, \Trond, x)$ un triplet préadmissible tel que $X_x(k_v)\neq
\emptyset$ pour toute place $v \in S$.  On cherche ici à relier l'existence de points
adéliques sur la courbe~$X_x$ à des conditions portant uniquement sur les
places de~$\Trond$.  Ceci n'est pas possible en général: d'autres
paramètres sont à prendre en compte, notamment l'obstruction de Brauer-Manin
verticale à l'existence d'un point rationnel sur~$X$, qui est également une
obstruction à l'existence d'un $x \in U(k)$ tel que $X_x(\A_k)\neq \emptyset$.
On va montrer que les hypothèses du théorème~\ref{ch1thprin} permettent d'exprimer
un tel lien à condition que~$B_0$ soit adéquatement choisi.

Pour $M \in \Mrond$, notons $\theta_M \in \kappa(M)$ l'image de~$t$ par le
morphisme \mbox{$k[t] \rightarrow \kappa(M)$} déduit de l'inclusion de~$M$ dans~$\A^1_k$.
Si~$E$ est une extension quadratique ou triviale de~$\kappa(M)$, on pose
$$
A_{E/\kappa(M)} = \Cores_{\kappa(M)(t)/k(t)}(E/\kappa(M), t-\theta_M) \in \tors{2}\Br(U) \rlap{\text{,}}
$$
\glossary{$A_{E/\kappa(M)}$}où $(\cdot,\cdot)$ désigne le symbole de Hilbert sur le corps $\kappa(M)(t)$.
Soit $B_0 \subset \Br(U)$\glossary{$B_0$} le sous-groupe engendré par les classes $A_{K_M/\kappa(M)}$
et $A_{L_M/\kappa(M)}$ pour $M \in \Mrond$ (cf.~lemme~\ref{ch1existekmd} pour la définition
de~$K_M$).

\bigskip
\begin{lemme}
\label{ch1invinfnul}
Soit $v \in \Omega_f$.
Supposons que pour tout $M \in \Mrond$, toute place de~$\kappa(M)$ divisant~$v$
soit totalement décomposée dans $L_MK_M$.
Alors $\inv_v A(x_v)=0$ pour tout $A\in B_0$ et tout $x_v \in U(k_v)$.
\end{lemme}

\bigskip
\begin{demo}
On peut supposer que $A=A_{E/\kappa(M)}$ pour un $M\in\Mrond$ et un $E\in\{L_M,K_M\}$.
On a alors
$$
\inv_v A(x_v)=\sum_{w|v} \inv_w (E/\kappa(M),x_v-\theta_M)\rlap{\text{,}}
$$
où la somme porte sur les places de~$\kappa(M)$ divisant~$v$.
L'image de
la classe de l'extension quadratique ou triviale $E/\kappa(M)$ dans $\kappa(M)_w^\star/\kappa(M)_w^{\star 2}$ est nulle par hypothèse,
d'où le résultat.
\end{demo}

\bigskip
\begin{proposition}
\label{ch1proppreadmadm}
Soit $(S,\Trond,x)$ un triplet préadmissible tel que le couple $(S,\Trond)$ soit
admissible et que
$$
\sum_{v \in S} \inv_v A(x)=0
$$
pour tout $A \in B_0$. Supposons de plus les extensions finies $L_M/k$ et $K_M/k$
non ramifiées hors de~$S$ pour tout $M \in \Mrond$.
Alors le triplet $(S,\Trond,x)$ est admissible.
\end{proposition}

\bigskip
\begin{demo}
On doit prouver que pour tout $M \in \Mrond$, la place~$w_M$ de~$\kappa(M)$
est totalement décomposée dans~$L_M$ et dans~$K_M$. Soit $E \in \{L_M,K_M\}$;
notons $A=A_{E/\kappa(M)}$.
Pour toute place $v \in \Omega$, on a
$$
\inv_v A(x) = \sum_{w|v} \inv_w (E/\kappa(M), x-\theta_M) \rlap{\text{,}}
$$
où la somme porte sur les places~$w$ de~$\kappa(M)$ divisant~$v$.
L'élément $x-\theta_M \in \kappa(M)$ est une unité en toute place~$w$ de~$\kappa(M)$
dont la trace n'appartient pas à~$S$, à l'exception des
places de $\ft{T'}M \cup \{w_M\}$,
en lesquelles $x-\theta_M$ est une uniformisante. Comme $E/\kappa(M)$ est non ramifiée en dehors
des places divisant~$S$ et que~$S$ contient les places dyadiques et les places
archimédiennes,
on en déduit que $\inv_v A(x) = 0$ pour tout $v \in \Omega \setminus (\ft{T}M \cup \{v_M\})$ et que
pour $v \in \ft{T}M \cup \{v_M\}$, $\inv_v A(x)=0$ si et seulement si
l'unique place de $\ft{T'}M \cup \{w_M\}$ divisant~$v$
est totalement décomposée dans~$E$ (cf.~\cite[Prop.~1.1.3]{cscrelle94}).
Cette dernière condition est satisfaite pour toute place
$v \in \ft{T}M$ puisque le couple $(\Trond,x)$ est admissible. L'hypothèse
de la proposition et la loi de réciprocité globale\index{loi de réciprocité globale} permettent d'en déduire que
$\inv_{v_M}A(x)=0$, et donc que~$w_M$ est bien totalement décomposée dans~$E$.
\end{demo}

\bigskip
\begin{proposition}
\label{ch1propexistadel}
Il existe un ensemble fini $S_0 \subset \Omega$ tel que pour tout
triplet admissible $(S, \Trond, x)$ vérifiant $S_0 \subset S$
et $X_x(k_v)\neq \emptyset$
pour tout $v \in S\setminus\Tinfty$, on ait $X_x(\A_k)\neq\emptyset$, en notant~$\Tinfty$
l'ensemble des places de~$S \setminus S_0$ en lesquelles~$x$ n'est pas entier.
\end{proposition}

\bigskip
\begin{demo}
La surface~$X$ est propre et lisse sur~$k$, donc projective, par un
théorème de Zariski (cf.~\cite{zarminmod}).
Il s'ensuit que le morphisme~$\pi$ est projectif et plat, d'où
l'existence d'un ensemble fini $S_0 \subset \Omega$ et d'un
morphisme projectif et plat $g \colon \Xtilde \rightarrow \P^1_{\Orond_{S_0}}$
dont la restriction au-dessus de~$\P^1_k$ est égale à~$\pi$.
Quitte à agrandir~$S_0$, on peut supposer le $\Orond_{S_0}$\nobreakdash-schéma $\Xtilde$
et les fibres de~$g$ au-dessus
du complémentaire de $\P^1_{\Orond_{S_0}} \cap \bigcup_{M \in \Mrond} \Mtilde$
dans $\P^1_{\Orond_{S_0}}$ lisses.

Le lemme \cite[Lemma~1.2]{csscrelle98} (assertions (b) et (c)) montre alors que,
quitte à agrandir encore~$S_0$, on a $X_x(k_v) \neq \emptyset$ pour toute place
$v \in \Omega \setminus S_0$ et tout $x \in U(k)$ pour lesquels
soit~$\xtilde$ ne rencontre au-dessus de~$v$
aucun~$\Mtilde$ pour $M \in \Mrond$,
soit il existe $M \in \Mrond$
et une composante irréductible~$Z$ de~$X_M$ de multiplicité~$1$
tels que~$\xtilde$ rencontre~$\Mtilde$ en une place~$w$ de~$\kappa(M)$ divisant~$v$
et totalement décomposée dans la fermeture algébrique de~$\kappa(M)$ dans~$\kappa(Z)$.
On conclut à l'aide du lemme suivant.
\end{demo}

\bigskip
\begin{lemme}
\label{ch1xmxrondm}
Pour tout $M \in \Mrond$, il existe une composante irréductible~$Z$ de la
fibre~$X_M$ de~$\pi$ en~$M$, de multiplicité~$1$ et telle que la fermeture algébrique
de~$\kappa(M)$ dans~$\kappa(Z)$ se plonge dans~$L_M K_M$.
\end{lemme}

\bigskip
\begin{demo}
Fixons~$M \in \Mrond$.
La suite exacte
$$
\myxyhook\xymatrix{
0 \ar[r] & \Erond^0_M \ar[r] & \Erond_M \ar[r] & \ff{F}M \ar[r] & 0
}
$$
induit une suite exacte
$$
\myxyhook\xymatrix{
H^1(L_M K_M, \Erond^0_M) \ar[r] & H^1(L_M K_M, \Erond_M) \ar[r] & H^1(L_M K_M, \ff{F}M) \rlap{\text{.}}
}
$$
L'image de~$[\Xrond]$ dans~$H^1(L_M K_M, \ff{F}M)$
est triviale (par définition de~$K_M$),
ce qui signifie que les composantes irréductibles
de $\Xrond_M \otimes_{\kappa(M)}L_MK_M$
sont géométriquement irréductibles sur~$L_MK_M$.  En effet, l'image dans $H^1(L_MK_M,\Erond_M)$
d'un élément de $H^1(L_MK_M,\Erond^0_M)$ représenté par un torseur~$T$
sous $\Erond^0_M \otimes_{\kappa(M)}L_MK_M$ est la classe du
torseur $\coprod_{\alpha \in \ff{F}M(L_M)} T(\alpha)$
sous $\Erond_M \otimes_{\kappa(M)}L_MK_M$, où $T(\alpha)$ désigne
le produit contracté sous~$\Erond^0_M$ de~$T$ et de l'image réciproque
de~$\alpha$ par $\Erond_M \rightarrow \ff{F}M$ (qui est un torseur sous~$\Erond^0_M$).

Notons $\Orond_M^\h$ l'anneau local hensélisé de~$\Orond_M$.
On peut relever l'extension finie $L_M K_M$ de~$\kappa(M)$ en une
$\Orond_M^\h$\nobreakdash-algèbre finie étale locale~$R$ de corps résiduel $L_M K_M$
(cf.~\cite[18.1.1]{ega44}).  On a vu que
la fibre spéciale du $R$\nobreakdash-schéma \mbox{$\Xrond \times_{\P^1_k} \Spec(R)$}
admet une composante irréductible géométriquement irréductible et de multiplicité~$1$.
Il en va nécessairement de même pour le $R$\nobreakdash-schéma \mbox{$X \times_{\P^1_k} \Spec(R)$}
d'après \cite[Lemma~1.1]{skodescent}, d'où le résultat.
\end{demo}

\subsection{Fin de la preuve}
\label{ch1parpreuve}

Maintenant que tous les outils nécessaires sont à notre disposition, passons à la preuve du théorème~\ref{ch1thprin}
proprement dite.
Soit $S_0 \subset \Omega$ un ensemble fini contenant les places archimédiennes,
un système de générateurs du groupe de classes de~$k$,
les places divisant un nombre premier inférieur ou égal au degré de~$\Mrond$
vu comme $k$\nobreakdash-schéma fini réduit et
les places finies de~$k$ au-dessus desquelles le morphisme
$$
\bigcup_{M \in \Mrond \cup \{\infty\}} \Mtilde \longrightarrow \Spec(\Orond)
$$
n'est pas étale, assez grand pour que les extensions $L_M/k$ et $K_M/k$ soient non
ramifiées hors de~$S_0$
pour tout $M \in \Mrond$, pour que
le sous-groupe de~$\Hgoth(U)$ image réciproque de $\{0,[\Xrond]\}$
par la flèche de droite de~(\ref{ch1sesgtsg}), fini d'après le théorème
de Mordell-Weil généralisé, soit inclus dans~$\Hgoth(\Urond_{\Orond_{S_0}})$,
pour qu'il existe un schéma abélien au-dessus de
$\left(\P^1_\Orond \setminus \bigcup_{M\in\Mrond}\Mtilde\right)\otimes_\Orond\Orond_{S_0}$
dont la fibre générique soit égale à~$E_\eta$
et pour que les conclusions du lemme~\ref{ch1symfams0} et de la proposition~\ref{ch1propexistadel}
soient satisfaites.
On a alors:

\bigskip
\begin{lemme}
\label{ch1bonnered}
Pour tout triplet préadmissible $(S,\Trond,x)$ avec $S_0\subset S$,
la courbe elliptique~$\Erond_x$ a bonne réduction hors de $T(x)\setminus \Tinfty$,
où~$\Tinfty$ désigne l'ensemble des places de $S\setminus S_0$ en lesquelles~$x$
n'est pas entier.
\end{lemme}

\bigskip
Soit $B_0 \subset \Br(U)$ le sous-groupe fini défini au paragraphe~\ref{ch1parexistadel}.
Soient \mbox{$S_1 \subset \Omega$} fini contenant~$S_0$ et $(x_v)_{v \in S_1}
\in \prod_{v \in S_1}U(k_v)$
satisfaisant aux hypothèses du théorème~\ref{ch1thprin}.
On va prouver l'existence
de $x \in U(k)$ arbitrairement proche de~$x_v$ pour $v \in S_1 \cap \Omega_f$
et arbitrairement grand aux places archimédiennes de~$k$,
tel que $x \in \RD$ si $\Mrond \neq \emptyset$, ou tel que $x \in \RDz$ et que~$x$
soit entier hors de~$S_1$.

Commençons par définir un ensemble fini $\Tinfty \subset \Omega$ disjoint de~$S_1$.
Nous allons prouver simultanément les deux conclusions du théorème, mais avec
des choix différents pour l'ensemble~$\Tinfty$.  Pour prouver l'assertion b) du théorème,
on pose $\Tinfty = \emptyset$. Pour prouver l'assertion a), on
pose $\Tinfty = \{v_\infty\}$, où $v_\infty$ est choisie comme suit.
Considérons le sous-groupe de $\Hgoth(U)$ engendré par $\Hgoth(\Urond_{\Orond_{S_1}})$
et par~$\SDz$.  Il est fini en vertu du théorème des unités de
Dirichlet et du corollaire~\ref{ch1finitude}
lorsque $\Mrond \neq \emptyset$.
Son intersection avec $\Hgoth(k)$ est donc elle aussi finie, et par conséquent incluse
dans~$\Hgoth(\Orond_{S'})$ pour un ensemble $S' \subset \Omega$ fini assez grand.
Le théorème de \smash{\v{C}ebotarev} fournit une infinité de places $v \in \Omega$ telles que
pour tout $M \in \Mrond$, toute place de~$\kappa(M)$ divisant~$v$ soit totalement décomposée
dans $L_M K_M$.  On choisit pour~$v_\infty$ une telle place hors de $S_1 \cup S'$
et l'on note $x_{v_\infty} \in k_{v_\infty}$ l'inverse d'une uniformisante.

\bigskip
Posons $S = S_1 \cup \Tinfty$\glossary{$S$, $S_1$, $\Tinfty$}.
Pour chaque place finie $v \in S$,
fixons un voisinage $v$\nobreakdash-adique arbitrairement petit $\Arond_v$ de
$x_v \in \P^1(k_v)$, suffisamment
petit pour que tout élément de $\Arond_v$
soit l'inverse d'une uniformisante si $v\in \Tinfty$,
pour que $\inv_v A(x)=\inv_v A(x_v)$ pour tout $x \in \Arond_v$ et tout $A \in B_0$,
et pour que $X_x(k_v) \neq \emptyset$ pour tout $x \in \Arond_v$ si $v\in S_1$.
Il est possible de satisfaire cette dernière condition grâce au théorème des fonctions
implicites et à
l'hypothèse que $X_{x_v}(k_v)\neq\emptyset$ pour $v \in S_1 \cap \Omega_f$.
Fixons de même un voisinage $v$\nobreakdash-adique arbitrairement petit $\Arond_v$
de $\infty \in \P^1(k_v)$ pour chaque place complexe $v \in \Omega$,
et un ouvert connexe non majoré arbitrairement petit $\Arond_v \subset U(k_v)$
pour chaque place réelle $v \in \Omega$,
suffisamment petit pour que $X_x(k_v)\neq\emptyset$ pour tout $x \in \Arond_v$.

\bigskip
\begin{lemme}
\label{ch1suminvax}
Soit $x \in U(k)$ appartenant à~$\Arond_v$ pour tout $v \in S$.
Alors $\sum_{v \in S} \inv_v A(x)=0$ pour tout $A \in B_0$.
\end{lemme}

\bigskip
\begin{demo}
On a $\inv_v A(x_v) = \inv_v A(x)$ pour $v \in S \cap \Omega_f$ par définition
des~$\Arond_v$.  Cette égalité vaut aussi pour~$v$ archimédienne,
puisque si~$v$ est réelle, $x_v$ et~$x$ appartiennent à la même composante
connexe de~$U(k_v)$, d'après l'hypothèse du théorème~\ref{ch1thprin} et
la définition de~$\Arond_v$.
Par conséquent
$$
\sum_{v \in S} \inv_v A(x) = \sum_{v \in S_1} \inv_v A(x_v) +
\sum_{v \in \Tinfty} \inv_v A(x_v) \rlap{\text{.}}
$$
Le premier terme est nul par hypothèse; il reste donc seulement à
vérifier que $\inv_v A(x_v)=0$ pour tout $v \in \Tinfty$,
or ceci résulte du lemme~\ref{ch1invinfnul}.
\end{demo}

\bigskip
Notons $\Lrond$\glossary{$\Lrond$} l'ensemble des entiers naturels~$n$ pour lesquels il
existe une famille~$\Trond$ et un point $x \in U(k)$ tels que le
triplet $(S, \Trond, x)$ soit admissible, que $x \in \Arond_v$ pour
tout $v \in S$ et que $\dim_{\ff{\F}2}\Sel_2(k,\Erond_x)=n$.

\bigskip
\begin{proposition}
\label{ch1propschinzel}
Admettons l'hypothèse de\index{hypothèse de Schinzel} Schinzel.  Soit $\Trond=(\ft{T'}M)_{M\in\Mrond}$ une famille
telle que le couple $(S, \Trond)$ soit préadmissible.  Supposons donné, pour
tout $M\in \Mrond$ et tout $v \in \ft{T}M$, un $x_v \in \A^1(\Orond_v)$ rencontrant
transversalement $\Mtilde$ en l'unique place de~$\ft{T'}M$ qui divise~$v$.  Alors il
existe $x \in U(k)$ arbitrairement proche de~$x_v$ pour $v \in T \cap \Omega_f$
et arbitrairement grand aux places archimédiennes de~$k$, tel que
le triplet $(S, \Trond, x)$ soit préadmissible.
\end{proposition}

\bigskip
\begin{demo}
Comme $\Pic(\Orond_S)=0$, il existe une famille de polynômes
irréductibles $(r_M(t))_{M \in \Mrond}$ de~$k[t]$ vérifiant~$r_M(M)=0$ et~$v(r_M)=0$
pour tout $M\in\Mrond$ et tout $v \in \Omega\setminus S$.
L'hypothèse de Schinzel fournit alors (cf.~\cite[Prop.~4.1]{cscrelle94}) un $x \in
U(k)$ arbitrairement proche de~$x_v$ pour $v \in T \cap \Omega_f$,
arbitrairement grand aux places archimédiennes et entier hors de~$T$, tel que
$\xtilde \cap \Mtilde \cap \P^1_{\Orond_T}$ soit le spectre d'un corps pour
tout $M \in \Mrond$.
(Remarquer que l'on peut supprimer le mot «~perhaps~» dans l'énoncé de
l'hypothèse $(\mathrm{H}_1)$ de~\cite{cscrelle94}: il suffit de choisir~$m$
assez grand à la fin de la preuve de \cite[Prop.~4.1]{cscrelle94}.)
Le triplet $(S,\Trond,x)$ sera bien préadmissible si l'on choisit~$x$
suffisamment proche de~$x_v$ pour $v \in T \setminus S$.
\end{demo}

\bigskip
La proposition~\ref{ch1propschinzel} appliquée à la famille $\Trond=(\emptyset)_{M\in\Mrond}$
prouve l'existence d'un $x \in U(k)$ appartenant à~$\Arond_v$ pour tout $v\in S$,
tel que le triplet $(S, \Trond, x)$ soit préadmissible.
Le triplet $(S, \Trond, x)$ est même admissible d'après la proposition~\ref{ch1proppreadmadm},
applicable grâce au lemme~\ref{ch1suminvax};
d'où $\Lrond \neq \emptyset$. L'ensemble $\Lrond$ possède donc un minimum.
Soient $\Trond=(\ft{T'}M)_{M\in\Mrond}$ et $x \in U(k)$ tels que le triplet $(S,\Trond,x)$ soit
admissible et que $\dim_{\ff{\F}2}\Sel_2(k,\Erond_x)=\min \Lrond$, avec $x \in \Arond_v$ pour tout $v \in S$.
On a $X_x(k_v) \neq \emptyset$ pour $v \in S_1$ par construction des voisinages~$\Arond_v$.
La conclusion de la proposition~\ref{ch1propexistadel} permet d'en déduire que~$x \in \RA$.

\bigskip
\begin{proposition}
\label{ch1propcentrale}
On a $x \in \RDz$.  De plus, si $\Tinfty \neq \emptyset$, on a $x \in \RD$.
\end{proposition}

\bigskip
La démonstration de la proposition~\ref{ch1propcentrale} va nous occuper jusqu'à la fin de ce paragraphe.

\bigskip
\begin{demo}
Reprenons les notations du paragraphe~\ref{ch1parsymfam}. Fixons
la famille de sous-espaces $(K_v)_{v \in S_0}$ et un isomorphisme $\tors{2}E_\eta \isoto (\Z/2)^2$.
L'injectivité de la restriction à l'image réciproque de $\{0,[\Xrond]\}$ par la flèche
de droite de~(\ref{ch1sesgtsg}) de la composée
$$
\SDz \subset \SG_2(\A^1_k,\Erond) \subset H^1(U, \tors{2}\Erond) \rightarrow
H^1(k, \tors{2}\Erond_x)
$$
est une conséquence immédiate du lemme~\ref{ch1evxiso} et de la définition de~$S_0$.

Soit $\alpha \in \Sel_2(k, \Erond_x)$.
La courbe elliptique $\Erond_x$ a bonne réduction hors de~$T(x)$ (lemme~\ref{ch1bonnered}).
Par conséquent $\Sel_2(k,\Erond_x) \subset \Hgoth(\Orond_{T(x)})$, et l'on peut donc
considérer l'image de~$\alpha$ par l'isomorphisme $\psi \colon \Hgoth(\Orond_{T(x)})
\isoto \Hgoth(\Urond_{\Orond_T})$ défini après le lemme~\ref{ch1evxiso}.
On a $\psi(\alpha) \in \SG_2(\A^1_k, \Erond)$
par le corollaire~\ref{ch1corpsinsg}.

\bigskip
\begin{lemme}
\label{ch1lemmepeuttuer}
Pour tout $e \in \tors{2}H^1(\A^1_k,\Erond) \setminus \TDz$, il existe $M \in \Mrond$
et une infinité de places $v \in \Omega \setminus S_0$
telles qu'il existe une place~$w$ de~$\kappa(M)$
totalement décomposée dans $L_MK_M$, divisant~$v$, non ramifiée et de degré résiduel~$1$
sur~$v$, telle que pour tout $x_v \in \A^1(\Orond_v)$ rencontrant transversalement~$\Mtilde$
en~$w$, l'image de~$e$ par la flèche $H^1(\A^1_k,
\Erond) \rightarrow H^1(k_v, \Erond_{x_v})$ d'évaluation en~$x_v$ soit non nulle.
\end{lemme}

\bigskip
\begin{demo}
Soient $\Yrond \rightarrow \A^1_k$ un torseur sous $\Erond_{\A^1_k}$ dont
la classe dans $H^1(\A^1_k, \Erond)$ soit égale à~$e$
et $Y \rightarrow \A^1_k$ un modèle propre et plat régulier minimal de~$\Yrond_\eta$
sur~$\A^1_k$.
On peut étendre
$Y \rightarrow \A^1_k$ en un morphisme propre et plat
$g \colon \Ytilde \rightarrow \A^1_V$ avec~$\Ytilde$ régulier, où~$V$ est un ouvert
dense de $\Spec(\Orond_{S_0})$
(cf.~\cite[8.8.2, 9.6.1 et~11.1.1]{ega43} et \cite[6.12.6]{ega42}).  Comme le modèle $Y \rightarrow \A^1_k$
est minimal, on peut supposer,
quitte à choisir~$V$ suffisamment petit,
qu'il existe un fermé $F \subset \Ytilde$ étale sur~$g(F)$
tel que la restriction de~$g$ à $\Ytilde \setminus F$ soit lisse et que
$g(F) \subset \A^1_V \cap \bigcup_{M \in \Mrond}\Mtilde$.

L'hypothèse que $e \not\in \TDz$ se traduit par l'existence d'un
$M \in \Mrond$ tel que $\delta_M(e)$ n'appartienne pas à $\{0,\delta_M([\Xrond])\}$.
Soient~$d$ l'image de~$e$ dans $\tors{2}H^1(\kappa(M),\ff{F}M)$
et $K_{M,d}/\kappa(M)$ une extension quadratique vérifiant
les conditions du lemme~\ref{ch1existekmd}.
Étant donné que l'extension $L_M K_M/L_M$ est quadratique
ou triviale et que le groupe $\ff{F}M(L_M K_M)=\ff{F}M(L_M)$ est cyclique, la suite exacte
d'inflation-restriction montre que
le noyau de la flèche de restriction $H^1(L_M, \ff{F}M) \rightarrow H^1(L_M K_M, \ff{F}M)$
est d'ordre au plus~$2$.
Notons~$N$ ce noyau. Si $N\neq 0$, alors $L_MK_M \neq L_M$,
de sorte que la définition de~$K_M$ entraîne que $\delta_M([\Xrond]) \in N \setminus\{0\}$;
le groupe~$N$
est donc dans tous les cas engendré par~$\delta_M([\Xrond])$.
Par conséquent, l'image de$~e$ dans $H^1(L_MK_M,\ff{F}M)$ n'est pas nulle.
L'extension quadratique $K_{M,d}/\kappa(M)$ ne se plonge donc pas dans $L_MK_M/\kappa(M)$.
On en déduit, à l'aide du théorème de \smash{\v{C}ebotarev},
l'existence d'une infinité de places finies
de $L_MK_M$ non ramifiées et de degré résiduel~$1$ sur~$k$,
inertes dans $L_M K_M K_{M,d}$
(cf.~\cite[Proposition~2.2]{harduke}). Soient~$w$ la trace sur~$\kappa(M)$ d'une telle place
et~$v \in \Omega_f$ la trace sur~$k$ de~$w$. Supposons que $v \in V$ et montrons que~$v$
convient.

La place~$w$ est totalement décomposée dans $L_M K_M$, inerte dans $K_{M,d}$,
non ramifiée et de degré résiduel~$1$ sur~$v$.  Supposons qu'il existe
$x_v \in \A^1(\Orond_v)$ rencontrant transversalement $\Mtilde$ en~$w$ et tel que
l'image de~$e$ dans $H^1(k_v,\Erond_{x_v})$ soit nulle, autrement dit tel
que $\Yrond_{x_v}(k_v) \neq \emptyset$, et aboutissons à une contradiction.

\bigskip
\begin{lemme}
\label{ch1ylyrondl}
Pour toute extension $\ell/k$ et tout $m \in \A^1(\ell)$,
$\fy{Y}m$ possède un $\ell$\nobreakdash-point régulier
si et seulement si $\Yrond_m(\ell)\neq\emptyset$.
\end{lemme}

\bigskip
\begin{demo}
Notons~$R$ le hensélisé de l'anneau local de~$\A^1_\ell$ en~$m$ et~$Q$
son corps des fractions.
D'après la propreté de $Y \rightarrow \A^1_k$, la régularité de~$Y$,
la lissité de $\Yrond \rightarrow \A^1_k$ et
l'injectivité de la flèche de restriction
$H^1(R, \Erond) \rightarrow H^1(Q, E_\eta)$,
les conditions de l'énoncé sont toutes les deux équivalentes
à l'existence d'un $Q$\nobreakdash-point sur~$\Yrond_\eta$.
(L'injectivité de la flèche de restriction découle de la suite spectrale de Leray
pour l'inclusion du point générique de~$R$ et le faisceau étale~$E_\eta$ sur $\Spec(Q)$,
une fois que l'on sait que $\Erond \times_{\P^1_k} \Spec(R)$ est un modèle de Néron
de $E_\eta \otimes_K Q$. Pour cela, cf.~\cite[7.2/2]{blr}.)
\end{demo}

\bigskip
Notons $z \in \A^1_V$ l'image du point fermé de $\Spec(\Orond_v)$ par
le morphisme $x_v \colon \Spec(\Orond_v) \rightarrow \A^1_V$.
Le lemme~\ref{ch1ylyrondl} montre que $\fy{Y}{x_v}(k_v) \neq \emptyset$. Comme $\fy{\Ytilde}{x_v}$
est propre sur $\Spec(\Orond_v)$, on en déduit que $\fy{\Ytilde}{x_v}(\Orond_v)\neq\emptyset$,
puis que~$\fy{\Ytilde}z$ possède un $\kappa(v)$\nobreakdash-point régulier, grâce au lemme suivant.

\bigskip
\begin{lemme}
Le schéma $\fy{\Ytilde}{x_v}$ est régulier.
\end{lemme}

\bigskip
\begin{demo}
Si~$S$ est un schéma et $s \in S$, notons $\ft{T}sS$ l'espace tangent à~$S$ en~$s$;
c'est par définition le $\kappa(s)$\nobreakdash-espace vectoriel dual de $\mgoth_s/\mgoth_s^2$.

Soit $p \in \fy{\Ytilde}{x_v}$.  On va montrer que $\fy{\Ytilde}{x_v}$ est régulier en~$p$,
c'est-à-dire que $\dim_{\kappa(p)}(\ft{T}p \fy{\Ytilde}{x_v})=\dim(\Orond_{\fy{\Ytilde}{x_v},p})$.
Comme la restriction de~$g$ à $\Ytilde \setminus F$ est lisse et que $\Spec(\Orond_v)$
est régulier, on peut supposer que $p \in \fy{\Ytilde}z \cap F$
(cf.~\cite[6.5.2]{ega42}).
Le morphisme $\ft{T}p F \rightarrow \ft{T}z \Mtilde
\otimes_{\kappa(z)}\kappa(p)$ induit par~$g$ est un isomorphisme,
puisque~$F$ est étale sur~$g(F)$.  Le carré commutatif
$$
\myxyhook\xymatrix{
\ft{T}p F \ar[d]^(0.44)\wr \ar@{^{ (}->}[r] & \ft{T}p \Ytilde \ar[d]^{\ft{T}p g} \\
\ft{T}z \Mtilde \otimes_{\kappa(z)} \kappa(p) \ar@{^{ (}->}[r] & \ft{T}z \A^1_V
\otimes_{\kappa(z)} \kappa(p)
}
$$
permet d'en déduire que l'image $I$
de~$\ft{T}pg$ contient $\ft{T}z \Mtilde\otimes_{\kappa(z)} \kappa(p)$.

Les sous-schémas fermés~$x_v$ et~$\Mtilde \otimes_\Orond \Orond_v$
de~$\P^1_{\Orond_v}$ se rencontrent transversalement en~$z$, par hypothèse.
Ceci se traduit par l'égalité $\ft{T}z x_v + \ft{T}z \Mtilde = \ft{T}z \A^1_V$
(moyennant le léger abus de notation consistant à identifier $\ft{T}z \A^1_V$
et $\ft{T}z \A^1_{\Orond_v}$), d'où l'on tire, d'après ce qui précède:
$\ft{T}zx_v\otimes_{\kappa(z)} \kappa(p) + I = \ft{T}z \A^1_V \otimes_{\kappa(z)}\kappa(p)$.
Le morphisme
$$
\myxyhook\xymatrix@C=6ex{
\ft{T}p \Ytilde \ar[r] & \left( \ft{T}z\A^1_V/\ft{T}zx_v\right) \otimes_{\kappa(z)}\kappa(p)
}
$$
induit par~$\ft{T}pg$
est donc surjectif, d'où une suite exacte de $\kappa(p)$\nobreakdash-espaces vectoriels
$$
\myxyhook\xymatrix{
0 \ar[r] & \ft{T}p \fy{\Ytilde}{x_v} \ar[r] & \ft{T}p \Ytilde \ar[r] &
\left( \ft{T}z\A^1_V/\ft{T}zx_v\right) \otimes_{\kappa(z)}\kappa(p) \ar[r] & 0 \rlap{\text{.}}
}
$$
On en déduit l'égalité $$\dim_{\kappa(p)}\ft{T}p\fy{\Ytilde}{x_v}=
\dim_{\kappa(p)}\ft{T}p \Ytilde - \dim_{\kappa(z)} \ft{T}z \A^1_V + \dim_{\kappa(z)} \ft{T}zx_v\rlap{\text{.}}$$
Les schémas $\A^1_V$, $\Ytilde$ et $x_v\simeq \Spec(\Orond_v)$ étant tous réguliers,
il en résulte que $\dim_{\kappa(p)}\ft{T}p\fy{\Ytilde}{x_v} =
\dim (\Orond_{\Ytilde,p}) - \dim (\Orond_{\A^1_V,z}) + \dim(\Orond_{x_v,z})$.
Enfin, grâce à la platitude de~$g$, le membre de droite de cette dernière équation
est égal à $\dim(\Orond_{\fy{\Ytilde}{x_v},p})$ (cf.~\cite[14.2.1]{ega43}, appliqué deux fois),
d'où le résultat recherché.
\end{demo}

\bigskip
Comme $\fy{\Ytilde}z$ est aussi la fibre en~$z$ de $\fy{\Ytilde}{\Mtilde}$, le lemme
de Hensel permet de relever un $\kappa(v)$\nobreakdash-point régulier de~$\fy{\Ytilde}z$ en un
$\kappa(M)_w$\nobreakdash-point régulier de~$\fy{Y}M$.  D'où
$\Yrond_M(\kappa(M)_w)\neq\emptyset$, compte tenu du lemme~\ref{ch1ylyrondl}.
L'image de~$e$ dans $H^1(\kappa(M)_w,\ff{F}M)$ est donc nulle.  La place~$w$
étant totalement décomposée dans~$L_M$, l'extension $L_M/\kappa(M)$ se plonge
dans $\kappa(M)_w/\kappa(M)$; par conséquent, vu la définition de~$K_{M,d}$,
il en va de même pour l'extension $K_{M,d}/\kappa(M)$, ce qui contredit
l'hypothèse selon laquelle~$w$ est inerte dans~$K_{M,d}$.
\end{demo}

\bigskip
Supposons que $\psi(\alpha) \not\in \SDz$.  On peut alors appliquer
le lemme~\ref{ch1lemmepeuttuer} à l'image~$e$ de~$\psi(\alpha)$ dans $\tors{2}H^1(\A^1_k,\Erond)$,
d'où
un point $M_0 \in \Mrond$, une place $v_0 \in \Omega \setminus T$ et une place~$w_0$
de~$\kappa(M_0)$
vérifiant la conclusion du lemme.
Soit $\Trond^+=(\ft{T'^+}M)_{M \in \Mrond}$ la famille définie par
$\ft{T'^+}M=\ft{T'}M$ pour $M \neq M_0$ et \mbox{$\ft{T'^+}{M_0}=\ft{T'}{M_0}\cup\{w_0\}$}.
Posons $S^\sharp=T$.
Soit $\Trond^\sharp=(\ft{T'^\sharp}M)_{M \in \Mrond}$ la famille
définie par $\ft{T'^\sharp}M=\emptyset$ pour tout~$M$.
Soit enfin $\Trond^{\sharp +}=(\ft{T'^{\sharp +}}M)_{M \in \Mrond}$ la famille
définie par $\ft{T'^{\sharp +}}M = \emptyset$ pour $M \neq M_0$
et $\ft{T'^{\sharp +}}{M_0}=\{w_0\}$.

Les couples $(S, \Trond^+)$ et $(S^\sharp, \Trond^{\sharp +})$ et le triplet
$(S^\sharp, \Trond^\sharp, x)$ sont admissibles, puisque~$w_0$ est totalement
décomposée dans $L_MK_M$.  Fixons $x_{v_0} \in \A^1(\Orond_{v_0})$ rencontrant
transversalement $\Mtildez$ en~$w_0$.  Un tel~$x_{v_0}$ existe car la
place~$w_0$ est non ramifiée et de degré résiduel~$1$ sur~$v_0$.  D'après la
proposition~\ref{ch1propschinzel}, il existe $x^+ \in U(k)$ arbitrairement
proche de~$x$ pour $v \in T\cap\Omega_f$ et arbitrairement grand aux places archimédiennes,
tel que le triplet $(S,\Trond^+, x^+)$ soit
préadmissible.  On peut notamment supposer que $x^+ \in \Arond_v$ pour tout $v
\in S$, auquel cas on a $\inv_v A(x^+)=\inv_v A(x)$ pour tout $A \in B_0$ et
tout $v \in S$.  La proposition~\ref{ch1proppreadmadm} permet d'en déduire que
le triplet $(S, \Trond^+, x^+)$ est admissible.  Par conséquent
$\dim_{\ff{\F}2}\Sel_2(k,\Erond_{x^+}) \in \Lrond$, d'où
\begin{equation}
\label{ch1selxselxplus}
\dim_{\ff{\F}2}\Sel_2(k,\Erond_x) \leq \dim_{\ff{\F}2}\Sel_2(k,\Erond_{x^+}) \rlap{\text{.}}
\end{equation}
Reprenons les notations du paragraphe~\ref{ch1parsymfam} associées
aux diverses familles que l'on vient de construire, en les assortissant
d'un $+$ et/ou d'un $\sharp$ en exposant.
On dispose ainsi d'un isomorphisme
$\psi^+ \colon \Hgoth(\Orond_{T^+(x^+)}) \isoto \Hgoth(\Urond_{\Orond_{T^+}})$,
d'espaces vectoriels $N_1$, $N_1^+$, $N_1^\sharp$, $N_1^{\sharp +}$,~etc.
Quitte à choisir~$x^+$ suffisamment proche de~$x$ aux places de $T \cap \Omega_f$ et
suffisamment grand aux places archimédiennes,
la proposition~\ref{ch1propindep} permet de supposer que
$\psi(N_1^\sharp)=\psi^+(N_1^{\sharp +})$ et que les accouplements
symétriques induits sur cet espace par les deux formes bilinéaires
$N_1^\sharp \times N_1^\sharp \rightarrow \Z/2$
et $N_1^{\sharp +} \times N_1^{\sharp +} \rightarrow \Z/2$
coïncident.

Notons $R=\psi(N^\sharp)$ et $R^+=\psi^+(N^{\sharp +})$. On a $N_1^\sharp=N^\sharp$
puisque $\ft{T}M^\sharp=\emptyset$ pour tout $M \in \Mrond$, d'où
$R=\psi^+(N_1^{\sharp +}) \subset R^+$.  La codimension de~$R$ dans~$R^+$
est au plus~$1$, comme le montre le lemme~\ref{ch1codimn1n}.
Notons $\phi \colon R^+ \times R^+ \rightarrow \Z/2$ la forme bilinéaire
symétrique donnée sur~$R^+$ et $\phi|_R$ sa restriction à $R \times R$.
La proposition~\ref{ch1noyauselmer} et le lemme~\ref{ch1bonnered}
montrent que les noyaux de $\phi$ et de $\phi|_R$ sont respectivement égaux
à $\psi^+(\Sel_2(k,\Erond_{x^+}))$ et à $\psi(\Sel_2(k,\Erond_x))$.
Par définition de~$v_0$, le $2$\nobreakdash-revêtement de~$\Erond_{x^+}$ déterminé
par la classe $\psi(\alpha)(x^+) \in H^1(k, \tors{2}\Erond_{x^+})$ ne possède
pas de $k_{v_0}$\nobreakdash-point.  Par conséquent $\psi(\alpha)$ n'appartient pas au noyau de~$\phi$,
bien qu'il appartienne au noyau de~$\phi|_R$.

Résumons la situation: on dispose d'un $\ff{\F}2$\nobreakdash-espace vectoriel de dimension finie~$R^+$,
d'un sous-espace $R \subset R^+$ de codimension au plus~$1$
et d'une forme bilinéaire symétrique~$\phi$ sur~$R^+$;
et l'on sait que le noyau de la restriction $\phi|_R$ n'est pas inclus dans
celui de~$\phi$. Il est alors automatique
que le noyau de~$\phi$ est inclus dans~$R$, et donc strictement inclus dans le
noyau de~$\phi|_R$.  Ceci contredit l'inégalité~(\ref{ch1selxselxplus}); ainsi
a-t-on prouvé, par l'absurde, que $\psi(\alpha) \in \SDz$.

La première partie de la proposition est maintenant établie.  Pour la seconde,
supposons que $\Tinfty \neq \emptyset$.  Il reste à démontrer
que $\psi(\alpha) \in \SD$.
Il suffit pour cela de vérifier que l'image de $\psi(\alpha)$ dans
$H^1(K_\infty^\sh, \tors{2}E_\eta)$ est nulle, autrement dit que
la valuation de $\psi(\alpha)$ au point $\infty \in \P^1_k$ est nulle.

\bigskip
\begin{lemme}
\label{ch1csqvinfty}
On a l'inclusion
$\SDz \cap \Hgoth(\Urond_{\Orond_T}) \subset \Hgoth(\Urond_{\Orond_{T \setminus \{v_\infty\}}})$.
\end{lemme}

\bigskip
\begin{demo}
Soit $a \in \SDz \cap \Hgoth(\Urond_{\Orond_T})$.
Comme l'application
$$
\Gm(\Urond_{\Orond_{S_0}})/2 \longrightarrow \prod_{M \in \Mrond} \Z/2
$$
est surjective (d'après la nullité de $\Pic(\A^1_{\Orond_{S_0}})$),
il existe $b \in \Hgoth(\Orond_T)$ et $c \in \Hgoth(\Urond_{\Orond_{S_0}})$
tels que $a=bc$.  Rappelons que tout élément de~$\Hgoth(k)$ appartenant au sous-groupe
de~$\Hgoth(K)$ engendré par~$\SDz$ et par $\Hgoth(\Urond_{\Orond_{S_1}})$ est de
valuation nulle en~$v_\infty$, par construction de~$v_\infty$.  En particulier,
comme $b=ac$, on a $b \in \Hgoth(\Orond_{T \setminus \{v_\infty\}})$, d'où
$a \in \Hgoth(\Urond_{\Orond_{T \setminus \{v_\infty\}}})$.
\end{demo}

\bigskip
L'intersection de~$\xtilde$ avec la section à l'infini au-dessus de~$v_\infty$
est transverse puisque $x \in \Arond_{v_\infty}$.
Le lemme~\ref{ch1csqvinfty} permet
d'en déduire l'égalité $$v_\infty(\psi(\alpha))=v_\infty(\psi(\alpha)(x))\rlap{\text{,}}$$ où
$v_\infty$ désigne dans le premier membre la valuation normalisée de~$K$ associée
au point $\infty \in \P^1_k$ et dans le second membre la valuation normalisée de~$k$
associée à la place~$v_\infty$.
La courbe~$\Erond_x$ a bonne réduction en~$v_\infty$ (lemme~\ref{ch1bonnered}).
On a par conséquent $v_\infty(\psi(\alpha)(x))=v_\infty(\alpha)=0$,
d'où le résultat.
Ainsi la proposition~\ref{ch1propcentrale} est-elle établie.
\end{demo}

\bigskip
Le théorème~\ref{ch1thprin} est maintenant prouvé.
En effet, $x$ est entier en dehors de~$S_1$ lorsque $\Tinfty=\emptyset$,
puisqu'il est entier en dehors de~$S$ (préadmissibilité du triplet $(S,\Trond,x)$)
et que $S=S_1$ si $\Tinfty=\emptyset$.

\section{\CondD{}, groupe~$\Drond$ et groupe de Brauer}
\label{ch1parcondd}

Les auteurs\index{condition (D)@\condD{}|(}
de~\cite{css} ont montré que la condition qu'ils notent~\cD{}
est étroitement liée à la nullité de la $2$\nobreakdash-torsion du groupe de Brauer
horizontal de la fibration $\pi \colon X \rightarrow C$ (cf.~\cite[§4]{css}).
Pour ce faire, ils ont défini en toute généralité (c'est-à-dire sans hypothèse
sur le type de réduction de~$E_\eta$) un groupe~$\Drond$, puis ont d'une part
étudié sa relation avec le groupe de Brauer horizontal et d'autre part
traduit leur \condD{}, sous certaines hypothèses, en termes du groupe~$\Drond$.

Nous allons ici rappeler la définition du groupe~$\Drond$ et préciser quelques-uns
des énoncés qui viennent d'être mentionnés.  Cela nous permettra de comparer notre
\condD{} avec le groupe~$\Drond$ en toute généralité
et avec la \condD{} de~\cite{css} sous les hypothèses de~\cite{css}.
Cela servira aussi à fixer les notations en vue de
l'utilisation de ces résultats au paragraphe~\ref{ch1parsecdesc}.

Reprenons les notations du paragraphe~\ref{ch1hypnot} et supposons que~$k$ soit
un corps de nombres.
Pour $M \in C$, notons~$D_M$\glossary{$D_M$, $\Delta_M$, $i_M$, $j$}
le groupe abélien libre sur l'ensemble des composantes
irréductibles de la fibre géométrique de~$\pi$ en~$M$ et~$\Delta_M$ le quotient de~$D_M$
par la classe de la fibre géométrique tout entière.
Ces groupes abéliens sont naturellement
munis d'une action continue du groupe de Galois absolu de~$\kappa(M)$.  On pourra donc
les considérer comme des faisceaux étales sur $\Spec(\kappa(M))$.
Notons enfin $i_M \colon \Spec(\kappa(M)) \rightarrow C$ (resp.~$j \colon \eta \rightarrow C$)
l'inclusion de $M\in\Mrond$ (resp.~du point générique de~$C$).

\bigskip
\begin{lemme}
\label{ch1lemmedeltam}
On a une suite exacte canonique
\begin{equation}
\label{ch1dsedeltam}
\myxyhook\xymatrix{
0 \ar[r] & \displaystyle \smash[b]{\prod_{M \in \Mrond} i_{M\star}\Delta_M} \ar[r] &
\PPic_{X/C} \ar[r] & j_\star \PPic_{X_\eta/\eta} \ar[r] & 0
}
\end{equation}
de faisceaux étales sur~$C$.
\end{lemme}

\bigskip
\begin{demo}
Il s'agit essentiellement de vérifier que pour tout $M \in C$, il existe une
suite exacte canonique de groupes abéliens
$$
\myxyhook\xymatrix{
0 \ar[r] & \Delta_M \ar[r] & \Pic(X_{\Orond_M^\sh}) \ar[r] & \Pic(X_{K_M^\sh}) \ar[r] & 0 \rlap{\text{.}}
}
$$
La régularité de~$X$ assure que la flèche de gauche est bien définie et que la flèche de droite
est surjective.  L'injectivité de la flèche de gauche est une conséquence de la propreté de~$\pi$.
Enfin, l'exactitude au milieu est évidente.
\end{demo}

\bigskip
\begin{lemme}
\label{ch1dh1h1}
On a canoniquement $H^1(C, j_\star \PPic_{X_\eta/\eta})=H^1(C,\Erond)/\langle[\Xrond]\rangle$.
\end{lemme}

\bigskip
\begin{demo}
La suite exacte
$$
\myxyhook\xymatrix{
0 \ar[r] & E_\eta \ar[r] & \PPic_{X_\eta/\eta} \ar[r] & \Z \ar[r] & 0
}
$$
induit un isomorphisme $H^1(K,\PPic_{X_\eta/\eta})=H^1(K,E_\eta)/\langle [\Xrond]\rangle$.
Comme les fibres de~$\pi$ sont réduites, on a \mbox{$X_M(K_M^\sh)\neq\emptyset$} pour tout $M\in C$.
La suite exacte ci-dessus induit donc aussi
un isomorphisme
$H^1(K_M^\sh,\PPic_{X_\eta/\eta})=H^1(K_M^\sh,E_\eta)$.
Les suites exactes de bas degré issues des
suites spectrales de Leray pour le morphisme~$j$ associées aux faisceaux étales
$\PPic_{X_\eta/\eta}$ et~$E_\eta$ permettent de conclure.
\end{demo}

\bigskip
Vu le lemme~\ref{ch1dh1h1}, la suite exacte~(\ref{ch1dsedeltam}) induit
une suite exacte
$$
\myxyhook\xymatrix{
H^1(C, \PPic_{X/C}) \ar[r] & \displaystyle\frac{H^1(C,\Erond)}{\langle[\Xrond]\rangle} \ar[r] &
\displaystyle \prod_{M\in \Mrond} H^2(\kappa(M),\Delta_M) \rlap{\text{.}}
}
$$
Posons alors:\glossary{$\Drond(C,X)$}
$$
\Drond(C,X)=\Ker\left(\frac{H^1(C,\Erond)}{\langle[\Xrond]\rangle} \longrightarrow
\prod_{M\in\Mrond}H^2(\kappa(M),\Delta_M)\right)\!\rlap{\text{.}}
$$

\bigskip
\begin{proposition}
Le groupe $\Drond(C,X)$ est égal au noyau de la flèche naturelle
$$
\frac{H^1(C,\Erond)}{\langle[\Xrond]\rangle} \longrightarrow
\prod_{M \in \Mrond}\frac{H^1(\kappa(M),\ff{F}M)}{\langle\delta_M^0([\Xrond])\rangle} \rlap{\text{,}}
$$
où $\delta_M^0 \colon H^1(C,\Erond)\rightarrow H^1(\kappa(M),\ff{F}M)$ est l'application
induite par le morphisme $\Erond \rightarrow i_{M\star} \ff{F}M$.
\end{proposition}

\bigskip
\begin{demo}
Cf.~\cite[Prop.~4.3.1]{css}.  En toute rigueur, il faut d'abord remarquer que l'on
peut supposer la fibration~$\pi$ relativement minimale (ce qui est évident),
puisque cette hypothèse est faite dans~\cite[§4.3]{css}. Le même commentaire s'applique
à plusieurs reprises ci-dessous mais nous ne le répéterons pas.
\end{demo}

\bigskip
Il en résulte immédiatement:

\bigskip
\begin{corollaire}
\label{ch1dtdd}
On a l'inclusion
$$
\Drond(C, X) \cap
\frac{\tors{2}H^1(C,\Erond)}{\langle[\Xrond]\rangle}
\subset
\frac{\TDC{C}}{\langle[\Xrond]\rangle}
$$
de sous-groupes de $H^1(C,\Erond)/\langle[\Xrond]\rangle$.
Si $L_M=\kappa(M)$ pour tout $M \in \Mrond$, cette inclusion est une égalité.
En particulier, on a alors une injection
canonique $\TDC{C}/\langle[\Xrond]\rangle \hookrightarrow \tors{2}\Drond(C,X)$,
qui est un isomorphisme si de plus $[\Xrond]$ n'est pas divisible par~$2$
dans $H^1(C,\Erond)$.
\end{corollaire}

\bigskip
\begin{corollaire}
\label{ch1dcompd}
Supposons que
$C=\P^1_k$ et que $\ff{F}M=\Z/2$ pour tout $M \in \Mrond$.
Sous ces hypothèses, une \condD{} est également définie dans~\cite[p.~583]{css}.
Alors, si~$\pi$ n'admet pas de section,
la \condD{} de~\cite{css} est satisfaite si et seulement si la nôtre l'est
et que le groupe~$E_\eta(K)$ est fini;
si~$\pi$ admet une section, la \condD{} de~\cite{css} est satisfaite si et seulement
si la nôtre l'est et que le groupe~$E_\eta(K)$ est de rang au plus~$1$.
\end{corollaire}

\bigskip
\begin{demo}
Le groupe
$\Drond(\P^1_k,X) \cap \tors{2}H^1(\P^1_k,\Erond)/\langle[\Xrond]\rangle$
est isomorphe au groupe noté
$\Drond^2(X/\P^1_k)$ dans~\cite{css}.
Par ailleurs, l'hypothèse $\ff{F}M=\Z/2$ entraîne que $L_M=\kappa(M)$
pour tout $M \in \Mrond$, d'où l'on déduit, grâce au corollaire~\ref{ch1dtdd}, que
$$\Drond(\P^1_k,X) \cap \tors{2}H^1(\P^1_k,\Erond)/\langle[\Xrond]\rangle=
\TDC{C}/\langle[\Xrond]\rangle\rlap{\text{.}}$$
Notre \condD{} est donc satisfaite si et seulement si le groupe
noté
$\Drond^2(X/\P^1_k)$ dans~\cite{css}
est nul.
Lorsque~$\pi$ n'admet pas de section, le premier énoncé de \cite[§4.7]{css}
permet de conclure. On constate tout de suite que la preuve de cet énoncé
s'adapte sans difficulté au cas où~$\pi$ possède une section,
avec la réserve que l'on doit autoriser le rang de~$E_\eta(K)$ à être
égal à~$1$.
\end{demo}

\bigskip
L'injection $\TDC{C}/\langle[\Xrond]\rangle \hookrightarrow \tors{2}\Drond(C,X)$
donnée par le corollaire~\ref{ch1dtdd}
n'existe malheureusement que sous l'hypothèse que $L_M=\kappa(M)$ pour tout $M \in \Mrond$.
Ainsi, à la différence de ce qui se passait dans la situation de~\cite{css},
le groupe $\Drond(C,X)$ ne suffit plus à exprimer la \condDC{C} dans le cas général
de réduction semi-stable.  C'est pourtant toujours $\Drond(C,X)$ qui est naturellement
lié au groupe de Brauer horizontal de~$X$, comme le montre la proposition suivante.

La suite spectrale de Leray pour la fibre générique de~$\pi$ et le faisceau étale~$\Gm$
fournit
une suite exacte
$$
\myxyhook\xymatrix{
\Br(K) \ar[r] & \Br_1(X_\eta) \ar[r] & H^1(K, \PPic_{X_\eta/\eta}) \ar[r] & H^3(K,\Gm) \rlap{\text{.}}
}
$$
On a $\Br_1(X_\eta)=\Br(X_\eta)$ puisque $X_\eta$ est une courbe lisse sur~$K$
(théorème de Tsen, cf.~\cite[Corollaire~1.2]{grothbr3}) et $H^3(K,\Gm)=0$ d'après l'hypothèse que~$k$ est un corps de nombres et la théorie du corps de classes (cf.~\cite[p.~241]{harduke}), d'où $\Br(X_\eta)/\Br(K)=H^1(K,\PPic_{X_\eta/\eta})$.
Ceci permet de considérer $\Drond(C,X)$ comme un sous-groupe de $\Br(X_\eta)/\Br(K)$,
grâce au lemme~\ref{ch1dh1h1}.

\bigskip
\begin{proposition}
\label{ch1dbrd}
On a l'inclusion $\Br_\hor(X) \subset \Drond(C,X)$ de\index{groupe de Brauer!horizontal} sous-groupes
de~$\Br(X_\eta)/\Br(K)$.  C'est une égalité si $C=\A^1_k$ ou si~$\pi$ admet une section.
\end{proposition}

\bigskip
\begin{demo}
Cf.~\cite[Th.~4.5.2]{css}.  Pour voir que l'inclusion est une égalité lorsque~$\pi$ admet une section,
cf.~les commentaires à la fin de \cite[§4.6]{css}.
\end{demo}

\bigskip
\begin{corollaire}
\label{ch1dbrdvert}
Si la \condDC{C} est satisfaite et si la classe $[\Xrond]$ n'est pas divisible par~$2$
dans $H^1(C,\Erond)$, le sous-groupe de torsion $2$\nobreakdash-primaire de $\Br(X)/\Br(k)$ est inclus
dans\index{groupe de Brauer!vertical} $\Br_\vert(X)/\Br(k)$.
\end{corollaire}

\bigskip
\begin{demo}
D'après le corollaire~\ref{ch1dtdd} et l'hypothèse sur~$[\Xrond]$,
le groupe $\tors{2}\Drond(C,X)$ est nul. La proposition~\ref{ch1dbrd} permet
alors de conclure.
\end{demo}

\bigskip
\begin{proposition}
\label{ch1dbrd2}
Si $C=\P^1_k$ et que $X_\infty$ est lisse,
$\Drond(C,X)$ est égal au sous-groupe $\Br^{\gnr/X}_\hor(X_{\A^1_k})$ de
$\Br(X_\eta)/\Br(K)$ engendré par les classes de~$\Br(X_{\A^1_k})$
géométriquement non ramifiées sur~$X$.
\end{proposition}

\bigskip
\begin{demo}
Cela résulte de \cite[Th.~4.5.1]{css} et de la surjectivité de la flèche
$$\Br(K) \longrightarrow \bigoplus_{M \in {\A^1_k}^{(1)}} H^1(\kappa(M),\Q/\Z)$$
induite par les résidus (cf.~\cite[Prop.~1.2.1]{cscrelle94}).
\end{demo}

\bigskip
On déduit enfin des propositions~\ref{ch1dbrd} et~\ref{ch1dbrd2}
et du corollaire~\ref{ch1dtdd}, compte tenu que le groupe
$\Br(X_{\A^1_k})$ est de torsion:

\bigskip
\begin{corollaire}
\label{ch1conddtrad}
Supposons que $[\Xrond]$ ne soit pas divisible par~$2$,
que $C=\P^1_k$, que $X_\infty$ soit lisse
et que $L_M=\kappa(M)$
pour tout $M \in \Mrond$.
Alors les conditions~\cDz{} et~\cD{} équivalent respectivement
aux inclusions
$\Br(X_{\A^1_k})\{2\} \subset \Br_\vert(X_{\A^1_k})$
et $\Br^{\gnr/X}(X_{\A^1_k})\{2\} \subset \Br_\vert(X_{\A^1_k})$.
\end{corollaire}

\bigskip
Il est possible de traduire les conditions~\cD{} et~\cDz{} en termes de groupe de Brauer
sans supposer que $L_M=\kappa(M)$ pour tout $M \in \Mrond$; cependant,
la traduction est alors nettement moins agréable et c'est pourquoi nous avons
choisi de\index{condition (D)@\condD{}|)} l'omettre.

\section{Applications à l'existence de points rationnels}
\label{ch1parapprat}

\subsection{Obstruction de Brauer-Manin}

Les conséquences du théorème~\ref{ch1thprin}
quant à l'existence de points rationnels sur les pinceaux semi-stables de courbes de genre~$1$
dont les jacobiennes ont leur $2$\nobreakdash-torsion rationnelle sont les suivantes.
Reprenons les notations du paragraphe~\ref{ch1hypnot}.

On suppose que~$k$ est un corps de nombres, que~$C$
est isomorphe à~$\P^1_k$ et que $\Mrond \neq \emptyset$.

\bigskip
\begin{theoreme}
\label{ch1appobm1}
Admettons l'hypothèse de Schinzel.
Alors l'adhérence de~$\RD$ dans~$C(\A_k)$ est égale à $\pi(X(\A_k)^{\Br_\vert})$.
\end{theoreme}

\bigskip
\begin{demo}
L'inclusion $\RD \subset \pi(X(\A_k)^{\Br_\vert})$ est triviale; de fait, tout
point adélique de~$X$ contenu dans une fibre de~$\pi$ au-dessus d'un point rationnel
de~$U$ est orthogonal à $\Br_\vert(X)$.
Soit $(\fp{P}v)_{v \in \Omega} \in X(\A_k)^{\Br_\vert}$.  Étant donnés un ensemble
$S \subset \Omega$ fini et pour chaque
$v \in S$, un voisinage $\Arond_v$ de $\pi(\fp{P}v)$ dans $C(k_v)$, on veut
montrer qu'il existe un élément de~$\RD$ appartenant à $\Arond_v$ pour tout $v
\in S$.

Le théorème des fonctions implicites et la propriété d'approximation faible sur~$\P^1_k$
entraînent l'existence d'un point $\infty \in C(k)$ au-dessus duquel
la fibre de~$\pi$ est lisse, tel que
$X_\infty(k_v)\neq\emptyset$ et $\infty \in \Arond_v$
pour toute place $v \in \Omega$ archimédienne,
et tel que pour toute place $v \in \Omega$ réelle, les points $\pi(\fp{P}v)$ et~$\infty$ appartiennent
à la même composante connexe de~$O(k_v)$, en notant~$O$ le plus grand ouvert
de~$C$ au-dessus duquel~$\pi$ est lisse.
Fixons un isomorphisme $C \isoto \P^1_k$ envoyant $\infty \in C(k)$ sur~$\infty \in \P^1(k)$
et~$\pi(\fp{P}v)$ dans la composante connexe non majorée de l'image de $U(k_v)$
pour~$v$ réelle, où $U = O \setminus \{\infty\}$.
On peut maintenant appliquer le théorème~\ref{ch1thprin}, ce qui
produit deux ensembles finis $S_0 \subset \Omega$ et $B_0 \subset \Br(U)$.
Quitte à agrandir $S$ et à poser
$\Arond_v = C(k_v)$ pour les places ainsi introduites, on peut supposer que~$S$ contient~$S_0$
et l'ensemble des places archimédiennes de~$k$.
D'après le théorème des fonctions implicites, la continuité de l'évaluation des classes
de~$\Br(X)$ sur $X(k_v)$ et la finitude de~$B_0$, il
existe $(\fp{P'}v)_{v \in \Omega} \in X(\A_k)^{B_0 \cap \Br(X)}$
arbitrairement proche de $(\fp{P}v)_{v \in \Omega}$ et tel que $\pi(\fp{P'}v) \in U(k_v)$
pour tout $v \in \Omega$. On peut notamment
supposer que $\pi(\fp{P'}v) \in \Arond_v \cap U(k_v)$ pour tout $v \in S$
et que~$\pi(\fp{P'}v)$ et~$\pi(\fp{P}v)$ appartiennent à la même composante connexe de~$U(k_v)$
pour $v \in \Omega$ réelle.
D'après le\index{lemme formel} «~lemme formel~»,
il existe $S_1 \subset \Omega$ fini contenant $S$
et $(Q_v)_{v \in S_1} \in \prod_{v \in S_1} X_U(k_v)$ tels que $Q_v=\fp{P'}v$ pour $v \in S$ et
que
$$\sum_{v \in S_1} \inv_v (\pi^\star A)(Q_v)=0$$
pour tout $A \in B_0$.
Posons $x_v=\pi(Q_v)$ pour $v \in S_1$.
Comme $(\pi^\star A)(Q_v)=A(x_v)$, le théorème~\ref{ch1thprin} permet maintenant de conclure.
\end{demo}

\bigskip
\begin{theoreme}
\label{ch1appobm2}
Admettons l'hypothèse de Schinzel et la finitude
des groupes de Tate-Shafarevich des courbes elliptiques~$\Erond_x$ pour $x \in U(k)$.
Supposons que la \condD{} soit vérifiée et que $X(\A_k)^{\Br_\vert}\neq \emptyset$.
Alors $X(k)\neq \emptyset$.
Si de plus~$\pi$ ne possède pas de section, l'ensemble~$X(k)$ est
Zariski-dense dans~$X$.
\end{theoreme}

\bigskip
\begin{demo}
Supposons que $X(\A_k)^{\Br_\vert}\neq\emptyset$.
Le théorème~\ref{ch1appobm1} permet d'en déduire que $\RD \neq \emptyset$.
Soit $x \in \RD$.
Puisque la \condD{} est satisfaite, le conoyau de la flèche naturelle $\Erond(U)/2 \rightarrow \SD$
est d'ordre~$\leq 2$ (cf.~paragraphe~\ref{ch1hypnot}).  L'appartenance de~$x$ à~$\RD$ et la commutativité du carré
$$
\myxyhook\xymatrix{
\Erond(U)/2 \ar[r] \ar[d] & \SD \ar[d] \\
\Erond_x(k)/2 \ar[r] & H^1(k, \tors{2}\Erond_x)
}
$$
permettent d'en déduire
que $\dim_{\ff{\F}2} (\tors{2}\Sha(k, \Erond_x)) \leq 1$.

L'accouplement\index{accouplement de Cassels-Tate}
de Cassels-Tate $\Sha(k,\Erond_x) \times \Sha(k,\Erond_x) \rightarrow \Q/\Z$
est non dégénéré d'après la finitude de $\Sha(k,\Erond_x)$.  Par ailleurs, il est alterné.
Comme tout groupe abélien fini\index{groupe de Tate-Shafarevich!finitude}
muni d'une forme bilinéaire alternée non dégénérée à valeurs dans~$\Q/\Z$,
le groupe $\Sha(k,\Erond_x)$ est isomorphe à $T \times T$ pour un certain
groupe~$T$.  L'entier $\dim_{\ff{\F}2} (\tors{2}\Sha(k, \Erond_x))$ est
donc pair; étant positif et $\leq 1$, il est nécessairement nul,
d'où $\tors{2}\Sha(k,\Erond_x)=0$. On a en particulier~$X_x(k) \neq \emptyset$
puisque $x \in \RA$, et donc $X(k)\neq\emptyset$.

Supposons de plus que~$\pi$ ne possède pas de section.
Notons $G \subset \SD$ le sous-groupe image réciproque
de $\{0,[\Xrond]\}$ par la flèche de droite de la suite exacte~(\ref{ch1sesgtsg}).
C'est un groupe d'ordre au moins~$8$ puisque la $2$\nobreakdash-torsion de~$E_\eta$ est
rationnelle et que~$\pi$ ne possède pas de section.
Par définition de~$\RD$, la flèche $G \rightarrow \Sel_2(k,\Erond_x)$ d'évaluation
en~$x$ est injective.
Le groupe $\Sel_2(k,\Erond_x)$ est donc lui aussi d'ordre au moins~$8$,
ce qui entraîne que la courbe elliptique $\Erond_x$ est de rang non nul, puisque
sa $2$\nobreakdash-torsion est rationnelle et que $\tors{2}\Sha(k,\Erond_x)=0$.
L'ensemble $X_x(k)$, en bijection avec $\Erond_x(k)$, est donc infini.
L'ensemble~$\RD$ étant non seulement non vide mais même infini d'après le théorème~\ref{ch1thprin},
on en déduit que les points rationnels de~$X$ sont Zariski-denses, en faisant varier $x \in \RD$.
\end{demo}

\bigskip
Lorsque $\Mrond = \emptyset$ et que $C=\P^1_k$, on peut bien sûr obtenir des
énoncés similaires, avec~$\RDz$ et la \condDz{} à la place de~$\RD$
et de la \condD{}.

\bigskip
Le théorème principal de~\cite{css} n'est autre que le
cas particulier du théorème~\ref{ch1appobm2}
où les hypothèses supplémentaires suivantes sont supposées satisfaites:
\begin{itemize}
\item pour tout $M \in \Mrond$, le groupe $\ff{F}M$ est d'ordre~$2$;
\item soit le rang de Mordell-Weil de la courbe elliptique $E_\eta/K$ est nul
et~$\pi$ ne possède pas de section, soit le rang est exactement~$1$ et~$\pi$
possède une section.
\end{itemize}
(Cf.~corollaire~\ref{ch1dcompd}.)

\subsection{Surfaces quartiques diagonales}
\label{ch1parquartiques}

Nous vérifions dans ce paragraphe que le théorème~\ref{ch1appobm2} permet de retrouver
les résultats de Swinnerton-Dyer \cite[§3]{sddiagquartic} sur l'arithmétique
des surfaces quartiques
diagonales définies sur~$\Q$ et nous en profitons pour les généraliser à tout corps de nombres.

\bigskip
\begin{theoreme}
\label{ch1thquartiques}
Soit~$k$ un corps de nombres. Admettons l'hypothèse de Schinzel et la finitude
des groupes de Tate-Shafarevich des courbes elliptiques sur~$k$.
Soit $X \subset \P^3_k$ la surface projective et lisse d'équation
$$
a_0 x_0^4 + a_1 x_1^4 + a_2 x_2^4 + a_3 x_3^4 = 0 \rlap{\text{,}}
$$
où $\uplet{a_0}{a_3} \in k^\star$.
Supposons que $X(\A_k)^{\Br_1} \neq \emptyset$.
Si $a_0a_1a_2a_3$ est un carré dans~$k$ mais pas une puissance
quatrième et que pour tous $i,j\in\{\uplet{0}{3}\}$ distincts,
ni $a_i a_j$ ni $-a_ia_j$ n'est un carré,
alors $X(k) \neq \emptyset$.
\end{theoreme}

\bigskip
Rappelons que $\Br_1(X)$ désigne le sous-groupe
des classes \emph{algébriques} de $\Br(X)$, c'est-à-dire tuées par extension
des scalaires de~$k$ à une clôture algébrique de~$k$.

\bigskip
\begin{remarque}
L'énoncé du théorème~\ref{ch1thquartiques} soulève naturellement la question
suivante: sous les hypothèses de ce théorème, a-t-on\index{groupe de Brauer!algébrique}
nécessairement $\Br_1(X)=\Br(X)$ ?
Le corollaire~\ref{ch1dbrdvert} permet, dans certains cas très particuliers, de prouver
l'égalité entre les sous-groupes de torsion $2$\nobreakdash-primaire de $\Br_1(X)/\Br(k)$ et de
$\Br(X)/\Br(k)$.
Il s'avère que le groupe $\Br_1(X)/\Br(k)$ est toujours fini d'ordre une puissance de~$2$
(cf.~\cite{bright}).  En revanche, sans condition sur les coefficients~$a_i$,
le groupe $\Br(X)/\Br(k)$ peut comporter des éléments de tout ordre; en effet,
la surface~$X$ est une surface~$K3$ de nombre de Picard géométrique~$20$, de sorte qu'il
existe un isomorphisme
$$
\Br(X \otimes_k \ksep) \simeq \bigoplus_\ell \left(\Q_\ell/\Z_\ell\right)^2\rlap{\text{,}}
$$
où~$\ell$ parcourt l'ensemble des nombres premiers et~$\ksep$ désigne une clôture algébrique de~$k$ (cf.~\cite[Corollaire~3.4]{grothbr2}).
\end{remarque}

\bigskip
\begin{demo}[ du théorème~\ref{ch1thquartiques}]%
La démonstration qui suit est essentiellement une transcription des
calculs de \cite[§2--§3]{sddiagquartic}.  Supposons que \mbox{$X(\A_k)\neq \emptyset$}.
D'après le théorème de Hasse-Minkowski,
la surface quadrique d'équation
$a_0 x_0^2 + a_1 x_1^2 + a_2 x_2^2 + a_3 x_3^2 = 0$
possède alors un point $k$\nobreakdash-rationnel.
Comme $a_0a_1a_2a_3 \in k^{\star 2}$, elle contient même une droite $k$\nobreakdash-rationnelle.
Par conséquent, sa trace sur tout hyperplan $k$\nobreakdash-rationnel admet un point $k$\nobreakdash-rationnel.
Il existe en particulier $r_1,r_2,r_3\in k$ non tous nuls tels que
$a_1 r_1^2 + a_2 r_2^2 + a_3 r_3^2 = 0$.
Remarquons qu'en vertu des hypothèses du théorème, aucun des~$r_i$ ne peut être nul.
Soit alors $\theta \in k^\star$ une racine
carrée de $a_0a_1a_2a_3$, et posons
\begin{align*}
A &= \theta r_2 x_0^2 + a_1 a_3 (r_3 x_1^2 - r_1 x_3^2) \rlap{\text{,}} &
B &= \theta r_3 x_0^2 - a_1 a_2 (r_2 x_1^2 + r_1 x_2^2) \rlap{\text{,}} \\
C &= a_3 \theta r_3 x_0^2 - a_1a_2a_3(r_2x_1^2-r_1x_2^2) \rlap{\text{,}} &
D &= -a_2\theta r_2 x_0^2 - a_1a_2a_3(r_3x_1^2+r_1x_3^2) \rlap{\text{,}}
\end{align*}
de sorte que
$a_0 x_0^4 + a_1 x_1^4 + a_2 x_2^4 + a_3 x_3^4 = (AD-BC)/(a_1^2 r_1^2 a_2 a_3)$.
Notons $\pi \colon X \rightarrow \P^1_k$ le morphisme
$[x_0:x_1:x_2:x_3]\mapsto [A:B]$.  (C'est évidemment une application rationnelle;
on peut vérifier à la main qu'elle est définie partout ou alors invoquer le lemme selon
lequel toute application rationnelle d'une surface~$K3$ vers~$\P^1_k$ est un morphisme.)
La fibre générique de~$\pi$ est la courbe de genre~$1$
définie par le système d'équations $B=tA$, $D=tC$ sur $\kappa(\P^1_k)=k(t)$.
Éliminant tour à tour $x_2^2$ et $x_3^2$ dans ces équations, on obtient le système équivalent
\begin{equation}
\label{ch1quarteq}
\left\{
\begin{aligned}
m_1 x_1^2 - m_3 x_3^2 &= (e_3 - e_1)x_0^2\rlap{\text{,}} \\
m_1 x_1^2 - m_2 x_2^2 &= (e_2 - e_1)x_0^2\rlap{\text{,}}
\end{aligned}
\right.
\end{equation}
où
$e_1=0$, $e_2=a_1a_2a_3^2r_1^2\theta f_3^2$, $e_3=-a_1a_3a_2^2r_1^2\theta f_2^2$,
$m_1=a_1^2a_2^2a_3^2r_1^2f_2f_3$, \mbox{$m_2=a_1^2a_2^2a_3^2r_1^3f_1f_3$},
$m_3=a_1^2a_2^2a_3^2r_1^3f_1f_2$,
$f_1=a_3t^2+a_2$,
$f_2=a_3r_2t^2-2a_3r_3t-a_2r_2$ et
$f_3=a_3r_3t^2+2a_2r_2t-a_2r_3$.
Il sera utile de remarquer que
$e_2-e_3=-a_2a_3a_1^2r_1^4\theta f_1^2$.

La jacobienne de~(\ref{ch1quarteq}) a pour équation de Weierstrass l'équation~(\ref{ch1expleqw}),
où les $e_i \in k[t]$ sont les polynômes que l'on vient de définir.
On vérifie tout de suite que les polynômes $f_1$, $f_2$ et~$f_3$ sont irréductibles
et deux à deux premiers entre eux. L'équation de Weierstrass~(\ref{ch1expleqw}) est donc minimale,
ce qui permet d'appliquer les résultats du paragraphe~\ref{ch1parexplicited}, dont on reprend
les notations et conventions.  Pour $i \in \{1,2,3\}$, soit $M_i \in \P^1_k$ le point
fermé en lequel s'annule~$f_i$.
La courbe elliptique~$E_\eta$ a réduction de type~$I_4$ en $M_1$, $M_2$ et $M_3$ et
bonne réduction partout ailleurs (cf.~commentaires au début du paragraphe~\ref{ch1parexplicited}).
D'autre part, la courbe~(\ref{ch1quarteq}) est un $2$\nobreakdash-revêtement de $E_\eta$;
sa classe dans $H^1(K,\tors{2}E_\eta)=\Hgoth(K)$ est la classe du couple $(m_1,m_2)$.
Notons-la~$\mgoth$.
La proposition~\ref{ch1explsg} montre que $\mgoth \in \SG_2(\P^1_k,\Erond)$,
de sorte que les fibres de~$\pi$ sont réduites.
Toutes les hypothèses du paragraphe~\ref{ch1hypnot} sont donc satisfaites.
Pour conclure, il reste seulement à vérifier que la \condD{} l'est aussi,
compte tenu du théorème~\ref{ch1appobm2} et de l'inclusion $\Br_\vert(X) \subset \Br_1(X)$,
qui est une conséquence du théorème de Tsen.

Soit $\mgoth' \in \SD$, représenté par un triplet $(m'_1,m'_2,m'_3)$.
Nous allons montrer que $\mgoth'$ appartient au sous-groupe engendré
par~$\mgoth$ et par les images des points d'ordre~$2$ de $E_\eta(K)$.
Vu la proposition~\ref{ch1explsg}, il existe $\alpha_1,\alpha_2,\alpha_3\in k^\star$
et $\epsilon_1,\epsilon_2,\epsilon_3 \in \{0,1\}$ tels que
$m'_i=\alpha_i f_j^{\epsilon_{\smash[b]{j}\vphantom{k}}}f_k^{\epsilon_k}$ pour toute permutation
$(i,j,k)$ de $(1,2,3)$.  Par ailleurs, à l'aide de la proposition~\ref{ch1expldlm},
on constate que pour toute permutation cyclique $(i,j,k)$ de $(1,2,3)$,
on a $\kappa(M_i)=k\!\left(\!\sqrt{-a_j a_k}\right)$ et $\gamma_{M_i}$ est égal à la classe
de $-a_i a_k \theta$ dans $\kappa(M_i)^\star/\kappa(M_i)^{\star 2}$.
Comme $\mgoth' \in \SD$, le corollaire~\ref{ch1expldeltai4} montre
que pour tout $i \in \{1,2,3\}$, l'image de $m'_i(M_i) \in \kappa(M_i)^\star$
dans $\kappa(M_i)^\star/\kappa(M_i)^{\star 2}$ appartient au sous-groupe
engendré par $\gamma_{M_i}$ et $m_i(M_i)$.
Il existe donc $\rho_1,\rho_2,\rho_3,\sigma_1,\sigma_2,\sigma_3\in\{0,1\}$
tels que $m'_i(M_i)(m_i(M_i))_{\vphantom{M_i}}^{\smash{\rho}_i}\gamma_{M_i}^{\sigma_i} \in \kappa(M_i)^{\star 2}$ pour tout
$i \in\{1,2,3\}$.
Passant aux normes de $\kappa(M_i)$ à~$k$ et remarquant que~$\gamma_{M_i}$
provient de $k^\star/k^{\star 2}$,
on en déduit que
\begin{equation}
\label{ch1quartrel}
N_{\kappa(M_i)/k}(f_j(M_i))^{\epsilon_j+\rho_i}
\sim N_{\kappa(M_i)/k}(f_k(M_i))^{\epsilon_k+\rho_i}
\end{equation}
pour toute permutation $(i,j,k)$ de $(1,2,3)$,
où le symbole $\sim$ désigne l'égalité dans $k^\star/k^{\star 2}$.
Un calcul facile montre que $N_{\kappa(M_i)/k}(f_j(M_i)) \sim -a_i a_j$
pour tous $i,j\in\{1,2,3\}$ distincts, de sorte que les relations~(\ref{ch1quartrel})
entraînent que $\epsilon_1=\epsilon_2=\epsilon_3=\rho_1=\rho_2=\rho_3$ puisqu'aucun des $-a_i a_j$
n'est un carré.
Quitte à remplacer $\mgoth'$ par $\mgoth\mgoth'$,
on peut donc supposer que $\epsilon_1=\epsilon_2=\epsilon_3=\rho_1=\rho_2=\rho_3=0$.
Que $m_i(M_i)\gamma_{M_i}^{\sigma_i} \in \kappa(M_i)^{\star 2}$
pour tout~$i$ signifie maintenant ceci:
pour toute permutation cyclique $(i,j,k)$ de $(1,2,3)$,
on a $\alpha_i \sim 1$ ou $\alpha_i \sim -a_ja_k$ ou $\alpha_i \sim a_i a_j \theta$
ou encore $\alpha_i \sim -a_i a_k \theta$.  Ces quatre valeurs possibles pour
la classe de $\alpha_i$ dans $k^\star/k^{\star 2}$ sont précisément les valeurs
obtenues lorsque~$\mgoth'$ est l'image de l'un des quatre points d'ordre~$2$ de $E_\eta(K)$;
on peut donc supposer que $\alpha_1 \sim 1$, quitte à translater~$\mgoth'$ par un
point d'ordre~$2$ bien choisi.  Comme $m'_1m'_2m'_3$ est un carré,
on a alors $\alpha_2 \sim \alpha_3$.
Compte tenu d'une part des valeurs possibles pour les classes de~$\alpha_2$ et de~$\alpha_3$
dans $k^\star/k^{\star 2}$ et d'autre part des hypothèses du théorème, on voit à présent
que nécessairement $\alpha_2\sim\alpha_3\sim 1$, d'où le résultat.
\end{demo}

\section{Secondes descentes et approximation faible}
\label{ch1parsecdesc}

Notons\index{seconde descente|(} $\RO$ l'ensemble des $x \in \RA$ tels que la classe de~$X_x$
dans $\Sha(k,\Erond_x)$ soit orthogonale à $\tors{2}\Sha(k,\Erond_x)$ pour
l'accouplement\index{accouplement de Cassels-Tate|(} de Cassels-Tate, et $\overline{\RO}$ son adhérence dans $C(\A_k)$
pour la topologie adélique.

On suppose que~$C$ est isomorphe à~$\P^1_k$ et que $\Mrond \neq \emptyset$.

\bigskip
\begin{theoreme}
\label{ch1thsecdesc}
Admettons l'hypothèse de Schinzel.  Supposons de plus que $L_M=\kappa(M)$
pour tout $M \in \Mrond$.  On a alors
$\pi(X(\A_k)^\Br) \subset \overline{\RO}$.
En particulier, si $X(\A_k)^\Br=X(\A_k)$,
l'ensemble~$\RO$ est dense dans~$\RA$.
\end{theoreme}

\bigskip
\begin{demo}
Soit $(\fp{P}v)_{v\in\Omega} \in X(\A_k)^\Br$.
Étant donnés un ensemble $S \subset \Omega$ fini et pour chaque $v \in S$,
un voisinage $\Arond_v$ de~$\pi(\fp{P}v)$ dans~$C(k_v)$, on veut montrer qu'il existe
un élément de~$\RO$ appartenant à~$\Arond_v$ pour tout $v \in S$.

Le théorème des fonctions implicites et la propriété d'approximation faible sur~$\P^1_k$
entraînent l'existence d'un point $\infty \in C(k)$ au-dessus duquel
la fibre de~$\pi$ est lisse, tel que
$X_\infty(k_v)\neq\emptyset$ et $\infty \in \Arond_v$
pour toute place $v \in \Omega$ archimédienne,
et tel que pour toute place $v \in \Omega$ réelle, les points $\pi(\fp{P}v)$ et~$\infty$ appartiennent
à la même composante connexe de~$O(k_v)$, en notant~$O$ le plus grand ouvert
de~$C$ au-dessus duquel~$\pi$ est lisse.
Fixons un isomorphisme $C \isoto \P^1_k$ envoyant $\infty \in C(k)$ sur~$\infty \in \P^1(k)$
et~$\pi(\fp{P}v)$ dans la composante connexe non majorée de l'image de $U(k_v)$
pour~$v$ réelle, où $U = O \setminus \{\infty\}$.
On peut maintenant appliquer le théorème~\ref{ch1thprin}, ce qui
produit deux ensembles finis $S_0 \subset \Omega$ et $B_0 \subset \Br(U)$.

Soit $s \colon \Br_\hor(X_{\A^1_k})\rightarrow \Br(X_{\A^1_k})$ une section
ensembliste de la projection naturelle.
Le corollaire~\ref{ch1dtdd} et la proposition~\ref{ch1dbrd} fournissent une
application canonique $i \colon \TDz \rightarrow \Br_\hor(X_{\A^1_k})$.
Notons $B_1 \subset \Br(X_{\A^1_k})$ le
sous-groupe engendré par $s(i(\TDz))$ et $B = B_1 + \pi^\star B_0 \subset \Br(X_U)$.
Le groupe~$B$ est fini en vertu du corollaire~\ref{ch1finitude}.
Quitte à agrandir $S$ et à poser
$\Arond_v = C(k_v)$ pour les places ainsi introduites, on peut supposer:
\begin{itemize}
\item que~$S$ contient $S_0$ et l'ensemble des places archimédiennes de~$k$;
\item que $\forall A \in B \cap \Br(X)$, $\forall v \in \Omega \setminus S$,
$\forall p \in X(k_v)$, $A(p)=0$
(en effet, comme~$B$ est fini, quitte à choisir~$S$ assez grand, on peut étendre~$X$
en un $\Orond_S$\nobreakdash-schéma propre sur lequel les classes de $B \cap \Br(X)$ sont
non ramifiées; tout $k_v$\nobreakdash-point de~$X$ s'étend alors en un $\Orond_v$\nobreakdash-point
de ce schéma, et la nullité de $\Br(\Orond_v)$ permet de conclure);
\item qu'il existe un morphisme propre et plat $\Xtilde \rightarrow \A^1_{\Orond_S}$
dont la restriction au-dessus de~$\A^1_k$ est égale à $X_{\A^1_k} \rightarrow \A^1_k$,
avec $\Xtilde$ régulier;
\item que $B_1 \subset \Br(\Xtilde)$.
\end{itemize}
D'après le théorème des fonctions implicites, la continuité de l'évaluation des classes
de~$\Br(X)$ sur $X(k_v)$ et la finitude de~$B$, il
existe $(\fp{P'}v)_{v \in \Omega} \in X(\A_k)^{B \cap \Br(X)}$
arbitrairement proche de $(\fp{P}v)_{v \in \Omega}$ et tel que $\pi(\fp{P'}v) \in U(k_v)$
pour tout $v \in \Omega$. On peut notamment
supposer que $\pi(\fp{P'}v) \in \Arond_v \cap U(k_v)$ pour tout $v \in S$
et que~$\pi(\fp{P'}v)$ et~$\pi(\fp{P}v)$ appartiennent à la même composante connexe de~$U(k_v)$
pour $v \in \Omega$ réelle.
D'après le\index{lemme formel} «~lemme formel~»,
il existe $S_1 \subset \Omega$ fini contenant~$S$
et $(Q_v)_{v \in S_1} \in \prod_{v \in S_1} X_U(k_v)$ tels que $Q_v=\fp{P'}v$ pour $v \in S$ et
que
\begin{equation}
\label{ch1sdeqs1}
\sum_{v \in S_1} \inv_v A(Q_v)=0
\end{equation}
pour tout $A \in B$.
Posons $x_v=\pi(Q_v)$ pour $v \in S_1$.
On a en particulier
$$\sum_{v \in S_1} \inv_v A(x_v)=0$$ pour tout $A \in B_0$,
d'où l'existence, grâce au théorème~\ref{ch1thprin}, de $x \in \RDz$ entier hors de~$S_1$
et arbitrairement proche de~$x_v$ pour tout $v \in S_1$.  On peut notamment supposer
que $x \in \Arond_v$ pour tout $v \in S$ et qu'il existe une famille $(Q'_v)_{v \in S_1} \in
\prod_{v\in S_1}X_x(k_v)$ arbitrairement proche
de $(Q_v)_{v \in S_1}$ (théorème des fonctions implicites) et en particulier telle que pour
tout $A \in B$ et tout $v \in S_1$,
on ait $\inv_v A(Q'_v)=\inv_v A(Q_v)$.  Choisissons arbitrairement
$Q'_v \in X_x(k_v)$ pour $v \in \Omega \setminus S_1$ (de tels~$Q'_v$
existent puisque $x \in \RA$).  On dispose à présent d'un point adélique $(Q'_v)_{v \in \Omega} \in X_x(\A_k)$.

Nous allons maintenant montrer que $x \in \RO$, ce qui conclura la preuve du théorème.
Fixons $\alpha \in \tors{2}\Sha(k,\Erond_x)$.
Si~$Z$ est une variété sur~$k$, notons
$$
\Ba(k,Z)=\Ker\left(\Br(Z) \longrightarrow \prod_{v \in \Omega} \Br(Z\otimes_k k_v)/\Br(k_v)
\right)\!\rlap{\text{.}}
$$
La suite spectrale de Leray (pour le morphisme structural $X_x \rightarrow \Spec(k)$
et le faisceau étale~$\Gm$) et l'égalité
$H^1(k,\PPic_{X_x/k})=H^1(k,\Erond_x)/\langle[X_x]\rangle$
fournissent un morphisme canonique
$r \colon \Br(X_x) \rightarrow H^1(k,\Erond_x)/\langle[X_x]\rangle$.
Notant $\langle \cdot, \cdot \rangle$ l'accouplement de Cassels-Tate sur $\Sha(k,\Erond_x)$,
on~a
$$
\langle \alpha, [X_x] \rangle = \sum_{v \in \Omega} \inv_v a(Q'_v)
$$
pour tout $a \in \Ba(k,X_x)$ tel que~$r(a)$ soit égal à la classe de~$\alpha$,
d'après un théorème de Manin (cf.~\cite[Theorem~6.2.3]{skotorsors}).

Comme $x \in \RDz$, il existe $b \in \TDz$ dont l'image par
la flèche $H^1(\A^1_k, \Erond) \rightarrow H^1(k, \Erond_x)$
d'évaluation en~$x$ soit égale à~$\alpha$.
Pour $A \in \Br(X_{\A^1_k})$, notons $A|_{X_x}$ l'image de~$A$ par la flèche
$\Br(X_{\A^1_k})\rightarrow \Br(X_x)$ déduite de l'inclusion
\mbox{$X_x \subset X_{\A^1_k}$}.

\bigskip
\begin{lemme}
On a $(s \circ i)(b)|_{X_x} \in \Ba(k,X_x)$
et la classe de~$\alpha$ dans $H^1(k,\Erond_x)/\langle[X_x]\rangle$
est égale à $r((s\circ i)(b)|_{X_x})$.
\end{lemme}

\bigskip
\begin{demo}
Le noyau de~$r$ est constitué des
classes constantes; de même pour l'application analogue à~$r$ au niveau
de chaque complété de~$k$.
La première assertion est donc une conséquence de la seconde et de l'appartenance de~$\alpha$
à $\Sha(k,\Erond_x)$.

Notons $j_U \colon U \rightarrow \A^1_k$, $j \colon \eta \rightarrow \A^1_k$
et $i_x\colon\Spec(k) \rightarrow \A^1_k$ les morphismes canoniques associés
à~$U$, $\eta$ et~$x$
et posons:
$$
P=\PPic_{X_{\A^1_k}/\A^1_k}
\quad;\quad \fp{P}U=\PPic_{X_U/U}
\quad;\quad \fp{P}\eta=\PPic_{X_\eta/\eta}
\quad;\quad \fp{P}x=\PPic_{X_x/k} \rlap{\text{,}}
$$
de sorte que $\fp{P}\eta=j^\star P$.
Le lemme~\ref{ch1lemmedeltam}, appliqué à $C=U$, montre
que $j_{U\star}\fp{P}U=j_{\star} \fp{P}\eta$,
d'où l'existence d'un morphisme $j_\star \fp{P}\eta \rightarrow i_{x\star}\fp{P}x$ rendant
le diagramme de faisceaux étales suivant commutatif:
$$
\myxyhook\xymatrix{
P \ar[dr] \ar[r] & j_\star \fp{P}\eta \ar[d] & \rlap{$\Erond=j_\star E_\eta$}\phantom{\Erond}
\ar[l] \ar[d] \\
& i_{x\star}\fp{P}x & i_{x\star} \Erond_x \ar[l]
}
$$
Passant en cohomologie, on en déduit le diagramme commutatif:
\begin{equation}
\label{ch1sddiag1}
\begin{aligned}
\myxyhook\xymatrix{
H^1(\A^1_k,P) \ar[dr]\ar[r] &
H^1(\A^1_k,j_\star \fp{P}\eta) \ar[d] & H^1(\A^1_k,\Erond)/\langle[\Xrond]\rangle \ar[l] \ar[d] \\
& H^1(k,\fp{P}x) & H^1(k,\Erond_x)\langle[X_x]\rangle \ar[l]
}
\end{aligned}
\end{equation}
Les flèches horizontales de droite sont des isomorphismes (cf.~lemme~\ref{ch1dh1h1}).

La régularité de~$X$ entraîne que $R^2 \pi_\star \Gm=0$, d'après un théorème d'Artin
(cf.~\cite[Cor.~3.2]{grothbr3}).  La suite spectrale de Leray pour le morphisme
$X_{\A^1_k} \rightarrow \A^1_k$ et le faisceau~$\Gm$ fournit donc une flèche
$\Br(X_{\A^1_k})\rightarrow H^1(\A^1_k,P)$ rendant le diagramme
\begin{equation}
\label{ch1sddiag2}
\begin{aligned}
\myxyhook\xymatrix{
\Br(X_\eta) \ar[r] & H^1(K,\fp{P}\eta) \\
\Br(X_{\A^1_k}) \ar[d] \ar[u] \ar[r] & H^1(\A^1_k,P) \ar[d]\ar[u] \\
\Br(X_x) \ar[r] & H^1(k, \fp{P}x)
}
\end{aligned}
\end{equation}
commutatif (les flèches verticales sont en effet induites par des morphismes entre
les suites spectrales qui donnent naissance aux flèches horizontales).
Combinant~(\ref{ch1sddiag1}) et~(\ref{ch1sddiag2}), on obtient le
diagramme commutatif:
$$
\myxyhook\xymatrix{
\Br(X_\eta)/\Br(K) \ar[rr]^\sim && H^1(K,\fp{P}\eta)
& \TDz \ar@{^{ (}->}[d] \\
\Br(X_{\A^1_k}) \ar[u] \ar[d] \ar[r] & H^1(\A^1_k, P) \ar[rd] \ar[ru] \ar[r] &
H^1(\A^1_k, j_\star \fp{P}\eta) \ar[d] \ar@{^{ (}->}[u] &
H^1(\A^1_k,\Erond)/\langle[\Xrond]\rangle \ar[l]_(0.54)\sim \ar[d] \\
\Br(X_x) \ar[rr]^r && H^1(k,\fp{P}x) & H^1(k,\Erond_x)/\langle[X_x]\rangle \ar[l]_(0.54)\sim
}
$$

La partie supérieure de ce diagramme montre que $(s\circ i)(b) \in \Br(X_{\A^1_k})$
s'envoie sur la classe de~$b$ dans $H^1(\A^1_k,\Erond)/\langle[\Xrond]\rangle$
par la composée des flèches horizontales de la ligne du milieu,
compte tenu de la définition du morphisme~$i$
(cf.~notamment les remarques qui précèdent la proposition~\ref{ch1dbrd}).
La partie inférieure du digramme permet d'en déduire le résultat voulu.
\end{demo}

\bigskip
Posons $A = (s\circ i)(b) \in \Br(X_{\A^1_k})$.
Nous avons maintenant établi que
\begin{equation}
\label{ch1sdeqct}
\langle \alpha, [X_x]\rangle=\sum_{v\in \Omega}\inv_v A(Q'_v)\rlap{\text{.}}
\end{equation}
Comme $A\in B_1$ et $\inv_v A(Q'_v)=\inv_v A(Q_v)$ pour tout $v \in S_1$,
on tire des équations~(\ref{ch1sdeqs1}) et~(\ref{ch1sdeqct})
l'égalité $\langle \alpha, [X_x]\rangle=\sum_{v \in \Omega\setminus S_1}\inv_v A(Q'_v)$.

\bigskip
\begin{lemme}
\label{ch1sdutila1}
On a $A(Q'_v)=0$ pour tout $v \in \Omega\setminus S_1$.
\end{lemme}

\bigskip
\begin{demo}
Soit $v \in \Omega\setminus S_1$.  Comme par hypothèse~$x$ est entier hors de~$S_1$,
il existe un morphisme $\xtilde \colon \Spec(\Orond_v) \rightarrow \A^1_{\Orond_S}$
rendant le carré
$$
\myxyhook\xymatrix{
\Spec(k_v) \ar[r]^(0.6){Q'_v} \ar[d] & \Xtilde \ar[d] \\
\Spec(\Orond_v) \ar@{.>}[ur] \ar[r]^(0.57)\xtilde & \A^1_{\Orond_S}
}
$$
commutatif (sans la flèche en pointillés).
La propreté de $\Xtilde \rightarrow \A^1_{\Orond_S}$
entraîne l'existence d'une flèche en pointillés telle que le diagramme reste commutatif.
Étant donné que $\Br(\Orond_v)=0$, on en déduit que la flèche
$\Br(\Xtilde)\rightarrow \Br(k_v)$ d'évaluation en~$Q'_v$
est nulle. En particulier $A(Q'_v)=0$, puisque $A \in B_1$ et $B_1 \subset \Br(\Xtilde)$.
\end{demo}

\bigskip
Ainsi avons-nous prouvé que $\langle \alpha, [X_x]\rangle=0$ pour
tout $\alpha \in \tors{2}\Sha(k,\Erond_x)$, autrement dit que $x \in \RO$; la démonstration du théorème~\ref{ch1thsecdesc}
est donc achevée.
\end{demo}

\bigskip
\begin{remarques}
(i) L'hypothèse $L_M=\kappa(M)$ pour tout~$M\in\Mrond$ est notamment satisfaite
lorsque $\ff{F}M=\Z/2$ pour tout~$M\in\Mrond$; le
théorème~\ref{ch1thsecdesc} s'applique donc à toutes les surfaces considérées
dans~\cite{css}.

(ii) L'hypothèse $L_M=\kappa(M)$ pour tout $M \in \Mrond$ n'est pas satisfaite
dans la situation considérée par Swinnerton-Dyer dans~\cite[§8.2]{sddiagquartic}.
\end{remarques}

\bigskip
Le résultat suivant, bien connu, précise le lien entre l'ensemble~$\RO$ et les
secondes descentes: si~$E$ est une courbe elliptique sur~$k$ et~$\alpha$
un élément de $\Sel_2(k,E)$, l'image de~$\alpha$
dans $\Sha(k,E)$ est orthogonale à $\tors{2}\Sha(k,E)$
pour l'accouplement de Cassels-Tate si et seulement si~$\alpha$ appartient à l'image
du morphisme $$\Sel_4(k,E) \longrightarrow \Sel_2(k,E)$$ induit par la flèche
$\tors{4}E \rightarrow \tors{2}E$ de multiplication par~$2$.
Notant $D \rightarrow E$ un $2$\nobreakdash-revêtement de~$E$ dont la classe dans $H^1(k,\tors{2}E)$
soit égale à~$\alpha$, cette condition équivaut encore à l'existence
d'un $4$\nobreakdash-revêtement $D' \rightarrow E$ et d'un morphisme $D' \rightarrow D$ tels que le triangle
$$
\myxyhook\xymatrix@R=2.5ex{
& D' \ar[dd] \ar[dl] \\
D \ar[dr] \\
& E
}
$$
soit commutatif et que $D'(\A_k)\neq\emptyset$.  Selon la terminologie consacrée, on dit
alors que le $2$\nobreakdash-revêtement~$D$ «~provient d'une seconde descente~».
La conclusion du théorème~\ref{ch1thsecdesc} signifie que si $X(\A_k)\neq\emptyset$
et si l'obstruction de Brauer-Manin
à l'approximation faible sur~$X$ s'évanouit, on peut trouver beaucoup de points rationnels de~$C$ au-dessus
desquels la fibre de~$\pi$ non seulement possède des points partout localement, mais de plus provient d'une
seconde descente.

\bigskip
Il existe un algorithme, dû à Cassels~\cite{casselssecond}, permettant de calculer
l'accouplement de Cassels-Tate sur deux éléments du groupe de $2$\nobreakdash-Selmer d'une
courbe elliptique dont on connaît une équation de Weierstrass.  À l'aide de cet
algorithme, il est parfois possible,
étant donné
$s \in \SD$, de trouver un ouvert $\Arond \subset \P^1(\A_k)$ tel
que $\pi(X(\A_k))\cap\Arond\neq\emptyset$ et
que pour tout $x \in \RA \cap \Arond$, la classe~$s(x) \in H^1(k,\tors{2}\Erond_x)$
appartienne au groupe de $2$\nobreakdash-Selmer de~$\Erond_x$ et ne soit pas orthogonale
à~$[X_x]$ pour l'accouplement de Cassels-Tate.  En particulier $X_x(k)=\emptyset$
pour tout $x \in \RA \cap \Arond$; la surface~$X$ présente alors un défaut
d'approximation faible.

Le théorème~\ref{ch1thsecdesc} montre, lorsque $L_M=\kappa(M)$ pour tout $M\in\Mrond$,
que \emph{tout défaut d'approximation faible
mis en évidence par cette méthode est expliqué par l'obstruction de Brauer-Manin},
si l'on admet l'hypothèse de Schinzel.  Notons que sa preuve n'est pas constructive;
l'explicitation d'une classe de~$\Br(X)$ responsable de l'obstruction n'est en rien
facilitée.

\bigskip
\begin{remarque}
Lorsque l'hypothèse $L_M=\kappa(M)$ n'est pas satisfaite, on peut malgré tout
appliquer le corollaire~\ref{ch1dtdd} pour $C=U$ et une variante de la
proposition~\ref{ch1dbrd} afin d'obtenir un morphisme canonique $i \colon \TDz
\rightarrow \Br_\hor(X_U)$. Le reste de la démonstration s'applique encore, excepté
le lemme~\ref{ch1sdutila1}. Se souvenant que le point~$x$
provenait d'un triplet $(S_1,\Trond,x)$ admissible (cf.~proposition~\ref{ch1propschinzel}),
on voit que $\inv_v A(Q'_v)=0$ pour tout $v \in \Omega\setminus T(x)$. Seules posent donc
problème les places de $T(x)\setminus S_1$.  En une place $v \in \ft{T}M\cup\{v_M\}$, il
est possible de relier $\inv_v A(Q'_v)$
aux résidus de~$(s\circ i)(b)$ en les composantes irréductibles de~$X_M$, mais on est alors confronté
à une autre difficulté, déjà présente lorsque $L_M=\kappa(M)$ pour tout $M \in \Mrond$
si l'on n'a pas pris soin de choisir~$(s\circ i)(b)$ non ramifié sur~$X_{\A^1_k}$: ces résidus
dépendent
de la section $s \colon \Br_\hor(X_U) \rightarrow \Br(X_U)$ choisie, et l'on ne peut donc
pas espérer qu'ils seront toujours\index{seconde descente|)} nuls\index{accouplement de Cassels-Tate|)} en~$v$.
\end{remarque}

\section{Application aux courbes elliptiques de rang élevé}
\label{ch1parrangel}

Soit $E$ une\index{courbe elliptique de rang élevé|(} courbe elliptique sur $k(t)$, où $k$ est un corps de nombres.
Il~existe un ouvert dense $U \subset \P^1_k$ au-dessus duquel~$E$ s'étend
en un schéma abélien \mbox{$\Erond \rightarrow U$}.  D'après un théorème de Silverman,
la flèche de spécialisation $E(k(t)) \rightarrow \Erond_x(k)$
est injective pour presque tout $x \in U(k)$ (cf.~\cite[Ch.~III,
Th.~11.4 et Ex.~3.16]{silvaec2});
l'expression «~presque tout~$x$~» signifie
bien sûr ici «~tout~$x$ sauf un nombre fini~».
En particulier, l'inégalité
$$
\rg(E(k(t))) \leq \rg(\Erond_x(k))
$$
vaut pour presque tout $x \in U(k)$.
Lorsqu'elle est stricte pour presque tout~$x$, la courbe elliptique~$E$ est
dite \emph{de rang élevé}.

À ce jour, aucun exemple de courbe elliptique de rang élevé
n'est connu.
Cassels et Schinzel~\cite{cassschin}
puis Rohrlich~\cite{rohrlich} ont fabriqué des courbes elliptiques isotriviales
sur~$\Q(t)$ dont ils montrent qu'elles sont de rang élevé en admettant
la conjecture de parité.
Dans la direction opposée, mais toujours sur~$\Q$,
Conrad, Conrad et Helfgott~\cite{cchrootnumbers} prouvent que toute courbe elliptique
de rang élevé sur~$\Q(t)$ est isotriviale, en admettant deux conjectures
de théorie analytique des nombres (la conjecture de Chowla et une conjecture
quantifiant le nombre de valeurs sans facteur carré prises par un polynôme de~$\Z[t]$)
ainsi qu'une conjecture de densité concernant les rangs des courbes elliptiques $\Erond_x$
lorsque~$x$ varie.  Le théorème suivant est également un résultat de non existence
de courbes elliptiques de rang élevé.  Il est lui aussi conditionnel, et par ailleurs
ne s'applique que sous certaines hypothèses sur les courbes elliptiques considérées;
en revanche, il présente l'avantage d'être valable sur tout corps de nombres
et de ne faire intervenir aucune conjecture concernant les courbes elliptiques.

\bigskip
\begin{theoreme}
\label{ch1rangeleve}
Admettons l'hypothèse de Schinzel.
Soient $k$ un corps de nombres et~$E$ une courbe elliptique sur~$k(t)$.
Notons $X \rightarrow \P^1_k$ un modèle propre et régulier
de~$E$ et $\Br^0(X)$ le noyau de la flèche
$\Br(X) \rightarrow \Br(\P^1_k)$ induite par la section nulle.
Si les points d'ordre~$2$ de~$E$ sont rationnels, si $E$ est à réduction de type~$I_2$
en tout point de mauvaise réduction et si $\tors{2}\Br^0(X)=0$, alors $E$ n'est
pas de rang élevé.
\end{theoreme}

\bigskip
\begin{demo}
Le morphisme $\pi \colon X \rightarrow \P^1_k$ satisfait aux
hypothèses du paragraphe~\ref{ch1hypnot}, dont on reprend les notations.
Le morphisme canonique $\Br^0(X) \rightarrow \Br_\hor(X)$ est
un isomorphisme; en effet, il est injectif puisque $\Br_\vert(X) \cap \Br^0(X)=0$,
surjectif puisque la suite exacte
$$
\xymatrix{
0 \ar[r] & \Br^0(X) \ar[r] & \Br(X) \ar[r] & \Br(\P^1_k) \ar[r] & 0
}
$$
est scindée. Il s'ensuit que $\tors{2}\Br_\hor(X)=0$.
Par ailleurs, comme la courbe elliptique~$E$ est à réduction de type~$I_2$ en tout point de~$\Mrond$,
on a \mbox{$L_M=\kappa(M)$} pour tout $M \in \Mrond$.
Le corollaire~\ref{ch1dtdd} et la proposition~\ref{ch1dbrd} permettent d'en déduire
l'existence d'une injection $\TD \hookrightarrow \tors{2}\Br_\hor(X)$,
compte tenu que \mbox{$[\Xrond]=0$}. Il en résulte que $\TD=0$, d'où finalement $\SD=E(k(t))/2$.

Supposons que la courbe elliptique~$E$ soit de rang élevé.
Il est clair que~$E$ ne peut pas être constante; autrement dit, on
a nécessairement $\Mrond \neq \emptyset$.
Le théorème~\ref{ch1appobm1} montre maintenant que l'ensemble $\RD$ est infini,
étant donné que $\pi(X(k))=\P^1(k)$.  Comme~$E$ est de rang élevé, il existe donc
$x \in \RD$ tel que $\rg(E(k(t)))<\rg(\Erond_x(k))$.  On a alors
$$
\begin{array}{ccll}
\rg(\Erond_x(k)) &=& \dim_{\ff{\F}2}(\Erond_x(k)/2)-2 & \qquad\!\!  \text{(car la $2$\nobreakdash-torsion de $\Erond_x$ est rationnelle)} \\
& \leq & \dim_{\ff{\F}2}\Sel_2(k, \Erond_x) -2 \\
& \leq & \dim_{\ff{\F}2}(\SD) -2 & \qquad\!\! \text{(car $x \in \RD$)} \\
& = & \dim_{\ff{\F}2}(E(k(t))/2) - 2 & \qquad\!\! \text{(car $\SD=E(k(t))/2$)} \\
& = & \rg(E(k(t))) & \qquad\!\! \text{(car la $2$\nobreakdash-torsion de $E$ est rationnelle),}
\end{array}
$$
d'où une contradiction.
\end{demo}

\bigskip
\begin{remarques}
(i) Le théorème~\ref{ch1rangeleve} ne concerne que les courbes elliptiques
non isotriviales, bien que cette restriction ne figure pas explicitement
dans son énoncé.  En effet, une courbe elliptique isotriviale sur~$k(t)$ ne peut
avoir réduction semi-stable sur~$\P^1_k$ à moins d'être constante:
l'invariant~$j$ admet un pôle en tout point de réduction multiplicative
(cf.~\cite[Ch.~IV, §9]{silvaec2}).

(ii) Le groupe $\Br^0(X)$ joue le rôle d'un «~groupe de Tate-Shafarevich~» pour la
courbe elliptique $E/k(t)$.
Il se plonge en tout cas dans le groupe\index{groupe de Tate-Shafarevich} de Tate-Shafarevich géométrique de~$E$
et l'on peut en principe calculer son sous-groupe de $2$\nobreakdash-torsion au moyen d'une $2$\nobreakdash-descente
(cf.~suite exacte~(\ref{ch1sesgtsg})).
De ce point de vue, la preuve du théorème~\ref{ch1rangeleve}
consiste essentiellement à exhiber une infinité de points $x \in U(k)$ pour lesquels
le groupe de $2$\nobreakdash-Selmer de~$\Erond_x$ ne soit pas plus gros que le «~groupe de $2$\nobreakdash-Selmer~»
de~$E$; comme ce dernier est par hypothèse réduit à $E(k(t))/2$, le rang de~$\Erond_x$
est nécessairement inférieur ou égal à celui de~$E$.

(iii) Il est bien sûr possible de déduire du théorème~\ref{ch1appobm1} un résultat
analogue au théorème~\ref{ch1rangeleve}
pour les courbes elliptiques à réduction semi-stable, sans restriction sur
le type de réduction, mais
la condition à imposer sur le groupe de Brauer est alors moins agréable à\index{courbe elliptique de rang élevé|)} exprimer.
\end{remarques}

\chapternotoc[Arithmétique des pinceaux semi-stables II]{Arithmétique des pinceaux semi-stables de courbes de genre~$1$ dont
les jacobiennes possèdent une section d'ordre~$2$}
\addcontentsline{toc}{chapter}{\chaptername\ \thechapter.\ Arithmétique des pinceaux semi-stables de courbes de genre~$1$ (seconde partie)}%

\section{Introduction}

À la suite des travaux de Colliot-Thélène, Skorobogatov et Swinnerton-Dyer~\cite{css}
concernant les points rationnels des surfaces munies d'un pinceau de courbes de genre~$1$
et de période~$2$
dont les jacobiennes ont leur $2$\nobreakdash-torsion entièrement rationnelle, Bender et
Swinnerton-Dyer~\cite{bsd} étudièrent la possibilité d'obtenir
des résultats similaires pour
les surfaces munies d'un pinceau de courbes de genre~$1$ et de période~$2$ dont les jacobiennes
ont \emph{au moins un} point d'ordre~$2$ rationnel.
Leur conclusion fut que les techniques de~\cite{css} s'appliquent encore dans ce
cadre, à condition d'effectuer des descentes par $2$\nobreakdash-isogénie simultanément sur
deux familles de courbes elliptiques au lieu d'une simple $2$\nobreakdash-descente complète
sur la fibration jacobienne du pinceau considéré.
Ils aboutirent ainsi à deux théorèmes, de domaines d'application disjoints;
nous renvoyons le lecteur à~\cite[§1]{bsd} pour leurs énoncés précis.
Colliot-Thélène~\cite{ctbsd} reformula ensuite ces théorèmes et leurs preuves en
termes d'obstruction de Brauer-Manin verticale.

À la fin de~\cite{bsd}, Swinnerton-Dyer expose quelques idées laissant espérer
que l'on puisse se servir du théorème~1 de~\cite{bsd} (ou du théorème~B de~\cite{ctbsd})
afin d'établir qu'une surface de del Pezzo de degré~$4$ admet un point rationnel
dès que l'obstruction de Brauer-Manin ne s'y oppose pas.  Il s'avère
que les résultats de~\cite[§6]{bsd} sont incorrects mais que la construction
proposée par Swinnerton-Dyer permet néanmoins de prouver que le principe
de Hasse vaut
sur de telles surfaces sous des conditions très générales,
quoique au prix d'un travail conséquent; c'est là le sujet du chapitre~3.
Cependant, les théorèmes de~\cite{bsd} et de~\cite{ctbsd} sont insuffisants
pour l'application en vue: ils contiennent plusieurs hypothèses techniques
qui ne seront pas satisfaites par les pinceaux considérés au chapitre~3,
notamment les conditions~(1) et~(2) de \cite[Theorem~B]{ctbsd}
(cf.~aussi les conditions~1 et~2 de~\cite{bsd} et l'avant-dernier paragraphe
de~\cite[p.~323]{bsd}; celui-ci contient un argument permettant d'affaiblir
quelque peu les conditions~1 et~2, mais pour le chapitre~3 il est nécessaire
de supprimer entièrement la condition~2).

Ce chapitre a pour but premier d'établir une version des théorèmes~A et~B
de~\cite{ctbsd} épurée de toute hypothèse parasite et
de laquelle il soit possible de déduire des résultats sur les surfaces de del Pezzo de degré~$4$.
À cette fin, nous emploierons systématiquement les propriétés des modèles
de Néron là où~\cite{bsd} et~\cite{ctbsd} font appel à des équations explicites.
D'autre part, nous généraliserons
au cas de réduction semi-stable arbitraire
les théorèmes de~\cite{ctbsd},
qui ne s'appliquent que lorsque la jacobienne de la fibre générique du pinceau
considéré est à réduction de type~$I_1$ ou~$I_2$ en chaque point de mauvaise réduction,
et nous supprimerons toute hypothèse
concernant le rang de Mordell-Weil de cette courbe elliptique, comme au chapitre~1.
Ces améliorations permettent notamment d'obtenir un seul énoncé
qui implique à la fois les théorèmes~A et~B de~\cite{ctbsd} (et
les théorèmes~1 et~2 de~\cite{bsd}), contribuant
ainsi à les clarifier.

Les idées générales qui sous-tendent la preuve du théorème principal de ce chapitre
sont bien sûr les mêmes que dans~\cite{bsd}, \cite{ctbsd}, \cite{css} ou que
dans le chapitre~1.  Néanmoins, leur mise en \oe{}uvre contient une innovation substantielle,
qui mérite d'être signalée: l'argument se déroule directement au niveau des groupes
de Selmer des courbes elliptiques concernées, sans plus passer par l'étude
d'un certain accouplement dont les noyaux à gauche et à droite sont isomorphes à
ces groupes de Selmer et dont la variation en famille est contrôlable.
Les sous-espaces totalement isotropes maximaux~$K_v$
(cf.~lemme~\ref{ch1existekv}, \cite[Proposition~1.1.2]{css}, \cite[Lemma~8]{bsd},
\cite[Proposition~1.1.2]{ctbsd}), dont le rôle semblait pour le moins mystérieux,
ont ainsi totalement disparu.
Cette modification de l'argument ne semble pas être une
opération de nature purement formelle; un indice en ce sens
est l'utilisation dans la preuve ci-dessous de réciprocités dans lesquelles
interviennent des algèbres de quaternions sur~$k(t)$ qui n'avaient été considérées
dans aucun des articles~\cite{bsd}, \cite{ctbsd}, \cite{css}
(cf.~paragraphe~\ref{ch2parrec} et preuve de la proposition~\ref{ch2pluslongueprop}).

Outre la simplification conceptuelle que cette innovation apporte,
elle représente un premier pas, quoique modeste,
vers l'obtention de résultats indépendants de l'hypothèse
de Schinzel.  Nombre de notions (principe de Hasse, obstruction de Brauer-Manin)
et d'outils (notamment, le «~lemme formel~» d'Harari) utiles dans l'étude
des questions d'existence de points rationnels possèdent un analogue
adapté aux questions d'existence de $0$\nobreakdash-cycles de degré~$1$.
Dans ce contexte, l'analogue de l'hypothèse de Schinzel est un théorème,
connu sous le nom d'astuce de Salberger (cf.~\cite[§3]{csscrelle98}).
Grâce à ce théorème, on sait démontrer inconditionnellement
l'existence d'un $0$\nobreakdash-cycle de degré~$1$
sur les surfaces
admettant un pinceau de coniques, dès que l'obstruction de Brauer-Manin ne
s'y oppose pas (cf.~\cite[§4]{csscrelle98}).
Si l'on cherche à combiner l'astuce de Salberger et les méthodes de~\cite{css}
afin d'étudier les $0$\nobreakdash-cycles de degré~$1$ sur des surfaces munies d'un
pinceau de courbes de genre~$1$, on est naturellement confronté à deux
difficultés, liées au fait que l'on doit comparer uniformément les groupes de Selmer de
courbes elliptiques~$\Erond_x$ pour divers points \emph{fermés}
(et non plus rationnels) $x\in U$,
où $\Erond \rightarrow U$ est une courbe elliptique relative et~$U$ un ouvert de~$\P^1_k$:
d'une part, si~$S$ est un ensemble fixé de places de~$k$, arbitrairement grand
mais indépendant de~$x$, on ne contrôle pas le groupe des $S$\nobreakdash-unités du corps~$\kappa(x)$,
et d'autre part, on ne contrôle pas non plus le groupe de classes de $S$\nobreakdash-idéaux de~$\kappa(x)$.
Le groupe des $S$\nobreakdash-unités intervient dans le diagramme~(\ref{ch1diagpsi}); quant au
groupe de classes de $S$\nobreakdash-idéaux, il est nécessaire qu'il soit nul pour que
l'on puisse \emph{définir} l'accouplement de la proposition~\ref{ch1estsymetrique}
(resp.~l'accouplement de \cite[Lemma~10]{bsd}, \cite[Proposition~1.4.3]{ctbsd}).
L'intérêt de la preuve que nous donnons dans ce chapitre est qu'en évitant l'emploi
d'un accouplement analogue à celui de la proposition~\ref{ch1estsymetrique},
elle fait tout simplement disparaître la seconde de ces difficultés; nous ne supposerons
à aucun moment qu'un quelconque ensemble fini de places de~$k$ contient un système
de générateurs du groupe de classes.

\section{Énoncé du résultat principal et applications}
\label{ch2enonceresul}

Soient~$k$\glossary{$k$, $C$, $\eta$, $X$, $\pi$}
un corps de caractéristique~$0$ et~$C$ un schéma de Dedekind connexe sur~$k$,
de point générique~$\eta$.
Supposons donnée une surface~$X$ lisse et géométriquement connexe sur~$k$, munie
d'un morphisme propre et plat $\pi \colon X \rightarrow C$ dont la fibre générique~$X_\eta$
est une courbe lisse de genre~$1$ sur~$K=\kappa(C)$ et dont toutes les fibres sont réduites.
Supposons de plus que la période de la courbe~$X_\eta$ divise~$2$, c'est-à-dire
que la classe de $H^1(K, E'_\eta)$ définie par le torseur~$X_\eta$ soit tuée par~$2$,
en notant~$E'_\eta$ la jacobienne de~$X_\eta$.  Supposons enfin que la courbe
elliptique~$E'_\eta$ soit à réduction semi-stable en tout point fermé de~$C$ et
qu'elle possède un point $K$\nobreakdash-rationnel~$P'$ d'ordre~$2$.

Notons\glossary{$E'_\eta$, $E''_\eta$, $P'$, $P''$} $E''_\eta$ la courbe elliptique quotient de~$E'_\eta$ par~$P'$,
$\phi'_\eta \colon E'_\eta \rightarrow E''_\eta$\glossary{$\phi'$, $\phi''$, $\phi'^0$,
$\phi''^0$}
la $2$\nobreakdash-isogénie associée
et $\phi''_\eta \colon E''_\eta \rightarrow E'_\eta$ la duale de~$\phi'_\eta$,
de sorte que $\phi'_\eta \circ \phi''_\eta = 2$ et $\phi''_\eta \circ \phi'_\eta = 2$.
Le noyau de $\phi''_\eta$ contient un unique point $K$\nobreakdash-rationnel non nul;
on le note $P'' \in E''_\eta(K)$.
Soient $\Erond'$\glossary{$\Erond'$, $\Erond''$, $\Erond'^0$, $\Erond''^0$} et~$\Erond''$ les modèles de Néron\index{modèle de Néron} respectifs de~$E'_\eta$ et
de~$E''_\eta$ sur~$C$.
D'après la propriété universelle des modèles de Néron, les isogénies~$\phi'_\eta$
et~$\phi''_\eta$ se prolongent en des morphismes de schémas en groupes
$\phi' \colon \Erond' \rightarrow \Erond''$ et $\phi'' \colon \Erond'' \rightarrow \Erond'$.
Ceux-ci induisent des morphismes surjectifs $\phi'^0 \colon \Erond'^0 \rightarrow \Erond''^0$
et $\phi''^0 \colon \Erond''^0 \rightarrow \Erond'^0$,
où $\Erond'^0 \subset \Erond'$ et~$\Erond''^0 \subset \Erond''$ désignent
les composantes neutres.

Notons $\Mrond \subset C$ l'ensemble des points fermés de mauvaise réduction
pour~$E'_\eta$.  C'est aussi l'ensemble des points fermés de mauvaise réduction
pour~$E''_\eta$, et la courbe elliptique~$E''_\eta$ est à
réduction semi-stable en ces points (cf.~\cite[7.3/7]{blr}).
Fixons un ouvert dense $U \subset C$ au-dessus duquel~$\pi$ est lisse.
Pour $M \in \Mrond$, notons respectivement~$\ff{F'}M$ et~$\ff{F''}M$\glossary{$\ff{F'}M$,
$\ff{F''}M$, $\phi'_M$, $\phi''_M$} les
$\kappa(M)$\nobreakdash-schémas en groupes finis étales $\Erond'_M/\Erond'^0_M$ et~$\Erond''_M/\Erond''^0_M$.
Les morphismes~$\phi'$ et~$\phi''$ induisent des morphismes $\phi'_M \colon \ff{F'}M
\rightarrow \ff{F''}M$ et $\phi''_M \colon \ff{F''}M \rightarrow \ff{F'}M$.

Posons\glossary{$\Mrond$, $\Mrond'$, $\Mrond''$}
$$
\phantom{\Mrond''}\llap{$\Mrond'$}=\bigensemble{M \in \Mrond}{\phi'_M\text{ est injectif}}
$$
et
$$
\Mrond''=\bigensemble{M \in \Mrond}{\phi''_M\text{ est injectif}}\!\rlap{\text{.}}
$$
La proposition~\ref{annpropppp} montre que $\Mrond' \cap \Mrond'' = \emptyset$
et que $\Mrond' \cup \Mrond'' = \Mrond$. Elle montre aussi que les
morphismes~$\phi'_M$ et~$\phi''_M$ s'insèrent dans des suites exactes
\begin{align}
\begin{aligned}
\smash[t]{\myxyhook\xymatrix@R=3.5ex@C=6ex{
0 \ar[r] & \ff{F'}M \ar[r]^{\phi'_M} & \ff{F''}M \ar[r] & \Z/2 \ar[r] & 0 \\
0 \ar[r] & \Z/2 \ar[r] & \ff{F''}M \ar[r]^{\phi''_M} & \ff{F'}M \ar[r] & 0
}}
\end{aligned}\\
\intertext{\smash{si $M \in \Mrond'$ et}}
\begin{aligned}
\smash[t]{\myxyhook\xymatrix@R=3.5ex@C=6ex{
0 \ar[r] & \Z/2 \ar[r] & \ff{F'}M \ar[r]^{\phi'_M} & \ff{F''}M \ar[r] & 0 \\
0 \ar[r] & \ff{F''}M \ar[r]^{\phi''_M} & \ff{F'}M \ar[r] & \Z/2 \ar[r] & 0
}}
\end{aligned}
\end{align}
\vskip-6.4pt\noindent{}si $M \in \Mrond''$.

Les suites exactes
$$
\myxyhook\xymatrix{
0 \ar[r] & \Z/2 \ar[r]^{P'} & E'_\eta \ar[r]^{\phi'_\eta} & E''_\eta \ar[r] & 0
}
$$
et
$$
\myxyhook\xymatrix{
0 \ar[r] & \Z/2 \ar[r]^{P''} & E''_\eta \ar[r]^{\phi''_\eta} & E'_\eta \ar[r] & 0
}
$$
induisent des suites exactes
$$
\myxyhook\xymatrix{
0 \ar[r] & E''_\eta(K)/\Im(\phi'_\eta) \ar[r] & K^\star/K^{\star 2} \ar[r]^(0.44){\alpha'} &
\tors{\phi'_\eta}H^1(K,E'_\eta) \ar[r] & 0
}
$$
et
$$
\myxyhook\xymatrix{
0 \ar[r] & E'_\eta(K)/\Im(\phi''_\eta) \ar[r] & K^\star/K^{\star 2} \ar[r]^(0.44){\alpha''} &
\tors{\phi''_\eta}H^1(K,E''_\eta) \ar[r] & 0 \rlap{\text{.}}
}
$$

\medskip\noindent{}Pour~$M \in C$, notons~$\Orond_M^\sh$ l'anneau strictement local de~$C$ en~$M$
et~$K_M^\sh$ son corps des fractions.
Rappelons que pour $\Erond \in \{\Erond', \Erond''\}$, si~$E_\eta$ désigne la fibre
générique de $\Erond \rightarrow C$,
le groupe $H^1(C,\Erond)$ s'identifie
au noyau du produit $H^1(K,E_\eta) \rightarrow \prod_{M \in C}H^1(K_M^\sh, E_\eta)$
(cf.~suite exacte~(\ref{ch1setsg})).
Notons respectivement $\SG_{\phi'}(C,\Erond')$
et $\SG_{\phi''}(C,\Erond'')$\glossary{$\SG_{\phi'}(C,\Erond')$, $\SG_{\phi''}(C,\Erond'')$}
les images réciproques des sous-groupes $\tors{\phi'}H^1(C,\Erond')
\subset \tors{\phi'_\eta}H^1(K,E'_\eta)$
et $\tors{\phi''}H^1(C,\Erond'')
\subset \tors{\phi''_\eta}H^1(K,E''_\eta)$ par~$\alpha'$ et~$\alpha''$,
de sorte que l'on obtient des suites exactes
\begin{equation}
\label{ch2sesgtsgp}
\myxyhook\xymatrix@C=4.53ex{
0 \ar[r] & E''_\eta(K)/\Im(\phi'_\eta) \ar[r] & \SG_{\phi'}(C,\Erond') \ar[r] &
\tors{\phi'}H^1(C,\Erond') \ar[r] & 0
}
\end{equation}
et
\begin{equation}
\label{ch2sesgtsgpp}
\myxyhook\xymatrix@C=4ex{
0 \ar[r] & E'_\eta(K)/\Im(\phi''_\eta) \ar[r] & \SG_{\phi''}(C,\Erond'') \ar[r] &
\tors{\phi''}H^1(C,\Erond'') \ar[r] & 0\rlap{\text{.}}
}
\end{equation}
Les groupes $\SG_{\phi'}(C,\Erond')$ et $\SG_{\phi''}(C,\Erond'')$ méritent qu'on les appelle
respectivement \emph{groupe de $\phi'$\nobreakdash-Selmer géométrique\index{groupe de Selmer!géométrique} de~$E'_\eta$} et \emph{groupe de $\phi''$\nobreakdash-Selmer
géométrique de~$E''_\eta$} (cf.~\cite[§4.2]{css}).

\bigskip
\begin{proposition}
\label{ch2propcompct}
On a
$\SG_{\phi'}(C,\Erond')=H^1(C \setminus \Mrond', \Z/2)$ et
$\SG_{\phi''}(C,\Erond'')=H^1(C \setminus \Mrond'', \Z/2)$
en tant que sous-groupes de $K^\star/K^{\star 2}=H^1(K,\Z/2)$.
\end{proposition}

\bigskip
\begin{demo}
Les deux assertions étant symétriques, il suffit de démontrer la première.
Vu la définition du groupe $\SG_{\phi'}(C,\Erond')$, il suffit de
prouver que $E''_\eta(K_M^\sh)/\Im(\phi'_\eta)=H^1(K_M^\sh,\Z/2)$
pour tout $M \in \Mrond'$ et $E''_\eta(K_M^\sh)/\Im(\phi'_\eta)=0$
pour tout $M \in C \setminus \Mrond'$.
Autrement dit, compte tenu que $H^1(K_M^\sh,\Z/2)=\Z/2$, il suffit
de prouver que l'application $E'_\eta(K_M^\sh) \rightarrow E''_\eta(K_M^\sh)$
induite par~$\phi'_\eta$ est surjective si et seulement si $M \not\in \Mrond'$,
pour $M \in C$. Ceci résulte du lemme~\ref{annlemmeisogekmsh}, de la
proposition~\ref{annpropppp}
et de la définition de~$\Mrond'$.
\end{demo}

\bigskip
Notons enfin $\SG_2(C,\Erond') \subset H^1(K, \tors{2}E'_\eta)$\glossary{$\SG_2(C,\Erond')$}
l'image réciproque de $H^1(C,\Erond') \subset H^1(K, E'_\eta)$
par le morphisme canonique $H^1(K,\tors{2}E'_\eta) \rightarrow H^1(K,E'_\eta)$.
Le morphisme~$\phi'$ induit un diagramme commutatif
\begin{equation}
\label{ch2hypdiags2spp}
\begin{aligned}
\myxyhook\xymatrix@C=5ex{
0 \ar[r] & E'_\eta(K)/2 \ar[d] \ar[r] & \SG_2(C,\Erond') \ar[d] \ar[r] & \tors{2}H^1(C,\Erond')
\ar[r] \ar[d] & 0 \\
0 \ar[r] & E'_\eta(K)/\Im(\phi''_\eta) \ar[r] & \SG_{\phi''}(C,\Erond'') \ar[r] &
\tors{\phi''}H^1(C,\Erond'') \ar[r] & 0
}
\end{aligned}
\end{equation}
dont les lignes sont exactes.

Comme les fibres de~$\pi$ sont réduites, on a~$X_M(K_M^\sh)\neq \emptyset$
pour tout $M \in C$.  La classe de $H^1(K,E'_\eta)$ définie par le torseur~$X_\eta$
appartient donc au sous-groupe $H^1(C,\Erond')$.  Cette classe correspond à
un faisceau représentable (cf.~\cite[Theorem~4.3]{milneec}), d'où l'existence
d'un torseur $\Xrond \rightarrow C$\glossary{$\Xrond$, $\Xrond''$} sous~$\Erond'$ dont la fibre générique est
égale à~$X_\eta$. Notons~$[\Xrond]$ sa classe dans $H^1(C,\Erond')$
et~$[\Xrond''] \in H^1(C,\Erond'')$ l'image de~$[\Xrond]$ par la flèche verticale de droite
du diagramme~(\ref{ch2hypdiags2spp}).

Pour $M \in \Mrond$, notons~$L'_M$\glossary{$L'_M$, $L''_M$, $\delta_M'$, $\delta_M''$}
(resp.~$L''_M$) l'extension quadratique ou
triviale minimale de~$\kappa(M)$ sur laquelle le $\kappa(M)$\nobreakdash-groupe $\ff{F'}M$
(resp.~$\ff{F''}M$) devient constant (cf.~proposition~\ref{annpropcyclpair}).
Les groupes $\ff{F'}M(L'_M)$ et $\ff{F''}M(L''_M)$ sont cycliques en vertu
de l'hypothèse de semi-stabilité (\emph{loco citato}).
Soient $\delta'_M \colon H^1(C,\Erond') \rightarrow H^1(L'_M,\ff{F'}M)$
et $\delta''_M \colon H^1(C,\Erond'') \rightarrow H^1(L''_M,\ff{F''}M)$ les composées
des flèches induites par les morphismes de faisceaux $\Erond' \rightarrow i_{M\star} \ff{F'}M$
et $\Erond'' \rightarrow i_{M\star} \ff{F''}M$, où~$i_M \colon \Spec(\kappa(M)) \rightarrow C$
désigne l'inclusion canonique, et des flèches de restriction $H^1(\kappa(M),\ff{F'}M)
\rightarrow H^1(L'_M,\ff{F'}M)$ et $H^1(\kappa(M),\ff{F''}M) \rightarrow H^1(L''_M,
\ff{F''}M)$.
Soient $\TDCp{C}$\glossary{$\TDCp{C}$, $\TDCpp{C}$}
le noyau de l'application composée
\begin{equation}
\label{ch2flechedeftdcp}
\myxyhook\xymatrix@C=8ex{
\tors{\phi'}H^1(C, \Erond') \ar[r]^(0.43){\prod \delta'_M} & \displaystyle
\prod_{M \in \Mrond''} H^1(L'_M, \ff{F'}M) \ar[r] &
\displaystyle \prod_{M \in \Mrond''} \frac{H^1(L'_M, \ff{F'}M)}{\langle \delta'_M([\Xrond]) \rangle}
}
\end{equation}
et $\TDCpp{C}$ le noyau de l'application composée
\begin{equation}
\label{ch2flechedeftdcpp}
\myxyhook\xymatrix@C=8ex{
\tors{\phi''}H^1(C, \Erond'') \ar[r]^(0.43){\prod \delta''_M} & \displaystyle
\prod_{M \in \Mrond'} H^1(L''_M, \ff{F''}M) \ar[r] &
\displaystyle
\prod_{M \in \Mrond'} \frac{H^1(L''_M, \ff{F''}M)}{\langle \delta''_M([\Xrond'']) \rangle}\rlap{\text{.}}
}
\end{equation}
On dira que\index{condition (D)@\condD{}} \emph{la \condDC{C} est satisfaite} si\glossary{$\SDCp{C}$,
$\SDCpp{C}$, $(\mathrm{D}/C)$} $\TDCp{C} \subset \{0,[\Xrond]\}$
et si~$\TDCpp{C}$ est engendré par~$[\Xrond'']$.
Notons enfin $\SDCp{C} \subset \SG_{\phi'}(C,\Erond')$ et
$\SDCpp{C} \subset \SG_{\phi''}(C,\Erond'')$ les images réciproques respectives
de~$\TDCp{C}$ et~$\TDCpp{C}$ par les flèches de droite des suites exactes~(\ref{ch2sesgtsgp})
et~(\ref{ch2sesgtsgpp}).

Supposons maintenant que~$k$ soit un corps de nombres et que~$C$ soit un
ouvert de~$\P^1_k$.  On note comme
précédemment~$\Omega$ l'ensemble des places de~$k$, $\Omega_f \subset \Omega$
l'ensemble de ses places finies et~$\A_k$ l'anneau des adèles de~$k$.
Un élément de~$K^\star$ peut être considéré comme une
fonction rationnelle non nulle sur $\P^1_\Orond$.  Par ailleurs, chaque place
finie de~$k$ définit un point de codimension~$1$ de~$\P^1_\Orond$ et donc une
valuation sur $\kappa(\P^1_\Orond)=K$.  On notera
$\SG_{\phi',S}(C,\Erond')$\glossary{$\SG_{\phi',S}(C,\Erond')$, $\SG_{\phi'',S}(C, \Erond'')$}
(resp.~$\SG_{\phi'',S}(C, \Erond'')$) le sous-groupe de
$\SG_{\phi'}(C,\Erond')$ (resp.~de~$\SG_{\phi''}(C,\Erond'')$) constitué des
classes appartenant au noyau de la flèche $K^\star/K^{\star 2} \rightarrow
(\Z/2)^{\Omega_f \setminus (S \cap \Omega_f)}$ induite par les valuations normalisées
associées aux places finies de $\Omega \setminus S$, pour $S \subset \Omega$.

Soit
$$
\RA = \bigensemble{x \in U(k)}{X_x(\A_k)\neq\emptyset}
$$
et\glossary{$\RA$, $\RDCS{C}{S}$}
pour $S \subset \Omega$, soit~$\RDCS{C}{S}$ l'ensemble des~$x \in \Rrond_A$ tels que tout élément
du groupe de $\phi'_x$\nobreakdash-Selmer de~$\Erond'_x$ appartienne à l'image de la composée
$$
\SDCp{C} \cap \SG_{\phi',S}(C,\Erond')
\subset \SG_{\phi'}(C,\Erond')=H^1(C \setminus \Mrond', \tors{\phi'}\Erond') \rightarrow
H^1(k, \tors{\phi'_x}\Erond'_x)
$$
(cf.~proposition~\ref{ch2propcompct}; la dernière flèche est l'évaluation en~$x$),
que tout élément
du groupe de $\phi''_x$\nobreakdash-Selmer de~$\Erond''_x$ appartienne à l'image de la composée
$$
\SDCpp{C} \cap \SG_{\phi'',S}(C,\Erond'')
\subset \SG_{\phi''}(C,\Erond'')=H^1(C \setminus \Mrond'', \tors{\phi''}\Erond'')
\rightarrow H^1(k, \tors{\phi''_x}\Erond''_x)\rlap{\text{,}}
$$
et que les restrictions des flèches
\begin{equation}
\label{ch2pourinjfp}
\SG_{\phi'}(C,\Erond')=H^1(C \setminus \Mrond', \tors{\phi'}\Erond') \rightarrow
H^1(k, \tors{\phi'_x}\Erond'_x)
\end{equation}
et
\begin{equation}
\label{ch2pourinjfpp}
\SG_{\phi''}(C,\Erond'')=H^1(C \setminus \Mrond'', \tors{\phi''}\Erond'')
\rightarrow H^1(k, \tors{\phi''_x}\Erond''_x)
\end{equation}
aux images réciproques respectives
de $\{0,[\Xrond]\}$ et de $\{0,[\Xrond'']\}$ par les flèches
\begin{equation}
\label{ch2injfp}
\SG_{\phi'}(C,\Erond')\rightarrow H^1(C,\Erond')
\end{equation}
et
\begin{equation}
\label{ch2injfpp}
\SG_{\phi''}(C,\Erond'')\rightarrow H^1(C,\Erond'')
\end{equation}
issues
des suites exactes~(\ref{ch2sesgtsgp}) et~(\ref{ch2sesgtsgpp}) soient injectives.
Pour $S=\Omega$, on note simplement\glossary{$\RDC{C}$, $\RDS{S}$, $\RDSz{S}$, $\RD$, $\RDz$} $\RDC{C}=\RDCS{C}{\Omega_f}$ l'ensemble $\RDCS{C}{S}$.

Lorsque~$C=\P^1_k$, on prend pour~$U$ le plus grand ouvert de~$\A^1_k$ au-dessus
duquel~$\pi$ est lisse et l'on
note d'une part $\RDS{S}$\glossary{$\TDp$, $\TDzp$, $\TDpp$, $\TDzpp$}\glossary{$\SDp$, $\SDzp$, $\SDpp$,
$\SDzpp$}\glossary{\cD{}, \cDz{}},
$\RD$, $\TDp$, $\TDpp$, $\SDp$, $\SDpp$ et~\cD{}
les ensembles $\RDCS{\P^1_k}{S}$, $\RDC{\P^1_k}$, $\TDCp{\P^1_k}$, $\TDCpp{\P^1_k}$,
$\SDCp{\P^1_k}$, $\SDCpp{\P^1_k}$ et la \condDC{\P^1_k}, et d'autre
part $\RDSz{S}$, $\RDz$, $\TDzp$, $\TDzpp$, $\SDzp$, $\SDzpp$ et~\cDz{}
les ensembles $\RDCS{\A^1_k}{S}$, $\RDC{\A^1_k}$, $\TDCp{\A^1_k}$,
$\TDCpp{\A^1_k}$, $\SDCp{\A^1_k}$, $\SDCpp{\A^1_k}$ et la \condDC{\A^1_k},
étant entendu que
$\RDCS{\A^1_k}{S}$, $\RDC{\A^1_k}$, $\TDCp{\A^1_k}$,
$\TDCpp{\A^1_k}$, $\SDCp{\A^1_k}$, $\SDCpp{\A^1_k}$ et la \condDC{\A^1_k} désignent
les ensembles et la condition obtenus en appliquant les définitions ci-dessus
avec $C=\A^1_k$ après
avoir restreint~$\pi$ au-dessus de $\A^1_k \subset \P^1_k$.

\bigskip
Le théorème principal de ce chapitre est le suivant.

\bigskip
\begin{theoreme}Admettons
\label{ch2thprin}
l'hypothèse de Schinzel. Supposons que $C=\P^1_k$ et que la fibre de~$\pi$
au-dessus du point $\infty \in \P^1(k)$ soit lisse.
Il existe alors un ensemble fini $S_0 \subset \Omega$ et un sous-groupe fini
$B_0 \subset \Br(U)$ tels que l'assertion suivante soit vérifiée.
Soient un ensemble $S_1 \subset \Omega$
fini contenant~$S_0$ et une famille $(x_v)_{v\in S_1}\in \prod_{v \in S_1} U(k_v)$.
Supposons que $X_{x_v}(k_v)\neq \emptyset$ pour tout $v \in S_1 \cap \Omega_f$ et que
$$
\sum_{v \in S_1} \inv_v A(x_v)=0
$$
pour tout $A \in B_0$. Supposons aussi que
pour toute place $v \in S_1$ réelle, on ait $X_\infty(k_v)\neq\emptyset$
et $x_v$ appartienne à la composante
connexe non majorée de~$U(k_v)$.  Alors
\begin{enumerate}
\item[a)] si $\Mrond' \neq \emptyset$ et $\Mrond''\neq\emptyset$,
il existe un élément de~$\RD$ arbitrairement
proche de~$x_v$ en chaque place $v \in S_1 \cap \Omega_f$ et arbitrairement
grand en chaque place archimédienne de~$k$;
\item[b)] il existe un élément de~$\RDSz{S_1}$ arbitrairement
proche de~$x_v$ en chaque place $v \in S_1 \cap \Omega_f$,
arbitrairement grand en chaque place archimédienne de~$k$ et entier hors de~$S_1$.
\end{enumerate}
\end{theoreme}

\bigskip
Voici les conséquences du théorème~\ref{ch2thprin} pour l'existence
et la Zariski-densité des points rationnels de~$X$.  Pour le restant de
ce paragraphe, supposons que $C=\P^1_k$ et que les ensembles $\Mrond'$ et~$\Mrond''$ ne soient
pas vides.

\bigskip
\begin{theoreme}
\label{ch2thrdbrvert}
Admettons l'hypothèse de Schinzel.  Alors l'adhérence de~$\RD$ dans~$C(\A_k)$
est égale à~$\pi(X(\A_k)^{\Br_\vert})$.
\end{theoreme}

\bigskip
\begin{demo}
Ce théorème se démontre à partir du théorème~\ref{ch2thprin} exactement comme
le théorème~\ref{ch1appobm1} se démontre à partir du théorème~\ref{ch1thprin}.
\end{demo}

\bigskip
\begin{theoreme}
\label{ch2thobmx}
Admettons l'hypothèse de Schinzel et la finitude des groupes de Tate-Shafarevich des
courbes elliptiques $\Erond'_x$ pour $x \in U(k)$. Supposons que la \condD{}
soit vérifiée et que $X(\A_k)^{\Br_\vert}\neq\emptyset$.
Alors $X(k)\neq\emptyset$.
Si de plus~$\pi$ ne possède pas de section, l'ensemble
$X(k)$ est Zariski-dense dans~$X$.
\end{theoreme}

\bigskip
\begin{demo}
Supposons que $X(\A_k)^{\Br_\vert}\neq\emptyset$. Le théorème~\ref{ch2thrdbrvert}
montre que l'ensemble~$\RD$ est infini.  Soit $x \in \RD$.
Compte tenu de l'équivalence $[\Xrond'']=0 \Leftrightarrow [\Xrond] \in \TDp$,
la \condD{} entraîne que
l'un des groupes $\TDp$ et $\TDpp$ est nul et que l'autre est d'ordre au plus~$2$.
D'après les carrés commutatifs
$$
\myxyhook\xymatrix{
\Erond''_U(U)/\Im(\phi'_U) \ar[d] \ar[r] & \SDp \ar[d] \\
\Erond''_x(k)/\Im(\phi'_x) \ar[r] & H^1(k, \tors{\phi'_x}\Erond'_x)
}
\quad\!\! \myxyhook\xymatrix@R=2.5ex{ \\ \text{et}} \quad\!\!
\myxyhook\xymatrix{
\Erond'_U(U)/\Im(\phi''_U) \ar[d] \ar[r] & \SDpp \ar[d] \\
\Erond'_x(k)/\Im(\phi''_x) \ar[r] & H^1(k, \tors{\phi''_x}\Erond''_x)
}
$$
et la définition de~$\RD$, l'un des groupes $\tors{\phi'_x}\Sha(k, \Erond'_x)$
et $\tors{\phi''_x}\Sha(k, \Erond''_x)$ est donc nul et l'autre est d'ordre au plus~$2$.
La suite exacte
$$
\myxyhook\xymatrix{
0 \ar[r] & \tors{\phi'_x}\Sha(k, \Erond'_x) \ar[r] & \tors{2}\Sha(k, \Erond'_x) \ar[r] & \tors{\phi''_x}\Sha(k,
\Erond''_x)
}
$$
permet d'en déduire que le groupe $\tors{2}\Sha(k,\Erond'_x)$
est d'ordre au plus~$2$.
La finitude\index{groupe de Tate-Shafarevich!finitude} de $\Sha(k,\Erond'_x)$ et les propriétés de l'accouplement de Cassels-Tate\index{accouplement de Cassels-Tate}
entraînent par ailleurs que cet ordre est un carré (cf.~preuve
du théorème~\ref{ch1appobm2}); par conséquent
$\tors{2}\Sha(k,\Erond'_x)=0$ et donc $X_x(k) \neq \emptyset$ puisque
$x \in \Rrond_A$.

Supposons de plus que~$\pi$ ne possède pas de section.
On vient de voir que
$\tors{\phi'_x}\Sha(k,\Erond'_x)=0$
(puisque $\tors{2}\Sha(k,\Erond'_x)=0$)
et $\dim_{\ff{\F}2}\tors{\phi''_x}\Sha(k,\Erond''_x)\leq 1$; d'où
$\dim_{\ff{\F}2}\tors{2}\Sha(k,\Erond''_x)\leq 1$, compte tenu de la suite exacte
$$
\myxyhook\xymatrix{
0 \ar[r] & \tors{\phi''_x}\Sha(k,\Erond''_x) \ar[r] & \tors{2}\Sha(k,\Erond''_x) \ar[r] &
\tors{\phi'_x}\Sha(k,\Erond'_x)\rlap{\text{.}}
}
$$
La suite exacte
$$
\myxyhook\xymatrix{
0 \ar[r] & \tors{\phi''_x}\Sha(k,\Erond''_x) \ar[r] & \Sha(k,\Erond''_x) \ar[r] &
\Sha(k,\Erond'_x)
}
$$
et la finitude de $\Sha(k,\Erond'_x)$
montrent
que le groupe $\Sha(k,\Erond''_x)$ est fini. Utilisant à nouveau les propriétés
de l'accouplement\index{accouplement de Cassels-Tate} de Cassels-Tate,
on obtient finalement que $\tors{2}\Sha(k,\Erond''_x)=0$ et
donc que $\tors{\phi''_x}\Sha(k,\Erond''_x)=0$.

Nous allons maintenant établir que le rang de la courbe elliptique $\Erond'_x$
n'est pas nul.  Comme les ensembles $X_x(k)$ et~$\Erond'_x(k)$ sont en bijection, cela
impliquera que $X_x(k)$ est infini puis que $X(k)$ est Zariski-dense
dans~$X$, étant donné que~$x$ a été choisi quelconque dans l'ensemble
infini~$\RD$.  Supposons que $\Erond'_x$ soit de rang nul et aboutissons à une
contradiction.  La nullité de $\tors{\phi''_x}\Sha(k,\Erond''_x)$
et du rang de~$\Erond'_x$ entraînent que
le groupe $\Sel_{\phi''_x}(k,\Erond''_x)$ est d'ordre~$2$.
Notons respectivement $G' \subset \SDp$ et $G'' \subset \SDpp$
les sous-groupes images réciproques de $\{0,[\Xrond]\}$ et de
$\{0,[\Xrond'']\}$ par les flèches $\SDp \rightarrow H^1(C,\Erond')$ et
$\SDpp \rightarrow H^1(C,\Erond'')$ issues des suites exactes~(\ref{ch2sesgtsgp})
et~(\ref{ch2sesgtsgpp}).
Par définition de~$\RD$, les flèches $G' \rightarrow \Sel_{\phi'_x}(k,\Erond'_x)$
et $G'' \rightarrow \Sel_{\phi''_x}(k,\Erond''_x)$ d'évaluation en~$x$
sont injectives. Le groupe~$G''$ est donc d'ordre au plus~$2$, d'où il résulte
que~$[\Xrond'']=0$.  On a alors $[\Xrond]\in\TDp$; comme~$\pi$ n'admet pas de section,
on en déduit que le groupe~$G'$ est d'ordre au moins~$4$. Il en va donc de même
pour $\Sel_{\phi'_x}(k,\Erond'_x)$.  Par ailleurs, les courbes elliptiques~$\Erond'_x$
et~$\Erond''_x$ ont même rang puisqu'elles sont isogènes.  La nullité
de $\tors{\phi'_x}\Sha(k,\Erond'_x)$ et du rang de~$\Erond'_x$ implique donc
que $\Sel_{\phi'_x}(k,\Erond'_x)$ est d'ordre au plus~$2$, d'où une contradiction.
\end{demo}

\bigskip
Le théorème suivant est celui dont nous nous servirons au troisième chapitre.
C'est un énoncé «~sur mesure~», ce qui explique sa forme quelque peu particulière.

\bigskip
\begin{theoreme}
\label{ch2thsurmesure}
Admettons l'hypothèse de Schinzel et la finitude des groupes de Tate-Shafarevich des
courbes elliptiques~$\Erond'_x$ pour $x \in U(k)$.  Supposons que la fibre de~$\pi$ au-dessus
du point $\infty \in \P^1(k)$ soit lisse, qu'elle possède un $k_v$\nobreakdash-point
pour toute place $v\in \Omega$ réelle et qu'il existe $x_0 \in U(k)$
appartenant à la composante connexe non majorée de~$U(k_v)$ pour toute place~$v$
réelle et
tel que $X_{x_0}(\A_k)\neq\emptyset$.  Soit $S\subset \Omega$ la réunion de l'ensemble des places archimédiennes
ou dyadiques de~$k$, de l'ensemble des places finies de mauvaise réduction pour
la courbe elliptique~$\Erond'_{x_0}$ et de l'ensemble des places finies~$v$ telles
que l'adhérence de $x_0\in \P^1(k_v)$ dans~$\P^1_{\Orond_v}$ rencontre celle
de $(\Mrond\cup\{\infty\}) \otimes_k k_v$. Supposons la \condE{} relative à~$S$ satisfaite:
\begin{itemize}
\medskip
\item[]\CondE{}\glossary{$(\mathrm{E})$}: l'image\index{condition (E)@\condE{}} de $\SDzp \cap \SG_{\phi',S}(\A^1_k, \Erond')$
dans $\tors{\phi'}H^1(\A^1_k,\Erond')$ par la seconde flèche de la suite exacte~(\ref{ch2sesgtsgp})
est incluse dans $\{0,[\Xrond]\}$
et l'image de \mbox{$\SDzpp \cap \SG_{\phi'',S}(\A^1_k, \Erond'')$}
dans $\tors{\phi''}H^1(\A^1_k,\Erond'')$
par la seconde flèche de la suite exacte~(\ref{ch2sesgtsgpp}) est incluse dans $\{0,[\Xrond'']\}$.
\end{itemize}
\medskip

\noindent{}Alors $X(k)\neq\emptyset$.
\end{theoreme}

\bigskip
\begin{demo}
Soient $S_0 \subset \Omega$ et $B_0 \subset \Br(U)$ les ensembles donnés par le théorème~\ref{ch2thprin}.
Soit $S_1 \subset \Omega$ fini contenant $S \cup S_0$ et contenant toutes les places $v \in \Omega$
pour lesquelles il existe $A \in B_0$ tel que $\inv_v A(x_0)\neq 0$.  Posons $x_v=x_0$ pour
tout $v \in S_1$.
On a alors $\sum_{v \in S_1} \inv_v A(x_v)=0$ pour tout $A \in B_0$
d'après la loi de réciprocité\index{loi de réciprocité globale} globale.
La conclusion du théorème~\ref{ch2thprin} permet d'en déduire l'existence de $x_1 \in \RDSz{S_1}$
arbitrairement proche de~$x_0$ en chaque place $v \in S_1 \cap \Omega_f$.
Pour $a \in \SG_{\phi',S_1}(\A^1_k,\Erond')$
(resp.~$a \in \SG_{\phi'',S_1}(\A^1_k,\Erond'')$)
et $z \in U(k)$, notons~$a(z)$ l'image de~$a$ par la composée
$$
\SG_{\phi',S_1}(\A^1_k,\Erond') \subset H^1(\A^1_k\setminus \Mrond', \tors{\phi'}\Erond')
\rightarrow H^1(k,\tors{\phi'_z}\Erond'_z)=k^\star/k^{\star 2}
$$
(resp.~$\SG_{\phi'',S_1}(\A^1_k,\Erond'') \subset H^1(\A^1_k\setminus \Mrond'', \tors{\phi''}\Erond'')
\rightarrow H^1(k,\tors{\phi''_z}\Erond''_z)=k^\star/k^{\star 2}$),
où l'inclusion est donnée par la proposition~\ref{ch2propcompct} et la seconde flèche
est l'évaluation en~$z$.
Quitte à choisir~$x_1$ suffisamment proche de~$x_0$ aux places finies de~$S_1$, on peut
supposer que la courbe elliptique~$\Erond'_{x_1}$ a bonne réduction
aux places de $S_1 \setminus S$ et que $v(a(x_1))=v(a(x_0))$
pour tout $v \in S_1 \setminus S$ et tout $a \in \SG_{\phi',S_1}(\A^1_k,\Erond')$
(resp.~tout $a \in \SG_{\phi'',S_1}(\A^1_k,\Erond'')$),
où~$v$ désigne l'application $k^\star/k^{\star 2} \rightarrow \Z/2$ induite par la valuation
normalisée associée à~$v$.

Prouvons maintenant que $x_1 \in \RDSz{S}$.  Étant donné que $x_1 \in \RDSz{S_1}$, il suffit pour cela
que tout $a \in \SDzp \cap \SG_{\phi',S_1}(\A^1_k,\Erond')$ tel que
$a(x_1) \in \Sel_{\phi'_{x_1}}(k,\Erond'_{x_1})$ appartienne à $\SG_{\phi',S}(\A^1_k,\Erond')$
(resp.~que tout $a \in \SDzpp \cap \SG_{\phi'',S_1}(\A^1_k,\Erond'')$ tel que
$a(x_1) \in \Sel_{\phi''_{x_1}}(k,\Erond''_{x_1})$ appartienne à $\SG_{\phi'',S}(\A^1_k,\Erond'')$).
Fixons un tel~$a$ et une place $v \in S_1 \setminus S$.
Il suffit de vérifier
que $v(a(x_0))=0$, puisque
l'adhérence de~$x_0$ dans~$\P^1_{\Orond_v}$
ne rencontre pas celle de $(\Mrond\cup\{\infty\}) \otimes_k k_v$.
Comme les courbes elliptiques~$\Erond'_{x_1}$
et~$\Erond''_{x_1}$ ont bonne réduction en~$v$ et que~$v$ n'est pas dyadique,
l'appartenance de~$a(x_1)$ au
groupe de Selmer entraîne que $v(a(x_1))=0$, d'où~$v(a(x_0))=0$.

De l'appartenance de~$x_1$ à~$\RDSz{S}$, de
la \condE{}
et de l'équivalence $[\Xrond'']=0 \Leftrightarrow [\Xrond] \in \tors{\phi'}H^1(\A^1_k,\Erond')$
résulte que l'un des deux groupes $\tors{\phi'_{x_1}}\Sha(k,\Erond'_{x_1})$
et $\tors{\phi''_{x_1}}\Sha(k,\Erond''_{x_1})$ est nul et que l'autre est d'ordre
au plus~$2$. La suite exacte
$$
\xymatrix{
0 \ar[r] & \tors{\phi'_{x_1}}\Sha(k,\Erond'_{x_1}) \ar[r] & \tors{2}\Sha(k,\Erond'_{x_1})
\ar[r] & \tors{\phi''_{x_1}}\Sha(k,\Erond''_{x_1})
}
$$
et la finitude du groupe $\Sha(k,\Erond'_{x_1})$ permettent de conclure,
exactement comme dans la preuve du théorème~\ref{ch2thobmx}.
\end{demo}

\bigskip
Comparons maintenant le théorème~\ref{ch2thobmx} avec les théorèmes~A et~B de~\cite{ctbsd}.

Sous les hypothèses du théorème~A, prenons pour~$X$ la surface notée~$X(m)$
dans~\cite[Th.~A]{ctbsd}.
La courbe~$X_\eta$ est naturellement un $2$\nobreakdash-revêtement du quotient de la courbe elliptique
d'équation
\begin{equation}
\label{ch2compbsdeq}
y^2=(x-c(t))(x^2-d(t))
\end{equation}
par le point de coordonnées $(c(t),0)$, d'où un choix canonique de $P' \in
E'_\eta(K)$ tel que la courbe elliptique $E''_\eta$ soit celle définie par
l'équation~(\ref{ch2compbsdeq}) et que le point $P'' \in E''_\eta(K)$ soit
celui de coordonnées $(c(t),0)$.  Les types de réduction de~$E'_\eta$
et~$E''_\eta$ sont les suivants, dans la notation de
Kodaira (cf.~\cite[(1.2.2)]{ctbsd} pour une équation de Weierstrass de~$E'_\eta$):

\begin{center}
\begin{tabular}{|c|c|c|}
\cline{2-3} \multicolumn{1}{c|}{} &
\vbox to .95em{}
\lower.1em\hbox to 13ex{\hfill$d=0$\hfill} &
\lower.1em\hbox to 13ex{\hfill$c^2-d=0$\hfill} \\
\hline \vbox to 2.8ex{}\raise.3ex\hbox{$E'_\eta$} &
\raise.3ex\hbox{$I_2$} & \raise.3ex\hbox{$I_1$} \\
\hline \vbox to 2.8ex{}\raise.3ex\hbox{$E''_\eta$} &
\raise.3ex\hbox{$I_1$} & \raise.3ex\hbox{$I_2$} \\
\hline
\end{tabular}
\end{center}

\smallskip
\noindent{}Par conséquent les ensembles~$\Mrond'$ et~$\Mrond''$ définis ci-dessus
coïncident respectivement avec les ensembles~$\Mrond''$ et~$\Mrond'$ de~\cite[p.~120]{ctbsd};
de plus $L'_M=L''_M=\kappa(M)$ pour tout $M\in\Mrond$, et l'on a
$\ff{F'}M=\Z/2$ pour tout $M \in \Mrond''$ et $\ff{F''}M=\Z/2$ pour tout $M \in \Mrond'$.
La proposition~\ref{ch2propcompct} montre que
les groupes~$\mathfrak{S}'$, $\mathfrak{S}''$, $\mathfrak{S}'_0$ et~$\mathfrak{S}''_0$
de~\cite[p.~120]{ctbsd} s'identifient respectivement aux groupes $\SG_{\phi''}(\A^1_k,\Erond'')$,
$\SG_{\phi'}(\A^1_k,\Erond')$, $\SG_{\phi''}(\P^1_k,\Erond'')$ et~$\SG_{\phi'}(\P^1_k,\Erond')$.
La commutativité du diagramme
\begin{equation}
\label{ch2compdiag1}
\begin{aligned}
\myxyhook\xymatrix{
\Z/2 \ar[d] \ar@{=}[r] & \tors{\phi'}\Erond' \ar[r] \ar[d] & \Erond' \ar[d] \\
\displaystyle \prod_{M \in \Mrond''} i_{M\star}\Z/2 \ar@{=}[r] &
\displaystyle \prod_{M \in \Mrond''} i_{M\star}\left(\tors{2}{}\ff{F'}M\right) \ar[r]^\sim &
\displaystyle \prod_{M \in \Mrond''} i_{M\star}\left(\ff{F'}M\right)
}
\end{aligned}
\end{equation}
entraîne celle du diagramme
\begin{equation}
\label{ch2compdiag2}
\begin{aligned}
\myxyin\myxyhook\xymatrix{
*+!!<0pt,\the\fontdimen22\textfont2>!<-4.3ex,0ex>{\SG_{\phi'}(\P^1_k, \Erond')}
\ar@{=}[r] & H^1(\P^1_k \setminus \Mrond', \Z/2)
\ar[r] \ar[d] &
H^1(\P^1_k \setminus \Mrond', \Erond') \ar[d] \\ &
\displaystyle \bigoplus_{M\in\Mrond''}H^1(\kappa(M), \Z/2) \ar[r]^\sim &
\displaystyle \bigoplus_{M\in\Mrond''}H^1(\kappa(M),\ff{F'}M)\rlap{\text{,}}
}\myxyout
\end{aligned}
\end{equation}
qui à son tour montre que la composée du morphisme
$$H^1(\P^1_k,\Erond')\xrightarrow{\;\;\prod \delta'_M\;}\prod_{M\in \Mrond''}H^1(L'_M,\ff{F'}M)$$
et de la flèche de droite de la suite exacte~(\ref{ch2sesgtsgp})
s'identifie à l'application notée~$\delta''_0$ dans~\cite[p.~120]{ctbsd}.
De même pour $\prod \delta''_M$ et~$\delta'_0$.
Par ailleurs, la classe $[\Xrond'']$ est nulle.

On déduit de ces considérations que les hypothèses (3.a) et (3.b) du théorème~A
de~\cite{ctbsd} impliquent la \condD{}. Le théorème~A est donc un cas particulier du
théorème~\ref{ch2thobmx}.
(En toute rigueur, pour que le théorème~\ref{ch2thobmx} implique formellement
le théorème~A de~\cite{ctbsd}, il aurait fallu détailler quelque peu sa conclusion.
Cependant, on vérifie tout de suite que les assertions supplémentaires
du théorème~A découlent directement de la preuve ci-dessus du théorème~\ref{ch2thobmx}.
Nul besoin pour cela de revenir à la preuve du théorème~\ref{ch2thprin}.)

Passons maintenant au théorème~B.  Prenons pour~$X$ la surface notée~$X$
dans~\cite[Th.~B]{ctbsd}. La courbe $X_\eta$ est un $2$\nobreakdash-revêtement de la courbe
elliptique définie par l'équation de Weierstrass~(\ref{ch2compbsdeq}), où~$c$ et~$d$
ont cette fois les significations données dans~\cite{ctbsd} juste avant l'énoncé du théorème~B.
Notons $P' \in E'_\eta(K)$ le point de coordonnées $(c(t),0)$. Les types de réduction
de~$E'_\eta$ et~$E''_\eta$ sont les suivants, dans la notation
de Kodaira (cf.~\cite[(1.2.2)]{ctbsd} pour une équation de Weierstrass de~$E''_\eta$):

\begin{center}
\begin{tabular}{|c|c|c|c|c|}
\cline{2-5} \multicolumn{1}{c|}{} &
\vbox to 1.12em{}
\raise.07em\hbox{\;$d_{01}=0$\;} &
\raise.07em\hbox{\;$d_{23}^2+4d_{24}d_{34}=0$\;} &
\raise.07em\hbox{\;$d_{04}^2-d_{02}d_{03}=0$\;} &
\raise.07em\hbox{\;$d_{14}^2-d_{12}d_{13}=0$\;} \\
\hline\vbox to 2.8ex{}
\raise.3ex\hbox{$E'_\eta$} &
\raise.3ex\hbox{$I_2$} &
\raise.3ex\hbox{$I_1$} &
\raise.3ex\hbox{$I_2$} &
\raise.3ex\hbox{$I_2$} \\
\hline\vbox to 2.8ex{}
\raise.3ex\hbox{$E''_\eta$} &
\raise.3ex\hbox{$I_4$} &
\raise.3ex\hbox{$I_2$} &
\raise.3ex\hbox{$I_1$} &
\raise.3ex\hbox{$I_1$} \\
\hline
\end{tabular}
\end{center}

\smallskip
\noindent{}Les ensembles~$\Mrond'$ et~$\Mrond''$ s'identifient donc
respectivement aux ensembles~$\Mrond'$ et~$\Mrond''$ du théorème~B, et l'on a
$\ff{F'}M=\Z/2$ pour tout $M \in \Mrond''$ et~$\ff{F''}M=\Z/2$ pour tout~$M
\in\Mrond'_2$.
La proposition~\ref{ch2propcompct} montre que
les groupes~$\mathfrak{S}'$, $\mathfrak{S}''$, $\mathfrak{S}'_0$ et~$\mathfrak{S}''_0$
de~\cite[p.~123]{ctbsd} s'identifient respectivement aux groupes $\SG_{\phi'}(\A^1_k,\Erond')$,
$\SG_{\phi''}(\A^1_k,\Erond'')$, $\SG_{\phi'}(\P^1_k,\Erond')$ et~$\SG_{\phi''}(\P^1_k,\Erond'')$.
Les extensions $L'_M/\kappa(M)$ sont toutes triviales, de même que
les extensions $L''_M/\kappa(M)$ pour $M \in \Mrond'_2 \cup \Mrond''$. Pour $M \in \Mrond'_1$,
l'extension $L''_M/\kappa(M)$ est celle obtenue en adjoignant à~$\kappa(M)$ une racine
carrée de l'élément noté $\delta''_M(-c)$ dans~\cite{ctbsd}, d'où une injection canonique
$$
\kappa(M)^\star/\langle\kappa(M)^{\star 2},\delta''_M(-c)\rangle
\longhookrightarrow L''^\star_M/L''^{\star 2}_M
= H^1(L''_M, \tors{2}{\ff{F''}M}) \longhookrightarrow H^1(L''_M, \ff{F''}M)\rlap{\text{.}}
$$
(La dernière flèche est injective car le $L''_M$\nobreakdash-groupe $\ff{F''}M
\otimes_{\kappa(M)} L''_M$ est
constant.) On vérifie de plus que $\delta'_M([\Xrond])$ (resp.~$\delta''_M([\Xrond''])$)
coïncide avec la classe notée~$\epsilon_M$ (resp.~$\delta''_M(d_{14}^2-d_{12}d_{13})$)
dans~\cite{ctbsd}, pour $M \in \Mrond''$ (resp.~pour $M \in \Mrond'_1$). Enfin,
le groupe~$\ff{F'}M$ étant trivial pour $M \in \Mrond'_2$, on a $\delta'_M([\Xrond])=0$
et donc $\delta''_M([\Xrond''])=0$ pour tout $M \in \Mrond'_2$.

Des diagrammes commutatifs analogues à~(\ref{ch2compdiag1}) et~(\ref{ch2compdiag2}) permettent
maintenant de voir que les hypothèses (3.a) et (3.b) du théorème~B de~\cite{ctbsd} impliquent
la \condD{}.  Le théorème~B est donc lui aussi un cas particulier du théorème~\ref{ch2thobmx}
(avec le même commentaire que pour le théorème~A).

\bigskip
Pour résumer, les théorèmes~A et~B de~\cite{ctbsd}, pris ensemble, correspondent au cas
particulier du théorème~\ref{ch2thobmx} où les hypothèses
supplémentaires suivantes sont supposées satisfaites:
\begin{itemize}
\item pour tout $M \in \Mrond$, le groupe $\ff{F'}M$ est d'ordre au plus~$2$;
\item les courbes elliptiques $E'_\eta$ et~$E''_\eta$ possèdent un seul point d'ordre~$2$
rationnel;
\item elles sont de rang de Mordell-Weil nul sur~$K$;
\item le morphisme~$\pi$ n'admet pas de section;
\item les hypothèses techniques (0.1), (0.5) et~(0.7) de~\cite{ctbsd} sont satisfaites.
\end{itemize}
Le théorème~A s'applique seulement lorsque $[\Xrond'']=0$ et le théorème~B
seulement lorsque $[\Xrond'']\neq0$. (En effet, si $[\Xrond'']=0$, la condition~(3.a)
du théorème~B entraîne que $[\Xrond]=0$, autrement dit que~$\pi$ possède une section.)

\section{Preuve du théorème~\ref{ch2thprin}}

Pour prouver le théorème~\ref{ch2thprin}, on peut évidemment supposer
la fibration $\pi \colon X \rightarrow \P^1_k$ relativement minimale.  On fera
dorénavant cette hypothèse, dont l'intérêt principal est qu'elle permet d'écrire
que $U=\A^1_k \setminus \Mrond$.
Les paragraphes~\ref{ch2premierparpreuve} à~\ref{ch2dernierparavantpreuve} contiennent
des définitions et résultats préliminaires à la preuve proprement dite du théorème~\ref{ch2thprin},
qui occupe le paragraphe~\ref{ch2parpreuve}. Le symbole~$S$ désigne pour le moment
un ensemble fini arbitraire de places de~$k$ contenant les places dyadiques et les places
archimédiennes; il sera précisé au paragraphe~\ref{ch2parpreuve}.

\subsection{Couples admissibles, préadmissibles}
\label{ch2premierparpreuve}
\label{ch2couplesadm}

Si $x \in \P^1_k$ est un point fermé, on note~$\xtilde$ son adhérence
schématique dans~$\P^1_{\Orond}$.  C'est un $\Orond$\nobreakdash-schéma fini
génériquement étale.  On définit de même l'adhérence~$\xtilde$ de~$x$
dans~$\P^1_{\Orond_v}$ lorsque~$x$ est un point fermé de~$\P^1_{k_v}$
et que $v \in \Omega_f$.  Si $\xtilde \cap \P^1_{\Orond_S}$ est
étale sur $\Spec(\Orond_S)$,
l'ensemble des points fermés de $\xtilde \cap \P^1_{\Orond_S}$
s'identifie à l'ensemble des places de~$\kappa(x)$ dont la trace
sur~$k$ n'appartient pas à~$S$.  On utilisera librement cette
identification par la suite.

Soit $\Trond$ une famille $(\ft{T'}M)_{M \in \Mrond}$, où $\ft{T'}M$ est un ensemble fini de
places finies de~$\kappa(M)$.
Notons $\ft{T}M \subset \Omega_f$ l'ensemble des traces sur $k$ des
places de $\ft{T'}M$.
On dit que la famille~$\Trond$
est \emph{préadmissible}\index{préadmissible!famille} si la condition suivante est satisfaite:
\begin{itemize}
\medskip
\item[(\refstepcounter{equation}\theequation{})]
les sous-ensembles $\ft{T}M \subset \Omega$ pour $M \in \Mrond$
sont deux à deux disjoints et disjoints de~$S$; le $\Orond_S$\nobreakdash-schéma
$\P^1_{\Orond_S} \cap \bigcup_{M \in \Mrond \cup \{\infty\}} \Mtilde$
est étale;
pour tout $M \in \Mrond$, l'application trace $\ft{T'}M \rightarrow \ft{T}M$ est bijective.
\end{itemize}
\medskip

Soit $x \in U(k)$. On dit que le couple $(\Trond, x)$ est
\emph{préadmissible}\index{préadmissible!couple} si la famille~$\Trond$ l'est et si de plus la
condition suivante est satisfaite:
\begin{itemize}
\medskip
\item[(\refstepcounter{equation}\theequation{})] $x$ est entier hors de~$S$ (\emph{i.e.}
$\xtilde \cap \widetilde{\infty} \cap \P^1_{\Orond_S}=\emptyset$);
pour tout $M \in \Mrond$, le schéma $\xtilde \cap \Mtilde \cap \P^1_{\Orond_S}$ est réduit, et son
ensemble sous-jacent est la réunion de $\ft{T'}M$ et
d'une place de $\kappa(M)$ hors de $\ft{T'}M$, que l'on note $w_M$.
\end{itemize}

\bigskip
Étant donné un tel couple $(\Trond, x)$, on notera $v_M$ la trace de $w_M$
sur $k$, et l'on posera\glossary{$v_M$, $w_M$, $\ft{T}M$, $T(x)$}
$$
T = S \cup \bigcup_{M \in \Mrond} \ft{T}M
$$
et $$T(x) = T \cup \ensemble{v_M}{M \in \Mrond}\rlap{\text{.}}$$

\bigskip
Les conditions de préadmissibilité sont de nature géométrique.  Nous voudrons aussi imposer
une condition arithmétique sur les triplets considérés afin d'assurer l'existence de
points locaux sur~$X_x$ (cf.~condition~(\ref{ch2conditionadmiss}) ci-dessous, et proposition~\ref{ch2exisadel}).

\bigskip
\begin{lemme}
\label{ch2existekmd}
Soit $M \in \Mrond$.
Pour tout $d \in \tors{2}H^1(\kappa(M), \ff{F'}M)$, il existe une extension quadratique ou triviale
minimale $K_{M,d}$ de~$\kappa(M)$ telle que~$d$ appartienne au noyau de
la restriction $H^1(\kappa(M), \ff{F'}M) \rightarrow H^1(L'_M K_{M,d}, \ff{F'}M)$
et que $L'_M K_{M,d}$ se plonge $L'_M$\nobreakdash-linéairement dans toute extension
$\ell/L'_M$ telle que l'image de~$d$ dans~$H^1(\ell,\ff{F'}M)$
soit nulle.
La même assertion avec~$\ff{F''}M$ et~$L''_M$ à la place de~$\ff{F'}M$ et~$L'_M$
est vraie.
\end{lemme}

\bigskip
\begin{demo}
Si le groupe $\ff{F'}M(L'_M)$ (resp.~$\ff{F''}M(L''_M)$) est d'ordre impair,
nécessairement $d=0$ et l'on peut donc prendre $K_{M,d}=\kappa(M)$.
Sinon, la démonstration est exactement la même que pour
le lemme~\ref{ch1existekmd}.
\end{demo}

\bigskip
On notera par la suite~$K_M$\glossary{$K_{M,d}$, $K_M$} le corps $K_{M,d}$ donné par le
lemme~\ref{ch2existekmd} en prenant pour~$d$ l'image de~$[\Xrond]$
dans $H^1(\kappa(M),\ff{F'}M)$.

\bigskip
La famille~$\Trond$ sera dite \emph{admissible}\index{admissible!famille} si elle est préadmissible
et que la condition suivante est vérifiée:
\begin{itemize}
\medskip
\item[(\refstepcounter{equation}\theequation{})\label{ch2conditionadmiss}]
pour tout $M \in \Mrond$, les places de~$\ft{T'}M$ sont totalement décomposées dans $L'_M L''_M K_M$.
\end{itemize}
\medskip

Si $x \in U(k)$, on dit enfin que le couple $(\Trond, x)$ est
\emph{admissible}\index{admissible!couple} s'il est préadmissible, si la famille~$\Trond$ est admissible et
si pour tout $M \in \Mrond$, la place $w_M$ est totalement décomposée dans $L'_M L''_M K_M$.

\subsection{Prélude à l'étude des groupes de Selmer en famille}
\label{ch2prelude}

Dans ce paragraphe, nous fixons un couple préadmissible $(\Trond,x)$ et nous faisons
l'hypothèse suivante:
\begin{itemize}
\medskip
\item[(\refstepcounter{equation}\theequation{})\label{ch2hypprelude}]
pour tout $M \in \Mrond$, la classe
du diviseur $\Mtilde \cap \A^1_{\Orond_S}$ dans $\Pic(\A^1_{\Orond_S})$ est nulle.
\end{itemize}

\vskip2ex
Cette hypothèse sera en tout cas satisfaite si~$S$ est suffisamment grand.

Posons $\Urond = \P^1_{\Orond} \setminus \left( \bigcup_{M \in \Mrond \cup
\{\infty\}} \Mtilde \right)$\glossary{$\Urond$}.
Par définition de la préadmissibilité,
le sous-schéma localement fermé
$\xtilde \cap \P^1_{\Orond_{T(x)}}$ de~$\P^1_{\Orond_T}$
est inclus dans l'ouvert $\Urond_{\Orond_T} = \Urond\otimes_\Orond\Orond_T = \Urond \cap
\P^1_{\Orond_T}$.

\bigskip
\begin{proposition}
\label{ch2propeviso}
Le morphisme $H^1(\Urond_{\Orond_T}, \Z/2) \rightarrow H^1(\Orond_{T(x)},\Z/2)$ induit par l'inclusion
de $\xtilde \cap \P^1_{\Orond_{T(x)}}$ dans $\Urond_{\Orond_T}$ est un isomorphisme.
\end{proposition}

\bigskip
\begin{demo}
Notons respectivement $\alpha \colon \Z^\Mrond \rightarrow \Pic(\A^1_{\Orond_T})$
et $\beta \colon \Z^\Mrond \rightarrow \Pic(\Orond_T)$ les applications $\Z$\nobreakdash-linéaires envoyant
$M\in \Mrond$ sur la classe du diviseur $\Mtilde \cap \A^1_{\Orond_T}$
dans $\Pic(\A^1_{\Orond_T})$
et sur la classe
de~$v_M$ dans $\Pic(\Orond_T)$.
Les inclusions de $\xtilde \cap \P^1_{\Orond_T}$ dans $\A^1_{\Orond_T}$
et de $\xtilde \cap \P^1_{\Orond_{T(x)}}$ dans~$\Urond_{\Orond_T}$
induisent le diagramme commutatif suivant:
$$
\myxyhook\xymatrix{
\Z^\Mrond \ar@{=}[d] \ar[r]^(0.375)\alpha & \Pic(\A^1_{\Orond_T}) \ar[d] \ar[r] & \Pic(\Urond_{\Orond_T}) \ar[r]
\ar[d]^\delta & 0 \\
\Z^\Mrond \ar[r]^(0.375)\beta & \Pic(\Orond_T) \ar[r]^(.46)\gamma & \Pic(\Orond_{T(x)}) \ar[r] & 0 \rlap{\text{.}}
}
$$
Le morphisme
$\Pic(\Orond_T) \rightarrow \Pic(\A^1_{\Orond_T})$ induit par
le morphisme structural du $\Orond_T$\nobreakdash-schéma $\A^1_{\Orond_T}$ est un isomorphisme.  La flèche verticale du milieu
du diagramme ci-dessus en est une rétraction; c'est donc aussi un isomorphisme. La commutativité du diagramme
et l'exactitude de ses lignes entraînent maintenant la bijectivité de~$\delta$.
Par ailleurs, l'hypothèse~(\ref{ch2hypprelude}) implique que $\alpha=0$, d'où l'on tire que $\beta=0$, ce qui signifie
encore que~$\gamma$ est bijective.

Considérons à présent le diagramme commutatif suivant:
$$
\myxyhook\xymatrix{
& 0 \ar[d] & 0 \ar[d] \\
0 \ar[r] & \Gm(\Orond_T)/2 \ar[d] \ar[r] & \Gm(\Orond_{T(x)})/2 \ar[d] \ar[r] & (\Z/2)^\Mrond \ar@{=}[d] \\
0 \ar[r] & H^1(\Orond_T,\Z/2) \ar[d] \ar[r] & H^1(\Orond_{T(x)},\Z/2) \ar[d] \ar[r] & (\Z/2)^\Mrond \\
& \tors{2}\Pic(\Orond_T) \ar[r]^\gamma & \tors{2}\Pic(\Orond_{T(x)}) \rlap{\text{.}}
}
$$
Les colonnes sont induites par la suite exacte de Kummer, la seconde ligne par la suite
spectrale de Leray associée à l'inclusion de $\Spec(\Orond_{T(x)})$ dans $\Spec(\Orond_T)$ et au
faisceau étale~$\Z/2$, et enfin la flèche de droite de la première ligne est induite par les valuations
normalisées de~$k$ associées aux places~$v_M$. La seconde ligne est exacte puisqu'elle provient
de la suite spectrale de Leray.  Comme~$\gamma$ est un isomorphisme, on en déduit que la première ligne
est exacte aussi.

Plaçons maintenant la première ligne du diagramme ci-dessus dans le diagramme
$$
\myxyhook\xymatrix{
0 \ar[r] & \Gm(\Orond_T)/2 \ar@{=}[d] \ar[r] & \Gm(\Urond_{\Orond_T})/2 \ar[d]^\epsilon \ar[r] &
(\Z/2)^\Mrond \ar@{=}[d] \ar[r] & 0 \\
0 \ar[r] & \Gm(\Orond_T)/2 \ar[r] & \Gm(\Orond_{T(x)})/2 \ar[r] & (\Z/2)^\Mrond \rlap{\text{,}}
}
$$
où~$\epsilon$ est l'évaluation en~$x$ et la flèche de droite de
la première ligne est induite par les valuations normalisées
de~$\kappa(\P^1_k)$ associées aux points $M \in \Mrond$.  Il résulte de l'hypothèse de
préadmissibilité sur le couple $(\Trond, x)$ que ce diagramme commute.
On a vu que la seconde ligne est exacte. La première l'est d'après l'hypothèse~(\ref{ch2hypprelude}).
Le lemme des cinq permet d'en déduire que la flèche~$\epsilon$ est bijective.

Insérons finalement la flèche de l'énoncé dans le diagramme commutatif
$$
\myxyhook\xymatrix{
0 \ar[r] & \Gm(\Urond_{\Orond_T})/2 \ar[d]^\epsilon \ar[r] & H^1(\Urond_{\Orond_T},\Z/2) \ar[d] \ar[r] &
\tors{2}\Pic(\Urond_{\Orond_T}) \ar[r] \ar[d]^\gamma & 0 \\
0 \ar[r] & \Gm(\Orond_{T(x)})/2 \ar[r] & H^1(\Orond_{T(x)},\Z/2) \ar[r] & \tors{2}\Pic(\Orond_{T(x)}) \ar[r]
& 0 \rlap{\text{,}}
}
$$
dont les lignes sont exactes.  On a vu que~$\gamma$ et~$\epsilon$ sont des isomorphismes.  Il n'y a plus qu'à
appliquer le lemme des cinq pour terminer la preuve.
\end{demo}

\bigskip
\begin{proposition}
\label{ch2h1u}
Pour tout ouvert dense~$V$ de~$\Spec(\Orond)$ sur lequel~$2$ est inversible,
le groupe $H^1(\Urond_V,\Z/2)$ s'identifie
au sous-groupe de~$K^\star/K^{\star 2}$ constitué des classes de fonctions
inversibles sur~$U$ dont la valuation sur chaque fibre de la
projection $\Urond_V\rightarrow V$ est paire.
\end{proposition}

\bigskip
\begin{demo}
La suite spectrale de Leray pour l'inclusion $j \colon U \rightarrow \Urond_V$ fournit
une suite exacte
\begin{equation}
\label{ch2h1use1}
\myxyhook\xymatrix{
0 \ar[r] & H^1(\Urond_V,\Z/2)\ar[r]&H^1(U,\Z/2)\ar[r]&H^0(\Urond_V,R^1j_\star\Z/2)\rlap{\text{.}}
}
\end{equation}
Compte tenu de la suite exacte de Kummer
et de la nullité de $R^1 j_\star \Gm$ (théorème de Hilbert~90),
le faisceau $R^1j_\star \Z/2$ est égal au conoyau
de la multiplication par~$2$ sur~$j_\star \Gm$. 
Grâce à la régularité de~$\Urond_V$,
on dispose par ailleurs de la suite exacte
des diviseurs de Weil
\begin{equation}
\label{ch2h1use2}
\myxyhook\xymatrix{
0 \ar[r] & \Gm \ar[r] & j_\star \Gm \ar[r] & \displaystyle\bigoplus_{v \in V^{(1)}} i_{v\star}\Z
\ar[r] & 0 \rlap{\text{,}}
}
\end{equation}
où $i_v \colon \Urond_v \rightarrow \Urond_V$ désigne l'immersion fermée canonique.
Celle-ci montre que le conoyau de la multiplication par~$2$ sur~$j_\star\Gm$
s'identifie naturellement
à $\bigoplus_{v\in V^{(1)}} i_{v\star}\Z/2$.
La suite exacte~(\ref{ch2h1use1}) se récrit donc comme suit:
\begin{equation}
\label{ch2h1use3}
\myxyhook\xymatrix{
0 \ar[r] & H^1(\Urond_V,\Z/2) \ar[r] & H^1(U,\Z/2) \ar[r] & \displaystyle\bigoplus_{v \in V^{(1)}}
\Z/2\rlap{\text{.}}
}
\end{equation}
On a $H^1(U,\Z/2)=\Gm(U)/2$ puisque~$U$ est un ouvert de~$\A^1_k$, et l'on vérifie tout
de suite que la flèche $\Gm(U)/2 \rightarrow \bigoplus_{v\in V^{(1)}}\Z/2$ donnée par
la suite~(\ref{ch2h1use3}) envoie la classe d'une fonction inversible sur la famille
de ses valuations modulo~$2$
aux points de codimension~$1$ de~$\Urond_V$ qui ne dominent pas~$V$.
La proposition est donc prouvée.
\end{demo}

\bigskip
Le lemme suivant ne présente aucune difficulté.  Nous l'énonçons en vue de son
utilisation répétée au paragraphe~\ref{ch2parpreuve}.

\bigskip
\begin{lemme}
\label{ch2transvcomm}
Pour~$v \in \Omega_f$, notons
$v \colon K^\star \rightarrow \Z$ la valuation normalisée associée au point de
codimension~$1$ de~$\P^1_\Orond$ défini par~$v$, et pour $M \in \Mrond \cup \{\infty\}$, notons
$\mathbf{v}_M \colon K^\star \rightarrow \Z$ la valuation normalisée associée au point~$M$.
Pour tout $M \in \Mrond\cup\{\infty\}$, tout $x \in U(k)$ et toute place $v \in
\Omega_f$ au-dessus de laquelle~$\xtilde$ rencontre transversalement~$\Mtilde$,
le diagramme
$$
\myxyin\myxyhook\xymatrix{
H^1(U,\Z/2) \ar[d] \ar[r] & H^1(K,Z/2)\ar@{=}[r] &
*+!!<0pt,\the\fontdimen22\textfont2>[r]!<5.7ex,0ex>{K^\star/K^{\star 2}}\ar[r]^(.43){v+\mathbf{v}_M}&\Z/2
\ar@{=}[d] \\
H^1(k,\Z/2) \ar@{=}[r] &
*+!!<0pt,\the\fontdimen22\textfont2>[r]!<7.7ex,0ex>{k^\star/k^{\star 2}}
\ar[rr]^(.39){v}&&\Z/2\rlap{\text{,}}
}\myxyout
$$
dans lequel la flèche verticale de gauche est l'évaluation en~$x$, est commutatif.
\end{lemme}

\bigskip
\begin{demo}
On a $H^1(U,\Z/2)=\Gm(U)/2$ puisque~$U$ est un ouvert de~$\A^1_k$.
Il suffit donc de démontrer que pour tout $f \in \Gm(U)$,
considérant~$f$ comme une fonction rationnelle sur~$\P^1_\Orond$
et~$f(x)$ comme une fonction rationnelle sur~$\xtilde$,
l'égalité $v(f)+\mathbf{v}_M(f)=v(f(x))$ vaut.
C'est évident si~$f$ est constante.  Par ailleurs, si~$v(f)=0$, autrement dit si~$f$
définit une fonction inversible sur~$\Urond_R$, où~$R$ désigne le localisé de~$\Orond$
en l'idéal premier défini par~$v$, cette égalité résulte immédiatement
de l'hypothèse de transversalité.
Le groupe~$\Gm(U)$ étant engendré par les sous-groupes~$\Gm(k)$ et~$\Gm(\Urond_R)$,
le lemme est donc prouvé.
\end{demo}

\subsection{Dualité locale pour les courbes $\Erond'_x$, $\Erond''_x$}
\label{ch2dualitelocale}

Les quelques notations suivantes nous seront utiles par la suite.
Pour $v \in \Omega$ et $x \in U(k_v)$, posons\glossary{$\fv{V}v$, $\fv{V'}v$,
$\fv{V''}v$, $W'_v$, $W''_v$}
$$
\fv{V}v=H^1(k_v,\Z/2) \quad;\quad
\fv{V'}v=H^1(k_v,\tors{\phi'_x}\Erond'_x) \quad;\quad
\fv{V''}v=H^1(k_v,\tors{\phi''_x}\Erond''_x)
$$
et
$$
W'_v=\Erond''_x(k_v)/\Im(\phi'_x) \quad;\quad
W''_v=\Erond'_x(k_v)/\Im(\phi''_x)\rlap{\text{.}}
$$

\medskip
\noindent{}Si $v \in \Omega_f$, posons de plus
$\ft{T}v=\Ker(\fv{V}v\rightarrow H^1(k_v^\nr,\Z/2))$,
$\ft{T'}v=\Ker(\fv{V'}v\rightarrow H^1(k_v^\nr,\tors{\phi'_x}\Erond'_x))$
et $\ft{T''}v=\Ker(\fv{V''}v\rightarrow H^1(k_v^\nr,\tors{\phi''_x}\Erond''_x))$,
où $k_v^\nr$ désigne une extension non ramifiée maximale de~$k_v$.
On a canoniquement $\fv{V}v=\fv{V'}v=\fv{V''}v$ et $\ft{T}v=\ft{T'}v=\ft{T''}v$.
Les suites exactes
$$
\myxyhook\xymatrix{
0 \ar[r] & \Z/2 \ar[r] & \Erond'_x \ar[r]^{\phi'_x} & \Erond''_x \ar[r] & 0
}
$$
et
$$
\myxyhook\xymatrix{
0 \ar[r] & \Z/2 \ar[r] & \Erond''_x \ar[r]^{\phi''_x} & \Erond'_x \ar[r] & 0
}
$$
montrent que~$W'_v$ (resp.~$W''_v$) s'identifie canoniquement à un sous-groupe
de~$\fv{V'}v$ (resp.~de~$\fv{V''}v$). L'accouplement de Weil
$\tors{\phi'_x}\Erond'_x \times \tors{\phi''_x}\Erond''_x \rightarrow \Z/2$
et l'injection canonique $H^2(k_v,\Z/2) \hookrightarrow \Z/2$ donnée par la théorie
du corps de classes local fournissent
un accouplement $\fv{V'}v \times \fv{V''}v \rightarrow \Z/2$,
non dégénéré d'après \cite[Ch.~I, Cor.~2.3 et Th.~2.13]{milneadt}.
Les sous-espaces~$W'_v$ et~$W''_v$ sont orthogonaux pour cet accouplement
(cf.~\cite[Lemma~3]{bsd}).

\bigskip
\begin{proposition}
\label{ch2ekvsurj}
Il existe un ensemble fini $S_0 \subset \Omega$ tel que pour toute
place finie $v \in \Omega \setminus S_0$,
tout $M \in \Mrond$ et tout $x \in U(k_v)$ tel que~$\xtilde$
rencontre transversalement~$\Mtilde\otimes_\Orond\Orond_v$
en une place~$w$ de~$\kappa(M)$ totalement décomposée dans~$L'_ML''_M$, les conditions
suivantes soient satisfaites:
\begin{enumerate}
\item Si $M \in \Mrond'$, les applications $H^1(k_v^\nr, \Erond'_x)\rightarrow
H^1(k_v^\nr,\Erond''_x)$ et $H^1(k_v,\Erond'_x)\rightarrow H^1(k_v,\Erond''_x)$
induites par~$\phi'_x$ sont injectives.
\item Si $M \in \Mrond''$, les applications
$\Erond'_x(k_v^\nr)\rightarrow\Erond''_x(k_v^\nr)$ et $\Erond'_x(k_v)\rightarrow
\Erond''_x(k_v)$ induites par~$\phi'_x$ sont surjectives.
\end{enumerate}
Le même énoncé reste vrai si l'on échange tous les~$'$ et les~$''$.
On a par conséquent $W'_v=\fv{V'}v$ et~$W''_v=0$ si $M \in \Mrond'$,
$W'_v=0$ et~$W''_v=\fv{V''}v$ si $M \in \Mrond''$.
\end{proposition}

\bigskip
\begin{demo}
La première des conditions de l'énoncé est équivalente
à la seconde une fois les~$'$ et les~$''$ échangés, d'après le théorème de
dualité locale pour les variétés abéliennes (cf.~\cite[Cor.~3.4]{milneadt}); le dual de
Pontrjagin du morphisme induit par~$\phi'_x$ est en effet le morphisme induit par~$\phi''_x$,
puisque~$\phi''_x$ est l'isogénie duale de~$\phi'_x$.
(En toute rigueur, le théorème de dualité locale s'applique lorsque $K=k_v$ mais pas lorsque~$K=k_v^\nr$;
cependant, l'injectivité de $H^1(k_v^\nr,\Erond'_x)\rightarrow H^1(k_v^\nr,\Erond''_x)$
se déduit facilement de l'injectivité de
$H^1(k_v,\Erond'_x)\rightarrow H^1(k_v,\Erond''_x)$.)
On peut donc se contenter d'établir d'une part l'énoncé $(P')$
obtenu en supprimant les deux dernières phrases de la proposition
et en remplaçant «~tout $M \in \Mrond$~» par «~tout $M \in \Mrond''$~»,
d'autre part l'énoncé $(P'')$ obtenu à partir de~$(P')$ en échangeant les~$'$ et les~$''$.

Nous prouvons ci-dessous $(P')$.  Le lecteur constatera que l'on obtient
une preuve de~$(P'')$ en échangeant simplement les~$'$ et les~$''$.

Prenons pour~$S_0$ un ensemble fini contenant les places archimédiennes et les
places dyadiques, suffisamment grand pour que
$\P^1_{\Orond_{S_0}}\cap\bigcup_{M\in\Mrond\cup\{\infty\}}\Mtilde$ soit étale
sur~$\Orond_{S_0}$.
D'autres conditions sur~$S_0$ seront précisées ci-dessous.
Soient $v \in \Omega \setminus S_0$, $M \in \Mrond''$ et $x \in U(k_v)$ tel que~$\xtilde$
rencontre transversalement $\Mtilde\otimes_\Orond\Orond_v$ en une place~$w$ de~$\kappa(M)$
totalement décomposée dans $L'_ML''_M$.
Notons respectivement $\Erondsouligne'_x$ et~$\Erondsouligne''_x$ les modèles
de Néron\index{modèle de Néron} de~$\Erond'_x$ et~$\Erond''_x$ au-dessus de $\xtilde=\Spec(\Orond_v)$.
Soient
$\Erondsouligne'^0_x \subset \Erondsouligne'_x$ et~$\Erondsouligne''^0_x
\subset \Erondsouligne''_x$ leurs composantes neutres
et $\souligne{F}'$ et~$\souligne{F}''$ les $\kappa(v)$\nobreakdash-groupes finis étales
fibres spéciales de $\Erondsouligne'_x/\Erondsouligne'^0_x$
et de $\Erondsouligne''_x/\Erondsouligne''^0_x$.

Supposons un instant que les courbes $\Erond'_x$ et~$\Erond''_x$ aient réduction
multiplicative, que le $k_v$\nobreakdash-point non nul du noyau de~$\phi'_x
\colon \Erond'_x \rightarrow \Erond''_x$ ne se spécialise pas sur la composante
neutre de la fibre spéciale de~$\Erondsouligne'_x$
et que les $\kappa(v)$\nobreakdash-groupes $\souligne{F}'$ et~$\souligne{F}''$ soient constants.
Le morphisme
$\souligne{F}'\rightarrow \souligne{F}''$ induit par~$\phi'_x$ s'insère alors dans
une suite exacte
\begin{equation}
\label{ch2ekvsurjse}
\myxyhook\xymatrix{
0 \ar[r] & \Z/2 \ar[r] & \souligne{F}' \ar[r] & \souligne{F}'' \ar[r] & 0
}
\end{equation}
(cf.~proposition~\ref{annpropppp}). Le lemme~\ref{annlemmeisogekmsh} suffirait à en déduire la
surjectivité de l'application $\Erond'_x(k_v^\nr) \rightarrow \Erond''_x(k_v^\nr)$
induite par~$\phi'_x$, mais nous allons de toute manière l'établir en même temps
que la surjectivité de $\Erond'_x(k_v) \rightarrow \Erond''_x(k_v)$.
On a un diagramme commutatif canonique
\begin{equation}
\label{ch2ekvsurjdiag}
\begin{aligned}
\myxyhook\xymatrix{
&& 0 \ar[d] & 0 \ar[d] \\
&& \Erondsouligne'^0_x \ar[d] \ar[r]^{\souligne{\phi}'^0_x}
& \Erondsouligne''^0_x \ar[r] \ar[d] & 0 \\
0 \ar[r] & \Z/2 \ar[d] \ar[r] & \Erondsouligne'_x \ar[d] \ar[r]^{\souligne{\phi}'_x}
& \Erondsouligne''_x \ar[d] \ar[r] & 0 \\
0 \ar[r] & i_\star \Z/2 \ar[r]\ar[d]
& i_\star \souligne{F}' \ar[d] \ar[r] & i_\star \souligne{F}''
\ar[d] \ar[r] & 0 \\
&0& 0 & 0
}
\end{aligned}
\end{equation}
de faisceaux étales sur~$\Spec(\Orond_v)$,
où~$i$ désigne l'inclusion du point fermé de~$\Spec(\Orond_v)$. Ses lignes et ses colonnes
sont exactes (cf.~lemme~\ref{annlemmesurje0}
pour la surjectivité de $\souligne{\phi}'^0_x$;
la surjectivité de~$\souligne{\phi}'_x$ résulte de l'exactitude
du reste du diagramme).
Compte tenu de la propriété universelle des modèles de Néron, l'exactitude de la
seconde ligne entraîne la surjectivité de l'application
$\Erond'_x(k_v^\nr)\rightarrow\Erond''_x(k_v^\nr)$ induite par~$\phi'_x$.

Passant aux sections globales dans~(\ref{ch2ekvsurjdiag}), on obtient
le diagramme commutatif
$$
\xymatrix{
\Erond'_x(k_v) \ar[d] \ar[r] & \Erond''_x(k_v) \ar[d] \ar[r] & H^1(\Orond_v,\Z/2) \ar[d] \\
\souligne{F}'(\kappa(v)) \ar[r] & \souligne{F}''(\kappa(v)) \ar[r] & H^1(\kappa(v),\Z/2)
\rlap{\text{,}}
}
$$
dont les lignes sont exactes.  La flèche verticale de droite
est injective (lemme de Hensel).  D'autre part,
les $\kappa(v)$\nobreakdash-groupes $\souligne{F}'$ et~$\souligne{F}''$ étant constants,
la flèche horizontale inférieure de gauche est surjective
(cf.~suite exacte~(\ref{ch2ekvsurjse})). Il~s'ensuit que l'application
$\Erond'_x(k_v) \rightarrow \Erond''_x(k_v)$ induite par~$\phi'_x$ est surjective.

Il reste donc seulement à établir les propriétés de~$\Erond'_x$ et~$\Erond''_x$
admises ci-dessus, lorsque~$S_0$ est assez grand.  Considérons seulement la courbe~$\Erond'_x$;
les raisonnements qui suivent s'appliqueront aussi bien à $\Erond''_x$ si l'on
échange tous les~$'$ et les~$''$, à quelques détails évidents près.
La courbe elliptique $E'_\eta/K$ possède une équation de Weierstrass
minimale de la forme
\begin{equation}
\label{ch2eqweier}
Y^2=(X-c)(X^2-d)
\end{equation}
avec $c,d\in k[t]$ telle que le point~$P'$ ait pour coordonnées $(X,Y)=(c,0)$. Son
discriminant est $\Delta'=16d(c^2-d)^2$.  Quitte à choisir l'ensemble~$S_0$ suffisamment
grand au début de la preuve de cette proposition, on peut supposer que les coefficients de~$c$ et de~$d$
sont des $S_0$\nobreakdash-entiers et que les coefficients dominants de~$d$ et de~$c^2-d$
sont des unités aux places de $\Omega \setminus S_0$.

Grâce à la minimalité de l'équation de Weierstrass,
l'ouvert de lissité sur~$\A^1_k$ du sous-schéma fermé de $\P^2_k \times_k \A^1_k$
défini par~(\ref{ch2eqweier}) s'identifie à $\Erond'^0_{\A^1_k}$
(cf.~\cite[Ch.~IV, §9, Cor.~9.1]{silvaec2}).
Il suffit donc de lire l'équation~(\ref{ch2eqweier}) modulo~$M$ pour déterminer
si~$P'$ se spécialise sur la composante neutre de~$\Erond'_M$ ou non, et pour
déterminer le type de réduction (bonne, multiplicative déployée, multiplicative non
déployée ou additive) de~$E'_\eta$ en~$M$.
Compte tenu des hypothèses sur~$E'_\eta$ et sur~$M$,
du lemme~\ref{annlemmeextreme},
de la proposition~\ref{annpropppp} et du corollaire~\ref{anncorcorfaux},
on voit ainsi que les polynômes~$d$ et~$c^2-d$ sont premiers entre eux,
que $(c^2-d)(M)=0$, que~$P'$ ne se spécialise pas sur~$\Erond'^0_M$
et que soit~$M$ est racine simple de $c^2-d$,
soit l'extension $L'_ML''_M/\kappa(M)$ est isomorphe
à~$\kappa(M)(\sqrt{2c(M)})$.

La courbe elliptique $\Erond'_x/k_v$ a pour équation de Weierstrass
\begin{equation}
\label{ch2eqweierx}
Y^2=(X-c(x))(X^2-d(x)) \rlap{\text{.}}
\end{equation}
D'après les hypothèses faites sur~$S_0$, les coefficients de l'équation~(\ref{ch2eqweierx})
sont des entiers $v$\nobreakdash-adiques.  (En effet, d'une part les coefficients
de~$c$ et de~$d$ sont des entiers $v$\nobreakdash-adiques, d'autre part~$x$ lui-même
en est un.) De plus, on a $v((c^2-d)(x))>0$ et~$v(d(x))=0$, ce qui entraîne
que l'équation de Weierstrass~(\ref{ch2eqweierx}) est minimale et que la courbe
elliptique $\Erond'_x/k_v$ est à réduction multiplicative, de type~$I_n$
avec $n=2v((c^2-d)(x))$.
Par le même raisonnement que ci-dessus, la minimalité de~(\ref{ch2eqweierx})
permet de voir que le $k_v$\nobreakdash-point non nul du noyau de~$\phi'_x$ ne se spécialise pas
sur~$\Erondsouligne'^0_x$.
En particulier, si $v((c^2-d)(x))=1$, les courbes elliptiques~$\Erond'_x$
et~$\Erond''_x$ sont respectivement à réduction de type~$I_2$ et~$I_1$
(cf.~proposition~\ref{annpropppp});
pour montrer que les $\kappa(v)$\nobreakdash-groupes $\souligne{F}'$ et~$\souligne{F}''$
sont constants, il suffit donc de vérifier que soit $v((c^2-d)(x))=1$,
soit $\Erond'_x$ est à réduction multiplicative
déployée (cf.~corollaire~\ref{anncorcorfaux}).
La place~$w$ étant
par hypothèse totalement décomposée dans~$L'_M L''_M$, il résulte
du calcul de cette extension que
soit~$M$ est racine simple de $c^2-d$, soit l'image de~$2c(M)$ dans~$\kappa(M)_w$
est un carré.
Dans le premier cas, on a $v((c^2-d)(x))=1$ puisque~$\xtilde$ rencontre
$\Mtilde\otimes_\Orond \Orond_v$ transversalement.
Dans le second, comme~$2c(M)$ est une unité $w$\nobreakdash-adique et que son image
dans~$\kappa(w)$ est égale
à celle de~$2c(x)$ dans~$\kappa(v)$, la courbe elliptique~$\Erond'_x$ est à
réduction multiplicative déployée.
\end{demo}

\subsection{Réciprocités et existence de points locaux}
\label{ch2parrec}

Pour\index{loi de réciprocité globale|(} $M \in \Mrond$, notons $\theta_M \in \kappa(M)$ l'image de~$t$ par le
morphisme \mbox{$k[t] \rightarrow \kappa(M)$} déduit de l'inclusion de~$M$ dans $\A^1_k=\Spec(k[t])$.
Si~$E$ est une extension quadratique ou triviale de~$\kappa(M)$, on pose
$$
A_{E/\kappa(M)} = \Cores_{\kappa(M)(t)/k(t)}(E/\kappa(M), t-\theta_M) \in \tors{2}\Br(U) \rlap{\text{,}}
$$
où\glossary{$A_{E/\kappa(M)}$} $(\cdot,\cdot)$ désigne le symbole de Hilbert sur le corps $\kappa(M)(t)$.

\bigskip
\begin{proposition}
\label{ch2invtotdec}
Soient un couple préadmissible $(\Trond,x)$, une place \mbox{$v\in\Omega \setminus S$},
un point $M \in \Mrond$ et une extension quadratique ou triviale $E/\kappa(M)$ non ramifiée en toute place
de~$\kappa(M)$ divisant~$v$. Alors
on a $\inv_v A_{E/\kappa(M)}(x)=0$ si et seulement si l'une des deux conditions suivantes
est vérifiée:
\begin{enumerate}
\item La place~$v$ appartient à $\ft{T}M \cup \{v_M\}$ et l'unique place de $\ft{T'}M \cup \{w_M\}$ divisant~$v$
est totalement décomposée dans~$E$.
\item La place~$v$ n'appartient pas à $\ft{T}M \cup \{v_M\}$.
\end{enumerate}
\end{proposition}

\bigskip
\begin{demo}
Considérons l'égalité
\begin{equation}
\label{ch2invaeinvw}
\inv_v A_{E/\kappa(M)}(x)=\sum_{w|v}\inv_w (E/\kappa(M), x-\theta_M)\rlap{\text{,}}
\end{equation}
où la somme porte sur les places~$w$ de~$\kappa(M)$ divisant~$v$.
Si~$w$ est une telle place,
l'image dans $H^1(\kappa(M)_w,\Z/2)=\kappa(M)_w^\star/\kappa(M)_w^{\star 2}$
de la classe de l'extension quadratique ou triviale $E/\kappa(M)$
est la classe d'une unité $w$\nobreakdash-adique, puisque~$w$ est par hypothèse non ramifiée dans~$E$.
Comme~$v$ n'est pas une place dyadique, on en déduit d'une part que
l'invariant $\inv_w (E/\kappa(M),x-\theta_M)$ est nul si \mbox{$x-\theta_M \in \kappa(M)$} est une unité
$w$\nobreakdash-adique, d'autre part que si l'image de $x-\theta_M$ dans~$\kappa(M)_w$ est une uniformisante,
cet invariant est nul si et seulement si~$w$ est totalement décomposée dans~$E$
(cf.~\cite[Prop.~1.1.3]{cscrelle94}).
En particulier, compte tenu de~(\ref{ch2invaeinvw}) et de l'hypothèse de préadmissibilité,
on a $\inv_v A_{E/\kappa(M)}(x)=0$ si $v \not\in \ft{T}M\cup\{v_M\}$
et $\inv_v A_{E/\kappa(M)}(x)=\inv_w (E/\kappa(M),x-\theta_M)$ si
$v \in \ft{T}M\cup\{v_M\}$ et que~$w$ est l'unique place de $\ft{T'}M\cup\{w_M\}$
divisant~$v$, auquel cas l'image de $x-\theta_M$ dans~$\kappa(M)_w$ est une
uniformisante. La proposition s'ensuit.
\end{demo}

\bigskip
\begin{proposition}
\label{ch2recprop}
Pour tout $M\in \Mrond$, tout couple préadmissible $(\Trond,x)$ et toute extension quadratique ou triviale
$E/\kappa(M)$ non ramifiée en toute place de~$\kappa(M)$ dont la trace sur~$k$ n'appartient pas à~$T$, les
conditions suivantes sont équivalentes:
\begin{enumerate}
\item La place $w_M$ de~$\kappa(M)$ associée à $(\Trond,x)$ est totalement décomposée dans~$E$.
\item On a $\displaystyle \sum_{v \in T} \inv_v A_{E/\kappa(M)}(x)=0$.
\end{enumerate}
\end{proposition}

\bigskip
\begin{demo}
La proposition~\ref{ch2invtotdec} montre que $\inv_v A_{E/\kappa(M)}(x)=0$
pour toute place~$v\in\Omega$ n'appartenant pas à $T \cup \{v_M\}$.
On a donc
$$
\sum_{v \in T} \inv_v A_{E/\kappa(M)}(x) = \inv_{v_M}A_{E/\kappa(M)}(x) \rlap{\text{,}}
$$
en vertu de la loi de réciprocité globale.
Une nouvelle application de la proposition~\ref{ch2invtotdec} permet maintenant de conclure.
\end{demo}

\bigskip
\begin{proposition}
\label{ch2reckmm}
Il existe un ensemble fini $S_0\subset \Omega$ contenant les places archimédiennes
tel que pour tout $v \in \Omega\setminus S_0$,
les assertions suivantes soient vérifiées.
Soient $m \in \SG_{\phi'}(\A^1_k,\Erond')$ et $M \in \Mrond''$.
Notons $d \in H^1(\kappa(M),\ff{F'}M)$ l'image de~$m$ par la flèche composée
\begin{equation}
\label{ch2reckmmeq}
\SG_{\phi'}(\A^1_k,\Erond')\rightarrow H^1(\A^1_k,\Erond')
\rightarrow H^1(\kappa(M),\ff{F'}M)
\end{equation}
et~$K_{M,d}/\kappa(M)$ une extension quadratique ou triviale
satisfaisant aux conditions du lemme~\ref{ch2existekmd}.
Si l'image de~$m$
dans $H^1(U_{k_v},\Z/2)$ appartient au sous-groupe $H^1(\Urond_{\Orond_v},\Z/2)$,
alors:
\begin{enumerate}
\item Toute place de~$\kappa(M)$ divisant~$v$ est non ramifiée dans~$K_{M,d}$.
\item Pour tout $x \in \A^1(\Orond_v)$ rencontrant transversalement $\Mtilde \otimes_\Orond
\Orond_v$ en une place~$w$ de~$\kappa(M)$ totalement décomposée dans~$L'_ML''_M$,
l'image de~$m$ par la flèche composée
\begin{equation}
\SG_{\phi'}(\A^1_k,\Erond')\rightarrow
H^1(\A^1_k,\Erond')\rightarrow H^1(k_v,\Erond'_x)
\end{equation}
est nulle si et seulement si~$w$ est totalement décomposée dans~$K_{M,d}$.
\end{enumerate}
Le même énoncé reste vrai si l'on échange tous les~$'$ et les~$''$.
\end{proposition}

\bigskip
\begin{demo}
Comme $\Pic(\A^1_k\setminus \Mrond')=0$,
les groupes $\SG_{\phi'}(\A^1_k,\Erond')$ et
$\Gm(\A^1_k\setminus \Mrond')/2$ sont canoniquement isomorphes
(cf.~proposition~\ref{ch2propcompct}). On utilisera librement cette identification ci-dessous.
Choisissons un point \mbox{$a \in U(k)$} et notons $r \colon \Gm(\A^1_k\setminus \Mrond')/2
\rightarrow \Gm(k)/2$
l'application d'évaluation en~$a$.  Comme~$r$ est une rétraction de la suite exacte
$$
\myxyhook\xymatrix{
0 \ar[r] & \Gm(k)/2 \ar[r] & \Gm(\A^1_k \setminus \Mrond')/2
\ar[r] & \displaystyle\prod_{M\in\Mrond'} \Z/2\rlap{\text{,}}
}
$$
son noyau est fini.
Pour $M \in \Mrond''$, soit $G_M \subset H^1(\kappa(M),\ff{F'}M)$
l'image de $\Ker(r)$ par la flèche~(\ref{ch2reckmmeq}).
Les groupes~$G_M$ et l'ensemble~$\Mrond''$ étant finis, il existe
un ensemble fini $S_0 \subset\Omega$ contenant
les places archimédiennes, tel que pour tout $v \in \Omega\setminus S_0$,
tout $M\in\Mrond''$ et toute place~$w$ de~$\kappa(M)$ divisant~$v$, le noyau
de la flèche de restriction $H^1(\kappa(M),\ff{F'}M)\rightarrow H^1(\kappa(M)_w^\nr,\ff{F'}M)$
contienne~$G_M$ (cf.~\cite[§6.1, Proposition~21]{serrecg}).
Quitte à agrandir~$S_0$, on peut supposer de plus
que~$\atilde \cap \P^1_{\Orond_{S_0}}$ est inclus dans~$\Urond_{\Orond_{S_0}}$,
que~$S_0$ contient les places dyadiques de~$k$ et les places finies
qui sont ramifiées dans l'une des extensions~$L'_M$ ou~$L''_M$ pour $M \in \Mrond$,
que le $\Orond_{S_0}$\nobreakdash-schéma $\P^1_{\Orond_{S_0}} \cap \bigcup_{M \in \Mrond \cup \{\infty\}}\Mtilde$ est étale
et que~$S_0$ contient l'ensemble du même nom donné par la proposition~\ref{ch2ekvsurj}.

Soient $v$, $m$, $M$ et~$d$ comme dans l'énoncé. L'existence du carré commutatif
$$
\myxyhook\xymatrix{
H^1(\Orond_v,\Z/2) \ar[d] & H^1(\Urond_{\Orond_v},\Z/2) \ar[d] \ar[l] \\
H^1(k_v^\nr,\Z/2) & \ar[l] H^1(U_{k_v^\nr},\Z/2)\rlap{\text{,}}
}
$$
dont la flèche verticale de gauche est nulle et dont les flèches horizontales
sont les applications d'évaluation en~$a$, entraîne que $r(m) \in \Ker(\Gm(k)/2
\rightarrow \Gm(k_v^\nr)/2)$.  Par conséquent, l'image de~$m$ dans
\mbox{$H^1((\A^1_k \setminus \Mrond')\otimes_k k_v^\nr,\Z/2)$} appartient à l'image
de~$\Ker(r)$.  Comme $M \in \Mrond''$,
la proposition~\ref{annpropppp} montre que le diagramme commutatif
$$
\myxyhook\xymatrix{
\Z/2\ar@{=}[r] \ar[d] & \tors{\phi'}\Erond' \ar[d] \ar[r] & \Erond' \ar[d] \\
i_{M\star} \Z/2 \ar[r]^(0.45)\sim & i_{M\star}\left(\tors{2}{}\ff{F'}M\right) \ar[r] & i_{M\star}
\ff{F'}M
}
$$
de faisceaux étales sur~$\P^1_k$ commute. On en déduit la commutativité du diagramme
$$
\myxyhook\xymatrix{
\SG_{\phi'}(\A^1_k,\Erond') \ar[r] \ar[d] &
H^1((\A^1_k\setminus\Mrond')\otimes_k k_v^\nr,\Z/2) \ar[d] \\
H^1(\A^1_k,\Erond') \ar[d] & H^1(\kappa(M)\otimes_k k_v^\nr,\Z/2) \ar[d] \\
H^1(\kappa(M),\ff{F'}M) \ar[r] & H^1(\kappa(M) \otimes_k k_v^\nr, \ff{F'}M) \rlap{\text{,}}
}
$$
où la flèche inférieure de la colonne de droite est induite par l'inclusion $\Z/2=\tors{2}{}\ff{F'}M
\subset \ff{F'}M$. Comme l'image de~$m$ par la flèche horizontale supérieure appartient
à l'image de~$\Ker(r)$, ce diagramme montre que pour toute place~$w$ de~$\kappa(M)$ divisant~$v$,
l'image de~$d$ dans $H^1(\kappa(M)_w^\nr, \ff{F'}M)$ appartient à l'image
de~$G_M$, et est donc nulle.
L'extension $L'_M/\kappa(M)$ étant non ramifiée en~$w$, elle se plonge dans~$\kappa(M)_w^\nr$.
L'extension $K_{M,d}/\kappa(M)$ se plonge donc elle aussi dans~$\kappa(M)_w^\nr$, par
définition de~$K_{M,d}$; autrement dit, elle est non ramifiée en~$w$.
La propriété~1 est établie.

Soit $x \in \A^1(\Orond_v)$ rencontrant transversalement $\Mtilde\otimes_\Orond\Orond_v$
en une place~$w$ de~$\kappa(M)$ totalement décomposée dans~$L'_ML''_M$.
Comme $M \in \Mrond''$, on a une suite exacte
\begin{equation}
\label{ch2fpfpp}
\myxyhook\xymatrix{
0 \ar[r] & \Z/2 \ar[r] & \ff{F'}M \ar[r] & \ff{F''}M \ar[r] & 0
}
\end{equation}
de $\kappa(M)$\nobreakdash-groupes (cf.~proposition~\ref{annpropppp}).
Les groupes~$\ff{F'}M$ et~$\ff{F''}M$ deviennent constants
après extension des scalaires de~$\kappa(M)$ à~$\kappa(M)_w$
puisque~$w$ est totalement décomposée dans~$L'_ML''_M$.
La flèche $H^1(\kappa(M)_w,\Z/2)\rightarrow H^1(\kappa(M)_w,\ff{F'}M)$
issue de la suite exacte~(\ref{ch2fpfpp}) est donc injective.
Par ailleurs, la conclusion de la proposition~\ref{ch2ekvsurj} montre
que la flèche $H^1(k_v,\Z/2)\rightarrow H^1(k_v,\Erond'_x)$ induite par
la suite exacte
$$
\myxyhook\xymatrix{
0 \ar[r] & \Z/2 \ar[r] & \Erond'_x \ar[r]^{\phi'_x} & \Erond''_x \ar[r] & 0
}
$$
est injective. Ces deux injections s'inscrivent dans le diagramme commutatif
$$
\myxyin\myxyhook\xymatrix{
& \Gm(k_v)/2 \ar@{^{ (}->}[r] & H^1(k_v,\Erond'_x) \\
*+!!<0pt,\the\fontdimen22\textfont2>!<-3.3ex,0ex>{\SG_{\phi'}(\A^1_k,\Erond')}
\ar@{=}[r] & \Gm(\A^1_k \setminus \Mrond')/2 \ar[d] \ar[r]
\ar[u] & H^1(\A^1_k \setminus \Mrond', \Erond') \ar[d] \ar[u] \\
& \Gm(\kappa(M)_w)/2 \ar@{^{ (}->}[r] & H^1(\kappa(M)_w,\ff{F'}M) \rlap{\text{,}}
}\myxyout\phantom{\SG_{\phi'}(\A^1_k,\Erond')==}
$$
où les flèches verticales supérieures (resp.~inférieures) sont les flèches d'évaluation
en~$x$ (resp.~en~$M$).
Étant donné que la place~$w$ est totalement décomposée dans~$L'_M$,
il résulte de la définition de~$K_{M,d}$ que
l'image de~$m$ dans $H^1(\kappa(M)_w,\ff{F'}M)$ est nulle si et seulement
si~$w$ est totalement décomposée dans~$K_{M,d}$. Compte tenu du diagramme ci-dessus,
il suffit donc,
pour conclure, de prouver que l'image de~$m$ dans~$\Gm(k_v)/2$ par l'évaluation
en~$x$ est nulle si et seulement si l'image de~$m$ dans~$\Gm(\kappa(M)_w)/2$ par l'évaluation
en~$M$ est nulle.  Comme l'image de~$m$ dans $H^1(U_{k_v},\Z/2)$ appartient
à $H^1(\Urond_{\Orond_v},\Z/2)$, la classe de~$m$ dans~$\Gm(\A^1_k\setminus \Mrond')/2$
est représentée par une fonction rationnelle~$f$ sur~$\A^1_{\Orond_v}$
inversible sur le complémentaire de l'adhérence de~\mbox{$\Mrond'\otimes_k k_v$}.
Notant encore~$M$ le point fermé de~$\A^1_{k_v}$ associé au point~$M\in\Mrond$
et à la place~$w$, les éléments~$f(M)$ de~$\kappa(M)_w$ et~$f(x)$ de~$k_v$
sont respectivement des unités $w$\nobreakdash-adiques et $v$\nobreakdash-adiques.  Ils ont même
réduction modulo~$v$ et sont donc simultanément des carrés (lemme de Hensel).
La propriété~2 est ainsi établie.

L'énoncé obtenu en échangeant les~$'$ et les~$''$ se prouve en appliquant la
même opération à la démonstration qui\index{loi de réciprocité globale|)} précède.
\end{demo}

\subsection{Finitude de $\SDzp$ et $\SDzpp$}
\label{ch2dernierparavantpreuve}

L'énoncé suivant servira pour la preuve de l'assertion a) du théorème~\ref{ch2thprin}.

\bigskip
\begin{proposition}
\label{ch2propfinitude}
Si $\Mrond' \neq \emptyset$ et $\Mrond'' \neq \emptyset$, les groupes
$\SDzp$ et~$\SDzpp$ sont finis.
\end{proposition}

\bigskip
\begin{demo}
Par symétrie, il suffit de considérer le groupe~$\SDzp$.  Soit $M \in \Mrond''$.
La flèche horizontale inférieure du carré commutatif canonique
$$
\xymatrix{
H^1(\A^1_k \setminus \Mrond', \tors{\phi'}\Erond') \ar[r] \ar[d] &
H^1(\A^1_k\setminus \Mrond', \Erond') \ar[d] \\
H^1(L'_M,\ff{\tors{\phi'_M}F'}M) \ar[r] & H^1(L'_M,\ff{F'}M)
}
$$
est injective puisque le $L'_M$\nobreakdash-groupe $\ff{F'}M \otimes_{\kappa(M)} L'_M$ est constant.
La flèche verticale de gauche s'identifie
à l'évaluation $H^1(\A^1_k \setminus \Mrond', \Z/2)
\rightarrow H^1(L'_M,\Z/2)$.  Il suffit de démontrer que le noyau de cette flèche d'évaluation
est fini, compte tenu de la proposition~\ref{ch2propcompct}.
Le sous-groupe \mbox{$H^1(k,\Z/2) \subset H^1(\A^1_k \setminus \Mrond',\Z/2)$} étant d'indice
fini, il reste seulement à vérifier que le noyau de la flèche
de restriction $H^1(k,\Z/2) \rightarrow H^1(L'_M,\Z/2)$
est fini, ce qui est élémentaire.
\end{demo}

\subsection{Fin de la preuve}
\label{ch2parpreuve}

\begin{proposition}
\label{ch2exisadel}
Il existe un ensemble $S_0 \subset \Omega$ tel que pour tout $S\subset \Omega$ fini contenant~$S_0$
et tout couple admissible $(\Trond,x)$ vérifiant $X_x(k_v)\neq\emptyset$ pour tout $v \in S
\setminus \Tinfty$,
on ait $X_x(\A_k)\neq\emptyset$, en notant~$\Tinfty$ l'ensemble des places de $S\setminus S_0$
en lesquelles~$x$ n'est pas entier.
\end{proposition}

\bigskip
\begin{demo}
Même démonstration que pour la proposition~\ref{ch1propexistadel}, en remplaçant
$\Erond^0_M$, $\Erond_M$, $\ff{F}M$ et~$L_MK_M$ dans le lemme~\ref{ch1xmxrondm}
par $\Erond'^0_M$, $\Erond'_M$, $\ff{F'}M$ et~$L'_ML''_MK_M$.
\end{demo}

\bigskip
Nous pouvons maintenant commencer la démonstration du
théorème~\ref{ch2thprin}.  Soit $S_0 \subset \Omega$ un ensemble fini
contenant les places archimédiennes, les places divisant un nombre premier inférieur
ou égal au degré de~$\Mrond$ vu comme $k$\nobreakdash-schéma fini réduit, les ensembles~$S_0$ donnés par les
propositions~\ref{ch2ekvsurj}, \ref{ch2reckmm} et~\ref{ch2exisadel} et les places finies de~$k$ au-dessus desquelles le morphisme
$$
\bigcup_{M\in\Mrond\cup\{\infty\}}\Mtilde \longrightarrow \Spec(\Orond)
$$
n'est pas étale, assez
grand pour que les extensions $L'_M/k$, $L''_M/k$ et~$K_M/k$ soient non
ramifiées hors de~$S_0$ pour tout $M\in\Mrond$,
pour que la condition~(\ref{ch2hypprelude}) soit satisfaite avec $S=S_0$,
pour que les courbes elliptiques~$E'_\eta$ et~$E''_\eta$ sur~$K$ s'étendent
en des schémas abéliens au-dessus de
$\left(\P^1_\Orond \setminus \bigcup_{M\in\Mrond}\Mtilde\right)\otimes_\Orond\Orond_{S_0}$,
et pour que les sous-groupes de $H^1(U,\Z/2)$ images réciproques
de~$\{0,[\Xrond]\}$ et de~$\{0,[\Xrond'']\}$ par les flèches~(\ref{ch2injfp}) et~(\ref{ch2injfpp}),
finis d'après le théorème de Mordell-Weil généralisé,
soient inclus dans $H^1(\Urond_{\Orond_{S_0}},\Z/2)$.
Soit $B_0 \subset \Br(U)$ le sous-groupe fini engendré par les classes
$A_{K_M/\kappa(M)}$, $A_{L'_M/\kappa(M)}$ et $A_{L''_M/\kappa(M)}$ pour $M\in\Mrond$
(cf.~paragraphe~\ref{ch2parrec}).
Soient \mbox{$S_1 \subset \Omega$}
fini contenant~$S_0$ et $(x_v)_{v\in S_1} \in \prod_{v \in S_1}U(k_v)$ satisfaisant aux hypothèses du
théorème~\ref{ch2thprin}.  On va prouver l'existence de $x \in U(k)$ arbitrairement proche de~$x_v$ pour
$v \in S_1 \cap \Omega_f$ et arbitrairement grand aux places archimédiennes de~$k$, tel que $x \in \RD$ si
$\Mrond'\neq \emptyset$ et $\Mrond''\neq\emptyset$, ou tel que $x \in \RDSz{S_1}$ et que~$x$ soit entier
hors de~$S_1$.

Commençons par définir un ensemble fini $\Tinfty \subset \Omega$ disjoint de~$S_1$. Nous allons prouver simultanément
les deux conclusions du théorème, mais avec des choix différents pour l'ensemble~$\Tinfty$. Pour prouver l'assertion
b) du théorème, on pose $\Tinfty=\emptyset$. Pour prouver l'assertion a), on pose $\Tinfty=\{v_\infty\}$, où~$v_\infty$
est choisie comme suit.  Lorsque $\Mrond'\neq\emptyset$ et $\Mrond''\neq\emptyset$, les
trois sous-groupes $H^1(\Urond_{\Orond_{S_1}}, \Z/2)$, $\SDzp$ et~$\SDzpp$ de $K^\star/K^{\star 2}$ sont finis
(le premier d'après le théorème des unités de Dirichlet et la finitude du groupe de classes d'un corps de nombres,
les deux autres d'après la proposition~\ref{ch2propfinitude}).
L'intersection du sous-groupe de $K^\star/K^{\star 2}$ qu'ils engendrent avec $k^\star/k^{\star 2}$ est donc
elle aussi finie, et par conséquent incluse dans $\Orond_{S'}^\star/\Orond_{S'}^{\star 2}$ pour un ensemble
$S' \subset \Omega$ fini assez grand. Le théorème de
\smash{\v{C}ebotarev}
fournit une infinité de places $v \in \Omega$
telles que pour tout $M \in \Mrond$, toute place de~$\kappa(M)$ divisant~$v$ soit totalement décomposée
dans~$L'_M L''_M K_M$. On choisit pour~$v_\infty$ une telle place hors de $S_1 \cup S'$ et l'on note
$x_{v_\infty}\in k_{v_\infty}$ l'inverse d'une uniformisante.

Posons $S=S_1\cup \Tinfty$\glossary{$S$, $S_1$, $\Tinfty$}.  On a, par définition de~$S_0$:

\bigskip
\begin{lemme}
\label{ch2bonnered}
Pour tout couple préadmissible $(\Trond,x)$,
les courbes elliptiques~$\Erond'_x$ et~$\Erond''_x$ ont
bonne réduction hors de~$T(x)\setminus \Tinfty$.
\end{lemme}

\bigskip
Pour chaque place finie $v \in S$, fixons un
voisinage $v$\nobreakdash-adique arbitrairement petit $\Arond_v$ de $x_v\in\P^1(k_v)$,
suffisamment petit pour que tout élément de~$\Arond_v$ soit l'inverse d'une
uniformisante si $v \in \Tinfty$, pour que $\inv_v A(x)=\inv_v A(x_v)$ pour
tout $x \in \Arond_v$ et tout $A \in B_0$, et
pour que $X_x(k_v)\neq\emptyset$ pour tout \mbox{$x \in \Arond_v$} si \mbox{$v \in S_1$}.
Il est possible de satisfaire cette dernière condition grâce au théorème des
fonctions implicites et à l'hypothèse que
$X_{x_v}(k_v)\neq\emptyset$ pour \mbox{$v \in S_1 \cap \Omega_f$}.
Fixons de même un voisinage $v$\nobreakdash-adique arbitrairement petit~$\Arond_v$
de $\infty \in \P^1(k_v)$ pour chaque place complexe $v \in \Omega$,
et un ouvert connexe non majoré arbitrairement petit $\Arond_v \subset U(k_v)$
pour chaque place réelle $v \in \Omega$,
suffisamment petit pour que $X_x(k_v)\neq\emptyset$ pour tout $x \in \Arond_v$
et pour que
les sous-groupes $\Erond'_{x}(k_v)/\Im(\phi''_{x})$
et $\Erond''_{x}(k_v)/\Im(\phi'_{x})$ de $H^1(k_v,\Z/2)$
ne dépendent pas de $x \in \Arond_v$ (cf.~preuve du lemme~\ref{ch1wvconst}).

\bigskip
\begin{lemme}
\label{ch2suminvax}
Soit $x \in U(k)$ appartenant à~$\Arond_v$ pour tout $v \in S$.
Alors $\sum_{v \in S} \inv_v A(x)=0$ pour tout $A \in B_0$.
\end{lemme}

\bigskip
\begin{demo}
Même démonstration que pour le lemme~\ref{ch1suminvax}.
\end{demo}

\bigskip
Si $(\Trond,x)$ est un couple préadmissible, on notera
$$\psi \colon H^1(\Orond_{T(x)},\Z/2) \longrightarrow H^1(\Urond_{\Orond_T},\Z/2)$$\glossary{$\psi$}
l'isomorphisme inverse de l'évaluation en~$x$ (cf.~proposition~\ref{ch2propeviso}).

\bigskip
\begin{proposition}
\label{ch2psiselsg}
Soit $(\Trond,x)$ un couple admissible.
Alors le groupe de $\phi'_x$\nobreakdash-Selmer de~$\Erond'_x$ et
le groupe de $\phi''_x$\nobreakdash-Selmer de~$\Erond''_x$ sont inclus dans $H^1(\Orond_{T(x)},\Z/2)$, et
l'on a
$$
\psi(\Sel_{\phi'_x}(k,\Erond'_x)) \subset \SG_{\phi'}(\A^1_k,\Erond')
$$
et
$$
\psi(\Sel_{\phi''_x}(k,\Erond''_x)) \subset \SG_{\phi''}(\A^1_k,\Erond'')
\rlap{\text{.}}
$$
\end{proposition}

\begin{demo}
La première assertion est une conséquence du lemme~\ref{ch2bonnered}.
La proposition~\ref{ch2ekvsurj} montre que tout élément du groupe de $\phi'_x$\nobreakdash-Selmer
de~$\Erond'_x$ est de valuation nulle en~$v_M$ pour tout $M \in \Mrond''$.
Comme~$\psi$ est à valeurs dans $H^1(\Urond_{\Orond_T},\Z/2)$,
le lemme~\ref{ch2transvcomm} et la proposition~\ref{ch2propcompct}
permettent d'en déduire que
$\psi(\Sel_{\phi'_x}(k,\Erond'_x))\subset\SG_{\phi'}(\A^1_k,\Erond')$.
On prouve l'autre inclusion de la même manière.
\end{demo}

\bigskip
Notons $\Lrond$\glossary{$\Lrond$} l'ensemble des couples $(m',m'')\in\N^2$ tels qu'il existe un
couple admissible $(\Trond,x)$ avec $x \in \Arond_v$ pour tout $v\in S$, tel
que:
\begin{equation}
\label{ch2listepropdeb}
\dim_{\ff{\F}2}\Sel_{\phi'_x}(k,\Erond'_x) = m' \rlap{\text{,}}
\end{equation}
\begin{equation}
\dim_{\ff{\F}2}\Sel_{\phi''_x}(k,\Erond''_x) = m'' \rlap{\text{,}}
\end{equation}
\begin{equation}
\label{ch2lincp}
\psi(\Sel_{\phi'_x}(k,\Erond'_x))\subset\SG_{\phi',S}(\A^1_k,\Erond') \rlap{\text{,}}
\end{equation}
\begin{equation}
\label{ch2lincpp}
\psi(\Sel_{\phi''_x}(k,\Erond''_x))\subset\SG_{\phi'',S}(\A^1_k,\Erond'') \rlap{\text{.}}
\end{equation}

\vspace{1mm}\noindent{}Munissons~$\Lrond$ de l'ordre partiel
$(m',m'')\leq(n',n'') \Longleftrightarrow m'\leq n'\text{ et }m''\leq n''$.

\bigskip
\begin{proposition}
\label{ch2propschinzel}
Admettons l'hypothèse de\index{hypothèse de Schinzel} Schinzel.  Soit $\Trond=(\ft{T'}M)_{M\in\Mrond}$ une famille
préadmissible.  Supposons donné, pour
tout $M\in \Mrond$ et tout $v \in \ft{T}M$, un $x_v \in \A^1(\Orond_v)$ rencontrant
transversalement $\Mtilde$ en l'unique place de~$\ft{T'}M$ qui divise~$v$.  Alors il
existe $x \in U(k)$ arbitrairement proche de~$x_v$ pour $v \in T \cap \Omega_f$
et arbitrairement grand aux places archimédiennes de~$k$, tel que
le couple $(\Trond, x)$ soit préadmissible.
\end{proposition}

\bigskip
\begin{demo}
Même démonstration que pour la proposition~\ref{ch1propschinzel}, à ceci près que l'on
utilise l'hypothèse~(\ref{ch2hypprelude}) au lieu de la nullité de $\Pic(\Orond_S)$.
\end{demo}

\bigskip
L'ensemble~$\Lrond$ n'est pas vide.  En effet, la
proposition~\ref{ch2propschinzel}, appliquée à la famille~$\Trond$
définie par $\ft{T'}M=\emptyset$ pour tout $M \in \Mrond$,
assure l'existence d'un $x \in U(k)$
appartenant à~$\Arond_v$ pour tout $v \in S$, tel que le couple $(\Trond,x)$
soit préadmissible.  Le lemme~\ref{ch2suminvax} et la
proposition~\ref{ch2recprop}, appliquée aux trois extensions $K_M/\kappa(M)$,
$L'_M/\kappa(M)$ et $L''_M/\kappa(M)$, montrent que ce couple est même
admissible.  Par ailleurs, les inclusions~(\ref{ch2lincp})
et~(\ref{ch2lincpp}) sont automatiques lorsque les ensembles~$\ft{T}M$ sont vides
(proposition~\ref{ch2psiselsg}).

Il existe donc un élément $(m',m'')\in\Lrond$ minimal.  Soient
$\Trond=(\ft{T'}M)_{M\in\Mrond}$ et $x \in U(k)$ avec $x \in \Arond_v$ pour
tout $v\in S$, tels que le couple $(\Trond,x)$ soit admissible et que les
propriétés~(\ref{ch2listepropdeb}) à~(\ref{ch2lincpp}) soient vérifiées.  On a
$X_x(k_v)\neq\emptyset$ pour tout $v \in S$ par construction des
voisinages~$\Arond_v$.  La conclusion de la proposition~\ref{ch2exisadel}
permet d'en déduire que $x \in \RA$.

\bigskip
\begin{proposition}
\label{ch2pluslongueprop}
On a $x \in \RDSz{S}$.  De plus, si $\Tinfty \neq \emptyset$, on a $x \in \RD$.
\end{proposition}

\bigskip
La démonstration de la proposition~\ref{ch2pluslongueprop} va nous occuper jusqu'à la fin de ce paragraphe.

\bigskip
\begin{demo}
L'injectivité des restrictions des flèches~(\ref{ch2pourinjfp}) et~(\ref{ch2pourinjfpp})
aux images réciproques de~$\{0,[\Xrond]\}$
et de~$\{0,[\Xrond'']\}$ par les flèches~(\ref{ch2injfp}) et~(\ref{ch2injfpp})
est une conséquence immédiate de la proposition~\ref{ch2propeviso} et de la définition
de~$S_0$.

Soit $\alpha \in \Sel_{\phi'_x}(k,\Erond'_x)$. Supposons que $\psi(\alpha) \not\in \SDSzp{S}$.
Vu l'hypothèse~(\ref{ch2lincp}), on a alors \mbox{$\psi(\alpha)\not\in\SDzp$}.

\bigskip
\begin{lemme}
\label{ch2lemmepeuttuerp}
Pour tout $e \in \tors{\phi'}H^1(\A^1_k,\Erond') \setminus \TDzp$, il existe $M \in \Mrond''$
et une infinité de places $v \in \Omega \setminus S_0$
telles qu'il existe une place~$w$ de~$\kappa(M)$
totalement décomposée dans $L'_ML''_MK_M$, divisant~$v$, non ramifiée et de degré résiduel~$1$
sur~$v$, telle que pour tout $x_v \in \A^1(\Orond_v)$ rencontrant transversalement~$\Mtilde$
en~$w$, l'image de~$e$ par la flèche $H^1(\A^1_k,
\Erond') \rightarrow H^1(k_v, \Erond'_{x_v})$ d'évaluation en~$x_v$ soit non nulle.
\end{lemme}

\bigskip
\begin{demo}
Soit $e \in \tors{\phi'}H^1(\A^1_k,\Erond') \setminus \TDzp$. Fixons un antécédent
$m \in \SG_{\phi'}(\A^1_k,\Erond')$ de~$e$ par la flèche de droite de la suite exacte~(\ref{ch2sesgtsgp}).
Par définition de~$\TDzp$, il existe $M \in \Mrond''$ tel que $\delta_M'(e)$
n'appartienne pas à $\{0,\delta_M'([\Xrond])\}$.  Soient~$d$ l'image de~$e$ dans
$\tors{2}H^1(\kappa(M),\ff{F'}M)$ et~$K_{M,d}/\kappa(M)$ une extension quadratique
vérifiant les conditions du lemme~\ref{ch2existekmd}.
Étant donné que l'extension $L'_MK_M/L'_M$ est quadratique ou triviale et que le groupe
$\ff{F'}M(L'_MK_M)=\ff{F'}M(L'_M)$ est cyclique, la suite exacte d'inflation-restriction
montre que le noyau de la flèche de restriction $H^1(L'_M,\ff{F'}M)\rightarrow
H^1(L'_MK_M,\ff{F'}M)$ est d'ordre au plus~$2$.
Notons~$N$ ce noyau. Si $N\neq 0$, alors $L'_MK_M \neq L'_M$,
de sorte que la définition de~$K_M$ entraîne que $\delta_M'([\Xrond]) \in N \setminus\{0\}$;
le groupe~$N$
est donc dans tous les cas engendré par~$\delta_M'([\Xrond])$.
Par conséquent, l'image de$~e$ dans $H^1(L'_MK_M,\ff{F'}M)$ n'est pas nulle.
L'extension quadratique $K_{M,d}/\kappa(M)$ ne se plonge donc pas dans $L'_MK_M/\kappa(M)$.
On en déduit, à l'aide du théorème de \smash{\v{C}ebotarev},
l'existence d'une infinité de places finies
de $L'_MK_M$ non ramifiées et de degré résiduel~$1$ sur~$k$,
inertes dans $L'_M K_M K_{M,d}$ (cf.~\cite[Proposition~2.2]{harduke}).
Soient~$w$ la trace sur~$\kappa(M)$ d'une telle place
et~$v \in \Omega_f$ la trace sur~$k$ de~$w$.
On peut supposer que~$v$ n'appartient pas à l'ensemble~$S_0$ de la proposition~\ref{ch2reckmm}
et que l'image de~$m$ dans $H^1(U_{k_v},\Z/2)$ appartient à $H^1(\Urond_{\Orond_v},\Z/2)$,
quitte à choisir~$v$ hors d'un certain ensemble fini.

La place~$w$ est totalement décomposée dans $L'_M K_M$, inerte dans $K_{M,d}$,
non ramifiée et de degré résiduel~$1$ sur~$v$. Elle est de plus totalement décomposée
dans~$L''_M$ car $L''_M/\kappa(M)$ se plonge dans $L'_M/\kappa(M)$, comme le montre
la suite exacte
$$
\myxyhook\xymatrix{
0 \ar[r] & \ff{F''}M \ar[r] & \ff{F'}M \ar[r] & \Z/2 \ar[r] & 0
}
$$
(cf.~proposition~\ref{annpropppp}).
La seconde partie de la proposition~\ref{ch2reckmm} permet donc de conclure.
\end{demo}

\bigskip
\begin{remarque}
Il n'était pas nécessaire de faire appel à la proposition~\ref{ch2reckmm} pour prouver
le lemme~\ref{ch2lemmepeuttuerp}; on aurait tout aussi bien pu reprendre la même démonstration
que pour le lemme~\ref{ch1lemmepeuttuer}.
\end{remarque}

\bigskip
La proposition~\ref{ch2psiselsg} et le lemme~\ref{ch2lemmepeuttuerp} montrent qu'il
existe un point \mbox{$M_0 \in \Mrond''$}, une place $v_0 \in \Omega \setminus T$ et une
place~$w_0$ de~$\kappa(M_0)$ totalement décomposée dans $L'_{M_0}L''_{M_0}K_{M_0}$,
divisant~$v_0$, non ramifiée et de degré résiduel~$1$ sur~$v_0$, tels que la condition
suivante soit satisfaite:
\begin{itemize}
\medskip
\item[(\refstepcounter{equation}\theequation{})\label{ch2condxvz}]
pour tout $x_{v_0}\in \A^1(\Orond_{v_0})$ rencontrant transversalement $\Mtildez$ en $w_0$,
l'image de $\psi(\alpha)$ dans $H^1(k_{v_0}, \Erond'_{x_{v_0}})$ est non nulle.
\end{itemize}
\medskip

Pour tout $M \in \Mrond$ et tout $\gamma$ appartenant à $\SG_{\phi',T}(\A^1_k,\Erond')$
(resp.~à $\SG_{\phi'',T}(\A^1_k,\Erond'')$), fixons une extension
quadratique ou triviale $K_{M,\gamma}/\kappa(M)$ vérifiant les conditions du
lemme~\ref{ch2existekmd} associées à l'image de~$\gamma$ dans $H^1(\kappa(M),\ff{F'}M)$
(resp.~dans $H^1(\kappa(M),\ff{F''}M)$).

Soit $\Trond^+=(\ft{T'^+}M)_{M \in \Mrond}$ la famille définie par
$\ft{T'^+}M=\ft{T'}M$ pour $M \neq M_0$ et $\ft{T'^+}{M_0}=\ft{T'}{M_0}\cup\{w_0\}$.
C'est une famille admissible.
Fixons $x_{v_0} \in \A^1(\Orond_{v_0})$ rencontrant
transversalement $\Mtildez$ en~$w_0$.  Un tel~$x_{v_0}$ existe car la
place~$w_0$ est non ramifiée et de degré résiduel~$1$ sur~$v_0$.  D'après la
proposition~\ref{ch2propschinzel}, il existe $x^+ \in U(k)$ arbitrairement
proche de~$x$ pour $v \in T\cap\Omega_f$ et arbitrairement grand aux places
archimédiennes, tel que le couple $(\Trond^+, x^+)$ soit
préadmissible.

\bigskip
\begin{lemme}
\label{ch2wvconst}
Pour tout $v \in \Omega$, les sous-groupes $\Erond'_{x_v}(k_v)/\Im(\phi''_{x_v})$
et $\Erond''_{x_v}(k_v)/\Im(\phi'_{x_v})$ de $H^1(k_v,\Z/2)$
sont des fonctions
localement constantes de $x_v \in (\P^1_k \setminus \Mrond)(k_v)$.
\end{lemme}

\bigskip
\begin{demo}
Démonstration similaire à celle du lemme~\ref{ch1wvconst}.
\end{demo}

\bigskip
Notons $\fv{V'}v$, $W'_v$, $\fv{V''}v$, $W''_v$, $\ft{T'}v$ et~$\ft{T''}v$
(resp.~$\fv{V'^+}v$, $W'^+_v$,
$\fv{V''^+}v$, $W''^+_v$, $\ft{T'^+}v$ et~$\ft{T''^+}v$)
les espaces définis au paragraphe~\ref{ch2dualitelocale} relativement
au point~$x$ (resp.~au point~$x^+$). Grâce au lemme~\ref{ch2wvconst} et à la
définition des voisinages~$\Arond_v$,
on peut supposer les conditions suivantes satisfaites, quitte à choisir~$x^+$ assez proche
de~$x$ aux places de $T\cap\Omega_f$ et assez grand aux places archimédiennes:
\begin{itemize}
\medskip
\item[(\refstepcounter{equation}\theequation{})]
$x \in \Arond_v$ pour tout $v \in S$;
\medskip
\item[(\refstepcounter{equation}\theequation{})\label{ch2hypa}]
pour tout $v \in T$ et tout $\gamma \in H^1(U,\Z/2)$,
l'image de $\gamma(x)$ dans $\fv{V'}v$ (resp.~$\fv{V''}v$) appartient à~$W'_v$ (resp.~$W''_v$) si
et seulement si l'image de $\gamma(x^+)$ dans $\fv{V'^+}v$ (resp.~$\fv{V''^+}v$) appartient à~$W'^+_v$
(resp.~$W''^+_v$);
\medskip
\item[(\refstepcounter{equation}\theequation{})\label{ch2hypb}]
pour tout
$\gamma \in \SG_{\phi',T}(\A^1_k,\Erond')$
(resp.~tout~$\gamma \in \SG_{\phi'',T}(\A^1_k, \Erond'')$),
tout \mbox{$M \in \Mrond$} et tout $v \in T$,
on a $\inv_v A_{K_{M,\gamma}/\kappa(M)}(x) = \inv_v A_{K_{M,\gamma}/\kappa(M)}(x^+)$.
\end{itemize}
\medskip

(Pour la condition~(\ref{ch2hypa}), on utilise la définition des voisinages~$\Arond_v$
pour~$v$ réelle et le fait que $H^1(U,\Z/2)=\Gm(U)/2$ est engendré par un sous-groupe
fini et par $\Gm(k)/2$.)

Notons $T^+$, $T^+(x^+)$, $w_M^+$ et $v_M^+$ pour $M \in \Mrond$ les données associées,
dans le paragraphe~\ref{ch2couplesadm}, au
couple $(\Trond^+, x^+)$
et $$\psi^+ \colon H^1(\Orond_{T^+(x^+)},\Z/2) \longrightarrow H^1(\Urond_{\Orond_{T^+}},\Z/2)$$
l'isomorphisme inverse de l'évaluation en~$x^+$.
On a $T^+=T \cup \{v_0\}$.

Comme la place~$w_0$ de~$\kappa(M_0)$ est totalement décomposée dans $L'_{M_0}L''_{M_0}K_{M_0}$,
l'admissibilité de la famille~$\Trond$ entraîne celle de~$\Trond^+$.
L'admissibilité de la famille $\Trond^+$,
la proposition~\ref{ch2invtotdec} (appliquée au couple $(\Trond^+,x^+)$,
aux places de $T^+ \setminus S$
et aux extensions $E \in \{L'_M,L''_M,K_M\}$),
le lemme~\ref{ch2suminvax} (appliqué au point~$x^+$)
et la proposition~\ref{ch2recprop} (appliquée au couple $(\Trond^+,x^+)$
et aux extensions $E \in \{L'_M,L''_M,K_M\}$)
montrent alors que le
couple $(\Trond^+, x^+)$ est admissible.   Nous allons maintenant établir
les inclusions~(\ref{ch2lincp}) et~(\ref{ch2lincpp}) relatives à~$x^+$ et les
inégalités
\begin{equation}
\label{ch2inegun}
\dim_{\ff{\F}2}(\Sel_{\phi'_{x^+}}(k, \Erond'_{x^+})) <
\dim_{\ff{\F}2}(\Sel_{\phi'_x}(k, \Erond'_x))
\end{equation}
et
\begin{equation}
\label{ch2inegdeux}
\dim_{\ff{\F}2}(\Sel_{\phi''_{x^+}}(k, \Erond''_{x^+})) \leq
\dim_{\ff{\F}2}(\Sel_{\phi''_x}(k, \Erond''_x)) \rlap{\text{;}}
\end{equation}
on aura alors abouti à une contradiction, compte tenu de la minimalité
du couple $(m',m'')$.

\bigskip
\begin{lemme}
\label{ch2v0compat}
Le diagramme
$$
\myxyin\myxyhook\xymatrix{
H^1(\Urond_{\Orond_{T^+}},\Z/2) \ar[d] \ar[r] & H^1(K,Z/2)\ar@{=}[r] &
*+!!<0pt,\the\fontdimen22\textfont2>[r]!<5.7ex,0ex>{K^\star/K^{\star 2}}\ar[r]^(.375){v_0}&\Z/2
\ar@{=}[d] \\
H^1(\Orond_{T^+(x^+)},\Z/2) \ar[r] & H^1(k,\Z/2)\ar@{=}[r] &
*+!!<0pt,\the\fontdimen22\textfont2>[r]!<5.7ex,0ex>{k^\star/k^{\star 2}}
\ar[r]^(.375){v_0+v_{M_0}^+}&\Z/2\rlap{\text{,}}
}\myxyout
$$
dans lequel la flèche verticale de gauche est l'évaluation en~$x^+$ et la flèche de droite
de la première ligne est induite par la valuation $\kappa(\P^1_\Orond)^\star\rightarrow \Z$
associée à la place~$v_0$, est commutatif.
\end{lemme}

\bigskip
\begin{demo}
Il suffit d'appliquer le lemme~\ref{ch2transvcomm} deux fois, d'abord avec \mbox{$v=v_0$},
puis avec~$v=v_{M_0}^+$, en remarquant que la flèche composée
$$
\myxyin\myxyhook\xymatrix{
H^1(\Urond_{\Orond_{T^+}},\Z/2) \ar[r] & H^1(K,Z/2)\ar@{=}[r] &
*+!!<0pt,\the\fontdimen22\textfont2>[r]!<5.7ex,0ex>{K^\star/K^{\star 2}}\ar[r]^(.425){\smash[t]{v_{M_0}^+}}&\Z/2
}\myxyout
$$
est nulle.
\end{demo}

\bigskip
\begin{lemme}
\label{ch2psiplusot}
L'image du groupe $\Sel_{\phi'_{x^+}}(k,\Erond'_{x^+})$ par~$\psi^+$ est incluse
dans $H^1(\Urond_{\Orond_T},\Z/2)$.
\end{lemme}

\bigskip
\begin{demo}
Vu le lemme~\ref{ch2v0compat} et la proposition~\ref{ch2h1u}, il suffit de prouver que tout élément
de $\Sel_{\phi'_{x^+}}(k,\Erond'_{x^+})$ est de valuation nulle en~$v_0$ et
en~$v_{M_0}^+$.  Ceci résulte de la conclusion de la proposition~\ref{ch2ekvsurj}.
\end{demo}

\bigskip
\begin{lemme}
\label{ch2psiplusotev}
L'image de $\psi^+(\Sel_{\phi'_{x^+}}(k,\Erond'_{x^+}))$ par la flèche
$H^1(\Urond_{\Orond_T},\Z/2) \rightarrow H^1(\Orond_{T(x)},\Z/2)$ d'évaluation
en~$x$ est incluse dans $\Sel_{\phi'_x}(k,\Erond'_x)$.
\end{lemme}

\bigskip
\begin{demo}
Soit $\beta \in \Sel_{\phi'_{x^+}}(k,\Erond'_{x^+})$.
Comme $\psi^+(\beta)\in H^1(\Urond_{\Orond_T},\Z/2)$, l'hypothèse~(\ref{ch2hypa})
assure que l'image de~$\psi^+(\beta)(x)$ dans~$\fv{V'}v$ appartient à~$W'_v$ pour tout~$v \in T$.
Pour $v \in \Omega \setminus T(x)$, on a $W'_v=\ft{T'}v$
puisque~$\Erond'_x$ a bonne réduction en~$v$ (lemme~\ref{ch2bonnered}),
or $\psi^+(\beta)(x) \in H^1(\Orond_{T(x)},\Z/2)$,
donc l'image de~$\psi^+(\beta)(x)$ dans~$\fv{V'}v$ appartient à~$W'_v$.
La conclusion de la proposition~\ref{ch2ekvsurj} montre par ailleurs que
pour tout $M \in \Mrond'$, l'image de~$\psi^+(\beta)(x)$
dans~$\fv{V'}{v_M}$ appartient à~$W'_{v_M}$.

Restent à considérer les places $v_M$ pour $M\in\Mrond''$.  Fixons $M \in \Mrond''$.
La proposition~\ref{ch2psiselsg} et le lemme~\ref{ch2psiplusot} montrent que
$\psi^+(\beta)\in\SG_{\phi',T}(\A^1_k,\Erond')$,
de sorte qu'une extension quadratique ou triviale
$K_{M,\psi^+(\beta)}/\kappa(M)$
vérifiant les conditions du
lemme~\ref{ch2existekmd} a été choisie précédemment.
La proposition~\ref{ch2reckmm}, appliquée à la place~$v_M^+$, au point~$x^+$
et à la classe $m=\psi^+(\beta)$,
montre que la place~$w_M^+$ de~$\kappa(M)$ est totalement décomposée dans $K_{M,\psi^+(\beta)}$;
les hypothèses de cette proposition sont satisfaites
parce que $\psi^+(\beta) \in H^1(\Urond_{\Orond_T},\Z/2)$
et que le couple $(\Trond^+,x^+)$ est admissible.
L'extension $K_{M,\psi^+(\beta)}/\kappa(M)$ est non ramifiée en toute place de~$\kappa(M)$
dont la trace sur~$k$ n'appartient pas à~$T$, d'après la première partie de
la proposition~\ref{ch2reckmm}.  On en déduit, grâce à la proposition~\ref{ch2recprop},
l'égalité
\begin{equation}
\label{ch2sominv1}
\sum_{v \in T^+} \inv_v A_{K_{M,\psi^+(\beta)}/\kappa(M)}(x^+)=0 \rlap{\text{.}}
\end{equation}
Si $M=M_0$, la proposition~\ref{ch2reckmm},
appliquée à~$v_0$ et~$x^+$, montre que la place~$w_0$ de~$\kappa(M)$
est totalement décomposée dans $K_{M,\psi^+(\beta)}$. Comme aucune place de~$\kappa(M)$
divisant~$v_0$ n'est ramifiée dans $K_{M,\psi^+(\beta)}$, la proposition~\ref{ch2invtotdec}
permet d'en déduire dans tous les cas que
\begin{equation}
\label{ch2invv0nul}
\inv_{v_0} A_{K_{M,\psi^+(\beta)}/\kappa(M)}(x^+)=0 \rlap{\text{.}}
\end{equation}
Des équations~(\ref{ch2sominv1}) et~(\ref{ch2invv0nul}) et de l'hypothèse~(\ref{ch2hypb}),
on tire l'égalité
\begin{equation}
\sum_{v \in T} \inv_v A_{K_{M,\psi^+(\beta)}/\kappa(M)}(x)=0 \rlap{\text{.}}
\end{equation}
Compte tenu de la proposition~\ref{ch2recprop}, cela entraîne que la place~$w_M$
de~$\kappa(M)$ est totalement décomposée dans $K_{M,\psi^+(\beta)}$.
Enfin, la proposition~\ref{ch2reckmm} permet d'en déduire
que l'image de~$\psi^+(\beta)(x)$ dans~$\fv{V'}{v_M}$
appartient à~$W'_{v_M}$.
\end{demo}

\bigskip
Au vu de la proposition~\ref{ch2propeviso} et des lemmes~\ref{ch2psiplusot} et~\ref{ch2psiplusotev},
on a exhibé une injection du groupe $\Sel_{\phi'_{x^+}}(k,\Erond'_{x^+})$ dans
$\Sel_{\phi'_x}(k,\Erond'_x)$.
Il résulte de l'hypothèse~(\ref{ch2condxvz}) que
$\psi(\alpha)(x^+) \not\in \Sel_{\phi'_{x^+}}(k,\Erond'_{x^+})$. Par conséquent, l'image de
cette injection ne contient pas~$\alpha$.
En particulier, ce n'est pas une bijection; l'inégalité~(\ref{ch2inegun})
est donc prouvée.

\bigskip
Intéressons-nous maintenant à l'inégalité~(\ref{ch2inegdeux}). On commence par
établir un analogue du lemme~\ref{ch2psiplusot}, dont la preuve sera néanmoins considérablement
plus complexe.

\bigskip
\begin{lemme}
\label{ch2psiplusotpp}
L'image du groupe $\Sel_{\phi''_{x^+}}(k,\Erond''_{x^+})$ par~$\psi^+$ est incluse
dans $H^1(\Urond_{\Orond_T},\Z/2)$.
\end{lemme}

\bigskip
\begin{demo}
D'après le lemme~\ref{ch2v0compat}
et la proposition~\ref{ch2h1u},
il suffit de prouver l'égalité $v_0(\beta)=v_{M_0}^+(\beta)$
dans~$\Z/2$ pour tout $\beta \in \Sel_{\phi''_{x^+}}(k,\Erond''_{x^+})$.
Fixons donc
$\beta \in \Sel_{\phi''_{x^+}}(k,\Erond''_{x^+})$ et considérons le cup-produit
$\psi(\alpha)(x^+) \cup \beta \in \tors{2}\Br(k)$ des classes $\psi(\alpha)(x^+),\;\!\!\beta\in
H^1(k,\Z/2)$.  La loi de réciprocité globale\index{loi de réciprocité globale} permet d'écrire que
\begin{equation}
\label{ch2psiabrecglob}
\sum_{v \in\Omega} \inv_v(\psi(\alpha)(x^+)\cup\beta)=0 \rlap{\text{.}}
\end{equation}

Montrons dans un premier temps que $\inv_v(\psi(\alpha)(x^+)\cup\beta)=0$
pour tout $v \in \Omega \setminus \{v_0, v_{M_0}^+\}$.
Pour tout $v \in \Omega$, cet invariant est égal à la valeur de
l'accouplement $\fv{V'^+}v \times \fv{V''^+}v \rightarrow \Z/2$ (cf.~paragraphe~\ref{ch2prelude})
sur l'image du couple $(\psi(\alpha)(x^+), \beta)$.
L'image de~$\beta$ dans~$\fv{V''^+}v$ appartient à $W''^+_v$ par hypothèse.  Comme
les sous-groupes $W'^+_v \subset \fv{V'^+}v$ et $W''^+_v \subset \fv{V''^+}v$ sont orthogonaux
pour cet accouplement, il suffit de prouver que l'image de $\psi(\alpha)(x^+)$
dans~$\fv{V'^+}v$ appartient à~$W'^+_v$ pour tout $v \in \Omega\setminus\{v_0,v_{M_0}^+\}$.
Pour $v \in T$, c'est une conséquence de l'hypothèse~(\ref{ch2hypa}).
Pour $v \in \Omega \setminus T^+(x^+)$, cela provient de ce que $W'^+_v=\ft{T'^+}v$
(lemme~\ref{ch2bonnered}).
Pour $v \in \ensemble{v_M^+}{M\in\Mrond'}$, la proposition~\ref{ch2ekvsurj}
montre que $W'^+_v=\fv{V'^+}v$.
Seules restent les places $v_M^+$
pour $M\in\Mrond''\setminus\{M_0\}$. Soit $M\in\Mrond''\setminus\{M_0\}$.
D'après les propositions~\ref{ch2psiselsg}
et~\ref{ch2reckmm}, la place~$w_M$ de~$\kappa(M)$ est totalement décomposée
dans~$K_{M,\psi(\alpha)}$.  On en déduit, à l'aide de
la proposition~\ref{ch2recprop}, applicable grâce à la première partie
de la proposition~\ref{ch2reckmm}, l'égalité
$$
\sum_{v \in T}\inv_v A_{K_{M,\psi(\alpha)}/\kappa(M)}(x)=0\rlap{\text{.}}
$$
Compte tenu
de l'hypothèse~(\ref{ch2hypb}),
de la proposition~\ref{ch2invtotdec} appliquée à la place $v=v_0$
et de ce que $M \neq M_0$, il
en résulte que
$$
\sum_{v \in T^+}\inv_v A_{K_{M,\psi(\alpha)}/\kappa(M)}(x^+)=0\rlap{\text{.}}
$$
Une nouvelle application de la proposition~\ref{ch2recprop} montre enfin que la
place~$w_M^+$ de~$\kappa(M)$ est totalement décomposée dans~$K_{M,\psi(\alpha)}$,
ce qui implique, d'après la proposition~\ref{ch2reckmm},
que l'image de $\psi(\alpha)(x^+)$ dans~$\fv{V'^+}{v^+_M}$ appartient bien à~$W'^+_{v_M^+}$.

Ainsi a-t-on prouvé que l'équation~(\ref{ch2psiabrecglob}) se réduit à
\begin{equation}
\label{ch2v0vm0plus}
\inv_{v_0} (\psi(\alpha)(x^+)\cup\beta) = \inv_{v_{M_0}^+} (\psi(\alpha)(x^+)\cup\beta)\rlap{\text{.}}
\end{equation}
Supposons momentanément que l'image de~$\psi(\alpha)(x^+)$
dans~$\fv{V'^+}{v_{M_0}^+}$ soit nulle.  La place~$w_{M_0}^+$ de~$\kappa(M_0)$
est alors totalement décomposée dans~$K_{M_0,\psi(\alpha)}$ (proposition~\ref{ch2reckmm}).
Il s'ensuit (proposition~\ref{ch2recprop}) que
\begin{equation}
\label{ch2sominv2}
\sum_{v \in T^+} \inv_v A_{K_{M_0,\psi(\alpha)}/\kappa(M_0)}(x^+)=0\rlap{\text{.}}
\end{equation}
On a par ailleurs, à nouveau grâce aux propositions~\ref{ch2reckmm} et~\ref{ch2recprop}:
\begin{equation}
\label{ch2sominv3}
\sum_{v \in T} \inv_v A_{K_{M_0,\psi(\alpha)}/\kappa(M_0)}(x)=0\rlap{\text{.}}
\end{equation}
De l'hypothèse~(\ref{ch2hypb}) et des deux équations~(\ref{ch2sominv2})
et~(\ref{ch2sominv3}) résulte l'égalité
$$
\inv_{v_0} A_{K_{M_0,\psi(\alpha)}/\kappa(M_0)}(x^+)=0\rlap{\text{.}}
$$
Les propositions~\ref{ch2invtotdec} et~\ref{ch2reckmm} permettent d'en déduire
que la place~$w_0$ est totalement décomposée dans~$K_{M_0,\psi(\alpha)}$ puis que
l'image de~$\psi(\alpha)(x^+)$ dans $H^1(k_{v_0},\Erond'_{x^+})$ est nulle,
contredisant ainsi l'hypothèse~(\ref{ch2condxvz}).
Nous avons donc prouvé, par l'absurde, que l'image de~$\psi(\alpha)(x^+)$
dans~$\fv{V'^+}{v_{M_0}^+}$ n'est pas nulle.  Compte tenu de l'hypothèse~(\ref{ch2condxvz}),
l'image de~$\psi(\alpha)(x^+)$ dans~$\fv{V'^+}v$ n'est donc nulle
pour aucun $v \in \{v_0,v_{M_0}^+\}$.  Montrons maintenant qu'elle
appartient à~$\ft{T'^+}v$
pour tout $v \in \{v_0,v_{M_0}^+\}$.
Le lemme~\ref{ch2transvcomm}, appliqué deux fois,
entraîne que $v_{M_0}^+(\psi(\alpha)(x^+))=v_{M_0}(\psi(\alpha)(x))=v_{M_0}(\alpha)$,
or $v_{M_0}(\alpha)=0$ d'après la proposition~\ref{ch2ekvsurj}.
Ceci prouve le résultat voulu pour $v=v_{M_0}^+$. La nullité de $v_0(\psi(\alpha)(x^+))$ se déduit
de celle de $v_{M_0}^+(\psi(\alpha)(x^+))$ grâce au lemme~\ref{ch2v0compat}.

En vertu de la formule du symbole modéré (cf.~\cite[Ch. XIV, Proposition~8]{serrecl}),
on déduit de l'appartenance de l'image de~$\psi(\alpha)(x^+)$ dans~$\fv{V'^+}v$
à~$\ft{T'^+}v \setminus \{1\}$
l'équivalence
$$
\inv_v(\psi(\alpha)(x^+)\cup\beta)=0 \Longleftrightarrow v(\beta)=0
$$
pour $v \in \{v_0,v_{M_0}^+\}$.  Vu l'équation~(\ref{ch2v0vm0plus}),
ceci achève la démonstration du lemme.
\end{demo}

\bigskip
\begin{lemme}
\label{ch2psiplusotevpp}
L'image de $\psi^+(\Sel_{\phi''_{x^+}}(k,\Erond''_{x^+}))$ par la flèche
$H^1(\Urond_{\Orond_T},\Z/2) \rightarrow H^1(\Orond_{T(x)},\Z/2)$ d'évaluation
en~$x$ est incluse dans $\Sel_{\phi''_x}(k,\Erond''_x)$.
\end{lemme}

\bigskip
\begin{demo}
On démontre ce lemme à partir du lemme~\ref{ch2psiplusotpp} exactement comme on
a démontré le lemme~\ref{ch2psiplusotev} à partir du lemme~\ref{ch2psiplusot};
il suffit d'échanger tous les~$'$ et les~$''$ et de remplacer
la référence au lemme~\ref{ch2psiplusot} par une référence au lemme~\ref{ch2psiplusotpp}.
\end{demo}

\bigskip
Il existe une injection naturelle
$\Sel_{\phi''_{x^+}}(k,\Erond''_{x^+}) \hookrightarrow \Sel_{\phi''_x}(k,\Erond''_x)$, d'après les lemmes~\ref{ch2psiplusotpp} et~\ref{ch2psiplusotevpp}
et la proposition~\ref{ch2propeviso}.
L'inégalité~(\ref{ch2inegdeux}) est donc établie. Les lemmes~\ref{ch2psiplusotev}
et~\ref{ch2psiplusotevpp} montrent de plus que
$$
\psi^+(\Sel_{\phi'_{x^+}}(k,\Erond'_{x^+})) \subset \psi(\Sel_{\phi'_x}(k,\Erond'_x))
$$
et que
$$
\psi^+(\Sel_{\phi''_{x^+}}(k,\Erond''_{x^+}))\subset
\psi(\Sel_{\phi''_x}(k,\Erond''_x))\rlap{\text{.}}
$$
Grâce aux hypothèses~(\ref{ch2lincp}) et~(\ref{ch2lincpp}), on en déduit que
$$\psi^+(\Sel_{\phi'_{x^+}}(k,\Erond'_{x^+}))\subset\SG_{\phi',S}(\A^1_k,\Erond')$$
et
$$\psi^+(\Sel_{\phi''_{x^+}}(k,\Erond''_{x^+}))\subset\SG_{\phi'',S}(\A^1_k,\Erond'')\rlap{\text{.}}$$
La minimalité du couple~$(m',m'')$ fournit maintenant une contradiction.

\bigskip
Ainsi avons-nous établi, par l'absurde, que $\psi(\alpha)\in\SDSzp{S}$
pour tout \mbox{$\alpha\in\Sel_{\phi'_x}(k,\Erond'_x)$}.  Échangeant ci-dessus
d'une part tous les~$'$ et les~$''$
(sauf bien sûr dans $\ft{T'}M$ et dans $(m',m'')$)
et remplaçant d'autre part le lemme~\ref{ch2lemmepeuttuerp} par le lemme suivant,
on obtient une preuve de l'assertion duale: $\phi(\alpha)\in\SDSzpp{S}$ pour tout
$\alpha\in\Sel_{\phi''_x}(k,\Erond''_x)$.
D'où finalement $x \in \RDSz{S}$.

\bigskip
\begin{lemme}
\label{ch2lemmepeuttuerpp}
Pour tout $e \in \tors{\phi''}H^1(\A^1_k,\Erond'') \setminus \TDzpp$, il existe $M \in \Mrond'$
et une infinité de places $v \in \Omega \setminus S_0$
telles qu'il existe une place~$w$ de~$\kappa(M)$
totalement décomposée dans $L'_ML''_MK_M$, divisant~$v$, non ramifiée et de degré résiduel~$1$
sur~$v$, telle que pour tout $x_v \in \A^1(\Orond_v)$ rencontrant transversalement~$\Mtilde$
en~$w$, l'image de~$e$ par la flèche $H^1(\A^1_k,
\Erond'') \rightarrow H^1(k_v, \Erond''_{x_v})$ d'évaluation en~$x_v$ soit non nulle.
\end{lemme}

\bigskip
\begin{demo}
L'énoncé du lemme~\ref{ch2lemmepeuttuerpp} s'obtient à partir de
celui du lemme~\ref{ch2lemmepeuttuerp} en échangeant les~$'$ et les~$''$; néanmoins,
la preuve du lemme~\ref{ch2lemmepeuttuerpp} ne se déduit pas formellement
de celle du lemme~\ref{ch2lemmepeuttuerp}.
L'opération naturelle que l'on peut appliquer à la fois à l'énoncé et à la preuve
du lemme~\ref{ch2lemmepeuttuerp} sans les invalider est celle consistant à
échanger les~$'$ et les~$''$ et à remplacer
respectivement $\delta'_M([\Xrond])$ et~$K_M$ par $\delta''_M([\Xrond''])$
et~$K_M''$, où~$K_M''$ désigne le corps $K_{M,d}$ donné par le
lemme~\ref{ch2existekmd} en prenant pour~$d$ l'image de~$[\Xrond'']$
dans $H^1(\kappa(M),\ff{F''}M)$.
Ceci prouve que le lemme~\ref{ch2lemmepeuttuerpp} devient vrai si l'on
remplace~$K_M$ par~$K''_M$ dans son énoncé. Pour conclure, il suffit donc
d'établir l'existence d'un plongement $\kappa(M)$\nobreakdash-linéaire
$K_M \hookrightarrow L'_ML''_MK''_M$ pour tout $M\in \Mrond'$.

Soit $M \in \Mrond'$. La suite
\begin{equation}
\label{ch2peuttuerppse}
\xymatrix{
0 \ar[r] & \ff{F'}M \ar[r] & \ff{F''}M \ar[r] & \Z/2 \ar[r] & 0
}
\end{equation}
est exacte (cf.~proposition~\ref{annpropppp}).  Il en résulte d'une part
que l'extension $L'_M/\kappa(M)$
se plonge dans $L''_M/\kappa(M)$ et d'autre part
que la flèche $$H^1(L''_M K''_M, \ff{F'}M) \longrightarrow H^1(L''_M K''_M, \ff{F''}M)$$ induite
par~$\phi'_M$ est injective,
compte tenu que le $L''_M$\nobreakdash-groupe $\ff{F''}M\otimes_{\kappa(M)}L''_M$ est constant.
On dispose donc d'un carré commutatif
$$
\myxyhook\xymatrix{
H^1(L'_M,\ff{F'}M) \ar[d] \ar[r] & H^1(L''_M K''_M, \ff{F'}M) \ar@{_{ (}->}[d] \\
H^1(L''_M,\ff{F''}M) \ar[r] & H^1(L''_M K''_M, \ff{F''}M)
}
$$
induit par~$\phi'_M$ et par les diverses flèches de restrictions.
L'image de $\delta_M'([\Xrond])$ par la flèche verticale de gauche
est égale à $\delta_M''([\Xrond''])$, or celle-ci appartient au noyau de la flèche
horizontale inférieure par définition de~$K''_M$. La classe
$\delta_M'([\Xrond])$ appartient donc au noyau de la flèche horizontale supérieure,
ce qui signifie que $L'_M K_M$ se plonge $L'_M$\nobreakdash-linéairement dans $L''_M K''_M$,
comme il fallait démontrer.
\end{demo}

\bigskip
Supposons maintenant que~$\Tinfty\neq\emptyset$.  Il reste à prouver que $x \in \RD$.
Il suffit pour cela de vérifier que pour tout $\alpha\in\Sel_{\phi'_x}(k,\Erond'_x)$
(resp.~$\alpha\in\Sel_{\phi''_x}(k,\Erond''_x)$), l'image de $\psi(\alpha)$
dans $H^1(K,\Z/2)=K^\star/K^{\star 2}$ est de valuation nulle au point $\infty \in \P^1_k$.

\bigskip
\begin{lemme}
Les groupes $\SDzp \cap H^1(\Urond_{\Orond_T},\Z/2)$
et $\SDzpp \cap H^1(\Urond_{\Orond_T},\Z/2)$
sont inclus dans
$H^1(\Urond_{\Orond_{T\setminus\{v_\infty\}}},\Z/2)$.
\end{lemme}

\bigskip
\begin{demo}
Soit $a \in \SDzp \cap H^1(\Urond_{\Orond_T},\Z/2)$.
La proposition~\ref{ch2h1u} et l'hypothèse~(\ref{ch2hypprelude}) pour~$S=S_1$
assurent
l'existence de $c \in H^1(\Urond_{\Orond_{S_1}},\Z/2)$
tel que pour tout $M \in \Mrond$, les images de~$a$ et de~$c$
dans $H^1(K_M^\sh,\Z/2)$ coïncident.
Posant $b=a-c$, on a nécessairement $b \in H^1(\Orond_T,\Z/2)$.
Il résulte de la construction de~$v_\infty$ que l'image de~$b$ dans $k^\star/k^{\star 2}$ est
de valuation nulle en~$v_\infty$, autrement dit que $b \in H^1(\Orond_{T\setminus
\{v_\infty\}},\Z/2)$.
On a alors $a \in H^1(\Urond_{\Orond_{T\setminus\{v_\infty\}}},\Z/2)$, comme annoncé.
L'inclusion de $\SDzpp \cap H^1(\Urond_{\Orond_T},\Z/2)$ dans
$H^1(\Urond_{\Orond_{T\setminus\{v_\infty\}}},\Z/2)$ se prouve exactement de la
même manière.
\end{demo}

\bigskip
Il n'y a plus qu'à appliquer le lemme~\ref{ch2transvcomm} à la place~$v_\infty$
pour conclure, compte tenu de ce que $v_\infty(\alpha)=0$
puisque~$v_\infty$ est une place de bonne réduction pour~$\Erond'_x$ et~$\Erond''_x$
(lemme~\ref{ch2bonnered}).
Ceci achève la démonstration de la proposition~\ref{ch2pluslongueprop}.
\end{demo}

\bigskip
Le théorème~\ref{ch2thprin} est maintenant prouvé.
En effet, $x$ est entier en dehors de~$S_1$ lorsque $\Tinfty=\emptyset$,
puisqu'il est entier en dehors de~$S$ (préadmissibilité du couple $(\Trond,x)$)
et que $S=S_1$ si $\Tinfty=\emptyset$.

\chapter{Principe de Hasse pour les surfaces de del Pezzo de degré~$4$}
\section{Introduction}

Ce chapitre est consacré aux deux conjectures suivantes.

\bigskip
\begin{conjectures}
\label{ch3introconj}
Soit~$k$ un corps de nombres.

(i) Soit~$X$ une surface de del Pezzo de degré~$4$ sur~$k$.
Si $X(\A_k)^\Br\neq\emptyset$, alors $X(k)\neq\emptyset$.

(ii) Soit $n \geq 5$. Toute intersection lisse de deux quadriques dans~$\P^n_k$ satisfait
au principe de Hasse.
\end{conjectures}

\bigskip
(La première de ces conjectures apparut d'abord
sous la forme d'une question dans~\cite{ctsannoordhoff};
la seconde en est une conséquence (cf.~\cite{harduke}).
Voir également \cite[§16]{ctsansdii} et \cite[§5.5 et §5.6]{cttoulouse}.
Par «~quadrique~», on entend
bien sûr «~hypersurface quadrique~».)

Rappelons qu'une \emph{surface de del\index{surface de del Pezzo} Pezzo}~$X$ sur~$k$ est une
surface projective et lisse sur~$k$, de faisceau anti-canonique ample; le \emph{degré} de~$X$,
compris entre~$1$ et~$9$,
est le nombre d'auto-intersection de son faisceau canonique.
L'étude systématique de l'arithmétique des surfaces de del Pezzo est justifiée par
le théorème de classification suivant, dû à Enriques, Manin et Iskovskikh~\cite{iskomin}:
toute surface rationnelle sur~$k$
est $k$\nobreakdash-birationnelle à une surface de del Pezzo ou à un fibré en coniques au-dessus d'une
conique.  Que les fibrés en coniques au-dessus d'une conique admettent un point
rationnel dès que l'obstruction de Brauer-Manin ne s'y oppose pas
est connu lorsqu'il y a au plus cinq fibres géométriques singulières
(cf.~\cite{ctsansdi}, \cite{ctsansdii}, \cite{ctsemtn}, \cite{salpezzoconique},
\cite{salskoweak}),
en toute généralité si l'on admet l'hypothèse de Schinzel
(cf.~\cite{ctsanschinzel}, \cite{cscrelle94}).
D'autre part,
on sait depuis Manin~\cite{maninratsurf} et Swinnerton-Dyer~\cite{swddelpezzo5} que
les surfaces de del Pezzo de degré~$\geq 5$ sur~$k$ satisfont toujours au principe de Hasse.
Les surfaces de del Pezzo de degré~$4$ constituent donc la classe de surfaces rationnelles
la plus simple pour laquelle il n'est pas connu, même en admettant l'hypothèse de Schinzel,
qu'il existe un point rationnel dès que l'obstruction de Brauer-Manin ne s'y oppose pas.
On sait en revanche qu'elles ne satisfont pas toujours au principe de Hasse.
Le premier contre-exemple fut trouvé par Birch et Swinnerton-Dyer~\cite{bsdcounter},
et tous les contre-exemples connus sont expliqués par l'obstruction de Brauer-Manin.

Les surfaces de del Pezzo de degré~$4$ sur~$k$ sont les surfaces lisses qui
s'écrivent comme l'intersection de deux
quadriques dans~$\P^4_k$ (cf.~\cite[p.~96]{maninratsurf});
ainsi les deux conjectures~\ref{ch3introconj} sont-elles très proches.
La question du principe de Hasse pour les intersections lisses
de deux quadriques dans~$\P^n_k$ possède une longue histoire, dont les grandes dates
furent les suivantes.
En~1959, Mordell~\cite{mordell12} prouva la conjecture~\ref{ch3introconj}~(ii)
pour $n \geq 12$ et $k=\Q$.
En~1964, Swinnerton-Dyer~\cite{swd10} établit la conjecture~\ref{ch3introconj}~(ii) pour $n \geq 10$
et~$k=\Q$, améliorant ainsi le résultat de Mordell
(cf.~\cite[Remark~10.5.2]{ctsansdii} pour une correction).
En~1971, Cook~\cite{cook} démontra le principe de Hasse
pour les intersections lisses de deux quadriques
définies par des formes quadratiques simultanément diagonales
dans~$\P^n_k$ avec $n \geq 8$ et~$k=\Q$.
En~1987, pour~$k$ arbitraire et par des
méthodes radicalement différentes, Colliot-Thélène, Sansuc et Swinnerton-Dyer
prouvèrent la conjecture~\ref{ch3introconj}~(ii)
pour $n \geq 8$ ainsi que dans quelques cas particuliers concernant
les intersections de deux quadriques
qui contiennent soit deux droites gauches conjuguées, soit une quadrique de dimension~$2$
(cf.~\cite{ctsansdi}, \cite{ctsansdii}).
Par la suite, Debbache obtint le cas des intersections lisses
de deux quadriques dans~$\P^7_k$ dont l'une est singulière
(non publié)
et Salberger celui des intersections lisses de deux quadriques
dans~$\P^n_k$ qui contiennent une conique définie sur~$k$, pour $n \geq 5$ (non publié).

Voici maintenant la liste des résultats connus au sujet de la conjecture~\ref{ch3introconj}~(i).
Ils concernent tous
des familles exceptionnelles de surfaces
de del Pezzo de degré~$4$.  Soit $X \subset \P^4_k$ une intersection lisse de deux
quadriques, de dimension~$2$.
Tout d'abord, le principe de Hasse vaut pour la surface~$X$
si elle contient deux droites gauches conjuguées;
en effet, en contractant ces droites, on obtient une surface
de del Pezzo de degré~$6$.  Plus généralement, dès que~$X$ n'est
pas une surface $k$\nobreakdash-minimale, une contraction permet de se ramener à une
surface de del Pezzo de degré~$>4$.
Lorsque~$X$ admet une structure de fibré en coniques sur~$\P^1_k$,
Salberger a démontré que $X(k)\neq\emptyset$ dès que $X(\A_k)^\Br\neq\emptyset$
(cf.~\cite{salpezzoconique}, notamment la remarque qui suit le théorème~(0.8)).
Ainsi le principe de Hasse en l'absence d'obstruction de Brauer-Manin est-il connu lorsque le groupe de Picard de~$X$ est de rang~$\geq 2$.
Enfin, en admettant l'hypothèse de Schinzel et la finitude des groupes de Tate-Shafarevich
des courbes elliptiques sur~$\Q$, Swinnerton-Dyer~\cite{swdegloffstein} a établi le principe
de Hasse pour la surface~$X$ lorsqu'elle est définie par un système d'équations de la forme
\begin{equation}
\label{ch3introeqsd}
\left\{
\begin{aligned}
a_0x_0^2 + a_1x_1^2 + a_2x_2^2 + a_3x_3^2 + a_4x_4^2 &= 0 \rlap{,} \\
b_0x_0^2 + b_1x_1^2 + b_2x_2^2 + b_3x_3^2 + b_4x_4^2 &= 0 \rlap{,}
\end{aligned}
\right.
\end{equation}
où $\uplet{a_0}{a_4},\uplet{b_0}{b_4}\in \Q$ sont «~suffisamment généraux~» en un sens
explicite (et qui implique notamment que $\Br(X)/\Br(\Q)=0$, tout en étant strictement
plus fort).
Ce résultat fut généralisé à tout corps de nombres
par Colliot-Thélène, Skorobogatov et Swinnerton-Dyer (cf.~\cite[§3.2]{css}).

\bigskip
Le but de ce chapitre est d'établir la conjecture~\ref{ch3introconj}~(ii) ainsi qu'une
grande partie de la conjecture~\ref{ch3introconj}~(i) (couvrant notamment le cas d'une surface
de del Pezzo de degré~$4$ «~suffisamment générale~»), en admettant l'hypothèse
de Schinzel et la finitude des groupes de Tate-Shafarevich des courbes elliptiques
sur les corps de nombres.
Signalons tout de suite qu'il y a un
espoir à moyen terme de débarrasser ces résultats de l'hypothèse
de Schinzel: étant donné que pour les surfaces de del Pezzo de degré~$4$, il
est équivalent de s'intéresser aux points rationnels
ou aux $0$\nobreakdash-cycles de degré~$1$
(en tout cas pour ce qui concerne les propriétés considérées ici; cf.~proposition~\ref{ch3findp4amerbrumer} ci-dessous), il suffirait pour cela de réussir à intégrer «~l'astuce de Salberger~» dans
la preuve du théorème principal du chapitre~2.  Pour plus de détails, le
lecteur pourra consulter l'introduction
du chapitre~2, où cette question est discutée.

Les hypothèses précises sous lesquelles nous établissons le principe de Hasse
pour les surfaces de del Pezzo de degré~$4$ étant de nature assez technique, nous
nous contentons d'énoncer ici les cas particuliers les plus simples.
Soit~$X$ une surface de del Pezzo de degré~$4$ sur~$k$.
Soient~$q_1$ et~$q_2$ des formes quadratiques homogènes en~$5$ variables, à coefficients
dans~$k$, telles que~$X$ soit isomorphe à la sous-variété de~$\P^4_k$ définie par
le système d'équations $q_1=q_2=0$.
Le polynôme homogène
$f(\lambda,\mu)=\det(\lambda q_1 + \mu q_2) \in k[\lambda,\mu]$,
est de degré~$5$.
Choisissons un corps de décomposition~$k'$ de~$f$,
numérotons les racines de~$f$ dans~$k'$
et notons $G \subset \mathfrak{S}_5$ le groupe de Galois de~$k'$ sur~$k$.
La classe de conjugaison du sous-groupe $G \subset \mathfrak{S}_5$
ne dépend que de la surface~$X$.

\bigskip
\begin{theoreme}[ (cf.~théorème~\ref{ch3thprindp4})]%
\label{ch3introthdp4}
Admettons l'hypothèse de Schinzel et la finitude des groupes de Tate-Shafarevich des
courbes elliptiques sur les corps de nombres. La surface~$X$ satisfait au principe
de Hasse dès que l'une des conditions suivantes est remplie:
\begin{itemize}
\item le sous-groupe $G \subset \mathfrak{S}_5$ est $3$\nobreakdash-transitif (\emph{i.e.} $G=\mathfrak{A}_5$
ou $G=\mathfrak{S}_5$);
\item le polynôme~$f$ admet exactement deux racines $k$\nobreakdash-rationnelles et d'autre part $\Br(X)/\Br(k)=0$;
\item le polynôme~$f$ est scindé et $\Br(X)/\Br(k)=0$.
\end{itemize}
\end{theoreme}

\bigskip
Notons que si la surface~$X$ est «~suffisamment générale~», on a $G=\mathfrak{S}_5$
(cf.~par exemple \cite[Theorem~7]{waterhouse}),
et en particulier~$G$ est un sous-groupe $3$\nobreakdash-transitif de~$\mathfrak{S}_5$.
Plus délicat, nous établirons que si~$X$ est une section hyperplane «~suffisamment générale~»
d'une intersection lisse de dimension~$3$
de deux quadriques dans~$\P^5_k$ \emph{arbitraire} et \emph{fixée},
on a encore $G=\mathfrak{S}_5$; à l'aide d'un argument de fibration standard, ceci nous
permettra de déduire du théorème~\ref{ch3introthdp4}:

\bigskip
\begin{theoreme}
\label{ch3introp5}
Admettons l'hypothèse de Schinzel et la finitude des groupes de Tate-Shafarevich des
courbes elliptiques sur les corps de nombres.
Soit \mbox{$n \geq 5$}.
Toute intersection lisse de deux quadriques dans $\P^n_k$ satisfait au principe de Hasse.
\end{theoreme}

\bigskip
Lorsque les formes quadratiques~$q_1$ et~$q_2$ sont simultanément diagonales,
c'est-à-dire lorsque la surface~$X$ est définie par un système d'équations de la
forme~(\ref{ch3introeqsd}), le polynôme~$f$ est scindé.
Le théorème~\ref{ch3introthdp4} peut donc être vu comme
une vaste généralisation du résultat de Colliot-Thélène,
Skorobogatov et Swinnerton-Dyer mentionné précédemment; de plus, même dans
le cas particulier où~$q_1$ et~$q_2$ sont simultanément diagonales, il l'améliore, puisqu'il
s'applique alors sans autre hypothèse que la nullité de $\Br(X)/\Br(k)$.
Tous les autres cas du théorème~\ref{ch3introthdp4} sont entièrement nouveaux.
Nous renvoyons le lecteur au paragraphe~\ref{ch3notations} pour des résultats
plus précis que ceux énoncés ici; en réalité, pour toute classe de conjugaison~$\Crond$
de sous-groupes de~$\mathfrak{S}_5$,
nous obtenons le principe de Hasse dès que~$X$ est «~suffisamment générale~»
(en un sens explicite) parmi
les surfaces de del Pezzo de degré~$4$ pour lesquelles~\mbox{$G \in \Crond$}.

\bigskip
Le point de départ de la preuve du théorème~\ref{ch3introthdp4} est
la construction que Swinnerton-Dyer expose dans \cite[§6]{bsd}, où il explique
comment associer à une surface de del Pezzo de degré~$4$ toute une famille de surfaces munies
d'un pinceau de courbes de genre~$1$ et de période~$2$ dont les
jacobiennes admettent une section d'ordre~$2$.  On peut espérer combiner cette construction
avec les méthodes du chapitre~2 pour démontrer
qu'une surface de del Pezzo de degré~$4$ admet un point rationnel dès que l'obstruction
de Brauer-Manin ne s'y oppose pas, sous l'hypothèse de Schinzel et la finitude des groupes
de Tate-Shafarevich.  L'effort de Swinnerton-Dyer en ce sens dans \cite[§6]{bsd}
ne fut que partiellement fructueux.
D'une part, il ne parvint pas à dégager de condition suffisante raisonnable
pour l'existence d'un point rationnel sur une surface de del Pezzo de degré~$4$
(cf.~\cite[Theorem~3]{bsd}), et d'autre part, il s'avère que les arguments de \cite[§6]{bsd}
contiennent plusieurs erreurs cruciales (de sorte que le «~Theorem~3~», notamment,
est faux).

\bigskip
Avant de donner plus de détails sur la preuve du théorème~\ref{ch3introthdp4},
décrivons rapidement l'organisation du chapitre.

Le paragraphe~\ref{ch3parobmvert} est consacré à un résultat d'ordre général
concernant l'obstruction de Brauer-Manin verticale pour des fibrations
au-dessus de~$\P^n_k$ (ou plutôt d'un ouvert de~$\P^n_k$ dont le complémentaire est
de codimension~$\geq 2$).  Nous y prouvons que s'il n'y a pas d'obstruction de Brauer-Manin
verticale à l'existence d'un point rationnel sur l'espace total d'une telle fibration
et que chaque fibre de codimension~$1$ est déployée par une extension abélienne du
corps de base, alors il existe une fibre au-dessus d'un point $k$\nobreakdash-rationnel contenant
un $k_v$\nobreakdash-point lisse pour toute place~$v$ de~$k$, si l'on admet l'hypothèse de Schinzel.
Un tel énoncé est connu lorsque~$n=1$
(cf.~\cite[Theorem~1.1]{csscrelle98}), et la preuve pour~$n$ quelconque consiste
à montrer que l'on peut supposer que $n=1$.  Un ingrédient de première importance
dans cette démonstration est un théorème récent d'Harari, qui établit un analogue de
\cite[Théorème~3.2.1]{harsmf} pour des fibrations en variétés géométriquement intègres
quelconques, avec la contrepartie que seul un sous-groupe fini du groupe de Brauer
de la fibre générique est pris en compte.

Le paragraphe~\ref{ch3sectiongeneralites} contient quelques rappels sur les pinceaux de quadriques
dans~$\P^n$ ainsi qu'une remarque qui nous conduira
à une reformulation plus conceptuelle
de la construction de Swinnerton-Dyer.
Au paragraphe~\ref{ch3notations}, nous énonçons la forme précise du théorème~\ref{ch3introthdp4};
sa preuve occupe les paragraphes~\ref{ch3parpremierpreuve} à~\ref{ch3pardernierpreuve}.
Dans \cite[§6]{bsd}, Swinnerton-Dyer espérait que sa méthode permette d'établir,
en toute généralité, qu'une surface de del Pezzo de degré~$4$ admet un point rationnel
dès que l'obstruction de Brauer-Manin ne s'y oppose pas
(sous l'hypothèse de Schinzel
et la finitude des groupes de Tate-Shafarevich). À l'aide des résultats du
paragraphe~\ref{ch3parconssd}, nous montrons au paragraphe~\ref{ch3parobstruction}
que cet espoir était \emph{prévisiblement} trop optimiste et qu'à moins d'avancées importantes
dans l'arithmétique des pinceaux de courbes de genre~$1$,
les techniques employées dans ce chapitre ne suffisent pas à l'étude des points rationnels
des surfaces de del Pezzo de degré~$4$ pour lesquelles~\mbox{$\Br(X)/\Br(k)\neq 0$}.

Enfin, les paragraphes~\ref{ch3monodr}
et~\ref{ch3parpreuve5} démontrent le théorème~\ref{ch3introp5}.

\bigskip
Voici maintenant quelques indications sur la preuve du théorème~\ref{ch3introthdp4}.
Soit~$X$ une surface de del Pezzo de degré~$4$ sur~$k$, vue comme intersection
de deux quadriques dans~$\P^4_k$.  Supposons qu'il n'y ait
pas d'obstruction de Brauer-Manin à l'existence d'un point rationnel sur~$X$.
Les résultats du chapitre~2 permettront de conclure quant à l'existence d'un
point rationnel sur~$X$ si l'on trouve une famille de sections hyperplanes
remplissant les conditions suivantes:
\begin{itemize}
\item[(i)] la famille doit être paramétrée par~$\P^1_k$;
\item[(ii)] sa fibre générique doit être une courbe lisse et géométriquement connexe
de genre~$1$ et de période~$2$ dont la
jacobienne admet un point d'ordre~$2$ rationnel;
\item[(iii)] la \condD{} du chapitre~2 (ou la \condE{}, cf.~paragraphe~\ref{ch2enonceresul}) doit être satisfaite;
\item[(iv)] il ne doit pas y avoir d'obstruction de Brauer-Manin verticale à l'existence d'un point rationnel sur l'espace
total de cette famille.
\end{itemize}
Une réduction standard (cf.~proposition~\ref{ch3findp4amerbrumer})
permet de ne s'intéresser qu'au cas où le polynôme~$f$
admet une racine $k$\nobreakdash-rationnelle, compte tenu qu'il est de degré impair.
Sous cette hypothèse, la construction
de Swinnerton-Dyer fournit une famille $\pi \colon Y \rightarrow B$ de sections hyperplanes de~$X$
remplissant la condition~(ii),
où~$B$ est un ouvert de~$\P^3_k$ dont le complémentaire est de codimension~$\geq 2$.
La projection naturelle $Y \rightarrow X$ admet une section rationnelle, de sorte que
la condition~(iv) est automatique.
Bien entendu, la condition~(i) est en défaut et la condition~(iii) n'a donc pas de sens
pour~$\pi$.

Se pose alors la question suivante: existe-t-il une droite $k$\nobreakdash-rationnelle
$D \subset \P^3_k$ incluse dans~$B$, telle que la famille $\pi^{-1}(D) \rightarrow D$
satisfasse aux conditions~(i) à~(iv)~? Pour les conditions~(i) et~(ii), il suffit
que~$D$ soit choisie hors d'un certain fermé de codimension~$1$
de l'espace des droites de~$\P^3_k$.  Voici comment forcer la condition~(iv) à être
satisfaite.  D'après les résultats du paragraphe~\ref{ch3parobmvert}, il
existe $b_0 \in B(k)$ tel que la fibre $\pi^{-1}(b_0)$ soit lisse et
admette un $k_v$\nobreakdash-point pour
toute place~$v$ de~$k$.  En particulier, si $b_0 \in D$, l'une
des fibres de la famille $\pi^{-1}(D) \rightarrow D$ au-dessus d'un point $k$\nobreakdash-rationnel
contient un $k_v$\nobreakdash-point lisse
pour toute place~$v$; ceci assure
que l'obstruction de Brauer-Manin verticale à l'existence d'un point rationnel sur $\pi^{-1}(D)$
s'évanouit.  On ne considérera donc par la suite que des droites~$D$ passant par~$b_0$.

Reste seulement la condition~(iii).
Au paragraphe~\ref{ch3parspecd}, nous définissons une version de la \condD{} adaptée à la fibration~$\pi$
tout entière; nous l'appelons «~\condD{} générique~».  Elle admet une formulation
dans laquelle interviennent:
la géométrie du lieu $\Delta \subset B$ des fibres singulières de~$\pi$,
la structure galoisienne de~$\Delta$ et enfin la structure galoisienne des fibres de~$\pi$ au-dessus des points
génériques des composantes irréductibles de~$\Delta$.
(L'expression «~structure galoisienne~» désigne la collection des objets que l'on peut associer à une variété (ou à un morphisme de variétés)
indépendamment de la nature du corps de base, mais qui sont triviaux lorsque
le corps de base est séparablement clos.  En l'occurrence, c'est d'une manière très délicate que la structure galoisienne
de~$\Delta$ apparaît dans la \condD{} générique;
il ne s'agit pas
simplement de l'action de $\Gal(\ksep/k)$ sur l'ensemble des composantes irréductibles
de $\Delta \otimes_k \ksep$ ni de la structure galoisienne de chaque composante irréductible
de~$\Delta$ considérée abstraitement comme $k$\nobreakdash-variété.)
Nous montrons ensuite
que si la \condD{} générique est satisfaite, il est possible de trouver des
droites $D\subset B$ passant par~$b_0$ pour lesquelles la fibration $\pi^{-1}(D) \rightarrow D$
vérifie la \condE{}.  L'argument employé nous oblige à nous contenter de la \condE{}
au lieu de la \condD{}, pour des raisons de finitude.  Rappelons que la \condE{} dépend d'un
point rationnel~$x$ de~$D$ et d'un ensemble fini~$S$ de places de~$k$. Il est absolument essentiel
ici que l'on puisse prévoir la taille de cet ensemble fini \emph{avant} de choisir~$D$.
Une telle prédiction est possible si l'on pose $x=b_0$
car l'ensemble~$S$ ne dépend que des propriétés de la fibre de~$\pi$
au-dessus de~$x$.  Ce sont ces considérations qui nous ont conduit à formuler
le théorème~\ref{ch2thsurmesure} tel qu'il est formulé,
et non en termes d'obstruction de Brauer-Manin
verticale;  il est à noter qu'une fois le point~$b_0$ fixé,
le «~lemme formel~» n'est plus jamais utilisé.

Ainsi suffit-il, pour démontrer le théorème~\ref{ch3introthdp4}, de traduire
la \condD{} générique en termes des propriétés de la surface de del Pezzo~$X$;
et c'est là la partie la plus délicate de l'argument.  En effet,
il s'avère qu'après extension des scalaires de~$k$ à une clôture algébrique~$\ksep$,
la \condD{} générique n'est \emph{jamais} satisfaite (bien que cette remarque n'ait rien
d'évident sur la définition de la \condD{} générique). Si l'on espère trouver des conditions suffisantes sur~$X$ pour que la
\condD{} générique soit satisfaite, il est donc indispensable de tenir compte
de la structure galoisienne de~$\Delta$.  (On voit ainsi que les arguments de~\cite[§6]{bsd}
ne pouvaient être corrects, étant purement géométriques.)

Au paragraphe~\ref{ch3parconssd}, nous formulons la construction de Swinnerton-Dyer
en termes abstraits, à l'aide uniquement des propriétés des pinceaux de quadriques
dans~$\P^n$.  Cette approche nous permet d'interpréter
les diverses composantes irréductibles de $\Delta \otimes_k \ksep$ (avec l'action
du groupe~$\Gal(\ksep/k)$).  Elle se révèle surtout particulièrement
efficace pour déterminer la structure des fibres de~$\pi$ au-dessus des points
génériques des composantes irréductibles de~$\Delta$ (notamment, elle
explique l'identité «~$\beta_{\Delta_5}=\epsilon_0$~» du lemme~\ref{ch3verifidentbeta},
qui peut sembler assez mystérieuse selon la manière dont on y arrive).
Néanmoins, les invariants
galoisiens fins de~$\Delta$ qui interviennent dans la \condD{} générique semblent actuellement
hors de portée de telles considérations.  Pour en donner une idée,
même les questions géométriques les plus simples au sujet de~$\Delta$ (par exemple, quels
sont
les degrés des composantes irréductibles de~$\Delta \otimes_k \ksep$~?)
demanderaient encore du travail pour être résolues par voie abstraite.

Il est possible d'étudier la géométrie de~$\Delta$ au moyen de calculs explicites,
car sur~$\ksep$, on peut toujours diagonaliser simultanément
les formes quadratiques qui s'annulent sur~$X$.
Pour ce qui est des questions tenant à la structure galoisienne de~$\Delta$, la situation est moins favorable:
il est totalement inconcevable de
travailler avec des équations explicites sur~$k$,
tant elles sont complexes.

Ce n'est qu'en combinant les trois techniques suivantes que nous parviendrons
à analyser la~\condD{} générique précisément.  Tout d'abord, pour certaines
questions, un argument de descente galoisienne est possible; on les résout
sur~$\ksep$ par le calcul, ce qui revient souvent à exhiber
des identités polynomiales à coefficients entiers
non évidentes, puis on suit l'action de $\Gal(\ksep/k)$
sur les formules obtenues.  Parmi les autres questions, celles qui concernent
les fibres singulières de~$\pi$ en codimension~$1$ sont traitées par voie abstraite.
Enfin, pour les questions restantes, des réductions permettent de supposer que le polynôme~$f$
admet trois racines $k$\nobreakdash-rationnelles, auquel cas certains calculs restent envisageables.

Il convient de signaler que la \condD{} générique se simplifie dans le cas où
le sous-groupe $G \subset \mathfrak{S}_5$ est $3$\nobreakdash-transitif, de sorte qu'une analyse
plus sommaire (mais devant toujours être effectuée sur~$k$ et non sur~$\ksep$) suffit si l'on ne s'intéresse
qu'au théorème~\ref{ch3introp5}.

\bigskip
\noindent\textbf{Notations.}\\
\noindent{}Si~$k$ est un corps de nombres, on note~$\Omega$ l'ensemble de
ses places et l'on pose $k_\Omega = \prod_{v \in \Omega} k_v$\glossary{$\Br_\nr(X)$, $k_\Omega$,
$X(k_\Omega)^{\Br_\nr}$}.
Si~$X$ est
une variété lisse et géométriquement intègre sur~$k$, on note $\Br_\nr(X)$
le groupe de Brauer\index{groupe de Brauer!non ramifié} non ramifié de~$X$ sur~$k$ et $X(k_\Omega)^{\Br_\nr}$
l'ensemble des $(\fp{P}v)_{v \in \Omega} \in X(k_\Omega)=\prod_{v \in \Omega}X(k_v)$
orthogonaux à $\Br_\nr(X)$, c'est\nobreakdash-à-\-dire tels que
$\sum_{v \in \Omega} \inv_v A(\fp{P}v)=0$ pour tout $A \in \Br_\nr(X)$
(cf.~\cite[§5.2]{skotorsors}).  Enfin, si l'on dispose de plus d'un
morphisme $X \rightarrow Y$ où~$Y$ est une variété intègre, on note
$\Br_{\nr,\vert}(X)=\Br_\nr(X) \cap \Br_\vert(X)$\glossary{$\Br_{\nr,\vert}(X)$, $X(k_\Omega)^{\Br_{\nr,\vert}}$}
et l'on définit $X(k_\Omega)^{\Br_{\nr,\vert}}$ de la manière évidente.

\sectionnotoc[Obstruction de Brauer-Manin verticale]{Obstruction de Brauer-Manin verticale et points adéliques dans les fibres d'un morphisme vers~$\P^n$}
\addcontentsline{toc}{section}{\protect\numberline{\thesection}Obstruction de Brauer-Manin verticale et points adéliques dans les fibres d'un morphisme vers~$\P^n$}
\label{ch3parobmvert}

L'objet de ce paragraphe est de démontrer le théorème suivant\index{groupe de Brauer!vertical}, qui généralise \cite[Theorem~1.1]{csscrelle98}.

\bigskip
\begin{theoreme}
\label{ch3thpointslocaux}
Soient un corps de nombres~$k$, un entier $n\geq 1$,
un fermé $E \subset \P^n_k$ de codimension $\geq 2$ et
une $k$\nobreakdash-variété géométriquement intègre~$X$
munie d'un morphisme projectif $f \colon X \rightarrow \P^n_k \setminus E$
dont la fibre générique est géométriquement intègre.
On note $X^0 \subset X$ l'ouvert de lissité de~$f$ et $X^\reg \subset X$ l'ouvert
des points réguliers de~$X$.

Supposons que pour tout point $m \in \P^n_k$ de codimension~$1$, la fibre~$X_m$
possède une composante irréductible~$\ft{Y}m$ de multiplicité~$1$ telle que la fermeture
algébrique de~$\kappa(m)$ dans~$\kappa(\fy{Y}m)$ soit une extension abélienne
de~$\kappa(m)$.

Soient~$U$ un ouvert dense de~$\P^n_k$,
$$\Rrond = \bigensemble{x \in (\P^n_k \setminus E)(k)}{X_x\text{ est géométriquement
intègre et }X_x^0(k_\Omega)\neq\emptyset}$$
et $\overline{\Rrond \cap U(k)}$ l'adhérence de~$\Rrond \cap U(k)$ dans~$\P^n(k_\Omega)$.
Admettons l'hypothèse de Schinzel. Alors
$$
f(X^\reg(k_\Omega)^{\Br_{\nr,\vert}}) \subset \overline{\Rrond \cap U(k)} \rlap{\text{.}}
$$
En particulier, s'il n'y a pas d'obstruction de Brauer-Manin verticale à
l'existence d'un point rationnel sur un modèle propre et lisse de~$X$,
l'ensemble~$\Rrond$ est Zariski-dense dans~$\P^n_k$.
\end{theoreme}

\bigskip
Le principe de la preuve du théorème~\ref{ch3thpointslocaux} est de se réduire au cas où~$X$ est une surface projective et lisse et où $n=1$,
qui est exactement la situation dans laquelle \cite[Theorem~1.1]{csscrelle98} s'applique.  Ainsi nous ne reprouvons pas
\cite[Theorem~1.1]{csscrelle98} mais l'invoquons à la fin de la démonstration.

Voici quelques corollaires d'intérêt général du théorème~\ref{ch3thpointslocaux}.

\bigskip
\begin{corollaire}
\label{ch3premcorthpl}
Soient~$X$ une variété projective, lisse et connexe sur un corps de nombres~$k$ et $f \colon X \rightarrow \P^n_k$ un morphisme de fibre générique géométriquement intègre.
Supposons que les fibres lisses de~$f$ au-dessus des points rationnels d'un ouvert dense de~$\P^n_k$ satisfassent à l'approximation faible
et que pour tout point $m \in \P^n_k$ de codimension~$1$, la fibre de~$f$ en~$m$
possède une composante irréductible de multiplicité~$1$ dans le corps des fonctions de laquelle la fermeture
algébrique de~$\kappa(m)$ soit une extension abélienne
de~$\kappa(m)$.
Admettons l'hypothèse de Schinzel.  Alors l'adhérence de~$X(k)$ dans $X(\A_k)$ est égale
à $X(\A_k)^{\Br_\vert}$.
\end{corollaire}

\bigskip
\begin{demo}
Notons~$\Omega$ l'ensemble des places de~$k$ et fixons une famille
$(\fp{P}v)_{v \in \Omega} \in X(\A_k)^{\Br_\vert}$, un ensemble fini $S\subset \Omega$ et un voisinage~$\Arond_v$
de~$\fp{P}v$ dans~$X(k_v)$ pour $v \in S$.  Nous allons exhiber un point rationnel de~$X$ dont l'image dans $\prod_{v \in S} X(k_v)$
appartienne à $\prod_{v \in S}\Arond_v$.
Le théorème des fonctions implicites, la continuité de l'évaluation des classes de~$\Br(X)$ sur~$X(k_v)$
et la finitude de $\Br_\vert(X)/\Br(k)$ (établie dans \cite[Lemma~3.1]{ctskodescenteouverte})
permettent de supposer,
quitte à remplacer $\fp{P}v$ par un autre point de~$\Arond_v$, que~$f$ est lisse en~$\fp{P}v$ pour tout $v \in S$.
Le théorème des fonctions implicites fournit alors, pour $v \in S$, un voisinage $\Urond_v$ de $f(\fp{P}v)$ dans $\P^n(k_v)$ et une section
analytique locale $\sigma_v \colon \Urond_v \rightarrow \Arond_v$ de l'application $\Arond_v \rightarrow \P^n(k_v)$ induite par~$f$.
Soit $U \subset \P^n_k$ un ouvert dense tel que pour tout $u \in U(k)$, la fibre~$X_u$ soit lisse et satisfasse à l'approximation faible.
D'après le théorème~\ref{ch3thpointslocaux}, il existe $u \in U(k)$ vérifiant $X_u(\A_k)\neq\emptyset$
et dont l'image dans $\prod_{v \in S}\P^n(k_v)$ appartienne à $\prod_{v \in S}\Urond_v$.
Comme~$X_u$ satisfait à l'approximation faible et que $X_u(\A_k)\neq\emptyset$, il existe $x \in X_u(k)$ arbitrairement proche de $\sigma_v(u)$ pour $v \in S$.
Si~$x$ est choisi suffisamment proche de $\sigma_v(u)$ pour $v \in S$, il appartiendra alors aux voisinages~$\Arond_v$,
puisque $\sigma_v(u)\in \Arond_v$.
\end{demo}

\bigskip
\begin{corollaire}
\label{ch3corsbgen}
Soit~$X$ une variété projective, lisse et connexe sur un corps de nombres~$k$.  Soient $n \geq 1$ et $f \colon X \rightarrow \P^n_k$ un morphisme dont la fibre générique
est une variété de Severi-Brauer\index{variété de Severi-Brauer}
généralisée au sens de~\cite[§2]{cscrelle94} (par exemple une conique, une quadrique de dimension~$2$ ou une variété de Severi-Brauer).
Admettons l'hypothèse de Schinzel.  Alors l'adhérence de~$X(k)$ dans $X(\A_k)$ est égale à $X(\A_k)^{\Br}$.
En particulier $X$ satisfait au principe de Hasse (resp.~à l'approximation faible) si l'obstruction de Brauer-Manin ne s'y oppose pas.
De plus $X$ satisfait à l'approximation faible faible.
\end{corollaire}

\bigskip
Rappelons qu'une variété~$X$ sur un corps de nombres~$k$ satisfait à
\emph{l'approximation faible faible} s'il existe un ensemble fini $S \subset
\Omega$ tel que l'image de
l'application diagonale \smash{$X(k) \rightarrow
\prod_{v \in \Omega \setminus S}X(k_v)$} soit dense, où~$\Omega$ désigne
l'ensemble des places de~$k$.

Rappelons d'autre part la notion de variété de
Severi-Brauer généralisée de~\cite[§2]{cscrelle94}: une variété~$X$ sur un
corps~$k$ est une \emph{variété de Severi-Brauer généralisée} s'il existe un
ensemble fini~$I$, une famille~$(k_i)_{i \in I}$ d'extensions finies
séparables de~$k$ et une famille $(V_i)_{i \in I}$ de variétés de
Severi-Brauer sur~$k$ telles que~$X$ soit isomorphe au produit des
restrictions des scalaires à la Weil à~$k$ des $k_i$\nobreakdash-variétés $V_i
\otimes_k k_i$ pour $i \in I$.

Le corollaire~\ref{ch3corsbgen}, dans le cas particulier où $n=1$, n'est autre que \cite[Theorem~4.2]{cscrelle94}.

\bigskip
\begin{demo}
Si $X(k)$ est dense dans $X(\A_k)^{\Br_\vert}$, il résulte de la finitude du
groupe $\Br_\vert(X)/\Br(k)$, démontrée dans \cite[Lemma~3.1]{ctskodescenteouverte},
que~$X$ satisfait à l'approximation faible faible.
Pour établir le corollaire~\ref{ch3corsbgen}, il suffit donc de prouver que l'hypothèse sur les fibres de codimension~$1$ de~$f$ dans le corollaire~\ref{ch3premcorthpl} est vérifiée.
Soit $m \in \P^n_k$ un point de codimension~$1$.  Le complété de l'anneau local de~$\P^n_k$ en~$m$ est isomorphe à $\kappa(m)[[t]]$ (théorème de Cohen, cf.~\cite[Theorem~5.5A]{hartshorne}).
Considérant l'image réciproque de~$f$ au-dessus du spectre de cet anneau, on voit alors que pour conclure,
il suffit de prouver que si~$\ell$ est un corps de caractéristique~$0$ et $V$ est un $\ell[[t]]$\nobreakdash-schéma propre et régulier
dont la fibre générique est une variété de Severi-Brauer généralisée, la fibre spéciale de~$V$ possède
une composante irréductible de multiplicité~$1$ dans le corps des fonctions de laquelle la fermeture
algébrique de~$\ell$ soit une extension abélienne
de~$\ell$.  La formation de la fermeture algébrique commute aux extensions régulières.
Quitte à remplacer~$\ell$ par la limite inductive des corps de fonctions de toutes les variétés géométriquement intègres sur~$\ell$, on peut donc
supposer~$\ell$ pseudo-algébriquement clos (c'est\nobreakdash-à-dire tel que toute variété géométriquement intègre sur~$\ell$ admette un point rationnel).
Soient $\uplet{A_1}{A_s} \in \Br(\ell((t)))$ les classes des algèbres simples centrales associées à la fibre générique de~$V$.
Pour $i \in \{\uplet{1}{s}\}$, le résidu de~$A_i$ est un élément de $H^1(\ell,\Q/\Z)$ et définit donc une extension cyclique $\ell_i/\ell$.  Soit $\ell'/\ell$ l'extension composée
des $\ell_i/\ell$ et soit~$A'_i$, pour $i \in \{\uplet{1}{s}\}$, l'image de~$A_i$ dans $\Br(\ell'((t)))$.
Par construction de~$\ell'$, le résidu de~$A'_i$ est nul (cf.~\cite[Proposition~1.1.1]{cscrelle94}); d'où $A'_i \in \Br(\ell'[[t]])$.
Le groupe $\Br(\ell'[[t]])$ est isomorphe à $\Br(\ell')$ (cf.~\cite[Corollaire~6.2]{grothbr1}), or celui-ci est nul car
la restriction des scalaires à la Weil à~$\ell$ d'une variété de Severi-Brauer sur~$\ell'$ possède nécessairement un point rationnel, le corps~$\ell$
étant pseudo-algébriquement clos.  D'où $A'_i=0$ pour tout $i \in \{\uplet{1}{s}\}$.
En particulier la fibre générique de~$V$ admet-elle un $\ell'((t))$\nobreakdash-point.  Comme~$V$ est propre sur $\ell[[t]]$ et régulier, il s'ensuit que la fibre spéciale de~$V$
admet un $\ell'$\nobreakdash-point lisse, d'où le résultat recherché puisque $\ell'/\ell$ est une extension abélienne.
\end{demo}

\bigskip
Les variétés fibrées en variétés de Severi-Brauer généralisées au-dessus de la droite projective, étudiées dans~\cite{cscrelle94},
sont des variétés rationnelles.  Il~n'en est pas de même, en général, des variétés fibrées en variétés de Severi-Brauer généralisées
au-dessus de l'espace projectif de dimension~$\geq 2$.
Il arrive notamment que le groupe de Brauer d'une telle variété comporte des classes transcendantes\index{groupe de Brauer!transcendant} (c'est-à-dire non algébriques).
Harari a montré que ces classes transcendantes peuvent jouer un rôle arithmétique; plus précisément,
il a exhibé une variété de dimension~$3$, fibrée en coniques
sur~$\P^2_\Q$, qui viole le principe de Hasse à cause d'une obstruction de Brauer-Manin créée par une classe transcendante du
groupe de Brauer (cf.~\cite{hartransc}).
Il y a donc une différence qualitative entre \cite[Theorem~1.1]{csscrelle98} et le théorème~\ref{ch3thpointslocaux}:
contrairement à \cite[Theorem~1.1]{csscrelle98},
le théorème~\ref{ch3thpointslocaux} s'applique dans des situations
où des classes transcendantes existent dans le groupe de Brauer et ne peuvent être ignorées.
Le corollaire~\ref{ch3corsbgen} est ainsi le premier résultat positif d'existence de points rationnels (ou d'approximation faible)
dont la preuve prenne effectivement en compte les classes transcendantes du groupe de Brauer.

Le reste de ce paragraphe est consacré à la démonstration du théorème~\ref{ch3thpointslocaux}.

\bigskip
\begin{demo}[ du théorème~\ref{ch3thpointslocaux}]%
Nous commençons par une réduction au cas projectif et lisse.

\bigskip
\begin{proposition}
\label{ch3propreduc}
Soit~$X'$ une $k$\nobreakdash-variété projective, lisse et géométriquement intègre munie d'un
morphisme $f' \colon X' \rightarrow \P^n_k$.
On suppose qu'il existe une application birationnelle $b \colon X \dashrightarrow X'$
telle que le diagramme
$$
\myxyhook\xymatrix{
X \ar@{-->}[r]^b \ar[d]^(.45)f & X' \ar[d]^(.45){f'} \\
\P^n_k \setminus E \ar@{^{ (}->}[r] & \P^n_k
}
$$
commute.  Alors pour tout point $m \in \P^n_k$ de codimension~$1$, la
fibre~$X'_m$ possède une composante irréductible~$\fy{Y}m$ de multiplicité~$1$ telle que
la fermeture algébrique de~$\kappa(m)$ dans~$\kappa(\fy{Y}m)$ soit une extension abélienne
de~$\kappa(m)$.  Notons~$X'^0$ l'ouvert de lissité de~$f'$ et
$$
\Rrond' = \bigensemble{x \in \P^n(k)}{X'_x\text{ est géométriquement
intègre et }X'^0_x(k_\Omega)\neq\emptyset}\!\rlap{\text{.}}
$$
Supposons que
l'on ait $f'(X'(\A_k)^{\Br_\vert}) \subset {\overline{\Rrond'\cap U(k)}}$
pour tout ouvert dense $U \subset \P^n_k$.
Alors on a l'inclusion $f(X^\reg(k_\Omega)^{\Br_{\nr,\vert}}) \subset \overline{\Rrond \cap U(k)}$
pour tout ouvert dense $U \subset \P^n_k$.
\end{proposition}

\bigskip
\begin{demo}
La première assertion résulte du lemme suivant.  Remarquons que dans la situation
de la proposition, le morphisme~$f'$ est dominant et sa fibre générique est automatiquement
lisse.

\bigskip
\begin{lemme}
\label{ch3lemmeinspi}
Soient $f \colon X \rightarrow S$ et $f' \colon X' \rightarrow S$ des morphismes
de schémas
dominants et de type fini avec~$S$ et~$X$ localement noethériens intègres, $X'$ régulier intègre,
$f'$ propre et de fibre générique lisse.
Soit~$s \in S$ un point de codimension~$1$ en lequel~$S$ est régulier.
Supposons qu'il existe une $S$\nobreakdash-application rationnelle de~$X$
vers~$X'$.
Si la fibre de~$f$ en~$s$ possède une composante irréductible de multiplicité
géométrique~$1$ (\emph{i.e.}~telle que la fibre soit géométriquement réduite au point générique
de cette composante)
dans le corps des fonctions de laquelle la fermeture algébrique de~$\kappa(s)$ est une
extension abélienne de~$\kappa(s)$,
alors il en va de même pour la fibre de~$f'$ en~$s$.
\end{lemme}

\bigskip
\begin{demo}
On peut supposer que~$S$ est un trait.
Les morphismes~$f$ et~$f'$ sont alors
plats, en vertu de \cite[III, 9.7]{hartshorne}.
Soit $Y \subset X_s$ une composante irréductible de la fibre spéciale de~$f$
satisfaisant à la condition de l'énoncé.
L'hypothèse de multiplicité géométrique~$1$ entraîne que~$X_s$ est lisse
au point générique~$y$ de~$Y$. Le morphisme~$f$, étant plat,
est donc lisse en~$y$.  Par conséquent, le schéma~$X$ est régulier en~$y$
et l'anneau local $\Orond_{X,y}$, de dimension~$1$ puisque~$f$ est plat
et que~$s$ est de codimension~$1$ dans~$S$, est donc un anneau de valuation discrète.
Notons $T=\Spec(\Orond_{X,y})$ et $U=\Spec(\kappa(X))$.
Par hypothèse, il existe un $S$\nobreakdash-morphisme $U \rightarrow X'$.
Le critère valuatif de propreté permet d'en déduire l'existence
d'un $S$\nobreakdash-morphisme $T \rightarrow X'$.
Celui-ci est nécessairement à valeurs dans l'ouvert de lissité de~$f'$, d'après
le lemme~\cite[3.6/5]{blr} appliqué à l'identité de~$X'$.
La fibre spéciale de~$f'$ contient donc un point lisse, disons~$c$,
dont le corps résiduel se plonge $\kappa(s)$\nobreakdash-linéairement dans~$\kappa(y)$.
Ce point appartient à une unique composante irréductible de~$X'_s$.
Notons-la~$Y'$. Elle est de multiplicité géométrique~$1$ puisque
l'ouvert de lissité de~$X'_s$, contenant~$c$, contient le point générique de~$Y'$.
Soit~$\ell$ la fermeture algébrique de~$\kappa(s)$ dans~$\kappa(Y')$.
Comme~$Y'$ est normal en~$c$, il existe un plongement $\kappa(s)$\nobreakdash-linéaire
de~$\ell$ dans $\Orond_{Y',c}$, donc dans~$\kappa(c)$, donc dans~$\kappa(y)$.
Il en résulte que l'extension $\ell/\kappa(s)$ est abélienne, ce qui termine de
prouver le lemme.
\end{demo}

\bigskip
\begin{remarque}
Le lemme~\ref{ch3lemmeinspi} et sa démonstration sont directement inspirés d'une preuve de Skorobogatov
(cf.~\cite[Lemma~1.1]{skodescent}).
\end{remarque}

\bigskip
Intéressons-nous maintenant à la seconde assertion. Supposons que
\begin{equation}
\label{ch3reduchypeq}
f'(X'(\A_k)^{\Br_\vert}) \subset {\overline{\Rrond'\cap U(k)}}
\end{equation}
pour tout ouvert dense $U \subset \P^n_k$.
Choisissons une famille $(\fp{P}v)_{v\in\Omega}\in X^\reg(k_\Omega)$
orthogonale à $\Br_{\nr,\vert}(X^\reg)$,
un ensemble fini~$S$ de places de~$k$,
un voisinage $v$\nobreakdash-adique $\Arond_v \subset \P^n(k_v)$ de~$f(\fp{P}v)$
pour chaque $v \in S$ et un ouvert dense $U \subset \P^n_k$.
On va maintenant prouver l'existence d'un $x \in \Rrond \cap U(k)$ appartenant à~$\Arond_v$
pour tout $v \in S$.

Soit $Y \subset X$ un ouvert dense sur lequel l'application rationnelle~$b$ soit définie
et induise un isomorphisme vers un ouvert~$Y'$ de~$X'$.
Quitte à rétrécir~$U$, on peut supposer que pour tout $x \in U$,
les conditions suivantes sont satisfaites:
l'ouvert~$\fy{Y'}x$ est dense dans~$X'_x$ (cf.~\cite[9.5.3]{ega43}), la variété~$X_x$ est géométriquement intègre sur~$k$
(cf.~\cite[9.7.7~(iv)]{ega43})
et la variété~$X'_x$ est lisse sur~$k$ (cf.~\cite[17.7.11]{ega44}; la fibre générique de~$f'$ est lisse car~$X'$
est lisse sur~$k$).
Notons que l'ouvert~$U$ est alors automatiquement disjoint de~$E$.

Le lemme suivant est une variation sur le lemme de\index{lemme de Nishimura} Nishimura.

\bigskip
\begin{lemme}
\label{ch3nishi}
Soient~$S$ un schéma, $X$ et~$Y$ deux $S$\nobreakdash-schémas
et \mbox{$f \colon X \dashrightarrow Y$} une $S$\nobreakdash-application rationnelle.
Supposons~$Y$ propre sur~$S$ et~$X$ localement noethérien régulier.
Soient~$s \in S$ et~$x$ un point rationnel de la fibre~$X_s$.
Alors il existe un point rationnel~$y$ de la fibre~$\fy{Y}s$
tel que l'égalité
$A(y)=(f^\star A)(x)$ d'éléments de~$\Br(\kappa(s))$ ait lieu pour
tout $A \in \Br(Y)$ tel que $f^\star A \in \Br(X)$.
\end{lemme}

\bigskip
Précisons le sens de l'énoncé: si~$U$ est un ouvert dense de~$X$ sur lequel~$f$
est définie, on dispose d'une classe $f^\star A \in \Br(U)$
pour tout $A \in \Br(Y)$.  Le morphisme de restriction $\Br(X) \rightarrow \Br(U)$
est une injection puisque~$X$ est régulier. On peut donc considérer~$\Br(X)$ comme
un sous-groupe de~$\Br(U)$ et évaluer $f^\star A$ en un point de~$X$
si $f^\star A$ appartient à ce sous-groupe.
Il est à noter que le point~$s$ n'est soumis à aucune condition; par exemple,
la fibre~$X_s$ pourrait tout à fait être disjointe du domaine de définition de~$f$.

\bigskip
\begin{demo}
Soit $U \subset X$ un ouvert dense sur lequel~$f$ soit définie.  Notons
encore~$f$ le $S$\nobreakdash-morphisme $U \rightarrow X$ déduit de~$f$.
Comme~$X$ est régulier, il existe un trait~$T$ et un morphisme $\phi \colon T \rightarrow X$
tels que $\phi(\eta)\in U$, $\phi(t)=x$ et que l'application $\kappa(x) \rightarrow \kappa(t)$
induite par~$\phi$ soit un isomorphisme,
où~$t$ et~$\eta$ désignent respectivement le point fermé et le point générique de~$T$.
(En effet, quitte à rétrécir~$U$, on peut supposer que $x \not\in U$ et que
chaque composante irréductible de $X \setminus U$ est de codimension~$1$ dans~$X$,
auquel cas $X \setminus U$ est un diviseur de Cartier, puisque~$X$ est régulier.
Une équation locale de $X \setminus U$ au voisinage de~$x$ détermine alors un vecteur
non nul $f_1 \in \mgoth_x/\mgoth_x^2$, que l'on peut
compléter en une base $(\uplet{f_1}{f_n})$; relevant arbitrairement $\uplet{f_2}{f_n}$
en $\uplet{e_2}{e_n} \in \mgoth_x$, on pose alors $T=\Spec(\Orond_x/(\uplet{e_2}{e_n}))$.
C'est un trait en vertu de \cite[0.16.5.6 et 0.17.1.7]{ega41}
et le morphisme canonique $\phi\colon T \rightarrow X$ remplit
bien les conditions voulues.) L'application rationnelle $f \circ \phi \colon T \dashrightarrow Y$
étant une $S$\nobreakdash-application rationnelle, elle est définie partout, compte
tenu que~$Y$ est propre sur~$S$ et que~$T$ est un trait (cf.~\cite[7.3.8]{ega2}).
Notons $\psi \colon T \rightarrow Y$ ce morphisme et posons $y=\psi(t)$;
c'est un point rationnel de~$\fy{Y}s$.  Soit $A \in \Br(Y)$ tel que $f^\star A \in \Br(X)$, de sorte
que l'on dispose d'une classe $\phi^\star f^\star A \in \Br(T)$.  Étant donné
que $(\phi^\star f^\star A)(t)=(f^\star A)(x)$ et que $(\psi^\star A)(t)=A(y)$,
il suffit, pour conclure, de prouver l'égalité $\phi^\star f^\star A=\psi^\star A$
dans $\Br(T)$.  Comme~$T$ est régulier, la flèche $j^\star \colon \Br(T) \rightarrow \Br(\eta)$
induite par l'inclusion $j \colon \eta \rightarrow T$
est injective.  Il suffit donc d'établir
que $j^\star \phi^\star f^\star A=j^\star \psi^\star A$,
c'est-à-dire que $(\phi\circ j)^\star f^\star A = (\psi \circ j)^\star A$; or ceci résulte
de ce que le triangle
$$
\xymatrix{
& \eta \ar[dr]^{\psi \circ j} \ar[dl]_{\phi \circ j} \\
U \ar[rr]^f && Y
}
$$
commute (par définition de~$\psi$).
\end{demo}

\bigskip
Appliquant le lemme~\ref{ch3nishi} à la $\P^n_{k_v}$\nobreakdash-application
rationnelle
$$X^\reg \otimes_k k_v \longdashrightarrow X' \otimes_k k_v$$
induite par~$b$
et au point $\fp{P}v \in X^\reg(k_v)$ pour chaque place $v \in \Omega$,
on obtient un point adélique $(Q_v)_{v\in\Omega} \in X'(\A_k)$ qui, par définition
du groupe de Brauer non ramifié, est orthogonal au groupe de Brauer vertical de~$f'$.
D'après l'inclusion~(\ref{ch3reduchypeq}), et compte tenu que $f'(Q_v)=f(\fp{P}v)$ pour tout $v \in \Omega$,
il existe donc $x \in \Rrond' \cap U(k)$ appartenant à~$\Arond_v$
pour tout $v \in S$.
Comme~$\fy{Y'}x$ est un ouvert dense de la variété lisse~$X'_x$, le théorème des fonctions implicites montre que $\fy{Y'}x(k_\Omega)\neq\emptyset$,
d'où il résulte que $X^0_x(k_\Omega)\neq\emptyset$.
La variété~$X_x$ étant géométriquement intègre sur~$k$, cela signifie que $x \in \Rrond$; ainsi la
proposition~\ref{ch3propreduc} est-elle démontrée.
\end{demo}

\bigskip
Il existe bien une variété~$X'$ et un morphisme~$f'$ vérifiant les conditions
de la proposition~\ref{ch3propreduc}.  On peut en effet choisir une compactification
projective arbitraire~$X'_0$ de la variété quasi-projective~$X$,
considérer l'adhérence~$X'_1$ dans $X'_0 \times_k \P^n_k$ du graphe de~$f$
puis trouver une variété projective, lisse et géométriquement intègre~$X'_2$
munie d'un morphisme birationnel $\sigma \colon X'_2 \rightarrow X'_1$, grâce à Hironaka.
Il suffit alors de poser $X'=X'_2$ et de prendre pour~$f'$ la composée de~$\sigma$
et de la projection $X'_1 \rightarrow \P^n_k$.
Pour démontrer le théorème, vu la proposition~\ref{ch3propreduc}, il est donc loisible de supposer
le fermé~$E$ vide et la variété~$X$ projective et lisse. Nous nous plaçons dorénavant sous
ces hypothèses.  Nous supposons de plus que $n \geq 2$, car si $n=1$, il suffit
maintenant d'appliquer \cite[Theorem~1.1]{csscrelle98} pour conclure.

\bigskip
La fibre générique de~$f$ est lisse, puisque~$X$ l'est, et elle est géométriquement intègre.
Quitte à rétrécir l'ouvert dense $U \subset \P^n_k$ apparaissant dans l'énoncé du théorème,
on peut donc supposer que toutes les fibres de~$f$ au-dessus de~$U$
sont lisses et géométriquement intègres.
Nous allons fixer un point $c \in U(k)$,
mais avant de préciser son choix, introduisons quelques notations qui en dépendent.
Soient~$\Delta$ la $k$\nobreakdash-variété des droites de~$\P^n_k$ passant par~$c$
et $\pi \colon D \rightarrow \Delta$ le fibré en droites projectives tautologique:
il existe un morphisme canonique $\sigma \colon D \rightarrow \P^n_k$ rendant le triangle
$$
\xymatrix{
& D \ar[dl]_\sigma \ar[dr]^\pi \\
\P^n_k \ar@{-->}[rr] && \Delta
}
$$
commutatif, où la flèche horizontale est l'application rationnelle associant à un point~$x$
de $\P^n_k \setminus\{c\}$ la droite contenant~$c$ et~$x$. Le morphisme~$\sigma$
n'est autre que l'éclatement de~$\P^n_k$ de centre~$c$.

\bigskip
\begin{proposition}
\label{ch3choixdec}
Quitte à bien choisir $c\in U(k)$, on peut supposer, pour prouver le théorème,
que l'application rationnelle $f' \colon X \dashrightarrow D$ induite
par~$\sigma^{-1} \circ f$ est un morphisme et
que la fibre générique du morphisme composé $\pi \circ f' \colon X \rightarrow \Delta$
est lisse et géométriquement intègre.
\end{proposition}

\bigskip
\begin{demo}
Notons~$\Gamma$ la grassmannienne des droites de~$\P^n_k$ et~$I$ la variété
d'incidence
$$
I = \bigensemble{(a,d) \in \P^n_k \times_k \Gamma}{a \in d}\!\rlap{\text{.}}
$$
Posons de plus
$$
X' = \bigensemble{(x,d) \in X \times_k \Gamma}{f(x) \in d}\!\rlap{\text{.}}
$$
Comme~$X$ est géométriquement irréductible et~$f$ dominant,
un théorème de Bertini (cf.~\cite[I, Th.~6.10]{jouanolou}) assure
que la fibre générique de la projection $X' \rightarrow \Gamma$
est géométriquement irréductible.
Il en va donc de même de la fibre
générique de la projection $X' \times_\Gamma I \rightarrow I$,
d'où l'existence d'un ouvert dense de~$I$ au-dessus duquel
les fibres de $X' \times_\Gamma I \rightarrow I$ sont toutes géométriquement
intègres.
Notant $p \colon I \rightarrow \P^n_k$ la projection canonique, on en déduit
l'existence d'un
ouvert dense $V \subset \P^n_k$ tel que la fibre de $X' \times_\Gamma I \rightarrow I$
au-dessus du point générique de $p^{-1}(c)$ soit géométriquement intègre
pour tout $c \in V$.  Quitte à rétrécir~$V$, on peut supposer~$f$ plat
au-dessus de~$V$ (cf.~\cite[6.9.1]{ega42}).

Fixons alors $c \in (U \cap V)(k)$.  La variété~$\Delta$ s'identifiant
à~$p^{-1}(c)$, il est légitime de noter $X' \times_\Gamma \Delta \rightarrow \Delta$
la restriction de $X' \times_\Gamma I \rightarrow I$ au-dessus de~$p^{-1}(c)$.
Celle-ci se factorise naturellement par~$D$, de sorte que l'on obtient le diagramme
commutatif suivant, dans lequel la flèche horizontale supérieure est la projection
canonique et le carré de gauche est cartésien.
$$
\xymatrix@C=8ex{
X \ar[d]^f & X'\times_\Gamma \Delta \ar[l]_{\sigma'} \ar[d] \ar[rd] \\
\P^n_k & D \ar[l]_\sigma \ar[r]^\pi & \Delta
}
$$
Par construction, la fibre générique de la flèche oblique $X' \times_\Gamma \Delta
\rightarrow \Delta$ dans le diagramme ci-dessus est géométriquement intègre.

Les éclatements commutent aux changements de base qui sont plats
au-dessus du sous-schéma fermé que l'on éclate.
Par définition de~$V$, le morphisme~$f$ est plat aux points de~$X_c$.
Comme par ailleurs~$\sigma$ est l'éclatement du point $c \in \P^n_k$
et
que le carré de gauche dans le diagramme ci-dessus est cartésien,
on en déduit que~$\sigma'$ est l'éclatement de la fibre~$X_c$ dans~$X$.
En particulier, étant donné que~$X$ est lisse et géométriquement intègre
et que~$X_c$ est lisse, la variété
$X' \times_\Gamma \Delta$ est elle-même lisse et géométriquement intègre
et $\sigma'^{-1}(X_c)$ est un fibré projectif (localement libre) sur~$X_c$
(cf.~\cite[II, 8.24]{hartshorne}).

La proposition~\ref{ch3propreduc}, appliquée à la variété $X' \times_\Gamma \Delta$,
au morphisme~$f \circ \sigma'$ et à l'application rationnelle $\sigma'^{-1}$, montre
que pour prouver le théorème pour $X \xrightarrow{f} \P^n_k$, il suffit
de l'établir
pour $X' \times_\Gamma \Delta \xrightarrow{f \circ \sigma'} \P^n_k$.
Quitte à remplacer~$X$ par $X' \times_\Gamma \Delta$, on peut donc supposer
que l'application rationnelle $\sigma^{-1}\circ f$ est définie partout
et que la fibre générique du morphisme $\pi \circ \sigma^{-1} \circ f$ est
géométriquement intègre.
L'hypothèse que $c \in U$ est en effet conservée puisque
la variété $(f \circ \sigma')^{-1}(c)$ admet une structure de fibré projectif (localement libre)
sur~$X_c$.  De plus, la fibre générique de $\pi \circ \sigma^{-1} \circ f$ est
automatiquement lisse puisque la variété~$X$ est elle-même lisse.
\end{demo}

\bigskip
Nous fixons dorénavant $c \in U(k)$ au moyen de la proposition~\ref{ch3choixdec} et
supposons les conclusions de cette proposition satisfaites.
Le diagramme commutatif suivant résume la situation.
$$
\xymatrix@C=7ex@R=7ex{
X \ar[r]^{f'} \ar[d]_f & D \ar[dl]_\sigma \ar[d]^\pi \\
\P^n_k \ar@{-->}[r] & \Delta
}
$$
Posons $g = \pi \circ f'$
et notons~$\eta$ le point générique de~$\Delta$.

\bigskip
\begin{lemme}
\label{ch3fggif}
Soit~$\delta\in\Delta$. La fibre générique
du morphisme \mbox{$f'_\delta \colon X_\delta \rightarrow D_\delta$} déduit
de~$f'$ par changement de base est lisse et géométriquement intègre.
\end{lemme}

\bigskip
\begin{demo}
La fibre générique de~$f'_\delta$ est égale à la fibre
de~$f$ au-dessus du point générique de~$\sigma(D_\delta)$.
Celui-ci appartient à~$U$ puisque $c \in U \cap \sigma(D_\delta)$
et que~$U$ est ouvert, ce qui prouve le lemme.
\end{demo}

\bigskip
Considérons maintenant le morphisme de $\kappa(\Delta)$\nobreakdash-variétés
$f'_\eta \colon X_\eta \rightarrow D_\eta$.
C'est un morphisme propre, sa fibre générique est lisse
et géométriquement intègre (cf.~lemme~\ref{ch3fggif}) et~$X_\eta$ est
lui-même géométriquement
intègre (cf.~conclusion de la proposition~\ref{ch3choixdec}).
Comme~$D_\eta$ est une droite projective
sur~$\eta$, ces hypothèses entraînent que~$f'_\eta$ est plat
(cf.~\cite[III, 9.7]{hartshorne}).  L'ensemble~$\Mrond$ des points de~$D_\eta$
au-dessus desquels la fibre de~$f'_\eta$ n'est pas lisse et géométriquement
intègre est donc un fermé strict de~$D_\eta$ (cf.~\cite[12.2.4 (iii)]{ega43}).
Pour chaque $m \in \Mrond$, notons $(\fy{Y}{m,i})_{i \in I_m}$ la famille des composantes
irréductibles de $f'^{-1}_\eta(m)$ et~$e_{m,i}$ la multiplicité de~$\fy{Y}{m,i}$
dans $f'^{-1}_\eta(m)$.
Notons~$F$ l'adhérence de~$\Mrond$ dans~$D$.  Pour $m \in \Mrond$, notons~$\mtilde$
l'adhérence de~$m$ dans~$D$ et pour $i \in I_m$, notons $\widetilde{\fy{Y}{m,i}}$ l'adhérence
de~$\fy{Y}{m,i}$ dans~$X$.  Munissons tous ces fermés de leur structure de sous-schéma
fermé réduit.

Le lemme suivant, appliqué aux morphismes $\widetilde{\fy{Y}{m,i}} \rightarrow \mtilde$,
fournit un ouvert dense $\fv{V}{m,i} \subset \mtilde$ et un $\fv{V}{m,i}$\nobreakdash-schéma fini étale
connexe~$R_{m,i}$ pour chaque \mbox{$m \in \Mrond$} et chaque $i \in I_m$.

\bigskip
\begin{lemme}
\label{ch3harrmi}
Soient~$S$ un schéma noethérien intègre et $f \colon X \rightarrow S$ un
morphisme de type fini dont la fibre générique est irréductible et géométriquement réduite.
Il existe un ouvert dense $V \subset S$ et un $V$\nobreakdash-schéma fini étale connexe~$R$
tels que pour tout $v \in V$ tel que~$R_v$ soit intègre,
le schéma~$X_v$ soit intègre et la fermeture
algébrique de~$\kappa(v)$ dans~$\kappa(X_v)$ soit isomorphe, comme $\kappa(v)$\nobreakdash-algèbre,
à $\kappa(R_v)$.
\end{lemme}

\bigskip
\begin{demo}
Quitte à rétrécir~$S$, on peut supposer toutes les fibres de~$f$ géométriquement réduites
(cf.~\cite[9.7.7]{ega43}).

Notons~$\eta$ le point générique de~$S$.
L'ensemble des points réguliers
de~$X_\eta$ est un ouvert de~$X_\eta$ (cf.~\cite[6.12.5]{ega42}), dense
dans~$X_\eta$ puisque~$X_\eta$ est intègre.
Le complémentaire~$X'$ de l'adhérence dans~$X$ de l'ensemble des points
singuliers de~$X_\eta$
est donc un ouvert de~$X$
dense dans chaque fibre de~$f$ au-dessus d'un voisinage de~$\eta$
(cf.~\cite[9.6.1 (ii)]{ega43}).
Comme les fibres de~$f$ sont réduites, on voit qu'il suffit de prouver le lemme
pour $X' \rightarrow S$ plutôt que pour $X \rightarrow S$. Autrement dit, on
peut supposer que $X'=X$, c'est-à-dire que la fibre générique de~$f$ est régulière.
On fait donc dorénavant cette hypothèse.

Notons~$K$ la fermeture algébrique
de~$\kappa(S)$ dans~$\kappa(X_\eta)$.  Comme $X_\eta$ est de type fini et géométriquement
réduit sur~$\eta$, l'extension $K/\kappa(S)$ est finie et séparable
(cf.~\cite[4.6.3]{ega42}).
Quitte à rétrécir~$S$, on peut donc supposer qu'il existe un $S$\nobreakdash-schéma
fini étale connexe~$R$ de fibre générique \mbox{$\Spec(K) \rightarrow \eta$}.
Le schéma~$X_\eta$ est normal, étant régulier; il existe donc un morphisme
canonique $X_\eta \rightarrow \Spec(K)$ dont la composée avec $\Spec(K)
\rightarrow \eta$ soit égale à~$f$.
Par conséquent, quitte à rétrécir encore~$S$, on peut supposer qu'il existe un
morphisme $f' \colon X \rightarrow R$
rendant le triangle
$$
\myxyhook\xymatrix{
X \ar[rr]^f \ar[dr]_(.4){f'} && S \\
& R \ar[ur]
}
$$
commutatif (cf.~\cite[8.8.2]{ega43}). De plus, $X_\eta$ étant géométriquement
irréductible sur~$K$, il existe un ouvert dense de~$R$ au-dessus duquel
les fibres de~$f'$ sont géométriquement irréductibles (cf.~\cite[9.7.7]{ega43}).
Comme le morphisme $R \rightarrow S$ est fermé, un dernier rétrécissement
de~$S$ permet de supposer que toutes
les fibres de~$f'$ sont géométriquement irréductibles.

Soit $v \in S$ tel que~$R_v$ soit intègre.  Le schéma~$X_v$ est à la fois
une fibre de~$f$ et une fibre de~$f'$.  Il est donc géométriquement réduit
sur~$\kappa(v)$ et géométriquement irréductible sur~$\kappa(R_v)$.
En particulier il est intègre et l'on a une
suite d'inclusions $\kappa(v) \subset \kappa(R_v) \subset \kappa(X_v)$,
où l'extension $\kappa(X_v)/\kappa(v)$ est séparable
(cf.~\cite[4.6.3]{ega42}) et où~$\kappa(R_v)$
est séparablement fermé dans~$\kappa(X_v)$ (cf.~\cite[4.6.1]{ega42}).
Il en résulte que la fermeture algébrique de~$\kappa(v)$ dans~$\kappa(X_v)$
est~$\kappa(R_v)$, comme il fallait démontrer.
\end{demo}

\bigskip
Traduisons en termes des~$e_{m,i}$ et des~$R_{m,i}$ les hypothèses sur le morphisme~$f$.

\bigskip
\begin{lemme}
\label{ch3exisabel}
Pour tout $m \in \Mrond$, il existe $i \in I_m$ tel que
$e_{m,i}=1$ et que le revêtement étale
$R_{m,i}\rightarrow \fv{V}{m,i}$ soit abélien.
\end{lemme}

\bigskip
\begin{demo}
Soit~$m \in \Mrond$ tel que $\sigma(m) \neq c$.
Comme~$\sigma$ est l'éclatement de centre~$c$,
la fibre de~$f'$ en~$m$ s'identifie
à la fibre de~$f$ en~$\sigma(m)$ et le point~$\sigma(m)$ est de codimension~$1$
dans~$\P^n_k$.
La conclusion recherchée résulte donc immédiatement
des hypothèses du théorème.

Nous allons maintenant établir que $c \not\in \sigma(\Mrond)$, ce qui terminera
la preuve du lemme.
Le seul point de~$\sigma^{-1}(c)$ qui appartienne à~$D_\eta$ est le point
générique de~$\sigma^{-1}(c)$.
Notons-le~$\xi$.
La variété~$X_\eta$ est projective, lisse et géométriquement intègre
(cf.~conclusion de la proposition~\ref{ch3choixdec})
et la fibre générique de \mbox{$f'_\eta \colon X_\eta \rightarrow D_\eta$} est géométriquement intègre
(cf.~lemme~\ref{ch3fggif}).
Par conséquent,
toute fibre lisse de~$f'_\eta$ est géométriquement intègre.
La fibre de~$f'_\eta$ en~$\xi$ est égale à la fibre générique du morphisme
$f^{-1}(c) \rightarrow \sigma^{-1}(c)$ induit par~$f'$.
Celle-ci est lisse; en effet, la variété $f^{-1}(c)$ est elle-même lisse, puisque $c \in U$.
Ainsi a-t-on prouvé que la fibre de~$f'_\eta$ en~$\xi$ est lisse et géométriquement
intègre, autrement dit que $\xi \not\in \Mrond$, d'où le résultat.
\end{demo}

\bigskip
Soit $\Delta^0 \subset \Delta$ un ouvert dense disjoint pour tout $m\in \Mrond$
du fermé $\pi\!\left(\mtilde \setminus \left(\bigcap_{i\in I_m} \fv{V}{m,i}\right)\right)$.
Quitte à rétrécir les ouverts~$\fv{V}{m,i} \subset \mtilde$, on peut supposer que
\begin{equation}
\label{ch3pplvm}
\fv{V}{m,i}=\pi^{-1}(\Delta^0) \cap \mtilde
\end{equation}
pour tout $m\in\Mrond$ et tout $i \in I_m$, ce qui n'a d'autre but que de simplifier
les notations.

\bigskip
\begin{lemme}
\label{ch3delta0petit}
Quitte à rétrécir l'ouvert $\Delta^0 \subset \Delta$, on peut supposer
que le morphisme \mbox{$\pi^{-1}(\Delta^0)\cap F \rightarrow \Delta^0$} induit par~$\pi$
est fini et étale et
que pour tout $\delta \in \Delta^0$, les propriétés suivantes sont vérifiées:
\begin{enumerate}
\item La variété~$X_\delta$ est lisse et géométriquement intègre sur~$\kappa(\delta)$.
\item Le morphisme $f'_\delta \colon X_\delta \rightarrow D_\delta$ est propre, plat,
de fibre générique géométriquement intègre.
\item Les fibres de~$f'_\delta$ au-dessus du complémentaire
de~$D_\delta \cap F$ dans~$D_\delta$ sont lisses et géométriquement intègres.
\item Soit $x \in D_\delta \cap F$.  Comme $\pi^{-1}(\Delta^0)\cap F \rightarrow \Delta^0$
est étale, il existe un unique $m \in \Mrond$ tel que $x \in \mtilde$.
Pour $i \in I_m$, notons $(\widetilde{\fy{Y}{m,i}})_x$ (resp.~$X_x$) la fibre en~$x$
de la composée $\widetilde{\fy{Y}{m,i}} \subset X \xrightarrow{f'} D$ (resp.~de~$f'$).
Alors
le morphisme canonique
$$
\coprod_{i \in I_m} \left(\widetilde{\fy{Y}{m,i}}\right){}_{\!\!x} \longrightarrow X_x
$$
induit une bijection entre composantes irréductibles, et pour tout $i \in I_m$,
la multiplicité dans~$X_x$ de chaque composante irréductible de~$(\widetilde{\fy{Y}{m,i}})_x$
est égale à $e_{m,i}$.
\end{enumerate}
\end{lemme}

\bigskip
\begin{demo}
L'assertion préliminaire résulte de la finitude de l'ensemble~$\Mrond$; celui-ci
s'identifie en effet à
la fibre générique du morphisme $F \rightarrow \Delta$ induit par~$\pi$.
La propriété~1 pour $\delta=\eta$ est satisfaite par hypothèse (cf.~conclusion
de la proposition~\ref{ch3choixdec}).  Elle est donc satisfaite pour tout~$\delta$
dans un ouvert dense de~$\Delta$.
La propriété~2
est mise pour mémoire (cf.~lemme~\ref{ch3fggif}; comme précédemment,
la platitude découle de la propriété~1 et de \cite[III, 9.7]{hartshorne}).
La propriété~3 est justiciable de \cite[17.7.11~(iii)]{ega44}, appliqué
au $\Delta$\nobreakdash-morphisme $X \setminus f'^{-1}(F) \rightarrow D\setminus F$
induit par~$f'$ (pour l'intégrité géométrique,
utiliser par exemple~\cite[9.7.7 et 9.5.2]{ega43}).
Le morphisme de la propriété~4 est la fibre en~$x$ de
la composée
$$
\xymatrix{
\displaystyle
\coprod_{i \in I_m} \widetilde{\fy{Y}{m,i}} \ar[r]^\alpha & (X_\mtilde)_{\mathrm{red}}
\ar[r]^(.57)\beta & X_\mtilde \rlap{\text{,}}
}
$$
où $(X_\mtilde)_{\mathrm{red}}$ désigne le sous-schéma fermé réduit de~$X_\mtilde$
d'espace topologique sous-jacent~$X_\mtilde$.
Considérons~$\alpha$ comme un morphisme de $\mtilde$\nobreakdash-schémas. Sa restriction
au-dessus du point générique de~$\mtilde$ est un morphisme birationnel.
Sa restriction au-dessus de chaque point d'un certain ouvert dense de~$\mtilde$
est donc un morphisme birationnel.
Quitte à rétrécir~$\Delta^0$, on peut supposer~$\Delta^0$
disjoint de l'image par~$\pi$ du complémentaire dans~$\mtilde$ de cet ouvert dense.
La première partie de la propriété~4 est alors établie.  L'assertion concernant
les multiplicités des composantes irréductibles résulte maintenant de \cite[9.8.6]{ega43}.
\end{demo}

\bigskip
Pour chaque $\delta \in \Delta$, on notera $\Br_\vert(X_\delta)$ le groupe
de Brauer vertical du morphisme $f'_\delta \colon X_\delta \rightarrow D_\delta$, c'est-à-dire le
sous-groupe de $\Br(X_\delta)$ constitué des classes dont la restriction à
la fibre générique de~$f'_\delta$ est l'image réciproque d'une classe
de $\Br(\kappa(D_\delta))$.

\bigskip
\begin{proposition}
\label{ch3brvertspecialise}
Le groupe $\Br_\vert(X_\eta)/g^\star \Br(\eta)$ est fini et il existe
un ensemble hilbertien $H \subset \Delta(k)$ tel que pour tout $\delta \in H$,
la flèche de spécialisation
$$
\Br_\vert(X_\eta)/g^\star \Br(\eta) \longrightarrow \Br_\vert(X_\delta)/g^\star \Br(k)
$$
soit bien définie et surjective.
\end{proposition}

\bigskip
Si $\uplet{A_1}{A_m} \in \Br_\vert(X_\eta)$ est une famille finie
dont les images dans $\Br_\vert(X_\eta)/g^\star \Br(\eta)$ engendrent ce groupe,
une flèche de spécialisation comme ci-dessus
est définie pour tout~$\delta$ appartenant à un ouvert dense $\Delta^0 \subset \Delta$
assez petit pour que les classes~$A_i$ appartiennent au groupe $f'^\star
\Br(\pi^{-1}(\Delta^0))$
(cf.~preuve ci-dessous).
Ce que la proposition affirme, c'est que cette flèche est surjective si~$\delta$
appartient à un sous-ensemble hilbertien bien choisi de~$\Delta^0$.

\bigskip
\begin{demo}
Rappelons que pour $m \in \Mrond$ et $i \in I_m$, on dispose
d'un revêtement étale connexe $R_{m,i} \rightarrow \pi^{-1}(\Delta^0) \cap \mtilde$.
Notons $$e_{m,i} \Res \colon H^1(\kappa(m),\Q/\Z) \longrightarrow H^1(\kappa(R_{m,i}),\Q/\Z)$$
la composée de la flèche de restriction associée à ce revêtement
et de la multiplication par~$e_{m,i}$.
Pour $m \in \Mrond$, posons
$$
K_m = \Ker\left( H^1(\kappa(m),\Q/\Z) \xrightarrow{\;\;\prod e_{m,i} \Res\;\;}
\prod_{i\in I_m} H^1(\kappa(R_{m,i}),\Q/\Z) \right)\!\rlap{\text{.}}
$$

\bigskip
\begin{lemme}
\label{ch3kmfini}
Pour tout $m \in \Mrond$, le groupe~$K_m$ est fini.
\end{lemme}

\bigskip
\begin{demo}
Cela résulte du lemme~\ref{ch3exisabel}.  En effet, si $i \in I_m$ est tel que $e_{m,i}=1$,
le noyau de la flèche $e_{m,i}\Res$ est déjà fini puisqu'il s'identifie
à l'ensemble des classes d'isomorphisme d'extensions cycliques de~$\kappa(m)$ qui se plongent
dans~$\kappa(R_{m,i})$ et que~$\kappa(R_{m,i})$ est une extension finie (et séparable) de~$\kappa(m)$.
\end{demo}

\bigskip
Fixons un point rationnel $\infty \in (D_\eta \setminus \Mrond)(\eta)$
et posons $W=D \setminus (F \cup \inftytilde)$, où $\inftytilde$ désigne
l'adhérence de~$\infty$ dans~$D$.  Pour chaque $m \in \Mrond$ et chaque $r \in K_m$,
choisissons une classe $A_{m,r} \in \Br(D_\eta \setminus \{m,\infty\})$
dont le résidu en~$m$ soit égal à~$r$; une telle classe existe d'après
la suite exacte de Faddeev (cf.~\cite[§1.2]{cscrelle94}), puisque
$D_\eta$ est isomorphe à~$\P^1_\eta$.
Soit $B \subset \Br(D_\eta \cap W)$ le sous-groupe engendré par la famille
$(A_{m,r})_{m\in \Mrond, r \in K_m}$; le lemme~\ref{ch3kmfini} montre que ce groupe est fini.
Soit $\Phi \subset H^1(\eta, \Q/\Z)$ l'image de~$B$ par la flèche
$\Br(D_\eta \cap W) \rightarrow H^1(\eta, \Q/\Z)$
qui à une classe de~$\Br(D_\eta \cap W)$ associe son résidu au point~$\infty$.
Soit enfin $B_0 \subset B$ le noyau de la surjection naturelle $B \rightarrow \Phi$.

\bigskip
\begin{lemme}
Le groupe $\Br_\vert(X_\eta)/g^\star\Br(\eta)$ est engendré par $f'^\star_\eta B_0$.
\end{lemme}

\bigskip
(Nous notons $f'^\star_\eta$ la flèche $\Br(\kappa(D_\eta)) \rightarrow \Br(\kappa(X_\eta))$ induite par
la fibre générique de~$f'_\eta$.  Ce léger abus ne cause pas d'ambiguïté puisque $\Br(X_\eta)$
est un sous-groupe de $\Br(\kappa(X_\eta))$. De tels raccourcis seront employés
à nouveau par la suite mais ne seront plus signalés.)

\bigskip
\begin{demo}
Vérifions d'abord que les classes de $f'^\star_\eta B_0$ sont
non ramifiées sur~$X_\eta$.  Elles sont évidemment non ramifiées sur $f'^{-1}_\eta(D_\eta \cap W)$.
Par définition de~$B_0$, elles sont non ramifiées sur $f'^{-1}_\eta(\infty)$.
De plus, pour tout $m \in \Mrond$ et tout $r \in K_m$, la classe $f'^\star_\eta A_{m,r}$ est
non ramifiée sur $f'^{-1}_\eta(m)$ d'après \cite[Proposition~1.1.1]{cscrelle94}, d'où l'assertion.

La proposition précitée montre par ailleurs que toute classe de $\Br_\vert(X_\eta)$
est l'image réciproque par~$f'_\eta$ d'une classe de $A \in \Br(D_\eta \cap W)$
dont le résidu en~$m$ appartient à~$K_m$ pour tout $m \in \Mrond$ et dont
le résidu au point~$\infty$ est trivial.  Il existe donc $A' \in B$ tel
que $A-A' \in \Br(D_\eta \setminus \{\infty\})$. Ce dernier groupe est
réduit à $\pi^\star \Br(\eta)$, d'après la suite exacte de Faddeev.  Par conséquent,
il existe $\gamma \in \Br(\eta)$ tel que $A=A'+\pi^\star \gamma$.  Il reste seulement à vérifier
que $A' \in B_0$, mais cela résulte de cette même égalité.
\end{demo}

\bigskip
Comme les groupes~$B$, $\Phi$ et~$K_m$ pour $m\in \Mrond$ sont finis,
quitte à rétrécir~$\Delta^0$, on peut supposer
que
$B \subset \Br(\pi^{-1}(\Delta^0) \cap W)$,
que $\Phi \subset H^1(\Delta^0,\Q/\Z)$
et que $K_m \subset H^1(\pi^{-1}(\Delta^0) \cap \mtilde,\Q/\Z)$ pour tout $m \in \Mrond$.
Par ailleurs, quitte à rétrécir encore~$\Delta^0$, on peut supposer le
morphisme $\pi^{-1}(\Delta^0) \cap (F \cup \inftytilde) \rightarrow \Delta^0$ induit par~$\pi$
fini et étale.

Notons $j_{\Delta^0} \colon \pi^{-1}(\Delta^0) \cap W \rightarrow \pi^{-1}(\Delta^0)$ l'immersion
ouverte canonique.  La suite spectrale de Leray
$$
H^p(\pi^{-1}(\Delta^0),
R^qj_{\Delta^0 \star}\Gm) \Longrightarrow H^{p+q}(\pi^{-1}(\Delta^0)\cap W, \Gm)
$$
pour le morphisme $j_{\Delta^0}$ et le faisceau étale~$\Gm$ induit une flèche
\begin{equation}
\label{ch3brouvert}
H^2(\pi^{-1}(\Delta^0),
j_{\Delta^0\star}\Gm) \longrightarrow \Br(\pi^{-1}(\Delta^0) \cap W) \rlap{\text{.}}
\end{equation}

\bigskip
\begin{lemme}
\label{ch3brh2}
Quitte à rétrécir~$\Delta^0$, on peut supposer que pour tout $m\in \Mrond$
et tout $r \in K_m$, la classe $A_{m,r}$ appartient à l'image de la flèche~(\ref{ch3brouvert}).
\end{lemme}

\bigskip
\begin{demo}
Notons $j_\eta \colon \pi^{-1}(\eta) \cap W \rightarrow \pi^{-1}(\eta)$ l'immersion
ouverte canonique.  Les deux lemmes suivants montrent que si~$\Delta^0$ est assez petit,
les restrictions des
classes $A_{m,r}$ à $\Br(\pi^{-1}(\eta)\cap W)$ appartiennent à l'image de la composée
$$
H^2(\pi^{-1}(\Delta^0), j_{\Delta^0\star}\Gm) \longrightarrow
H^2(\pi^{-1}(\eta), j_{\eta\star}\Gm) \longrightarrow
\Br(\pi^{-1}(\eta) \cap W)\rlap{\text{,}}
$$
où la première flèche est la flèche de restriction à la fibre générique de~$\pi$
et la seconde est induite par
la suite spectrale de Leray associée à~$j_\eta$.

\bigskip
\begin{lemme}
Soient $C$ une courbe lisse sur un corps parfait
et $j \colon U \rightarrow C$
une immersion ouverte.
La flèche $H^2(C, j_\star \Gm) \rightarrow \Br(U)$ induite
par la suite spectrale de Leray est un isomorphisme.
\end{lemme}

\bigskip
\begin{demo}
Soient~$\cbarre$ un point géométrique de~$C$ et $V=U \times_C \Spec(\Orond_{C,\cbarre})$.
Pour tout entier~$q$, la tige en~$\cbarre$ du faisceau étale $R^q j_\star \Gm$
est égale à $H^q(V,\Gm)$.  Pour $q=1$, ce groupe s'identifie au groupe de Picard de~$V$
(théorème de Hilbert~90).  Pour $q=2$, c'est par définition le groupe de Brauer de~$V$.
Le schéma~$V$ peut être un trait strictement local,
le spectre d'un corps de degré de transcendance~$\leq 1$
sur un corps algébriquement clos, ou le schéma vide.  Dans les trois cas,
les groupes de Picard et de Brauer de~$V$ sont nuls (cinq de ces assertions
sont triviales, la sixième résulte du théorème de Tsen).

Ainsi a-t-on prouvé que $R^q j_\star \Gm=0$ pour $q \in \{1,2\}$.
La suite spectrale de Leray permet d'en déduire le résultat voulu.
\end{demo}

\bigskip
\begin{lemme}
Soient~$S$ un schéma intègre, $X \rightarrow S$
un morphisme de schémas quasi-compact et quasi-séparé
et $j \colon U \rightarrow X$ une immersion ouverte quasi-compacte.
Pour tout $S$\nobreakdash-schéma $V$, notons $j_V \colon U_V \rightarrow X_V$
le morphisme déduit de~$j$ par changement de base.
Notons de plus~$I$ l'ensemble ordonné filtrant des ouverts
affines non vides de~$S$ et~$\eta$ le point générique de~$S$.
Alors, pour tout entier~$n \geq 0$, le morphisme canonique
$$
\limind_{V \in I}\, H^n(X_V, j_{V\star}\Gm) \longrightarrow H^n(X_\eta, j_{\eta\star}\Gm)
$$
est un isomorphisme.
\end{lemme}

\bigskip
\begin{demo}
C'est un cas particulier de \cite[Exp.~VII, Cor.~5.8]{sga42},
compte tenu de ce que pour tout $\fv{V}0 \in I$
et tout $V \in I$ tel que $V \subset \fv{V}0$ ou pour $V=\eta$,
notant~$f$ l'inclusion de~$X_V$ dans~$X_{\fv{V}0}$,
la flèche canonique de faisceaux étales sur~$X_V$
$$
f^\star j_{\fv{V}0\star}\Gm \longrightarrow j_{V \star}\Gm
$$
est un isomorphisme.  Cette dernière affirmation est triviale lorsque $V \in I$.
Pour $V=\eta$, le point délicat est de vérifier qu'étant donnés un schéma~$W$, un
morphisme $W \rightarrow X_\eta$
étale et de présentation finie et un $U_\eta$\nobreakdash-morphisme $W \times_{X_\eta} U_\eta
\rightarrow \Gmsur{U_\eta}$, on peut étendre tous ces objets au-dessus d'un ouvert
dense de~$\fv{V}0$; ceci résulte de \cite[8.8.2]{ega43} et \cite[17.7.8~(ii)]{ega44}.
\end{demo}

\bigskip
Le morphisme $\eta \rightarrow \Delta^0$ induit un morphisme entre les suites spectrales
de Leray associées à~$j_\eta$ et à~$j_{\Delta^0}$, d'où la commutativité du carré
$$
\myxyhook\xymatrix{
H^2(\pi^{-1}(\Delta^0),j_{\Delta^0 \star}\Gm) \ar[d] \ar[r] &
\Br(\pi^{-1}(\Delta^0) \cap W) \ar[d] \\
H^2(\pi^{-1}(\eta),j_{\eta\star}\Gm) \ar[r] & \Br(\pi^{-1}(\eta)\cap W)\rlap{\text{,}}
}
$$
dont la flèche verticale de droite est injective. Ceci termine de prouver le lemme~\ref{ch3brh2}.
\end{demo}

\bigskip
Nous sommes maintenant en position de définir l'ensemble hilbertien~$H$.
Rappelons que pour tout $m \in \Mrond$, le morphisme $\pi^{-1}(\Delta^0) \cap \mtilde
\rightarrow \Delta^0$ est par hypothèse un revêtement étale connexe
(cf.~lemme~\ref{ch3delta0petit}),
et l'on dispose par ailleurs pour chaque $i \in I_m$ d'un revêtement étale
connexe \mbox{$R_{m,i} \rightarrow \pi^{-1}(\Delta^0) \cap \mtilde$} (cf.~équation~(\ref{ch3pplvm}));
d'où, par composition, un revêtement étale connexe $R_{m,i} \rightarrow \Delta^0$.
D'autre part, chaque élément de~$\Phi$ définit un revêtement étale connexe (cyclique)
de~$\Delta^0$ puisque $\Phi \subset H^1(\Delta^0,\Q/\Z)$.
Soit~$H$ l'ensemble des $\delta \in \Delta^0(k)$ tels que la fibre en~$\delta$
de chacun des revêtements étales connexes de~$\Delta^0$ évoqués ci-dessus
soit connexe.

Soit $\delta \in H$.  Notons $\phi \colon D_\delta \cap W \rightarrow \pi^{-1}(\Delta^0) \cap W$
l'immersion fermée canonique et $\phi^\star B_0$ l'image de~$B_0$
par $\phi^\star \colon \Br(\pi^{-1}(\Delta^0)\cap W)
\rightarrow \Br(D_\delta \cap W)$.
Nous allons maintenant démontrer que le groupe $\Br_\vert(X_\delta)/g^\star\Br(k)$
est engendré par $f'^\star_\delta\phi^\star B_0$, ce qui terminera la
preuve de la proposition.

Pour la commodité du lecteur, nous résumons dans le lemme ci-dessous les propriétés
du morphisme~$f'_\delta$.  Elles découlent immédiatement des conclusions du
lemme~\ref{ch3delta0petit}, de la définition des revêtements~$R_{m,i}$
et de la définition de~$H$.

\bigskip
\begin{lemme}
\label{ch3descfpdelta}
Le morphisme $f'_\delta \colon X_\delta \rightarrow D_\delta$ est un morphisme propre, plat, de fibre
générique géométriquement intègre, entre variétés projectives, lisses et géométriquement intègres.
Il est lisse au-dessus de $D_\delta \setminus \ff{F}\delta$.
L'application $\Mrond \rightarrow D_\delta$, $m \mapsto D_\delta \cap \mtilde$ définit une
bijection entre~$\Mrond$ et l'ensemble des points fermés de~$D_\delta$ au-dessus desquels
la fibre de~$f'_\delta$ est singulière. Pour tout $m \in \Mrond$, les composantes
irréductibles de $f'^{-1}_\delta(D_\delta \cap \mtilde)$ sont exactement
les fermés $(\widetilde{\fy{Y}{m,i}})_{D_\delta \cap \mtilde}$ pour $i \in I_m$.
Pour $m \in \Mrond$ et $i \in I_m$, la fermeture algébrique de $\kappa(D_\delta \cap \mtilde)$
dans $\kappa((\widetilde{\fy{Y}{m,i}})_{D_\delta\cap \mtilde})$ est isomorphe,
comme $\kappa(D_\delta \cap \mtilde)$\nobreakdash-algèbre, à $\kappa((R_{m,i})_\delta)$,
et la multiplicité de
$(\widetilde{\fy{Y}{m,i}})_{D_\delta \cap \mtilde}$ dans $f'^{-1}_\delta(D_\delta \cap \mtilde)$
est~$e_{m,i}$.
\end{lemme}

\bigskip
Vérifions que les classes de $f'^\star_\delta \phi^\star B_0$
sont non ramifiées sur~$X_\delta$.
Elles sont évidemment non ramifiées sur $f'^{-1}_\delta(D_\delta \cap W)$,
puisque \mbox{$B_0 \subset \Br(\pi^{-1}(\Delta^0) \cap W)$}.
Qu'elles soient non ramifiées aux autres points de codimension~$1$ de~$X_\delta$
résulte de la description ci-dessus des fibres singulières de~$f'_\delta$,
des deux lemmes suivants et de \cite[Proposition~1.1.1]{cscrelle94}.

\bigskip
\begin{lemme}
\label{ch3compatev}
Pour tout $m \in \Mrond$, tout $r \in K_m$ et tout $s \in \Mrond \cup \{\infty\}$,
le résidu de $\phi^\star A_{m,r}$
au point fermé $D_\delta \cap \stilde$ de~$D_\delta$ est égal à l'image du résidu
de~$A_{m,r}$ en~$s$
par la flèche de restriction
\begin{equation}
\label{ch3restrh1}
H^1(\pi^{-1}(\Delta^0) \cap \stilde, \Q/\Z) \longrightarrow H^1(D_\delta \cap \stilde, \Q/\Z)
\rlap{\text{.}}
\end{equation}
\end{lemme}

\vspace{-3mm}
(Le résidu de~$A_{m,r}$ en~$s$ appartient à $H^1(\pi^{-1}(\Delta^0) \cap \stilde,\Q/\Z)$
car il appartient à~$K_m$.)

\bigskip
\begin{demo}
Les variétés considérées sont résumées dans le diagramme suivant, qui permet en outre
de donner des noms aux flèches qui n'en possèdent pas encore.
$$
\xymatrix@C=8ex{
\pi^{-1}(\Delta^0) \cap W \ar[r]^(.55){j_{\Delta^0}} & \pi^{-1}(\Delta^0) &
\pi^{-1}(\Delta^0) \cap \stilde \ar[l]_(.52){i_{\Delta^0}} \\
D_\delta \cap W \ar[u]_\phi \ar[r]^(.55){j_\delta} & D_\delta \ar[u]_i &
D_\delta \cap \stilde \ar[l]_(.52){i_\delta} \ar[u]
}
$$
Tous les morphismes de ce diagramme sont des immersions fermées, excepté $j_{\Delta^0}$
et~$j_\delta$, qui sont des immersions ouvertes.

Les schémas $\pi^{-1}(\Delta^0) \cap \stilde$ et $D_\delta \cap \stilde$ sont normaux,
puisqu'ils sont par hypothèse étales respectivement au-dessus de~$\Delta^0$ et de~$\delta$.
Toute composante connexe d'un ($\pi^{-1}(\Delta^0) \cap \stilde$\hspace{0.15ex})-schéma étale
(resp.~d'un ($D_\delta \cap \stilde$\hspace{0.15ex})-schéma étale) est donc aussi une composante irréductible.
Cette remarque permet de définir un morphisme
$j_{\Delta^0 \star}\Gm \rightarrow i_{\Delta^0 \star} \Z$
(resp.~$j_{\delta\star}\Gm \rightarrow i_{\delta \star} \Z$)
de faisceaux étales sur~$\pi^{-1}(\Delta^0)$
(resp.~sur~$D_\delta$) en associant à une fonction rationnelle sur
un $\pi^{-1}(\Delta^0)$\nobreakdash-schéma étale connexe
(resp.~sur un $D_\delta$\nobreakdash-schéma étale connexe)
la famille de ses valuations sur les composantes irréductibles
de l'image réciproque de $\pi^{-1}(\Delta^0)\cap \stilde$ (resp.~de $D_\delta\cap\stilde$).
Les morphismes que l'on vient de définir s'inscrivent dans le carré
\begin{equation}
\begin{aligned}
\label{ch3carrecommute}
\xymatrix{
j_{\Delta^0 \star} \Gm \ar[d] \ar[r] & i_{\Delta^0\star}\Z \ar[d] \\
i_\star j_{\delta \star} \Gm \ar[r] & i_\star i_{\delta \star} \Z
}
\end{aligned}
\end{equation}
de faisceaux étales sur~$\pi^{-1}(\Delta^0)$, où les flèches verticales
sont les flèches canoniques.

\bigskip
\begin{souslemme}
\label{ch3souslemmecarrecom}
Le carré~(\ref{ch3carrecommute}) est commutatif.
\end{souslemme}

\bigskip
\begin{demo}
Le point clé est que les sous-variétés fermées $\pi^{-1}(\Delta^0) \cap \stilde$
et~$D_\delta$ de~$\pi^{-1}(\Delta^0)$ se rencontrent transversalement
au point~$D_\delta \cap \stilde$.  En déduire la commutativité du carré~(\ref{ch3carrecommute})
est une vérification aisée, compte tenu que la transversalité de l'intersection de deux
sous-schémas fermés réguliers d'un schéma régulier est préservée par image réciproque
par un morphisme étale.  Pour établir le point clé, il suffit d'utiliser l'hypothèse
que $\pi^{-1}(\Delta^0) \cap \stilde$ est étale sur~$\Delta^0$.  Celle-ci
entraîne que tout vecteur tangent à~$\Delta^0$ se relève en un vecteur
tangent à $\pi^{-1}(\Delta^0) \cap \stilde$.  Tout vecteur tangent à $\pi^{-1}(\Delta^0)$
en $D_\delta \cap \stilde$ s'écrit donc comme la somme d'un vecteur
tangent à $\pi^{-1}(\Delta^0) \cap \stilde$ et d'un vecteur dont l'image par~$\pi$
est nulle, autrement dit, d'un vecteur tangent à~$D_\delta$.  Ceci signifie, par définition,
que $\pi^{-1}(\Delta^0) \cap \stilde$ et $D_\delta$ se rencontrent transversalement
en~$D_\delta \cap \stilde$.
\end{demo}

\bigskip
Appliquant le foncteur $H^2(\pi^{-1}(\Delta^0),-)$ au carré~(\ref{ch3carrecommute}),
on obtient la face de gauche du diagramme suivant. Les flèches obliques de droite
sont les flèches de résidu.
Comme $\pi^{-1}(\Delta^0) \cap \stilde$ et $D_\delta \cap \stilde$ sont des schémas
normaux, il résulte de la suite spectrale de Leray pour l'inclusion du point générique
que les groupes $H^q(\pi^{-1}(\Delta^0),\Q)$ et $H^q(\pi^{-1}(D_\delta \cap \stilde\mkern1.2mu),\Q)$
sont nuls pour tout $q \geq 1$.
Par conséquent, les flèches naturelles
$H^1(\pi^{-1}(\Delta^0), \Q/\Z) \rightarrow H^2(\pi^{-1}(\Delta^0), \Z)$
et
$H^1(D_\delta \cap \stilde, \Q/\Z) \rightarrow H^2(D_\delta \cap \stilde, \Z)$
sont des isomorphismes; ce sont leurs inverses qui apparaissent dans le diagramme ci-dessous.
$$
\myxyin\myxyhook\xymatrix{
H^2(\pi^{-1}(\Delta^0), j_{\Delta^0 \star}\Gm) \ar@<-2ex>[dd] \ar[rr]
\save!<-3.1mm,0mm>\ar+(15,-11.17)\restore &&
*+!!<0pt,\the\fontdimen22\textfont2>!<7.6ex,0ex>{\Br(\pi^{-1}(\Delta^0) \cap W)}
\ar@<-7.6ex>@{-}+(0,-13.1)
\save+<0mm,-15.5mm>\ar@<-7.6ex>+(0,-10.1)\restore
\save!<-7.6ex,0ex>\ar+(15,-11.17)\restore
\\
&
*+!!<0pt,\the\fontdimen22\textfont2>!<19ex,0ex>{H^2(\pi^{-1}(\Delta^0) \cap \stilde, \Z)}
\ar[dd]<-19ex>
\save-<12.83mm,0ex>\ar+<4.64mm,0ex>^(.45)\sim\restore
&
*+!!<0pt,\the\fontdimen22\textfont2>!<21.5ex,0ex>{H^1(\pi^{-1}(\Delta^0) \cap \stilde, \Q/\Z)}
\save+<-14.24mm,0mm>\ar@{^{ (}->}+(4.9,0)\restore \ar@<-22ex>[dd] &
*+!!<0pt,\the\fontdimen22\textfont2>!<26.5ex,0ex>{H^1(\kappa(s), \Q/\Z)} \\
*+!!<0pt,\the\fontdimen22\textfont2>!<2ex,0ex>{H^2(D_\delta, j_{\delta\star}\Gm)}
\save!<-3.1mm,0mm>\ar+(15,-11.17)\restore
\ar@{-}+(12.5,0)
\save+<14.95mm,0mm>\ar@{-}+(35.54,0)
\save+<38.01mm,0mm>\ar+(10.3,0)
\restore\restore &&
*+!!<0pt,\the\fontdimen22\textfont2>!<7.6ex,0ex>{\Br(D_\delta \cap W)}
\save!<-7.6ex,0mm>\ar+(15,-11.17)\restore
\\
&
*+!!<0pt,\the\fontdimen22\textfont2>!<19ex,0ex>{H^2(D_\delta\cap\stilde, \Z)}
\save-<17.31mm,0ex>\ar+<9.73mm,0ex>^\sim\restore &
*+!!<0pt,\the\fontdimen22\textfont2>!<24ex,0ex>{H^1(D_\delta \cap \stilde, \Q/\Z)}
\ar@{=}[r] &
*+!!<0pt,\the\fontdimen22\textfont2>!<29ex,0ex>{H^1(\kappa(D_\delta \cap \stilde), \Q/\Z)}
}\myxyout
$$
Ce diagramme est commutatif. En effet, la face de gauche l'est grâce au
sous-lemme~\ref{ch3souslemmecarrecom} et les faces horizontales le sont par définition
des flèches de résidu.  Pour tout $m \in \Mrond$ et tout $r \in K_m$,
la classe $A_{m,r} \in \Br(\pi^{-1}(\Delta^0) \cap W)$ provient par hypothèse d'un
élément de $H^2(\pi^{-1}(\Delta^0), j_{\Delta^0\star}\Gm)$
(cf.~conclusion du lemme~\ref{ch3brh2}).  La commutativité
du diagramme ci-dessus
permet donc de conclure.
\end{demo}

\bigskip
Pour $m \in \Mrond$ et $i \in I_m$, notons
$$
e_{m,i} \Res_\delta \colon H^1(D_\delta \cap \mtilde,\Q/\Z) \longrightarrow
H^1((R_{m,i})_\delta,\Q/\Z)
$$
la composée de la flèche de restriction associée au revêtement
étale connexe $(R_{m,i})_\delta \rightarrow D_\delta \cap \mtilde$
et de la multiplication par~$e_{m,i}$.

\bigskip
\begin{lemme}
\label{ch3kmsuiteex}
Pour tout $m \in \Mrond$, la suite
$$
\xymatrix{
K_m \ar[r] & H^1(D_\delta \cap \mtilde, \Q/\Z) \ar[rr]^(.465){\prod e_{m,i} \Res_\delta} &&
\displaystyle \prod_{i\in I_m} H^1((R_{m,i})_\delta,\Q/\Z) \rlap{\text{,}}
}
$$
dans laquelle la première flèche est la restriction de~(\ref{ch3restrh1}) à~$K_m$,
est exacte.
\end{lemme}

\bigskip
\begin{demo}
Soit $m \in \Mrond$.
Que cette suite forme un complexe résulte de la définition de~$K_m$ et de la commutativité des carrés
$$
\xymatrix@C=10ex{
H^1(\pi^{-1}(\Delta^0)\cap \mtilde, \Q/\Z) \ar[d] \ar[r]^(.54){e_{m,i}\Res} &
H^1(R_{m,i}, \Q/\Z) \ar[d] \\
H^1(D_\delta \cap \mtilde, \Q/\Z) \ar[r]^(.54){e_{m,i}\Res_\delta} & H^1((R_{m,i})_\delta,\Q/\Z)
}
$$
pour $i \in I_m$.  Considérons maintenant une classe de $H^1(D_\delta \cap \mtilde,
\Q/\Z)$ appartenant au noyau de la seconde flèche de ce complexe
et choisissons une extension cyclique $\ell/\kappa(D_\delta \cap \mtilde)$ la
représentant (rappelons que \mbox{$D_\delta \cap \mtilde$} est le spectre d'un corps).
Il existe un $i_0 \in I_m$
tel que $e_{m,i_0}=1$ (cf.~lemme~\ref{ch3exisabel}).
L'appartenance de la classe de~$\ell$ au noyau de la restriction
$$H^1(D_\delta \cap \mtilde,\Q/\Z) \longrightarrow H^1((R_{m,i_0})_\delta,\Q/\Z)$$
signifie que l'extension $\ell/\kappa(D_\delta \cap \mtilde)$
se plonge dans $\kappa((R_{m,i_0})_\delta)/\kappa(D_\delta \cap \mtilde)$.
Autrement dit, le morphisme $(R_{m,i_0})_\delta \rightarrow D_\delta \cap \mtilde$
se factorise
en
$(R_{m,i_0})_\delta \rightarrow \Spec(\ell) \rightarrow D_\delta \cap \mtilde$.
Comme $(R_{m,i_0})_\delta$ est connexe,
cette factorisation s'étend en une factorisation
$$R_{m,i_0} \longrightarrow L \longrightarrow \pi^{-1}(\Delta^0) \cap \mtilde$$
du revêtement étale connexe $R_{m,i_0} \rightarrow \pi^{-1}(\Delta^0) \cap \mtilde$,
où $L \rightarrow \pi^{-1}(\Delta^0) \cap \mtilde$ est un revêtement cyclique
dont la fibre en~$\delta$ s'identifie à~$\Spec(\ell)$.
Il reste à vérifier
que la classe~$[L]$ de~$L$ dans $H^1(\pi^{-1}(\Delta^0) \cap \mtilde, \Q/\Z)$ appartient à~$K_m$.
Soit donc $i \in I_m$.  Notons $L_i \rightarrow \pi^{-1}(\Delta^0) \cap \mtilde$
un revêtement cyclique connexe dont la classe~$[L_i]$
dans $H^1(\pi^{-1}(\Delta^0) \cap \mtilde,\Q/\Z)$ soit égale à $e_{m,i}[L]$.
L'appartenance de la classe de $\ell$ au noyau de $e_{m,i}\Res_\delta$
se traduit par l'existence d'une factorisation du morphisme
$(R_{m,i})_\delta \rightarrow D_\delta \cap \mtilde$ par
$(L_i)_\delta$.  Comme $(R_{m,i})_\delta$ est connexe, ceci entraîne
que le morphisme $R_{m,i} \rightarrow \pi^{-1}(\Delta^0) \cap \mtilde$
se factorise par~$L_i$,
ce qui à son tour
implique que~$[L_i]$ appartient au noyau de la flèche de restriction
$H^1(\pi^{-1}(\Delta^0) \cap \mtilde,\Q/\Z) \rightarrow H^1(R_{m,i},\Q/\Z)$,
autrement dit que~$[L]$ appartient au noyau de $e_{m,i}\Res$, comme il fallait démontrer.
\end{demo}

\bigskip
Le lemme~\ref{ch3descfpdelta}
et la proposition \cite[Proposition~1.1.1]{cscrelle94} montrent que toute classe
de $\Br_\vert(X_\delta)$ est l'image réciproque par $f'_\delta$ d'une
classe $A \in \Br(D_\delta \cap W)$ dont le résidu en $D_\delta \cap \mtilde$
appartient au noyau de $e_{m,i}\Res_\delta$ pour tout $m \in \Mrond$ et tout $i \in I_m$
et dont le résidu en $D_\delta \cap \inftytilde$ est nul.
D'après les lemmes~\ref{ch3kmsuiteex} et~\ref{ch3compatev},
il existe $A' \in B$ tel que $\phi^\star A'$ et~$A$ aient mêmes résidus
en $D_\delta \cap \mtilde$ pour tout $m \in \Mrond$.
On a alors $A - \phi^\star A' \in
\Br(D_\delta \setminus (D_\delta \cap \inftytilde))$.  Ce dernier groupe
est réduit à $\pi_\delta^\star \Br(k)$ d'après la suite exacte de Faddeev; d'où l'existence
de $\gamma \in \Br(k)$ tel que $A = \phi^\star A' + \pi_\delta^\star \gamma$.  En particulier,
comme~$A$ est non ramifiée au point $D_\delta \cap \inftytilde$, il en va de même
pour $\phi^\star A'$.  Autrement dit, et compte tenu du lemme~\ref{ch3compatev},
l'image de~$A'$ dans~$\Phi$ appartient au noyau
de la flèche
$H^1(\Delta^0,\Q/\Z)\rightarrow H^1(\delta,\Q/\Z)$
d'évaluation en~$\delta$.
La restriction de celle-ci au sous-groupe $\Phi \subset H^1(\Delta^0,\Q/\Z)$
est injective par définition de~$H$.  L'image de~$A'$ dans~$\Phi$ est donc nulle,
ce qui signifie encore que $A' \in B_0$.  On a alors bien établi que toute classe
de $\Br_\vert(X_\delta)$ s'écrit comme la somme d'une classe constante et de
l'image réciproque par $f'_\delta$ d'une classe de~$\phi^\star B_0$,
et la proposition~\ref{ch3brvertspecialise} est donc démontrée.
\end{demo}

\bigskip
Nous pouvons maintenant terminer la preuve du théorème~\ref{ch3thpointslocaux}.
Soient $(\fp{P}v)_{v \in \Omega} \in \prod_{v \in \Omega} X(k_v)$ une famille
orthogonale à $\Br_\vert(X)$ et $B\subset\Br_\vert(X_\eta)$
un sous-groupe fini engendrant $\Br_\vert(X_\eta)/g^\star\Br(\eta)$.
(Rappelons que dans $\Br_\vert(X)$, le symbole «~vert~» fait référence au
morphisme~$f$ (ou au morphisme~$f'$, c'est équivalent), alors que dans $\Br_\vert(X_\delta)$
pour $\delta \in \Delta$ (et notamment pour $\delta=\eta$), il fait
référence au morphisme $f'_\delta$.)

La variété~$\Delta$ est $k$\nobreakdash-isomorphe à~$\P^{n-1}_k$.
Par ailleurs, chaque fibre de~$g$
contient un ouvert géométriquement intègre d'après le lemme~\ref{ch3fggif}.
Toutes les conditions sont donc réunies pour appliquer le théorème suivant, dû à Harari,
au morphisme $g \colon X \rightarrow \Delta$ et au sous-groupe $B \subset \Br(X_\eta)$.

\bigskip
\begin{theoreme}[ {(cf.~\cite[Théorème~1]{harprepub})}]%
\label{ch3thharari}
Soient~$k$ un corps de nombres et~$X$ une variété projective, lisse et géométriquement connexe sur~$k$, munie d'un morphisme $f \colon X \rightarrow \P^n_k$ de fibre générique
géométriquement intègre.  Supposons que chaque fibre de~$f$ au-dessus du complémentaire
d'un fermé de codimension~$\geq 2$ contienne un ouvert géométriquement intègre.
Soient~$\eta$ le point générique de~$\P^n_k$
et $B \subset \Br(X_\eta)$ un sous-groupe fini.
Notons $B' \subset \Br(X)$ le groupe $(B+f^\star\Br(\eta))\cap \Br(X)
\subset \Br(\kappa(X))$ et~$U \subset \P^n_k$ un ouvert dense au-dessus
duquel~$f$ soit lisse, suffisamment petit pour que $B\subset \Br(X_U)$.
Pour tout $(\fp{P}v)_{v\in\Omega}\in X(\A_k)^{B'}$
et tout ensemble hilbertien $H \subset \P^n(k)$,
il existe $\theta \in H \cap U$
et une famille $(Q_v)_{v\in\Omega} \in X_\theta(\A_k)^{B_\theta}$
arbitrairement proche de~$(\fp{P}v)_{v \in \Omega}\in X(\A_k)$,
où $B_\theta \subset \Br(X_\theta)$
désigne l'image de~$B$ par la flèche de restriction $\Br(X_U)\rightarrow \Br(X_\theta)$.
\end{theoreme}

\bigskip
Comme
$(B+g^\star\Br(\eta)) \cap \Br(X) \subset \Br_\vert(X)$,
il existe $\delta \in H$
et une famille $(Q_v)_{v \in \Omega} \in X_\delta(\A_k)$ arbitrairement proche
de $(\fp{P}v)_{v \in \Omega}$ et orthogonale à l'image de~$B$ par spécialisation en~$\delta$,
où~$H$ est l'intersection de l'ensemble hilbertien donné par la
proposition~\ref{ch3brvertspecialise} et de celui déterminé
par les revêtements étales connexes $R_{m,i}\rightarrow \Delta^0$
pour $m \in \Mrond$ et $i \in I_m$.
Vu la définition de~$H$, la famille $(Q_v)_{v\in \Omega}$ est alors orthogonale
à $\Br_\vert(X_\delta)$.
Des conclusions du lemme~\ref{ch3delta0petit},
du lemme~\ref{ch3exisabel} et de la définition des revêtements~$R_{m,i}$
résultent les propriétés suivantes: la variété~$X_\delta$ est projective, lisse, géométriquement
intègre, le morphisme
$f'_\delta \colon X_\delta \rightarrow D_\delta$ est plat, sa fibre générique est géométriquement
intègre et pour tout $m \in D_\delta$, la fibre $f'^{-1}_\delta(m)$ contient une composante
irréductible de multiplicité~$1$ dans le corps des fonctions de laquelle la fermeture
algébrique de~$\kappa(m)$ est une extension abélienne de~$\kappa(m)$.
On peut donc appliquer \cite[Theorem~1.1]{csscrelle98} au morphisme~$f'_\delta$
et au point adélique $(Q_v)_{v \in \Omega}$, et conclure quant à l'existence
de $x \in D_\delta(k) \cap \sigma^{-1}(U)$ arbitrairement proche de $(f'(Q_v))_{v\in\Omega}$ et
tel que $f'^{-1}(x)$ soit lisse et possède un point adélique.
Le point $\sigma(x) \in \P^n(k)$
appartient alors à $\Rrond\cap U(k)$ et est arbitrairement proche de la famille $(f(\fp{P}v))_{v\in\Omega}$.
\end{demo}

\bigskip
\begin{remarque}
L'hypothèse de Schinzel\index{hypothèse de Schinzel} n'a servi qu'à appliquer \cite[Theorem~1.1]{csscrelle98}.
Avec une version adaptée du théorème d'Harari,
on obtiendrait donc, suivant exactement la même preuve,
des résultats inconditionnels sur les $0$\nobreakdash-cycles de degré~$1$ généralisant
ceux de \cite{csscrelle98} à des fibrations au-dessus de~$\P^n_k$.
\end{remarque}

\bigskip
\begin{remarque}
On pourrait être tenté de prouver le théorème en raisonnant
comme dans \cite[Theorem~1.1]{csscrelle98}.  Il s'agirait,
une fois la réduction au cas projectif et lisse accomplie, d'appliquer
directement le\index{lemme formel} «~lemme formel~» à une famille finie
$\uplet{A_1}{A_m} \in \Br_\vert(X_\eta)$
engendrant $\Br_\vert(X_\eta)/g^\star\Br(\eta)$, ce qui fournirait un ensemble
fini~$S$ de places de~$k$.  Notant~$H$ l'ensemble hilbertien donné par
la proposition~\ref{ch3brvertspecialise} et $V \subset \Delta$ un ouvert
dense assez petit pour que les classes~$A_i$ soient non ramifiées sur $g^{-1}(V)$,
on choisirait ensuite un $\delta \in H \cap V$ arbitrairement
proche d'une famille fixée $(\delta_v)_{v\in S} \in \prod_{v\in S}\Delta(k_v)$, à l'aide de la
version d'Ekedahl du théorème d'irréductibilité de Hilbert\index{théorème d'irréductibilité de Hilbert} (cf.~\cite{ekedahl}).
L'hypothèse de Schinzel fournirait enfin un $x \in D_\delta(k)$ pour lequel on pourrait
espérer que $f'^{-1}(x)$ possède un point adélique.  L'existence d'un $k_v$\nobreakdash-point
de~$f'^{-1}(x)$ serait facile à assurer pour $v \in S$
(proximité de~$\delta$ et de~$\delta_v$)
 et pour toute place~$v$
de bonne réduction pour~$f'^{-1}(x)$, quitte à imposer à~$S$ d'être suffisamment grand
(estimations de Lang-Weil et lemme de Hensel).
On voudrait alors se servir d'un argument de réciprocité
basé sur la conclusion du «~lemme formel~»
pour obtenir l'existence d'un $k_v$\nobreakdash-point aux places de mauvaise réduction
hors de~$S$ (places que l'on peut même choisir hors d'un ensemble fini arbitrairement
grand \emph{dépendant de $\delta$}).
Le problème de cette approche est que l'on n'a aucun contrôle sur les places $v \not\in S$
\emph{de bonne réduction pour $f'^{-1}(x)$} pour lesquelles la réduction de~$\delta$ modulo~$v$
appartient à la réduction modulo~$v$ du fermé $\Delta \setminus V$.
C'est un obstacle à l'argument de réciprocité dont il vient d'être question;
l'emploi d'un théorème plus fort tel que le théorème~\ref{ch3thharari} semble
ainsi inévitable.
\end{remarque}

\section{Généralités sur les pinceaux de quadriques dans~$\P^n$}
\label{ch3sectiongeneralites}

Nous\index{pinceau de quadriques!généralités, revêtement associé|(} commençons par rappeler quelques propriétés élémentaires et bien connues
des pinceaux de quadriques dans~$\P^n$ puis nous définissons des
revêtements de l'espace projectif dual~$(\P^n)^\star$ qui leur sont naturellement
associés.  Ceux-ci joueront un rôle central dans les preuves des théorèmes~\ref{ch3introthdp4}
et~\ref{ch3introp5} (cf.~notamment le paragraphe~\ref{ch3monodr}).

Dans tout ce paragraphe, nous fixons un corps~$k$ de caractéristique différente de~$2$,
une clôture séparable~$\ksep$ de~$k$,
un $k$\nobreakdash-espace vectoriel~$V$ de dimension finie et deux formes quadratiques $q_1$
et~$q_2$ sur~$V$.  L'espace projectif~$\P(V)$ des droites de~$V$
sera considéré comme une variété sur~$k$.  Notons~$n$ sa dimension et supposons
que $n \geq 2$.
Pour $t \in \P^1_k$ de coordonnées homogènes $t=[\lambda:\mu]$, soit $Q_t \subset \P(V)$
la quadrique projective définie par l'équation
$$
\lambda q_1 + \mu q_2 = 0 \rlap{\text{.}}
$$
Toutes les propriétés de~$q_1$ et~$q_2$ auxquelles on va maintenant s'intéresser
ne dépendent en fait que de la famille de quadriques $(Q_t)_{t\in\P^1_k}$.
Celle-ci définit un pinceau si les formes~$q_1$ et~$q_2$ ne sont pas proportionnelles.

Notons $f(\lambda,\mu)=\det(\lambda q_1 + \mu q_2)$ le déterminant
de $\lambda A + \mu B$, où~$A$ et~$B$ sont les matrices de~$q_1$ et de~$q_2$
dans une base quelconque de~$V$.  C'est un polynôme homogène en~$(\lambda, \mu)$,
nul ou de degré~$n+1$, bien défini à multiplication par un élément de~$k^\star$ près.  En particulier,
ses racines et leurs multiplicités respectives sont bien définies.

Soit $X \subset \P(V)$ l'intersection des quadriques~$(Q_t)_{t \in \P^1_k}$, autrement
dit la sous-variété de~$\P(V)$ définie par les équations $q_1=q_2=0$.

\bigskip
\begin{proposition}
\label{ch3genpropeq}
Les propriétés suivantes sont équivalentes:
\begin{itemize}
\item[(i)] La variété~$X$ est lisse sur~$k$ et purement de codimension~$2$ dans~$\P(V)$.
\item[(ii)] Le polynôme homogène $f(\lambda,\mu) \in k[\lambda,\mu]$ est non nul et séparable.
\item[(iii)] Les formes quadratiques~$q_1$ et~$q_2$ ne sont pas proportionnelles et pour tout $t \in \P^1_k$, le lieu singulier de la quadrique~$Q_t$ est
disjoint de~$X$.
\end{itemize}

\smallskip
\noindent{}De plus, elles impliquent:
\begin{itemize}
\item[(iv)] Parmi les formes quadratiques $\lambda q_1 + \mu q_2$ sur $V \otimes_k \ksep$
avec $(\lambda,\mu)\in \ksep^2 \setminus \{(0,0)\}$,
toutes sont de rang $\geq n$
et au moins une est de rang~$n+1$.
\end{itemize}
\end{proposition}

\bigskip
\begin{lemme}
\label{ch3genracsimple}
Soit $t \in \P^1(k)$.
Le point~$t$ est racine simple du polynôme homogène $f(\lambda,\mu)$ si et seulement si
la quadrique~$Q_t$ possède un unique point singulier et que celui-ci n'appartient pas à~$X$.
\end{lemme}

\bigskip
Nous allons simultanément prouver le lemme et la proposition.

\bigskip
\begin{demo}
Les propriétés~(i) à~(iii) étant invariantes par extension des scalaires,
on peut supposer~$k$ algébriquement clos.
Soient $(\lambda_0,\mu_0)\in k^2 \setminus\{(0,0)\}$ une racine de~$f$
et~$r$ le rang de $\lambda_0 q_1 + \mu_0 q_2$.
Choisissons une base de~$V$ dans laquelle la matrice de $\lambda_0 q_1 + \mu_0 q_2$
soit diagonale, les~$r$ premiers coefficients diagonaux étant non nuls;
on note cette matrice~$J$.
Pour $i \in \{\uplet{1}{n+1}\}$, soit~$D_i$ (resp.~$D'_i$)
la matrice obtenue en remplaçant la $i$\nobreakdash-ème colonne de~$J$
par la $i$\nobreakdash-ème colonne de la matrice de~$q_1$ (resp.~de~$q_2$).
Vu la forme de~$J$, on a $\det(D_i)=\det(D'_i)=0$ pour tout~$i$ si~$r<n$
et pour tout $i\leq n$ si~$r=n$.  Compte tenu des égalités
$$
\frac{\partial f}{\partial \lambda}(\lambda_0, \mu_0)=\sum_{i=1}^{n+1}\det(D_i)
\qquad\text{ et }\qquad
\frac{\partial f}{\partial \mu}(\lambda_0, \mu_0)=\sum_{i=1}^{n+1}\det(D'_i) \rlap{\text{,}}
$$
il s'ensuit que $(\lambda_0,\mu_0)$ est une racine multiple de~$f$ si et seulement
si $\det(D_{n+1})=\det(D'_{n+1})=0$ ou~$r<n$.
Ceci prouve déjà que (ii)$\Rightarrow$(iv).

Notons $t=[\lambda_0:\mu_0]$.
Si $r<n$ et $r>0$, le lieu singulier de~$Q_t$ contient une droite;
celle-ci rencontre nécessairement~$X$.  Si $r=n$, la condition $\det(D_{n+1})=\det(D'_{n+1})=0$
équivaut à l'appartenance à~$X$ de l'unique point singulier de~$Q_t$.
On a donc aussi prouvé le lemme, et par suite, l'équivalence entre~(ii) et~(iii).

Il reste à établir que~(i)$\Leftrightarrow$(iii). Notons $\ft{T}mM$ l'espace
tangent à une variété~$M$ en un point~$m$.
Soient $t \in \P^1(k)$ et $x \in X(k)$.
Choisissons un $u \in \P^1(k) \setminus \{t\}$.
Si la propriété~(i) est satisfaite, le sous-espace $\ft{T}x X \subset \P(V)$ est de codimension~$2$.
Étant donné que $\ft{T}x Q_u$ est de codimension au plus~$1$ dans~$\P(V)$
et que $\ft{T}x X = \ft{T}x Q_t \cap \ft{T}x Q_u$,
on en déduit que~$\ft{T}x Q_t$
est de codimension au moins~$1$, ce qui signifie que~$x$ est un point régulier de~$Q_t$.
D'où (i)$\Rightarrow$(iii).
Supposons enfin que la propriété~(i) ne soit pas satisfaite et montrons que~(iii)
ne l'est pas non plus. Le cas où~$X$ n'est pas purement de codimension~$2$ étant trivial,
on peut supposer que~$X$ possède un point singulier $x \in X(k)$.
Le noyau de la matrice jacobienne en~$x$ du système $q_1=q_2=0$ est alors non nul,
or la donnée d'un vecteur non nul de ce noyau équivaut exactement à la donnée
d'un $t \in \P^1(k)$ tel que la quadrique~$Q_t$ soit singulière en~$x$,
d'où le résultat.
\end{demo}

\bigskip
Supposons dorénavant que~$X$ est lisse sur~$k$ et purement de codimension~$2$ dans~$\P(V)$.
Comme on vient de le voir, le polynôme~$f$
possède alors $n+1$ racines deux à deux distinctes dans~$\P^1(\ksep)$.
Notons-les $\uplet{t_0}{t_n}$.
La propriété~(iv) ci-dessus montre que chacune des quadriques~$Q_{t_i}$ possède un
unique $\ksep$\nobreakdash-point singulier, que l'on note~$\fp{P}i$.  Les $\ksep$\nobreakdash-points $\uplet{\fp{P}0}{\fp{P}n}$
de~$\P(V)$ sont globalement stables sous l'action du groupe de Galois de~$\ksep$ sur~$k$;
on n'hésitera donc pas à identifier leur ensemble à une réunion de points fermés
de~$\P(V)$.

\bigskip
\begin{proposition}
\label{ch3gensimdiag}
Supposons que $k=\ksep$.  Pour chaque $i \in \{\uplet{0}{n}\}$, soit
\mbox{$v_i \in V\setminus\{0\}$} un vecteur dont l'image dans~$\P(V)$ soit égale à~$\fp{P}i$.
Alors $(\uplet{v_0}{v_n})$ est l'unique base de~$V$, à permutation des coordonnées
près et à homothétie près (indépendamment sur chacun des vecteurs de base),
dans laquelle les formes quadratiques~$q_1$ et~$q_2$ soient simultanément diagonales.
\end{proposition}

\bigskip
(On dit qu'une forme quadratique est \emph{diagonale} dans une base donnée
si sa matrice dans cette base est diagonale, ce qui équivaut encore à ce que la
base soit orthogonale pour la forme quadratique en question.)

\bigskip
\begin{demo}
Prouvons d'abord que les vecteurs $\uplet{v_0}{v_n}$ forment une base
de~$V$ dans laquelle~$q_1$ et~$q_2$ sont diagonales.
Il suffit pour cela de vérifier d'une part que les~$v_i$ sont deux à deux orthogonaux
à la fois pour~$q_1$ et pour~$q_2$ et d'autre part qu'aucun des~$v_i$ n'est isotrope
à la fois pour~$q_1$ et pour~$q_2$.  La seconde assertion résulte de
la proposition~\ref{ch3genpropeq}, propriété~(iii).  Pour la première,
il suffit de montrer
que pour tous $i,j\in \{\uplet{0}{n}\}$ distincts, il existe
$t,t' \in \P^1(k)$ distincts
tels que les vecteurs~$v_i$ et~$v_j$ soient orthogonaux pour les deux formes quadratiques
$\lambda q_1 + \mu q_2$ et $\lambda' q_1 + \mu' q_2$,
où $t=[\lambda:\mu]$ et $t'=[\lambda':\mu']$.
Il est immédiat que cette
condition est satisfaite pour $(t,t')=(t_i,t_j)$, puisque~$v_i$ (resp.~$v_j$) est alors
orthogonal à tout~$V$ pour $\lambda q_1 + \mu q_2$ (resp.~$\lambda' q_1 + \mu' q_2$).

Établissons maintenant l'unicité.  Supposons que dans une base de~$V$ l'on
puisse écrire $q_1=\sum_{i=0}^n a_i x_i^2$ et $q_2=\sum_{i=0}^n b_i x_i^2$.
Les racines de~$f$ sont alors les $t_i=[-b_i:a_i]\in\P^1(k)$
pour $i \in \{\uplet{0}{n}\}$. La droite engendrée
par le $i$\nobreakdash-ème vecteur de base est évidemment un point singulier de~$Q_{t_i}$,
ce qui signifie que ce vecteur est colinéaire à~$v_i$.
\end{demo}

\bigskip
\begin{corollaire}
\label{ch3genpiengendrent}
Aucun hyperplan de~$\P(V)$ ne contient simultanément tous les $\ksep$\nobreakdash-points $\uplet{\fp{P}0}{\fp{P}n}$.
\end{corollaire}

\bigskip
\begin{corollaire}
\label{ch3gencorsimdiag}
Le polynôme~$f$ est scindé si et seulement s'il
existe une base de~$V$ dans laquelle les formes quadratiques~$q_1$ et~$q_2$ sont simultanément
diagonales.
\end{corollaire}

\bigskip
\begin{demo}
En effet, pour tout $i \in \{\uplet{0}{n}\}$, le point $t_i \in \P^1(\ksep)$ est
$k$\nobreakdash-rationnel si et seulement si le point~$\fp{P}i$ l'est, puisque ce dernier est
l'unique point singulier de~$Q_{t_i}$.
\end{demo}

\bigskip
Intéressons-nous maintenant à la trace du pinceau $(Q_t)_{t\in \P^1_k}$ sur
un hyperplan variable de~$\P(V)$.  Les hyperplans de~$\P(V)$ sont paramétrés
par l'espace projectif dual~$\P(V^\star)$.
Posons
$$
Z = \bigensemble{(t,L) \in \P^1_k \times_k \P(V^\star)}{ Q_t \cap L \text{ n'est pas lisse}}\!\rlap{\text{,}}
$$
où l'intersection est à prendre au sens schématique.  Ce sous-ensemble de \mbox{$\P^1_k
\times_k \P(V^\star)$} est un fermé (cf.~\cite[12.1.7]{ega43}); munissons-le de sa structure
de sous-schéma fermé réduit.

Les fibres géométriques de la projection naturelle $p \colon Z \rightarrow \P(V^\star)$ sont
particulièrement simples à décrire. Pour $L \in \P(V^\star)$, posons
$f_L(\lambda,\mu)=\det((\lambda q_1 + \mu q_2)|_L)$, avec la notation évidente
pour la restriction d'une forme quadratique au sous-espace vectoriel de~$V$ associé à~$L$.
C'est un polynôme homogène; la fibre géométrique de~$p$ en~$L$ s'identifie à
l'ensemble de ses racines (avec multiplicités).

\bigskip
\begin{proposition}
\label{ch3genpfiniplat}
Le morphisme $p \colon Z \rightarrow \P(V^\star)$ est fini et plat, de degré~$n$.
\end{proposition}

\bigskip
\begin{demo}
Supposons d'abord le morphisme~$p$ fini.  Il résulte alors de la description de ses
fibres géométriques qu'elles sont toutes de degré~$n$, ce qui entraîne
que~$p$ est un morphisme plat (lemme de Nakayama).
Il suffit donc de vérifier que~$p$ est fini.  Pour cela,
il suffit de vérifier qu'il est quasi-fini,
puisqu'il est évidemment propre. Par l'absurde, supposons qu'il existe $L \in \P(V^\star)$
tel que $Q_t \cap L$ ne soit lisse pour aucun $t \in \P^1_k$.
Considérant successivement les points $t=t_i$ pour $i\in\{\uplet{0}{n}\}$,
on trouve que l'hyperplan~$L$ doit contenir~$\fp{P}i$ pour tout~$i$, ce qui contredit
le corollaire~\ref{ch3genpiengendrent}.
\end{demo}

\bigskip
Ainsi avons-nous associé à tout pinceau de quadriques dans~$\P(V)$
un revêtement plat de~$\P(V^\star)$.  Nous y reviendrons au paragraphe~\ref{ch3monodr},
où nous étudierons la monodromie de ces revêtements.  Pour l'instant, limitons-nous à
la proposition suivante, fondamentale pour l'étude des points
rationnels sur les surfaces de del Pezzo de degré~$4$.

\bigskip
\begin{proposition}
\label{ch3genzkrat}
La variété~$Z$ est $k$\nobreakdash-rationnelle.
\end{proposition}

\bigskip
C'est en réalité la description explicite d'une certaine équivalence birationnelle entre~$Z$
et un espace projectif qu'il nous faudra connaître.  Voici comment la construire.

\bigskip
Pour $(t,L) \in Z$ tel que $t \not\in \{\uplet{t_0}{t_n}\}$, la quadrique
$Q_t \cap L \subset L$ possède un unique point singulier.  On définit une application rationnelle
$Z \dashrightarrow \P(V)$ en associant ce point au couple~$(t,L)$.

Pour $x \in \P(V) \setminus (X \cup \{\uplet{\fp{P}0}{\fp{P}n}\})$, il existe un unique $t \in \P^1_k$
tel que \mbox{$x \in Q_t$}. Le point~$x$ est alors régulier sur~$Q_t$. On définit une application
rationnelle $\P(V) \dashrightarrow Z$ en associant à~$x$ le couple $(t, \ft{T}x Q_t)$, où~$\ft{T}x Q_t$
désigne l'hyperplan tangent à~$Q_t$ en~$x$.

\bigskip
On vérifie tout de suite que ces deux applications rationnelles sont bien inverses l'une
de l'autre.  La seconde est définie sur
\mbox{$\P(V) \setminus (X \cup \{\uplet{\fp{P}0}{\fp{P}n}\})$}, qui est un ouvert de~$\P(V)$ dont le
complémentaire est de codimension~$2$.

\bigskip
Voici un résultat parfois utile lorsque l'on cherche à restreindre le revêtement $Z \rightarrow \P(V^\star)$
au-dessus d'hyperplans spécifiques de~$\P(V^\star)$.

\bigskip
\begin{proposition}
\label{ch3genlp0}
Soit $z=(t,L)\in Z$.
Supposons que l'un des points~$\fp{P}i$, disons~$\fp{P}0$, soit $k$\nobreakdash-rationnel.
S'il existe un point singulier de~$Q_t \cap L$ qui soit un point régulier de~$Q_t$
et qui appartienne à l'hyperplan de~$\P(V)$ contenant $\{\uplet{\fp{P}1}{\fp{P}n}\}$,
alors $\fp{P}0 \in L$.
\end{proposition}

\bigskip
\begin{demo}
On peut supposer que le corps~$k$ est algébriquement clos et que le point~$z$ est $k$\nobreakdash-rationnel.
Reprenons alors les notations de l'énoncé
de la proposition~\ref{ch3gensimdiag}.
L'hyperplan de~$V$ engendré par $(\uplet{v_1}{v_n})$ est orthogonal
au vecteur~$v_0$ pour toute forme quadratique $\lambda q_1 + \mu q_2$ avec $\lambda,\mu \in k$.
Par ailleurs, si~$x$ est un point singulier de~$Q_t \cap L$ qui soit un point régulier de~$Q_t$,
l'hyperplan de~$V$ associé à~$L$ est l'orthogonal de la droite vectorielle définie par~$x$
pour la forme quadratique $\lambda q_1 + \mu q_2$, où $t=[\lambda:\mu]$. La proposition s'ensuit.
\end{demo}

\bigskip
Pour conclure ce paragraphe, mentionnons une condition nécessaire pour qu'une intersection
de deux quadriques dans~$\P(V)$ possède une composante irréductible qui soit une sous-variété
linéaire (c'est-à-dire de degré~$1$) de~$\P(V)$.

\bigskip
\begin{proposition}
\label{ch3gensousvarlin}
Ne supposons plus~$X$ lisse sur~$k$.
Si~$n \geq 3$ et si~$X$ possède une composante irréductible qui est une sous-variété
linéaire de~$\P(V)$, alors le polynôme~$f$ n'admet aucune racine simple.
\end{proposition}

\bigskip
\begin{demo}
On peut supposer~$X$ purement de codimension~$2$ dans~$\P(V)$, les autres
cas étant triviaux. On peut d'autre part supposer~$k$ algébriquement clos.
Soit une composante irréductible~$I \subset X$ de degré~$1$
dans~$\P(V)$ et une racine simple $t_0 \in \P^1(k)$ de~$f$.
D'après le lemme~\ref{ch3genracsimple}, la quadrique~$Q_{t_0}$ est un cône quadratique possédant un unique point singulier,
disons~$\fp{P}0$, et celui-ci n'appartient pas à~$X$.
Il existe donc un hyperplan $L \subset \P(V)$ contenant~$I$ mais ne contenant pas~$\fp{P}0$.
Comme $I \subset X \subset Q_{t_0}$, on a alors $I \subset Q_{t_0} \cap L$;
autrement dit, la quadrique $Q_{t_0} \cap L$, qui est lisse et de codimension~$1$
dans~$L$ puisque $\fp{P}0\not\in L$ et que~$n \geq 3$, contient un sous-espace linéaire
de codimension~$1$ dans~$L$. C'est une\index{pinceau de quadriques!généralités, revêtement associé|)} contradiction.
\end{demo}

\section{Surfaces de del Pezzo de degré~$4$}
\subsection{Notations, énoncés des résultats}
\label{ch3notations}

Soit~$X$ une surface de del Pezzo de degré~$4$ sur un corps~$k$ de caractéristique~$0$.
Soient~$\ksep$ une clôture algébrique de~$k$ et~$\Gamma$ le groupe de Galois
de~$\ksep$ sur~$k$.

Il existe un plongement $X \subset \P^4_k$ permettant de voir~$X$ comme
l'intersection de deux quadriques dans~$\P^4_k$.
Soient~$q_1$ et~$q_2$ deux formes quadratiques homogènes en cinq variables
telles que la surface~$X$ soit définie par le système $q_1=q_2=0$.
Reprenons alors les notations du paragraphe~\ref{ch3sectiongeneralites} associées
à ces deux formes. 
Elles engendrent un pinceau $(Q_t)_{t\in\P^1_k}$ de quadriques dans~$\P^4_k$.
Comme~$X$ est une surface lisse, il
existe exactement cinq valeurs du paramètre $t \in \P^1(\ksep)$ pour lesquelles la
quadrique~$Q_t$ n'est pas lisse (cf.~proposition~\ref{ch3genpropeq}).  On
les note~$\uplet{t_0}{t_4}$\glossary{$\uplet{t_0}{t_4}$, $\uplet{\fp{P}0}{\fp{P}4}$, $\Srond$}
et l'on note $\fp{P}i \in \P^4(\ksep)$ l'unique
point singulier de~$Q_{t_i}$.

Le groupe~$\Gamma$ opère naturellement sur $\{\uplet{t_0}{t_4}\}$.
Notons~$\Srond$ l'ensemble de ses orbites.
Nous considérerons les éléments de~$\Srond$ comme des points fermés de~$\P^1_k$,
ce qui confère notamment un sens au symbole~$\kappa(t)$ pour~$t \in \Srond$.
La longueur d'une orbite s'interprète comme le degré sur~$k$ du point fermé correspondant.
Soit $t \in \Srond$.  La quadrique $Q_t \subset \P^4_k$, définie sur~$\kappa(t)$,
possède un unique point singulier.  Soit~$L$ un hyperplan $\kappa(t)$\nobreakdash-rationnel
de~$\P^4_k$ ne contenant pas ce point.  Considérons alors la quadrique $Q_t \cap L \subset L$;
elle est lisse, et comme~$L$ est de dimension impaire, toutes les formes quadratiques non nulles sur
le sous-espace vectoriel de~$V$ associé à~$L$ qui s'annulent sur $Q_t \cap L$ auront le même discriminant.
Celui-ci
ne dépend pas de l'hyperplan~$L$ choisi,
comme on le voit tout de suite en
considérant une base dans laquelle~$Q_t$ est définie par une forme
quadratique diagonale.
On note ce discriminant $\epsilon_t \in \kappa(t)^\star/\kappa(t)^{\star 2}$\glossary{$\epsilon_t$, $\epsilon_0$}.
Une autre façon de définir~$\epsilon_t$ est la suivante: si~$L$ est un
hyperplan $\kappa(t)$\nobreakdash-rationnel de~$\P^4_k$ contenant le point singulier de~$Q_t$
et tangent à~$Q_t$ selon une droite, l'intersection~$Q_t \cap L$ est géométriquement
une réunion de deux plans, et~$\epsilon_t$ est précisément la classe dans
$\kappa(t)^\star/\kappa(t)^{\star 2}$
de l'extension quadratique ou triviale minimale de~$\kappa(t)$ qui permute ces deux plans.

\bigskip
\begin{exemple}
\label{ch3exemplesimdiag1}
Supposons que les formes quadratiques~$q_1$ et~$q_2$ soient simultanément diagonalisables
sur~$k$, autrement dit (cf.~corollaire~\ref{ch3gencorsimdiag}) que~$\Gamma$
agisse trivialement sur $\{\uplet{t_0}{t_4}\}$.  On peut alors écrire
\vskip-1.5em
\begin{align*}
q_1 &= a_0 x_0^2 + a_1 x_1^2 + a_2 x_2^2 + a_3 x_3^2 + a_4 x_4^2 \rlap{\text{,}}\\
q_2 &= b_0 x_0^2 + b_1 x_1^2 + b_2 x_2^2 + b_3 x_3^2 + b_4 x_4^2
\end{align*}
avec $\uplet{a_0}{a_4},\uplet{b_0}{b_4}\in k$.
Les cinq formes dégénérées dans le pinceau
$\lambda q_1 + \mu q_2$ sont les $a_i q_2 - b_i q_1$ pour $i \in \{\uplet{0}{4}\}$,
ce qui permet de choisir $t_i=[-b_i:a_i] \in \P^1(k)$.
Posant $d_{ij}=a_ib_j-a_jb_i$ et $\epsilon_i = \epsilon_{t_i}$,
on voit alors que les $\epsilon_i \in k^\star/k^{\star 2}$ sont donnés par la formule suivante:
$$
\epsilon_i = \prod_{\genfrac{}{}{0pt}{}{0 \leq j \leq 4}{j \neq i}} d_{ij} \rlap{\text{.}}
$$
\end{exemple}

Lorsque l'on supposera que l'un des~$t_i$ est $k$\nobreakdash-rationnel, on conviendra toujours
que~$t_0$ l'est, quitte à renuméroter les~$t_i$. On posera
alors
$\epsilon_0=\epsilon_{t_0}$
et l'on notera $\Srond' \subset \Srond$ l'ensemble des orbites
de~$\Gamma$ sur~$\{\uplet{t_1}{t_4}\}$. On aura à considérer l'hypothèse suivante:\glossary{$\Srond'$}
\begin{itemize}
\medskip
\item[(\refstepcounter{equation}\theequation{})\label{ch3enonceshyp}]
Pour tout $t \in \Srond'$ de degré au plus~$3$ sur~$k$,
on a $\epsilon_t \neq 1$ dans $\kappa(t)^\star/\kappa(t)^{\star 2}$.
\end{itemize}
\medskip

\bigskip
Les résultats que nous démontrons au sujet du principe de Hasse pour
les surfaces de del Pezzo de degré~$4$ sont les suivants.

\bigskip
\begin{theoreme}
\label{ch3thprindp4}
Supposons que~$k$ soit un corps de nombres.
Admettons l'hypothèse de Schinzel et
la finitude des groupes de Tate-Shafarevich des courbes elliptiques sur les corps
de nombres. Alors le principe de Hasse vaut pour~$X$ dans chacun des cas
suivants:
\begin{itemize}
\item[(i)] le groupe~$\Gamma$ agit $3$\nobreakdash-transitivement sur $\{\uplet{t_0}{t_4}\}$;
\item[(ii)] l'un des~$t_i$ est $k$\nobreakdash-rationnel et~$\Gamma$ agit $2$\nobreakdash-transitivement
sur les quatre autres;
\item[(iii)] exactement deux des~$t_i$ sont $k$\nobreakdash-rationnels et $\Br(X)/\Br(k)=0$;
\item[(iv)] tous les~$t_i$ sont $k$\nobreakdash-rationnels et $\Br(X)/\Br(k)=0$ (cas simultanément diagonal);
\item[(v)] le point~$t_0$ est $k$\nobreakdash-rationnel,
la condition~(\ref{ch3enonceshyp}) est satisfaite, on a $\Br(X)/\Br(k)=0$ et enfin
soit $\epsilon_0=1$ dans
$k^\star/k^{\star 2}$, soit il existe $t \in \Srond'$ tel que l'image de~$\epsilon_0$
dans $\kappa(t)^\star/\kappa(t)^{\star 2}$ soit distincte de~$1$ et de~$\epsilon_t$.
\end{itemize}
\end{theoreme}

\bigskip
Nous verrons que l'hypothèse~(v) est l'hypothèse la plus générale sous laquelle la méthode permet
directement de prouver le principe de Hasse avec les outils dont on dispose à l'heure actuelle.
Des progrès dans l'étude des points rationnels sur les pinceaux de courbes de genre~$1$
auraient cependant des répercussions immédiates pour les surfaces de del Pezzo de degré~$4$.
Notamment, si l'on savait supprimer la «~\condD{}~» dans le théorème~\ref{ch2thobmx} du
chapitre précédent (quitte à considérer toute l'obstruction de Brauer-Manin au lieu de
sa seule partie verticale), la méthode employée dans le présent chapitre prouverait
en toute généralité l'existence d'un point rationnel sur~$X$ dès que l'obstruction de Brauer-Manin
ne s'y oppose pas (toujours en admettant l'hypothèse de Schinzel et la
finitude des groupes de Tate-Shafarevich).

C'est la preuve du principe de Hasse lorsque l'hypothèse~(v) est satisfaite qui occupera
l'essentiel des paragraphes à venir.  Nous en déduirons le principe de Hasse
sous chacune des autres hypothèses à l'aide de considérations \emph{ad hoc}.
Les hypothèses~(i) à~(iv) sont les énoncés les plus simples auxquels nous aboutirons;
cependant, ils sont loin de couvrir la totalité des cas que l'hypothèse~(v) permet
d'obtenir.  Celle-ci entraîne notamment que pour chaque décomposition
possible de l'ensemble $\{\uplet{t_0}{t_4}\}$ en orbites sous~$\Gamma$,
le principe de Hasse vaut dès que les~$\epsilon_t$ sont «~suffisamment généraux~» en un
certain sens, qui peut être rendu explicite.

Le théorème~\ref{ch3thprindp4} généralise~\cite[Proposition~3.2.1]{css},
où le cas simultanément diagonal (correspondant ici à l'hypothèse~(iv))
était traité avec quelques hypothèses de généricité sur les coefficients strictement plus
fortes\footnote{Le lecteur attentif remarquera que telles quelles, les hypothèses de
\cite[Proposition~3.2.1]{css} n'impliquent pas que $\Br(X)/\Br(k)=0$.
Il s'avère que cette proposition est erronée;
c'est l'antépénultième phrase de la démonstration qui est en cause.
On obtient un énoncé correct en remplaçant la seconde hypothèse par la suivante:
«~If moreover
the six classes $-d_{12}d_{13}$, $-d_{21}d_{23}$,
$-d_{10}d_{14}$, $-d_{20}d_{24}$, $-d_{30}d_{34}$,
$\prod_{i\neq 4}d_{4i}$ are independent in~$k^\star/k^{\star 2}$,
then \condD{} holds for~$\pi$.~»
Le théorème~\ref{ch3thbr} ci-dessous montre que cette condition est bien strictement plus
forte que $\Br(X)/\Br(k)=0$ (cf.~notamment l'exemple~\ref{ch3exemplesimdiag2}).
}
que
la condition $\Br(X)/\Br(k)=0$.
Il est remarquable que
le théorème~\ref{ch3thprindp4} établisse le principe de Hasse
dans le cas simultanément diagonal ainsi que dans le cas où exactement deux des~$t_i$ sont $k$\nobreakdash-rationnels
sous la seule hypothèse que \mbox{$\Br(X)/\Br(k)=0$},
alors que pour chacune des autres décompositions possibles de l'ensemble $\{\uplet{t_0}{t_4}\}$
en orbites sous~$\Gamma$, une hypothèse supplémentaire, quoique légère, semble indispensable pour que
la preuve fonctionne.

\subsection{Groupe de Brauer des surfaces de del Pezzo de degré~$4$}
\label{ch3parpremierpreuve}

Conservons\index{groupe de Brauer!d'une surface de del Pezzo de degré~$4$|(} les notations du paragraphe~\ref{ch3notations} et supposons que~$k$
soit un corps de nombres.
Nous allons entièrement déterminer la structure du groupe $\Br(X)/\Br(k)$
en termes des invariants~$\epsilon_t$ associés à la surface~$X$. Soulignons
que ceux-ci se lisent très facilement sur les équations
des formes quadratiques~$q_1$ et~$q_2$.

\bigskip
\begin{theoreme}
\label{ch3thbr}
Le groupe $\Br(X)/\Br(k)$ est isomorphe à $(\Z/2)^m$
avec $m={\max(0,n-d-1)}$,
où
$$
n=\Card\bigensemble{t \in \Srond}{\epsilon_t \neq 1 \text{ dans }
\kappa(t)^\star/\kappa(t)^{\star 2}}
$$
et~$d$ est
la dimension du sous-$\Z/2$\nobreakdash-espace vectoriel de $k^\star/k^{\star 2}$
engendré par les normes $N_{\kappa(t)/k}(\epsilon_t)$ pour $t \in \Srond$.
\end{theoreme}

\bigskip
\begin{corollaire}
\label{ch3brcor541}
Si le groupe~$\Gamma$ agit transitivement sur $\{\uplet{t_0}{t_4}\}$, ou
s'il agit trivialement sur l'un des~$t_i$ et transitivement sur les quatre
autres, alors $\Br(X)/\Br(k)=0$.
\end{corollaire}

\bigskip
\begin{demo}
La transitivité de l'action de~$\Gamma$ sur les cinq~$t_i$
entraîne que $n \leq 1$, d'où $n-d-1 \leq 0$.
Si maintenant l'un des~$t_i$ est $k$\nobreakdash-rationnel, disons~$t_0$, et que~$\Gamma$ agit transitivement
sur les quatre autres, alors $n \leq 2$ et la seule possibilité
pour que $n-d-1>0$ est que~$n=2$ et~$d=0$; mais alors
$\epsilon_0$ devrait être trivial puisque~$d=0$ et non trivial puisque~$n=2$.
\end{demo}

\bigskip
La proposition suivante sera un sous-produit de la preuve du théorème~\ref{ch3thbr}.

\bigskip
\begin{proposition}
\label{ch3brproduit}
On a l'égalité $\prod_{t \in \Srond} N_{\kappa(t)/k}(\epsilon_t)=1$ dans $k^\star/k^{\star 2}$.
\end{proposition}

\bigskip
\begin{corollaire}
\label{ch3brtroisgroupes}
Le groupe $\Br(X)/\Br(k)$ est isomorphe à l'un des trois groupes suivants: le groupe
trivial, $\Z/2$ ou $\Z/2 \times \Z/2$.
\end{corollaire}

\bigskip
\begin{demo}
Il s'agit de prouver que $n-d\leq 3$. Si $n \leq 3$, c'est trivial.
Si $n=4$, l'ensemble~$\Srond$ contient un point rationnel $t \in \Srond$
tel que $\epsilon_t \neq 1$, d'où $d \geq 1$ et donc $n-d \leq 3$.
Si enfin $n=5$, la proposition~\ref{ch3brproduit} montre que
les $\epsilon_t$ ne peuvent pas être tous égaux, sous peine d'être
alors tous égaux à~$1$. D'où $d \geq 2$ puis $n-d \leq 3$.
\end{demo}

\bigskip
\begin{remarque}
Les corollaires~\ref{ch3brcor541} et~\ref{ch3brtroisgroupes}
étaient déjà connus (cf.~respectivement \cite[Theorem~3.19]{ctsansdi} et~\cite{sdbrauergroup}).
Nous ne prétendons d'ailleurs à aucune originalité quant à la preuve du théorème~\ref{ch3thbr}:
nous avons simplement poussé quelque peu les arguments de \cite{ctsansdi}.
\end{remarque}

\bigskip
\begin{exemple}[ (suite de l'exemple~\ref{ch3exemplesimdiag1})]%
\label{ch3exemplesimdiag2}
Si les formes quadratiques~$q_1$ et~$q_2$ sont simultanément diagonalisables sur~$k$
et qu'aucun des $\epsilon_i \in k^\star/k^{\star 2}$ pour $i \in \{\uplet{0}{4}\}$
n'est égal à~$1$, alors $\Br(X)/\Br(k)=0$ si et seulement s'il n'existe pas d'autre
relation non triviale entre les~$\epsilon_i$ que $\epsilon_0\epsilon_1\epsilon_2\epsilon_3
\epsilon_4=1$.
\end{exemple}

\bigskip
\begin{demo}[ du théorème~\ref{ch3thbr}]%
Posons $K=\kappa(X)$ et fixons une clôture algébrique~$\Ksep$ de~$K$
puis un plongement
$k$\nobreakdash-linéaire $\ksep \hookrightarrow \Ksep$.

\bigskip
\begin{lemme}
\label{ch3brisopic}
La flèche canonique
$H^1(k, \Pic(X_\ksep)) \rightarrow H^1(K,\Pic(X_\Ksep))$
est un isomorphisme.
\end{lemme}

\bigskip
\begin{demo}
Comme~$X$ est une surface de del Pezzo,
les $\Z$\nobreakdash-modules $\Pic(X_\ksep)$ et $\Pic(X_\Ksep)$ sont libres de type fini
et l'application naturelle $\Pic(X_\ksep) \rightarrow \Pic(X_\Ksep)$ est bijective.
Cette application est de plus compatible aux actions respectives de~$\Gal(\ksep/k)$
et~$\Gal(\Ksep/K)$ sur~$\Pic(X_\ksep)$ et~$\Pic(X_\Ksep)$ relativement
au morphisme canonique $\Gal(\Ksep/K)\rightarrow \Gal(\ksep/k)$.
Si~$K\ksep$ désigne l'extension composée de~$K$ et~$\ksep$ dans~$\Ksep$,
le sous-groupe $\Gal(\Ksep/K\ksep)
\subset \Gal(\Ksep/K)$ agit donc trivialement sur~$\Pic(X_\Ksep)$.
Il en résulte que $H^1(K\ksep, \Pic(X_\Ksep))=0$ (puisque $\Pic(X_\Ksep)$ est sans torsion), d'où
un isomorphisme
canonique $H^1(K,\Pic(X_\Ksep))=H^1(\Gal(K\ksep/K),\Pic(X_{K\ksep}))$.
Enfin, comme~$k$ est algébriquement fermé dans~$K$,
le morphisme canonique $\Gal(K\ksep/K) \rightarrow \Gal(\ksep/k)$ est
un isomorphisme, ce qui permet de considérer $\Pic(X_{K\ksep})$ comme
un $\Gal(\ksep/k)$-module.  Compte tenu de ce qui précède, l'application
canonique $\Pic(X_\ksep) \rightarrow \Pic(X_{K\ksep})$ est un isomorphisme
de $\Gal(\ksep/k)$-modules; elle induit donc un isomorphisme en cohomologie.
\end{demo}

\bigskip
\begin{lemme}
\label{ch3brisobr}
La flèche canonique $\Br(X)/\Br(k) \rightarrow \Br(X_K)/\Br(K)$ est un isomorphisme.
\end{lemme}

\bigskip
\begin{demo}
Notons $\Br_1(X)=\Ker(\Br(X)\rightarrow \Br(X_\ksep))$
et $\Br_1(X_K)=\Ker(Br(X_K)\rightarrow \Br(X_\Ksep))$.
Les suites spectrales de Hochschild-Serre pour les extensions~$\ksep/k$
et~$\Ksep/K$ et leurs groupes multiplicatifs respectifs fournissent le diagramme
commutatif
$$
\xymatrix@C=3.5ex{
\Br(k) \ar[d] \ar[r] & \Br_1(X) \ar[d] \ar[r]^(.38){\delta_k} &
H^1(k, \Pic(X_\ksep)) \ar[r] \ar[d]^\gamma & H^3(k, \Gm) \ar[r] \ar[d] & H^3(X,\Gm) \ar[d] \\
\Br(K) \ar[r] & \Br_1(X_K)
\ar[r]^(.38){\delta_K} & H^1(K, \Pic(X_\Ksep)) \ar[r] & H^3(K, \Gm) \ar[r] & H^3(X_K,\Gm)
\rlap{\text{,}}
}
$$
dont les lignes sont exactes.
On a $\Br_1(X)=\Br(X)$ et $\Br_1(X_K)=\Br(X_K)$ puisque la surface~$X$ est rationnelle.
D'après la théorie
du corps de classes,
l'hypothèse que~$k$ est un corps de nombres entraîne l'annulation de $H^3(k,\Gm)$,
d'où la surjectivité de~$\delta_k$.
La flèche $\delta_K$ est également surjective car
le $K$\nobreakdash-point tautologique de~$X$ définit une rétraction
de $H^3(K,\Gm)\rightarrow H^3(X_K,\Gm)$.  Enfin, le lemme~\ref{ch3brisopic}
montre que~$\gamma$ est un isomorphisme. Le résultat s'ensuit.
\end{demo}

\bigskip
L'ensemble~$\Srond$ et les entiers~$n$ et~$d$ ne varient pas
lorsqu'on étend les scalaires de~$k$ à~$K$,
puisque~$k$ est algébriquement fermé dans~$K$.
Grâce au lemme~\ref{ch3brisobr}, il suffit donc de connaître la conclusion du théorème pour la
$K$\nobreakdash-variété~$X_K$ afin de l'obtenir pour la $k$\nobreakdash-variété~$X$.
Étant donné que $X(K) \neq \emptyset$, on a ainsi établi:

\emph{Pour prouver le théorème~\ref{ch3thbr} (et la proposition~\ref{ch3brproduit}),
on peut remplacer l'hypothèse que~$k$ est un corps de nombres par celle
que~$X(k) \neq \emptyset$.}

C'est ce que nous faisons désormais.

\bigskip
\begin{lemme}
Supposons~$X(k) \neq \emptyset$.  Alors la surface~$X$ est $k$\nobreakdash-birationnelle
à un fibré en coniques $\pi \colon C \rightarrow \P^1_k$ (\emph{i.e.}~$C$ est propre, lisse
et géométriquement connexe sur~$k$ et la fibre générique de~$\pi$ est une conique)
dont les fibres au-dessus
de~$\P^1_k \setminus \Srond$ sont lisses et qui vérifie:
pour tout $t \in \Srond$, la fibre~$\pi^{-1}(t)$ est intègre
et la fermeture algébrique de~$\kappa(t)$ dans~$\kappa(\pi^{-1}(t))$
est une extension quadratique ou triviale de~$\kappa(t)$
dont la classe dans $\kappa(t)^\star/\kappa(t)^{\star 2}$
est égale à~$\epsilon_t$.
\end{lemme}

\bigskip
(On appelle \emph{conique} une
courbe propre, lisse, géométriquement connexe et de genre~$0$.)

\bigskip
\begin{demo}
Choisissons $x \in X(k)$. Comme la surface~$X$ est lisse et de codimension~$2$ dans~$\P^4_k$,
il existe un pinceau d'hyperplans de~$\P^4_k$ tangents à~$X$ en~$x$.
Celui-ci découpe sur~$X$ un pinceau de courbes de genre arithmétique~$1$, singulières en~$x$.
On définit $\pi \colon C \rightarrow \P^1_k$ comme
un modèle propre et régulier relativement minimal de ce pinceau.
Il n'est pas difficile de vérifier que~$\pi$ satisfait aux conditions de l'énoncé;
pour les détails, cf.~\cite[p.~61]{ctsansdi}.
\end{demo}

\bigskip
Rappelons un lemme bien connu sur le groupe de Brauer des coniques (on le démontre
par exemple en écrivant la suite spectrale de Hochschild-Serre).

\bigskip
\begin{lemme}
\label{ch3coniquebienconnu}
Soit~$D$ une conique sur un corps~$K$.
La flèche naturelle $\Br(K) \rightarrow \Br(D)$ est surjective et son noyau
est engendré par la classe de l'algèbre de quaternions associée à~$D$.
\end{lemme}

\bigskip
Les groupes~$\Br(C)$ et~$\Br(X)$ sont isomorphes (invariance birationnelle du groupe de Brauer,
cf.~\cite[Corollaire~7.5]{grothbr3}).
De plus, on a $\Br(C)=\Br_\vert(C)$ car la fibre générique de~$\pi$ est une conique
(cf.~lemme~\ref{ch3coniquebienconnu}).
Il reste donc seulement à déterminer la structure du groupe $\Br_\vert(C)/\Br(k)$.

Notons $A \in \Br(\kappa(\P^1_k))$ la classe de l'algèbre de quaternions
définie par la fibre générique de~$\pi$.  Il résulte de la description des
fibres de~$\pi$ que $A \in \Br(\P^1_k \setminus \Srond)$ et que le
résidu de~$A$ en~$t \in \Srond$ est égal à $\epsilon_t \in
\kappa(t)^\star/\kappa(t)^{\star 2} = H^1(\kappa(t),\Z/2) \subset H^1(\kappa(t),\Q/\Z)$.
Notamment, si $A \in \Br(k)$, alors $n=0$. Par ailleurs, dans ce cas,
la surface~$C$ s'écrit comme le produit de~$\P^1_k$ et d'une conique sur~$k$,
ce qui montre que $\Br(C)/\Br(k)=0$ (cf.~\cite[Proposition~2.1.4]{cscrelle94}).
L'énoncé du théorème est donc vérifié lorsque $A\in\Br(k)$;
on suppose dorénavant que l'image de~$A$ dans $\Br(\kappa(\P^1_k))/\Br(k)$ est non nulle.

De la suite exacte de Faddeev (cf.~\cite[§1.2]{cscrelle94}) se déduit la
suite exacte
\begin{equation}
\label{ch3brsexacte}
\xymatrix@C=3ex{
0 \ar[r] & \Br(k) \ar[r] & \Br(\P^1_k \setminus \Srond) \ar[r]^(.41)\alpha &
\displaystyle \bigoplus_{t \in \Srond} H^1(\kappa(t),\Q/\Z)
\ar[r]^(.58)\beta & H^1(k,\Q/\Z) \rlap{\text{,}}
}
\end{equation}
où~$\alpha$ est le produit des résidus aux points de~$\Srond$
et~$\beta$ est la somme des corestrictions de~$\kappa(t)$ à~$k$.
La nullité de l'image de $A \in \Br(\P^1_k \setminus \Srond)$ dans $H^1(k,\Q/\Z)$
est équivalente à l'égalité $\prod_{t\in \Srond}N_{\kappa(t)/k}(\epsilon_t)=1$,
ce qui prouve la proposition~\ref{ch3brproduit}.

Soit~$N$ le noyau du produit
$$
\bigoplus_{t \in \Srond} H^1(\kappa(t),\Q/\Z) \longrightarrow
\left( \bigoplus_{t \in \Srond} H^1(\kappa(t)(\sqrt{\epsilon_t}\;\!),\Q/\Z) \right)
\times H^1(k,\Q/\Z)
$$
des flèches de restriction et de~$\beta$.
D'après \cite[Proposition~1.1.1]{cscrelle94} et la suite exacte~(\ref{ch3brsexacte}),
le groupe~$N$ s'identifie au sous-groupe de $\Br(\P^1_k \setminus \Srond)/\Br(k)$ constitué
des classes dont l'image réciproque par~$\pi$ est non ramifiée sur~$C$. On a donc une
suite exacte
$$
\xymatrix@C=2.7ex{
\!0 \ar[r] & \Ker\!\left(\Br(\kappa(\eta))/\Br(k) \xrightarrow{\;\!\!\pi_\eta^\star\;\!\!\!}
\Br(C_\eta)/\Br(k)\right)\! \ar[r]^(.88)\alpha & N \ar[r] &
\Br_\vert(C)/\Br(k) \ar[r] & 0 \rlap{\text{,}}
}
$$
où~$\eta$ désigne le point générique de~$\P^1_k$.
Le lemme~\ref{ch3coniquebienconnu} et l'hypothèse selon laquelle l'image de~$A$
dans $\Br(\eta)/\Br(k)$ est non nulle entraînent que le groupe de gauche
est égal à~$\Z/2$. Le groupe~$N$ est par ailleurs clairement un $\Z/2$\nobreakdash-espace
vectoriel de dimension~$n-d$; le théorème~\ref{ch3thbr} est donc
démontré.
\end{demo}

\bigskip
\begin{remarque}
On a supposé que~$k$ est un corps de nombres uniquement
pour pouvoir affirmer que $H^3(k,\Gm)=0$ (cf.~preuve du lemme~\ref{ch3brisobr}).
Le théorème~\ref{ch3thbr} et la proposition~\ref{ch3brproduit} sont donc
vrais sous cette hypothèse plus faible, qui est par exemple satisfaite lorsque~$k$
est le corps des fonctions rationnelles en une variable sur un corps de nombres
(cf.~\cite[p.~241]{harduke}).
Sans hypothèse sur~$k$ (autre que celle de caractéristique nulle, énoncée
au début du paragraphe~\ref{ch3notations} et de laquelle on ne se départira pas),
on peut néanmoins déduire de la preuve du théorème~\ref{ch3thbr} d'une
part que
la proposition~\ref{ch3brproduit} reste vraie,
d'autre part que
le groupe $\Br(X)/\Br(k)$ est un $\Z/2$\nobreakdash-espace vectoriel de dimension
\mbox{$\leq n-d-1$}.
En particulier, les corollaires~\ref{ch3brcor541} et~\ref{ch3brtroisgroupes}
sont vrais sans\index{groupe de Brauer!d'une surface de del Pezzo de degré~$4$|)} hypothèse sur~$k$.
\end{remarque}

\subsection{La construction de Swinnerton-Dyer}
\label{ch3parconssd}

Les notations sont les mêmes que dans le paragraphe~\ref{ch3notations}
mais nous supposons de plus que~$t_0$ est $k$\nobreakdash-rationnel.
Comme annoncé dans l'introduction, nous allons construire
une famille de sections hyperplanes de~$X$ paramétrée par une variété $k$\nobreakdash-rationnelle
et telle que la jacobienne de la fibre générique possède un point rationnel d'ordre~$2$.
Nous étudierons ensuite quelques-unes de ses propriétés, en particulier concernant le lieu de ses
fibres singulières.  Certains des résultats ci-dessous seront reprouvés par le
calcul au paragraphe~\ref{ch3calcexpl}; il en va notamment ainsi de
la proposition~\ref{ch3conssd2rat}.  Nous avons cependant jugé les preuves
abstraites concernées suffisamment éclairantes pour être incluses malgré
la redondance.

\bigskip
Considérons d'abord la famille $\pi \colon C \rightarrow (\P^4_k)^\star$
de toutes les sections hyperplanes de~$X$, où~$(\P^4_k)^\star$ désigne l'espace
projectif dual; autrement dit, on pose\glossary{$C$, $\pi$}
$$
C=\bigensemble{(x,L) \in X \times_k (\P^4_k)^\star}{x \in L}
$$
et l'on définit~$\pi$ comme la seconde projection. Cette famille est bien évidemment paramétrée
par une variété $k$\nobreakdash-rationnelle, mais on peut montrer que la
$2$\nobreakdash-torsion rationnelle de la jacobienne de la fibre générique de~$\pi$ est toujours triviale.
Le phénomène sous-jacent, que nous préciserons par la suite
(cf.~preuve de la proposition~\ref{ch3conssd2rat}),
est le suivant: pour $L \in (\P^4_k)^\star$
tel que $X \cap L$ soit lisse, la jacobienne de~$X \cap L$ possède d'autant plus de points
d'ordre~$2$ rationnels que le polynôme homogène~$f_L(\lambda,\mu)$ introduit au
paragraphe~\ref{ch3sectiongeneralites} possède de racines rationnelles; or celui-ci
est même irréductible si~$L$ est générique (cf.~théorème~\ref{ch3propmonod}).  On peut forcer l'existence d'une
racine rationnelle de ce polynôme en étendant les scalaires à un corps de rupture,
ce qui revient, pour~$L$ générique, à considérer la famille obtenue à partir de~$\pi$
par le changement de base $Z \rightarrow (\P^4)^\star$ défini au
paragraphe~\ref{ch3sectiongeneralites}.  Par chance, la famille obtenue sera encore
paramétrée par une variété $k$\nobreakdash-rationnelle (cf.~proposition~\ref{ch3genzkrat}).
Cependant, il ne suffit pas que le polynôme~$f_L(\lambda,\mu)$ ait une racine
dans un corps pour que la jacobienne de~$X \cap L$ admette un point d'ordre~$2$ défini
sur ce corps.  Il suffit en revanche
que le polynôme~$f_L(\lambda,\mu)$ ait deux racines
dans un corps pour que la jacobienne de~$X \cap L$ admette un point d'ordre~$2$ défini
sur ce corps.
On pourrait alors envisager d'utiliser le changement de
base $Z \times_{(\P^4_k)^\star} Z \rightarrow (\P^4_k)^\star$ plutôt
que $Z \rightarrow (\P^4_k)^\star$, mais la variété
$Z \times_{(\P^4_k)^\star} Z$ n'a pas de raison de posséder de composante
irréductible $k$\nobreakdash-rationnelle (ou même rationnelle)
autre que celle de la diagonale.
Le même problème se pose si l'on considère directement le changement
de base de~$(\P^4_k)^\star$ à la variété des points d'ordre~$2$ des jacobiennes
des fibres de~$\pi$.

Une autre manière de forcer l'existence d'une racine rationnelle du polynôme~$f_L(\lambda,\mu)$
est de ne considérer que les hyperplans~$L$ contenant $\fp{P}0 \in \P^4(k)$.
(Le point~$\fp{P}0$ est $k$\nobreakdash-rationnel puisque~$t_0$ l'est par hypothèse.)
Pour un tel hyperplan, la quadrique $Q_{t_0} \cap L$ est en effet singulière,
ce qui signifie que~$t_0$ est racine de~$f_L(\lambda,\mu)$. La famille
ainsi obtenue est bien paramétrée par une variété $k$\nobreakdash-rationnelle,
même un espace projectif; mais là encore, une
seule racine ne suffit pas à assurer l'existence d'un point d'ordre~$2$ dans
la jacobienne.  Rien n'empêche toutefois de combiner les deux approches qui viennent
d'être décrites.  Telle est l'idée que nous allons dès à présent mettre en \oe{}uvre.

Soit $H \subset \P^4_k$\glossary{$H$, $\Lambda$, $Z$, $\rho$} l'unique hyperplan $k$\nobreakdash-rationnel contenant $\uplet{\fp{P}1}{\fp{P}4}$.
Posons $$\Lambda = \bigensemble{L \in (\P^4_k)^\star}{\fp{P}0 \in L}$$
et
$$
Z=\bigensemble{(t,L)\in \P^1_k \times_k \Lambda}{ Q_t \cap L \cap H \text{ n'est pas lisse}}
$$
et notons $\rho \colon Z \rightarrow \Lambda$ le morphisme induit par la seconde projection.
Comme $\fp{P}0 \not\in H$ (cf.~corollaire~\ref{ch3genpiengendrent}), il
existe un isomorphisme canonique $\Lambda=H^\star$
faisant correspondre à un hyperplan de~$\P^4_k$ contenant~$\fp{P}0$ sa trace sur~$H$. Le morphisme
$Z \xrightarrow{\;\rho\;} \Lambda$ s'interprète, \emph{via} cet isomorphisme, comme le revêtement
de~$H^\star$ associé
au pinceau de quadriques $(Q_t \cap H)_{t \in \P^1_k}$ de~$H$
par la construction générale du paragraphe~\ref{ch3sectiongeneralites}.
En particulier, le morphisme~$\rho$ est fini et plat, de degré~$3$
(cf.~proposition~\ref{ch3genpfiniplat})
et la variété~$Z$ est $k$\nobreakdash-rationnelle
(cf.~proposition~\ref{ch3genzkrat}).
Plus précisément, si $E=(H \cap X) \cup \{\uplet{\fp{P}1}{\fp{P}4}\}$\glossary{$E$, $H^0$,
$\sigma$, $\pi_\Lambda$, $\pi_Z$, $\pi_{H^0}$}
et $H^0 = H \setminus E$,
les commentaires qui suivent la proposition~\ref{ch3genzkrat} définissent un morphisme
birationnel $\sigma \colon H^0 \rightarrow Z$.
Notons $\pi_\Lambda$, $\pi_Z$ et~$\pi_{H^0}$ les morphismes déduits de~$\pi$ par changement de base
de~$(\P^4_k)^\star$ à~$\Lambda$, $Z$ et à~$H_0$.
Le diagramme commutatif suivant, dont tous les carrés sont cartésiens, résume la situation.
(La flèche oblique est induite par la première projection.)
\begin{equation}
\label{ch3conssddiag}
\begin{aligned}
\myxyin\myxyhook\xymatrix@R=2.5em@C=6ex{
*+!!<0pt,\the\fontdimen22\textfont2>!<-5mm,6mm>{X} & \ar[l]!(5,-6) C \ar[d]^\pi & \ar@{_{ (}->}[l] C_\Lambda \ar[d]^{\pi_\Lambda} & \ar[l] C_Z
\ar[d]^{\pi_Z} & C_{H^0} \ar[l] \ar[d]^{\pi_{H^0}} \\
& (\P^4_k)^\star & \ar@{_{ (}->}[l] \Lambda & \ar[l]_(.47)\rho Z & \ar[l]_\sigma H^0
}\myxyout
\end{aligned}
\end{equation}

\bigskip
\begin{proposition}
\label{ch3conssdfibrespi}
Toutes les fibres de~$\pi$ (et donc de~$\pi_\Lambda$, $\pi_Z$ et~$\pi_{H^0}$)
sont des courbes géométriquement connexes de genre arithmétique~$1$.
\end{proposition}

\bigskip
\begin{demo}
La surface~$X$ est géométriquement intègre puisque c'est une intersection
complète lisse dans~$\P^4_k$ de dimension~$\geq 1$. De plus, elle n'est
contenue dans aucun hyperplan de~$\P^4_k$ (cf.~\cite[Lemma~1.3 (i)]{ctsansdi}).
Chaque section hyperplane de~$X$ est donc
une intersection de dimension~$1$ de deux quadriques dans un
espace projectif de dimension~$3$.  En particulier, les fibres de~$\pi$ sont des courbes
géométriquement connexes de genre arithmétique~$1$.
\end{demo}

\bigskip
\begin{proposition}
\label{ch3conssdreguliers}
Les variétés~$C$ et~$C_\Lambda$ sont lisses et géométriquement connexes sur~$k$,
de dimensions respectives~$5$ et~$4$.
\end{proposition}

\bigskip
\begin{demo}
Étant donné que $\fp{P}0 \not\in X$, les projections naturelles $C \rightarrow X$
et $C_\Lambda \rightarrow X$ sont des fibrés projectifs (localement
libres) de dimensions relatives respectivement~$3$ et~$2$.
L'assertion résulte donc de la lissité et de la connexité géométrique
de la surface~$X$ sur~$k$.
\end{demo}

\bigskip
\begin{proposition}
\label{ch3conssd2rat}
La fibre générique de~$\pi_{H^0}$ est une courbe lisse et géométriquement connexe
de genre~$1$ dont la jacobienne possède un point d'ordre~$2$ rationnel.
\end{proposition}

\bigskip
\begin{demo}
Pour établir que la fibre générique de~$\pi_{H^0}$ est une courbe lisse et géométriquement
connexe de genre~$1$, il suffit de montrer que la fibre générique de~$\pi_\Lambda$
vérifie ces propriétés.  Compte tenu de la proposition~\ref{ch3conssdfibrespi},
il suffit de vérifier que la fibre générique de~$\pi_{\Lambda}$ est lisse,
mais cela découle de la régularité de~$C_\Lambda$ (cf.~proposition~\ref{ch3conssdreguliers}).

Les morphismes~$\pi_{H^0}$ et~$\pi_Z$ ont bien évidemment même fibre générique.
Les lemmes suivants vont nous permettre de montrer que la jacobienne de la fibre générique
de~$\pi_Z$ possède un point d'ordre~$2$ rationnel.

\bigskip
\begin{lemme}
\label{ch3conssdlemmerescub}
Soient~$k$ un corps de caractéristique différente de~$2$ et de~$3$ et~$C \subset \P^3_k$
une courbe lisse, intersection de deux quadriques.  Soit $(Q_t)_{t \in \P^1_k}$
le pinceau que ces deux quadriques engendrent.
Les valeurs de~$t$ pour
lesquelles la quadrique $Q_t \subset \P^3_k$ n'est pas lisse définissent
un $k$\nobreakdash-schéma étale de degré~$4$ (cf.~proposition~\ref{ch3genpropeq}), que l'on
note~$T$.  Soient respectivement~$A$ et~$B$
les anneaux des schémas affines~$T$ et $\tors{2}J \setminus \{0\}$,
où~$J$ désigne la jacobienne de~$C$.
Alors la $k$\nobreakdash-algèbre~$B$ est isomorphe à la résolvante\index{résolvante cubique|(} cubique de~$A$.
\end{lemme}

\bigskip
(Si~$A$ est une algèbre étale de degré~$4$ sur un corps~$k$, la \emph{résolvante cubique}~$B$
de~$A$
est une $k$\nobreakdash-algèbre étale de degré~$3$ qu'il est possible de construire
canoniquement à partir de~$A$ mais que nous nous contentons ici de définir à isomorphisme
près.
L'action par conjugaison du groupe symétrique~$\mathfrak{S}_4$
sur l'ensemble de ses trois sous-groupes de $2$\nobreakdash-Sylow fournit un morphisme
de groupes
$\mathfrak{S}_4 \rightarrow \mathfrak{S}_3$ (si l'on accepte de numéroter les $2$\nobreakdash-Sylow).
Celui-ci induit une application $\rho \colon
H^1(k,\mathfrak{S}_4) \rightarrow H^1(k,\mathfrak{S}_3)$
entre ensembles pointés de cohomologie galoisienne non abélienne.
Pour tout~$n$, l'ensemble $H^1(k,\mathfrak{S}_n)$ classifie les $k$\nobreakdash-algèbres étales
de degré~$n$ à isomorphisme près.  La classe d'isomorphisme
de~$B$ est par définition l'image par~$\rho$ de la classe
d'isomorphisme de~$A$.
Voir \cite{knustignol} pour plus de détails sur ces notions
et pour la comparaison
avec le cas classique où~$A$ est un corps.  Voir aussi \cite[p.~10]{birchsd}
pour des formules explicites tout à fait générales\footnote{On peut aussi
définir directement la résolvante cubique de~$A$ au moyen de ces formules, au moins lorsque~$A$
possède un élément primitif, ce qui est toujours le cas si~$k$ est infini.  L'un des inconvénients
de cette approche est qu'il n'est pas clair
que la classe d'isomorphisme de l'algèbre obtenue ne dépend pas de l'élément primitif choisi.
Cela dit, pour la situation qui nous intéresse, cela n'a aucune importance puisque l'on dispose
d'un élément primitif privilégié.}
permettant d'exprimer~$B$
en fonction du polynôme minimal d'un élément primitif de~$A$.

Si l'on part d'un polynôme séparable $f \in k[t]$ de degré~$\leq 4$, on appellera
\emph{résolvante cubique de~$f$} tout\footnote{Nous prenons cette précaution
car dans la littérature se trouvent réellement plusieurs définitions explicites
et contradictoires
du polynôme «~résolvante cubique~» associé à un polynôme de degré~$4$.} polynôme non nul~$g\in k[t]$
de degré~$\leq 3$
tel que la $k$\nobreakdash-algèbre $k[t]/(g(t)) \times k^{3-\deg(g)}$
soit isomorphe à la résolvante cubique de la $k$\nobreakdash-algèbre $k[t]/(f(t)) \times k^{4-\deg(f)}$.

Une propriété importante de la résolvante cubique est la suivante. Notons $(k_i)_{i\in I}$
et $(k'_j)_{j \in J}$ des familles finies d'extensions finies séparables de~$k$
telles que $A=\prod_{i\in I}k_i$ et $B=\prod_{j\in J}k'_j$.
Si~$K/k$ est une extension galoisienne dans laquelle~$k_i/k$
se plonge pour tout~$i \in I$, alors
$k'_j/k$ se plonge aussi dans~$K/k$ pour tout~$j\in J$.
Cette propriété
entraîne notamment que~$\Spec(B)$ possède un point $k$\nobreakdash-rationnel
si $\Spec(A)$ en possède deux.
Elle résulte immédiatement de la définition
cohomologique de la résolvante cubique, compte tenu de la remarque suivante:
si~$K/k$ est une extension galoisienne, le «~noyau~» de l'application pointée de restriction
$H^1(k,\mathfrak{S}_n) \rightarrow H^1(K,\mathfrak{S}_n)$ classifie à isomorphisme
près les $k$\nobreakdash-algèbres
étales~$A$ de degré~$n$ telles que $A \otimes_k K$ soit $K$\nobreakdash-isomorphe à~$K^n$.)

\bigskip
\begin{demo}[ du lemme~\ref{ch3conssdlemmerescub}]%
Soient~$q_1$ et~$q_2$ des formes quadratiques associées à des quadriques distinctes
dans le pinceau considéré.
Le $k$\nobreakdash-schéma~$T$ est défini par l'équation homogène
$\det(\lambda q_1 + \mu q_2)=0$, le choix d'une base
étant sous-entendu (cf.~paragraphe~\ref{ch3sectiongeneralites}).
Notons $\PPic_{C/k}^2$ la composante de degré~$2$ du foncteur de Picard relatif de~$C$
sur~$k$.
Il est bien connu d'une part que~$\PPic_{C/k}^2$ est $k$\nobreakdash-birationnel
à la courbe affine d'équation
$$
y^2=\det(\lambda q_1 + q_2)
$$
(cf.~par exemple~\cite{siksek}),
d'autre part que la jacobienne d'un modèle propre et lisse de la
courbe de genre~$1$ définie par cette équation
a pour équation de Weierstrass $y^2=g(\lambda)$,
où~$g(\lambda)$ est une résolvante\index{résolvante cubique|)} cubique (unitaire de degré~$3$) du polynôme $\det(\lambda q_1 + q_2)$
(cf.~\cite[p.~10]{birchsd}).
Les courbes~$C$ et~$\PPic_{C/k}^2$ ayant même jacobienne, le lemme s'ensuit.
\end{demo}

\bigskip
\begin{lemme}
\label{ch3conssdmemez}
Pour $(t,L)\in \P^1_k \times_k \Lambda$, on a $(t,L)\in Z \cup (\{t_0\} \times\Lambda)$ si et seulement
si $Q_t \cap L$ n'est pas lisse.
\end{lemme}

\bigskip
\begin{demo}
La question étant de nature géométrique, on peut supposer~$k$ algébriquement clos.
D'après la proposition~\ref{ch3gensimdiag}, on peut alors choisir une base
dans laquelle~$q_1$ et~$q_2$ sont simultanément diagonales.  Un calcul explicite
très facile utilisant le critère jacobien permet de conclure.
\end{demo}

\bigskip
Nous pouvons maintenant terminer la démonstration de la proposition~\ref{ch3conssd2rat}.
Notons~$L$ l'hyperplan $\kappa(Z)$\nobreakdash-rationnel
de~$\P^4_k$ défini par l'image du point générique de~$Z$
par la seconde projection $Z \rightarrow \Lambda$.
D'après le lemme~\ref{ch3conssdlemmerescub} et les remarques qui le suivent, tout ce
qu'il reste à vérifier est que le polynôme homogène
\mbox{$f_L(\lambda,\mu) \in \kappa(Z)[\lambda,\mu]$}
s'annule en au moins deux points distincts de~$\P^1(\kappa(Z))$; or
il s'annule en $t_0 \in \P^1(k)$ puisque $\fp{P}0 \in L$, et le lemme~\ref{ch3conssdmemez}
montre qu'il s'annule aussi en l'image dans~$\P^1(\kappa(Z))$ du point
générique de~$Z$ par la première projection~\mbox{$Z \rightarrow \P^1_k$}.
\end{demo}

\bigskip
Nous aurons besoin de quelques informations supplémentaires sur la fibration $\pi_{H^0}$.

\bigskip
\begin{proposition}
\label{ch3conssdplatcm}
Les morphismes~$\pi$, $\pi_\Lambda$, $\pi_Z$ et~$\pi_{H^0}$ sont plats
et la variété~$C_{H^0}$ est de Cohen-Macaulay.
\end{proposition}

\bigskip
\begin{demo}
Il résulte des propositions~\ref{ch3conssdfibrespi} et~\ref{ch3conssdreguliers}
et de \cite[15.4.2]{ega43} que le morphisme~$\pi$ est plat et que ses fibres
sont de Cohen-Macaulay.
Les morphismes~$\pi_\Lambda$, $\pi_Z$ et~$\pi_{H^0}$ sont
donc eux aussi plats et à fibres de Cohen-Macaulay (cf.~\cite[6.7.2]{ega42}).
Cela suffit à assurer que la variété~$C_{H^0}$ est de
Cohen-Macaulay, puisque~$H^0$ est lui-même de Cohen-Macaulay (cf.~\cite[6.3.5]{ega42}).
\end{demo}

\bigskip
\begin{proposition}
\label{ch3conssdgeomint}
La variété~$C_{H^0}$ est géométriquement intègre sur~$k$.
\end{proposition}

\bigskip
\begin{demo}
On peut supposer~$k$ algébriquement clos.
Le morphisme~$\pi_{H^0}$ est propre et plat
(cf.~proposition~\ref{ch3conssdplatcm}), donc ouvert et fermé.
Par conséquent, chaque composante connexe de~$C_{H^0}$ rencontre
la fibre générique de~$\pi_{H^0}$.  Celle-ci étant
connexe (cf.~proposition~\ref{ch3conssdfibrespi}), toute composante
connexe de~$C_{H^0}$ contient la fibre générique de~$\pi_{H^0}$;
d'où la connexité de~$C_{H^0}$.

Il existe une unique composante irréductible de~$C_{H^0}$ qui rencontre
la fibre générique de~$\pi_{H^0}$, puisque cette dernière
est irréductible (cf.~proposition~\ref{ch3conssd2rat}).  Notons-la~$I_1$.
Supposons qu'il existe une composante irréductible~$I_2$ de~$C_{H^0}$ distincte de~$I_1$.
Interprétant la dimension d'une variété irréductible par
le degré de transcendance du corps résiduel de son point générique,
on voit tout de suite
que $\dim(I_1)=4$ et que $\dim(I_2) \leq \dim(\pi_{H^0}(I_2))+1$,
compte tenu que~$H^0$ est de dimension~$3$ et que chaque fibre de~$\pi$
est de dimension~$1$.
Comme $I_1 \neq I_2$, la composante irréductible~$I_2$ ne rencontre
pas la fibre générique de~$\pi_{H^0}$, par définition de~$I_1$.
On a donc $\dim(\pi_{H^0}(I_2))<3$, d'où $\dim(I_2)<\dim(I_1)$.
Comme la variété~$C_{H^0}$
est connexe et de Cohen-Macaulay (cf.~proposition~\ref{ch3conssdplatcm}),
elle est équidimensionnelle (cf.~\cite[0.16.5.4]{ega41}), ce qui est
contradictoire avec l'inégalité précédente.  Ainsi a-t-on prouvé que~$C_{H^0}$ est irréductible.

Il reste seulement à établir que la variété irréductible~$C_{H^0}$ est réduite.
Comme elle est de Cohen-Macaulay, elle vérifie la propriété~$(\mathrm{S}_1)$,
et il suffit donc de s'assurer que
son anneau local au point générique est réduit (critère $(\mathrm{R}_0)+(\mathrm{S}_1)$,
cf.~\cite[5.8.5]{ega42}); mais ceci est clair
puisque la fibre générique
de~$\pi_{H^0}$ est lisse (cf.~proposition~\ref{ch3conssd2rat}).
\end{demo}

\bigskip
Intéressons-nous maintenant au lieu des fibres singulières de~$\pi_\Lambda$.
En d'autres termes, il s'agit de déterminer quels hyperplans de~$\P^4_k$ contenant~$\fp{P}0$
sont tangents à~$X$.

\bigskip
\begin{lemme}
\label{ch3conssdxhlisse}
La variété~$X \cap H$ est une courbe lisse et géométriquement connexe de genre~$1$.
\end{lemme}

\bigskip
\begin{demo}
On peut supposer~$k$ algébriquement clos puis choisir
une base dans laquelle~$q_1$ et~$q_2$ sont simultanément diagonales
(cf.~proposition~\ref{ch3gensimdiag}).  Il est alors facile de déduire par le calcul
la lissité de~$X \cap H$ de celle de~$X$ (en utilisant au choix
la proposition~\ref{ch3genpropeq} ou le critère jacobien).
Que cette variété soit une courbe géométriquement connexe de genre arithmétique~$1$
est déjà connu puisque c'est une fibre de~$\pi$ (cf.~proposition~\ref{ch3conssdfibrespi}).
\end{demo}

\bigskip
Soit $L \in \Lambda$. Le lemme~\ref{ch3conssdmemez} montre que
l'ensemble des racines de~$f_L$ est égal à la réunion de~$\{t_0\}$ et de l'ensemble
des racines de~$f_{L \cap H}$, où $f_{L\cap H}(\lambda,\mu)=\det((\lambda q_1 + \mu q_2)|_{L \cap H})$.
Le polynôme~$f_L$ est donc non nul et séparable si et seulement si $f_{L \cap H}(t_0) \neq 0$ et
que~$f_{L\cap H}$ est non nul et séparable.
Comme les variétés~$X \cap L$ et~$X \cap L \cap H$ sont purement de codimension~$2$
respectivement dans~$L$ et dans~$L \cap H$ (proposition~\ref{ch3conssdfibrespi} et lemme~\ref{ch3conssdxhlisse}),
il en résulte, grâce à la proposition~\ref{ch3genpropeq}, que~$X \cap L$ est singulier si et seulement
si \mbox{$Q_{t_0} \cap L \cap H$} est singulier ou $X \cap L \cap H$ est singulier.
Autrement dit, la fibre de~$\pi_\Lambda$ en~$L$ est singulière si et seulement
si l'hyperplan \mbox{$L \cap H$} de~$H$ est tangent à \mbox{$Q_{t_0} \cap H$} ou à \mbox{$X \cap H$}.
Selon la terminologie usuelle, cela signifie que le lieu des fibres singulières de~$\pi_\Lambda$
s'identifie, \emph{via} l'isomorphisme canonique $\Lambda=H^\star$, à la réunion des
variétés duales $(Q_{t_0} \cap H)^\star \subset H^\star$ et $(X \cap H)^\star \subset H^\star$.
Posons $Q=(Q_{t_0}\cap H)^\star$ et $R=(X \cap H)^\star$\glossary{$Q$, $R$}.

\bigskip
\begin{proposition}
\label{ch3conssdrqhypirr}
Les sous-variétés~$Q$ et~$R$ de~$\Lambda$ sont des hypersurfaces géométriquement irréductibles distinctes
(de degrés respectifs~$2$ et~$8$, mais cela ne nous servira pas).
\end{proposition}

\bigskip
\begin{demo}
Un calcul facile montre que la variété duale d'une hypersurface quadrique lisse est une hypersurface quadrique lisse
(choisir une base dans laquelle la quadrique est donnée par une forme quadratique diagonale et exprimer
le discriminant de sa trace sur un hyperplan variable).
Ceci prouve la partie du lemme portant sur~$Q$.  Intéressons-nous maintenant à~$R$.

De manière générale, si $Y \subset H$ est une sous-variété lisse et géométriquement intègre,
le lieu des hyperplans $h \in H^\star$ tels que~$h$ soit tangent à~$Y$ est un fermé géométriquement
irréductible de~$H^\star$.  En effet, si l'on note ce lieu~$Y^\star$ et que l'on pose
$$\Ytilde=\bigensemble{(y,h)\in Y \times_k H^\star}{\ft{T}y Y \subset h}\!\rlap{\text{,}}$$
la seconde projection induit un morphisme propre $\Ytilde \rightarrow H^\star$ dont l'image
est égale à~$Y^\star$, et la variété $\Ytilde$ est géométriquement intègre car la première projection
$\Ytilde \rightarrow Y$ en fait un fibré projectif (localement libre) sur~$Y$.
Comme $\Ytilde \rightarrow Y$ est de dimension relative \mbox{$\dim(H)-\dim(Y)-1$}, on voit
de plus que $Y^\star$ est une hypersurface de~$H$ si et seulement si le morphisme $\Ytilde \rightarrow
Y^\star$ est génériquement fini, où l'on munit $Y^\star$ de sa structure de sous-variété fermée réduite de~$H^\star$.

Dans la situation qui nous intéresse, c'est-à-dire pour $Y=X \cap H$, le morphisme $\Ytilde \rightarrow Y^\star$ est génériquement
fini car~$Y$ est une courbe et que ce n'est pas une droite (cf.~\cite[p.~174]{kleimanvanc}; le «~codéfaut de dualité~» ne peut
valoir que~$1$).  La sous-variété $R \subset H^\star$ est donc bien une hypersurface.
Que~$Q$ et~$R$ soient des sous-variétés distinctes découle de la propriété de réciprocité, toujours
valide en caractéristique nulle (critère de Monge-Segre, cf.~\cite[Theorem~4]{kleimanvanc}),
selon laquelle on a $Q^\star = Q_{t_0} \cap H$ et $R^\star = X \cap H$ \emph{via} l'identification
canonique~\mbox{$H^{\star\star}=H$}.

Lorsque~$Y^\star$ est une hypersurface, il existe une formule pour son degré
en termes du degré du morphisme $\Ytilde \rightarrow Y^\star$ et des
caractéristiques d'Euler-Poincaré de~$Y$, d'une section hyperplane lisse de~$Y$ et d'une section de~$Y$
par un sous-espace linéaire de~$H$ de codimension~$2$ (cf.~\cite[Proposition~5.7.2]{katzsga}).
Cette formule permet de vérifier l'assertion facultative que l'hypersurface~$R$ est de degré~$8$, compte tenu du lemme~\ref{ch3conssdxhlisse}
et de ce que le morphisme $\Ytilde \rightarrow Y^\star$ est birationnel, étant génériquement fini
et~$k$ étant de caractéristique nulle (cf.~\cite[Proposition~15]{kleimanvanc}).
\end{demo}

\bigskip
\begin{proposition}
\label{ch3conssdi1i2}
La fibre géométrique de~$\pi_\Lambda$ au-dessus du point générique de~$Q$ est réduite et possède deux
composantes irréductibles; la classe dans $\kappa(Q)^\star/\kappa(Q)^{\star 2}$
de l'extension quadratique ou triviale minimale de~$\kappa(Q)$ par laquelle se
factorise l'action du groupe de Galois absolu de~$\kappa(Q)$ sur ces deux composantes
est égale à~$\epsilon_0$.

La fibre géométrique de~$\pi_\Lambda$ au-dessus du point générique de~$R$
est intègre.
\end{proposition}

\bigskip
\begin{remarque}
La proposition ci-dessus révèle que
chaque fibre de~$\pi_\Lambda$ au-dessus d'un point de codimension~$1$ de~$\Lambda$
est déployée par une extension \emph{constante} des scalaires, puisque $\epsilon_0 \in k^\star/k^{\star 2}$.
Cette propriété était prévisible.  En effet, si~$k$ est algébriquement clos,
la surface~$X$ contient une droite (et même seize), et toute droite de~$X$ définit naturellement
une section rationnelle du morphisme~$\pi_\Lambda$.  Une telle section est nécessairement définie en codimension~$1$
et rencontre une unique composante irréductible
de chaque fibre de~$\pi_\Lambda$ au-dessus de son domaine de définition,
car la variété~$C$ est lisse (cf.~proposition~\ref{ch3conssdreguliers}).
\end{remarque}

\bigskip
\begin{demo}[ de la proposition~\ref{ch3conssdi1i2}]%
Quitte à étendre les scalaires de~$k$ à~$\kappa(Q)$ et à~$\kappa(R)$,
et dans le seul but de simplifier les notations, on peut
se contenter de prouver les assertions pour la fibre de~$\pi_\Lambda$ au-dessus
de tout point rationnel de~$Q \cup R$ appartenant à un ouvert dense suffisamment
petit (et dont le lecteur vérifiera immédiatement qu'il provient par extension
des scalaires d'un ouvert défini sur le corps~$k$ initial).

Soit donc~$L$ un hyperplan $k$\nobreakdash-rationnel de~$\P^4_k$ correspondant à un point de~$Q$.
La proposition~\ref{ch3conssdrqhypirr} montre que $Q \cap R \neq Q$.  On peut donc supposer que~$L$ n'est pas
tangent à~$X \cap H$, ce qui revient à supposer le polynôme $f_{L \cap H}$ non nul et séparable
(cf.~lemme~\ref{ch3conssdxhlisse} et proposition~\ref{ch3genpropeq}).
D'après la définition de~$Q$, l'hyperplan~$L$ est tangent à~$Q_{t_0}$ en
un point lisse de~$Q_{t_0}$.  Comme~$Q_{t_0}$ est une quadrique de dimension~$3$
admettant un unique point singulier, cela entraîne que $Q_{t_0} \cap L$ est
géométriquement une réunion de deux plans distincts.  Ceux-ci sont permutés par l'extension
quadratique $k(\sqrt{\epsilon_0}\;\!)/k$, par définition de~$\epsilon_0$.
Pour terminer de prouver la première partie de la proposition,
il suffit de vérifier que pour $t \in \P^1(\ksep) \setminus \{t_0\}$, la trace
de~$Q_t$ sur chacun des deux plans contenus dans $Q_{t_0} \cap L$ est une
conique lisse dans ce plan.  Si tel n'était pas le cas, la courbe~$X \cap L$
contiendrait géométriquement une droite, ce qui contredirait
la proposition~\ref{ch3gensousvarlin} puisque le polynôme~$f_L$
admet une racine simple (et même deux) d'après
l'hypothèse que~$f_{L \cap H}$ est
non nul et séparable
et la description des racines de~$f_L$ en fonction
de celles de~$f_{L \cap H}$ (cf.~le paragraphe qui suit le lemme~\ref{ch3conssdxhlisse}).

Soit maintenant~$L$ un hyperplan $k$\nobreakdash-rationnel de~$\P^4_k$ correspondant à un point de~$R$.
La proposition~\ref{ch3conssdrqhypirr} montre que $Q \cap R \neq R$.  On peut donc supposer que~$L$
n'est pas tangent à~$Q_{t_0} \cap H$, ce qui revient à supposer que \mbox{$f_{L\cap H}(t_0) \neq 0$}.
Par ailleurs, on peut supposer que~$L$ ne contient aucun des~$\fp{P}i$ pour $i \in \{\uplet{1}{4}\}$.

\bigskip
\begin{lemme}
\label{ch3conssdlxrangtrois}
Pour tout $(\lambda,\mu)\in \ksep^2 \setminus \{(0,0)\}$,
la forme quadratique $(\lambda q_1 + \mu q_2)|_{L \otimes_k \ksep}$
sur le $\ksep$\nobreakdash-espace vectoriel de dimension~$4$ associé à $L \otimes_k \ksep$
est de rang~$\geq 3$.
\end{lemme}

\bigskip
\begin{demo}
Restreindre une forme quadratique à un hyperplan ne peut abaisser son rang qu'au plus de~$2$.
La conclusion du lemme est donc satisfaite si $[\lambda:\mu] \not\in \{\uplet{t_0}{t_4}\}$.
Pour $[\lambda:\mu]=t_i$ avec $i \in \{\uplet{1}{4}\}$, la forme quadratique
$(\lambda q_1 + \mu q_2)|_{L \otimes_k \ksep}$ est non dégénérée puisque $\fp{P}i \not\in L$;
elle est donc de rang~$4$.
Enfin, pour $[\lambda:\mu]=t_0$, son rang est égal à~$3$ puisque l'hyperplan~$L$ n'est pas tangent
à~$Q_{t_0}$ selon une droite, n'étant pas tangent à~$Q_{t_0} \cap H$.
\end{demo}

\bigskip
Si la fibre de~$\pi_\Lambda$ au-dessus du point de~$\Lambda$ défini par~$L$
n'était pas géométriquement
intègre, la variété $X \cap L$, qui est une intersection de deux quadriques dans~$L$ purement
de codimension~$2$ dans~$L$ (proposition~\ref{ch3conssdfibrespi}), contiendrait géométriquement
soit deux composantes irréductibles de degré~$2$, soit une composante irréductible de degré~$2$
et de multiplicité~$2$, soit au moins une composante irréductible de degré~$1$.
Les deux premières alternatives sont exclues à cause du lemme~\ref{ch3conssdlxrangtrois}
(cf.~\cite[Lemma~1.7 et Lemma~1.10]{ctsansdi}).  La dernière l'est à cause de la proposition~\ref{ch3gensousvarlin},
compte tenu que~$t_0$ est une racine simple de~$f_L$, vu la description des racines de~$f_L$ en fonction
de celles de~$f_{L \cap H}$ et l'hypothèse que $f_{L\cap H}(t_0)\neq 0$.
\end{demo}

\bigskip
\begin{proposition}
\label{ch3conssdrhor}
Le lieu de branchement de $\rho \colon Z \rightarrow \Lambda$ est égal à~$R$.
La fibre de~$\rho$ au-dessus du point générique de~$R$ comporte un point double et un point de multiplicité~$1$.
\end{proposition}

\bigskip
\begin{demo}
D'après la définition de~$\rho$ et la description de ses fibres géométriques (cf.~les remarques
qui précèdent la proposition~\ref{ch3genpfiniplat}), un point de~$\Lambda$ appartient au lieu
de branchement de~$\rho$ si et seulement si le polynôme~$f_{L \cap H}$ n'est pas séparable,
où~$L$ désigne l'hyperplan de~$\P^4_k$ correspondant au point considéré.
Cette condition équivaut bien à l'appartenance de ce point à~$R$ (cf.~les remarques
qui précèdent la proposition~\ref{ch3conssdrqhypirr}).

Pour prouver la seconde assertion, on peut supposer~$k$ algébriquement clos,
auquel cas~$\fp{P}1$ et~$\fp{P}2$ sont des points rationnels de~$H$.
Le sous-espace linéaire de~$H^\star$ constitué des hyperplans de~$H$ qui contiennent~$\fp{P}1$ et~$\fp{P}2$
est de dimension~$1$.  Son intersection avec l'hypersurface~$R$ est donc non vide: il existe
un hyperplan $L \subset \P^4_k$ contenant~$\fp{P}0$, $\fp{P}1$ et~$\fp{P}2$ et tangent à~\mbox{$X \cap H$}.
Le polynôme~$f_{L \cap H}$ s'annule alors en~$t_1$ et en~$t_2$. En particulier, il possède deux racines
distinctes; autrement dit, la fibre de~$\rho$ au-dessus du point de~$R$ défini par~$L$
possède au moins deux points géométriques.
D'après la semi-continuité inférieure du nombre géométrique de points dans
les fibres d'un morphisme fini et plat (cf.~\cite[15.5.1~(i)]{ega43}),
la fibre de~$\rho$ au-dessus du point générique de~$R$ possède nécessairement elle aussi au moins deux points
géométriques,
ce qui prouve la proposition.
\end{demo}

\bigskip
Combinant les propositions~\ref{ch3conssdi1i2} et~\ref{ch3conssdrhor} et la description
explicite du morphisme~$\sigma$,
nous pouvons maintenant essentiellement déterminer le lieu des fibres singulières de~$\pi_{H^0}$, au moins
en codimension~$1$, ainsi que la structure des fibres correspondantes.

\bigskip
\begin{theoreme}
\label{ch3conssdtheoreme}
(i) Il existe un unique point de codimension~$1$ de~$H^0$ au-dessus duquel
la fibre de~$\pi_{H^0}$ est singulière et géométriquement intègre.

(ii) Pour tout $t \in \Srond'$, la fibre géométrique de~$\pi_{H^0}$ au-dessus du
point générique de~$Q_t \cap H$ est réduite et possède deux composantes
irréductibles; la classe dans $\kappa(Q_t \cap H)^\star/\kappa(Q_t \cap H)^{\star 2}$
de l'extension quadratique ou triviale minimale par laquelle se factorise l'action
du groupe de Galois absolu de \mbox{$\kappa(Q_t \cap H)$} sur ces deux composantes
est égale à~$\epsilon_t$.

(iii) Pour tout $h \in H^0$ de codimension~$1$ dans~$H^0$ n'appartenant à~$Q_t$
pour aucun $t \in \Srond'$ et tel que la fibre géométrique de~$\pi_{H^0}$
au-dessus de~$h$ ne soit pas intègre, celle-ci est réduite et possède deux composantes irréductibles;
la classe dans $\kappa(h)^\star/\kappa(h)^{\star 2}$ de l'extension
quadratique ou triviale minimale de~$\kappa(h)$ par laquelle se factorise
l'action du groupe de Galois absolu de~$\kappa(h)$ sur ces deux composantes
est égale à l'image de~$\epsilon_0$.
\end{theoreme}

\bigskip
\begin{demo}
Notons $\tau \colon Z \dashrightarrow H$ l'application rationnelle inverse de~$\sigma$.
La restriction du morphisme $\sigma \colon H^0 \rightarrow Z$
à l'ouvert $H^0 \setminus (\bigcup_{t \in \Srond'} (Q_t \cap H^0))$
est une immersion ouverte, comme il résulte de la description de~$\sigma$ et de~$\tau$
(cf.~proposition~\ref{ch3genzkrat}).
On en déduit, vu que~$\rho$ est fini et plat, que les seuls points de
codimension~$1$ de~$H^0$ dont les images par~$\rho \circ \sigma$ ne sont pas des points de codimension~$1$ de~$\Lambda$
sont les points génériques de $Q_t \cap H$ pour $t \in \Srond'$.  Notons ceux-ci $(\xi_t)_{t \in \Srond'}$.
Les seuls points de codimension~$1$ de~$H^0$ au-dessus desquels la fibre de~$\pi_{H^0}$ est susceptible
d'être singulière sont donc les $\xi_t$ pour $t \in \Srond'$
et les antécédents des points génériques de~$Q$ et de~$R$.
Compte tenu de la proposition~\ref{ch3conssdi1i2}, ceci prouve l'assertion~(iii).
Par ailleurs, supposant~(ii) connue, ceci prouve aussi que les seuls points de codimension~$1$
de~$H^0$ au-dessus desquels la fibre de~$\pi_{H^0}$ est singulière et géométriquement intègre
sont les antécédents du point générique de~$R$.  Pour en déduire l'assertion~(i), il
suffit de vérifier qu'il existe un unique tel antécédent dans~$H^0$, ce qui résulte
du lemme suivant.

\bigskip
\begin{lemme}
\label{ch3conssdramifsigma}
L'unique point de la fibre de~$\rho$ au-dessus du point générique de~$R$ en lequel le morphisme~$\rho$ soit étale (resp.~ramifié)
(cf.~proposition~\ref{ch3conssdrhor}) appartient (resp.~n'appartient pas) à l'image de~$\sigma$.
\end{lemme}

\bigskip
\begin{demo}
On peut supposer~$k$ algébriquement clos.
Notons~$z=(t,L)\in Z$ le point considéré.
L'hypersurface $R \subset \Lambda$ n'étant contenue dans aucun hyperplan de~$\Lambda$,
on a $\fp{P}i \not\in L$ et par suite $t \neq t_i$ pour tout $i \in \{\uplet{1}{4}\}$.
Il en résulte que~$z$ appartient au domaine de définition de~$\tau$ (cf.~la description
explicite de~$\tau$).  Dans ces conditions, le seul antécédent possible de~$z$ par~$\sigma$
est~$\tau(z)$. (En effet, si $Z^0$ désigne le domaine de définition de~$\tau$,
le morphisme $\sigma^{-1}(Z^0) \rightarrow H$ induit par~$\tau \circ \sigma$ coïncide
avec l'inclusion $H^0 \subset H$ sur un ouvert dense de~$H^0$, donc sur~$H^0$ tout entier.)
Le point~$z$ appartient donc à l'image
de~$\sigma$ si et seulement si $\tau(z) \in H^0$, si et seulement si $\tau(z) \not\in E$.
Rappelons que~$\tau(z)$ est par définition l'unique point singulier de~$Q_t \cap L \cap H$.
Celui-ci ne peut être l'un des~$\fp{P}i$ puisque $\fp{P}i \not\in L$ pour $i \in\{\uplet{1}{4}\}$.
On a donc $\tau(z) \not\in E$ si et seulement si l'unique point singulier de~$Q_t \cap L \cap H$
n'appartient pas à~$X \cap H$, ce qui équivaut, d'après le lemme~\ref{ch3genracsimple},
à ce que~$t$ soit racine simple du polynôme~$f_{L\cap H}$, autrement dit à ce que~$\rho$ soit
étale en~$z$ (cf.~la description des fibres géométriques de~$\rho$, immédiatement avant la proposition~\ref{ch3genpfiniplat}).
\end{demo}

\bigskip
Il reste à démontrer l'assertion~(ii).  On raisonne comme dans la preuve de la proposition~\ref{ch3conssdi1i2}.
Soit $t \in \Srond'$.  Notons $L \in \Lambda$ l'image par~$\rho \circ \sigma$ de~$\xi_t$.
C'est le point générique de la variété duale $(Q_t \cap H)^\star \subset \Lambda$.
Un calcul simple montre
que $(Q_t \cap H)^\star \not\subset (Q_{t_0} \cap H)^\star$
(si $Q_{t_0} \cap H$ et $Q_t \cap H$ ont pour équations respectives
$a_1 x_1^2 + \cdots + a_4 x_4^2=0$ et $b_2x_2^2 + \cdots + b_4 x_4^2=0$,
alors $(Q_{t_0} \cap H)^\star$ est la quadrique d'équation $\lambda_1^2/a_1 + \cdots +
\lambda_4^2/a_4=0$ et $(Q_t \cap H)^\star$ est défini par le système d'équations
$\lambda_1=\lambda_2^2/b_2+\cdots+\lambda_4^2/b_4=0$).
L'hyperplan~$L$ n'est donc pas tangent à~$Q_{t_0} \cap H$,
d'où $f_{L \cap H}(t_0) \neq 0$.  En particulier, $t_0$ est racine simple du polynôme~$f_L$
(cf.~la description des racines de~$f_L$ en fonction de celles de~$f_{L \cap H}$,
après le lemme~\ref{ch3conssdxhlisse}), ce qui entraîne que la courbe~$X \cap L$ ne contient pas de droite,
même géométriquement (cf.~proposition~\ref{ch3gensousvarlin}).  Comme $X \cap L = Q_{t_0} \cap Q_t \cap L$
et que~$Q_t \cap L$ est géométriquement une réunion de deux plans permutés par
l'extension quadratique ou triviale définie par~$\epsilon_t$, le résultat voulu s'ensuit.
\end{demo}

\bigskip
Enfin, nous aurons besoin d'une propriété de régularité pour~$C_{H^0}$
et d'une propriété permettant de propager à~$C_{H^0}$ l'hypothèse que l'obstruction
de Brauer-Manin à l'existence d'un point rationnel sur~$X$ s'évanouit, lorsque~$k$ est un corps
de nombres.

\bigskip
\begin{proposition}
\label{ch3conssdregularite}
L'image par~$\pi_{H^0}$ du lieu singulier de la variété~$C_{H^0}$ est un fermé de codimension~\mbox{$\geq 2$}
dans~$H^0$.
\end{proposition}

\bigskip
\begin{demo}
Notons $S \subset H^0$ l'image par~$\pi_{H^0}$ du lieu singulier de~$C_{H^0}$;  c'est
un fermé puisque~$\pi_{H^0}$ est propre.
Nous avons déjà remarqué que la restriction de~$\sigma$ 
à $H^0 \setminus (\bigcup_{t \in \Srond'} (Q_t \cap H^0))$
est une immersion ouverte.  C'est en particulier un morphisme étale.
Le lemme~\ref{ch3conssdramifsigma} permet d'en déduire que la restriction
à $H^0 \setminus (\bigcup_{t \in \Srond'} (Q_t \cap H^0))$
du morphisme $\rho \circ \sigma \colon H^0 \rightarrow \Lambda$
est étale, d'où il résulte, compte tenu de la lissité de la variété~$C_\Lambda$ (cf.~proposition~\ref{ch3conssdreguliers}),
que $S \subset \bigcup_{t \in \Srond'} (Q_t \cap H^0)$.

Pour démontrer la proposition, on peut supposer~$k$ algébriquement clos.
D'après ce que l'on vient d'établir, si le fermé~$S$ n'était pas de codimension~$\geq 2$, il contiendrait $Q_{t_i} \cap H^0$
pour un $i \in \{\uplet{1}{4}\}$.
Par l'absurde, et quitte à renuméroter les~$t_i$, supposons
que $Q_{t_1} \cap H^0 \subset S$.
Soit alors $D \subset H$ la droite de~$H$ passant par~$\fp{P}2$ et~$\fp{P}3$.

\bigskip
\begin{lemme}
\label{ch3conssdlemmede}
On a $D \cap E = \{\fp{P}2, \fp{P}3\}$.
\end{lemme}

\bigskip
\begin{demo}
C'est évident sur les équations si l'on choisit une base dans laquelle~$q_1$ et~$q_2$
sont simultanément diagonales (cf.~proposition~\ref{ch3gensimdiag}), compte tenu
de la lissité de~$X$.
\end{demo}

\bigskip
Comme $D \cap Q_{t_1} \neq \emptyset$ et que~$Q_{t_1}$ ne contient ni~$\fp{P}2$ ni~$\fp{P}3$
(puisque ces points n'appartiennent pas à~$X$ et qu'ils appartiennent déjà respectivement
à~$Q_{t_2}$ et~$Q_{t_3}$), le lemme~\ref{ch3conssdlemmede}
montre que $D \cap (Q_{t_1} \cap H^0) \neq \emptyset$, d'où $D \cap S \neq \emptyset$.
Le lieu singulier de~$C_{H^0}$ rencontre donc $\pi_{H^0}^{-1}(D \cap H^0)$,
ce qui entraîne que $\pi_{H^0}^{-1}(D \cap H^0)$ n'est pas une variété régulière
(en effet, elle est à la fois de codimension~$2$ et définie par deux équations dans~$C_{H^0}$, qui
est de Cohen-Macaulay d'après la proposition~\ref{ch3conssdplatcm}).

Soit~$D' \subset \Lambda$ la droite formée des hyperplans $L \in \Lambda$ qui contiennent~$\fp{P}1$ et~$\fp{P}4$.

\bigskip
\begin{lemme}
\label{ch3conssdrhosigmad}
La restriction de~$\rho \circ \sigma$ à~$D \cap H^0$ est injective
et son image est incluse dans~$D'$.
\end{lemme}

\bigskip
\begin{demo}
On peut se contenter de raisonner sur les $k$\nobreakdash-points puisque le corps~$k$ est
algébriquement clos.
Soit~$h \in (D \cap H^0)(k)$. Il existe un unique \mbox{$t \in \P^1(k)$} tel que $h \in Q_t$,
et~$h$ est un point régulier de~$Q_t \cap H$.
Notons~$L$ l'hyperplan de~$\P^4_k$ contenant~$\fp{P}0$ et tangent à $Q_t \cap H$ en~$h$, de sorte que le
point de~$\Lambda$ associé à~$L$ est par définition~$\rho(\sigma(h))$.
On peut retrouver le point \mbox{$h \in (D \cap H^0)(k)$} à partir de la seule donnée de l'hyperplan~$L$
puisque c'est l'unique point d'intersection de~$L$
et de~$D$.
La restriction de~$\rho \circ \sigma$ à~$D \cap H^0$ est donc bien injective. 
Que son image soit incluse dans~$D'$ découle de la proposition~\ref{ch3genlp0}.
\end{demo}

\bigskip
\begin{lemme}
\label{ch3conssdpilambdadp}
La variété $\pi^{-1}_\Lambda(D')$ est régulière.
\end{lemme}

\bigskip
\begin{demo}
La projection naturelle $\pi^{-1}_\Lambda(D') \rightarrow X$ permet d'identifier
$\pi^{-1}_\Lambda(D')$ à la variété obtenue en faisant éclater~$X \cap K$ dans~$X$,
où~$K$ désigne le plan de~$\P^4_k$ contenant~$\fp{P}0$, $\fp{P}1$ et~$\fp{P}4$.
Comme~$X$ et~$X \cap K$ sont lisses (pour~$X \cap K$, c'est une vérification exactement
analogue à celle du lemme~\ref{ch3conssdxhlisse}), le lemme s'ensuit (cf.~\cite[Theorem~8.24]{hartshorne}).
\end{demo}

\bigskip
D'après le lemme~\ref{ch3conssdrhosigmad}, on dispose d'un morphisme injectif $D \cap H^0 \rightarrow D'$
induit par~$\rho \circ \sigma$.  Comme tout morphisme injectif entre ouverts de~$\P^1_k$, celui-ci
est une immersion ouverte. La variété $\pi_{H^0}^{-1}(D\cap H^0)$ s'identifie donc à un
ouvert de~$\pi^{-1}_\Lambda(D')$;
le lemme~\ref{ch3conssdpilambdadp} permet d'en déduire qu'elle est régulière, d'où une contradiction.
\end{demo}

\bigskip
\begin{proposition}
\label{ch3conssdadmetsect}
La restriction de la projection naturelle $C_{H^0} \rightarrow X$ (cf.~diagramme~(\ref{ch3conssddiag}))
à l'ouvert de lissité de~$C_{H^0}$ sur~$k$
admet une section rationnelle.
\end{proposition}

\bigskip
\begin{demo}
Notons $p \colon \P^4_k \setminus \{\fp{P}0\} \rightarrow H$ la projection sur~$H$ depuis~$\fp{P}0$.
La surface~$X$ n'est contenue dans aucun cône fermé de sommet~$\fp{P}0$ et de dimension~$2$,
puisque $\fp{P}0 \not\in X$.  En particulier, la trace sur~$X$ de~$p^{-1}(E)$ est un fermé
strict de~$X$, ce qui signifie que le morphisme~$p$ induit une application
rationnelle $r \colon X \dashrightarrow H^0$.

Montrons maintenant que~$r$
se factorise en $X \dashrightarrow C_{H^0} \xrightarrow{\pi_{H^0}} H^0$,
où la première flèche est une section rationnelle de la projection
$C_{H^0} \rightarrow X$. Il suffit pour cela de vérifier que pour
tout $x \in X$ tel que $p(x) \in H^0$, l'hyperplan de~$(\P^4_k)^\star$
défini par~$\rho(\sigma(p(x)))$ contient~$x$.
Il contient~$\fp{P}0$ puisque $\rho(\sigma(p(x)))$ appartient à~$\Lambda$
et il contient~$p(x)$ par définition du morphisme~$\sigma$.
Comme les points~$\fp{P}0$, $x$ et~$p(x)$ sont alignés, il contient
aussi~$x$.

Il reste à vérifier que l'image~$\xi$ du point générique de~$X$ par cette section rationnelle
appartient à l'ouvert de lissité de~$C_{H^0}$ sur~$k$. Comme $X \subset Q_{t_0} \setminus \{\fp{P}0\}$,
on a $p(X) \subset p(Q_{t_0} \setminus \{\fp{P}0\})=Q_{t_0} \cap H$ et
donc $p(X)=Q_{t_0} \cap H$ puisque le morphisme $X \rightarrow H$ induit par~$p$
est quasi-fini et que $Q_{t_0}\cap H$ est irréductible et de dimension~$2$.
Il en résulte que $\pi_{H^0}(\xi)$ est le point générique de $Q_{t_0}\cap H$;
c'est donc un point de~$H^0$ de codimension~$1$ et la proposition~\ref{ch3conssdregularite}
permet d'en déduire le résultat voulu.
\end{demo}

\subsection{Calculs explicites}
\label{ch3calcexpl}

Nous rendons ici explicite la construction du paragraphe~\ref{ch3parconssd}, dont nous conservons les hypothèses
et notations.
L'intérêt de cette entreprise est qu'elle nous permettra d'obtenir des informations
subtiles de nature non géométrique et qui semblent actuellement hors de portée d'arguments
théoriques (cf.~sous-lemmes~\ref{ch3verifdppf1d6}, \ref{ch3verifdppf1d0}, \ref{ch3verifdppf12d6}
et~\ref{ch3verifdppf12d0}). Certains des calculs ci-dessous sont inspirés
de ceux de Swinnerton-Dyer (cf.~\cite[p.~325]{bsd}); d'autres sont entièrement originaux (notamment
les propositions~\ref{ch3explicitepropp6} et~\ref{ch3explicitepsdeq}).

\subsubsection{Situation générale}

Soit $e=(\uplet{e_0}{e_4})$ la base canonique de~$k^5$.  Le point $\fp{P}0 \in \P^4_k$ est $k$\nobreakdash-rationnel
par hypothèse.  Une transformation linéaire permet de supposer qu'il a pour coordonnées
homogènes $[1:0:0:0:0]$ et que l'hyperplan~$H$ est défini par l'équation $x_0=0$,
où $[x_0:x_1:x_2:x_3:x_4]$ désignent les coordonnées homogènes de~$\P^4_k$.
Les hypothèses sur~$\fp{P}0$ et~$H$ équivalent à ce que le vecteur~$e_0$ soit orthogonal
à $\uplet{e_1}{e_4}$ pour les formes quadratiques~$q_1$ et~$q_2$ (cf.~proposition~\ref{ch3gensimdiag}), autrement dit,
à ce que la variable~$x_0$ soit simultanément diagonalisée pour~$q_1$ et~$q_2$ dans la base~$e$.

Par construction de la fibration~$\pi_{H^0}$, ses fibres s'écrivent comme des intersections de deux quadriques dans~$\P^3_k$
avec deux variables simultanément diagonalisées.  Nous allons maintenant expliciter une telle écriture.
Soient $[0:y_1:y_2:y_3:y_4]$ les coordonnées homogènes d'un point $h \in H^0(k)$.  Posons $y=\sum_{i=1}^4 y_i e_i$
et $q=q_1(y)q_2-q_2(y)q_1$.
Notons respectivement~$\phi_1$, $\phi_2$ et~$\phi$ les formes bilinéaires symétriques associées
aux formes quadratiques~$q_1$, $q_2$ et~$q$.
La quadrique projective d'équation $q(x)=0$ est
l'unique quadrique du pinceau $(Q_t)_{t \in \P^1_k}$ qui contienne~$h$. L'hyperplan tangent en~$h$ à cette quadrique
est l'espace projectif associé à l'orthogonal~$T$ de~$y$ pour~$q$.
Vu la définition de~$\pi_{H^0}$, tout ce qu'il nous reste à faire est d'exprimer une base
$(\uplet{f_0}{f_3})$ de l'hyperplan~$T$ telle que les vecteurs~$f_0$ et~$f_1$ soient orthogonaux entre
eux et à~$f_2$ et~$f_3$.

Posons $\gamma_i = \phi(e_i,y)$
et $\delta_{ij}=\phi_1(e_i,y)\phi_2(e_j,y)-\phi_1(e_j,y)\phi_2(e_i,y)$ pour $i,j\in\{\uplet{1}{4}\}$,
de sorte que
\begin{equation}
\label{ch3explicitemu}
\begin{aligned}
\phi_1(\gamma_je_i-\gamma_ie_j,y) &= \delta_{ij}q_1(y)\rlap{\text{,}} \\
\phi_2(\gamma_je_i-\gamma_ie_j,y) &= \delta_{ij}q_2(y)
\end{aligned}
\end{equation}
pour tous~$i$ et~$j$.
Quitte à permuter les vecteurs $\uplet{e_1}{e_4}$, on peut supposer que~$y_1 \neq 0$.
Quitte à permuter ensuite les vecteurs $\uplet{e_2}{e_4}$, on peut supposer
que~$\gamma_4 \neq 0$.
En effet, si l'on avait $\gamma_2=\gamma_3=\gamma_4=0$,
la relation $\sum_{i=1}^4 \gamma_i y_i = q(y) =0$ et l'hypothèse que~$y_1 \neq 0$
permettraient d'en déduire que $\gamma_i=0$ pour tout $i \in \{\uplet{1}{4}\}$,
autrement dit que $e_i \in T$ pour tout $i \in \{\uplet{0}{4}\}$,
puisque l'on a toujours $e_0 \in T$; mais
ceci est impossible puisque~$T$ est un hyperplan.
Posons maintenant
$$
\begin{aligned}
f_0&=e_0\text{,}\;\;&
f_1&=y\text{,}\;\;&
f_2&=\delta_{24}y + \gamma_2 e_4 - \gamma_4 e_2\text{,}\;\;&
f_3&=\delta_{34}y + \gamma_3 e_4 - \gamma_4 e_3\text{.}
\end{aligned}
$$
Les vecteurs~$f_1$, $f_2$ et~$f_3$ appartiennent à l'hyperplan de~$k^5$ engendré
par $\uplet{e_1}{e_4}$;
compte tenu des hypothèses sur~$e_0$ et~$\uplet{e_1}{e_4}$, cela entraîne qu'ils sont orthogonaux à~$f_0$
pour les deux formes quadratiques~$q_1$ et~$q_2$.
Par ailleurs, les vecteurs~$f_2$ et~$f_3$ sont orthogonaux à~$f_1$ pour~$q_1$ et~$q_2$
d'après les relations~(\ref{ch3explicitemu}).
Comme~$f_1=y$, il en résulte notamment que $f_i \in T$ pour tout $i \in \{\uplet{0}{3}\}$.
Enfin, l'hypothèse $y_1 \gamma_4 \neq 0$ assure que la famille $f=(\uplet{f_0}{f_3})$ est libre;
c'est donc une base de~$T$.
Ainsi la famille~$f$ satisfait-elle bien à toutes les conditions voulues.
Les matrices dans~$f$ des restrictions à~$T$ des formes quadratiques~$q_1$ et~$q_2$
sont donc
\begin{equation}
\label{ch3explicitematricesres}
\begin{pmatrix}
\alpha_0 & 0 & 0 & 0 \\
0 & \alpha_1 & 0 & 0 \\
0 & 0 & \alpha_2 & \alpha_4 \\
0 & 0 & \alpha_4 & \alpha_3
\end{pmatrix}
\;\qquad\text{ et }\;\qquad
\begin{pmatrix}
\beta_0 & 0 & 0 & 0 \\
0 & \beta_1 & 0 & 0 \\
0 & 0 & \beta_2 & \beta_4 \\
0 & 0 & \beta_4 & \beta_3
\end{pmatrix}\!\rlap{\text{,}}
\end{equation}
où $\alpha_i=q_1(f_i)$ et $\beta_i=q_2(f_i)$ pour $i\in\{0,1,2,3\}$,
$\alpha_4=\phi_1(f_2,f_3)$ et $\beta_4=\phi_2(f_2,f_3)$.
On peut exprimer les~$\alpha_i$ comme suit:
\begin{equation}
\label{ch3explicitealphaigen}
\begin{aligned}
\alpha_0 &= q_1(e_0)\rlap{\text{,}}\\
\alpha_1 &= q_1(y)\rlap{\text{,}}\\
\alpha_2 &= q_1(\gamma_2 e_4 - \gamma_4 e_2) - \delta_{24}^2 q_1(y)\rlap{\text{,}}\\
\alpha_3 &= q_1(\gamma_3 e_4 - \gamma_4 e_3) - \delta_{34}^2 q_1(y)\rlap{\text{,}}\\
\alpha_4 &= \phi_1(\gamma_2 e_4 - \gamma_4 e_2, \gamma_3 e_4 - \gamma_4 e_3) - \delta_{24}\delta_{34}q_1(y)\rlap{\text{;}}
\end{aligned}
\end{equation}
on obtient des expressions similaires pour les~$\beta_i$ en remplaçant~$q_1$ par~$q_2$ et~$\phi_1$ par~$\phi_2$.
Les deux premières formules ci-dessus sont mises pour mémoire.  Les trois suivantes se déduisent
tout de suite des relations~(\ref{ch3explicitemu}) et de la définition des~$\alpha_i$.

La fibre de~$\pi_{H^0}$ en~$h$
s'identifie à l'intersection des deux quadriques de~$\P(T)$ définies par les matrices symétriques~(\ref{ch3explicitematricesres}).
Celle-ci est lisse si et seulement si
\begin{equation}
\label{ch3explicitelisse}
\dgras_{01}(\dgras_{04}^2-\dgras_{02}\dgras_{03})(\dgras_{14}^2-\dgras_{12}\dgras_{13})(\dgras_{23}^2+4\dgras_{24}\dgras_{34})
\neq 0\rlap{\text{,}}
\end{equation}
où l'on a posé $\dgras_{ij}=\alpha_i \beta_j - \alpha_j \beta_i$ pour $i,j\in\{\uplet{0}{4}\}$
(cf.~proposition~\ref{ch3genpropeq}).
Cette équation fournit une description explicite du lieu des fibres singulières de~$\pi_{H^0}$
au-dessus de l'ouvert $y_1 \gamma_4 \neq 0$ de~$H^0$.  Comparons-la maintenant avec la conclusion
du théorème~\ref{ch3conssdtheoreme}.

\bigskip
\begin{proposition}
\label{ch3expliciteabstth}
On a $\dgras_{01}=0$ si et seulement si $h \in Q_{t_0}$,
$\dgras_{01}(\dgras_{04}^2-\dgras_{02}\dgras_{03})=0$ si et seulement si $\rho(\sigma(h))\in Q$,
$\dgras_{14}^2-\dgras_{12}\dgras_{13}=0$ si et seulement si $h \in \bigcup_{t\in\Srond'}Q_t$
et enfin $(\dgras_{14}^2-\dgras_{12}\dgras_{13})(\dgras_{23}^2+4\dgras_{24}\dgras_{34})=0$ si et seulement si $\rho(\sigma(h))\in R$.
\end{proposition}

\bigskip
\begin{demo}
Posons $q_0=\alpha_0 q_2 - \beta_0 q_1$.  Comme $\alpha_0=q_1(e_0)$ et $\beta_0=q_2(e_0)$,
la quadrique projective d'équation $q_0(x)=0$ n'est autre que~$Q_{t_0}$, d'où la première assertion.
Pour la seconde, remarquons que $\rho(\sigma(h)) \in Q$ si et seulement si la restriction de~$q_0$
à $T\cap \langle \uplet{e_1}{e_4}\rangle$ est dégénérée, vu la définition de~$Q$ et le lemme~\ref{ch3conssdxhlisse}.
Compte tenu que $T \cap \langle \uplet{e_1}{e_4} \rangle = \langle f_1, f_2, f_3 \rangle$ et que la
matrice dans $(f_1,f_2,f_3)$ de la restriction de~$q_0$ à ce sous-espace est
$$
\begin{pmatrix}
\dgras_{01} & 0 & 0 \\
0 & \dgras_{02} & \dgras_{04} \\
0 & \dgras_{04} & \dgras_{03}
\end{pmatrix}\!\rlap{\text{,}}
$$
il s'ensuit que
$\dgras_{01}(\dgras_{04}^2-\dgras_{02}\dgras_{03})=0$ si et seulement si $\rho(\sigma(h))\in Q$.

Comme l'hyperplan $\P(T) \cap H \subset H$ est tangent en un point lisse (à savoir~$h$) à la quadrique
projective d'équation $q(x)=0$ dans~$H$, le rang de la restriction de~$q$ à $T \cap \langle \uplet{e_1}{e_4} \rangle$ est égal à~$r-2$,
où~$r$ désigne le rang de la restriction de~$q$ à $\langle \uplet{e_1}{e_4} \rangle$.
La matrice de la restriction de~$q$ à $T \cap \langle \uplet{e_1}{e_4}\rangle$ dans la base~$(f_1,f_2,f_3)$ étant égale à
$$
\begin{pmatrix}
0 & 0 & 0 \\
0 & \dgras_{12} & \dgras_{14} \\
0 & \dgras_{14} & \dgras_{13}
\end{pmatrix}\!\rlap{\text{,}}
$$
on en déduit que $r<4$ si et seulement si $\dgras_{14}^2-\dgras_{12}\dgras_{13}=0$.
La condition $r<4$ est par ailleurs équivalente à $h \in \bigcup_{t \in \Srond'} Q_t$,
d'où la conclusion recherchée.

Enfin, on a $\rho(\sigma(h)) \in R$ si et seulement si $\P(T) \cap H \cap X$ n'est pas lisse
(par définition de~$R$), si et seulement si le discriminant du polynôme~$f_{\P(T) \cap H}$ est nul
(cf.~proposition~\ref{ch3genpropeq}).
La dernière assertion de la proposition s'en déduit en calculant ce polynôme,
qui est un déterminant, dans la
base $(f_1,f_2,f_3)$ de $T \cap \langle \uplet{e_1}{e_4}\rangle$.
\end{demo}

\bigskip
Rappelons maintenant comment s'exprime la jacobienne de la fibre de~$\pi_{H^0}$ en~$h$, lorsque
celle-ci est lisse.  Posons
$c=4\dgras_{04}\dgras_{14}-2\dgras_{02}\dgras_{13}-2\dgras_{03}\dgras_{12}$
et $d=4\dgras_{01}^2(\dgras_{23}^2+4\dgras_{24}\dgras_{34})$\glossary{$c$, $d$, $c^2-d$}.
Il est utile de remarquer que
$c^2-d=16(\dgras_{04}^2-\dgras_{02}\dgras_{03})(\dgras_{14}^2-\dgras_{12}\dgras_{13})$.

\bigskip
\begin{proposition}
\label{ch3explicitegeneralites}
Supposons la fibre de~$\pi_{H^0}$ en~$h$ lisse et notons~$E'$ sa jacobienne.
La courbe elliptique~$E'$ a pour équation de Weierstrass
$Y^2=(X-c)(X^2-d)$.  De plus, si $\phi' \colon E' \rightarrow E''$ désigne le quotient
de~$E'$ par le point de coordonnées $(X,Y)=(c,0)$ et $\phi'' \colon E'' \rightarrow E'$ l'isogénie
duale de~$\phi'$,
la courbe elliptique~$E''$ a pour équation
de Weierstrass $Y^2=(X+2c)(X^2-4(c^2-d))$
et la fibre de~$\pi_{H^0}$ en~$h$ est canoniquement un $2$\nobreakdash-revêtement de~$E'$, déterminant
une classe de $H^1(k,\tors{2}E')$ dont l'image dans $H^1(k,\tors{\phi''}E'')=k^\star/k^{\star 2}$
par la flèche induite par $\phi' \colon \tors{2}E' \rightarrow \tors{\phi''}E''$ est égale à la classe de $\dgras_{14}^2
-\dgras_{12}\dgras_{13}$
dans $k^\star/k^{\star 2}$.
\end{proposition}

\bigskip
\begin{demo}
Ces généralités sont décrites en détail dans~\cite[§3]{bsd}.
\end{demo}

\subsubsection{Cas simultanément diagonal}

Supposons les formes quadratiques~$q_1$ et~$q_2$ simultanément diagonales dans la base~$e$.
On pose alors $a_i=q_1(e_i)$, $b_i=q_2(e_i)$ et $d_{ij}=a_i b_j - a_j b_i$
pour $i,j \in \{\uplet{0}{4}\}$, conformément aux notations introduites dans l'exemple~\ref{ch3exemplesimdiag1}.
Définissons des polynômes $\uplet{p_0}{p_6}\in k[\uplet{y_1}{y_4}]$ par les formules suivantes
(pour $i, j\in\{1,2,3,4\}$ distincts fixés, nous convenons de noter~$k$ et~$\ell$
les entiers uniquement déterminés par les conditions $1 \leq k < \ell \leq 4$ et $\{i,j,k,\ell\}=\{1,2,3,4\}$;
de même pour~$i'$, $j'$, $k'$ et~$\ell'$):\glossary{$\uplet{p_0}{p_6}$}
\begin{align*}
p_i&=\sum_{j=1}^4 d_{ij} y_j^2 \quad\text{ pour } i \in \{0,1,2,3,4\}\rlap{\text{,}}\\
p_5&=\sum_{1 \leq i < j \leq 4} d_{ij}^2 d_{0k} d_{0\ell} y_i^2 y_j^2\rlap{\text{,}}\\
p_6&=\sum_{\genfrac{}{}{0pt}{}{\scriptstyle 1 \leq i < j \leq 4}{\scriptstyle 1 \leq i' < j' \leq 4}}
d_{ij}^2 d_{i'j'}^2
\left(
d_{k'k}d_{\ell\ell'} + d_{\ell'k}d_{\ell k'}
\right)
y_i^2 y_j^2 y_{i'}^2 y_{j'}^2\rlap{\text{.}}
\end{align*}
Les équations~(\ref{ch3explicitealphaigen}) se simplifient et fournissent les expressions
suivantes pour les~$\alpha_i$ et les~$\beta_i$:
\begin{equation}
\label{ch3explicitealphaisim}
\begin{aligned}
\alpha_0 &= a_0\rlap{\text{,}}&
\beta_0 &= b_0\rlap{\text{,}}\\
\alpha_1 &= q_1(y)\rlap{\text{,}}&
\beta_1 &= q_2(y)\rlap{\text{,}}\\
\alpha_2 &= a_4 \gamma_2^2 + a_2 \gamma_4^2 - \delta_{24}^2 q_1(y)\rlap{\text{,}}&
\beta_2 &= b_4 \gamma_2^2 + b_2 \gamma_4^2 - \delta_{24}^2 q_2(y)\rlap{\text{,}}\\
\alpha_3 &= a_4 \gamma_3^2 + a_3 \gamma_4^2 - \delta_{34}^2 q_1(y)\rlap{\text{,}}&
\beta_3 &= b_4 \gamma_3^2 + b_3 \gamma_4^2 - \delta_{34}^2 q_2(y)\rlap{\text{,}}\\
\alpha_4 &= a_4 \gamma_2\gamma_3 - \delta_{24}\delta_{34}q_1(y)\rlap{\text{,}}&
\beta_4 &= b_4 \gamma_2\gamma_3 - \delta_{24}\delta_{34}q_2(y)\rlap{\text{.}}
\end{aligned}
\end{equation}
On a de plus $\gamma_i=-y_ip_i$ et $\delta_{ij}=d_{ij}y_iy_j$ pour tous $i,j \in \{\uplet{1}{4}\}$.

\bigskip
\begin{proposition}
\label{ch3expliciteegalitessim}
Les égalités suivantes ont lieu:
\vspace{-9pt}\begin{align*}
\dgras_{01}&=p_0\rlap{\text{,}} &
\dgras_{14}^2-\dgras_{12}\dgras_{13} &= y_1^2\gamma_4^2 \prod_{i=1}^4 p_i\rlap{\text{,}}\\
\dgras_{04}^2-\dgras_{02}\dgras_{03} &= -y_1^2 \gamma_4^2 p_5\rlap{\text{,}}&
\dgras_{23}^2+4\dgras_{24}\dgras_{34} &= y_1^4 \gamma_4^4 p_6\rlap{\text{.}}
\end{align*}
\end{proposition}

\vspace{-9pt}
%\bigskip
\begin{demo}
Ces égalités sont purement formelles, au sens où elles valent dans l'anneau des polynômes à coefficients entiers
en les~$a_i$, les~$b_i$ et les~$y_i$.  Quelques explications semblent néanmoins appropriées (excepté pour l'égalité
$\dgras_{01}=p_0$,
qui est une trivialité).
Introduisons, à la suite de Swinnerton-Dyer, les invariants fondamentaux
\begin{align*}
\Theta_{\alpha\alpha}&=\alpha_4^2-\alpha_2\alpha_3\rlap{\text{,}}&
\Theta_{\alpha\beta}&=2\alpha_4\beta_4-\alpha_2\beta_3-\alpha_3\beta_2\rlap{\text{,}}&
\Theta_{\beta\beta}&=\beta_4^2-\beta_2\beta_3
\end{align*}
du couple de formes quadratiques considéré.
Si l'on admet les égalités
\begin{align}
\label{ch3expliciteinvfond1}
\Theta_{\alpha\alpha} &= -y_1^2 \gamma_4^2 \sum_{1 \leq i < j \leq 4} d_{ij}^2 a_k a_\ell y_i^2 y_j^2\rlap{\text{,}} \\
\label{ch3expliciteinvfond2}
\Theta_{\alpha\beta} &= -y_1^2 \gamma_4^2 \sum_{1 \leq i < j \leq 4} d_{ij}^2 \left( a_k
b_\ell + a_\ell b_k\right) y_i^2 y_j^2\rlap{\text{,}} \\
\label{ch3expliciteinvfond3}
\Theta_{\beta\beta} &= -y_1^2 \gamma_4^2 \sum_{1 \leq i < j \leq 4} d_{ij}^2
b_k b_\ell y_i^2 y_j^2\rlap{\text{,}}
\end{align}
(cf.~\cite[p.~326]{bsd}), les formules
\begin{align*}
\dgras_{04}^2-\dgras_{02}\dgras_{03}&=\alpha_0^2 \Theta_{\beta\beta} - \alpha_0\beta_0\Theta_{\alpha\beta}
+ \beta_0^2\Theta_{\alpha\alpha}\rlap{\text{,}}\\
\dgras_{14}^2-\dgras_{12}\dgras_{13}&=\alpha_1^2\Theta_{\beta\beta}-\alpha_1\beta_1\Theta_{\alpha\beta}
+\beta_1^2\Theta_{\alpha\alpha}\rlap{\text{,}}\\
\dgras_{23}^2+4\dgras_{24}\dgras_{34}&=\Theta_{\alpha\beta}^2-4\Theta_{\alpha\alpha}\Theta_{\beta\beta}
\end{align*}
(cf.~\cite[p.~324]{bsd}) permettent de conclure.  Il reste donc à établir~(\ref{ch3expliciteinvfond1})
et~(\ref{ch3expliciteinvfond2}), l'égalité~(\ref{ch3expliciteinvfond3}) étant symétrique
de~(\ref{ch3expliciteinvfond1}).  Nous allons détailler la preuve de~(\ref{ch3expliciteinvfond1});
celle de~(\ref{ch3expliciteinvfond2}) est similaire. À l'aide de~(\ref{ch3explicitealphaisim}), on
peut écrire que
\begin{align*}
\Theta_{\alpha\alpha} &=
(a_4\gamma_2\gamma_3-\delta_{24}\delta_{34}q_1(y))^2 - (a_4\gamma_2^2+a_2\gamma_4^2-\delta_{24}^2q_1(y))
(a_4\gamma_3^2+a_3\gamma_4^2
\\
&\hspace*{9.2cm}-\delta_{34}^2q_1(y)) \\
&=q_1(y)\left((a_2\delta_{34}^2  + a_3\delta_{24}^2)\gamma_4^2 + a_4(\delta_{34}\gamma_2-\delta_{24}\gamma_3)^2\right)
\\&\hspace*{5.83cm}- \gamma_4^2(a_2 a_3\gamma_4^2  + a_3 a_4\gamma_2^2  + a_2 a_4\gamma_3^2) \rlap{\text{.}}
\end{align*}
L'identité $\delta_{34}\gamma_2+\delta_{42}\gamma_3+\delta_{23}\gamma_4=0$ permet de simplifier l'expression ci-dessus;
on trouve ainsi que
\begin{align*}
-\Theta_{\alpha\alpha}/\gamma_4^2 &=
a_2 a_3\gamma_4^2  + a_3 a_4\gamma_2^2  + a_2 a_4\gamma_3^2-
q_1(y)(a_2 \delta_{34}^2 + a_3\delta_{24}^2 + a_4\delta_{23}^2)\rlap{\text{.}}
\end{align*}
Pour vérifier que la formule que l'on vient d'obtenir pour $-\Theta_{\alpha\alpha}/\gamma_4^2$ coïncide avec
celle à laquelle on veut aboutir, il est utile de remarquer qu'elles sont toutes deux invariantes par permutation des
indices~$\{2,3,4\}$ et homogènes de degré~$3$ en les~$y_i^2$.
Il suffit donc de calculer les coefficients,
dans $-\Theta_{\alpha\alpha}/\gamma_4^2$, des monômes suivants : $y_2^2 y_3^2 y_4^2$,
$y_2^4 y_3^2, y_2^6$, $y_1^2 y_2^2 y_3^2$, $y_1^2 y_2^4$, $y_1^4 y_2^2$, $y_1^6$.
Il est immédiat que les coefficients de~$y_2^4 y_3^2$, $y_2^6$, $y_1^2y_2^4$ et~$y_1^6$ sont nuls
et que celui de~$y_1^4y_2^2$ est bien~$d_{12}^2a_3a_4$.  Que le coefficient
de~$y_1^2 y_2^2 y_3^2$ soit égal à~$d_{23}^2a_1a_4$ et que celui de~$y_2^2 y_3^2 y_4^2$ soit nul
résulte des identités $a_1 d_{32}+a_2 d_{13}+a_3 d_{21}=0$
et $a_2^2 d_{34}^2 + a_3^2 d_{24}^2 + a_4^2 d_{23}^2 = 2(a_2 a_3 d_{42} d_{43} + a_3 a_4 d_{23} d_{24} + a_2 a_4 d_{32} d_{34})$.
\end{demo}

\bigskip
Notons $\Delta \subset H^0$ l'intersection de~$H^0$ et de l'hypersurface de~$H$ définie par l'équation
$$
\prod_{i=0}^6 p_i = 0 \rlap{\text{.}}
$$

\bigskip
\begin{proposition}
\label{ch3explicitedelta}
L'ensemble des $h \in H^0$ au-dessus desquels la fibre de~$\pi_{H^0}$ est singulière
est égal à l'ensemble sous-jacent à~$\Delta$.
\end{proposition}

\bigskip
\begin{demo}
Soit $h \in H^0$, de coordonnées homogènes $[0:y_1:\cdots:y_4]$.
Compte tenu de l'invariance des polynômes $\uplet{p_0}{p_6}$ par permutation
des indices~$1$ à~$4$, les commentaires qui suivent les équations~(\ref{ch3explicitemu})
permettent de supposer que $y_1 \gamma_4 \neq 0$ pour établir que la fibre de~$\pi_{H^0}$
en~$h$ est lisse si et seulement si $h \not\in \Delta$.
Sous cette hypothèse, on a vu que la fibre de~$\pi_{H^0}$ en~$h$ est lisse
si et seulement si la relation~(\ref{ch3explicitelisse}) est satisfaite.
Celle-ci équivaut bien à ce que $h \not\in \Delta$, d'après la
proposition~\ref{ch3expliciteegalitessim}.
\end{demo}

\bigskip
\begin{proposition}
\label{ch3expliciteirred}
Les polynômes $p_i \in k[\uplet{y_1}{y_4}]$ pour $i \in \{\uplet{0}{6}\}$ sont irréductibles et premiers entre
eux deux à deux.
\end{proposition}

\bigskip
\begin{demo}
La lissité de~$X$ entraîne que $d_{ij}\neq 0$ pour $i,j$ distincts.
Les polynômes $\uplet{p_0}{p_4}$ sont donc des formes quadratiques de rang~$\geq 3$,
ce qui assure leur irréductibilité.
Il résulte des propositions~\ref{ch3expliciteabstth} et~\ref{ch3expliciteegalitessim}
que l'image par~$\rho\circ\sigma$ de chaque point générique de l'hypersurface
de~$H$ définie par l'équation $p_0p_5=0$ appartient à~$Q$.  Les fibres de~$\pi_{H^0}$
au-dessus de ces points ne sont donc pas géométriquement intègres (cf.~proposition~\ref{ch3conssdi1i2}
et \cite[12.2.4~(viii)]{ega43}).  Comme par ailleurs l'hypersurface d'équation
$p_1p_2p_3p_4=0$ coïncide avec $\bigcup_{t \in \Srond'} Q_t \cap H$, les deux premières assertions du théorème~\ref{ch3conssdtheoreme}
permettent d'en déduire que le polynôme~$p_6$ est irréductible.

L'irréductibilité de~$p_5$ est prouvée dans \cite[p.~326]{bsd}; pour la commodité du lecteur, nous reproduisons l'argument ici.
Considérons l'action du groupe $G=(\Z/2)^4$ sur $k[\uplet{y_1}{y_4}]$ donnée
par $(\uplet{\epsilon_1}{\epsilon_4})\cdot y_i=(-1)^{\epsilon_i} y_i$.
Étant donné que~$p_5$ est invariant sous cette action, le groupe~$G$ agit naturellement sur l'ensemble~$F$ des facteurs irréductibles de~$p_5$
à multiplication par une unité près.
Si cette action n'était pas transitive,
il existerait $f,g \in k[\uplet{y_1}{y_4}]$
non constants,
invariants sous~$G$ à multiplication par une unité près, tels que \mbox{$p_5=fg$}.
Les polynômes~$f$ et~$g$ seraient nécessairement invariants sous~$G$
puisque~$p_5$ n'est divisible par~$y_i$ pour aucun $i \in \{\uplet{1}{4}\}$;
ils appartiendraient donc au sous-anneau $k[\uplet{y_1^2}{y_4^2}]$.
Quitte à échanger~$f$ et~$g$, on pourrait supposer qu'il
existe $i,j \in \{\uplet{1}{4}\}$ distincts tels que~$f$ soit de degré non nul en~$y_i$ et en~$y_j$.
Comme~$p_5$ est de degré total~$4$ et de degré~$2$ en~$y_i$ et en~$y_j$, que~$g$ est non constant et que les degrés de~$f$ en~$y_i$ et en~$y_j$ sont pairs
(puisque $f \in k[\uplet{y_1^2}{y_4^2}]$) et non nuls, on conclurait que le coefficient de $y_i^2 y_j^2$ dans~$p_5$ est nul, ce qu'il n'est pas.
L'action de~$G$ sur~$F$ est donc transitive.
Il en résulte que~$p_5$ est irréductible, compte tenu que le coefficient dominant de~$p_5$ vu comme polynôme en~$y_1$ à coefficients dans $k[y_2,y_3,y_4]$
est lui-même un polynôme irréductible, étant une forme quadratique de rang~$3$.

Enfin, pour vérifier que les polynômes $\uplet{p_0}{p_6}$ sont premiers entre eux deux à deux, sachant qu'ils sont irréductibles,
il suffit de remarquer que
les~$p_i$ pour $i \in \{\uplet{0}{4}\}$ sont deux à deux non proportionnels
et que \mbox{$\deg(p_i)<\deg(p_5)<\deg(p_6)$} pour tout $i\in\{\uplet{0}{4}\}$.
\end{demo}

\bigskip
\begin{remarque}
Il est possible d'établir l'irréductibilité de~$p_6$ en travaillant directement
sur sa définition.
Inversement, une variante du théorème~\ref{ch3propmonod} permettrait sans doute
d'obtenir l'irréductibilité de~$p_5$ par voie abstraite.
\end{remarque}

\bigskip
Voici enfin quelques identités dont nous aurons besoin par la suite.  Notons $c_6(y_i^m)$
le coefficient de~$y_i^m$ dans~$p_6$, où l'on considère~$p_6$ comme un polynôme en~$y_i$
à coefficients dans $k[(y_\ell)_{\ell\neq i}]$.

\bigskip
\begin{proposition}
\label{ch3explicitepropp6}
Les égalités suivantes ont lieu:
\begin{gather}
\label{ch3expliciteeqp1p6}
\begin{split}
16 p_1^3 y_2^2 y_3^2 y_4^2 d_{34}^2 d_{23}^2 d_{24}^2 d_{21} d_{31} d_{41} =
(2 c_6(y_1^4) y_1^2 + c_6(y_1^2))^2 - 4c_6(y_1^4)p_6\rlap{\text{,}}
\end{split} \\
\label{ch3expliciteeqp2p6}
16 p_2^3 y_1^2 y_3^2 y_4^2 d_{34}^2 d_{13}^2 d_{14}^2 d_{12} d_{32} d_{42} =
(2 c_6(y_2^4) y_2^2 + c_6(y_2^2))^2 - 4c_6(y_2^4)p_6\rlap{\text{,}} \\
\label{ch3expliciteeqp1p2p0}
4d_{01}d_{02}p_1p_2 = (d_{01}p_2 + d_{02}p_1)^2 - d_{12}^2 p_0^2 \rlap{\text{.}}
\end{gather}
\end{proposition}

\vspace{-9pt}
%\bigskip
\begin{demo}
La troisième égalité n'est autre que l'identité $d_{01}p_2 - d_{02}p_1=-d_{12}p_0$ élevée au carré.
Les deux premières étant symétriques,
il reste seulement à établir~(\ref{ch3expliciteeqp1p6}).  C'est à nouveau une égalité formelle, au sens où elle vaut
dans l'anneau des polynômes à coefficients entiers en les~$a_i$, $b_i$, $y_i$;
tout logiciel de calcul symbolique permet de la vérifier.
Alternativement, il est possible de la démontrer à la main en moins de trois pages,
comme me le prouva par l'exemple Swinnerton-Dyer peu de temps après
que je lui eus envoyé une version préliminaire du manuscrit
(lettre de Swinnerton-Dyer à l'auteur, 21~octobre~2005).
\end{demo}

\subsubsection{Cas presque simultanément diagonal}

Nous supposons ici que les formes quadratiques~$q_1$ et~$q_2$ sont «~presque simultanément diagonales~» dans la
base~$e$, c'est-à-dire que pour ces deux formes, les vecteurs $e_0$, $e_1$ et~$e_2$ sont deux à deux orthogonaux et orthogonaux
à~$e_3$ et~$e_4$. Soient $a_i=q_1(e_i)$ et $b_i=q_2(e_i)$ pour $i\in\{\uplet{0}{4}\}$, $a_5=\phi_1(e_3,e_4)$,
$b_5=\phi_2(e_3,e_4)$ et $d_{ij}=a_ib_j-a_jb_i$ pour $i,j\in\{\uplet{0}{5}\}$.
Les équations~(\ref{ch3explicitealphaigen}) deviennent ici:
\begin{equation}
\label{ch3explicitealphaipsim}
\begin{aligned}
\alpha_0 &= a_0\rlap{\text{,}}\\
\alpha_1 &= q_1(y)\rlap{\text{,}}\\
\alpha_2 &= a_4 \gamma_2^2 + a_2 \gamma_4^2 - \delta_{24}^2 q_1(y)\rlap{\text{,}}\\
\alpha_3 &= a_4 \gamma_3^2 + a_3 \gamma_4^2 - 2a_5\gamma_3\gamma_4 - \delta_{34}^2 q_1(y)\rlap{\text{,}}\\
\alpha_4 &= a_4 \gamma_2\gamma_3 - a_5\gamma_2\gamma_4 - \delta_{24}\delta_{34}q_1(y)\rlap{\text{,}}
\end{aligned}
\end{equation}
et de même en remplaçant~$\alpha_i$ et~$a_i$ respectivement par~$\beta_i$ et~$b_i$ pour tout~$i$, et~$q_1$ par~$q_2$.
Posons $p_i=a_i q_2 - b_i q_1$ pour $i\in\{0,1\}$.

\bigskip
\begin{proposition}
\label{ch3explicitepsdeq}
Le polynôme~$p_0$ divise $d_{01}^2 c + 4 d_{02}(d_{03}d_{04}-d_{05}^2)y_1^2\gamma_4^2p_1^3$ dans $k[\uplet{y_1}{y_4}]$.
\end{proposition}

\bigskip
Remarquons que $\epsilon_0 \in k^\star/k^{\star 2}$ est représenté
par $d_{01}d_{02}(d_{03}d_{04}-d_{05}^2)$.

\bigskip
\begin{demo}
Une fois remplacés dans la définition de~$c$ les~$\alpha_i$ et les~$\beta_i$ par les membres
de droite des équations~(\ref{ch3explicitealphaipsim}) (et des équations analogues pour les~$\beta_i$) et développée l'expression
obtenue,
l'identité $\delta_{34}\gamma_2-\delta_{24}\gamma_3=-\delta_{23}\gamma_4$ permet
de voir que~$c$ est divisible par~$\gamma_4^2$, et plus précisément que
$$
\frac{c}{2\gamma_4^2}=
-2\gamma_3\gamma_4(d_{02}f_5+d_{05}f_2)
-2d_{05}\gamma_2^2 f_5 +
 \sum_{\{i,j,k\}=\{2,3,4\}} d_{0k} \gamma_i^2 f_j
\mod p_0
\rlap{\text{,}}
$$
où $f_i=a_i q_2 - b_i q_1$ pour $i\in\{\uplet{2}{5}\}$.
Le résultat voulu s'en déduit aisément, compte tenu des relations suivantes:
\begin{eqnarray*}
d_{01}f_i &=& d_{0i}p_1\!\mod p_0 \quad \text{ pour tout } i\in\{\uplet{2}{5}\}\rlap{\text{,}}\\
d_{01}\gamma_2 &=& -d_{02}y_2 p_1\!\mod p_0\rlap{\text{,}}\\
d_{01}\gamma_3 &=& -(d_{03}y_3 + d_{05}y_4)p_1\!\mod p_0\rlap{\text{,}}\\
d_{01}\gamma_4 &=& -(d_{05}y_3 + d_{04}y_4)p_1\!\mod p_0\rlap{\text{,}}\\
d_{01}y_1^2 &=& -d_{02}y_2^2 - d_{03}y_3^2 - d_{04}y_4^2 - 2d_{05}y_3y_4\!\mod p_0\rlap{\text{.}}
\end{eqnarray*}
\end{demo}

\vspace{-10pt}

\subsection{Spécialisation de la \condD{}}
\label{ch3parspecd}

À toute droite de~$H$ suffisamment générale est associée un pinceau de courbes
de genre~$1$: la restriction de~$\pi_{H^0}$ au-dessus de cette droite.  Nous
avons vu que la fibration jacobienne de
ce pinceau admet une section d'ordre~$2$ (cf.~proposition~\ref{ch3conssd2rat}).  Cette fibration étant par ailleurs
semi-stable, comme nous le prouverons ci-dessous, nous pouvons lui appliquer
les résultats du chapitre~2, et notamment le théorème~\ref{ch2thsurmesure}, dans
l'énoncé duquel une «~\condE{}~» a été définie.
L'objet de ce paragraphe est d'exhiber beaucoup
de droites de~$H$ pour lesquelles le pinceau associé satisfait à la \condE{}, lorsque~$k$
est un corps de nombres.
Il sera pour cela nécessaire de supposer qu'une certaine forme générique
de la~\condD{} est vérifiée.   Occupons-nous donc d'abord de définir celle-ci précisément.

Dans tout ce paragraphe, on suppose donnés un point $h_0 \in H^0(k)$
au-dessus duquel la fibre de~$\pi_{H^0}$ est lisse
et un hyperplan $k$\nobreakdash-rationnel $\Pi \subset H$ ne contenant pas~$h_0$\glossary{$h_0$, $\Pi$, $E'$}.

Soit~$E'$ la jacobienne de la fibre générique de~$\pi_{H^0}$.

\bigskip
\begin{proposition}
\label{ch3specdss}
La courbe elliptique~$E'$ est à réduction semi-stable en tout point de codimension~$1$
de~$H^0$.  Elle possède un unique point rationnel d'ordre~$2$, et si l'on note~$E''$
le quotient de~$E'$ par ce point, la courbe elliptique~$E''$ possède elle aussi un unique point rationnel d'ordre~$2$.
\end{proposition}

\bigskip
\begin{demo}
Il est déjà connu que~$E'$ possède au moins un point rationnel d'ordre~$2$ (cf.~proposition~\ref{ch3conssd2rat}).
Toute courbe elliptique $2$\nobreakdash-isogène à~$E'$ possède donc un point rationnel d'ordre~$2$.
Ainsi, pour prouver la proposition, on peut supposer~$k$ algébriquement clos, auquel
cas les formes quadratiques~$q_1$ et~$q_2$ sont simultanément diagonalisables
(cf.~proposition~\ref{ch3gensimdiag}).
Soit $h \in H^0$ un point de codimension~$1$.
Reprenons les notations du paragraphe~\ref{ch3calcexpl};
une transformation linéaire dans~$\P^4_k$ permet de supposer que~$q_1$ et~$q_2$ sont simultanément diagonales
dans la base~$e$, que le point~$\fp{P}0$ a pour coordonnées homogènes $[1:0:0:0:0]$
et que $(y_1\gamma_4)(h)\neq 0$, comme expliqué immédiatement après les équations~(\ref{ch3explicitemu}).
Notons $r \colon k[\uplet{y_1}{y_4}] \rightarrow \Orond_{H^0,h}$ l'application
qui à $f(\uplet{y_1}{y_4})$ associe $f(1,y_2/y_1,y_3/y_1,y_4/y_1)$; c'est un morphisme
d'anneaux.  Nous avons vu que la courbe elliptique~$E'$
a pour équation de Weierstrass
$$
Y^2=(X-r(c))(X^2-r(d))
$$
(cf.~proposition~\ref{ch3explicitegeneralites}, appliquée au point générique
de~$H^0$; en toute rigueur, dans le paragraphe~\ref{ch3calcexpl}, le point~$h$ était
supposé $k$\nobreakdash-rationnel, mais rien n'empêche d'étendre les scalaires de~$k$ à~$\kappa(H)$
avant d'appliquer cette proposition).
Notons~$v$ la valuation normalisée associée à l'anneau de valuation
discrète~$\Orond_{H^0,h}$.
Compte tenu que~$v(r(y_1))=v(r(\gamma_4))=0$,
les propositions~\ref{ch3expliciteegalitessim} et~\ref{ch3expliciteirred} montrent
que l'un de~$v(r(d))$ et de~$v(r(c^2-d))$ est nul.
Il en résulte d'une part que l'équation de Weierstrass ci-dessus est minimale
(en tant qu'équation de Weierstrass à coefficients dans~$\Orond_{H^0,h}$)
et d'autre part, comme son discriminant est~$16r(d)r(c^2-d)^2$,
que la courbe elliptique~$E'$ est à réduction semi-stable en~$h$.
Les propositions~\ref{ch3expliciteegalitessim} et~\ref{ch3expliciteirred} montrent de
plus que ni~$r(d)$ ni~$r(c^2-d)$ ne sont des carrés dans~$\kappa(H)$,
ce qui prouve la seconde assertion,
étant donné que le quotient de~$E'$ par le
point de coordonnées $(X,Y)=(r(c),0)$ a pour équation de
Weierstrass $Y^2=(X+2r(c))(X^2 - 4r(c^2-d))$ (cf.~proposition~\ref{ch3explicitegeneralites}).
\end{demo}

\bigskip
Soient $\phi' \colon E' \rightarrow E''$ le morphisme canonique et $\phi'' \colon E''\rightarrow E'$ l'isogénie
duale\glossary{$\phi'$, $\phi''$, $\tgoth'$, $\tgoth''$, $\mgoth''$}.
Notons $\tgoth' \in \kappa(H)^\star/\kappa(H)^{\star 2}$ (resp.~$\tgoth'' \in \kappa(H)^\star/\kappa(H)^{\star 2}$) l'image
de l'unique point rationnel d'ordre~$2$ de~$E''$ (resp.~de~$E'$)
par la flèche $E''(\kappa(H)) \rightarrow H^1(\kappa(H),\Z/2)$
(resp.~$E'(\kappa(H)) \rightarrow H^1(\kappa(H),\Z/2)$) induite par
la suite exacte
$$
\xymatrix{
0 \ar[r] & \Z/2 \ar[r] & E' \ar[r]^{\phi'} & E'' \ar[r] & 0
}
$$
(resp.
$$
\xymatrix{
0 \ar[r] & \Z/2 \ar[r] & E'' \ar[r]^{\phi''} & E' \ar[r] & 0\rlap{\text{).}}
}
$$
Notons de plus $\mgoth'' \in \kappa(H)^\star/\kappa(H)^{\star 2}$ l'image
par la flèche $H^1(\kappa(H), \tors{2}E') \rightarrow H^1(\kappa(H), \tors{\phi''}E'')$
induite par $\phi'\colon \tors{2}E'\rightarrow \tors{\phi''}E''$ de la classe
de la fibre générique de~$\pi_{H^0}$, vue comme $2$\nobreakdash-revêtement de~$E'$
(cf.~proposition~\ref{ch3explicitegeneralites}).

Soit $\Delta \subset H^0$\glossary{$\Delta$, $\Delta'$, $\Delta''$}
l'ensemble des points de~$H^0$ au-dessus desquels la fibre de~$\pi_{H^0}$
est singulière.  C'est un fermé purement de codimension~$1$ dans~$H^0$ d'après
la proposition~\ref{ch3explicitedelta}, que l'on peut appliquer après extension des scalaires
de~$k$ à~$\ksep$ et diagonalisation simultanée de~$q_1$ et~$q_2$.
(Cette propriété résulte aussi des considérations du paragraphe~\ref{ch3parconssd}.)
Si~$h$ est le point générique d'une composante irréductible de~$\Delta$, on peut donc considérer
le modèle de Néron de~$E'$ au-dessus de $\Spec(\Orond_{H^0,h})$. Deux cas se présentent alors,
selon que l'unique point rationnel
d'ordre~$2$ de~$E'$ se spécialise sur la composante neutre de la fibre spéciale
du modèle de Néron ou non.
Notons~$\Delta'$ (resp.~$\Delta''$)
la réunion des composantes irréductibles de~$\Delta$ pour lesquelles le point rationnel d'ordre~$2$
de~$E'$ se spécialise (resp.~ne se spécialise pas) sur la composante neutre.

\bigskip
\begin{proposition}
\label{ch3specdtptppmpp}
Les classes $\mgoth''$ et $\tgoth''$
appartiennent au sous-groupe $\Gm(H \setminus (\Delta'' \cup \Pi))/2$
de~$\kappa(H)^\star/\kappa(H)^{\star 2}$.
De même, la classe $\tgoth'$ appartient
au sous-groupe $\Gm(H \setminus (\Delta' \cup \Pi))/2$.
\end{proposition}

\bigskip
\begin{demo}
Comme $\Pi$ est un hyperplan, il suffit de vérifier que les classes $\mgoth''$ et $\tgoth''$
appartiennent (resp.~que la classe~$\tgoth'$ appartient)
aux noyaux des flèches \mbox{$\kappa(H)^\star/\kappa(H)^\star \rightarrow \Z/2$} induites
par les valuations normalisées aux points de codimension~$1$ de~$H \setminus \Delta''$
(resp.~de~\mbox{$H \setminus \Delta'$}).  Fixons un tel point~$h$.
Le $k$\nobreakdash-schéma $\Spec(\Orond_{H,h})$ est un schéma de Dedekind connexe, puisque c'est un trait;
vu la proposition~\ref{ch3specdss}, nous sommes donc en situation d'appliquer
la proposition~\ref{ch2propcompct} aux modèles de Néron des courbes elliptiques~$E'$ et~$E''$
au-dessus de $\Spec(\Orond_{H,h})$.
Pour conclure, il reste à remarquer que~$\mgoth''$ et~$\tgoth''$ appartiennent
(resp.~$\tgoth'$ appartient) au groupe de $\phi''$\nobreakdash-Selmer géométrique
de~$E''$ (resp.~au groupe de $\phi'$\nobreakdash-Selmer géométrique de~$E'$) relativement
à $\Spec(\Orond_{H,h})$.  Ceci est évident pour les classes~$\tgoth'$ et~$\tgoth''$ puisqu'elles
proviennent de points rationnels; quant à~$\mgoth''$, cela résulte de ce que la fibre
de~$\pi_{H^0}$ au-dessus de~$h$ est géométriquement réduite (cf.~théorème~\ref{ch3conssdtheoreme}).
\end{demo}

\bigskip
\begin{proposition}
\label{ch3specdi1i2}
Pour tout $h \in H^0$ de codimension~$1$, la surface $C_{H^0} \times_{H^0} \Spec(\Orond_{H,h})$
est un modèle propre et régulier minimal de la fibre générique de~$\pi_{H^0}$ au-dessus
de~$\Spec(\Orond_{H,h})$ et la courbe~$\pi_{H^0}^{-1}(h)$ est
de type~$I_0$, $I_1$ ou~$I_2$.
\end{proposition}

\bigskip
Nous disons qu'une courbe est \emph{de type~$I_0$} si elle est propre, lisse, géométriquement
connexe et de genre~$1$;
\emph{de type~$I_1$} si elle est propre, géométriquement intègre,
rationnelle et qu'elle possède un unique point singulier, qui est un point double ordinaire;
\emph{de type~$I_2$} si elle est propre et qu'elle est géométriquement réunion de deux courbes
rationnelles lisses se rencontrant transversalement en deux points distincts.

\bigskip
\begin{demo}
Que cette surface soit régulière résulte de la proposition~\ref{ch3conssdregularite};
sa connexité découle de la proposition~\ref{ch3conssdgeomint}.
Notons temporairement~$S$ le spectre du hensélisé strict de~$\Orond_{H,h}$
et $C_S = C_{H^0} \times_{H^0} S$.  La surface~$C_S$ est elle aussi régulière.
La fibre spéciale de la projection $\pi_S \colon C_S \rightarrow S$ est réduite
(cf.~théorème~\ref{ch3conssdtheoreme}).  Elle admet donc un point rationnel lisse;
celui-ci se relève en une section de~$\pi_S$,
d'où un isomorphisme entre la fibre générique de~$\pi_S$ et la courbe elliptique
$E' \times_{\kappa(H)} \kappa(S)$.  Il s'ensuit qu'après éventuellement une ou plusieurs
contractions (dont on note~$m$ le nombre) de composantes irréductibles de la fibre spéciale de~$\pi_S$,
la surface~$C_S$ devient isomorphe à un modèle propre et régulier
minimal de la courbe elliptique~$E' \times_{\kappa(H)} \kappa(S)$ au-dessus de~$S$.  Étant donné que celle-ci est à réduction
semi-stable (cf.~proposition~\ref{ch3specdss}) et que la fibre spéciale de~$\pi_S$
possède au plus deux composantes irréductibles (cf.~théorème~\ref{ch3conssdtheoreme}), on en déduit que si la conclusion de la proposition
est en défaut, alors~$m=1$, la courbe elliptique~$E'$ est à réduction de type~$I_1$ en~$h$,
la fibre spéciale de~$\pi_S$ admet deux composantes irréductibles
et soit l'une de ces composantes irréductibles est de genre arithmétique~$1$
(cas où~$C_S$ est obtenu en éclatant un point régulier de la fibre spéciale d'un modèle propre et régulier
minimal de $E' \times_{\kappa(H)} \kappa(S)$), soit l'une de ces composantes irréductibles est
de multiplicité~$2$ (cas où~$C_S$ est obtenu en éclatant le point singulier).
Une telle situation ne peut arriver; en effet, d'une part nous savons que la fibre spéciale
de~$\pi_S$
est réduite, d'autre part, comme c'est une intersection de deux quadriques dans~$\P^3$,
si elle est réductible, alors ses composantes irréductibles sont des
courbes de degré~$\leq 2$ et donc de genre arithmétique nul.
\end{demo}

\bigskip
Notons $\Irr(Z)$ l'ensemble des composantes irréductibles d'une variété~$Z$.
Pour $m \in \Irr(\Delta)$, la proposition~\ref{ch3specdi1i2} montre que
la fibre géométrique de~$\pi_{H^0}$ au-dessus du point générique de~$m$
possède au plus deux composantes irréductibles et au plus deux points singuliers.
La plus petite extension de~$\kappa(m)$ par laquelle se factorise l'action
du groupe de Galois absolu de~$\kappa(m)$ sur ces deux composantes irréductibles
(resp.~sur ces deux points singuliers) est donc une extension quadratique ou triviale;
notons~$\beta_m$\glossary{$\beta_m$, $\gamma_m$}
(resp.~$\gamma_m$) sa classe dans $\kappa(m)^\star/\kappa(m)^{\star 2}$.

\bigskip
\begin{definition}
\label{ch3defcdg}
Nous dirons\index{condition (D)@\condD{}!générique} que \emph{la condition~\cDgp{} est satisfaite}
si le noyau de la flèche naturelle\glossary{\cDg{}, \cDgp{}, \cDgpp{}}
$$
\Gm(H \setminus (\Delta' \cup \Pi))/2 \longrightarrow \prod_{m \in \Irr(\Delta'')}
\kappa(m)^\star/\langle \kappa(m)^{\star 2}, \beta_m \rangle
$$
est engendré par~$\tgoth'$.
Nous dirons que \emph{la condition~\cDgpp{} est satisfaite}
si le noyau de la flèche naturelle
$$
\Gm(H \setminus (\Delta'' \cup \Pi))/2 \longrightarrow \prod_{m \in \Irr(\Delta')}
\kappa(m)^\star/\langle \kappa(m)^{\star 2}, \beta_m, \gamma_m \rangle
$$
est engendré par~$\mgoth''$ et~$\tgoth''$.
Nous dirons enfin que \emph{la condition~\cDg{}} (ou: la \emph{\condD{} générique}) \emph{est satisfaite} si les conditions~\cDgp{}
et~\cDgpp{} le sont.
\end{definition}

\bigskip
Nous sommes maintenant en position d'énoncer le résultat principal de ce paragraphe.
Notons~$\Rrond$\glossary{$\Rrond$} l'ensemble des points $p \in \Pi(k)$ tels que les conditions suivantes soient vérifiées:
\begin{itemize}
\item[(i)] le point~$p$ n'appartient pas à~$\Delta$;
\item[(ii)] la droite~$D$ de~$H$ passant par~$p$ et~$h_0$ est incluse dans~$H^0$;
\item[(iii)] la surface $\pi_{H^0}^{-1}(D)$ est lisse et géométriquement connexe sur~$k$;
\item[(iv)] le morphisme $\pi_{H^0}^{-1}(D) \rightarrow D$ induit par~$\pi_{H^0}$ est propre, plat, à fibres réduites,
de fibre générique une courbe lisse et géométriquement connexe de genre~$1$ et de période $\leq 2$
dont la jacobienne est une courbe elliptique à réduction semi-stable sur~$D$ admettant un unique point rationnel d'ordre~$2$;
\item[(v)] il existe un $k$\nobreakdash-isomorphisme $\tau \colon D \isoto \P^1_k$ tel que $\tau(p)=\infty$,
tel que la\index{condition (E)@\condE{}} \condE{} du théorème~\ref{ch2thsurmesure} soit satisfaite
pour la famille $$\tau \circ \pi_{H^0} \colon \pi_{H^0}^{-1}(D) \longrightarrow \P^1_k$$
relativement à l'ensemble
de places~$S$ que l'énoncé du théorème~\ref{ch2thsurmesure} associe à cette
famille et au point $x_0=\tau(h_0) \in \P^1(k)$,
et tel que $\tau(h_0)$ n'appartienne à la composante connexe non minorée de $U(k_v)$ pour aucune place~$v$ réelle,
où $U=\tau(D \setminus (D \cap (\Delta \cup \Pi)))$.
\end{itemize}

\bigskip
\begin{proposition}
\label{ch3specdcamarche}
Supposons que~$k$ soit un corps de nombres.
Si la condition~\cDg{} est satisfaite,
l'ensemble~$\Rrond$ est dense dans $\Pi(\A_k)$ pour la topologie adélique.
\end{proposition}

\bigskip
\begin{demo}
Choisissons, pour chaque place réelle~$v$ de~$k$, un élément $a_v \in k^\star$ qui soit négatif en~$v$ et
positif en toute autre place réelle.
Soit $S \subset \Omega$ un ensemble fini arbitrairement grand de places de~$k$ contenant les places archimédiennes,
les places dyadiques, un système de générateurs du groupe de classes de~$k$,
les places finies de mauvaise réduction pour la jacobienne de~$\pi_{H^0}^{-1}(h_0)$,
les places finies~$v$ telles que l'adhérence de $h_0 \in H \subset \P^4_k$ dans $\P^4_{\Orond_v}$
rencontre celle de $\Delta \cup \Pi$ et
enfin les places finies en lesquelles au moins l'un des~$a_v$ pour~$v$ réelle n'est pas une unité.
Nous allons exhiber un élément de~$\Rrond$ arbitrairement proche de $k_v$\nobreakdash-points
fixés de~$\Pi$ pour $v \in S$.

Si~$Z$ est un sous-schéma de~$\P^4_k$, notons~$\Ztilde$ son adhérence
dans~$\P^4_{\Orond_S}$, où~$\Orond_S$ désigne l'anneau des $S$\nobreakdash-entiers de~$k$.
Pour $m \in \Irr(\Delta'')$ (resp.~$m \in \Irr(\Delta')$),
soit $G_m \subset \kappa(m)^\star/\kappa(m)^{\star 2}$ le sous-groupe
engendré par $\beta_m$ (resp.~$\beta_m$ et~$\gamma_m$) et par l'image
de la flèche
$$\Gm\!\left(\Htilde \setminus \left(\widetilde{\Delta'\cup\Pi}\right)\right)\!/2
\longrightarrow
\kappa(m)^\star/\kappa(m)^{\star 2}$$
(resp. $$\Gm\!\left(\Htilde \setminus \left(\widetilde{\Delta''\cup\Pi}\right)\right)\!/2
\longrightarrow
\kappa(m)^\star/\kappa(m)^{\star 2}\rlap{\text{)}}$$
d'évaluation au point générique de~$m$.
Les groupes~$G_m$ sont finis en vertu du théorème des unités de Dirichlet.

Notons $q \colon H \setminus \{h_0\} \rightarrow \Pi$ le morphisme de projection
sur~$\Pi$ depuis~$h_0$ et munissons $\Delta \subset H$ de sa structure de sous-schéma
fermé réduit.  Le morphisme $\Delta \rightarrow \Pi$ induit par~$q$
est fini (puisque $h_0 \not\in \Delta$)
et génériquement étale
(puisque~$\Delta$ est réduit).
Il existe donc un ouvert dense $\Pi^0 \subset \Pi$ tel que
le morphisme $q^{-1}(\Pi^0) \cap \Delta \rightarrow \Pi^0$ induit par~$q$ soit fini étale.
Quitte à rétrécir~$\Pi^0$,
on peut supposer que $\Pi^0 \cap q(E)=\emptyset$
et donc que $q^{-1}(\Pi^0) \subset H^0$,
puisque~$E$ est de codimension~$\geq 2$ dans~$H$.
On peut aussi supposer que $\Pi^0 \cap \Delta = \emptyset$
puisque~$\Delta$ ne contient pas d'hyperplan de~$H$ (cf.~proposition~\ref{ch3explicitedelta},
que l'on peut appliquer après extension des scalaires de~$k$ à~$\ksep$; au lieu d'utiliser cette
propriété non triviale, on pourrait tout aussi bien prendre pour hypothèse
que $\Pi \not\subset\Delta$ dans tout le paragraphe).

\bigskip
\begin{lemme}
Quitte à rétrécir l'ouvert $\Pi^0 \subset \Pi$, on peut supposer
que $\pi_{H^0}^{-1}(D)$ est une surface lisse et connexe pour tout
$p \in \Pi^0$, notant~$D$ la droite de~$H$ passant par~$p$ et~$h_0$.
\end{lemme}

\bigskip
\begin{demo}
Soit $E' \subset H$ la réunion de~$E$ et de l'image par~$\pi_{H^0}$ du lieu singulier
de~$C_{H^0}$.
Posons
$$
I=\bigensemble{(p,h)\in \Pi\times_k (H\setminus E')}{p, h \text{ et } h_0 \text{ sont alignés dans } H}\!\rlap{\text{.}}
$$
La seconde projection $I \rightarrow H \setminus E'$ permet d'identifier~$I$ à l'image réciproque
de $H \setminus E'$ dans le $H$\nobreakdash-schéma obtenu en faisant éclater~$h_0$.
Notons $C_{H\setminus E'}=\pi_{H^0}^{-1}(H\setminus E')$ et
$C_I=C_{H\setminus E'} \times_{H \setminus E'} I$.
Comme le morphisme~$\pi_{H^0}$ est plat (cf.~proposition~\ref{ch3conssdplatcm}) et que les éclatements commutent à de tels changements
de base,
la première projection $C_I \rightarrow C_{H \setminus E'}$ permet d'identifier~$C_I$
à l'image réciproque de $C_{H\setminus E'}$ dans le $C_{H^0}$-schéma obtenu en faisant éclater $\pi_{H^0}^{-1}(h_0)$.
Comme $C_{H\setminus E'}$ et $\pi_{H^0}^{-1}(h_0)$ sont des variétés lisses et connexes ($C_{H\setminus E'}$ est
lisse par définition de~$E'$, connexe d'après la proposition~\ref{ch3conssdgeomint}),
il en résulte que~$C_I$ est lisse et connexe (cf.~\cite[II, 8.24]{hartshorne}).
En particulier, c'est une variété régulière et irréductible, d'où l'on tire
que la fibre générique de la composée des projections
$C_I \rightarrow I \rightarrow \Pi$ est régulière et irréductible.
La proposition s'ensuit, étant donné que~$q(E')$ est un fermé strict de~$\Pi$
(cf.~proposition~\ref{ch3conssdregularite}).
\end{demo}

\bigskip
Pour $m \in \Irr(\Delta)$, chaque $g \in G_m \setminus \{1\}$ définit
un revêtement double connexe de~$m$.  Quitte à rétrécir~$\Pi^0$, on peut
supposer que pour tout $m \in \Irr(\Delta)$ et tout $g \in G_m \setminus \{1\}$,
le revêtement double de~$m$ défini par~$g$ est étale au-dessus de $q^{-1}(\Pi^0)\cap m$.
Nous noterons alors $g(M) \in \kappa(M)^\star/\kappa(M)^{\star 2}$ la classe
de la fibre en~$M$ du revêtement associé à $g\in G_m$,
pour $M \in q^{-1}(\Pi^0)\cap m$.

Quitte à rétrécir encore~$\Pi^0$, on peut supposer que pour tout $m \in \Irr(\Delta)$,
si la fibre de~$\pi_{H^0}$ au-dessus du point générique
de~$m$ est de type~$I_1$ (resp.~$I_2$), alors toute fibre de~$\pi_{H^0}$ au-dessus
de $q^{-1}(\Pi^0)\cap m$ est géométriquement intègre
(resp.~est réduite et possède géométriquement deux composantes irréductibles,
qui sont des courbes rationnelles lisses se rencontrant en deux points distincts).
(Cela résulte de \cite[9.7.7 et 9.9.5]{ega43}, que l'on applique
d'une part au morphisme obtenu à partir de $\pi_{H^0}^{-1}(m) \rightarrow m$ par le changement de base
de degré~$2$ défini par~$\beta_m$ suivi de la restriction à une composante irréductible
de l'espace total, d'autre part, dans le cas~$I_2$, au morphisme obtenu par restriction de $\pi_{H^0}^{-1}(m)\rightarrow m$
au lieu de non lissité, changement de base de degré~$2$ défini par~$\gamma_m$ puis
restriction à une composante irréductible de l'espace total.)

Enfin, il résulte du lemme~\ref{ch3harrmi}
qu'un dernier rétrécissement de~$\Pi^0$ permet
de supposer que pour tout $m \in \Irr(\Delta)$ tel que la fibre de~$\pi_{H^0}$ au-dessus
du point générique de~$m$ soit de type~$I_2$ et tout $M \in q^{-1}(\Pi^0) \cap m$,
la plus petite extension de~$\kappa(M)$ par laquelle se factorise l'action du groupe de Galois
absolu de~$\kappa(M)$ sur les deux composantes irréductibles
(resp.~sur les deux points singuliers) de la fibre géométrique de~$\pi_{H^0}$ en~$m$
est une extension quadratique ou triviale dont la classe dans $\kappa(M)^\star/\kappa(M)^{\star 2}$
est égale à~$\beta_m(M)$ (resp.~$\gamma_m(M)$).

Nous disposons maintenant d'un certain nombre de revêtements étales connexes
de~$\Pi^0$: les revêtements connexes $q^{-1}(\Pi^0) \cap m \rightarrow \Pi^0$ pour
$m \in \Irr(\Delta)$ (ils sont étales parce que $q^{-1}(\Pi^0) \cap \Delta \rightarrow \Pi^0$ l'est
par hypothèse) et les revêtements obtenus en composant ceux-ci
avec les revêtements étales doubles connexes de $q^{-1}(\Pi^0) \cap m$ définis
par les éléments de $G_m\setminus\{1\}$.
D'après le théorème d'irréductibilité de Hilbert\index{théorème d'irréductibilité de Hilbert} avec approximation faible
(cf.~\cite{ekedahl}), il existe $p \in \Pi^0(k)$ arbitrairement proche
de $k_v$\nobreakdash-points fixés de~$\Pi$ aux places $v \in S$, tel que les fibres en~$p$ de tous
ces revêtements de~$\Pi^0$ soient intègres.

Il reste à vérifier que $p \in \Rrond$. Notons~$D$ la droite de~$H$ passant par~$p$ et~$h_0$.
Les conditions~(i) et~(ii) de la définition de~$\Rrond$ sont satisfaites par construction de~$\Pi^0$;
de même pour la lissité et la connexité de la surface $\pi_{H^0}^{-1}(D)$. Sachant que cette surface
est lisse et connexe, sa connexité géométrique résultera de~(iv).
Il suffit donc d'établir les propriétés~(iv) et~(v).
Le morphisme $\pi_{H^0}^{-1}(D) \rightarrow D$ est évidemment propre.
Il est plat d'après la proposition~\ref{ch3conssdplatcm}.
Ses fibres sont réduites par construction de~$\Pi^0$: plus précisément, nous savons
que ses fibres géométriques singulières sont intègres ou sont des réunions de deux
courbes rationnelles lisses.  Comme $h_0 \not\in \Delta$, le point générique de~$D$ n'appartient
pas à~$\Delta$; la fibre générique de ce morphisme est donc lisse.  C'est une courbe géométriquement
connexe de genre~$1$ d'après la proposition~\ref{ch3conssdfibrespi},
de période~$\leq 2$ d'après la proposition~\ref{ch3explicitegeneralites}.
Notons~$E'_0$ la jacobienne de la fibre générique du morphisme $\pi_{H^0}^{-1}(D) \rightarrow D$.
Le lemme suivant termine de prouver que la propriété~(iv) est satisfaite.

\bigskip
\begin{lemme}
\label{ch3specde0semi}
La courbe elliptique~$E'_0$ est à réduction semi-stable en tout point fermé de~$D$.
Elle possède un unique point rationnel d'ordre~$2$, et si l'on note~$E''_0$ le
quotient de~$E'_0$ par ce point, la courbe elliptique~$E''_0$ possède elle aussi
un unique point rationnel d'ordre~$2$.
De plus, l'ensemble des points fermés de~$D$
de mauvaise réduction pour~$E'_0$ au-dessus desquels le point rationnel d'ordre~$2$ de~$E'_0$
se spécialise (resp.~ne se spécialise pas) sur la composante neutre du modèle de Néron
est égal à~$D \cap \Delta'$ (resp.~$D \cap \Delta''$).
\end{lemme}

\bigskip
\begin{demo}
Compte tenu qu'il est déjà connu que~$E'_0$ admet un point d'ordre~$2$ rationnel
(cf.~proposition~\ref{ch3explicitegeneralites}), on peut s'autoriser une extension des
scalaires pour démontrer le lemme.  On peut donc supposer les formes quadratiques~$q_1$ et~$q_2$
simultanément diagonalisables (cf.~proposition~\ref{ch3gensimdiag}). Soit $h \in D \cap \Delta$. Dans les arguments
ci-dessous, le point~$h$ sera parfois soumis à une contrainte, parfois quelconque.
Reprenons les notations du paragraphe~\ref{ch3calcexpl}; une transformation linéaire
de~$\P^4_k$ permet de supposer que~$q_1$ et~$q_2$ sont simultanément diagonales dans la base~$e$,
que le point~$\fp{P}0$ a pour coordonnées homogènes $[1:0:0:0:0]$ et que $(y_1\gamma_4)(h)\neq 0$.
(La même remarque que dans la preuve
de la proposition~\ref{ch3specdss} vaut ici: si l'on veut que le point~$h$ soit rationnel, comme
cela avait été supposé au paragraphe~\ref{ch3calcexpl}, il suffit d'étendre les scalaires de~$k$ à~$\kappa(h)$
avant d'appliquer les résultats de ce paragraphe. Nous ne referons plus cette remarque à l'avenir.)
Notons $r \colon k[\uplet{y_1}{y_4}]\rightarrow \Orond_{D,h}$ le morphisme d'anneaux qui à $f(\uplet{y_1}{y_4})$
associe le germe de fonction induit par $f(1,y_2/y_1,y_3/y_1,y_4/y_1)$
et~$v$ la valuation discrète normalisée associée à l'anneau de valuation discrète $\Orond_{D,h}$.
La courbe elliptique~$E'_0$ a pour équation de Weierstrass
$$
Y^2=(X-r(c))(X^2-r(d))
$$
(cf.~proposition~\ref{ch3explicitegeneralites}).
Compte tenu que $v(r(y_1))=v(r(\gamma_4))=0$, la proposition~\ref{ch3expliciteegalitessim}
montre que
\begin{equation}
\label{ch3specvrd}
v(r(d))=2v(r(p_0))+v(r(p_6))
\end{equation}
et
\begin{equation}
\label{ch3specdvrcd}
v(r(c^2-d))=\sum_{i=1}^5 v(r(p_i))\rlap{\text{.}}
\end{equation}
Les hypersurfaces de~$H^0$ d'équation $p_i=0$ pour $i\in\{\uplet{0}{6}\}$ sont
les composantes irréductibles de~$\Delta$ (cf.~proposition~\ref{ch3explicitedelta}).
Par ailleurs, comme par construction de~$\Pi^0$ le schéma $D \cap \Delta$ est réduit,
le point~$h$ appartient à une unique composante irréductible de~$\Delta$ et celle-ci
rencontre~$D$ transversalement en~$h$.  Il s'ensuit qu'il existe un unique $i(h) \in \{\uplet{0}{6}\}$
tel que $v(r(p_{i(h)}))\neq 0$ et que l'on a alors $v(r(p_{i(h)}))=1$.
Les équations~(\ref{ch3specvrd}) et~(\ref{ch3specdvrcd}) permettent d'en déduire
que l'un de~$v(r(d))$ et de~$v(r(c^2-d))$ est nul, ce qui assure
que la courbe elliptique~$E'_0$ est à réduction semi-stable en~$h$
(cf.~preuve de la proposition~\ref{ch3specdss} pour les détails).

Choisissons temporairement le point~$h$ dans
la composante irréductible de~$\Delta$
d'équation $p_6=0$.
On a alors $i(h)=6$ et donc $v(r(d))=1$, vu l'équation~(\ref{ch3specvrd}), ce qui
entraîne que~$r(d)$ n'est pas un carré dans~$\kappa(D)$.  Nous avons ainsi prouvé que la courbe elliptique~$E'_0$
admet un unique point rationnel d'ordre~$2$.
Choisissons maintenant le point~$h$ dans la composante irréductible de~$\Delta$ d'équation $p_5=0$.
On a alors $i(h)=5$ et donc $v(r(c^2-d))=1$, vu l'équation~(\ref{ch3specdvrcd}), ce qui
entraîne que~$r(c^2-d)$ n'est pas un carré dans~$\kappa(D)$.  Nous avons ainsi prouvé que la courbe elliptique~$E''_0$
admet un unique point rationnel d'ordre~$2$, étant donné qu'elle a pour équation
de Weierstrass $Y^2=(X+2r(c))(X^2-4r(c^2-d))$.

Revenons à la situation où le point~$h$ est quelconque.
Le point rationnel d'ordre~$2$ de~$E'_0$
se spécialise sur la composante neutre du modèle de Néron de~$E'_0$ au-dessus de $\Spec(\Orond_{D,h})$
si et seulement si $v(r(d))>0$, autrement dit, si et seulement si~$h$ appartient à l'hypersurface
d'équation $p_0p_6=0$.  Le même calcul, effectué au niveau des points génériques des composantes
irréductibles de~$\Delta$, toujours en supposant~$q_1$ et~$q_2$ simultanément diagonales,
montre que $\Delta'$ est l'hypersurface de~$H^0$ définie par l'équation $p_0p_6=0$
(cf.~preuve de la proposition~\ref{ch3specdss} pour les détails).  Ceci démontre la dernière assertion du lemme.
\end{demo}

\bigskip
\begin{lemme}
\label{ch3specde0i1i2}
La surface $\pi^{-1}_{H^0}(D)$ est un modèle propre et régulier minimal de la fibre de~$\pi_{H^0}^{-1}$
au-dessus du point générique de~$D$.  Les fibres singulières du morphisme $\pi_{H^0}^{-1}(D) \rightarrow D$
sont de type~$I_1$ ou~$I_2$.
\end{lemme}

\bigskip
\begin{demo}
Nous savons déjà que cette surface est propre, lisse et connexe et que les fibres géométriques
réductibles de $\pi_{H^0}^{-1}(D) \rightarrow D$ sont des réunions de deux courbes rationnelles lisses.
Sachant que~$E'_0$ est à réduction semi-stable en tout point fermé de~$D$ (cf.~lemme~\ref{ch3specde0semi}),
il s'ensuit que cette surface est relativement minimale au-dessus de~$D$
et que ses fibres singulières sont de type~$I_1$ ou~$I_2$ (cf.~fin de la preuve de la proposition~\ref{ch3specdi1i2}).
\end{demo}

\bigskip
Soient $\phi'_0 \colon E'_0 \rightarrow E''_0$ le morphisme canonique et $\phi''_0 \colon E''_0 \rightarrow E'_0$
l'isogénie duale. Notons \mbox{$\tgoth'_0 \in \kappa(D)^\star/\kappa(D)^{\star 2}$}
(resp.~$\tgoth''_0 \in \kappa(D)^\star/\kappa(D)^{\star 2}$) l'image de l'unique point rationnel d'ordre~$2$ de~$E''_0$
(resp.~de~$E'_0$) par la flèche $E''_0(\kappa(D)) \rightarrow H^1(\kappa(D),\Z/2)$
(resp.~$E'_0(\kappa(D)) \rightarrow H^1(\kappa(D),\Z/2)$) induite par la suite exacte
$$
\xymatrix{
0 \ar[r] & \Z/2 \ar[r] & E'_0 \ar[r]^{\phi'_0} & E''_0 \ar[r] & 0
}
$$
(resp.
$$
\xymatrix{
0 \ar[r] & \Z/2 \ar[r] & E''_0 \ar[r]^{\phi''_0} & E'_0 \ar[r] & 0\rlap{\text{).}}
}
$$
Notons de plus $\mgoth''_0 \in \kappa(D)^\star/\kappa(D)^{\star 2}$
l'image par la flèche $H^1(\kappa(D), \tors{2}E'_0) \rightarrow H^1(\kappa(D),\tors{\phi''_0}E''_0)$
induite par $\phi'_0 \colon \tors{2}E'_0 \rightarrow \tors{\phi''_0}E''_0$ de la classe
de la fibre de~$\pi_{H^0}$ au-dessus du point générique de~$D$, vue comme $2$\nobreakdash-revêtement
de~$E'_0$ (cf.~proposition~\ref{ch3explicitegeneralites}).

\bigskip
\begin{lemme}
\label{ch3specdgothspec}
L'image de~$\tgoth'$ (resp.~$\tgoth''$, $\mgoth''$) par la flèche
$$
\Gm(H \setminus (\Delta \cup \Pi))/2 \longrightarrow \kappa(D)^\star/\kappa(D)^{\star 2}
$$
d'évaluation au point générique de~$D$ est égale à $\tgoth'_0$ (resp.~$\tgoth''_0$, $\mgoth''_0$).
\end{lemme}

\bigskip
\begin{demo}
Notons~$h$ le point générique de~$D$ et
reprenons les notations du paragraphe~\ref{ch3calcexpl},
de sorte que $(y_1\gamma_4)(h)\neq 0$.
Soit $r \colon k[\uplet{y_1}{y_4}] \rightarrow \Orond_{H,h}$ le morphisme
d'anneaux qui à $f(\uplet{y_1}{y_4})$ associe $f(1,y_2/y_1,y_3/y_1,y_4/y_1)$.
Il résulte de la proposition~\ref{ch3explicitegeneralites} que~$\tgoth'$
(resp.~$\tgoth''$, $\mgoth''$)
est égal à la classe de $r(d)$ (resp.~$r(c^2-d)$, $r(\dgras_{14}^2-\dgras_{12}\dgras_{13})$)
dans $\Orond_{H,h}^\star/\Orond_{H,h}^{\star 2}$
et que~$\tgoth'_0$ (resp.~$\tgoth''_0$, $\mgoth''_0$) est égal à l'image
de $r(d)$ (resp.~$r(c^2-d)$, $r(\dgras_{14}^2-\dgras_{12}\dgras_{13})$)
par la flèche naturelle $\Orond_{H,h}^\star \rightarrow \kappa(h)^\star/\kappa(h)^{\star 2}$,
d'où le lemme.
\end{demo}

\bigskip
\begin{lemme}
\label{ch3specdexistetau}
Il existe un $k$\nobreakdash-isomorphisme $\tau \colon D \isoto \P^1_k$ tel que $\tau(p)=\infty$
et qui se prolonge en un $\Orond_S$\nobreakdash-isomorphisme $\Dtilde \isoto \P^1_{\Orond_S}$.
\end{lemme}

\bigskip
\begin{demo}
C'est essentiellement une conséquence de l'hypothèse que $\Pic(\Orond_S)=0$.
Celle-ci entraîne en effet que le groupe $\PGL_{n+1}(\Orond_S)$ agit transitivement
sur $\P^n(\Orond_S)$ pour tout $n \geq 1$.
De manière générale, pour toute droite $L \subset \P^4_k$
et tout point $\ell \in L(k)$,
il existe un $\Orond_S$\nobreakdash-isomorphisme $\Ltilde \isoto \P^1_{\Orond_S}$
envoyant~$\ell$ sur~$\infty$.  En effet, la transitivité de l'action
de $\PGL_5(\Orond_S)$ sur~$\P^4(\Orond_S)$, puis
de $\PGL_4(\Orond_S)$ sur~$\P^3(\Orond_S)$, puis
de $\PGL_3(\Orond_S)$ sur~$\P^2(\Orond_S)$, permet
de supposer que la droite~$L$ est incluse dans l'hyperplan d'équation $x_0=0$,
puis qu'elle est incluse dans le sous-espace linéaire d'équation $x_0=x_1=0$,
puis qu'elle a pour équation $x_0=x_1=x_2=0$,
où $[x_0:\cdots:x_4]$ sont les coordonnées homogènes de~$\P^4_k$; il est alors
évident que~$\Ltilde$ est $\Orond_S$\nobreakdash-isomorphe à~$\P^1_{\Orond_S}$
et l'on conclut grâce à la transitivité de l'action de $\PGL_2(\Orond_S)$ sur~$\P^1(\Orond_S)$.
\end{demo}

\bigskip
Fixons~$\tau$ comme dans le lemme~\ref{ch3specdexistetau}.
Pour chaque place réelle $v \in \Omega$, on peut supposer, quitte à remplacer~$\tau$ par sa composée avec l'automorphisme $t \mapsto a_v t$ de~$\P^1_k$,
que~$\tau(h_0)$ n'appartient pas à la composante connexe non minorée de $\tau(D \setminus (D \cap (\Delta \cup \Pi)))(k_v)$.

Nous savons que
$\pi_{H^0}^{-1}(D) \xrightarrow{\;\tau \circ \pi_{H^0}} \P^1_k$
est une famille
de courbes de genre~$1$ satisfaisant aux hypothèses générales du chapitre~2.  Nous pouvons donc reprendre les
notations de ce chapitre relatives à cette famille, notamment $\Mrond$, $\Mrond'$, $\Mrond''$,
$\SG_{\phi',S}(\A^1_k, \Erond')$, $\SG_{\phi'',S}(\A^1_k, \Erond'')$, $L'_M$, $L''_M$, $\ff{F'}M$, $\ff{F''}M$, $\delta_M'$,
$\delta_M''$, $[\Xrond]$, $[\Xrond'']$.  (Ceci rend la notation~$\phi'$, $\phi''$ conflictuelle; le contexte
permettra facilement de lever toute ambiguïté.)
Soient respectivement $\widetilde{\Mrond'}$ et $\widetilde{\Mrond''}$ les adhérences
de~$\Mrond'$ et~$\Mrond''$ dans~$\A^1_{\Orond_S}$.  Dans le lemme suivant, nous identifions le corps des fonctions
de~$\P^1_k$ à~$\kappa(D)$ \emph{via}~$\tau$.

\bigskip
\begin{lemme}
\label{ch3specdconde}
Pour que la \condE{} du théorème~\ref{ch2thsurmesure} soit satisfaite relativement à l'ensemble de places
que l'énoncé de ce théorème associe au point $x_0=\tau(h_0)$, il suffit que le noyau de la flèche naturelle
$$
\Gm\!\left(\A^1_{\Orond_S} \setminus \widetilde{\Mrond'}\right)\!/2
\longrightarrow \prod_{M \in \Mrond''} \kappa(M)^\star/\langle\kappa(M)^{\star 2},\beta_m(M)\rangle
$$
soit inclus dans le sous-groupe de $\kappa(D)^\star/\kappa(D)^{\star 2}$ engendré par~$\tgoth'_0$
et que le noyau de la flèche naturelle
$$
\Gm\!\left(\A^1_{\Orond_S} \setminus \widetilde{\Mrond''}\right)\!/2
\longrightarrow \prod_{M \in \Mrond'} \kappa(M)^\star/\langle\kappa(M)^{\star 2},
\beta_m(M),\gamma_m(M)\rangle
$$
soit inclus dans le sous-groupe de $\kappa(D)^\star/\kappa(D)^{\star 2}$ engendré par~$\tgoth''_0$ et~$\mgoth''_0$.
\end{lemme}

\bigskip
\begin{demo}
Étant donné que~$\tau$ se prolonge en un $\Orond_S$\nobreakdash-isomorphisme $\Dtilde\isoto \P^1_{\Orond_S}$,
que $\tau(D \cap (\Delta \cup \Pi))=\Mrond\cup\{\infty\}$ et
que l'adhérence de~$h_0$
dans~$\P^4_{\Orond_S}$ ne rencontre pas celle de~$\Delta \cup \Pi$ (par définition de~$S$),
l'adhérence de~$\tau(h_0)$ dans~$\P^1_{\Orond_S}$
ne rencontre pas celle de $\Mrond\cup\{\infty\}$. Il s'ensuit que l'ensemble de places que l'énoncé
du théorème~\ref{ch2thsurmesure} associe au point~$\tau(h_0)$ est inclus dans~$S$, et il suffit donc de
vérifier la \condE{} relative à~$S$.

Les groupes $\SG_{\phi',S}(\A^1_k, \Erond')$ et
$\SG_{\phi'',S}(\A^1_k, \Erond'')$ s'identifient respectivement à
$\Gm\!\left(\A^1_{\Orond_S} \setminus \widetilde{\Mrond'}\right)\!/2$
et $\Gm\!\left(\A^1_{\Orond_S} \setminus \widetilde{\Mrond''}\right)\!/2$,
d'après la proposition~\ref{ch2propcompct} et l'hypothèse selon laquelle $\Pic(\Orond_S)=0$.
Les groupes~$\ff{F'}M$ pour $M \in \Mrond$ sont tous d'ordre~$\leq 2$
puisque les fibres géométriques de $\pi_{H^0}^{-1}(D) \rightarrow D$ comportent au plus deux composantes irréductibles.
Par conséquent, les extensions $L'_M/\kappa(M)$ pour $M \in \Mrond$ sont toutes triviales
et le groupe~$H^1(L'_M,\ff{F'}M)$ se plonge naturellement dans $\kappa(M)^\star/\kappa(M)^{\star 2}$.
Pour tout $M \in \Mrond''$, l'image de~$\delta_M'([\Xrond])$
dans $\kappa(M)^\star/\kappa(M)^{\star 2}$ est égale à la classe de la
fermeture algébrique de~$\kappa(M)$ dans une composante irréductible de $\pi_{H^0}^{-1}(\tau^{-1}(M))$.
Cette dernière n'est autre que~$\beta_m(M)$, par hypothèse; d'où un plongement naturel
$$\prod_{M \in \Mrond''} H^1(L'_M,\ff{F'}M)/\langle \delta_M'([\Xrond])\rangle\longhookrightarrow
\prod_{M \in \Mrond''} \kappa(M)^\star/\langle \kappa(M)^{\star 2}, \beta_m(M)\rangle\rlap{\text{.}}$$
De plus, un diagramme commutatif analogue à~(\ref{ch2compdiag2}) (précisément, celui obtenu en remplaçant~$\P^1_k$
par~$\A^1_k$) assure la commutativité du carré
\begin{equation}
\label{ch3specdcarre1}
\begin{aligned}
\myxyhook\xymatrix{
\SG_{\phi',S}(\A^1_k, \Erond') \ar[d]^\wr \ar[r] & \displaystyle\prod_{M\in\Mrond''}\frac{H^1(L'_M,\ff{F'}M)}
{\langle\delta_M'([\Xrond])\rangle} \ar@{_{ (}->}[d] \\
\Gm\!\left(\A^1_{\Orond_S} \setminus \widetilde{\Mrond'}\right)\!/2 \ar[r] &
\displaystyle\prod_{M\in \Mrond''}\kappa(M)^\star / \langle\kappa(M)^{\star 2},\beta_m(M)\rangle \rlap{\text{,}}
}
\end{aligned}
\end{equation}
dont on vient de construire les flèches verticales, dont la flèche horizontale supérieure est
la composée de la flèche $\SG_{\phi',S}(\A^1_k,\Erond')\rightarrow \tors{\phi'}H^1(\A^1_k,\Erond')$
issue de la suite exacte~(\ref{ch2sesgtsgp}) et de la flèche~(\ref{ch2flechedeftdcp}) et dont la flèche horizontale
inférieure est celle apparaissant dans l'énoncé de ce lemme.
On prouve de la même manière l'existence d'un carré commutatif
\begin{equation}
\label{ch3specdcarre2}
\begin{aligned}
\myxyhook\xymatrix{
\SG_{\phi'',S}(\A^1_k, \Erond'') \ar[d]^\wr \ar[r] & \displaystyle\prod_{M\in\Mrond'}\frac{\tors{2}H^1(L''_M,\ff{F''}M)}
{\langle\delta_M''([\Xrond''])\rangle} \ar@{_{ (}->}[d] \\
\Gm\!\left(\A^1_{\Orond_S} \setminus \widetilde{\Mrond''}\right)\!/2 \ar[r] &
\displaystyle\prod_{M\in \Mrond'}L''^\star_M / \langle L''^{\star 2}_M,\beta_m(M)\rangle \rlap{\text{.}}
}
\end{aligned}
\end{equation}
La seule différence est que les groupes~$\ff{F''}M$ pour $M \in \Mrond$ ne sont pas tous d'ordre~$\leq 2$
et que les extensions $L''_M/\kappa(M)$ ne sont pas toutes triviales, mais il n'en reste pas moins vrai
que~$\tors{2}H^1(L''_M,\ff{F''}M)=H^1(L''_M,\ff{\tors{2}F''}M)=L''^\star_M/L''^{\star 2}_M$.

Le lemme~\ref{ch3specde0i1i2} montre que pour tout $M \in \Mrond$, après extension des scalaires de~$\kappa(M)$ à l'extension
quadratique ou triviale définie par~$\beta_m(M)$, la fibre en~$M$ d'un
modèle propre et régulier minimal
de~$E'_0$ au-dessus de~$D$ s'identifie à~$\pi_{H^0}^{-1}(M)$.  Il en résulte, compte tenu du
corollaire~\ref{anncorsingrat}, que pour $M \in \Mrond'$, l'image dans
$\kappa(M)^\star/\langle\kappa(M)^{\star 2},\beta_m(M)\rangle$
de la classe de l'extension quadratique ou triviale $L''_M/\kappa(M)$
est égale à l'image de~$\gamma_m(M)$. Le noyau de la flèche horizontale inférieure du carré~(\ref{ch3specdcarre2})
coïncide donc avec le noyau de la seconde flèche de l'énoncé de ce lemme.

Nous avons maintenant établi que la condition du lemme est satisfaite si et seulement
si les noyaux des flèches horizontales supérieures des carrés commutatifs~(\ref{ch3specdcarre1}) et~(\ref{ch3specdcarre2}),
qui s'identifient à $\SDzp \cap \SG_{\phi',S}(\A^1_k, \Erond')$
et à $\SDzpp \cap \SG_{\phi'',S}(\A^1_k, \Erond'')$, sont respectivement inclus
dans les sous-groupes de~$\kappa(D)^\star/\kappa(D)^{\star 2}$ engendrés par~$\tgoth'_0$
et par $\tgoth''_0$ et~$\mgoth''_0$.  Cette dernière condition implique bien la \condE{}, puisque
l'image de~$\mgoth''_0$ dans le groupe $\tors{\phi''}H^1(\A^1_k,\Erond'')$ n'est autre que
la classe~$[\Xrond'']$.
\end{demo}

\bigskip
Le lemme~\ref{ch3specde0semi} montre que $\tau(D\cap \Delta')=\Mrond'$ et $\tau(D\cap \Delta'')=\Mrond''$,
ce qui permet de définir des morphismes
$\Gm(H \setminus (\Delta' \cup \Pi))/2 \rightarrow \Gm(\A^1_k \setminus \Mrond')/2$
et
$\Gm(H \setminus (\Delta'' \cup \Pi))/2 \rightarrow \Gm(\A^1_k \setminus \Mrond'')/2$
par restriction à $\tau^{-1}(\A^1_k \setminus \Mrond')$
et à~$\tau^{-1}(\A^1_k\setminus\Mrond'')$.

\bigskip
\begin{lemme}
\label{ch3specdrestiso}
Les morphismes de restriction
$\Gm(H \setminus (\Delta' \cup \Pi))/2 \rightarrow \Gm(\A^1_k \setminus \Mrond')/2$
et
$\Gm(H \setminus (\Delta'' \cup \Pi))/2 \rightarrow \Gm(\A^1_k \setminus \Mrond'')/2$
sont des isomorphismes.
Ils induisent par restriction des isomorphismes
$$\Gm\!\left(\Htilde \setminus \left(\widetilde{\Delta'\cup\Pi}\right)\right)\!/2
\longisoto \Gm\!\left(\A^1_{\Orond_S} \setminus \widetilde{\Mrond'}\right)\!/2$$
et
$$\Gm\!\left(\Htilde \setminus \left(\widetilde{\Delta''\cup\Pi}\right)\right)\!/2
\longisoto
\Gm\!\left(\A^1_{\Orond_S} \setminus \widetilde{\Mrond''}\right)\!/2 \rlap{\text{.}}
$$
\end{lemme}

%\bigskip
\begin{demo}
Intéressons-nous au premier de ces morphismes seulement, la preuve pour le second étant symétrique.
Le schéma $D \cap m$ est par hypothèse intègre pour tout $m\in\Irr(\Delta')$.
Il en résulte tout d'abord que
l'inclusion \mbox{$D\cap \Delta' \subset \Delta'$} induit une bijection entre composantes irréductibles,
ce qui entraîne que l'ensemble sous-jacent à~$\Mrond'$ s'identifie à~$\Irr(\Delta')$
et permet ainsi de définir le diagramme
$$
\xymatrix{
0 \ar[r] & \Gm(k)/2 \ar@{=}[d] \ar[r] & \Gm(H \setminus (\Delta'\cup\Pi))/2 \ar[d] \ar[r] &
(\Z/2)^{\Irr(\Delta')} \ar@{=}[d] \ar[r] & 0 \\
0 \ar[r] & \Gm(k)/2 \ar[r] & \Gm(\A^1_k \setminus \Mrond')/2 \ar[r] & (\Z/2)^{\Mrond'} \ar[r] & 0 \rlap{\text{,}}
}
$$
dont les flèches horizontales de droite sont induites par les valuations normalisées.
Il en résulte ensuite que la droite~$D$ rencontre transversalement chacune des composantes irréductibles de~$\Delta'$,
de sorte que ce diagramme est commutatif. Ses lignes étant évidemment exactes, le lemme des cinq permet de conclure quant à
la première assertion du lemme.

Considérons maintenant le carré
$$
\xymatrix{
\Gm(H \setminus (\Delta'\cup\Pi))/2 \ar[r] \ar[d]^\wr & (\Z/2)^{\Omega \setminus S} \ar@{=}[d] \\
\Gm(\A^1_k \setminus \Mrond')/2 \ar[r] & (\Z/2)^{\Omega \setminus S} \rlap{\text{,}}
}
$$
dans lequel les flèches horizontales sont induites par les valuations normalisées aux points
de codimension~$1$ de~$\A^1_{\Orond_S}$ et de~$\Htilde$ associés aux places de~$\Omega \setminus S$.
Les noyaux des flèches horizontales sont précisément
les sous-groupes
$\Gm\!\left(\Htilde \setminus \left(\widetilde{\Delta'\cup\Pi}\right)\right)\!/2$
et $\Gm\!\left(\A^1_{\Orond_S} \setminus \widetilde{\Mrond'}\right)\!/2$, puisque $\Pic(\Orond_S)=0$;
tout ce qu'il reste à faire pour établir la seconde assertion du lemme est donc de prouver
que ce carré est commutatif.
Notons $\xi \in \Dtilde$ le point générique de la fibre en~$v$ du morphisme
structural $\Dtilde \rightarrow \Spec(\Orond_S)$.
Comme~$\tau$ s'étend en un $\Orond_S$\nobreakdash-isomorphisme $\Dtilde\isoto \P^1_{\Orond_S}$,
il suffit de vérifier que pour toute fonction $f \in \Gm(H \setminus (\Delta' \cup \Pi))$ et toute
place $v \in \Omega \setminus S$, la valuation de~$f$ au point générique de la fibre en~$v$ du morphisme
structural $\Htilde \rightarrow \Spec(\Orond_S)$
est égale à la valuation en~$\xi$ de la
restriction de~$f$ à $D \setminus (D \cap (\Delta' \cup \Pi))$. On peut supposer la première de ces
deux valuations nulle, grâce à l'hypothèse $\Pic(\Orond_S)=0$.  Il suffit alors pour conclure de prouver que~$\xi$
n'appartient pas à $\widetilde{\Delta'\cup\Pi}$; mais si tel était le cas, la fibre en~$v$
de~$\Dtilde$ serait incluse dans $\widetilde{\Delta'\cup\Pi}$ et en particulier l'adhérence de~$h_0$
dans~$\Htilde$ rencontrerait celle de~$\Delta'\cup\Pi$, ce qui est exclu par définition de~$S$.
\end{demo}

\bigskip
Pour tout $m \in \Irr(\Delta)$,
on dispose d'une flèche d'évaluation $G_m \rightarrow \kappa(M)^\star/\kappa(M)^{\star 2}$,
où $M=D \cap m$,
puisque les revêtements de~$m$ définis par les éléments de~$G_m$ sont par hypothèse étales au-dessus
d'un voisinage de~$M$.  Cette flèche est injective par construction de~$p$;
les flèches
$$G_m/\langle \beta_m(M)\rangle \rightarrow \kappa(M)^{\star}/\langle\kappa(M)^{\star 2},\beta_m(M)\rangle$$
(pour $m \in \Irr(\Delta'')$)
et
$$G_m/\langle \beta_m(M),\gamma_m(M)\rangle \rightarrow \kappa(M)^{\star}/\langle\kappa(M)^{\star 2},\beta_m(M),\gamma_m(M)\rangle$$
(pour $m \in \Irr(\Delta')$)
qu'elle induit sont donc elles aussi injectives.
Celles-ci s'insèrent dans les carrés commutatifs
$$
\myxyhook\xymatrix{
\Gm\!\left(\Htilde \setminus \left(\widetilde{\Delta'\cup\Pi}\right)\right)\!/2 \ar[r] \ar[d]^\wr &
\displaystyle\prod_{m\in \Irr(\Delta'')} G_m/\langle \beta_m(M)\rangle \ar@{_{ (}->}[d] \\
\Gm\!\left(\A^1_{\Orond_S} \setminus \widetilde{\Mrond'}\right)\!/2 \ar[r] &
\displaystyle\prod_{M\in\Mrond''}\kappa(M)^\star/\langle\kappa(M)^{\star 2},\beta_m(M)\rangle
}
$$
et
$$
\myxyhook\xymatrix{
\Gm\!\left(\Htilde \setminus \left(\widetilde{\Delta''\cup\Pi}\right)\right)\!/2 \ar[r] \ar[d]^\wr &
\displaystyle\prod_{m\in \Irr(\Delta')} G_m/\langle \beta_m(M),\gamma_m(M)\rangle \ar@{_{ (}->}[d] \\
\Gm\!\left(\A^1_{\Orond_S} \setminus \widetilde{\Mrond''}\right)\!/2 \ar[r] &
\displaystyle\prod_{M\in\Mrond'}\kappa(M)^\star/\langle\kappa(M)^{\star 2},\beta_m(M),\gamma_m(M)\rangle\rlap{\text{,}}
}
$$
dont les flèches verticales de gauche sont les isomorphismes donnés par le lemme~\ref{ch3specdrestiso}.
Combinant ces deux carrés, le lemme~\ref{ch3specdgothspec} et le lemme~\ref{ch3specdconde}, on voit
maintenant que si la condition~\cDg{} est satisfaite, alors la \condE{} relative à
l'ensemble de places associé au point~$\tau(h_0)$ l'est aussi, ce qui achève la preuve de la proposition~\ref{ch3specdcamarche}.
\end{demo}

\bigskip
Concluons ce paragraphe avec la remarque suivante.

\bigskip
\begin{proposition}
\label{ch3specdtoujoursnoyau}
Le noyau de la flèche apparaissant dans la définition de la condition~\cDgp{} (resp.~\cDgpp{})
contient toujours~$\tgoth'$ (resp.~$\mgoth''$ et~$\tgoth''$).
\end{proposition}

\bigskip
\begin{demo}
Soient $m \in \Irr(\Delta'')$ et~$h$ le point générique de~$m$.
Nous avons déjà vu (cf.~preuve de la proposition~\ref{ch3specdtptppmpp}) que~$\tgoth'$ appartient
au groupe de $\phi'$\nobreakdash-Selmer géométrique
de~$E'$ relativement au trait $\Spec(\Orond_{H,h})$.
La courbe elliptique~$E'$ ayant réduction de type~$I_1$ ou~$I_2$ en~$h$
(cf.~proposition~\ref{ch3specdi1i2}), le diagramme commutatif~(\ref{ch2compdiag2})
montre que l'image de~$\tgoth'$ dans $\kappa(m)^\star/\kappa(m)^{\star 2}$ est nulle.
Soient maintenant $m \in \Irr(\Delta')$ et~$h$ le point générique de~$m$.
Les mêmes arguments que ceux que l'on vient d'employer prouvent d'une part que l'image de~$\tgoth''$
dans $\kappa(m)^\star/\langle \kappa(m)^{\star 2}, \gamma_m \rangle$ est nulle,
compte tenu que~$\gamma_m$ est la classe de l'extension quadratique ou triviale
minimale de~$\kappa(h)$ sur laquelle le groupe des composantes connexes de la fibre
spéciale du modèle de Néron de~$E''$ au-dessus de $\Spec(\Orond_{H,h})$ devienne constant
(cf.~corollaire~\ref{anncorsingrat}), et d'autre part que l'image de~$\mgoth''$
dans $\kappa(m)^\star/\kappa(m)^{\star 2}$ est égale à~$\beta_m$
(le raisonnement est exactement le même que celui qui précède le diagramme~(\ref{ch3specdcarre1})).
\end{demo}

\subsection{Vérification de la \condD{} générique}
\label{ch3ssverification}

Les hypothèses du paragraphe précédent sont toujours en vigueur.  Supposons de plus que~$k$
soit un corps de nombres.
Nous allons prouver que la condition~\cDg{} est satisfaite
sous l'hypothèse~(v) du théorème~\ref{ch3thprindp4}; la démonstration de la proposition ci-dessous occupera tout le paragraphe~\ref{ch3ssverification}.

\bigskip
\begin{proposition}
\label{ch3verifprin}
Supposons que la condition~(\ref{ch3enonceshyp}) soit satisfaite, que $\Br(X)/\Br(k)=0$ et
que soit $\epsilon_0=1$ dans $k^\star/k^{\star 2}$, soit il existe $t \in \Srond'$ tel que
l'image de~$\epsilon_0$ dans $\kappa(t)^\star/\kappa(t)^{\star 2}$ soit distincte de~$1$
et de~$\epsilon_t$.  Alors\index{condition (D)@\condD{}!générique} la condition~\cDg{} est satisfaite.
\end{proposition}

\bigskip
Il est possible d'étudier certaines questions liées à la structure galoisienne de~$X$
en appliquant sur~$\ksep$ les calculs explicites du paragraphe~\ref{ch3calcexpl} pour le cas simultanément diagonal
et en suivant l'action de~$\Gamma$ sur les coefficients des polynômes ainsi obtenus.
Ce type de raisonnement sera utilisé à plusieurs reprises dans la preuve de la proposition~\ref{ch3verifprin}.  Nous le détaillons
ici une fois pour toutes.

Pour $i \in \{\uplet{0}{4}\}$, soit $e_i \in \ksep^5$ un vecteur non nul dont l'image dans~$\P^4(\ksep)$ soit égale à~$\fp{P}i$.
Quitte à multiplier les~$e_i$ par des scalaires non nuls, on peut supposer le vecteur~$e_0$ invariant sous~$\Gamma$
et la famille $(\uplet{e_1}{e_4})$ globalement stable sous~$\Gamma$.
D'après la proposition~\ref{ch3gensimdiag}, la famille $e=(\uplet{e_0}{e_4})$ est une base de~$\ksep^5$
et les formes quadratiques~$q_1$ et~$q_2$ sont simultanément diagonales dans~$e$.  Reprenons donc les notations
du paragraphe~\ref{ch3calcexpl}, cas simultanément diagonal, associées à~$e$ et au corps~$\ksep$;
d'où notamment des éléments $\uplet{a_0}{a_4},\uplet{b_0}{b_4}\in\ksep$.
Notons $\chi \colon \Gamma \rightarrow \mathfrak{S}_4$ le morphisme de groupes
tel que $\gamma e_i = e_{\chi(\gamma)(i)}$ pour tout $\gamma\in \Gamma$ et
tout $i \in \{\uplet{1}{4}\}$.
La base~$e$ détermine un isomorphisme de $\ksep$\nobreakdash-algèbres
\begin{equation}
\label{ch3verifiso}
\Gamma\!\left(H \otimes_k \ksep, \;\bigoplus_{d \geq 0} \Orond(d)\!\right) \verylongisoto \ksep[\uplet{y_1}{y_4}] \rlap{\text{,}}
\end{equation}
d'où une action du groupe~$\Gamma$ sur la $k$\nobreakdash-algèbre $\ksep[\uplet{y_1}{y_4}]$, par transport de structure.
L'action de~$\Gamma$ sur les monômes unitaires est donnée par~$\chi$; plus précisément,
on a $\gamma y_i = y_{\chi(\gamma)(i)}$ pour tout $\gamma \in \Gamma$ et tout $i\in\{\uplet{1}{4}\}$.
Par ailleurs, les polynômes $q_1(y)=a_1y_1^2 + \cdots + a_4y_4^2$ et $q_2(y)=b_1y_1^2 + \cdots + b_4y_4^2$ sont invariants sous~$\Gamma$
puisque $q_1,q_2 \in \Gamma(H,\Orond(2))$; il en résulte que $\gamma a_i=a_{\chi(\gamma)(i)}$ et $\gamma b_i=b_{\chi(\gamma)(i)}$
pour tout $\gamma \in \Gamma$ et tout $i \in \{\uplet{1}{4}\}$.  Enfin, on a $a_0,b_0\in k$ puisque~$e_0$
est invariant sous~$\Gamma$.
Ainsi a-t-on prouvé:

\bigskip
\begin{lemme}
\label{ch3veriflemmeinvariant}
Soit $f \in \Z[\uplet{a_0}{a_4},\uplet{b_0}{b_4},\uplet{y_1}{y_4}]$
un polynôme invariant par l'action
de $\chi(\Gamma) \subset \mathfrak{S}_4$ sur les indices $\{\uplet{1}{4}\}$,
où les~$a_i$ et les~$b_i$ désignent ici des indéterminées.
L'image de~$f$ dans $\ksep[\uplet{y_1}{y_4}]$ par le morphisme d'anneaux
qui envoie les indéterminées $\uplet{a_0}{a_4},\uplet{b_0}{b_4}$ sur
les éléments de~$\ksep$ précédemment notés $\uplet{a_0}{a_4},\uplet{b_0}{b_4}$
est invariante sous~$\Gamma$ et définit donc un élément de
$\Gamma\!\left(H, \;\bigoplus_{d \geq 0} \Orond(d)\!\right)$ \emph{via}
l'isomorphisme~(\ref{ch3verifiso}).
\end{lemme}

\bigskip
\begin{demo}[ de la proposition~\ref{ch3verifprin}]%
Nous allons d'abord récrire la condition~\cDg{} en tenant compte des résultats obtenus
aux paragraphes~\ref{ch3parconssd} et~\ref{ch3calcexpl}.
Posons $\Delta_0 = Q_{t_0} \cap H^0$\glossary{$\Delta_0$, $\Delta_{1234}$, $\Delta_6$}
et $\Delta_{1234}=\bigcup_{t \in \Srond'}Q_t \cap H^0$.
Soit $\Delta_6 \subset H^0$ l'hypersurface intègre contenant l'unique point de codimension~$1$
de~$H^0$ au-dessus duquel la fibre de~$\pi_{H^0}$ est singulière et géométriquement
intègre (cf.~théorème~\ref{ch3conssdtheoreme}). Les variétés~$\Delta_0$ et~$\Delta_6$ sont
géométriquement intègres (pour~$\Delta_6$, cela découle par exemple de l'assertion d'unicité
dans le théorème~\ref{ch3conssdtheoreme} (i), que l'on applique cette fois après extension
des scalaires de~$k$ à~$\ksep$); la variété $\Delta_{1234}\otimes_k \ksep$ possède quatre
composantes irréductibles, sur lesquelles~$\Gamma$ agit comme sur $\{\uplet{t_1}{t_4}\}$.
Dans une base de~$\ksep^5$ dans laquelle~$q_1$ et~$q_2$ sont simultanément diagonales,
les hypersurfaces $\Delta_0 \otimes_k \ksep$ et $\Delta_{1234}\otimes_k \ksep$ de~$H^0$
sont respectivement définies par les équations~$p_0=0$ et $p_1p_2p_3p_4=0$, dans
la notation du paragraphe~\ref{ch3calcexpl}.
Les propositions~\ref{ch3expliciteabstth}, \ref{ch3expliciteegalitessim}
et~\ref{ch3conssdi1i2} montrent que la
fibre géométrique de~$\pi_{H^0}$ au-dessus du point générique de l'hypersurface d'équation $p_6=0$
(cf.~proposition~\ref{ch3expliciteirred})
est singulière et intègre; vu la définition de~$\Delta_6$, cela signifie que $p_6=0$ est l'équation
de l'hypersurface $\Delta_6 \otimes_k \ksep$.
On peut maintenant déduire des propositions~\ref{ch3explicitedelta} et~\ref{ch3expliciteirred}
que les composantes irréductibles de~$\Delta_{1234}$ ainsi que~$\Delta_0$ et~$\Delta_6$
sont des composantes irréductibles de~$\Delta$, que~$\Delta$ possède une unique autre composante
irréductible et que celle-ci est géométriquement irréductible.  Notons-la $\Delta_5$ et munissons-la
de sa structure de sous-schéma fermé réduit de~$\Delta$.

\bigskip
\begin{lemme}
\label{ch3verifidentdelta}
On a $\Delta'= \Delta_0 \cup \Delta_6$ et $\Delta''= \Delta_{1234} \cup \Delta_5$.
\end{lemme}

\bigskip
\begin{demo}
Pour démontrer ce lemme, on peut supposer les formes quadratiques~$q_1$ et~$q_2$ simultanément diagonalisables sur~$k$,
quitte à étendre les scalaires.  Soit~$h$ le point générique
d'une composante irréductible de~$\Delta$.
Reprenons les notations du paragraphe~\ref{ch3calcexpl};
comme dans la preuve de la proposition~\ref{ch3specdss}, on peut supposer que~$q_1$ et~$q_2$ sont simultanément
diagonales dans~$e$ et que $(y_1\gamma_4)(h)\neq 0$, de sorte que la courbe elliptique~$E'$ a pour équation
de Weierstrass $Y^2=(X-r(c))(X^2-r(d))$, où $r \colon k[\uplet{y_1}{y_4}]\rightarrow \Orond_{H,h}$
est l'application qui à $f(\uplet{y_1}{y_4})$ associe $f(1,y_2/y_1,y_3/y_1,y_4/y_1)$.
Nous avons vu, au cours de la preuve de la proposition~\ref{ch3specdss}, que cette équation est minimale en tant
qu'équation de Weierstrass à coefficients dans~$\Orond_{H,h}$.
Par conséquent, on a $h \in \Delta'$ si et seulement si~$r(d)$ n'est pas inversible,
si et seulement si $d(h)=0$.  La proposition~\ref{ch3expliciteegalitessim} montre
que $d=4y_1^4\gamma_4^4 p_0^2 p_6$, d'où le lemme.
\end{demo}

\bigskip
Pour $t \in \Srond'$, notons $\Delta_t=Q_t \cap H^0$\glossary{$\Delta_5$, $\Delta_t$}.  C'est une $k$\nobreakdash-variété intègre. La fermeture algébrique
de~$k$ dans~$\kappa(\Delta_t)$ s'identifie naturellement à~$\kappa(t)$.  Il est donc légitime de considérer~$\epsilon_t$
comme un élément de $\kappa(\Delta_t)^\star/\kappa(\Delta_t)^{\star 2}$ pour $t \in \Srond'$; de même,
pour $i \in \{0,5,6\}$, comme~$\Delta_i$ est géométriquement intègre sur~$k$, on peut identifier
le groupe $k^\star/k^{\star 2}$ à son image dans $\kappa(\Delta_i)^\star/\kappa(\Delta_i)^{\star 2}$.

\bigskip
\begin{lemme}
\label{ch3verifidentbeta}
Les~$\beta_m$ pour $m \in \Irr(\Delta)$ sont donnés par les
égalités suivantes:
\begin{align*}
\beta_{\Delta_0}&=\epsilon_0\text{;} &
\beta_{\Delta_t}&=\epsilon_t \text{ pour tout } t \in \Srond'\text{;} &
\beta_{\Delta_5}&=\epsilon_0\text{;} &
\beta_{\Delta_6}&=1\text{.}
\end{align*}
Par ailleurs, on a $\gamma_{\Delta_6}=1$.
\end{lemme}

\bigskip
\begin{demo}
Il est clair que $\beta_{\Delta_6}=1$ et $\gamma_{\Delta_6}=1$ puisque la fibre
de~$\pi_{H^0}$ au-dessus du point générique de~$\Delta_6$ est de type~$I_1$
(cf.~proposition~\ref{ch3specdi1i2}, compte tenu que cette fibre est géométriquement irréductible
par définition de~$\Delta_6$).  Les autres égalités résultent des assertions~(ii) et~(iii)
du théorème~\ref{ch3conssdtheoreme}.
\end{demo}

\bigskip
Afin de simplifier les notations, posons $\gamma_0=\gamma_{\Delta_0}$\glossary{$\gamma_0$}.
En combinant les lemmes~\ref{ch3verifidentdelta} et~\ref{ch3verifidentbeta}, nous obtenons:

\bigskip
\begin{lemme}
\label{ch3verifrecritd}
La condition~\cDgp{} est équivalente à ce que le noyau de la flèche naturelle
$$
\Gm(H \setminus (\Delta_0 \cup \Delta_6 \cup \Pi))/2 \xrightarrow{\,\,\,\theta'}
\kappa(\Delta_5)^\star/\langle \kappa(\Delta_5)^{\star 2}, \epsilon_0\rangle \times
\prod_{t\in \Srond'} \kappa(\Delta_t)^\star/\langle \kappa(\Delta_t)^{\star 2}, \epsilon_t\rangle
$$
soit engendré par~$\tgoth'$.
La condition~\cDgpp{} est équivalente à ce que le noyau de la flèche naturelle
$$
\Gm(H \setminus (\Delta_{1234} \cup \Delta_5 \cup \Pi))/2 \xrightarrow{\,\,\,\,\,\theta''\,\,}
\kappa(\Delta_0)^\star/\langle \kappa(\Delta_0)^{\star 2}, \epsilon_0, \gamma_0 \rangle
\times \kappa(\Delta_6)^\star/\kappa(\Delta_6)^{\star 2}
$$
soit engendré par~$\mgoth''$ et~$\tgoth''$.
\end{lemme}

\bigskip
Pour $i \in \{0,5,6\}$ ou $i\in\Srond'$, fixons une fonction rationnelle $f_i \in \kappa(H)^\star$ dont le diviseur
soit égal à $[\Delta_i]-\deg(\Delta_i)[\Pi]$ et posons $f_{1234}=\prod_{t \in \Srond'}f_t$.
Le lemme suivant précise la proposition~\ref{ch3specdtptppmpp}.

\bigskip
\begin{lemme}
\label{ch3verifgoth}
La classe~$\tgoth'$ (resp.~$\tgoth''$, $\mgoth''$)
est le produit d'un élément de $k^\star/k^{\star 2}$ et de la classe de~$f_6$
(resp.~$f_{1234}f_5$, $f_{1234}$).
\end{lemme}

\bigskip
\begin{demo}
On peut supposer les formes quadratiques~$q_1$ et~$q_2$ simultanément diagonalisables sur~$k$, quitte à étendre
les scalaires, et raisonner comme dans la
preuve du lemme~\ref{ch3verifidentdelta}, dont on reprend les notations.
La classe $\tgoth'$ (resp.~$\tgoth''$, $\mgoth''$) est égale à l'image
dans $\kappa(H)^\star/\kappa(H)^{\star 2}$ de $r(d)$
(resp.~$r(c^2-d)$, $r(\dgras_{14}^2-\dgras_{12}\dgras_{13})$)
(cf.~proposition~\ref{ch3explicitegeneralites} pour l'expression de~$\mgoth''$),
qui à son tour est égale à l'image
de $r(p_6)$ (resp.~$r(-p_1p_2p_3p_4p_5)$, $r(p_1p_2p_3p_4)$),
comme le montre la proposition~\ref{ch3expliciteegalitessim}. Le résultat voulu s'ensuit.
\end{demo}

\bigskip
\begin{lemme}
\label{ch3verifdpaf0}
Aucune classe de la forme $\alpha f_0$ pour $\alpha \in k^{\star}/k^{\star 2}$
n'appartient au noyau de~$\theta'$.
\end{lemme}

\bigskip
\begin{demo}
On peut supposer~$k$ algébriquement clos (au sens où d'une part l'hypothèse que~$k$ est un corps de nombres ne va jouer ici
aucun rôle et d'autre part, le lemme pour~$k$ algébriquement clos implique
le lemme pour~$k$ arbitraire; la réciproque n'est pas claire \emph{a priori}).
Sous cette hypothèse, l'ensemble~$\Srond'$ s'identifie à $\{t_1,t_2,t_3,t_4\}$ et l'on
a $\epsilon_t=1$ pour tout $t \in \Srond'$, de sorte qu'il suffit de prouver par exemple que
la fonction de~$\kappa(\Delta_1)^\star$ induite par~$f_0$ n'est pas un carré, où l'on
a posé $\Delta_1=\Delta_{t_1}$.
Vu les définitions de~$\Delta_0$ et de~$\Delta_1$, on a $\Delta_0 \cap \Delta_1=X \cap H^0$.
Cette variété est une courbe lisse et connexe (cf.~proposition~\ref{ch3conssdxhlisse}); notons~$\xi$
son point générique.  Étant donné que $\xi \not\in \Pi$ (cf.~\cite[Lemma~1.3 (i)]{ctsansdi}),
les fonctions rationnelles~$f_0$ et~$f_1$ sont respectivement des équations locales pour~$\Delta_0$
et~$\Delta_1$ au voisinage de~$\xi$.  On a donc
$\Orond_{\Delta_i,\xi}=\Orond_{H,\xi}/(f_i)$ pour $i \in \{0,1\}$,
d'où $\Orond_{\Delta_1,\xi}/(f_0)=\Orond_{H,\xi}/(f_0,f_1)=\Orond_{\Delta_0\cap\Delta_1,\xi}$.
Comme $\Delta_0\cap\Delta_1$ est intègre, ceci prouve que $\Orond_{\Delta_1,\xi}/(f_0)$ est un corps.
Autrement dit, la fonction de~$\kappa(\Delta_1)^\star$ induite par~$f_0$ s'annule à l'ordre~$1$
sur le diviseur $\Delta_0 \cap \Delta_1$; ce n'est donc pas un carré.
\end{demo}

\bigskip
\begin{lemme}
\label{ch3verifdpk}
Le noyau de~$\theta'$ ne rencontre pas $(k^\star/k^{\star 2}) \setminus \{1\}$.
\end{lemme}

\bigskip
\begin{demo}
La flèche naturelle
$k^\star/k^{\star 2} \rightarrow \kappa(\Delta_5)^\star/\kappa(\Delta_5)^{\star 2}$ est injective
puisque $\Delta_5$ est géométriquement intègre sur~$k$.
Vu la définition de~$\theta'$, ceci prouve déjà que
les seules classes de $k^\star/k^{\star 2}$ susceptibles d'appartenir au noyau de~$\theta'$
sont la classe triviale et~$\epsilon_0$.  Supposons donc que $\theta'(\epsilon_0)=1$.
Considérant les autres facteurs
du but de~$\theta'$, on voit que pour tout $t \in \Srond'$, l'image de~$\epsilon_0$ dans
$\kappa(\Delta_t)^\star/\kappa(\Delta_t)^{\star 2}$ doit appartenir au sous-groupe $\{1,\epsilon_t\}$.
L'hypothèse de la proposition~\ref{ch3verifprin} permet d'en déduire que $\epsilon_0=1$, ce qui
termine de prouver le lemme.
\end{demo}

\bigskip
Étant donné que le groupe $\Gm(H\setminus(\Delta_0 \cup \Delta_6 \cup \Pi))/2$
est engendré par les classes de~$f_0$ et de~$f_6$ et par le sous-groupe $k^\star/k^{\star 2}$,
les lemmes~\ref{ch3verifrecritd},
\ref{ch3verifgoth}, \ref{ch3verifdpaf0} et~\ref{ch3verifdpk} et la proposition~\ref{ch3specdtoujoursnoyau} prouvent ensemble que
la condition~\cDgp{} est satisfaite.

\bigskip
\begin{lemme}
\label{ch3verifdppk}
Le noyau de~$\theta''$ ne rencontre pas $(k^\star/k^{\star 2})\setminus\{1\}$.
\end{lemme}

\bigskip
\begin{demo}
La flèche naturelle $k^\star/k^{\star 2}\rightarrow \kappa(\Delta_6)^\star/\kappa(\Delta_6)^{\star 2}$ est
en effet injective, la variété~$\Delta_6$ étant géométriquement intègre sur~$k$.
\end{demo}

\bigskip
Si~$J$ est une partie de~$\Srond'$, on appellera \emph{longueur de~$J$} le degré sur~$k$ du sous-schéma réduit de~$\Srond'$
associé et l'on notera $f_J=\prod_{t \in J}f_t$.

\bigskip
\begin{lemme}
\label{ch3verifdpplong1}
Aucune classe de la forme $\alpha f_J$ pour $\alpha \in k^\star/k^{\star 2}$
et $J \subset \Srond'$ de longueur~$1$ n'appartient au noyau de~$\theta''$.
\end{lemme}

\bigskip
\begin{demo}
Supposons qu'il existe $J \subset \Srond'$ de longueur~$1$ et $\alpha \in k^\star/k^{\star 2}$
tels que $\alpha f_J \in \Ker(\theta'')$.  L'ensemble~$J$ est un singleton et son unique élément
est $k$\nobreakdash-rationnel; quitte à renuméroter les~$t_i$, on peut donc supposer que $t_1 \in \P^1(k)$
et que $J=\{t_1\}$.  Notons alors $f_1=f_{t_1}=f_J$, $\Delta_1=\Delta_{t_1}$ et $\epsilon_1=\epsilon_{t_1}$.
Fixons des vecteurs $e_0,e_1\in k^5\setminus\{0\}$ dont les images dans~$\P^4(k)$ soient
respectivement égales à~$\fp{P}0$ et~$\fp{P}1$ et posons $a_i=q_1(e_i)$ et $b_i=q_2(e_i)$ pour $i\in\{0,1\}$
et $d_{01}=a_0b_1-a_1b_0$.

La seule condition imposée jusqu'à présent sur la fonction rationnelle~$f_1$ est que son diviseur soit égal
à $[\Delta_1]-\deg(\Delta_1)[\Pi]$.  Nous allons maintenant la normaliser.
Soit $\ell \in \Gamma(H,\Orond(1)) \setminus \{0\}$ une forme linéaire s'annulant sur~$\Pi$.
La forme quadratique $a_1 q_2 - b_1 q_1$ définit, par restriction, un
élément de $\Gamma(H, \Orond(2))$, que l'on note~$p_1$.
La fonction rationnelle $p_1/\ell^2 \in \kappa(H)$ a pour diviseur $[\Delta_1]-\deg(\Delta_1)[\Pi]$.
Quitte à multiplier~$f_1$ par une constante, on peut donc supposer que $f_1=p_1/\ell^2$.

\bigskip
\begin{souslemme}
\label{ch3verifdppf1d6}
La classe dans $\kappa(\Delta_6)^\star/\kappa(\Delta_6)^{\star 2}$ de la fonction rationnelle
induite par~$f_1$ est égale à la classe de $d_{01}\epsilon_1$.
\end{souslemme}

\bigskip
\begin{demo}
Choisissons des vecteurs $e_2,e_3,e_4\in\ksep^5$ tels que la famille
$e=(\uplet{e_0}{e_4})$ soit une base de~$\ksep^5$ stable sous~$\Gamma$ dans laquelle
les formes quadratiques~$q_1$ et~$q_2$ sont simultanément diagonales
et reprenons les notations introduites immédiatement avant le lemme~\ref{ch3veriflemmeinvariant}.
Remarquons que les notations~$a_0$, $a_1$, $b_0$, $b_1$, $d_{01}$ et~$p_1$ définies ci-dessus
coïncident avec celles du paragraphe~\ref{ch3calcexpl}.
Définissons des éléments $\rho\in\ksep^\star$, $A \in \Gamma(H\otimes_k\ksep,\Orond(3))$,
$B \in \Gamma(H\otimes_k\ksep,\Orond(6))$
et $C \in \Gamma(H\otimes_k\ksep, \Orond(4))$ par les formules:
\begin{equation}
\label{ch3verifformulesrhoabc}
\begin{aligned}
\rho &= (d_{34}d_{23}d_{24})^2d_{21}d_{31}d_{41} \rlap{\text{,}} \\
A &= 4y_2y_3y_4 \rlap{\text{,}} \\
B &= 2c_6(y_1^4)y_1^2 + c_6(y_1^2) \rlap{\text{,}}\\
C &= -4c_6(y_1^4) \rlap{\text{.}}
\end{aligned}
\end{equation}
(Rappelons que $c_6(y_i^m)$ désigne
le coefficient de~$y_i^m$ dans~$p_6$, où~$p_6$ est considéré comme un polynôme en~$y_i$
à coefficients dans $k[(y_\ell)_{\ell\neq i}]$.)
On vérifie tout de suite qu'en tant qu'éléments de $\Z[\uplet{a_0}{a_4},\uplet{b_0}{b_4},\uplet{y_1}{y_4}]$,
les membres de droite des égalités ci-dessus et de l'égalité définissant~$p_6$ sont
invariants par toute permutation des indices $\{2,3,4\}$.
Compte tenu que le vecteur~$e_1$ est invariant sous~$\Gamma$, le lemme~\ref{ch3veriflemmeinvariant}
permet d'en déduire que~$\rho$, $A$, $B$, $C$ et~$p_6$ appartiennent respectivement
à $k$, $\Gamma(H,\Orond(3))$, $\Gamma(H,\Orond(6))$, $\Gamma(H,\Orond(4))$
et $\Gamma(H,\Orond(8))$.
D'après la proposition~\ref{ch3explicitepropp6}, l'égalité
$$
(A p_1)^2 \rho p_1 = B^2 + Cp_6
$$
d'éléments de $\Gamma(H,\Orond(12))$ a lieu. En restreignant cette égalité à $\Gamma(\Delta_6,\Orond(12))$ et
en la divisant par~$\ell^{12}$, on voit que la classe dans $\kappa(\Delta_6)^\star/\kappa(\Delta_6)^{\star 2}$
de la fonction rationnelle induite par~$f_1$ est égale à la classe de~$\rho$ dans $k^\star/k^{\star 2}$.

Il reste seulement à vérifier que $d_{01}\epsilon_1\rho$ est un carré dans~$k$.
Pour cela, on peut s'autoriser à remplacer~$k$ par une extension finie de degré impair.
En particulier, on peut supposer que~$\Gamma$ n'agit pas transitivement sur $\{e_2,e_3,e_4\}$;
quitte à permuter ces trois vecteurs, on peut donc supposer~$e_2$ invariant sous~$\Gamma$.
Pour obtenir une expression de~$\epsilon_1$ en fonction des~$d_{ij}$, il faut trouver un hyperplan
$k$\nobreakdash-rationnel de~$H$ ne contenant pas~$\fp{P}1$, exprimer une base $k$\nobreakdash-rationnelle de l'espace vectoriel associé
puis calculer le déterminant de la matrice de $a_1q_2-b_1q_1$ dans cette base.
Posons donc $f_0=e_0$, $f_1=e_2$, $f_2=e_3+e_4$ et $f_3=d_{34}(e_3-e_4)$.
Les vecteurs~$f_i$ sont invariants sous~$\Gamma$ puisqu'ils sont invariants par permutation
des indices~$3$ et~$4$.  La famille $f=(\uplet{f_0}{f_3})$ est clairement libre, et la matrice
dans~$f$ de la restriction au sous-espace engendré par~$f$ de la forme quadratique $a_1q_2-b_1q_1$
a pour déterminant $4d_{10}d_{12}d_{34}^2d_{13}d_{14}$.  Le quotient de cet élément de~$k$
par $d_{01}\rho$ est égal à $4(d_{23}d_{24})^{-2}$: c'est bien le carré d'un élément de~$k$
puisque $d_{23}d_{24}$ est invariant par permutation des indices~$3$ et~$4$.
\end{demo}

\bigskip
\begin{souslemme}
\label{ch3verifdppf1d0}
La classe dans $\kappa(\Delta_0)^\star/\kappa(\Delta_0)^{\star 2}$ de la fonction rationnelle
induite par~$f_1$ est égale à la classe de $d_{01}\epsilon_0\gamma_0$.
\end{souslemme}

\bigskip
\begin{demo}
Pour prouver ce sous-lemme, on peut remplacer~$k$ par une extension finie de degré impair,
ce qui permet de supposer que l'un de $t_2$, $t_3$ et~$t_4$ est $k$\nobreakdash-rationnel.
Quitte à renuméroter
les~$t_i$, on peut supposer que $t_2\in\P^1(k)$.
Une transformation linéaire de~$\P^4_k$ permet de supposer de plus les formes quadratiques~$q_1$ et~$q_2$
presque simultanément diagonales dans la base canonique~$e$ de~$k^5$.
Reprenons alors les notations du
paragraphe~\ref{ch3calcexpl}, cas presque simultanément diagonal,
en prenant pour~$h$ le point générique de~$\Delta_0$, de sorte que $(y_1\gamma_4)(h)\neq 0$.
Soit $r \colon k[\uplet{y_1}{y_4}]\rightarrow \Orond_{H,h}$ le morphisme d'anneaux
qui à $f(\uplet{y_1}{y_4})$ associe $f(1,y_2/y_1,y_3/y_1,y_4/y_1)$.
Comme on l'a vu au cours de la preuve de la proposition~\ref{ch3specdss},
l'équation \mbox{$Y^2=(X-r(c))(X^2-r(d))$} est une équation de Weierstrass minimale
à coefficients dans l'anneau de valuation discrète $\Orond_{H,h}$ pour
la courbe elliptique~$E'$.  On peut donc déterminer la classe~$\gamma_0$
en lisant sur cette équation la plus petite
extension de~$\kappa(h)$ sur laquelle les pentes des tangentes de la cubique
réduite en son point singulier deviennent rationnelles
(cf.~lemmes~\ref{annlemmeextreme} et~\ref{annlemmetoresing}).
On voit ainsi que~$\gamma_0$ coïncide avec la classe
de~$r(-c)$ dans $\kappa(h)^\star/\kappa(h)^{\star 2}$,
compte tenu que~$r(d)$ s'annule en~$h$ (cf.~proposition~\ref{ch3expliciteabstth}).
Il suffit donc pour conclure
de prouver que la classe de~$f_1$ est égale à celle de $-d_{01}\epsilon_0 r(c)$.
La proposition~\ref{ch3explicitepsdeq} et la remarque qui la suit
montrent que
$$
-d_{01}\epsilon_0 r(c) = \left(2 \frac{\epsilon_0}{d_{01}}r(y_1\gamma_4p_1)\right)^2 r(p_1) \mod \mgoth_{H,h} \rlap{\text{.}}
$$
Comme $r(p_1)=(r(\ell))^2 f_1$, ceci achève la démonstration.
\end{demo}

\bigskip
\begin{souslemme}
\label{ch3verifdppgamma0pasconst}
La classe $\gamma_0 \in \kappa(\Delta_0)^\star/\kappa(\Delta_0)^{\star 2}$ n'appartient
pas au sous-groupe $k^\star/k^{\star 2}$.
\end{souslemme}

\bigskip
\begin{demo}
Compte tenu du sous-lemme~\ref{ch3verifdppf1d0}, il suffit
de prouver que la fonction de $\kappa(\Delta_0 \otimes_k \ksep)^\star$ induite par~$f_1$ n'est
pas un carré.  L'argument est exactement le même que celui employé dans
la preuve du lemme~\ref{ch3verifdpaf0}.
\end{demo}

\bigskip
L'injectivité de la flèche naturelle $k^\star/k^{\star 2}\rightarrow \kappa(\Delta_6)^\star/\kappa(\Delta_6)^{\star 2}$,
le sous-lemme~\ref{ch3verifdppf1d6}
et l'hypothèse que $\theta''(\alpha f_1)=0$
entraînent que $\alpha=d_{01}\epsilon_1$, d'où l'on déduit, grâce au sous-lemme~\ref{ch3verifdppf1d0},
que la classe $\epsilon_0\epsilon_1\gamma_0 \in \kappa(\Delta_0)^\star/\kappa(\Delta_0)^{\star 2}$
appartient au sous-groupe engendré par~$\epsilon_0$ et~$\gamma_0$.
Vu le sous-lemme~\ref{ch3verifdppgamma0pasconst}, il en résulte que $\epsilon_1 \in \{1,\epsilon_0\}$.
L'hypothèse~(\ref{ch3enonceshyp}) assure que $\epsilon_1\neq 1$; on a donc $\epsilon_0\epsilon_1=1$
et $\epsilon_0 \neq 1$.
Notons~$d'$ la dimension du sous-$\Z/2$\nobreakdash-espace vectoriel de~$k^\star/k^{\star 2}$
engendré par les normes $N_{\kappa(t)/k}(\epsilon_t)$ pour $t\in\Srond'$.
L'égalité $\epsilon_0\epsilon_1=1$ et la proposition~\ref{ch3brproduit}
fournissent la relation $\prod_{t \in \Srond'\setminus\{t_1\}}\epsilon_t=1$,
d'où résulte l'inégalité $d'<\Card(\Srond')$.
Soient~$n$ et~$d$ les entiers définis dans l'énoncé du théorème~\ref{ch3thbr}.
Comme~$\Srond'$ contient un point $k$\nobreakdash-rationnel,
l'hypothèse~(\ref{ch3enonceshyp})
entraîne que $\epsilon_t \neq 1$ pour tout $t\in\Srond'$.
D'autre part, on a vu que $\epsilon_0 \neq 1$. On
a donc $n=\Card(\Srond)=\Card(\Srond')+1$.
La proposition~\ref{ch3brproduit} montre par ailleurs que $d=d'$;
d'où finalement $d<n-1$.  Cette inégalité contredit la nullité de $\Br(X)/\Br(k)$ d'après le théorème~\ref{ch3thbr}.
\end{demo}

\bigskip
\begin{lemme}
\label{ch3verifdpplong2}
Aucune classe de la forme $\alpha f_J$ pour $\alpha \in k^\star/k^{\star 2}$
et $J \subset \Srond'$ de longueur~$2$ n'appartient au noyau de~$\theta''$.
\end{lemme}

\bigskip
\begin{demo}
La preuve de ce lemme est similaire à celle du lemme~\ref{ch3verifdpplong1}.
Supposons qu'il existe $J \subset \Srond'$ de longueur~$2$ et $\alpha \in k^\star/k^{\star 2}$
tels que $\alpha f_J \in \Ker(\theta'')$.  Quitte à renuméroter les~$t_i$, on peut supposer
que l'ensemble des points géométriques de~$J$ est $\{t_1,t_2\}$.
Fixons des vecteurs $\uplet{e_0}{e_4} \in \ksep^5 \setminus \{0\}$ dont les images dans~$\P^4(\ksep)$ soient
respectivement égales à $\uplet{\fp{P}0}{\fp{P}4}$, avec la condition supplémentaire que $e_0 \in k^5$ et
que les ensembles $\{e_1,e_2\}$ et $\{e_3,e_4\}$ soient stables sous~$\Gamma$.
Pour $i,j \in\{0,1,2\}$, posons $a_i=q_1(e_i)$, $b_i=q_2(e_i)$, $d_{ij}=a_ib_j-a_jb_i$.

Soit $\ell \in \Gamma(H,\Orond(1)) \setminus \{0\}$ une forme linéaire s'annulant sur~$\Pi$.
La forme homogène $(a_1 q_2-b_1 q_1)(a_2 q_2-b_2q_1)$ est invariante sous~$\Gamma$ et définit
donc un élément de $\Gamma(H,\Orond(4))$, que l'on note~$p_{12}$.
La fonction rationnelle $p_{12}/\ell^4 \in \kappa(H)$ a pour diviseur
$\sum_{t \in J}([\Delta_t] - \deg(\Delta_t)[\Pi])$. Quitte à multiplier l'un des~$f_t$
pour $t \in J$ par une constante, on peut donc supposer que $f_J=p_{12}/\ell^4$.

\bigskip
\begin{souslemme}
\label{ch3verifdppf12d6}
La classe dans $\kappa(\Delta_6)^\star/\kappa(\Delta_6)^{\star 2}$ de la fonction rationnelle
induite par~$f_J$ est égale à la classe de $d_{01}d_{02}\prod_{t\in J}N_{\kappa(t)/k}(\epsilon_t)$.
\end{souslemme}

\bigskip
(La formule $d_{01}d_{02}\prod_{t\in J}N_{\kappa(t)/k}(\epsilon_t)$ définit bien un élément de~$k^\star$
puisque $d_{01}d_{02}$ est invariant sous~$\Gamma$.  La même remarque vaudra pour le
prochain sous-lemme.)

\bigskip
\begin{demo}
Définissons~$\rho$, $A$, $B$ et~$C$ par les formules~(\ref{ch3verifformulesrhoabc})
et~$\rho'$, $A'$, $B'$ et~$C'$ par les formules obtenues à partir de~(\ref{ch3verifformulesrhoabc})
en échangeant les indices~$1$ et~$2$.

En tant qu'éléments de $\Z[\uplet{a_0}{a_4},\uplet{b_0}{b_4},\uplet{y_1}{y_4}]$,
les formules définissant les polynômes $AA'$, $\rho\rho'$, $BB'$, $CC'$, $CB'^2+BC'^2$ et~$p_6$ sont
invariantes si l'on permute les indices~$1$ et~$2$ ou les indices~$3$ et~$4$.
Compte tenu que~$\Gamma$ agit sur $\uplet{e_1}{e_4}$ à travers le sous-groupe engendré
par ces deux transpositions, le lemme~\ref{ch3veriflemmeinvariant}
permet d'en déduire que tous ces polynômes homogènes déterminent
des sections sur~$H$ des faisceaux~$\Orond(d)$ pour divers entiers~$d$.
D'après la proposition~\ref{ch3explicitepropp6}, on a alors l'égalité
$$
(AA'p_{12})^2 \rho\rho' p_{12} = (BB')^2 + (CB'^2+C'B^2+CC'p_6)p_6
$$
dans $\Gamma(H,\Orond(24))$.  (Considérer le produit des équations~(\ref{ch3expliciteeqp1p6})
et~(\ref{ch3expliciteeqp2p6}).)  En la restreignant à $\Gamma(\Delta_6,\Orond(24))$ et en la divisant
par~$\ell^{24}$, on voit que la classe dans $\kappa(\Delta_6)^\star/\kappa(\Delta_6)^{\star 2}$ de la fonction
rationnelle induite par~$f_J$ est égale à la classe de $\rho\rho'$ dans $k^\star/k^{\star 2}$.

Il reste seulement à vérifier que $\rho\rho'd_{01}d_{02}\prod_{t\in J}N_{\kappa(t)/k}(\epsilon_t)$ est un carré dans~$k$.
Si~$J$ comporte deux points rationnels, on procède exactement comme dans la preuve du sous-lemme~\ref{ch3verifdppf1d6}
pour calculer~$\epsilon_{t_1}$; en échangeant les indices~$1$ et~$2$, on obtient alors le calcul de la classe~$\epsilon_{t_2}$.
Lorsque~$J$ est un singleton, disons $J=\{t\}$, on identifie
l'extension quadratique $\kappa(t)/k$ à la sous-extension de~$\ksep/k$ engendrée par
les coordonnées de $t_1 \in \P^1(\ksep)$, de sorte que la famille~$f$ du sous-lemme~\ref{ch3verifdppf1d6}
est formée de vecteurs $\kappa(t)$\nobreakdash-rationnels et permet donc encore une fois de calculer
la classe $\epsilon_t \in \kappa(t)^\star/\kappa(t)^{\star 2}$
et de conclure comme pour le sous-lemme~\ref{ch3verifdppf1d6}.
\end{demo}

\bigskip
\begin{souslemme}
\label{ch3verifdppf12d0}
La classe dans $\kappa(\Delta_0)^\star/\kappa(\Delta_0)^{\star 2}$ de la fonction rationnelle
induite par~$f_J$ est égale à la classe de $d_{01}d_{02}$.
\end{souslemme}

\bigskip
\begin{demo}
Le lemme~\ref{ch3veriflemmeinvariant} permet d'interpréter l'équation~(\ref{ch3expliciteeqp1p2p0})
comme une égalité dans $\Gamma(H, \Orond(4))$, compte tenu
que $d_{01}d_{02}$, $d_{01}p_2+d_{02}p_1$ et $d_{12}^2$ sont invariants par permutation des indices~$1$ et~$2$
et par permutation des indices~$3$ et~$4$.  Restreignant cette égalité à $\Gamma(\Delta_0,\Orond(4))$ et
divisant ses deux membres par~$\ell^4$, on obtient le résultat voulu.
\end{demo}

\bigskip
L'injectivité de la flèche naturelle $k^\star/k^{\star 2}\rightarrow \kappa(\Delta_6)^\star/\kappa(\Delta_6)^{\star 2}$,
le sous-lemme~\ref{ch3verifdppf12d6} et l'hypothèse que $\theta''(\alpha f_J)=0$
entraînent que $$\alpha=d_{01}d_{02}\prod_{t\in J}N_{\kappa(t)/k}(\epsilon_t)\rlap{\text{,}}$$ d'où l'on déduit, grâce au sous-lemme~\ref{ch3verifdppf12d0},
que la classe $\prod_{t \in J}N_{\kappa(t)/k}(\epsilon_t)\in \kappa(\Delta_0)^\star/\kappa(\Delta_0)^{\star 2}$
appartient au sous-groupe engendré par~$\epsilon_0$ et~$\gamma_0$.
Vu le sous-lemme~\ref{ch3verifdppgamma0pasconst}, il en résulte que $\prod_{t \in J}N_{\kappa(t)/k}(\epsilon_t) \in \{1,\epsilon_0\}$,
ce que la proposition~\ref{ch3brproduit} permet de reformuler en disant qu'au moins l'une des deux
égalités
\begin{equation}
\label{ch3verifprodtjstj}
\prod_{t \in J}N_{\kappa(t)/k}(\epsilon_t)=1 \qquad \text{ et } \qquad
\prod_{t \in \Srond' \setminus J}N_{\kappa(t)/k}(\epsilon_t)=1
\end{equation}
a lieu.  On a donc $d'<\Card(\Srond')$,
en notant~$d'$ la dimension du sous-$\Z/2$\nobreakdash-espace vectoriel de~$k^\star/k^{\star 2}$
engendré par les normes $N_{\kappa(t)/k}(\epsilon_t)$ pour $t\in\Srond'$.
Soient~$n$ et~$d$ les entiers définis dans l'énoncé du théorème~\ref{ch3thbr}.
Comme~$\Srond'$ contient un point de degré~$\leq 2$,
l'hypothèse~(\ref{ch3enonceshyp})
entraîne que $\epsilon_t \neq 1$ pour tout $t\in\Srond'$,
d'où il résulte que $n=\Card(\Srond')$ si $\epsilon_0=1$ et $n=\Card(\Srond')+1$ sinon.
On a donc $d'<n-1$ si $\epsilon_0\neq 1$.
Par ailleurs, si $\epsilon_0=1$, la proposition~\ref{ch3brproduit} montre que les deux égalités~(\ref{ch3verifprodtjstj})
sont équivalentes; elles sont donc toutes les deux satisfaites, ce qui entraîne
que $d'<\Card(\Srond')-1$ et donc à nouveau $d'<n-1$.
Enfin, la proposition~\ref{ch3brproduit} montre que $d=d'$.  L'inégalité $d<n-1$ a donc toujours lieu,
ce qui contredit la nullité de $\Br(X)/\Br(k)$ d'après le théorème~\ref{ch3thbr}.
\end{demo}

\bigskip
Étant donné que le groupe $\Gm(H \setminus (\Delta_{1234} \cup \Delta_5 \cup \Pi))/2$ est engendré par les classes
des fonctions~$f_t$ pour $t \in \Srond'$ et de~$f_5$ et par le sous-groupe $k^\star/k^{\star 2}$,
les lemmes~\ref{ch3verifrecritd}, \ref{ch3verifgoth}, \ref{ch3verifdppk}, \ref{ch3verifdpplong1}, et~\ref{ch3verifdpplong2}
et la proposition~\ref{ch3specdtoujoursnoyau}
prouvent ensemble que
la condition~\cDgpp{} est satisfaite.
\end{demo}

\subsection{Preuve du théorème~\ref{ch3thprindp4}}
\label{ch3pardernierpreuve}

Nous admettons désormais
l'hypothèse de Schinzel et la finitude des groupes de Tate-Shafarevich des courbes elliptiques sur les corps de nombres.

\bigskip
Commençons par prouver que la surface~$X$ satisfait au principe de Hasse sous
l'hypothèse~(v) du théorème~\ref{ch3thprindp4}.
Comme le point~$t_0$ est supposé $k$\nobreakdash-rationnel,
on peut se servir des résultats du paragraphe~\ref{ch3parconssd}.
On dispose notamment
d'une variété~$C_{H^0}$ géométriquement intègre sur~$k$ (cf.~proposition~\ref{ch3conssdgeomint}),
d'une variété~$H^0$ isomorphe au complémentaire dans~$\P^3_k$ d'un fermé de codimension~$\geq 2$
et d'un morphisme projectif $\pi_{H^0} \colon C_{H^0} \rightarrow H^0$
dont la fibre générique est géométriquement intègre (cf.~proposition~\ref{ch3conssd2rat})
et dont la fibre au-dessus de tout point $m\in H^0$ de codimension~$1$ possède une composante
irréductible~$Y$ de multiplicité~$1$ dans le corps des fonctions de laquelle la fermeture algébrique de~$\kappa(m)$
est une extension abélienne de~$\kappa(m)$ (cf.~théorème~\ref{ch3conssdtheoreme}).
On peut donc appliquer le théorème~\ref{ch3thpointslocaux} à ce morphisme.

Supposons que $X(k_\Omega)\neq\emptyset$.
D'après la proposition~\ref{ch3conssdadmetsect},
il existe un ouvert dense $U \subset X$ et un $k$\nobreakdash-morphisme $f \colon U \rightarrow (C_{H^0})^\reg$,
où $(C_{H^0})^\reg$ désigne l'ouvert de lissité de~$C_{H^0}$ sur~$k$.
Étant donné que $\Br(X)/\Br(k)=0$, tout $k_\Omega$\nobreakdash-point de~$U$ est
orthogonal à~$\Br_\nr(U)$.
Comme $f^\star \Br_\nr((C_{H^0})^\reg) \subset \Br_\nr(U)$
(cf.~\cite[Théorème~7.4]{grothbr3}),
le lemme~\ref{ch3nishi}, appliqué au morphisme $f\otimes_k k_v$
pour chaque place $v \in \Omega$ (avec $S=Y=(C_{H^0})^\reg\otimes_k k_v$),
permet d'en déduire que
l'image par~$f$ de tout $k_\Omega$\nobreakdash-point de~$U$ est
orthogonale à $\Br_\nr((C_{H^0})^\reg)$.
Le~théorème des fonctions implicites assure par ailleurs l'existence d'un $k_\Omega$\nobreakdash-point de~$U$;
on a ainsi montré que $(C_{H^0})^\reg(k_\Omega)^{\Br_\nr}\neq\emptyset$.
D'après le théorème~\ref{ch3thpointslocaux} appliqué au morphisme~$\pi_{H^0}$
et à l'ouvert $H^0 \setminus \Delta$,
il existe donc $h_0 \in (H^0 \setminus \Delta)(k)$ tel que la fibre $\pi_{H^0}^{-1}(h_0)$ admette
un $k_v$\nobreakdash-point pour tout $v \in \Omega$.

Soit $h_1 \in (H^0 \setminus \Delta)(k)$ un point distinct de~$h_0$ mais suffisamment proche
de~$h_0$ aux places réelles de~$k$ pour que, notant~$D_0$ la droite de~$H$ passant par~$h_0$ et~$h_1$,
les points~$h_0$ et~$h_1$ appartiennent à la même composante
connexe de $(D_0 \setminus (D_0 \cap \Delta))(k_v)$ pour toute place~$v$ réelle.
Comme $\pi_{H^0}^{-1}(h_0)$ est lisse et possède un $k_v$\nobreakdash-point pour toute place~$v$ réelle
et que le morphisme~$\pi_{H^0}$ est plat (cf.~proposition~\ref{ch3conssdplatcm}),
le théorème des fonctions implicites permet de supposer, quitte à choisir~$h_1$ suffisamment proche
de~$h_0$, que $\pi_{H^0}^{-1}(h_1)$ possède un $k_v$\nobreakdash-point pour toute place~$v$ réelle.

Soit $\Pi \subset H$ un hyperplan $k$\nobreakdash-rationnel contenant~$h_1$ mais ne contenant pas~$h_0$.
Nous sommes maintenant dans la situation du paragraphe~\ref{ch3parspecd}, où une «~condition~\cDg~» a été
définie (cf.~définition~\ref{ch3defcdg}).  Celle-ci est satisfaite en vertu de la proposition~\ref{ch3verifprin}.
D'après la proposition~\ref{ch3specdcamarche}, il existe donc un point $p \in \Pi(k)$ arbitrairement
proche de~$h_1$ aux places réelles de~$k$ et un $k$\nobreakdash-isomorphisme $\tau \colon D \isoto \P^1_k$
vérifiant toutes les conditions énoncées immédiatement avant la proposition~\ref{ch3specdcamarche},
où~$D$ désigne la droite de~$H$ passant par~$h_0$ et~$p$.
Quitte à choisir~$p$ suffisamment proche de~$h_1$, on peut supposer que~$h_0$ et~$p$ appartiennent
à la même composante connexe de $(D \setminus (D \cap \Delta))(k_v)$ pour toute place~$v$ réelle
et que la fibre $\pi_{H^0}^{-1}(p)$ possède un $k_v$\nobreakdash-point pour toute place~$v$ réelle.
Il résulte alors de la condition~(v) précédant l'énoncé de la proposition~\ref{ch3specdcamarche}
que $\tau(h_0)$ appartient à la composante connexe non majorée
de $\tau(D \setminus (D \cap (\Delta \cup \Pi)))(k_v)$ pour toute place~$v$ réelle.

Toutes les hypothèses du théorème~\ref{ch2thsurmesure} sont satisfaites
pour la famille $\tau \circ \pi_{H^0} \colon \pi_{H^0}^{-1}(D) \rightarrow \P^1_k$
et le point $\tau(h_0)$.  Sa conclusion l'est donc aussi; il en résulte que la variété $C_{H^0}$
admet un point $k$\nobreakdash-rationnel.  L'existence d'un morphisme $C_{H^0} \rightarrow X$ permet de
conclure quant à l'existence d'un point $k$\nobreakdash-rationnel sur la surface~$X$.

\bigskip
Reste à prouver le principe de Hasse pour~$X$ sous chacune des hypothèses~(i) à~(iv)
du théorème~\ref{ch3thprindp4}.  Les deux propositions suivantes nous seront utiles pour éliminer
un certain nombre de cas particuliers.
Rappelons que la surface~$X$ admet un point rationnel
si elle admet un point sur une extension finie de~$k$ de degré impair; sous diverses formes plus générales,
ce résultat fut démontré par Amer, Brumer et Coray (cf.~\cite{brumer}, \cite{coray}).

\bigskip
\begin{proposition}
\label{ch3findp4amerbrumer}
Soit $k'/k$ une extension finie de degré impair.
Supposons que
la $k'$\nobreakdash-variété
$X \otimes_k k'$ admette un point rationnel si l'obstruction de Brauer-Manin
ne s'y oppose pas. Alors
la $k$\nobreakdash-variété~$X$ admet un point rationnel si l'obstruction de Brauer-Manin
ne s'y oppose pas.
\end{proposition}

\bigskip
\begin{demo}
Soit $(\fp{P}v)_{v \in \Omega} \in X(\A_k)^\Br$. Notons~$\Omega'$ l'ensemble des places de~$k'$
et $t\colon \Omega' \rightarrow \Omega$ l'application trace.
Pour $w \in \Omega'$, posons $\fp{P'}w=\fp{P}{t(w)}$. Le point adélique $(\fp{P'}w)_{w \in \Omega'}$
de la $k'$\nobreakdash-variété $X \otimes_k k'$
est orthogonal au groupe de Brauer, comme il résulte des égalités
$$
\sum_{w \in \Omega'} \inv_w A(\fp{P'}w) = \sum_{v \in \Omega} \sum_{w|v} \inv_w A(\fp{P}v)
= \sum_{v \in \Omega} \inv_v (\Cores_{k'/k}A)(\fp{P}v) \rlap{\text{,}}
$$
valables pour tout $A \in \Br(X\otimes_k k')$.
Vu l'hypothèse de la proposition,
on a alors $X(k')\neq\emptyset$,
d'où l'on déduit que $X(k) \neq \emptyset$ grâce au théorème d'Amer, Brumer et Coray.
\end{demo}

\bigskip
\begin{proposition}
\label{ch3findp4incondi}
Supposons qu'il existe $t \in \Srond$ de degré impair sur~$k$ tel que $\epsilon_t=1$.
Alors la surface~$X$ admet un point rationnel si l'obstruction de Brauer-Manin ne s'y
oppose pas.  En particulier, sous cette hypothèse, le principe de Hasse vaut pour~$X$ dès que $\Br(X)/\Br(k)=0$.
\end{proposition}

\bigskip
\begin{demo}
La proposition~\ref{ch3findp4amerbrumer} permet de supposer que le point~$t$ est $k$\nobreakdash-rationnel.
L'hypothèse $\epsilon_t=1$ et la $k$\nobreakdash-rationalité de~$t$ entraînent que
la surface~$X$ est fibrée en coniques au-dessus d'une conique, avec quatre fibres géométriques singulières (cf.~\cite[Corollary~3.4]{kunsktsf}).
Il est connu que de telles surfaces admettent un point rationnel
dès que l'obstruction de Brauer-Manin ne s'y oppose pas (cf.~\cite[Théorème~2]{ctsemtn}, \cite{salpezzoconique}).
\end{demo}

\bigskip
Voici maintenant quelques conditions suffisantes pour que l'hypothèse~(v) du théorème~\ref{ch3thprindp4} soit
satisfaite.

\bigskip
\begin{proposition}
\label{ch3findp4cor1rat}
Supposons que~$t_0$ et~$t_1$ soient $k$\nobreakdash-rationnels.
Si la condition~(\ref{ch3enonceshyp}) est satisfaite et que $\Br(X)/\Br(k)=0$,
l'hypothèse~(v) du théorème~\ref{ch3thprindp4} est satisfaite.
\end{proposition}

\bigskip
\begin{demo}
Supposons, par l'absurde, que $\epsilon_0 \neq 1$ et que pour tout $t \in \Srond'$,
l'image de~$\epsilon_0$ dans $\kappa(t)^\star/\kappa(t)^{\star 2}$ appartienne
au sous-groupe $\{1,\epsilon_t\}$.  Choisissant $t=t_1$, on en déduit
que $\epsilon_0=\epsilon_{t_1}$.
Cette égalité et la proposition~\ref{ch3brproduit} montrent que l'entier noté~$d$
dans l'énoncé du théorème~\ref{ch3thbr} est $<\Card(\Srond)-1$.
Par ailleurs, la condition~(\ref{ch3enonceshyp}) et l'hypothèse que $\epsilon_0\neq 1$ entraînent que l'entier noté~$n$
dans l'énoncé du théorème~\ref{ch3thbr} est égal à $\Card(\Srond)$,
compte tenu que tout point de~$\Srond'$ est de degré~$\leq 3$.
Il résulte donc du théorème~\ref{ch3thbr} que $\Br(X)/\Br(k)\neq 0$, ce qui contredit l'hypothèse de la proposition.
\end{demo}

\bigskip
\begin{proposition}
\label{ch3findp4cortransit}
Supposons que~$t_0$ soit $k$\nobreakdash-rationnel.
Si le groupe~$\Gamma$ agit transitivement sur $\{\uplet{t_1}{t_4}\}$ et si
$$
\epsilon_0 \not\in \Ker\left(k^\star/k^{\star 2} \rightarrow \kappa(t)^\star/\kappa(t)^{\star 2}\right)\setminus\{1\} \rlap{\text{,}}
$$
où~$t$ désigne l'unique élément de~$\Srond'$, l'hypothèse~(v) du théorème~\ref{ch3thprindp4} est satisfaite.
\end{proposition}

\bigskip
\begin{lemme}
\label{ch3findp4lemmecor}
Supposons que~$t_0$ soit $k$\nobreakdash-rationnel.
Si $\epsilon_0 \neq 1$, il existe $t \in \Srond'$ tel que l'image de~$\epsilon_0$
dans $\kappa(t)^\star/\kappa(t)^{\star 2}$ soit distincte de~$\epsilon_t$.
\end{lemme}

\bigskip
\begin{demo}
Supposons que pour tout $t \in \Srond'$, l'image de~$\epsilon_0$ dans $\kappa(t)^\star/\kappa(t)^{\star 2}$
soit égale à~$\epsilon_t$.  On a alors $N_{\kappa(t)/k}(\epsilon_t)=\epsilon_0^{[\kappa(t):k]}$
pour tout $t \in \Srond$.  Compte tenu que $\sum_{t\in \Srond}[\kappa(t):k]=5$
et que $\prod_{t\in\Srond} N_{\kappa(t)/k}(\epsilon_t)=1$ (cf.~proposition~\ref{ch3brproduit}),
il en résulte que $\epsilon_0=1$.
\end{demo}

\bigskip
\begin{demo}[ de la proposition~\ref{ch3findp4cortransit}]%
La condition~(\ref{ch3enonceshyp}) est vide lorsque~$\Gamma$ agit transitivement sur $\{\uplet{t_1}{t_4}\}$.
De même, le groupe $\Br(X)/\Br(k)$ est automatiquement nul (cf.~corollaire~\ref{ch3brcor541}).
Notant~$t$ l'unique point de~$\Srond'$, le lemme~\ref{ch3findp4lemmecor} montre que l'image de~$\epsilon_0$ dans
$\kappa(t)^\star/\kappa(t)^{\star 2}$ est distincte de~$\epsilon_t$.  Elle est non triviale par hypothèse,
d'où le résultat.
\end{demo}

\bigskip
\begin{proposition}
\label{ch3findp4corprimit}
Supposons que~$t_0$ soit $k$\nobreakdash-rationnel.
Si le groupe~$\Gamma$ agit transitivement et primitivement sur $\{\uplet{t_1}{t_4}\}$
(ce qui est par exemple le cas s'il agit $2$\nobreakdash-transitivement), l'hypothèse~(v) du théorème~\ref{ch3thprindp4} est satisfaite.
\end{proposition}

\bigskip
\begin{demo}
Sous cette hypothèse, si~$t$ désigne l'unique élément de~$\Srond'$, l'extension quartique $\kappa(t)/k$
ne contient aucune sous-extension quadratique, de sorte que l'on peut appliquer
la proposition~\ref{ch3findp4cortransit}.
\end{demo}

\bigskip
Nous pouvons à présent terminer la démonstration du théorème~\ref{ch3thprindp4}.

Supposons l'hypothèse~(iv) de ce théorème satisfaite.
S'il existe $t \in \Srond$ tel que $\epsilon_t=1$, la proposition~\ref{ch3findp4incondi} permet
de conclure.  Sinon, la condition~(\ref{ch3enonceshyp}) est satisfaite et la proposition~\ref{ch3findp4cor1rat}
montre donc que l'hypothèse~(v) du théorème l'est aussi, d'où le résultat.

Supposons l'hypothèse~(iii) du théorème~\ref{ch3thprindp4} satisfaite.
Tout point de~$\Srond$ est alors de degré impair sur~$k$.  S'il existe $t \in \Srond$ tel
que $\epsilon_t=1$, la proposition~\ref{ch3findp4incondi} permet donc de conclure. Sinon,
quitte à renuméroter les~$t_i$, on peut supposer que~$t_0$ et~$t_1$ sont $k$\nobreakdash-rationnels;
la condition~(\ref{ch3enonceshyp}) est alors satisfaite et l'hypothèse~(v) l'est donc aussi, d'après
la proposition~\ref{ch3findp4cor1rat}.

Si l'hypothèse~(ii) du théorème~\ref{ch3thprindp4} est satisfaite,
on peut supposer que~$t_0$ est $k$\nobreakdash-rationnel,
quitte à renuméroter les~$t_i$.
Il suffit alors d'appliquer la proposition~\ref{ch3findp4corprimit}
pour conclure.

Enfin, si l'hypothèse~(i) du théorème~\ref{ch3thprindp4} est satisfaite, l'ensemble~$\Srond$ est un singleton.
Notons~$t$ son unique élément.  Comme l'extension $\kappa(t)/k$ est de degré~$5$, le théorème d'Amer, Brumer et Coray
montre qu'il suffit d'établir le principe de Hasse pour la surface $X \otimes_k \kappa(t)$; or celle-ci
est justiciable d'un cas déjà établi du théorème~\ref{ch3thprindp4} puisque l'hypothèse~(ii) relative à $X\otimes_k\kappa(t)$
est satisfaite.

\subsection{Groupe de Brauer et obstruction à la méthode}
\label{ch3parobstruction}

Le théorème~\ref{ch3thprindp4} ne prédit la validité du principe de Hasse que pour des surfaces
de del Pezzo de degré~$4$ vérifiant $\Br(X)/\Br(k)=0$.  De fait,
cette hypothèse est nécessaire pour que la \condD{} générique soit satisfaite
(cf.~proposition~\ref{ch3verifprin}; il est facile de voir que
l'implication apparaissant dans l'énoncé de cette
proposition est une équivalence).
L'objet de ce paragraphe est d'établir
\emph{a priori} (c'est-à-dire indépendamment de tout calcul)
l'existence, dès que $\Br(X)/\Br(k)\neq 0$,
d'une obstruction à la méthode employée pour prouver le théorème~\ref{ch3thprindp4}.
Plus précisément,
nous allons montrer que toute classe non constante de~$\Br(X)$ donne naissance à une classe
de $\Br(C_{H^0})$ non verticale par rapport à~$\pi_{H^0}$.
Il est vraisemblable que si \mbox{$D \subset H^0$} est une droite suffisamment générale
et si $A \in \Br(C_{H^0})$ est une classe non verticale,
la restriction de~$A$ à~$\pi^{-1}_{H^0}(D)$ ne soit jamais verticale.
À tout le moins, on ne sait pas trouver systématiquement de droite~$D$
mettant cette propriété en défaut; en particulier ne sait-on pas trouver de droite~$D$
telle que la \condD{} associée au pinceau $\pi^{-1}_{H^0}(D) \rightarrow D$ soit
satisfaite (cf.~corollaire~\ref{ch1dbrdvert}; un énoncé analogue au corollaire~\ref{ch1dbrdvert} vaut dans la situation
du chapitre~2).

Les notations sont celles du paragraphe~\ref{ch3parconssd}.
Le morphisme sous-entendu dans $\Br_\vert(C_\Lambda)$ et $\Br_\hor(C_\Lambda)$
(resp.~dans $\Br_\vert(C_{H^0})$ et $\Br_\hor(C_{H^0})$)
sera toujours~$\pi_\Lambda$ (resp.~$\pi_{H^0}$).

\bigskip
\begin{theoreme}
L'application\index{groupe de Brauer} naturelle $\Br(X)/\Br(k) \rightarrow \Br_{\hor}(C_{H^0})$
(cf.~diagramme~(\ref{ch3conssddiag})) est\index{groupe de Brauer!horizontal} injective.
\end{theoreme}

\bigskip
\begin{demo}
Notons $\alpha \colon \Br(X)/\Br(k) \rightarrow \Br_\hor(C_\Lambda)$ l'application induite
par la projection $C_\Lambda\rightarrow X$.
Comme le morphisme $H^0 \xrightarrow{\;\rho\circ\sigma\;\!} \Lambda$ est dominant,
la flèche $\Br(C_\Lambda) \rightarrow \Br(C_{H^0})$ issue du diagramme~(\ref{ch3conssddiag})
induit une application $\beta\colon \Br_\hor(C_\Lambda) \rightarrow \Br_\hor(C_{H^0})$.
La composée
$$
\xymatrix{
\Br(X)/\Br(k) \ar[r]^(.535)\alpha & \Br_\hor(C_\Lambda) \ar[r]^(.48)\beta & \Br_\hor(C_{H^0})
}
$$
est égale à l'application dont on veut montrer qu'elle est injective.

\bigskip
\begin{lemme}
\label{ch3obsvertvert}
La restriction de~$\beta$ au sous-groupe de $2$\nobreakdash-torsion de $\Br_\hor(C_\Lambda)$
est injective.
\end{lemme}

\bigskip
\begin{demo}
Comme la variété~$C_\Lambda$ est régulière et connexe (cf.~proposition~\ref{ch3conssdreguliers}),
son groupe de Brauer s'injecte dans celui de son corps de fonctions.  Il suffit donc de s'intéresser
aux groupes de Brauer des corps de fonctions des variétés envisagées
(cf.~proposition~\ref{ch3conssdgeomint} pour l'intégrité de~$C_{H^0}$);
ceux-ci s'inscrivent dans le carré commutatif suivant:
$$
\xymatrix{
\Br(\kappa(C_{\Lambda})) \ar[r]^{\gamma'} & \Br(\kappa(C_{H^0})) \\
\Br(\kappa(\Lambda)) \ar[u]_\delta \ar[r]^\gamma & \Br(\kappa(H^0)) \ar[u]_{\delta'}\rlap{\text{.}}
}
$$
Étant donné $A \in \Br(\kappa(C_\Lambda))$ tel que $\gamma'(A)\in\Im(\delta')$ et $2A \in \Im(\delta)$, on doit alors prouver
que $A \in \Im(\delta)$.  Soit $B \in \Br(\kappa(H^0))$ tel que $\delta'(B)=\gamma'(A)$.
Comme $\kappa(C_{H^0}) = \kappa(C_\Lambda) \otimes_{\kappa(\Lambda)} \kappa(H^0)$,
on a $\Cores_{\kappa(C_{H^0})/\kappa(C_\Lambda)}(\delta'(B))=\delta(\Cores_{\kappa(H^0)/\kappa(\Lambda)}(B))$.
Par ailleurs, l'extension $\kappa(C_{H^0})/\kappa(C_{\Lambda})$ étant de degré~$3$,
on a $\Cores_{\kappa(C_{H^0})/\kappa(C_\Lambda)}(\gamma'(A))=3A$.
Combinant les trois égalités qui précèdent, on obtient
$A = \delta(\Cores_{\kappa(H^0)/\kappa(\Lambda)}(B)) - 2A$, d'où $A \in \Im(\delta)$.
\end{demo}

\bigskip
\begin{lemme}
\label{ch3obsbvertlambda}
L'application~$\alpha$ est injective.
\end{lemme}

\bigskip
\begin{demo}
La projection $C_\Lambda \rightarrow X$ est un fibré projectif (localement libre).  Par conséquent,
l'application $\Br(X) \rightarrow \Br(C_\Lambda)$ qu'elle induit est un isomorphisme.
L'injectivité de~$\alpha$ équivaut donc à la nullité
du groupe $\Br_\vert(C_\Lambda)/\Br(k)$; c'est cette dernière propriété que nous allons démontrer.
Rappelons qu'il existe une hypersurface quadrique lisse $Q \subset \Lambda$ telle que les
fibres de~$\pi_\Lambda$ au-dessus des points de codimension~$1$ de $\Lambda\setminus Q$ soient géométriquement intègres
(cf.~propositions~\ref{ch3conssdrqhypirr} et~\ref{ch3conssdi1i2}).
Il en résulte que toute classe de $\Br_\vert(C_\Lambda)$ est l'image réciproque par~$\pi_\Lambda$
d'une classe de $\Br(\Lambda \setminus Q)$ (cf.~\cite[Proposition~1.1.1, Theorem~1.3.2]{cscrelle94}
et la proposition~\ref{ch3conssdreguliers}).
Pour conclure, il reste seulement à montrer que si $A \in \Br(\Lambda \setminus Q)$ est
telle que $\pi_\Lambda^\star A \in \Br(C_\Lambda)$, alors
$\pi_\Lambda^\star A$ appartient à l'image de $\Br(k)$.

Comme~$Q$ est géométriquement intègre sur~$k$, la flèche de restriction
$$\Br(X)/\Br(k) \longrightarrow \Br(X \otimes_k \kappa(Q))/\Br(\kappa(Q))$$ est
injective (l'argument est exactement le même que celui employé dans les preuves
des lemmes~\ref{ch3brisopic} et~\ref{ch3brisobr}; le point clé est que le groupe de Picard géométrique
de~$X$ est un $\Z$\nobreakdash-module libre de type fini); pour prouver la proposition, on peut donc
supposer que $Q(k)\neq \emptyset$, quitte à étendre les scalaires de~$k$ à~$\kappa(Q)$.
Soit alors $q \in Q(k)$.  Choisissons une droite $k$\nobreakdash-rationnelle $D \subset \Lambda$
passant par~$q$, tangente à~$Q$ en~$q$ mais non incluse dans $Q \cup R$,
et posons $D^0=D \setminus\{q\}$.
(Une telle droite existe car $Q \cup R$ ne contient pas de plan.)
Comme $D \not\subset Q \cup R$, la fibre de~$\pi_\Lambda$ au-dessus du point
générique de~$D$ est irréductible.  Le morphisme~$\pi_\Lambda$ étant propre et plat
et la variété~$C_\Lambda$ étant de Cohen-Macaulay, on en déduit que la surface
$\pi_\Lambda^{-1}(D^0)$ est irréductible (cf.~preuve de la proposition~\ref{ch3conssdgeomint}).
Le morphisme composé $\pi_\Lambda^{-1}(D^0) \hookrightarrow C_\Lambda \rightarrow X$
est birationnel; en effet, si $L^0 \subset \P^4_k$ désigne l'hyperplan correspondant à~$q$,
on peut définir un inverse birationnel de ce morphisme sur l'ouvert dense $X \setminus (X \cap L^0) \subset X$
en associant à~$x$ le couple formé de~$x$ et de l'unique hyperplan de~$D$ contenant~$x$.
Il en résulte que la flèche horizontale du triangle commutatif
$$
\xymatrix{
\Br(C_\Lambda) \ar[dr] \\
\Br(X) \ar[r] \ar[u]_{\wr} & \Br(\pi_\Lambda^{-1}(D^0))
}
$$
est injective.  La flèche de restriction
$\Br(C_\Lambda) \rightarrow \Br(\pi_\Lambda^{-1}(D^0))$ est donc elle aussi injective.  Celle-ci s'insère dans le diagramme
commutatif
$$
\myxyhook\xymatrix{
\Br(C_\Lambda) \ar@{^{ (}->}[d] \ar@{^{ (}->}[rd] \\
\Br(\pi_\Lambda^{-1}(\Lambda \setminus Q)) \ar[r] & \Br(\pi_\Lambda^{-1}(D^0)) \\
\Br(\Lambda\setminus Q) \ar[u]_{\pi_\Lambda^\star} \ar[r] & \Br(D^0) \ar[u]_{\pi_\Lambda^\star} \rlap{\text{,}}
}
$$
où les flèches dépourvues d'étiquette sont les flèches de restriction; on a en effet $D^0 \subset \Lambda \setminus Q$
puisque~$Q$ est une quadrique, que~$D$ est tangente à~$Q$ en~$q$ et que $D \not\subset Q$.
La flèche naturelle $\Br(k)\rightarrow \Br(D^0)$
est un isomorphisme
puisque la variété $D^0$ est $k$\nobreakdash-isomorphe à~$\A^1_k$. Grâce à cette remarque, il est maintenant évident sur le diagramme ci-dessus que
si $A \in \Br(\Lambda \setminus Q)$ est telle
que $\pi_\Lambda^\star A \in \Br(C_\Lambda)$, alors $\pi_\Lambda^\star A$ appartient à l'image de~$\Br(k)$.
\end{demo}

\bigskip
On déduit des lemmes~\ref{ch3obsvertvert} et~\ref{ch3obsbvertlambda} que la restriction de $\beta \circ \alpha$
au sous-groupe de $2$\nobreakdash-torsion de $\Br(X)/\Br(k)$ est injective.
Compte tenu que $\Br(X)/\Br(k)$ est tué par~$2$ (cf.~théorème~\ref{ch3thbr} et la remarque qui suit sa démonstration),
ceci prouve le théorème.
\end{demo}

\section[Intersections de deux quadriques dans~$\P^n$ pour $n\geq 5$]{Principe de Hasse pour les intersections de deux quadriques dans~$\P^n$ avec $n \geq 5$}
\subsection{Un résultat de monodromie}
\label{ch3monodr}

Reprenons les notations du paragraphe~\ref{ch3sectiongeneralites} et supposons la variété~$X$
lisse sur~$k$ et purement de codimension~$2$ dans~$\P(V)$.  Nous avons défini
un morphisme $p \colon Z \rightarrow \P(V^\star)$, fini et plat de degré $n=\dim \P(V)$
(cf.~proposition~\ref{ch3genpfiniplat}).

\bigskip
\begin{theoreme}
\label{ch3propmonod}
Le\index{pinceau de quadriques!monodromie du revêtement associé|(} revêtement $p \colon Z \rightarrow \P(V^\star)$ est irréductible et génériquement étale.
Son groupe de monodromie est isomorphe au groupe symétrique~$\mathfrak{S}_n$.
\end{theoreme}

\bigskip
(Le \emph{groupe de monodromie} d'un revêtement irréductible et génériquement étale est par définition le groupe
de Galois d'une clôture galoisienne de sa fibre générique.)

\bigskip
Nous donnons deux démonstrations de cette proposition.

\bigskip
\begin{premdemo}
Soit $U \subset \P(V^\star)$ le plus grand ouvert au-dessus duquel le morphisme~$p$ soit étale.
Comme il est bien connu (cf.~\cite[Exp.~V]{sga1}), la proposition peut se reformuler ainsi: l'ouvert~$U$ est non vide et
si~$b$ est un point géométrique de~$U$,
le morphisme canonique du groupe fondamental $\pi_1(U,b)$ vers le groupe symétrique sur $p^{-1}(b)$ est surjectif
(cette propriété ne dépend pas du point~$b$ considéré).
C'est en prenant pour~$b$ le point générique géométrique de~$U$ que l'on fait le
lien avec l'énoncé de la proposition en termes de clôture galoisienne de la fibre générique de~$p$.
La clé de la preuve va être de choisir le point base~$b$ plus judicieusement, en tenant compte de la géométrie de la situation.

Rappelons que l'on dispose d'une famille $\uplet{t_0}{t_n} \in \P^1(\ksep)$
et d'une famille $\uplet{\fp{P}0}{\fp{P}n}$ de $\ksep$\nobreakdash-points de~$\P(V)$ qui ne sont pas contenus dans
un même hyperplan (cf.~corollaire~\ref{ch3genpiengendrent}).  Pour $i \in \{\uplet{0}{n}\}$, notons
$\fp{P}i^\star$ le $\ksep$\nobreakdash-point de~$\P(V^\star)$ associé à l'hyperplan de~$\P(V)$ contenant $\ensemble{\fp{P}j}{j \in \{\uplet{0}{n}\}
\setminus \{i\}}$.  Pour chaque $i \in \{\uplet{0}{n}\}$,
l'inclusion $Z \subset \P^1_k \times_k \P(V^\star)$ permet d'identifier la fibre $p^{-1}(\fp{P}i^\star)$ à l'ensemble des couples
$(t_j,\fp{P}i^\star)$ pour $j \in \{\uplet{0}{n}\}\setminus\{i\}$.  Comme les~$t_j$ sont deux à deux distincts,
il en résulte que $p^{-1}(\fp{P}i^\star)$ est de cardinal~$n$ pour tout~$i$; les points
$\uplet{\fp{P}0^\star}{\fp{P}n^\star}$ appartiennent donc tous à~$U(\ksep)$, ce qui prouve au passage que $U\neq\emptyset$.

Posons $b=\fp{P}0^\star$.
Nous allons maintenant montrer que l'image du morphisme canonique
$\pi_1(U,b) \rightarrow \Aut(\pi^{-1}(b))$ (où le symbole «~$\Aut$~» désigne simplement le groupe
symétrique) contient toutes les transpositions.

Si~$S$ est un schéma, notons $\FEt(S)$ la catégorie des $S$\nobreakdash-schémas finis étales. Notons~$\Ens$
la catégorie des ensembles finis.
Pour $i \in\{\uplet{0}{n}\}$, notons $\ff{F}i \colon \FEt(U) \rightarrow \Ens$ le foncteur fibre en~$\fp{P}i^\star$:
par définition, si $S \rightarrow U$ est un morphisme fini étale, $\ff{F}i(S)$ est sa fibre en~$\fp{P}i^\star$.
Pour $i,j\in\{\uplet{0}{n}\}$ distincts, notons $D_{ij}$ la $\ksep$\nobreakdash-droite de $\P(V^\star)$ passant par~$\fp{P}i^\star$
et~$\fp{P}j^\star$ et $\ff{F}{ij} \colon \FEt(U \cap D_{ij}) \rightarrow \Ens$ le foncteur fibre en~$\fp{P}i^\star$.
Compte tenu que $U \cap D_{ij}=U \cap D_{ji}$ et que ce schéma est connexe,
les foncteurs~$\ff{F}{ij}$ et~$\ff{F}{ji}$ sont (non canoniquement) isomorphes (cf.~\cite[Exp.~V, §7]{sga1}).
Choisissons, indépendamment pour chaque couple d'entiers $i,j \in \{\uplet{0}{n}\}$ distincts,
un isomorphisme $u_{ij} \colon \ff{F}{ij} \isoto \ff{F}{ji}$.  Des isomorphismes $v_{ij} \colon \ff{F}i \isoto \ff{F}j$
s'en déduisent par composition avec le foncteur de restriction $\FEt(U) \rightarrow \FEt(U \cap D_{ij})$.

\bigskip
\begin{lemme}
\label{ch3monodlemme}
Soient $i,j\in\{\uplet{0}{n}\}$ distincts. La bijection
$$p^{-1}(\fp{P}i^\star) \longrightarrow p^{-1}(\fp{P}j^\star)$$ induite par~$v_{ij}$
applique $(t_\ell,\fp{P}i^\star)$ sur $(t_\ell,\fp{P}j^\star)$ pour tout $\ell \in \{\uplet{0}{n}\}\setminus\{i,j\}$
et $(t_j,\fp{P}i)^\star$ sur $(t_i,\fp{P}j^\star)$.
\end{lemme}

\bigskip
\begin{demo}
Notons $f \colon p^{-1}(\fp{P}i^\star) \rightarrow p^{-1}(\fp{P}j^\star)$ cette application.
Il suffit de prouver que $f((t_\ell,\fp{P}i^\star))=(t_\ell,\fp{P}j^\star)$
pour tout $\ell\in\{\uplet{0}{n}\}\setminus\{i,j\}$, puisque~$f$ est bijective.
Fixons donc un tel entier~$\ell$.
La droite~$D_{ij}$ s'identifie au pinceau des hyperplans de~$\P(V)$
contenant $\ensemble{\fp{P}s}{s \in \{\uplet{0}{n}\}\setminus \{i,j\}}$.
En particulier, les hyperplans de~$D_{ij}$ contiennent tous~$\fp{P}\ell$;
il en résulte que l'on peut définir une section~$s$ du $\ksep$\nobreakdash-morphisme
$p^{-1}(U \cap D_{ij}) \rightarrow U \cap D_{ij}$
induit par~$p$ en posant $s(L)=(t_\ell,L)$.
Considérons~$s$ comme un morphisme dans la catégorie \mbox{$\FEt(U \cap D_{ij})$}.
La fonctorialité de l'isomorphisme~$u_{ij}$ fournit alors un carré commutatif
$$
\xymatrix{
\{\fp{P}i^\star\} \ar[r] \ar[d] & p^{-1}(\fp{P}i^\star) \ar[d]^f \\
\{\fp{P}j^\star\} \ar[r] & p^{-1}(\fp{P}j^\star) \rlap{\text{,}}
}
$$
où les flèches horizontales sont induites par~$s$.  Celui-ci montre
bien que $f((t_\ell,\fp{P}i^\star))=(t_\ell,\fp{P}j^\star)$, vu la définition de~$s$.
\end{demo}

\bigskip
Pour $i,j\in\{\uplet{1}{n}\}$ distincts, notons~$w_{ij}$ l'automorphisme de~$\ff{F}0$
défini comme la composée
$$
\xymatrix@C=6ex{
\ff{F}0 \ar[r]^{v_{0i}} & \ff{F}i \ar[r]^{v_{ij}} & \ff{F}j \ar[r]^{v_{j0}} & \ff{F}0 \rlap{\text{.}}
}
$$
Le lemme~\ref{ch3monodlemme} montre que l'automorphisme de $p^{-1}(b)$ induit par~$w_{ij}$
n'est autre que la transposition $(i\ j)$, si l'on identifie $p^{-1}(b)$ à
$\{\uplet{1}{n}\}$ par l'application $(t_i,\fp{P}0^\star)\mapsto i$.
Par définition du groupe fondamental, les auto\-morphismes~$w_{ij}$ sont des éléments
de $\pi_1(U,b)$. Nous venons de prouver que l'ensemble de leurs images dans $\Aut(p^{-1}(b))$
est exactement l'ensemble des transpositions de~$p^{-1}(b)$; ainsi le morphisme $\pi_1(U,b)\rightarrow \Aut(p^{-1}(b))$
est-il bien surjectif.
\end{premdemo}

\bigskip
\begin{remarque}
Lorsque~$k$ est de caractéristique~$0$, on voit facilement que l'on peut supposer que $k=\C$, auquel cas
la démonstration ci-dessus s'interprète topologiquement comme suit.  Pour simplifier les notations, prenons $i=1$
et~$j=2$. Considérons alors la courbe connexe $T=D_{01} \cup D_{12} \cup D_{02}$, que l'on pourrait qualifier
de «~triangle~».  Comme les sommets de ce triangle appartiennent à~$U$, il existe un lacet~$\gamma$ sur~$T(\C)$,
d'origine~$\fp{P}0^\star$, évitant le lieu de branchement de~$p$ et «~faisant une fois le tour de~$T$~»
(dans un sens ou dans l'autre; pour fixer les idées, mettons que~$\gamma$ visite le point~$\fp{P}1^\star$
avant le point~$\fp{P}2^\star$).
Les fibres fermées de~$p$ sont canoniquement des sous-ensembles
de~$\P^1(\C)$.  De ce point de vue, la fibre $p^{-1}(\fp{P}0^\star)$ s'identifie à
l'ensemble $\{\uplet{t_1}{t_n}\}$; voyons ce qu'il en advient lorsqu'on parcourt le lacet~$\gamma$.
Dans un premier temps, si l'on va de~$\fp{P}0^\star$ à~$\fp{P}1^\star$ tout en ne quittant pas la droite~$D_{01}$,
les~$t_i$ pour $i>2$ restent fixes tandis que le $n$\nobreakdash-ème point de la fibre, initialement égal à~$t_1$,
se déplace progressivement pour finalement atteindre~$t_0$.  Dans un second temps, les~$t_i$
pour $i>2$ ou $i=0$ sont fixes et le point restant passe de~$t_2$ à~$t_1$.
Enfin, le parcours du troisième côté de~$T$ laisse les~$t_i$ pour $i>2$ ou $i=1$ fixes et déplace le
point restant de~$t_0$ vers~$t_2$.  Le lacet~$\gamma$ agit donc bien par la transposition attendue.

Au premier abord, il pourrait sembler surprenant que l'action de~$\gamma$ sur $p^{-1}(\fp{P}0^\star)$ ne dépende
pas du choix de ce lacet;
il existe pourtant une infinité de chemins deux à deux non homotopes
permettant de relier~$\fp{P}i^\star$ à~$\fp{P}j^\star$ sur $(U \cap D_{ij})(\C)$, pour chaque couple $(i,j)$.
(Ce choix correspond à celui des isomorphismes~$u_{ij}$ dans la preuve algébrique.)
L'explication est la suivante.
Bien que le revêtement plat $p^{-1}(D_{ij}) \rightarrow D_{ij}$
soit ramifié au-dessus de $R \cap D_{ij}$, où $R=\P(V^\star)\setminus U$,
il ne s'agit que de ramification factice:
la courbe $p^{-1}(D_{ij})$ n'est pas normale, et le revêtement de~$D_{ij}$ obtenu après normalisation
de $p^{-1}(D_{ij})$ est étale (et même constant).  Ce phénomène a lieu car l'hypersurface $R \subset \P(V^\star)$
est tangente à la droite~$D_{ij}$ en tout point d'intersection.
\end{remarque}

\bigskip
\begin{secdemo}
Les fibres géométriques de la première projection $Z \rightarrow \P^1_k$ sont irréductibles; en effet, d'un
point de vue ensembliste, la fibre
au-dessus de $t \in \P^1(\ksep)$ coïncide avec la variété duale de la quadrique lisse~$Q_t$
si $t \not\in \{\uplet{t_0}{t_n}\}$ ou avec l'ensemble des hyperplans de~$\P(V)$
passant par~$\fp{P}i$ si $t=t_i$.
Par conséquent, la variété~$Z$ est géométriquement irréductible.

Si $n=2$, il n'y a plus rien à démontrer. Supposons donc $n \geq 3$.
Il résulte alors de la description des fibres géométriques de~$p$
(cf.~les remarques qui précèdent la proposition~\ref{ch3genpfiniplat})
et de la proposition~\ref{ch3genpropeq} que
le lieu de branchement $R \subset \P(V^\star)$ de~$p$ s'identifie à l'ensemble des hyperplans de~$\P(V)$
tangents à~$X$, c'est-à-dire à la variété duale de~$X$ puisque~$X$ est lisse.
Ceci montre déjà que $R \neq \P(V^\star)$, donc que~$p$ est génériquement étale,
et que~$R$ est un fermé géométriquement irréductible
(cf.~\cite[Proposition~3.1.4]{katzsga}; la variété $X$ est elle-même géométriquement irréductible
car c'est une intersection complète lisse de dimension~$\geq 1$ dans~$\P(V)$).
Notons $U=\P(V^\star) \setminus R$.  Si~$R$ était de codimension~$\geq 2$ dans~$\P(V^\star)$,
la variété $U\otimes_k \ksep$ serait simplement connexe
(cf.~\cite[Exp.~XI, Proposition~1.1 et Exp.~X, Corollaire~3.3]{sga1})
et $p^{-1}(U\otimes_k \ksep)$ serait donc un revêtement constant de $U\otimes_k \ksep$, contredisant ainsi l'irréductibilité
géométrique de~$Z$.  Le fermé $R \subset \P(V^\star)$ est donc une hypersurface irréductible.

Comme~$R$ est une hypersurface, le pinceau des hyperplans de~$\P(V)$ contenant $\uplet{\fp{P}2}{\fp{P}n}$
rencontre~$R$.  Autrement dit, il existe $r \in R(\ksep)$ tel que $p^{-1}(r)$ contienne
les points $(t_i,r)$ pour $i \in \{\uplet{2}{n}\}$.  En particulier, comme les~$t_i$ sont deux à deux distincts,
cette fibre géométrique possède
au plus un point double.  Compte tenu de la semi-continuité inférieure du nombre géométrique de
points dans les fibres d'un morphisme fini et plat (cf.~\cite[15.5.1~(i)]{ega43}), il en résulte
que la fibre de~$p$ au-dessus du point générique de~$R$ possède au plus un point double.

Notons~$K$ le corps des fonctions de la grassmannienne des droites de~$\P(V^\star)$
et~$D$ la droite générique de~$\P(V^\star)$ (c'est-à-dire la fibre générique du fibré en droites
tautologique sur la grassmannienne des droites de~$\P(V^\star)$).
Posons $Z_D=Z \times_{\P(V^\star)} D$, où le produit fibré
est relatif à~$p$ et au morphisme canonique $i \colon D \rightarrow \P(V^\star)$.
Comme~$p$ est dominant et que~$Z$ est géométriquement irréductible sur~$k$,
$Z_D$ est géométriquement irréductible sur~$K$ (cf.~\cite[I, Th.~6.10]{jouanolou}).
Soit $Z'_D \rightarrow Z_D$ la normalisation de~$Z_D$.
On a vu que les fibres de~$p$ au-dessus des points de codimension~$1$ possèdent au plus un point double.
Compte tenu que l'image du morphisme~$i$ ne contient aucun point de codimension~$>1$,
il en résulte que les fibres de la projection $Z_D \rightarrow D$ possèdent au plus un point double.
Il en va nécessairement de même pour la projection $Z'_D \rightarrow D$. Ainsi disposons-nous
d'une courbe régulière~$Z'_D$, géométriquement irréductible sur~$K$, et d'un $K$\nobreakdash-morphisme $Z'_D \rightarrow D$
génériquement étale de degré~$n$ dont les fibres ont au plus un point double, où~$D$ est $K$\nobreakdash-isomorphe à~$\P^1_K$.
Il est bien connu, au moins en caractéristique~$0$, que le groupe de monodromie d'un tel revêtement
est isomorphe à~$\mathfrak{S}_n$ (cf.~\cite[Proposition~2.2]{benderwitt} pour une preuve en toute caractéristique).
Comme le morphisme $\pi_1(i^{-1}(U),b) \rightarrow \pi_1(U,i(b))$ induit par~$i$ est
compatible aux actions respectives de ces groupes sur la fibre~$p^{-1}(i(b))$, où~$b$ désigne un
point géométrique de~$i^{-1}(U)$, le résultat voulu s'en\index{pinceau de quadriques!monodromie du revêtement associé|)} déduit.
\end{secdemo}

\subsection{Preuve du théorème~\ref{ch3introp5}}
\label{ch3parpreuve5}

Le corps~$k$ est maintenant un corps de nombres.  Rappelons que le théorème~\ref{ch3introp5}
affirme que toute intersection lisse de deux quadriques dans~$\P^n_k$ pour $n \geq 5$
satisfait au principe de Hasse, si l'on admet l'hypothèse de Schinzel et la finitude des
groupes de Tate-Shafarevich des courbes elliptiques sur les corps de nombres.

\bigskip
\begin{demo}[ du théorème~\ref{ch3introp5}]%
Soit $X \subset \P^5_k$ une intersection lisse de deux quadriques.
Il est clair que l'on peut supposer~$X$ purement de codimension~$2$ dans~$\P^5_k$,
compte tenu du théorème de Hasse-Minkowski. Notons alors
$$
C=\bigensemble{(x,L)\in X \times_k (\P^5_k)^\star}{x \in L}
$$
et munissons ce fermé de $X \times_k (\P^5_k)^\star$ de sa structure de sous-schéma fermé
réduit.  La première projection $C \rightarrow X$ fait de~$C$ un fibré projectif (localement libre)
sur~$X$.  La variété~$C$ est donc lisse et géométriquement intègre sur~$k$.

Notons $f \colon C \rightarrow (\P^5_k)^\star$ la seconde projection
et $p \colon Z \rightarrow (\P^5_k)^\star$ le revêtement plat de degré~$5$
que l'on a associé à~$X$ au paragraphe~\ref{ch3sectiongeneralites}.
Le théorème~\ref{ch3propmonod} montre qu'il existe un ensemble hilbertien~$H$ de points $k$\nobreakdash-rationnels de $(\P^5_k)^\star$
au-dessus desquels la fibre de~$p$ est le spectre d'une extension finie de~$k$
dont la clôture galoisienne a pour groupe de Galois~$\mathfrak{S}_5$.
Comme la fibre générique de~$f$ est lisse (la variété~$C$ étant lisse sur~$k$), quitte à rétrécir~$H$,
on peut supposer que les fibres de~$f$ au-dessus de~$H$ sont lisses.
D'après le théorème~\ref{ch3thprindp4} (hypothèse~(i)), le principe de Hasse vaut pour les fibres de~$f$ au-dessus de~$H$.

Comme~$X$ est une intersection complète lisse de dimension~$\geq 3$ dans un espace
projectif, toutes les fibres de~$f$ sont géométriquement intègres
(cf.~\cite[Remark~7.5]{fultonlaz}).  (On aurait aussi pu déduire cette propriété de la proposition~\ref{ch3genpropeq},
du corollaire~\ref{ch3genpiengendrent} et de \cite[Lemma~1.11]{ctsansdi}.)
Nous pouvons donc appliquer le théorème suivant
au morphisme $f^{-1}(U) \rightarrow U$ induit
par~$f$, où~$U$ est un ouvert de~$(\P^5_k)^\star$ isomorphe à~$\A^5_k$;
l'hypothèse (Sect) est satisfaite en vertu
du théorème de Tsen (cf. la remarque~1 précédant le~§3 de~\cite{skofibration}).

\bigskip
\begin{theoreme}
Soit $f \colon X \rightarrow Y$ un morphisme projectif entre variétés géométriquement intègres sur un corps de nombres~$k$.
Supposons que les fibres de~$f$ soient géométriquement intègres, qu'il existe un entier $n \geq 1$ tel que~$Y$ soit le complémentaire
dans~$\A^n_k$ d'un fermé de~$\A^n_k$ de codimension~$\geq 2$, que
le principe de Hasse vaille pour les fibres de~$f$ au-dessus d'un sous-ensemble hilbertien de~$Y(k)$
et que la condition suivante soit satisfaite:
\begin{itemize}
\medskip
\item[]{\rm(Sect)} Il existe un ouvert dense~$U$ de l'espace des droites de~$\A^n_k$ tel que pour toute droite $D \in U$, si~$\ksep$ désigne une clôture algébrique de~$k$
et~$\eta$ le point générique de $D \otimes_k \ksep \subset \A^n_\ksep$, la $\kappa(\eta)$\nobreakdash-variété $X \times_Y \eta$ admet un point rationnel lisse.
\end{itemize}
\medskip

\noindent{}Alors le principe de Hasse vaut pour~$X$.
\end{theoreme}

\bigskip
\begin{demo}
Ce théorème est prouvé par Skorobogatov dans~\cite{skofibration}.
La technique sous-jacente, due à Colliot-Thélène, Sansuc et Swinnerton-Dyer, est expliquée en détail dans l'introduction de~\cite{ctsansdi}.
\end{demo}

\bigskip
Ainsi le principe de Hasse vaut-il pour~$C$, et donc pour~$X$.

Le théorème~\ref{ch3introp5} est à présent démontré pour $n=5$. Le cas général s'en déduit par
récurrence, en suivant exactement l'argument que l'on vient d'employer pour établir le
cas~$n=5$, à ceci près que l'on remplace chaque occurrence du chiffre~$5$ par~$n$
(bien entendu, l'ensemble~$H$ est maintenant inutile).
\end{demo}

\chapter*{Annexe\markboth{Annexe}{Annexe}}
\makeatletter
\def\theenonce{A.\@arabic\c@enonce}
\def\theequation{A.\@arabic\c@equation}
\setcounter{enonce}{0}
\setcounter{equation}{0}
\makeatother
\addcontentsline{toc}{chapter}{Annexe}

Un certain nombre de lemmes généraux au sujet des courbes elliptiques à réduction
semi-stable sur le corps des fonctions d'un trait hensélien sont utilisés de manière
répétée tout au long du texte.  Il a semblé plus commode de les rassembler
ici, en dépit de leur disparité.

Dans toute l'annexe, on fixe un anneau de valuation discrète~$R$ de
corps résiduel~$\kappa$ et de corps des fractions~$K$ et l'on suppose~$\kappa$
parfait\footnote{Cette hypothèse est sans doute inutile; elle permet néanmoins
d'appliquer la classification de Kodaira-Néron telle qu'on la trouve dans la
littérature.}
de caractéristique $p\neq 2$ (éventuellement nulle).
Si $E$ est une courbe elliptique sur~$K$,
on note $\Erond$ son modèle de Néron\index{modèle de Néron} au-dessus de~$R$,
$\Erond^0 \subset \Erond$ la composante
neutre de~$\Erond$ et $F$ le $\kappa$\nobreakdash-schéma en groupes fini étale
des composantes connexes de la fibre spéciale de~$\Erond$,
c'est-à-dire la fibre spéciale de $\Erond/\Erond^0$.
De même, si~$E'$ et~$E''$ sont des courbes elliptiques sur~$K$, on désigne
par $\Erond'$, $\Erond'^0$, $F'$, $\Erond''$, $\Erond''^0$, $F''$ les objets
correspondants.  Soit enfin~$\kappabarre$ une clôture algébrique de~$\kappa$.

\bigskip
\begin{lemme}
\label{annlemmesurje0}
Soient~$E'$ et~$E''$ des courbes elliptiques sur~$K$
et $\phi' \colon E' \rightarrow E''$ une isogénie.
Si~$p$ ne divise pas $\deg(\phi')$, le morphisme $\phi'^0 \colon \Erond'^0
\rightarrow \Erond''^0$ induit par $\phi'$ est un épimorphisme de faisceaux étales
sur $\Spec(R)$.
\end{lemme}

\bigskip
\begin{demo}
Il suffit de vérifier que~$\phi'^0$ est un morphisme de schémas surjectif et étale.
Il est surjectif d'après \cite[7.3/6]{blr}.  Pour qu'il soit étale, il
suffit que le morphisme $\Erond' \rightarrow \Erond''$ induit par~$\phi'$ le
soit, et ceci résulte de \cite[7.3/2~(b) et~7.3/5]{blr}.
\end{demo}

\bigskip
\begin{lemme}
\label{annlemmeisogekmsh}
Supposons~$R$ strictement hensélien.
Soient~$E'$ et~$E''$ des courbes elliptiques sur~$K$
et $\phi' \colon E' \rightarrow E''$ une isogénie dont le degré n'est pas divisible par~$p$.
Alors l'application $E'(K) \rightarrow E''(K)$ induite par~$\phi'$ est surjective si et
seulement si l'application $F'(\kappa) \rightarrow F''(\kappa)$ induite par~$\phi'$ l'est.
\end{lemme}

\bigskip
\begin{demo}
Comme~$R$ est strictement hensélien, les lignes du diagramme commutatif suivant sont
exactes:
$$
\xymatrix{
0 \ar[r] & \Erond'^0(R) \ar[r] \ar[d] & \Erond'(R) \ar[r] \ar[d] & F'(\kappa) \ar[r] \ar[d] & 0 \\
0 \ar[r] & \Erond''^0(R) \ar[r] & \Erond''(R) \ar[r] & F''(\kappa) \ar[r] & 0 \rlap{\text{.}}
}
$$
La flèche verticale de gauche est surjective d'après le lemme~\ref{annlemmesurje0},
ce qui permet de conclure, compte tenu que $\Erond'(R)=E'(K)$ et $\Erond''(R)=E''(K)$.
\end{demo}

\bigskip
\begin{proposition}
\label{annpropcyclpair}
Soit~$E$ une courbe elliptique sur~$K$, à réduction multiplicative.
Le groupe $F(\kappabarre)$ est cyclique.
Si de plus $\tors{2}E(K) \simeq \Z/2 \times \Z/2$, le groupe $F(\kappabarre)$ est cyclique
d'ordre pair.
Dans tous les cas, si le $\kappa$\nobreakdash-groupe $F$ n'est pas constant,
il existe une extension quadratique $\ell/\kappa$,
unique à isomorphisme près, telle que le $\ell$\nobreakdash-groupe $F\otimes_\kappa \ell$ soit constant.
Le groupe $\Gal(\ell/\kappa)$ agit alors sur~$F(\ell)$ par multiplication par~$-1$
(et en particulier $F(\kappa)$ est d'ordre~$\leq 2$).
\end{proposition}

\bigskip
\begin{demo}
Notons~$\sbarre$ le point géométrique de $\Spec(R)$ défini par~$\kappabarre$
et $E^\star$ un modèle propre et régulier minimal de~$E$ au-dessus de~$R$.
Que le groupe $F(\kappabarre)$ soit cyclique est une conséquence bien
connue de l'hypothèse de réduction multiplicative;
plus précisément, il résulte de cette hypothèse que la fibre~$E^\star_\sbarre$
de $E^\star$ au-dessus de~$\sbarre$
est réduite et qu'il existe
un générateur $\alpha \in F(\kappabarre)$
et une bijection
canonique et $\Gal(\kappabarre/\kappa)$\nobreakdash-équivariante $m$
de~$F(\kappabarre)$ sur l'ensemble des composantes irréductibles de~$E^\star_\sbarre$
tels que $\{m(\alpha),m(0),m(-\alpha)\}$ soit
exactement l'ensemble des composantes irréductibles de~$E^\star_\sbarre$
rencontrant~$m(0)$.

Faute de référence
satisfaisante,
nous donnons ici une preuve de cette dernière assertion.
(Le c\oe{}ur de la démonstration se trouve au bas de~\cite[p.~105]{neronihes}.)
L'ouvert de lissité de $E^\star$ sur $\Spec(R)$ est canoniquement
isomorphe à~$\Erond$ (cf.~\cite[Ch.~IV, §6, Th.~6.1]{silvaec2}). La fibre spéciale
de~$E^\star$ étant réduite (grâce à l'hypothèse de réduction multiplicative,
cf.~\cite[Ch~IV, §9, Th.~8.2]{silvaec2}), il en résulte une bijection canonique
et $\Gal(\kappabarre/\kappa)$\nobreakdash-équivariante
entre l'ensemble des composantes irréductibles de
la fibre
de $E^\star$ au-dessus de~$\sbarre$
et~$F(\kappabarre)$.
Pour vérifier qu'elle remplit la condition
voulue,
on peut supposer~$R$ strictement
hensélien, puisque la formation du modèle propre
et régulier minimal et la formation du modèle de Néron
commutent tous deux aux changements de base étales (cf.~
\cite[Chapter~9, Proposition~3.28]{liu} et \cite[1.2/2]{blr}, respectivement).
Vu la structure de la fibre spéciale de~$E^\star$ dans le cas de
réduction multiplicative (cf.~\cite[Ch~IV, §9, Th.~8.2]{silvaec2}),
il est possible de numéroter les composantes
irréductibles $\uplet{C_0}{C_{n-1}}$ de la
fibre spéciale de $E^\star$ de telle sorte que~$C_0$
soit la composante rencontrant la section nulle,
que $C_i$ rencontre transversalement $C_{i+1}$ en un unique point $z_i \in E^\star$
pour tout $i\in \{\uplet{0}{n-1}\}$, où l'on a posé $C_n=C_0$,
que les~$z_i$ soient deux à deux distincts
et qu'il n'y ait aucune autre intersection entre les~$C_i$ que celles que l'on
vient de spécifier.
Si $n \leq 2$, le groupe $F(\kappabarre)$ est trivial ou isomorphe à~$\Z/2$ et
il n'y a rien à prouver. Supposons donc que $n>2$.
Il reste seulement à établir que $C_1+C_i=C_{i+1}$ pour tout $i \in \{\uplet{0}{n-1}\}$,
la somme étant calculée dans $F(\kappabarre)$.
Compte tenu que la fibre spéciale de~$E^\star$ est réduite et que~$R$ est strictement
hensélien, il existe un point $a \in E(K)$ qui se spécialise sur~$C_1$.
D'après la propriété universelle caractérisant~$E^\star$,
l'automorphisme de~$E$ de translation par~$a$ s'étend en un automorphisme~$\tau$ de~$E^\star$.
La restriction de~$\tau$ à $\Erond \subset E^\star$ est l'automorphisme de~$\Erond$ déduit de
la translation par~$a$ par la propriété universelle du modèle de Néron; en particulier,
si l'on munit l'ensemble des composantes irréductibles de la fibre spéciale de $E^\star$
de la loi de groupe de $F(\kappabarre)$, l'application~$\tau$ induit un
endomorphisme du groupe $\{\uplet{C_0}{C_{n-1}}\}$.
Comme~$C_0$ est le neutre de ce groupe et que $\tau(C_0)=C_1$, il suffit, pour conclure,
de montrer que $\tau(C_i)=C_{i+1}$ pour tout $i \in \{\uplet{0}{n-1}\}$.
Utilisant à nouveau la propriété universelle du modèle de Néron, celle de $E^\star$
et la comparaison entre~$E^\star$ et~$\Erond$, on voit que
la multiplication par~$-1$ sur~$\Erond$ s'étend en un automorphisme
$\sigma \colon E^\star \rightarrow E^\star$. Celui-ci stabilise~$C_0$
puisqu'il stabilise~$\Erond^0$; étant
donné que~$C_0$ est $\kappa$\nobreakdash-isomorphe à~$\P^1_\kappa$,
que~$C_0$ est ensemblistement la réunion de~$\Erond^0_\kappa$ et de~$\{z_{n-1},z_0\}$,
que~$\sigma$ fixe un point de~$\Erond^0_\kappa$ (le neutre du groupe),
que la restriction de~$\sigma$ à~$\Erond^0_\kappa$ n'est pas l'identité (c'est
la multiplication par~$-1$) et que l'identité est le seul automorphisme de~$\P^1_\kappa$
fixant trois points, il est nécessaire que $\sigma(z_0)=z_{n-1}$.
Par ailleurs, on a $\tau \sigma = \sigma \tau^{-1}$ puisque cette égalité vaut
après restriction à la fibre générique de~$E^\star$. D'autre part, on a
$\tau(\{z_{n-1},z_0\})=\{z_0,z_1\}$ puisque $\tau(C_0)=C_1$ et que~$\tau$ stabilise
le lieu singulier de la fibre spéciale de~$E^\star$.  Si l'on avait $\tau(z_0)=z_0$,
les relations $\tau \sigma = \sigma \tau^{-1}$ et $\sigma(z_0)=z_{n-1}$
entraîneraient que $\tau(z_{n-1})=z_{n-1}$, d'où $z_{n-1} \in \{z_0,z_1\}$, ce qui
contredirait l'hypothèse selon laquelle $n>2$.
On a donc $\tau(z_0)=z_1$.   Le résultat voulu s'en déduit par des considérations
purement combinatoires; en effet,
$\tau$ agit sur le graphe
dont les sommets sont les~$z_i$ et dont les arêtes sont les~$C_i$,
or le seul automorphisme de ce graphe envoyant~$C_0$ sur~$C_1$ et~$z_0$ sur~$z_1$
est la rotation évidente.

Compte tenu que la propriété que deux composantes irréductibles de la fibre de~$E^\star$
au-dessus de~$\sbarre$ se rencontrent est préservée par l'action de $\Gal(\kappabarre/\kappa)$
et que~$m(0)$ est invariant sous $\Gal(\kappabarre/\kappa)$,
l'assertion que l'on vient de démontrer suffit à assurer que $\Gal(\kappabarre/\kappa)$
agit soit trivialement sur~$F(\kappabarre)$, soit par multiplication par~$-1$ à travers
le groupe de Galois d'une extension quadratique $\ell/\kappa$, nécessairement
unique à isomorphisme près.

Comme $p\neq 2$, la flèche
de spécialisation $\tors{2}E(K) \rightarrow \tors{2}\Erond(\kappa)$
est injective (cf.~\cite[7.3/3]{blr}).
Il s'ensuit
que $F(\kappabarre)$ est d'ordre pair si $\tors{2}E(K)\simeq \Z/2\times\Z/2$,
compte tenu que $F(\kappabarre)=\Erond(\kappabarre)/\Erond^0(\kappabarre)$ et
que $\Erond^0(\kappabarre)\simeq\kappabarre^\star$.
\end{demo}

\bigskip
Rappelons que si~$E$ est une courbe elliptique sur~$K$ à réduction multiplicative,
la fibre spéciale~$\Erond^0_\kappa$ de~$\Erond^0$ est un tore, et l'on dit que~$E$ est à réduction
multiplicative déployée lorsque ce tore est lui-même déployé.
Les tores de dimension~$1$ sur~$\kappa$ sont classifiés par le groupe $H^1(\kappa,\Z/2)$,
de sorte qu'un tel tore est soit déployé sur~$\kappa$, soit
déployé sur une extension quadratique,
unique à isomorphisme près.  Ceci fournit une seconde extension quadratique ou triviale
de~$\kappa$ naturellement associée à toute courbe elliptique à réduction multiplicative sur~$K$
(la première étant donnée par la proposition~\ref{annpropcyclpair}).
Il est faux que ces deux extensions soient isomorphes en général, mais ce n'est pas
non plus loin d'être vrai:

\bigskip
\begin{proposition}
\label{annpropcorfauxrat}
Soit~$E$ une courbe elliptique sur~$K$, à réduction de type~$I_n$ avec $n>2$.
L'extension quadratique ou triviale minimale $\ell/\kappa$ telle que
le $\ell$\nobreakdash-groupe $F \otimes_\kappa \ell$ soit constant (cf.~proposition~\ref{annpropcyclpair})
est isomorphe à l'extension quadratique ou triviale minimale de~$\kappa$
qui déploie le tore~$\Erond^0_\kappa$.
\end{proposition}

\bigskip
\begin{corollaire}
\label{anncorcstdep}
Soit~$E$ une courbe elliptique sur~$K$, à réduction multiplicative.
Le $\kappa$\nobreakdash-groupe $F$ est constant si et seulement si~$E$ est à réduction
multiplicative déployée ou à réduction de type~$I_1$ ou~$I_2$.
\end{corollaire}

\bigskip
\begin{demo}
Si le type de réduction de~$E$ est~$I_1$ ou~$I_2$, il est évident que~$F$ ne peut
être que constant.  Sinon, il suffit d'appliquer la proposition~\ref{annpropcorfauxrat}.
\end{demo}

\bigskip
Avant de prouver la proposition~\ref{annpropcorfauxrat},
rappelons, sans démonstration, une propriété très élémentaire des tores de dimension~$1$.

\bigskip
\begin{lemme}
\label{anntoreelem}
Soit $T$ un tore sur~$\kappa$, $T \hookrightarrow \P^1_\kappa$ une $\kappa$\nobreakdash-immersion
ouverte et $x \in \P^1_\kappa \setminus T$.
L'extension $\kappa(x)/\kappa$ est la plus petite extension
de~$\kappa$ qui déploie le tore~$T$.
\end{lemme}

\bigskip
\begin{demo}[ de la proposition~\ref{annpropcorfauxrat}]%
La preuve de la proposition~\ref{annpropcyclpair} montre
que pour toute sous-extension quadratique ou triviale $\ell/\kappa$
de $\kappabarre/\kappa$,
lorsque $n>2$,
le groupe $\Gal(\kappabarre/\ell)$ agit trivialement sur~$F(\kappabarre)$
si et seulement si les points singuliers de la fibre spéciale~$E^\star_\kappa$ de~$E^\star$
sont tous $\ell$\nobreakdash-rationnels, si et seulement si l'un d'entre eux est $\ell$\nobreakdash-rationnel.
En particulier, notant~$\ell/\kappa$ l'extension quadratique ou triviale
minimale telle que le $\ell$\nobreakdash-groupe $F \otimes_\kappa \ell$ soit constant,
tous les points singuliers de~$E^\star_\kappa$
ont un corps résiduel $\kappa$\nobreakdash-isomorphe à~$\ell$.
Par ailleurs, la composante irréductible de~$E^\star_\kappa$ contenant~$\Erond^0_\kappa$ est
$\kappa$\nobreakdash-isomorphe à~$\P^1_\kappa$ et le complémentaire de~$\Erond^0_\kappa$ dans
cette composante est formé de points singuliers de~$E^\star_\kappa$.
Le lemme~\ref{anntoreelem} permet donc de conclure.
\end{demo}

\bigskip
Le lemme suivant est extrêmement bien connu, mais une démonstration semble
plus facile à fournir qu'une référence.

\bigskip
\begin{lemme}
\label{annlemmeextreme}
Soit~$E$ une courbe elliptique sur~$K$, à réduction multiplicative.
Choisissons une équation de Weierstrass minimale de~$E$ et notons~$C$
la courbe sur~$\kappa$ définie par cette équation (c'est donc une cubique plane à point
double).
Alors la courbe~$E$ est à réduction multiplicative déployée si et seulement si
les pentes des directions tangentes à~$C$ en son point singulier sont $\kappa$\nobreakdash-rationnelles.
\end{lemme}

\bigskip
\begin{demo}
Comme l'équation de Weierstrass choisie est minimale, la $\kappa$\nobreakdash-variété $\Erond^0_\kappa$
est isomorphe à l'ouvert de lissité~$C^0$ de~$C$
(cf.~\cite[Ch.~IV, §9, Cor.~9.1]{silvaec2}).  D'autre part, si $\pi \colon \Ctilde \rightarrow C$
désigne la normalisation de~$C$, le complémentaire de $\pi^{-1}(C^0)$ dans~$\Ctilde$
est un $\kappa$\nobreakdash-schéma fini de degré~$2$, déployé si
et seulement si les pentes des directions
tangentes à~$C$ en son point singulier sont $\kappa$\nobreakdash-rationnelles.
Le lemme~\ref{anntoreelem} permet de conclure, compte tenu que~$\Ctilde$
est $\kappa$\nobreakdash-isomorphe à~$\P^1_\kappa$.
\end{demo}

\bigskip
Comparons maintenant les groupes de composantes connexes des fibres spéciales
des modèles de Néron de deux courbes elliptiques $2$\nobreakdash-isogènes
dans le cas de réduction multiplicative.

\bigskip
\begin{proposition}
\label{annpropppp}
Soient~$E'$ et~$E''$ des courbes elliptiques sur~$K$, à réduction multiplicative,
et $\phi' \colon E' \rightarrow E''$ une isogénie de degré~$2$.
Notons $P'$ le $K$\nobreakdash-point non nul de $\Ker(\phi')$
et $P''$ le $K$\nobreakdash-point non nul de $\Ker(\phi'')$, où~$\phi''$ est l'isogénie
duale de~$\phi'$. Les conditions suivantes sont équivalentes:
\begin{enumerate}
\item Le point~$P'$ se spécialise dans $\Erond'^0(\kappa)$.
\item Le point~$P''$ ne se spécialise pas dans $\Erond''^0(\kappa)$.
\item Le morphisme $F' \rightarrow F''$ induit par~$\phi'$ est injectif et a pour
conoyau $\Z/2$.
\item Le morphisme $F'' \rightarrow F'$ induit par~$\phi''$ est surjectif et a
pour noyau $\Z/2$.
\end{enumerate}
\end{proposition}

\bigskip
\begin{demo}
Comme $p \neq 2$, la courbe elliptique~$E'$ admet une équation de Weierstrass
minimale de la forme $y^2=(x-c)(x^2-d)$ avec $c,d\in R$, telle que
le point~$P'$ ait pour coordonnées $(x,y)=(c,0)$.
Son discriminant est $\Delta'=16d(c^2-d)^2$.  L'hypothèse de réduction multiplicative
se traduit par la condition que l'un de~$d$ et de $c^2-d$
est inversible dans~$R$ et que l'autre ne l'est pas. Il en résulte
notamment que l'équation de Weierstrass
$y^2=(x+2c)(x^2-4(c^2-d))$ est minimale.  Celle-ci définit une courbe elliptique isomorphe
à~$E''$,
son discriminant vaut $\Delta''=4096 d^2 (c^2-d)$ et le point~$P''$ a pour
coordonnées $(x,y)=(-2c,0)$.  Comme les équations de Weierstrass considérées sont
minimales, les ouverts de lissité des $R$\nobreakdash-schémas projectifs $W'$ et~$W''$ qu'elles définissent
sont respectivement isomorphes à~$\Erond'^0$ et~$\Erond''^0$;
un point rationnel de~$E'$ (resp.~de~$E''$) se spécialise donc
sur $\Erond'^0$ (resp.~$\Erond''^0$) si et seulement si la section de~$W'$ (resp.~$W''$)
qu'il définit ne rencontre pas le point singulier de la fibre spéciale de~$W'$ (resp.~$W''$).
En particulier, la condition~1 (resp.~la condition~2) équivaut à ce que $d\not\in R^\star$
(resp.~$c^2-d \in R^\star$), d'où l'équivalence entre~1 et~2.

Notons~$v$ la valuation normalisée de~$K$ associée à~$R$.  Les groupes $F'(\kappabarre)$
et $F''(\kappabarre)$ sont cycliques, d'ordres respectifs $v(\Delta')$ et $v(\Delta'')$
(cf.~\cite[Step~2, p.~366]{silvaec2}).  On voit sur les expressions
de~$\Delta'$ et~$\Delta''$ que
$d \not\in R^\star$ si et seulement si $v(\Delta'')=2v(\Delta')$,
si et seulement si $v(\Delta') \neq 2v(\Delta'')$.
L'observation suivante permet donc de conclure:
si~$A$ et~$B$ sont deux groupes cycliques tels
que $\Card(B)=2\Card(A)$ et si $u \colon A \rightarrow B$
et $v \colon B \rightarrow A$ sont des morphismes de groupes vérifiant $uv=2$ et $vu=2$,
on a nécessairement $\Ker(u)=0$, $\Coker(u)=\Z/2$, $\Ker(v)=\Z/2$ et $\Coker(v)=0$.
\end{demo}

\bigskip
\begin{corollaire}
\label{anncorcorfaux}
Soient~$E'$ et~$E''$ des courbes elliptiques sur~$K$, à réduction multiplicative,
et $\phi' \colon E' \rightarrow E''$ une isogénie de degré~$2$.  Si $E'$
n'est pas à réduction de type $I_1$ ou~$I_2$, alors l'extension quadratique ou triviale
minimale $\ell/\kappa$ qui déploie le tore $\Erond'^0_\kappa$ est aussi la plus petite
extension de~$\kappa$ telle que les deux $\ell$\nobreakdash-groupes $F'\otimes_\kappa\ell$
et $F''\otimes_\kappa\ell$ soient constants.
\end{corollaire}

\bigskip
\begin{demo}
Vu la proposition~\ref{annpropcorfauxrat}, il suffit de vérifier
que pour toute extension $\ell/\kappa$, si le $\ell$\nobreakdash-groupe $F' \otimes_\kappa\ell$
est constant, il en va de même de $F'' \otimes_\kappa\ell$.
D'après la proposition~\ref{annpropppp}, l'un des deux $\kappa$\nobreakdash-groupes~$F'$ et~$F''$ est isomorphe
à un sous-groupe de l'autre.  Si $F'' \hookrightarrow F'$, l'assertion est évidente.
Supposons donc que $F' \hookrightarrow F''$ et que $F' \otimes_\kappa \ell$ soit
constant.
Le groupe $F'(\ell)$ est \mbox{d'ordre~$>2$}, compte tenu de l'hypothèse sur~$E'$;
il en résulte que $F''(\ell)$ est lui aussi d'ordre~$>2$.
La proposition~\ref{annpropcyclpair} montre alors
que le $\ell$\nobreakdash-groupe
$F'' \otimes_\kappa \ell$ est nécessairement constant.
\end{demo}

\bigskip
\begin{corollaire}
\label{anncorsingrat}
Mêmes notations que dans la proposition~\ref{annpropppp}.
Soient $E'^\star$ un modèle propre et régulier minimal de~$E'$ au-dessus de~$R$
et $x \in E'^\star$ un point en lequel $E'^\star \rightarrow \Spec(R)$ n'est pas lisse.
Si~$E'$ est à réduction de type~$I_1$ ou~$I_2$ et que~$P'$ se spécialise
dans $\Erond'^0(\kappa)$, alors l'extension quadratique ou triviale minimale $\ell/\kappa$
telle que le $\ell$\nobreakdash-groupe $F''\otimes_\kappa\ell$ soit constant
est $\kappa$\nobreakdash-isomorphe à~$\kappa(x)$.
\end{corollaire}

\bigskip
\begin{demo}
Notons~$n$ l'ordre du groupe $F'(\kappabarre)$.
D'après la proposition~\ref{annpropppp},
les courbes elliptiques~$E'$ et~$E''$ sont respectivement à réduction de type~$I_n$ et~$I_{2n}$.
Si $n=1$, alors $F''=\Z/2$ et $\ell=\kappa$, et d'autre
part la fibre spéciale de $E'^\star$ contient un unique point singulier, qui est
rationnel; la conclusion du corollaire est donc satisfaite.
Supposons maintenant que $n=2$.
L'extension $\ell/\kappa$ est alors la plus petite extension qui
déploie le tore $\Erond''^0_\kappa$
(cf.~proposition~\ref{annpropcorfauxrat}).
Comme les tores $\Erond'^0_\kappa$ et $\Erond''^0_\kappa$
sont isomorphes (deux $\kappa$\nobreakdash-tores isogènes
de dimension~$1$ étant nécessairement isomorphes), le lemme suivant permet de conclure.
\end{demo}

\bigskip
\begin{lemme}
\label{annlemmetoresing}
Soient~$E$ une courbe elliptique sur~$K$, à réduction de type~$I_2$.
Soient $E^\star$ un
modèle propre et régulier minimal de~$E$ au-dessus de~$R$ et $x \in E^\star$ un point
en lequel $E^\star \rightarrow \Spec(R)$ n'est pas lisse.
L'extension quadratique ou triviale minimale de~$\kappa$ qui déploie le tore $\Erond^0_\kappa$
est $\kappa$\nobreakdash-isomorphe à~$\kappa(x)$.
\end{lemme}

\bigskip
\begin{demo}
Vu la structure de la fibre spéciale de~$E^\star$, et compte tenu que $\Erond$
est isomorphe à l'ouvert de lissité de $E^\star$ sur~$R$
(cf.~\cite[Ch.~IV, §6, Th.~6.1]{silvaec2}),
il suffit d'appliquer le lemme~\ref{anntoreelem}.
\end{demo}

\bigskip
Voici enfin quelques propriétés spécifiques aux courbes elliptiques
dont tous les points d'ordre~$2$ sont rationnels.

\bigskip
\begin{proposition}
\label{annpropekshsurdeux}
Soit~$E$ une courbe elliptique sur~$K$ telle
que $\tors{2}E(K) \simeq \Z/2 \times \Z/2$.
Si~$R$ est hensélien
(resp.~strictement hensélien), la flèche de spécialisation $E(K)/2 \rightarrow F(\kappa)/2$
est surjective (resp.~bijective).
\end{proposition}

\bigskip
\begin{demo}
Considérons la suite exacte
\begin{equation}
\label{annse}
\xymatrix{
0 \ar[r] & \Erond^0 \ar[r] & \Erond \ar[r] & i_\star F \ar[r] & 0
}
\end{equation}
de faisceaux étales sur~$\Spec(R)$, où $i \colon \Spec(\kappa) \rightarrow \Spec(R)$
désigne l'immersion fermée canonique. Si~$R$ est strictement hensélien,
cette suite reste exacte par passage aux sections globales, d'où il résulte (compte
tenu que $E(K)=\Erond(R)$) que le noyau de la flèche de spécialisation
$E(K)/2 \rightarrow F(\kappa)/2$ est un quotient de $\Erond^0(R)/2$; or
ce dernier groupe est nul (cf.~lemme~\ref{annlemmesurje0}).

Il reste à établir la surjectivité de la flèche $E(K)/2 \rightarrow F(\kappa)/2$ lorsque~$R$
est hensélien.  Nous allons en fait montrer que
$s \colon E(K) \rightarrow F(\kappa)$ est surjective.  On peut évidemment supposer
que~$E$ a mauvaise réduction.
Comme $p\neq 2$, la flèche
de spécialisation $\tors{2}E(K) \rightarrow \tors{2}\Erond(\kappa)$
est injective (cf.~\cite[7.3/3]{blr}).  D'autre part, le groupe $\tors{2}\Erond^0(\kappa)$
est d'ordre~$\leq 2$ puisque $\Erond^0(\kappabarre^\star)$
est isomorphe à~$\kappabarre$ ou à~$\kappabarre^\star$.
Compte tenu que $\tors{2}E(K)$ est d'ordre~$4$, il en résulte
que~$s$ est surjective si~$F(\kappa)$ est d'ordre~$\leq 2$.
Supposons maintenant $F(\kappa)$ d'ordre~$>2$.
Si~$E$ est à réduction multiplicative,
la proposition~\ref{annpropcyclpair} montre que le $\kappa$\nobreakdash-groupe~$F$ est
constant; comme il est constant et d'ordre~$>2$, la courbe elliptique~$E$ est alors à réduction
multiplicative déployée (cf.~corollaire~\ref{anncorcstdep}).
Ainsi le groupe $\Erond^0_\kappa$ est-il nécessairement isomorphe à~$\Gm$ ou à~$\Ga$.
Dans les deux cas, on a $H^1(\kappa,\Erond^0_\kappa)=0$ (théorème de Hilbert~90 et
sa version additive, cf.~\cite[Ch.~X, §1]{serrecl}), d'où
l'on déduit, puisque~$R$ est hensélien, que $H^1(R, \Erond^0)=0$
(par exemple à l'aide de la représentabilité des éléments
de ce groupe par des torseurs, cf.~\cite[Theorem~4.3]{milneec}).
L'exactitude de~(\ref{annse}) est donc préservée par passage aux sections
globales, ce qui prouve le résultat voulu.
\end{demo}

\bigskip
\begin{remarque}
Voici un exemple montrant que la conclusion de la proposition~\ref{annpropekshsurdeux}
peut être en défaut si l'on suppose seulement que $\tors{2}E(K) \neq 0$.
Soit~$k$ un corps de caractéristique~$\neq 2$ sur lequel il existe une conique
sans point rationnel, disons $ax^2+by^2=1$ avec $a,b \in k^\star$.
Notons~$C$ cette conique et
considérons la courbe elliptique~$E$ sur $k((t))$ définie
par l'équation de Weierstrass $y^2=(x+a)(x^2-bt^2)$.
Elle est à réduction de type~$I_2$, possède un point rationnel d'ordre~$2$, et l'on
vérifie sans peine (par un calcul d'éclatement) que la composante connexe non neutre
de la fibre spéciale de son modèle de Néron est $k$\nobreakdash-birationnelle à~$C$, de sorte
qu'elle ne possède pas de point rationnel.  La flèche de spécialisation
$E(K) \rightarrow F(\kappa)$ est donc nulle, alors que $F(\kappa)=\Z/2$.
\end{remarque}

\bigskip
\begin{lemme}
\label{annlemmeuniquespec}
Soit~$E$ une courbe elliptique sur~$K$, à réduction multiplicative et telle
que $\tors{2}E(K) \simeq \Z/2 \times \Z/2$. Il existe un unique
point de $\tors{2}E(K) \setminus \{0\}$ qui se spécialise dans $\Erond^0(\kappa)$.
\end{lemme}

\bigskip
\begin{demo}
Comme la flèche de spécialisation $\tors{2}E(K) \rightarrow \tors{2}\Erond(\kappa)$
est injective
et que $\tors{2}\Erond^0(\kappa)=\Z/2$ (puisque $\Erond^0$ est un tore de dimension~$1$),
l'unicité est claire.  Quant à l'existence, il suffit de remarquer
que la composée $\tors{2}E(K) \rightarrow \tors{2}\Erond(\kappa) \rightarrow \tors{2}F(\kappa)$
ne peut être injective, le groupe~$\tors{2}F(\kappa)$ étant d'ordre~$2$
(cf.~proposition~\ref{annpropcyclpair}).
\end{demo}

\bigskip
\begin{proposition}
\label{annpropcaracteksh}
Supposons~$R$ strictement hensélien.
Soit~$E$ une courbe elliptique sur~$K$, à réduction multiplicative et telle
que \mbox{$\tors{2}E(K) \simeq \Z/2 \times \Z/2$}.
Le sous-groupe $E(K)/2 \hookrightarrow H^1(K,\tors{2}E)$ est égal à l'image
du morphisme $$H^1(K,\Z/2) \longrightarrow H^1(K,\tors{2}E)$$ induit par l'inclusion
$\Z/2 \subset \tors{2}E$ de l'unique point de $\tors{2}E(K)\setminus\{0\}$
qui se spécialise dans $\Erond^0(\kappa)$.
\end{proposition}

\bigskip
\begin{demo}
Notons $\phi \colon E \rightarrow E''$ le quotient de~$E$ par l'unique
point~$P$ de $\tors{2}E(K)\setminus\{0\}$ qui se spécialise dans $\Erond^0(\kappa)$
et $\phi'' \colon E'' \rightarrow E$ l'isogénie duale.
Il résulte de la proposition~\ref{annpropppp} et du lemme~\ref{annlemmeisogekmsh}
que l'application $E''(K) \rightarrow E(K)$ induite par~$\phi''$ est surjective.
Le diagramme commutatif
$$
\xymatrix{
0 \ar[r] & \tors{2}E \ar[d] \ar[r] & E \ar[r]^2 \ar[d]^\phi & E \ar[r] \ar@{=}[d] & 0 \\
0 \ar[r] & \Z/2 \ar[r] & E'' \ar[r]^{\phi''} & E \ar[r] & 0 \rlap{\text{,}}
}
$$
dont les lignes sont exactes, permet d'en déduire la nullité de la composée
$$E(K)/2 \longhookrightarrow H^1(K,\tors{2}E) \longrightarrow H^1(K,\Z/2)\rlap{\text{,}}$$
où la seconde
flèche est induite par la flèche $\tors{2}E \rightarrow \Z/2$ de quotient par~$P$,
d'où le résultat.
\end{demo}

\bigskip
\begin{remarque}
On se gardera de croire que la proposition~\ref{annpropcaracteksh} affirme
quoi que ce soit au sujet
des classes dans $E(K)/2$ des points d'ordre~$2$ de~$E$.
Ces classes sont toutes nulles dès que le type de réduction de~$E$ est~$I_n$
avec $n>2$ (sous l'hypothèse que~$R$ est strictement hensélien).
\end{remarque}

\backmatter
\let\oldbibliography\thebibliography
\renewcommand{\thebibliography}[1]{%
  \oldbibliography{#1}%
  \setlength{\itemsep}{1.55pt}%
}
\renewcommand{\refname}{Bibliographie}
\bibliographystyle{custbibstyle}
\bibliography{book}
\printnotation
\def\indexname{Index terminologique}
\onecolindex
\printindex
\end{document}